\documentclass[twoside,a4paper,11pt]{amsbook}

\usepackage[T1]{fontenc}
\usepackage{mathpazo} 

\usepackage[utf8]{inputenc}
\usepackage{fullpage}
\usepackage{amscd}
\usepackage{amsmath}
\usepackage{amssymb}
\usepackage{amsthm}
\usepackage{color}
\usepackage{graphicx}
\usepackage{hyperref}
\usepackage[noabbrev,capitalize]{cleveref}
\usepackage{mathrsfs}
\usepackage{mathtools}
\usepackage{pict2e}
\usepackage{stackrel}
\usepackage[T1]{fontenc}
\usepackage{tikz}
\usepackage{tikz-cd}
\usetikzlibrary{decorations.pathmorphing,decorations.markings,arrows,calc,shapes.geometric,arrows.meta,positioning}
\usepackage[utf8]{inputenc}
\usepackage{wasysym}
\usepackage[all]{xy}
\usepackage{xcolor}
\usepackage{xy}
\usepackage{ipaex-type1}
\usepackage{epigraph}
\usepackage{tabularx}

\hypersetup{
    colorlinks,
    linkcolor={red!50!black},
    citecolor={blue!50!black},
    urlcolor={blue!80!black}
}

\usepackage{palatino}
\usepackage{mathpazo} 
\usepackage{eulervm}
\allowdisplaybreaks

\newtheorem{lemma}{Lemma}[section]
\newtheorem{theorem}[lemma]{Theorem}
\newtheorem{Corollary}[lemma]{Corollary}
\newtheorem{Proposition}[lemma]{Proposition}

\newtheorem*{Notation}{Notation}

\newtheorem{theoremintro}{Theorem}

\theoremstyle{definition}
\newtheorem{Definition}[lemma]{Definition}
\newtheorem{Remark}[lemma]{\sc Remark}
\newtheorem{Counterexample}[lemma]{\sc Counter-Example}
\newtheorem{Example}[lemma]{\sc Example}

\def\colim{\mathop{\mathrm{colim}}}

\newcommand{\palg}{\PP\textsf{-}\mathsf{alg}}

\newcommand{\upartialop}{\mathsf{upOp}}
\newcommand{\partialop}{\mathsf{pOp}}
\newcommand{\upartialcoop}{\mathsf{upCoop}}
\newcommand{\partialcoop}{\mathsf{pCoop}}

\newcommand{\smod}{\mathbb{S}\textsf{-}\mathsf{mod}}
\newcommand{\kmod}{\mathbb{K}\textsf{-}\mathsf{mod}}

\newcommand{\ac}{\scriptstyle \textrm{!`}}

\newcommand{\qi}{\xrightarrow{ \,\smash{\raisebox{-0.65ex}{\ensuremath{\scriptstyle\sim}}}\,}}
\newcommand{\lqi}{\xleftarrow{ \,\smash{\raisebox{-0.65ex}{\ensuremath{\scriptstyle\sim}}}\,}}

\newcommand{\Cobar}{\Omega}

\newcommand{\C}{\mathcal{C}}
\newcommand{\I}{\mathcal{I}}

\newcommand{\kk}{\mathbb{K}}

\newcommand{\PP}{\mathcal{P}}
\newcommand{\Q}{\mathcal{Q}}
\newcommand{\D}{\mathcal{D}}
\newcommand{\G}{\mathcal{G}}
\newcommand{\ucom}{u\mathcal{C}om}
\newcommand{\ucomd}{u\mathcal{C}om^{*}}

\newcommand{\Assc}{c\mathcal{A}ss}
\newcommand{\Liec}{c\mathcal{L}ie}
\newcommand{\cP}{(\PP,d_\PP, \Theta_\PP)}

\input cyracc.def
\font\tencyr=wncysc10
\def\cyr{\tencyr\cyracc}
\def\diracComb{\mbox{\cyr SH}}

\begin{document}
\begin{titlepage}
  \begin{center}
    \large
    UNIVERSITÉ PARIS XIII - SORBONNE PARIS NORD 

    École Doctorale Sciences, Technologies, Santé Galilée 

    \vspace{1cm}
    \huge
    \hrulefill
    \vspace{0.5cm}

    \textbf{Théorie de l'intégration des algèbres de Lie à homotopie près absolues courbées.}  

    \vspace{0.5cm}
    \Large
    \textbf{The integration theory of curved absolute homotopy Lie algebras.} 

    \hrulefill

    \vspace{1.5cm}

    \large
    THÈSE DE DOCTORAT\\
    présentée par\\
    \smallskip
    \Large
    \textbf{Víctor Roca i Lucio}\\ 
    \smallskip
    \large
    Laboratoire Analyse, Géométrie et Applications\\ 
    \bigskip
    pour l'obtention du grade de\\
    DOCTEUR EN MATHEMATIQUES 

    \vspace{1cm}

    \bigskip
    \normalsize
    soutenue le 4 juillet 2022 devant le jury d'examen composé de : 
    \bigskip

    \begin{tabularx}{\textwidth}{X @{} l}
      \textbf{Campos Ricardo}, CNRS--Université de Montpellier \dotfill  Examinateur\\
      \textbf{Cirici Joana}, Université de Barcelone \dotfill  Examinatrice\\
      \textbf{Ginot Grégory}, Université Paris XIII \dotfill  Examinateur\\
      \textbf{Horel Geoffroy}, Université Paris XIII \dotfill  Examinateur\\
      \textbf{Moerdijk Ieke}, Université d'Utrecht \dotfill  Président du jury\\ 
      \textbf{Vallette Bruno}, Université Paris XIII \dotfill  Directeur de thèse\\   
      \\
      \textbf{Moerdijk Ieke}, Université d'Utrecht \dotfill  Rapporteur\\
      \textbf{Willwacher Thomas}, Ecole Polytechnique Fédérale de Zurich \dotfill Rapporteur\\
     \end{tabularx}  
  \end{center}
 
\end{titlepage}

\newpage

\phantom{Page} 
\setcounter{page}{2}

\newpage

\epigraph{À Marie.}{}

\newpage

\tableofcontents

\newpage

\phantom{Page} 

\newpage

\section*{Remerciements--Agradecimientos.}
Les mauvaises langues disent que les remerciements d'une thèse sont l'endroit où le jeune mathématicien assouvi ses pulsions littéraires en devenant, l'espace de quelques lignes, un écrivain médiocre. Ce sera malgré tout un plaisir de les entendre jaser, après être fraichement sorti de cette sombre grotte dans laquelle la peste nous a tous terré. J'en suis sorti avec un désir de vengeance, tel le comte de Monte--Cristo, et non content d'infliger ce pavé au monde, je compte aussi lui infliger ces remerciements. Lisez donc, à votre guise, soit des remerciements sincères, soit une liste détaillée des complices sans qui il aurait sans doute été impossible de perpétrer tout ceci.

\medskip

En primer lugar, debo dar las gracias a mi madre. Sin su apoyo constante y su entusiasmo particular, sin la libertad que me dio desde pequeño para que me dedicase a lo que quisiese, en fin, sin todos sus sacrificios, jamás habría llegado tan lejos. También es menester agradecer el apoyo de mi familia, a pesar de todos nuestros \textit{forcejeos}. 

\medskip 

Le premier remerciement d'ordre mathématique (mais évidemment personnel aussi), je le dois sans nul doute à mon directeur, Bruno Vallette. Au-delà des discussions mathématiques toujours intéressantes, j'ai beaucoup apprécié sa profonde humanité, son sens de l'humour (parfois douteux), sa franchise et sa patience. Surtout, je tiens sincèrement à le remercier pour m'avoir appris à \textit{faire} des mathématiques, et m'avoir initié à ce fabuleux métier. 

\medskip

Je tiens aussi à remercier Geoffroy Horel et Brice Le Grignou. Discuter avec eux m'a beaucoup appris, je suis grandement reconnaissant envers leur générosité. Je remercie aussi tous les membres de l'équipe de Topologie Algébrique. Les multiples groupes de travail, séminaires et autres activités, toujours dans une ambiance chaleureuse et conviviale, malgré les effluves d'amiante, m'ont énormément aidé à progresser. Ce groupe a su garder une activité remarquable pendant les temps troubles qui ont accompagné une bonne partie de cette thèse.
Bien sûr, merci aux doctorants du LAGA, et tout particulièrement ceux de mon équipe (Hugo, Jean-Michel, Guillaume, Coline, et les autres), pour nombre de discussions surtout absurdes, parfois sensées. Sans eux, les pauses café auraient été bien plus ternes. Je remercie Roberta pour m'avoir faire découvrir la poésie des mathématiques, ainsi que tous ses amis dévoués. Il faut aussi mentionner la grande hospitalité de Jussieu pendant toute cette thèse, que ce soit au couloir des doctorants à l'IMJ ou au LPSM. 

\medskip

Une pensée pour tous les autres qui ont participé avec bienveillance à ma formation de mathématicien. Depuis le lycée, avec les cours de mathématiques avancées et ceux de philosophie, dont je garde un beau souvenir, jusqu'au master, et notamment Benoît Stroh avec qui j'ai découvert en TER les joies de l'algèbre pure. 

\medskip

Je remercie également tous les membres du jury et les rapporteurs, qui, pour une raison parfaitement obscure à mes yeux, ont acceptés de s'infliger la lecture de ces 188 pages. 

\medskip

Enfin, je ne peux ignorer les autres amis qui m'ont accompagné pendant toutes ces années. Un grand merci à Emilien, Diego, Bianca, Anthony, et Brivael sans qui ce long voyage aurait été infiniment plus solitaire et ennuyeux. Penser au chemin parcouru ensemble donne parfois le vertige, mais la route a été belle. Sans oublier Thibault, Dylan, Pierre, Thomas et leurs invectives trop familières ou familiales. Un mot pour Adrien, ce fin esthète avec qui il fait bon de boire une bière ou d'avoir des crampes. La liste complète serait trop longue, mais merci à tous ceux qui ont rempli de charme cette vie parisienne depuis ses débuts. Hâte de découvrir avec vous de quoi l'avenir est fait. Finalement, un silence pudique pour Marie : tu ne sais que trop bien tout ce que je te dois.

\section*{Prologue.}
Ce prologue se veut un dictionnaire des principaux concepts qui apparaissent dans cette thèse, à l'usage des non-initiés et des curieux. Nous nous efforcerons d'ébaucher ces concepts dans un but clair de concision.

\medskip

On peut concevoir la mathématique comme l'étude des objets mathématiques. En dehors de la sempiternelle controverse sur leur ontologie, un fait marquant est que les objets mathématiques sont le plus souvent munis de \textit{structures}. Par exemple, l'objet le plus élémentaire de la mathématique, l'ensemble des nombres entiers $\mathbb{Z} = \{\cdots, -1, 0, 1, \cdots \}~$, est muni de l'\textit{addition} et de la \textit{multiplication}. Multiplier deux nombres n'est pas une opération quelconque : elle satisfait des propriétés bien précises. Notamment, elle est linéaire par rapport à l'addition, elle est \textit{associative}, \textit{commutative} et \textit{unitaire}. Considérer de façon générale un ensemble quelconque avec des opérations qui ressembleraient à une addition et une multiplications, sans s'attarder sur les exemples \textit{concrets} tels que $\mathbb{Z}$, est ce qui donne lieu à la notion d'\textit{anneau}. Cette approche peut être vue comme le point de départ de l'algèbre moderne, un des sujets les plus féconds du dernier siècle. 

\medskip

Les mathématiciens ont par la suite décidé de séparer ces opérations des ensembles sur lesquelles elles habitent. En effet, on peut définir un objet appelé \textit{opérade}, constitué uniquement des opérations abstraites et des relations qu'elles satisfont. Leur force réside en ce que pour un ensemble, avoir une certaine structure algébrique revient à recevoir une certaine "participation" de la part de l'opérade associée à cette structure. Ce sont des objets universels qui encodent la faune algébrique dans toute sa diversité. Un jeu de billard commun en algèbre est de produire des constructions universelles qui associent à une structure algébrique une autre, sans faire de choix. Par exemple, à partir de toute algèbre de Lie on peut construire une algèbre associative, que l'on appelle son algèbre enveloppante universelle. Les opérades permettent d'éclairer ces procédés : établir des constructions universelles entre différents types d'algèbres revient à établir une relation entre les opérades qui les encodent. Ces objets fournissent donc de puissants outils pour étudier ce qu'on pourrait appeler l'\textit{algèbre universelle}, qui étudie les relations entre les différents types de structures algébriques.

\medskip

Un autre concept clef sera celui de \textit{cogèbre}. Multiplier est une opération qui prend deux éléments en entrée, disons $a$ et $b$, et vous donne un seul élément "$a$ fois $b$", que l'on dénote par $a.b$. Une cogèbre est la donnée d'un procédé dual : en partant d'un seul élément, vous voulez obtenir un ensemble de paires d'éléments. L'exemple le plus simple est encore donné par $\mathbb{Z}$. En effet, non seulement on est capable ici de vous multiplier des nombres, mais on sait aussi les \textit{diviser}. Concrètement, car il y a bien un moment où tout algébriste veut du sens, cela s'instancie de la façon suivante : 

\[
\begin{tikzcd}[column sep=4pc,row sep=0.5pc]
\Delta: \mathbb{Z} \arrow[r]
&\mathbb{Z} \times \mathbb{Z} \\
a \arrow[r,mapsto]
&\{(b,c) ~~ \text{tel que} ~~ b.c = a\}~.
\end{tikzcd}
\]

On associe à un nombre $a$ le sous-ensemble de $\mathbb{Z} \times \mathbb{Z}$ constitué de toutes les paires $(b,c)$ telles que $b.c =a$. Cette opération de \textit{décomposition} n'est pas anodyne non plus, elle satisfait des propriétés que l'on appelle la \textit{coassociativité} et la \textit{counitalité}, et qui peuvent très bien s'écrire en termes d'arbres en renversant la description arboricole de la relation d'associativité et d'unitalité. Nous invitons le lecteur à aller lire le très beau article de survol sur les cogèbres de T. Fox \cite{ekwcfox} pour approfondir ce sujet. 

\medskip 

Un dernier aspect fondamental du domaine est la notion d'\textit{homotopie}. Elle est née dans le giron de la topologie pour étudier les espaces de façon plus simple, car la relation "être les mêmes", appelée \textit{homéomorphisme}, s'est rapidement avérée trop rigide. Ainsi, on a considéré des espaces topologiques à déformation continue près. On appelle l'existence d'une telle deformation continue une \textit{homotopie}, et le but du jeu revient donc à savoir quand deux espaces sont \textit{homotopes}. Le changement de paradigme est assez conséquent : l'analogue de "l'égalité" pour les structures algébriques peut lui aussi être remplacé par des notions plus faibles d'équivalences. Cela donne lieu à un sujet que l'on peut appeler \textit{algèbre homotopique}. Via les \textit{invariants algébriques}, c'est-à-dire, via la construction d'une structure algébrique à partir de la donnée d'un espace topologique, on peut relier ces deux notions: l'homotopie des espaces et l'homotopie des algèbres. Un des but du domaine fut de trouver un invariant algébrique \textit{fidèle}, autrement dit, pour lequel deux on aurait que deux espaces sont homotopes si et seulement si leurs invariants algébriques sont homotopes. Cette question, en partie résolue, reste une de motivations de cette thèse.

\medskip

Comme vous le soupçonnez peut-être, la théorie des opérades rélève pleinement sa force quand on l'utilise pour faire de l'algèbre homotopique. Un de ses aspects fondamentaux est la dualité de Koszul au sens large, qui joue un rôle fondamental pour comprendre la théorie de l'homotopie des algèbres encodées par une opérade. A partir d'une opérade, on peut lui en associer une autre, sa duale de Koszul. Cette dualité s'instancie alors en une dualité entre certaines algèbres encodées par la première et certaines cogèbres encodées par la seconde. Ces cogèbres, munies elles aussi d'une relation d'equivalence propre, nous permettent d'y \textit{lire} l'information homotopique de nos algèbres de départ. Cette dualité a une place prédominante dans l'étude qui nous concerne.

\medskip

Nous avions promis de la concision, alors nous nous en tiendrons à cette ébauche pour l'instant. Cette thèse s'inscrit dans le vaste monde de l'algèbre homotopique et de la théorie des opérades, et vise à y apporter sa modeste contribution. On pourra y voir des résultats concernant cette dualité de Koszul, une extension de celle-ci dans le cas \textit{courbé}, et finalement la construction d'invariants algébriques pour les espaces topologiques en termes de d'\textit{algèbres de Lie à homotopie près courbées}. On peut comprendre ces algèbres comme une version \textit{homotopique} des algèbres de Lie classiques, munies de courbure. Nous saurions gré au lecteur taquin de ne nous demander davantage d'explications pour l'instant, et nous l'invitons à poursuivre son chemin jusqu'au chapitre suivant, où il trouvera les détails qui brillent ici par leur absence. 

\section*{Résumé--Abstract}
Le but principal de cette thèse est de développer la théorie de l'intégration des algèbres de Lie à homotopie près absolues courbées. Dans un premier chapitre, nous développons les outils du calcul opéradique nécessaire: nous codons les cogèbres non conilpotentes via des opérades et nous introduisons la notion duale d'algèbres sur une cooperade, dites "algèbres absolues".
Les premières sections de ce chapitre rappellent l'état de l'art et les trois dernières sont faites de résultats originaux. Dans le deuxième chapitre, nous développons la nouvelle théorie du calcul opéradique courbé, qui permet d'encoder les types d'algèbres courbés avec des operades courbées. Finalement, dans le troisièmes chapitre, nous développons la théorie de l'intégration des algèbres de Lie à homotopie près absolues courbées grâce aux outils introduits dans les deux chapitres précédents. Notre nouvelle approche de l'intégration nous permet d'obtenir plusieurs applications à la théorie de la déformation et à la théorie de l'homotopie rationnelle.

\medskip

The main goal of this thesis is to develop the integration theory of curved homotopy Lie algebras. In the first chapter, we develop the operadic calculus needed: we encode non-necessarily conilpotent coalgebras with operads and introduce their dual notion of an algebra over a cooperad, which we call "absolute algebras". The first sections of this chapter are a state of the art, and the last three sections are made of original results. In the second chapter, we develop the new theory of curved operadic calculus, which allows us to encode curved types of algebras with curved operads. Finally, in the third chapter, we develop the integration theory of curved absolute homotopy Lie algebras using the tools introduced in the first two chapters. This new approach allows us to obtain applications to deformation theory and rational homotopy theory.

\newpage

\phantom{Page} 

\newpage

\chapter{Absolute algebras, contramodules and duality squares.}\label{Chapter 1}

\epigraph{"Le silence éternel des décompositions infinies m'effraie."}{Blaise Pascal, précurseur de la notion de cogèbre.}

\section*{Introduction}
Let $A$ be a vector space. The data of an associative algebra structure amounts to the data of a "multiplication table" of elements of $A$: given two elements, it states which element in $A$ is their product. One can compile this "multiplication table" into a single morphism
\[
\left\{
\begin{tikzcd}[column sep=3pc,row sep=0pc]
\gamma_A: \displaystyle \bigoplus_{n \geq 1} A^{\otimes n} \arrow[r]
&A \\
a_1 \otimes \cdots \otimes a_n \arrow[r,mapsto]
&\gamma_A(a_1 ,\cdots,a_n)~,
\end{tikzcd}
\right.
\]
which assigns to any $n$-tuple $(a_1, \cdots, a_n)$ the value of their product $\gamma_A(a_1,\cdots,a_n)$. The conditions that an associative algebra structure has to satisfy are encoded by the fact that $\gamma_A$ defines a structure of an \textit{algebra over a monad}. But in this classical algebraic framework, infinite collections of elements do not have an assigned value by this "multiplication table". In order to given a value to expressions like $\sum_{n \geq 1} a_{1,n} \otimes \cdots \otimes a_{n,n}$, the classical approach used so far has been to add the data of a topology on $A$, so that the values assigned to any partial sum \textit{converge} to a well-defined element in $A$. This solution has several issues, specially when one tries to mix \textit{homotopical algebra} and \textit{topology}. For instance, neither the category of topological abelian groups nor the category of topological modules over a topological ring are abelian. This type of problems is the main motivation for the recent approaches of \cite{scholze} and \cite{barwick}.

\medskip

Absolute algebras are a new type of algebraic structures where infinite sums of operations have well-defined images \textit{by definition}. For instance, the data of an \textit{absolute algebra structure} on a vector space $A$ amounts to the data of a "transfinite multiplication table", where any series of elements in $A$ is assigned a value in $A$. This information compiles into a single structural map
\[
\left\{
\begin{tikzcd}[column sep=2.5pc,row sep=-0.5pc]
\gamma_A: \displaystyle \prod_{n \geq 1} A^{\otimes n} \arrow[r]
&A \\
\displaystyle \sum_{n \geq 1} a_{1,n} \otimes \cdots \otimes a_{n,n} \arrow[r,mapsto]
&\displaystyle \gamma_A \left(\sum_{n \geq 1} a_{1,n} \otimes \cdots \otimes a_{n,n} \right)~.
\end{tikzcd}
\right.
\]
The bright point is that this type of structure is defined as an algebra over a monad, hence the category of absolute associative algebras enjoys many desirable properties. For example, it is complete and cocomplete. Any absolute algebra structure gives an algebra structure in the classical sense by restricting the structural map to finite sums. There are many "absolute analogues" of standard types of algebraic structures: absolute Lie algebras, absolute $\mathcal{A}_\infty$-algebra, absolute $\mathcal{L}_\infty$-algebras, etc. This stems from the fact that this type of structures are defined as \textit{algebras over a cooperad}. 

\medskip

Recall that most of classical types of algebraic structures can be successfully encoded by an operad, see \cite{LodayVallette12}. Now, given an operad $\mathcal{P}$, we get a category of dg $\mathcal{P}$-algebras. The operadic calculus provides us with a powerful tool to study the homotopy category of these dg $\mathcal{P}$-algebras. The key element is the Koszul duality for operads, that allows one to construct a Koszul dual conilpotent cooperad $\mathcal{P}^{\hspace{1pt}\ac}$. This data gives a Bar-Cobar adjunction relative to the canonical twisting morphism $\kappa$ 
\[
\begin{tikzcd}[column sep=5pc,row sep=5pc]
\mathsf{dg}~\mathcal{P}\text{-}\mathsf{alg} \arrow[r,"\mathrm{B}"{name=B},shift left=1.1ex] 
&\mathsf{dg}~\mathcal{P}^{\hspace{1pt}\ac}\text{-}\mathsf{coalg} \arrow[l,"\Omega"{name=C},shift left=1.1ex] \arrow[phantom, from=C, to=B, , "\dashv" rotate=90]
\end{tikzcd}
\]
between the category of dg $\mathcal{P}$-algebras and the category of dg $\mathcal{P}^{\hspace{1pt}\ac}$-coalgebras. Contrary to popular belief, any coalgebra over a cooperad is, by definition, conilpotent. Hence we omit this adjective when possible. Using this adjunction, one can transfer the model structure where weak-equivalences are given by quasi-isomorphism onto the category of conilpotent $\mathcal{P}^{\hspace{1pt}\ac}$-coalgebras and obtain a Quillen equivalence. This approach, first developed by V. Hinich for dg Lie algebras in \cite{Hinich01} and by K. Lefèvre-Hasegawa in \cite{LefevreHasegawa03} for dg associative algebras, was extended to all types of dg $\mathcal{P}$-algebras in \cite{Vallette14}. Conilpotent coalgebraic types of structures are given by decomposition maps, which assign finite sums to any element. Let $C$ be a vector space. For example, the data of a non-counital \textit{conilpotent} coassociative coalgebra structure on $C$ is equivalent to the data of a "decomposition table", compiled into a map
\[
\left\{
\begin{tikzcd}[column sep=3pc,row sep=-0.5pc]
\Delta_C: C \arrow[r]
&\displaystyle \bigoplus_{n \geq 1} C^{\otimes n}\\
c \arrow[r,mapsto]
&\displaystyle \Delta_C(c) = \sum_{i \in \mathrm{I}} c_1^{(i)} \otimes \cdots \otimes c_n^{(i)}~,
\end{tikzcd}
\right.
\]
where $\mathrm{I}$ is a finite set. It is \textit{conilpotent} precisely because this sum is finite, and therefore $\Delta_C$ lands on the direct sum instead of the product. But most coassociative coalgebras that appear in nature are not conilpotent, e.g: $\kk$ with its diagonal. 

\medskip

In the same way as algebras are "the Koszul dual notion" to conilpotent coalgebras, since they both are "finitary types of structures", absolute types of algebras are the "Koszul dual notion" to non-conilpotent coalgebras. The infinite decompositions of elements in these non-conilpotent coalgebras are reflected in the "transfinite multiplication tables" of absolute algebras. The notion of an algebra over a cooperad was introduced in \cite{grignoulejay18} in order to study the homotopy theory of non necessarily conilpotent coalgebras. Since these types of structures are encoded as \textit{coalgebras over an operad}, the reasonable thing to do is to look at what an algebra over a cooperad looks like. Starting from a cooperad, one can construct a monad by considering a dual version of the Schur functor. Algebras over a cooperad are defined as algebras over its associated monad. Given an operad $\mathcal{P}$ and its Koszul dual cooperad $\mathcal{P}^{\hspace{1pt}\ac}$, the authors of \textit{loc.cit} construct a \textit{complete Bar-Cobar adjunction}
\[
\begin{tikzcd}[column sep=5pc,row sep=5pc]
\mathsf{dg}~\mathcal{P}\text{-}\mathsf{coalg} \arrow[r,"\widehat{\Omega}"{name=B},shift left=1.1ex] 
&\mathsf{dg}~\mathcal{P}^{\hspace{1pt}\ac}\text{-}\mathsf{alg}~, \arrow[l,"\widehat{\mathrm{B}}"{name=C},shift left=1.1ex] \arrow[phantom, from=C, to=B, , "\dashv" rotate=-90]
\end{tikzcd}
\]
and show that in some cases, the homotopy theory of these non necessarily conilpotent coalgebra can be recovered from the homotopy theory of the dual absolute algebras with a transferred model structure along this complete Bar-Cobar adjunction.

\medskip

The first sections of this chapter can be considered as a survey on the theory developed in \cite{grignoulejay18}. We introduce all the required operadic notions in order to explain their main results in Section \ref{Section: Complete Bar-Cobar}. The last three sections contain the main original results of this chapter, which are stated in \cite{lucio2022contra}. The first one is to relate the existing Bar-Cobar adjunction with the new complete Bar-Cobar adjunction introduced before. We show that there are duality functors that intertwine both of these adjunctions in a duality square of commuting adjunctions.

\begin{theoremintro}[Duality square, Theorem \ref{thm: magical square}]
There exists a square of adjunctions 
\[
\begin{tikzcd}[column sep=5pc,row sep=5pc]
\left(\mathsf{dg}~\mathcal{P}\text{-}\mathsf{alg}\right)^{\mathsf{op}} \arrow[r,"\mathrm{B}^{\mathsf{op}}"{name=B},shift left=1.1ex] \arrow[d,"(-)^\circ "{name=SD},shift left=1.1ex ]
&\left(\mathsf{dg}~\mathcal{P}^{\hspace{1pt} \mathsf{\ac}} \text{-}\mathsf{coalg}\right)^{\mathsf{op}} \arrow[d,"(-)^*"{name=LDC},shift left=1.1ex ] \arrow[l,"\Omega^{\mathsf{op}}"{name=C},,shift left=1.1ex]  \\
\mathsf{dg}~\mathcal{P}\text{-}\mathsf{coalg} \arrow[r,"\widehat{\Omega}"{name=CC},shift left=1.1ex]  \arrow[u,"(-)^*"{name=LD},shift left=1.1ex ]
&\mathsf{dg}~\mathcal{P}^{\hspace{1pt} \mathsf{\ac}} \text{-}\mathsf{alg}~, \arrow[l,"\widehat{\mathrm{B}}"{name=CB},shift left=1.1ex] \arrow[u,"(-)^\vee"{name=TD},shift left=1.1ex] \arrow[phantom, from=SD, to=LD, , "\dashv" rotate=0] \arrow[phantom, from=C, to=B, , "\dashv" rotate=-90]\arrow[phantom, from=TD, to=LDC, , "\dashv" rotate=0] \arrow[phantom, from=CC, to=CB, , "\dashv" rotate=-90]
\end{tikzcd}
\] 
which commutes in the following sense: right adjoints going from the top right corner to the bottom left corner are naturally isomorphic.
\end{theoremintro}

There the functor $(-)^\circ$ is a generalization of the Sweedler dual functor constructed in \cite{Sweedler69}. The above theorem admits a much more general formulation. It also holds for any curved twisting morphism, as it will be shown in Section \ref{Section: Curved duality square}. Furthermore, in the appropriate cases, this square becomes a square of Quillen adjunctions. This results plays an essential role in Chapter \ref{Chapter 3}, where we develop the integration theory of curved absolute $\mathcal{L}_\infty$-algebras.

\medskip

We then treat examples and develop the theory in those particular cases of interest. The first example is that of \textit{contramodules}, which appear as a particular case of this definition. Contramodules over coassociative coalgebras where first introduced by S. Eilenberg and J. C. Moore in \cite{eilenbergmoore65} but later somewhat forgotten until they were extensively studied by L. Positselski, see for instance \cite{positselski2021contramodules}. Since a coassociative coalgebra is a cooperad concentrated in arity one, we show that contramodules are a particular example of absolute algebras. Engulfing this theory provides us with illuminating examples and counterexamples that shed a light into what is to be expected of this new type of algebraic structures. After, we treat \textit{in extenso} the cases of dg absolute associative algebras and dg absolute Lie algebras. We show that \textit{nilpotent} associative algebras are a particular examples of absolute associative algebras, and describe various new constructions that can be performed on this category.

\medskip

Finally, we apply this new framework to Lie theory. In \cite{campos2020lie}, the authors proved the following theorem: two nilpotent Lie algebras are isomorphic as Lie algebras if and only if their universal enveloping algebras are isomorphic as associative algebras. For proving this statement, they compared the deformation complexes of $\mathcal{C}_\infty$ and $\mathcal{A}_\infty$-coalgebras. In this last section, we show how to reinterpret this result on deformation complexes in the more general context of absolute Lie algebras and their universal enveloping absolute algebras. By doing so, we were able to generalize their theorem as follows.

\begin{theoremintro}[Theorem \ref{thm: iso envelopantes Lie absolues} and Theorem \ref{thm: isos envelopantes absolues de L infinies}]
Let $\kk$ be a field of characteristic zero and let $\mathfrak{g}$ and $\mathfrak{h}$ be two complete graded absolute Lie algebras (resp. two complete minimal absolute $\mathcal{L}_\infty$-algebras). They are isomorphic as complete graded absolute Lie algebras (resp. as complete minimal absolute $\mathcal{L}_\infty$-algebras) if and only if their universal enveloping absolute algebras (resp. absolute $\mathcal{A}_\infty$-algebras) are isomorphic.
\end{theoremintro}

The first thing to mention on this result is that is follows from a more general statement concerning dg absolute Lie algebras and absolute $\mathcal{L}_\infty$-algebras. Indeed, we show that two of them are linked by a zig-zag a weak-equivalences if and only if their universal enveloping constructions are. Since these weak equivalences are in particular quasi-isomorphisms, the above result follows when the differential is zero. The second thing is that class of complete graded absolute Lie algebras includes nilpotent graded Lie algebras without degree restrictions. In this particular case, we compute explicitly the universal enveloping absolute algebra, which is given by the completed tensor algebra modulo the standard relation. Analogously, nilpotent $\mathcal{L}_\infty$-algebras in the sense of \cite{Getzler09} are examples of absolute $\mathcal{L}_\infty$-algebras. Therefore the above theorem hold for minimal nilpotent $\mathcal{L}_\infty$-algebras as well.

\subsection{Conventions.}
Let $\mathbb{K}$ be a ground field of characteristic $0$. The ground category is the symmetric monoidal category $(\mathsf{gr}\textsf{-}\mathsf{mod}, \otimes, \mathbb{K})$ of graded $\mathbb{K}$-modules, with the tensor product $\otimes$ of graded modules given by 
\[
(A\otimes B)_n \coloneqq \displaystyle \bigoplus_{p+q=n}~A_p \otimes~ B_q ~.
\]
The tensor is taken over the base field $\mathbb{K}$, which will be implicit from now on. The isomorphism $\tau_{A,B}: A \otimes B \longrightarrow B \otimes A$ is given by the Koszul sign rule $\tau (a \otimes b) \coloneqq (-1)^{|a|.|b|} \hspace{1pt} b\otimes a$ on homogeneous elements. We work with the \textit{homological} degree convention, that is, the degree of the (pre)-differentials considered is $-1$. The suspension of a graded module $V$ is denoted by $sV$, given by $(sV)_{p} \coloneqq V_{p-1}$~. 

\section{Monoidal definition of operads}
We recall the definition of an $\mathbb{S}$-module, of the composition product of $\mathbb{S}$-modules and the monoidal definition of an operad. For a comprehensive exposition, see \cite[Chapter 5]{LodayVallette12}. We omit the adjective \textit{graded}, which will be implicit from now on until Section \ref{Section: Complete Bar-Cobar}.

\begin{Definition}[$\mathbb{S}$-module]
A $\mathbb{S}$\textit{-module} $M \coloneqq \{ M(n) \}_{n \geq 0}$ is the data of a right $\mathbb{K}[\mathbb{S}_n]$-module for all $n \geq 0$. A \textit{morphism} of $\mathbb{S}$-module $f: M \longrightarrow N$ is given by a family of morphisms $\{f(n): M(n) \longrightarrow N(n)\}_{n \geq 0}$ of right $\mathbb{K}[\mathbb{S}_n]$-modules.
\end{Definition}

They form a category, denoted by $\mathbb{S}\textsf{-}\mathsf{mod}$. This category can be endowed with several monoidal products; we are going to focus on the composition product that defines operads. 

\begin{Definition}[Composition product of $\mathbb{S}$-module]
The composition product of two $\mathbb{S}$-modules $M$ and $N$, denoted by $M \circ N$, is defined as:  

\[
M \circ N (n) \coloneqq \bigoplus_{k\geq 0} M(k) \otimes_{\mathbb{S}_k} \left( \bigoplus_{i_1 + \cdots + i_k = n} \mathsf{Ind}_{\mathbb{S}_{i_1} \times \cdots \times \mathbb{S}_{i_k}}^{\mathbb{S}_n} (N(i_1) \otimes \cdots \otimes N(i_k))\right)~,
\]

where the second sum runs through all $k$-uples $(i_1, \cdots, i_k) \in \mathbb{N}^k$ such that $i_1 + \cdots + i_k = n$. The first sum over $k$ is finite if and only if $N(0) = 0$.
\end{Definition}

This endows with a monoidal structure the category of $\mathbb{S}$-modules, where the unit is given by $\I \coloneqq (0, \mathbb{K}.\mathrm{id}, 0, \cdots)$. We denote this monoidal category by $(\mathbb{S}\textsf{-}\mathsf{mod}, \circ, \I)$. 

\begin{Notation}
Since we already use $\circ$ for the composition product of $\mathbb{S}$-modules, we will sometimes denote the composition of two morphisms by $\cdot$ when the context might be ambiguous.
\end{Notation}

\begin{Definition}[Operad]\label{def: operad}
An \textit{operad} $(\PP, \gamma, \eta)$ is the data of a monoid in the monoidal category $(\mathbb{S}\textsf{-}\mathsf{mod}, \circ, \I)$.  An \textit{morphism} $f: (\PP, \gamma_\PP, \eta_\PP) \longrightarrow (\Q, \gamma_\Q, \eta_\Q)$ of operads is given by a morphism of $\mathbb{S}$-modules $f: \PP \longrightarrow \Q$ a that commutes with the monoid structure.
\end{Definition}

It is well known \cite[Proposition 5.3.1]{LodayVallette12} that the data of the composition map $\gamma: \PP \circ \PP \longrightarrow \PP$ is equivalent to the data of a family of linear composition maps:

\[
\gamma(i_1, \cdots, i_k): \PP(k) \otimes \PP(i_1) \otimes \cdots \otimes \PP(i_k) \longrightarrow \PP(n)~,
\]

satisfying conditions corresponding to the equivariance and to the associativity, which we call May's axioms. 

\begin{enumerate}
\medskip

\item Each map $\gamma(i_1, \cdots, i_k)$ is $\mathbb{S}_{i_1} \times \cdots \times \mathbb{S}_{i_k}$-equivariant map where the action on $\PP(n)$ is given by restricting the action of $\mathbb{S}_n$ to $\mathbb{S}_{i_1} \times \cdots \times \mathbb{S}_{i_k}$.

\medskip

\item The map:

\[
\bigoplus_{i_1 + \cdots + i_k = n} \gamma(i_1, \cdots, i_k) : \PP(k) \otimes \left( \bigoplus_{i_1 + \cdots + i_k = n} \PP(i_1) \otimes \cdots \otimes \PP(i_k) \right) \longrightarrow \PP(n)
\]

factors through the tensor product over $\mathbb{S}_k$. 

\medskip

\item For any set of indices $(j_{1,1}, \cdots, j_{1,i_1},j_{2,1}, \cdots, j_{2,i_2}, \cdots, j_{n,1}, \cdots, j_{n,i_n})$ , the following square of linear composition maps is commutative:

\[
\begin{tikzcd}[column sep=4pc,row sep=2pc]
\PP(n) \otimes \PP(r_1) \otimes \cdots \otimes \PP(r_n) \arrow[dddr,bend left = 35]
& \\
\PP(n) \otimes \PP(i_1) \otimes \PP(j_{1,1}) \otimes \cdots \otimes \PP(j_{1,i_1}) \otimes \PP(i_2) \otimes \PP(j_{2,1}) \otimes \cdots \otimes \PP(j_{n,i_n}) \arrow[u] \arrow[d,"\cong"]
& \\
\PP(n) \otimes \PP(i_1) \otimes \cdots \otimes \PP(i_n) \otimes \PP(j_{1,1}) \otimes \cdots \otimes \PP(j_{n,i_n})\arrow[d]
& \\
\PP(m) \otimes \PP(j_{1,1}) \otimes \cdots \otimes \PP(j_{n,i_n}) \arrow[r]
&\PP(l)~,
\end{tikzcd}
\]

where $r_k = j_{k,1} + \cdots + j_{k,i_k}$ for $k = 1, \cdots , n~,$ where $m = i_1 + \cdots + i_n ~,$ and where $l = r_1 + \cdots + r_n ~.$

\medskip

\item An element $\eta(1) \coloneqq \mathrm{id}_\PP \in \PP(1)$ such that $\gamma(n): \PP(1) \otimes \PP(n) \to \PP(n)$ on $(\mathrm{id}_\PP;-)$ is the identity of $\PP(n)$ and $\gamma(1, \cdots , 1): \PP(n) \otimes \PP(1) \otimes \cdots \otimes \PP(1)$ on $(-;\mathrm{id}_\PP, \cdots , \mathrm{id}_\PP)$ is also the identity of $\PP(n)$, for all $n \geq 0$. 
\end{enumerate}

\medskip

To any $\mathbb{S}$-module one can associate an analytic functor in the category $\kmod$ via the \textit{Schur realization functor}. It is given by

\[
\begin{tikzcd}[column sep=4pc,row sep=0.5pc]
\mathscr{S} : \smod \arrow[r]
&\mathsf{End}(\kmod) \\
M \arrow[r,mapsto]
&\mathscr{S}(M)(-) \coloneqq \bigoplus_{n \geq 0} M(n) \otimes_{\mathbb{S}_n} (-)^{\otimes n}~,
\end{tikzcd}
\]

where $\mathscr{S}(M)$ is called the \textit{Schur functor} associated to $M ~.$

\begin{lemma}
The Schur realization functor defines a strong monoidal functor between the monoidal categories $(\mathbb{S}\textsf{-}\mathsf{mod}, \circ, \I)$ and $(\mathsf{End}(\kmod), \circ,\mathsf{Id})$, where the monoidal product on endofuctors is given by composition. 
\end{lemma}

\begin{proof}
It is a straightforward computation since the tensor product distributes over coproducts.
\end{proof}

The functor $\mathscr{S}$ is strong monoidal, hence it preserves monoids. For any operad $\PP$, its Schur functor $\mathscr{S}(\PP)$ is canonically endowed with a monad structure. 

\begin{Definition}[$\PP$-algebras]
Let $(\PP, \gamma, \eta)$ be an operad. A $\PP$-\textit{algebra} $(A, \gamma_A)$ amounts to the data of an algebra over the monad $\mathscr{S}(\PP)$. The $\PP$-algebra structure $\gamma_A$ on the $\mathbb{K}$-module $A$ is therefore given by a map 

\[
\gamma_A: \bigoplus_{n \geq 0} \PP(n) \otimes_{\mathbb{S}_n} A^{\otimes n} \longrightarrow A~,
\]

which satisfies associativity and unitary conditions.
\end{Definition}

\begin{Remark}
The data of the map $\gamma_A$ is equivalent to the data of a family of maps 
\[
\left\{\gamma_A^n: \PP(n) \otimes_{\mathbb{S}_n} A^{\otimes n} \longrightarrow A \right\}
\]
for all $n \geq 0$. The associativity and unitary conditions of $\gamma_A$ can be rewritten in terms of these $\gamma_A^n$.
\end{Remark}

\begin{Example}
The classical examples include the operads $\mathcal{C}om$ and $\mathcal{A}ss$, which encode respectively non-unital commutative algebras and non-unital associative algebras as its algebras. See \cite[Chapter 13]{LodayVallette12} for more examples.
\end{Example} 

\begin{lemma}
Let $(\PP, \gamma, \eta)$ be an operad. The category of $\PP$-algebras is presentable. The free $\PP$-algebra on a $\mathbb{K}$-module $V$ is given by $\mathscr{S}(\PP)(V)$ together with the map 
\[
\mathscr{S}(\gamma)(V) : \mathscr{S}(\PP) \circ \mathscr{S}(\PP)(V) \cong \mathscr{S}(\PP \circ \PP)(V) \longrightarrow \mathscr{S}(\PP)(V)~.
\]
\end{lemma}

\begin{proof}
The ground category $\kmod$ is presentable and the endofunctor $\mathscr{S}(\PP)(-)$ is clearly accessible.
\end{proof}

\begin{Definition}[Endomorphism operad] 
To any $\mathbb{K}$-module $A$, one can associate an operad called the \textit{endomorphism operad}. It is given by the $\mathbb{S}$-module:

\[
\mathrm{End}_A(n) \coloneqq \mathrm{Hom}(A^{\otimes n},A)~.
\]

together with the composition maps given by the composition of functions, where the unit is $\mathrm{id}_A \in \mathrm{End}_A(1)$.  
\end{Definition} 

\begin{lemma}\label{lemma: endomorphisms = P alg}
Let $(\PP, \gamma, \eta)$ be an operad. The data of a $\PP$-algebra $(A,\gamma_A)$ is equivalent to a morphism of operads $\Gamma_A: \PP \longrightarrow \mathrm{End}_A$.
\end{lemma}

\begin{proof}
The maps $\gamma_A^n: \PP(n) \otimes_{\mathbb{S}_n} A^{\otimes n} \longrightarrow A$ transpose under the tensor-hom adjunction into $\mathbb{S}_n$-equivariant maps $\Gamma_A^n: \PP(n) \longrightarrow \mathrm{Hom}(A^{\otimes n},A)$. The compatibility conditions of a $\PP$-algebra and those of an operad morphism are sent to each other under this adjunction. 
\end{proof}

\section{Categorical dual of operads and cooperads}
If one wants a dual notion to that of operads, one can consider two possibilities: either consider comonoids in the category of $\mathbb{S}$-modules or dualize May's axioms. Unlike with operads, these two notions do not coincide. Cooperads in the sense of \cite{GetzlerJones94} or in \cite{LodayVallette12} are defined as comonoids in category of $\mathbb{S}$-modules. We introduce the notion of a May cooperad. Classical cooperads induce May cooperads, but the converse is not true. This phenomena appear due to the arity $0$ and arity $1$ components, which creates infinite sums in the decompositions of operations. 

\begin{Definition}[Cooperad]\label{def: cooperad}
A \textit{cooperad} $(\C,\Delta,\epsilon)$ is the data of a comonoid in $(\smod, \circ, \I)~.$ 
\end{Definition}

\begin{Remark}
Notice that the structural map

\[
\Delta: \C(n) \longrightarrow \bigoplus_{k\geq 0} \C(k) \otimes_{\mathbb{S}_k} \left( \bigoplus_{i_1 + \cdots + i_k = n} \mathsf{Ind}_{\mathbb{S}_{i_1} \times \cdots \times \mathbb{S}_{i_k}}^{\mathbb{S}_n} (\C(i_1) \otimes \cdots \otimes \C(i_k))\right)
\]

lands on a direct sum. Therefore an element in $\C(n)$ can only be decomposed into a \textit{finite sum} of elements.
\end{Remark}

If $\C$ is a cooperad, then $\mathscr{S}(\C)$ has a canonical comonad structure because the Schur realization functor is strong monoidal. This allows us to define a well-behaved category of coalgebras over a cooperad.

\begin{Definition}[$\C$-coalgebra]\label{def: C-coalgebra}
Let $(\C,\Delta,\epsilon)$ be a cooperad. A $\C$-\textit{coalgebra} $(C,\Delta_C)$ amounts to the data of a coalgebra over the comonad $\mathscr{S}(\C)$. The structural map $\Delta_C$ is therefore given by

\[
\Delta_C: C \longrightarrow \bigoplus_{n \geq 0} \C(n) \otimes_{\mathbb{S}_n} C^{\otimes n}~,
\]

and it satisfies coassociative and counitary conditions.
\end{Definition}

\begin{Remark}
The notion of a coalgebra over a cooperad can only encapsulate types of coalgebras which are conilpotent in some sense. The reason is that the structural map
\[
\Delta_C: C \longrightarrow \bigoplus_{n \geq 0} \mathcal{C}(n) \otimes_{\mathbb{S}_n} C^{\otimes n}~,
\]
lands on the direct sum of a non-completed tensor product. Thus decompositions of elements in $C$ are supported in finitely many arities and in finitely many steps of the coradical filtration of $\mathcal{C}$.
\end{Remark}

\begin{Example}
The linear dual $\mathcal{A}ss^*$ of $\mathcal{A}ss$ has a cooperad structure structure given by the dual maps of the operad structure. The category of $\mathcal{A}ss^*$-coalgebras is equivalent to the category of conilpotent non-counital coassociative coalgebras. Likewise, the category of $\mathcal{C}om^*$-coalgebras is equivalent to the category of conilpotent non-counital cocommutative coalgebras.
\end{Example}

\begin{lemma}
Let $(\C,\Delta,\epsilon)$ be a cooperad. The category of $\C$-coalgebras is presentable. The cofree $\C$-coalgebra on a $\mathbb{K}$-module $V$ is given by $\mathscr{S}(\C)(V)$ together with the map 
\[
\mathscr{S}(\Delta_\C)(V): \mathscr{S}(\C)(V) \longrightarrow \mathscr{S}(\C \circ \C)(V) \cong \mathscr{S}(\C) \circ \mathscr{S}(\C)(V)~.
\]
\end{lemma}

\begin{proof}
The ground category $\kmod$ is presentable and the endofunctor $\mathscr{S}(\C)(-)$ is clearly accessible.
\end{proof}

\begin{Remark}
The coassociative coalgebra $\mathscr{S}(\mathcal{A}ss^*)(V)$ is isomorphic to $\overline{\mathcal{T}}^c(V)$, the cofree conilpotent non-counital coassociative coalgebra, which \textit{is not} cofree in general. (Only satisfies the universal property with respect to other conilpotent coassociative coalgebras).
\end{Remark}

Let $(\PP, \gamma, \eta)$ be an operad such that $\PP(n)$ is a finite dimensional $\kk$-module for all $n \geq 0$. The linear dual of the composition map $\gamma^*: \PP^* \longrightarrow (\PP \circ \PP)^*$ lands on the linear dual of the composition product $(\PP \circ \PP)^*$. This linear dual is given by the \textit{complete composition product} of $\mathbb{S}$-modules. 

\begin{Definition}[Complete composition product]
Let $M$ and $N$ be two $\mathbb{S}$-modules, the \textit{complete composition product} $M ~\widehat{\circ}~ N$ is defined by: 

\[
M ~\widehat{\circ}~ N (n) \coloneqq \prod_{k \geq 0} \left( M(k) \otimes \left(\prod_{i_1 + \cdots + i_k = n} \mathsf{Coind}_{\mathbb{S}_{i_1} \times \cdots \times \mathbb{S}_{i_k}}^{\mathbb{S}_n} (N(i_1) \otimes \cdots \otimes N(i_k))\right)\right)^{\mathbb{S}_k}~,
\]

where the second product runs over all $k$-uples $(i_1, \cdots, i_k) \in \mathbb{N}^k$ such that $i_1 + \cdots + i_k = n$. When $k$ is fixed, there are only finitely many partitions of $n$ with $k$ terms, therefore the second product is finite and isomorphic to a direct sum.
\end{Definition}

\begin{lemma}
Let $M,N$ be an $\mathbb{S}$-module and let $M(n)$ be a finite dimensional $\kk$-module. There is an isomorphism of $\mathbb{S}$-modules $(M \circ N)^* \cong M ~\widehat{\circ}~ N ~.$
\end{lemma}

\begin{proof}
The main ingredient of the proof is the isomorphism shown in Lemma \ref{lemmaindandcoind}:

\[
\mathrm{Hom}(\mathsf{Ind}_H^G(V),\kk) \cong \mathsf{Coind}_H^G(V^*)~,
\]
\vspace{0.1pc}

where $G = \mathbb{S}_n$ and $H = \mathbb{S}_{i_1} \times \cdots \times \mathbb{S}_{i_k}~.$ The rest is a straightforward computation.
\end{proof}

\begin{Remark}\label{rmk: may cooperad = cooperad quand C(0) =0}
Since $\mathbb{K}$ is a field of characteristic 0, one can identify the coinvariants with in invariants in the definition above. In this case, there is also a natural isomorphism between $\mathsf{Ind}$ and $\mathsf{Coind}$. This gives a natural inclusion $M \circ N \hookrightarrow M ~\widehat{\circ}~ N~.$ Furthermore, if $N(0)=0$, there are finitely many non-trivial $k$-uples $(i_1, \cdots, i_k) \in \mathbb{N}^k$ such that $i_1 + \cdots + i_k = n~,$ implying that the natural inclusion gives an isomorphism $M \circ N \cong M ~\widehat{\circ}~ N~.$
\end{Remark}

Contrary to popular belief, if there are no assumptions on the arity $0$ component, the complete composite product $\widehat{\circ}$ does not form a monoidal product (tensor products do not distribute over infinite products). Therefore the linear dual of the composition map of an operad of finite dimension arity-wise $\gamma^*: \PP^* \longrightarrow \PP^* ~ \widehat{ \circ} ~ \PP^* $ can not satisfy a coassociativity axiom since the product $\widehat{\circ}$ \textit{is not monoidal}. 

\begin{Definition}[May cooperad]\label{May cooperad}
A \textit{May cooperad} structure on an $\mathbb{S}$-module $\C$ is the data of a family of decomposition maps:

\[
\Delta(i_1,\cdots,i_k): \C(n) \longrightarrow \C(k) \otimes \C(i_1) \otimes \cdots \otimes \C(i_k)
\]

for all $k$-tuples $(i_1,\cdots,i_k) \in \mathbb{N}$ such that $i_1 + \cdots + i_k =n$ and a counit map $\epsilon: \C \longrightarrow \I$ satisfying the categorical dual of May's axioms for operads. 

\begin{enumerate}
\medskip 

\item Each map $\Delta(i_1, \cdots, i_k)$ is $\mathbb{S}_{i_1} \times \cdots \times \mathbb{S}_{i_k}$-equivariant map where the action on $\C(n)$ is given by restricting the action of $\mathbb{S}_n$ to $\mathbb{S}_{i_1} \times \cdots \times \mathbb{S}_{i_k}~.$

\medskip

\item The map: 

\[
\bigoplus_{i_1 + \cdots + i_k = n} \Delta(i_1, \cdots, i_k) : \C(n) \longrightarrow \C(k) \otimes \left( \bigoplus_{i_1 + \cdots + i_k = n} \C(i_1) \otimes \cdots \otimes \C(i_k) \right)
\]

factors through the invariants of $\mathbb{S}_k$ on the right.

\medskip

\item For any set of indices $(j_{1,1}, \cdots, j_{1,i_1},j_{2,1}, \cdots, j_{2,i_2}, \cdots, j_{n,1}, \cdots, j_{n,i_n})$, the following square is commutative:

\[
\begin{tikzcd}[column sep=4pc,row sep=2pc]
\C(n) \otimes \C(r_1) \otimes \cdots \otimes \C(r_n) \arrow[d]
& \\
\C(n) \otimes \C(i_1) \otimes \C(j_{1,1}) \otimes \cdots \otimes \C(j_{1,i_1}) \otimes \C(i_2) \otimes \C(j_{2,1}) \otimes \cdots \otimes \C(j_{n,i_n})  \arrow[d,"\cong"]
& \\
\C(n) \otimes \C(i_1) \otimes \cdots \otimes \C(i_n) \otimes \C(j_{1,1}) \otimes \cdots \otimes \C(j_{n,i_n})
& \\
\C(m) \otimes \C(j_{1,1}) \otimes \cdots \otimes \C(j_{n,i_n}) \arrow[u]
&\C(l)~, \arrow[l]  \arrow[uuul,bend right = 35]
\end{tikzcd}
\]

where $r_k = j_{k,1} + \cdots + j_{k,i_k}$ for $k = 1, \cdots , n~,$ where $m = i_1 + \cdots + i_n ~,$ and where $l = r_1 + \cdots + r_n ~.$

\medskip

\item The following diagrams commute with the counit $\epsilon: \C \longrightarrow \I$ :
\[
\begin{tikzcd}[column sep=3pc,row sep=3pc]
\C(n) \arrow[r,"\Delta(n)"] \arrow[dr,"\mathrm{id}_{\C(n)}",swap] 
&\C(1) \otimes \C(n) \arrow[d,"\epsilon \hspace{1pt} \otimes \hspace{1pt}\mathrm{id}_{\C(n)}"]
&\C(n) \arrow[dr,"\mathrm{id}_{\C(n)}",swap] \arrow[r,"\Delta (1 \mathrm{,} \cdots \mathrm{,} 1)"]
&\C(n) \otimes \C(1) \otimes \cdots \otimes \widehat{\C}(1) \arrow[d,"\mathrm{id}_{\C(n)} \otimes \epsilon^{\otimes n}"] \\
&\C(n) \cong \mathbb{K} \otimes \C(n) 
&
&\C(n) \otimes \mathbb{K}^{\otimes n} \cong \C(n)
\end{tikzcd}
\]
\end{enumerate}

\medskip

A \textit{morphism} of May cooperads $f: (\C,\{\Delta_\C(i_1,\cdots,i_k)\},\epsilon_\C) \longrightarrow (\D,\{\Delta_\D(i_1,\cdots,i_k)\},\epsilon_\D)$ is an $\mathbb{S}$-morphism $f: \C \longrightarrow \D$ commuting with the decomposition maps and the counits. 
\end{Definition}

\begin{lemma}\label{lemma: May induit un delta sur le chapeau}
The data of decomposition maps $\{\Delta_\C(i_1,\cdots,i_k)\}$ on an $\mathbb{S}$-module $\C$ satisfying axioms $(1)$ and $(2)$ is equivalent to a morphism of $\mathbb{S}$-modules $\Delta: \C \longrightarrow \C~\widehat{\circ}~\C~.$ 
\end{lemma}

\begin{proof}
Suppose that the maps $\Delta(i_1,\cdots,i_k)$ satisfy axioms $(1)$ and $(2)$, then

\[
\bigoplus_{i_1 + \cdots + i_k = n} \Delta(i_1, \cdots, i_k) : \C(n) \longrightarrow \C(k) \otimes \left( \bigoplus_{i_1 + \cdots + i_k = n} \mathsf{Coind}_{\mathbb{S}_{i_1} \times \cdots \times \mathbb{S}_{i_k}}^{\mathbb{S}_n}\left( \C(i_1) \otimes \cdots \otimes \C(i_k)\right) \right)^{\mathbb{S}_k}
\]

are $\mathbb{S}_n$-equivariant. These maps assemble into an unique map $\Delta: \C \longrightarrow \C~\widehat{\circ}~\C$ by the universal property of the product. The other way around is also straightforward.
\end{proof}

\begin{lemma}\label{pété}
Let $\C$ be an $\mathbb{S}$-module such that $\C(0) = 0$. The data of a cooperad structure on $\C$ is equivalent to the data of a May cooperad structure on $\C$.
\end{lemma}

\begin{proof}
Since $\C(0)=0$, we have that $\C \circ \C \cong \C~ \widehat{\circ} ~\C$. By Lemma \ref{lemma: May induit un delta sur le chapeau}, the data of the family of maps $\{\Delta_\C(i_1,\cdots,i_k)\}$ satisfying axioms $(1)$ and $(2)$ is equivalent to a morphism of $\mathbb{S}$-modules $\C \longrightarrow \C \circ \C$. In this situation, the morphism $\C \longrightarrow \C \circ \C$ is coassociative with respect to the composition product if and only if the family of maps $\{\Delta_\C(i_1,\cdots,i_k)\}$ satisfies axiom $(3)$. The data of $\epsilon: \C \longrightarrow \I$ and the counit axiom $(4)$ are equivalent as well.
\end{proof}

\begin{Proposition}\label{comparaisoncoop}
Let $\C$ be an $\mathbb{S}$-module. 

\medskip

\begin{enumerate}
\item Any cooperad structure $(\Delta, \epsilon)$ on $\C$ induces a May cooperad structure $(\{\Delta(i_1,\cdots,i_k)\},\epsilon)$ by first composing $\Delta$ with the natural inclusion $\C \circ \C \hookrightarrow \C ~\widehat{\circ}~\C$ and then splitting the composite into $\{\Delta(i_1,\cdots,i_k)\}$ using Lemma \ref{lemma: May induit un delta sur le chapeau}.

\medskip

\item Let $(\{\Delta(i_1,\cdots,i_k)\},\epsilon)$ be a May cooperad structure on $\C$. It induces a cooperad structure on $\C$ if and only if the map $\Delta: \C \longrightarrow \C ~\widehat{\circ}~\C$ given by Lemma \ref{lemma: May induit un delta sur le chapeau} factors through $\C \circ \C~.$
\end{enumerate}
\end{Proposition}

\begin{proof}
The proof is straightforward, completely analogue to the proof of Lemma \ref{pété}.
\end{proof}

\begin{Definition}[Strict May cooperad]
A May cooperad $(\C,\{\Delta(i_1,\cdots,i_k)\},\epsilon)$ is a \textit{strict May cooperad} if the decomposition map $\Delta: \C \longrightarrow \C ~\widehat{\circ}~\C$ does not factor through $\C \circ \C~.$
\end{Definition}

\begin{Example}
Let $\ucom$ be the operad encoding unital commutative algebras. Although it is finitely dimensional in every arity, its linear dual $\ucomd$ cannot be endowed with a cooperad structure. The presence of the unit $u \in \ucomd(0)$ creates infinite sums in the 2-leveled decomposition of any operation in $\ucomd$. Hence the map 

\[
\Delta_{\ucomd} : \ucomd \longrightarrow \ucomd ~\widehat{\circ}~ \ucomd
\]
\vspace{0.1pc}

does not factor through $\ucomd \circ \ucomd$. The same phenomenon can be observed for $u\mathcal{A}ss^*$, the linear dual of the operad encoding unital associative algebras.
\end{Example}  

\section{Partial operads and partial cooperads}\label{Section: partial (co)operads}
In this section we give the definition of a partial operad, introduced by M. Markl in \cite{Markl96}, and the dual definition of a partial cooperad. Unital partial operads are equivalent to operads defined as monoids, whereas counital partial cooperads are equivalent to May cooperads.

\begin{Definition}
Let $n,k \in \mathbb{N}$, for all $1 \leq i \leq n$, we define maps:
\begin{equation}\label{morph: rho i}
\circ_i: \mathbb{S}_n \times \mathbb{S}_k \longrightarrow \mathbb{S}_{n+k-1}~.
\end{equation}
For $(\tau,\sigma) \in \mathbb{S}_n \times \mathbb{S}_k$, $\sigma \circ_i \tau$ is given by the unique permutation in $\mathbb{S}_{n+k-1}$ which first acts as $\tau$ on the set $\{1,\cdots,n+k-1\} - \{i,\cdots, i+k-1\}$ and sends $\{i,\cdots, i+k-1\}$ to $\{\tau(i), \cdots, \tau(i) + k -1 \}$; then acts as $\sigma$ on the set $\{\tau(i), \cdots, \tau(i) + k -1 \}$.
\end{Definition}

\begin{Definition}[Partial operad]\label{def: partialoperad}
A \textit{partial operad} $(\PP, \{\circ_i\})$ is the data of a graded $\mathbb{S}$-module $\PP$ endowed with a family of \textit{partial composition maps} of degree $0$:
\[
\circ_i : \PP(n) \otimes \PP(k) \longrightarrow \PP(n+k-1)
\]
for $1 \leq i \leq n$, subject to the following conditions.

\begin{enumerate}

\item It satisfies the \textit{sequential} axiom: for $1 \leq i \leq n, 1 \leq j \leq k$, the following diagram commutes
\[
\begin{tikzcd}[column sep=4.5pc,row sep=3pc]
\PP(n) \otimes \PP(k) \otimes \PP(m) \arrow[r,"\circ_i ~\otimes~ \mathrm{id}_{\PP(m)}"] \arrow[d,"\mathrm{id}_{\PP(n)}~ \otimes~ \circ_j",swap] 
&\PP(n+k-1) \otimes \PP(m) \arrow[d,"\circ_{i+j-1} "] \\
\PP(n) \otimes \PP(k+m-1) \arrow[r,"\circ_i"]
&\PP(n+k+m-2)~.
\end{tikzcd}
\]
\item It satisfies the \textit{parallel} axiom: for $1 \leq i < j \leq n$, the following diagram commutes

\[
\begin{tikzcd}[column sep=4.5pc,row sep=3pc]
\PP(n) \otimes \PP(k) \otimes \PP(m) \arrow[d,"\circ_i  ~\otimes~ \mathrm{id}_{\PP(m)}",swap] \arrow[r,"\circ_j ~\otimes~ \mathrm{id}_{\PP(k)}"]
&\PP(n+m-1) \otimes \PP(k) \arrow[d,"\circ_i"] \\
\PP(n+k-1) \otimes \PP(m) \arrow[r,"\circ_{j+k-1}"]
&\PP(n+k+m-2)~.
\end{tikzcd}
\]

\item The maps
\[
\circ_i : \PP(n) \otimes \PP(k) \longrightarrow \PP(n+k-1)
\]
satisfy the following condition: let $(\tau,\sigma) \in \mathbb{S}_n \times \mathbb{S}_k$, then the following diagram commutes
\[
\begin{tikzcd}[column sep=3pc,row sep=3pc]
\PP(n) \otimes \PP(k) \arrow[r,"\circ_i "] \arrow[d,"\tau ~\otimes~\sigma ",swap]
&\PP(n+k-1) \arrow[d,"\tau ~\circ_i~ \sigma " ] \\
\PP(n) \otimes \PP(k) \arrow[r,"\circ_{\tau(i)} "]
&\PP(n+k-1)~.
\end{tikzcd}
\]
\end{enumerate} 
\end{Definition}

\begin{Remark}
Partial operads are sometimes called Markl operads in the literature, see \cite{Markl96}.
\end{Remark}

\begin{Notation}
The maps $\circ_i$ should be denoted $\circ_i^{n,k}$ since they are given for each pair $(n,k)$. For clarity, we simply denote them $\circ_i$ when there is no ambiguity. 
\end{Notation}

\begin{Definition}[Unital partial operad]\label{def: unital partial operad}
A \textit{unital partial operad} $(\mathcal{P},\{\circ_i\},\eta)$ is the data of a partial operad $(\mathcal{P},\{\circ_i\})$ together with a morphism $\eta: \I \longrightarrow \PP$, such that $\circ_1(\eta(\mathrm{id}),-)$ and $\circ_i(-,\eta(\mathrm{id}))$ are the identity on $\PP(n)$.
\end{Definition}

\begin{Proposition}\label{prop: unital partial operads are operads}
The category of operads defined as monoids is equivalent to the category of unital partial operads. 
\end{Proposition}

\begin{proof}
Let $(\PP,\{\circ_i\},\eta)$ be a unital partial operad, the map 

\[
\gamma(i_1, \cdots, i_k): \PP(k) \otimes \PP(i_1) \otimes \cdots \otimes \PP(i_k) \longrightarrow \PP(n)~,
\]
\vspace{0.1pc}

given by iterating $k$ times the partial composition maps $(- \circ_1 (- \circ_2 ( \cdots (\circ_k -)))))$ lands on $\PP(k + i_1 + \cdots i_k -k)$. The axioms of a partial operad ensure that $\gamma(i_1, \cdots, i_k)$ satisfies May's axioms, hence it endows $\PP$ with an operad structure. Conversely, given an operad $(\PP,\gamma,\eta)$ with a unit, partial compositions maps by plugging the unit everywhere $\mu \circ_i \nu \coloneqq \gamma(\mu; \mathrm{id}, \cdots,  \nu, \cdots, \mathrm{id})$, except at the $i^{\text{th}}$-place. Mays's axioms ensure that it satisfies the axioms of a unital partial operad. 
\end{proof}

\begin{Remark}
Note that the presence of the unit is crucial for the converse of the proof. One could also consider a non-unital version of an operad, that is an $\mathbb{S}$-module $\PP$ endowed with a map $\gamma: \PP \circ \PP \longrightarrow \PP$ that is associative. But this type of structure would not be equivalent to a partial operad structure. In fact, one can always construct such a map $\gamma$ from the partial composition map $\{\circ_i\}$, but not the other way around. In a partial operad, one knows where the operations are composed, hence it is a finer structure. 
\end{Remark}

\begin{Definition}[Augmented partial operad]
An \textit{augmented} partial operad $(\PP, \{\circ_i\}, \eta, \xi)$ is the data of a unital partial operad $(\PP, \{\circ_i\}, \eta)$ equipped with a morphism of unital partial operads $\xi: \PP \longrightarrow \I$ such that $\xi \circ \eta = \mathrm{id}$. 
\end{Definition}

Since $\xi: \PP \longrightarrow \I$ is a morphism of operads, the operadic structure of $\PP$ splits in $\PP \cong \overline{\PP} \oplus \I$, where $\overline{\PP} \coloneqq \text{Ker}~\xi~$. This entails that the unit $\mathrm{id}$ in $\PP$ cannot be written as the partial composition of other operations in $\PP$. 

\begin{Proposition}[Unital extension]\label{prop: add a unit to a partial operad}
Let $(\PP, \{\circ_i\})$ be a partial operad. Consider

\[
\PP^{u} \coloneqq (\PP \oplus \I, \{\circ_i\}, \eta, \xi)~,
\]

where $\eta$ is given by the inclusion $\I \hookrightarrow \PP \oplus \I$ and $\xi$ by the projection $\PP \oplus \I \twoheadrightarrow \I$.
This data forms an augmented partial operad. This defines a fully faithful functor from the category of partial operads to the category of unital partial operads which is left adjoint of the forgetful functor:  

\[
\begin{tikzcd}[column sep=4pc,row sep=2pc]
\partialop \arrow[r,"(-)^u ",""{name=A, below},bend left] & \upartialop \arrow[l,bend left,"\mathrm{U}",""{name=B,above}] \arrow[phantom, from=B, to=A, "\dashv" rotate=-90] 
\end{tikzcd}~.
\]

Its essential image is the category of augmented partial operads, hence this adjunction restricts to an equivalence between the category of partial operads and the category of augmented partial operads.
\end{Proposition}

\begin{proof}
It is a straightforward computation. 
\end{proof}

\begin{Definition}[Partial cooperad]\label{def: partialcoop}
A \textit{partial cooperad} $(\C, \{\Delta_i\})$ is the data of a graded $\mathbb{S}$-module $\C$ together with \textit{partial decomposition maps}: 
\[
\Delta_i: \C(n+k-1) \longrightarrow \C(n) \otimes \C(k)
\]
for $1 \leq i \leq n$, subject to the following conditions.

\begin{enumerate}
\item It satisfies the \textit{sequential} axiom: for $1 \leq i \leq n, 1 \leq j \leq k$, the following diagram commutes
\[
\begin{tikzcd}[column sep=4.5pc,row sep=3pc]
\C(n+k+m-2) \arrow[r,"\Delta_{i+j-1}"] \arrow[d,"\Delta_i ",swap]
&\C(n+k-1) \otimes \C(m) \arrow[d,"\Delta_i  \hspace{1pt} \otimes  \hspace{1pt} \mathrm{id}_{\C(m)}"] \\
\C(n) \otimes \C(k+n-1) \arrow[r,"\mathrm{id}_{\C(n)}  \hspace{1pt} \otimes  \hspace{1pt} \Delta_j "]
&\C(n) \otimes \C(k) \otimes \C(m)~.
\end{tikzcd}
\]
\item It satisfies the \textit{parallel} axiom: for $1 \leq i < j \leq n$, the following diagram commutes
\[
\begin{tikzcd}[column sep=4.5pc,row sep=3pc]
\C(n+k+m-2) \arrow[r,"\Delta_{i+k-1}"] \arrow[d,"\Delta_i ",swap]
&\C(n+k-1) \otimes \C(m) \arrow[d,"\Delta_i \hspace{1pt} \otimes  \hspace{1pt} \mathrm{id}_{\C(m)}"] \\
\C(n+m-1) \otimes \C(k) \arrow[r,"\Delta_j  \hspace{1pt} \otimes  \hspace{1pt} \mathrm{id}_{\C(k)}"]
&\C(n) \otimes \C(k) \otimes \C(m)~.
\end{tikzcd}
\]
\item The maps
\[
\Delta_i: \C(n+k-1) \longrightarrow \C(n) \otimes \C(k)
\]
satisfy the following condition: let $(\tau,\sigma) \in \mathbb{S}_n \times \mathbb{S}_k$, then the following diagram commutes
\[
\begin{tikzcd}[column sep=3pc,row sep=3pc]
\C(n+k-1) \arrow[r,"\Delta_{\tau(i)} "] \arrow[d,"\tau ~\circ_i~\sigma ",swap]
&\C(n) \otimes \C(k)  \arrow[d,"\tau ~ \otimes ~ \sigma "]\\
\C(n+k-1) \arrow[r,"\Delta_i "]
&\C(n) \otimes \C(k)~.
\end{tikzcd}
\] 
\end{enumerate} 
\end{Definition}

\begin{Notation}
The maps $\Delta_i$ should be denoted $\Delta_i^{n,k}$ since they are given for each pair $n,k \in \mathbb{N}$. Again, for clarity, we simply denote them $\Delta_i$ when there is no ambiguity. 
\end{Notation}

\begin{Definition}[Counital partial cooperad]
A \textit{counital partial cooperad} $(\C, \{\Delta_i\},\epsilon)$ is the data of a partial cooperad $(\C, \{\Delta_i\})$ and a morphism of $\mathbb{S}$-modules $\epsilon: \C \longrightarrow \I$ such that:

\[
\begin{tikzcd}
\C(n) \arrow[r,"\Delta_i"] \arrow[rd,"\cong",swap]
&\C(n) \otimes \C(1) \arrow[d,"\mathrm{id}~\otimes~\epsilon"]
&\C(n) \arrow[r,"\Delta_1"] \arrow[rd,"\cong",swap]
&\C(1) \otimes \C(n) \arrow[d,"\epsilon~\otimes~\mathrm{id}"] \\
&\C(n) \cong \C(n) \otimes \kk
&
&\C(n) \cong  \kk \otimes \C(n)~. \\
\end{tikzcd}
\]
\end{Definition}

\begin{Proposition}\label{prop: May coop = counital partial coop}
The category of counital partial cooperads is equivalent to the category of May cooperads. 
\end{Proposition}

\begin{proof}
The proof is completely dual to the proof of Proposition \ref{prop: unital partial operads are operads}. To construct the maps $\Delta(i_1, \cdots, i_k)$ we iterate the partial decompositions maps: 

\[
\Delta(i_1, \cdots, i_k) \coloneqq \Delta_1^{k, i_1} \circ \cdots \circ \Delta_k^{n-i_k+1,i_k}
\]

for each $k$-tuple $(i_1, \cdots, i_k)$ such that $i_1 + \cdots + i_k = n$. The axioms of a partial cooperad ensure that these maps satisfy the dual of May's axioms. The data of the counit on both sides is equivalent. 

\medskip

Conversely, the presence of the counit in the May cooperad allows us to define partial composition maps:

\[
\Delta_i^{n,k} \coloneqq (\epsilon, \cdots, \epsilon, \mathrm{id}, \epsilon, \cdots, \epsilon) \circ \Delta(1, \cdots, 1, k, 1, \cdots, 1)
\]

by killing all the decompositions outside of the $i^{\text{th}}$-place. May's dual axioms on the maps $\Delta(i_1, \cdots, i_k)$ ensure that these partial decomposition maps satisfy the axioms of a counitary partial cooperad. 
\end{proof}

\begin{Remark}
As in the case of operads, a partial cooperad structure $\{\Delta_i\}$ allows us to construct maps $\Delta(i_1, \cdots, i_k)$, but without a counit it is impossible to go the other way around. Hence a partial cooperad structure is finer than a May cooperad structure in the case where there is no counit. 
\end{Remark}

\begin{Definition}[Coaugmented partial cooperad]
A \textit{coaugmented} partial cooperad $(\C, \{\Delta_i\}, \epsilon, \nu)$ is the data of a counital partial cooperad $(\C, \{\Delta_i\}, \epsilon)$ equipped with a morphism of counital partial cooperads $\nu: \I \longrightarrow \C$ such that $\epsilon \circ \nu = \mathrm{id}$. 
\end{Definition}

Since $\nu : \I \longrightarrow \C$ is a morphism of partial cooperads, the structure of $\C$ splits in $\C \cong \overline{\C} \oplus \I$, where $\overline{\C} \coloneqq \text{Ker}~\epsilon~$. This entails that the counit $\nu(\mathrm{id})$ has no non-trivial partial decompositions in the partial cooperad $\C$.

\begin{Proposition}[Counital extension]\label{prop: add a counit to a partial cooperad}
Let $(\C, \{\Delta_i\})$ be a partial cooperad. Consider

\[
\C^{u} \coloneqq (\C \oplus \I, \{\Delta_i\}, \epsilon, \nu)~,
\]

where $\nu$ is given by the inclusion $\I \hookrightarrow \C \oplus \I$ and $\epsilon$ by the projection $\C \oplus \I \twoheadrightarrow \I$.
It forms an coaugmented partial cooperad. This defines a fully faithful functor from the category of partial cooperads to the category of counital partial cooperads which is right adjoint of the forgetful functor:  

\[
\begin{tikzcd}[column sep=4pc,row sep=2pc]
\partialcoop \arrow[r,"(-)^u ",""{name=A, below},bend left] & \upartialcoop \arrow[l,bend left,"\mathrm{U}",""{name=B,above}] \arrow[phantom, from=A, to=B, "\dashv" rotate=90] 
\end{tikzcd}~.
\]

Its essential image is exactly the category of coaugmented partial cooperads, hence this adjunction restricts to an equivalence of categories between the category of partial cooperads and the category of coaugmented partial cooperads.
\end{Proposition}

\begin{proof}
It is a straightforward computation. 
\end{proof}

\section{Filtrations on partial (co)operads}\label{Section: Filtration on partial (co)operads}
The category of unital partial cooperads and the category of partial operads can both be encoded as algebras over monads. The first monad $\mathscr{T}$ is called the \textit{tree monad}, and the second monad $\overline{\mathscr{T}}$ is called the \textit{reduced tree monad}. Let us briefly recall their definitions.

\medskip

We define the endofunctor $\overline{\mathscr{T}}$ on the category of $\mathbb{S}$ as follows: for an $\mathbb{S}$-module $M$, a basis of $\overline{\mathscr{T}}(M)(n)$ is given by non-planar rooted trees $\tau$ of arity $n$, where the vertices of $\tau$ are labeled with elements of $M$ in the following way. A vertex $v \in \text{Ver}(\tau)$ with $k$ incoming edges is labeled with an element of $M(k)$. The leaves of the rooted trees are labeled with numbers $\{1,\cdots,n\}$, the action of $\mathbb{S}_n$ is given by permuting the labels. This endofunctor has a monad structure given by the substitution of trees and where the unit is given by the projection of the corollas labeled by $M$ onto $M$. If we add the trivial tree $|$ to $\overline{\mathscr{T}}$, then we get $\mathscr{T}$. See \cite[Section 5.6]{LodayVallette12} for more details.

\begin{Proposition}
There are isomorphisms of categories between:

\medskip

\begin{enumerate}
\item The category of algebras over the tree monad $\mathscr{T}$ and the category of unital partial operads.

\medskip

\item The category of algebras over the reduced tree monad $\overline{\mathscr{T}}$ and the category of partial operads.
\end{enumerate}
\end{Proposition}

\begin{proof}
Let $(\PP,\{\circ_i\})$ be a partial operad. One can construct a map $\gamma_\PP: \overline{\mathscr{T}}(\PP) \longrightarrow \PP$ by sending any rooted tree $\tau$ labeled with elements of $\PP$ to their corresponding iteration of compositions in $\PP$ along $\tau$. This maps satisfies the associative and unital axioms of an algebra over the reduced tree monad. 

\medskip

Conversely, by restricting an algebra map $\gamma_\PP: \overline{\mathscr{T}}(\PP) \longrightarrow \PP$ to all the rooted trees with two vertices, one gets the partial composition of a partial operad. If $\PP$ is unital of unit $\eta$, one sends $|$ to $\eta(\mathrm{id})$. The other way around, restricting $\gamma_\PP$ to $|$ gives the unit map $\eta$.
\end{proof}

\begin{Remark}
Note that $\mathscr{T}(M)(n)$ is given by the direct sum over rooted trees, hence only finite families appear. If $M(0)=0$, this does not make a difference since there are only finitely many rooted trees that can appear. On the other hand, when $M$ contains elements of arity $0$, there are infinitely many rooted trees because of the "corks" in $M(0)$.
\end{Remark}

The tree monad has a natural weight grading $\mathscr{T}^{(\omega)}$ given by the number of internal edges $\omega$ of the rooted trees. Notice that: 
\[
\mathscr{T}^{(-1)}(M) = |~, \quad \mathscr{T}^{(0)}(M) = M~, \quad \text{and}~~ \mathscr{T}^{(1)}(M) = M \circ_{(1)} M~.
\]

\begin{Definition}[Canonical filtration on a partial operad] 
Let $(\PP,\{\circ_i\})$ be a partial operad and let $\gamma_\PP: \overline{\mathscr{T}}(\PP) \longrightarrow \PP$ be its $\overline{\mathscr{T}}$-algebra structure. Its \textit{canonical filtration} is the decreasing filtration given by
\[ 
\mathscr{F}_{\omega} \PP \coloneqq \mathrm{Im}\left(\gamma_\PP^{(\geq \omega)}: \overline{\mathscr{T}}^{(\geq \omega)}(\PP) \longrightarrow \PP \right)
\]
for all $\omega \geq 0$, where $\overline{\mathscr{T}}^{(\geq \omega)}(\PP)$ denotes the elements of weight greater or equal to $\omega~.$ Each $\mathscr{F}_{\omega} \PP$ defines an operadic ideal of $\PP$ since they are stable under partial composition. Notice that: 
\[
\PP = \mathscr{F}_{0} \PP \supseteq \mathscr{F}_{1} \PP \supseteq \cdots \supseteq \mathscr{F}_{\omega} \PP \supseteq \cdots.
\]
\end{Definition}

\begin{Remark}
If we were to define an analogue filtration for unital partial operads which are not augmented it would be trivial. Indeed, suppose $(\PP,\{\circ_i\},\eta)$ is a unital partial operad such that the unit $\eta(\mathrm{id})$ can be written as a non-trivial composition of operations in $\PP$. Any operation is equal to its composition with $\eta(\mathrm{id})$ by definition, thus we can obtain any operation as the composition of an arbitrary number of operations in $\PP$. This implies that $\mathscr{F}_{\omega} \PP = \PP$ for all $\omega \geq 0~.$ A good example of the phenomenon is $\ucom$, where the unit $\mathrm{id} \in \ucom(1)$ can be written as the composite of $\mu_2 \in \ucom(2)$ and $\mu_0 \in \ucom(0)$. 
\end{Remark}

\begin{Definition}[Nilpotent partial operad]
Let $(\PP,\{\circ_i\})$ be a partial operad. It is said to be \textit{nilpotent} if there exists an $\omega \geq 1$ such that 
\[
\PP/\mathscr{F}_{\omega} \PP \cong \PP~.
\]
The partial operad is said to be $\omega_0$\textit{-nilpotent} if $\omega_0$ is the smallest integer such that the above isomorphism exits.
\end{Definition}

\begin{Definition}[Completion of a partial operad]\label{def: completion functor for operads}
Let $(\PP,\{\circ_i\})$ be a partial operad, its \textit{completion} $\widehat{\PP}$ is given by the following limit
\[
\widehat{\PP} \coloneqq \lim_{\omega} \PP/\mathscr{F}_{\omega} \PP 
\]
taken in the category of partial operads. 
\end{Definition}

Any morphism of partial operads $f: (\PP,\{\circ_i\}_\PP) \longrightarrow (\Q,\{\circ_i\}_\Q)$ is continuous with respect to the canonical filtration (i.e $f(\mathscr{F}_{\omega} \PP) \subseteq \mathscr{F}_{\omega} \Q$). Thus the completion of partial operads is functorial. For every $\omega \geq 1$, there are projection morphisms of partial operads $\varphi_\omega: \PP \twoheadrightarrow \PP/\mathscr{F}_{\omega} \PP$ which induce a canonical morphism of partial operads: 

\[
\varphi: \PP \longrightarrow \widehat{\PP}. 
\]
\begin{Definition}[Complete partial operad]\label{def: complete operads}
Let $(\PP,\{\circ_i\})$ be a partial operad. It is \textit{complete} if the canonical morphism $\varphi: \PP\longrightarrow \widehat{\PP}$ is an isomorphism of partial operads. 
\end{Definition} 

\begin{Example}
Any nilpotent partial operad is complete. Any complete partial operad is the limit of a tower of nilpotent partial operads.
\end{Example}

\begin{Example}
Any partial operad $(\PP,\{\circ_i\})$ generated by operations of arity greater or equal to $2$ is complete. Any operation of arity $n$ is the composite of at most $n-1$ generators. For binary operads, the canonical filtration is given by $\mathscr{F}_{\omega} \PP = \{\PP(n)\}_{n \geq \omega +1}$. This includes the \textit{three graces} $\mathcal{L}ie$, $\mathcal{A}ss$ and $\mathcal{C}om$, as well as many well known binary operads. Other examples of complete partial operads are $\mathcal{L}_{\infty}, \mathcal{A}_{\infty}$ and $\mathcal{C}_{\infty}$, the partial operads encoding $\mathcal{L}_{\infty}$-algebras, $\mathcal{A}_{\infty}$-algebras and $\mathcal{C}_{\infty}$-algebras respectively.
\end{Example} 

\begin{Definition}[Completed tree endofunctor]
The \textit{completed reduced tree endofunctor} $\overline{\mathscr{T}}^\wedge$ is the endofunctor of the category of $\mathbb{S}$-modules given by
\[
\overline{\mathscr{T}}^\wedge(M) \coloneqq \lim_{\omega} \overline{\mathscr{T}}(M)/\overline{\mathscr{T}}^{(\geq \omega)}(M)~,
\]
that is, the completion of the tree monad with respect to this weigh filtration. 
\end{Definition}

\begin{Remark}
The reduced complete tree endofunctor is isomorphic to the \textit{product} over rooted trees where vertices are labeled by elements of $M$ in the following way: if $\tau$ is a rooted tree and $v$ is one of its vertices with $k$ incoming edges, then $v$ must be labeled with an element of $M(k)$. In particular, if $M(0)=0$ and $M(1)=0$, then there are only a finite number of rooted trees in each arity and one has that 
\[
\overline{\mathscr{T}}^\wedge(M) \cong \overline{\mathscr{T}}(M)~.
\]
\end{Remark}

In order to compare the notion of a (possibly counital) partial cooperad with the notion of a cooperads defined as a comonoid, one needs to introduce filtrations.

\begin{lemma}\label{Existence of delta C}
Let $(\C,\{\Delta_i\})$ be a partial cooperad. The partial decomposition maps $\{\Delta_i\}$ induce a morphism of $\mathbb{S}$-modules

\[
\Delta_\C: \C \longrightarrow \overline{\mathscr{T}}^{\wedge}(\C)~,
\]

where $\overline{\mathscr{T}}^{\wedge}(\C)$ denotes the $\mathbb{S}$-module of the completed reduced tree endofunctor on $\C$.
\end{lemma}

\begin{proof}
The map $\Delta_\C$ is defined as follows: it maps to an operation $\mu \in \C(n)$ the sum of all possible rooted trees associated to the partial decompositions of this operation. This sum might not be finite, but since the target is complete, this gives a well defined map of $\mathbb{S}$-modules.
\end{proof}

\begin{Definition}[Coradical filtration of a partial cooperad]\label{coradical filtration coop}
Let $(\C, \{ \Delta_i \})$ be a partial cooperad. The \textit{coradical filtration} of $\C$ is the increasing filtration given by 
\[
\mathscr{R}_\omega \C \coloneqq \mathrm{Ker}\left(\Delta_\C^{(\geq \omega)}: \C \longrightarrow (\overline{\mathscr{T}}^{\wedge})^{(\geq \omega)}(\C) \right)
\]
for all $\omega \geq 0$, where $(\overline{\mathscr{T}}^{\wedge})^{(\geq \omega)}$ denotes the elements of weight greater or equal to $\omega$. Each $\mathscr{R}_\omega \C$ defines a partial sub-cooperad of $\C$ made of operations in $\C$ which admit non-trivial $\omega$-iterated partial decompositions. Notice that 
\[
0 = \mathscr{R}_{0} \C \subseteq \mathscr{R}_{1} \C \subseteq \cdots \subseteq \mathscr{R}_{\omega} \C \subseteq \cdots.
\]
\end{Definition} 

\begin{Definition}[Conilpotent partial cooperad]\label{conilpartcoop}
Let $(\C, \{\Delta_i \})$ be a partial cooperad. It is \textit{conilpotent} if the canonical map from the colimit
\[ 
\psi: \colim_{\omega} \mathscr{R}_{\omega} \C \longrightarrow \C 
\]
is an isomorphism of partial cooperads. This is equivalent to the fact that the coradical filtration is exhaustive. 
\end{Definition}

The underlying endofunctor of the reduced tree monad $\overline{\mathscr{T}}$ has a comonad structure $(\Delta_{\overline{\mathscr{T}}},\epsilon_{\overline{\mathscr{T}}})$, where $\Delta_{\overline{\mathscr{T}}}$ is given by partitioning trees and $\epsilon_{\overline{\mathscr{T}}}$ by the inclusion of the corollas. See \cite[Section 5.8.8]{LodayVallette12} for more details. The endofunctor $\overline{\mathscr{T}}$ together with its comonad structure will be denoted $\overline{\mathscr{T}}^c$, and called \textit{the reduced tree comonad}.

\begin{Remark}
The same coproduct $\Delta_{\overline{\mathscr{T}}}$ cannot be extended to the underlying endofunctor of the tree monad $\mathscr{T}$; the presence of the trivial rooted tree $|$ creates infinite sums and does \textbf{not} induce a comonad structure on $\mathscr{T}$. 
\end{Remark}

\begin{Proposition}\label{Prop: conil coop = coalgebra over tree comonad}
The category of coalgebras over the reduced tree comonad $\overline{\mathscr{T}}^c$ is equivalent to the category of conilpotent partial cooperads.
\end{Proposition}

\begin{proof}
Let $(\C, \Delta_\C)$ be a coalgebra over $\overline{\mathscr{T}}^c$. Consider the map

\[
\Delta_\C^{(2)}: \C \longrightarrow \C \circ_{(1)} \C
\]

obtained by composition the structural map $\Delta_\C$ with the projection into weigh $1$ elements. This induces partial decomposition maps $\{\Delta_i\}$ on $\C$. One can show that the coassociativity and counitary axioms of $\Delta_\C$ imply that $\{\Delta_i\}$ endow $\C$ with a partial cooperad structure. Furthermore, since the reduced tree comonad is defined as a direct sum over rooted trees, this partial cooperad structure is indeed conilpotent. 

\medskip

On the other hand, given a partial cooperad $(\C,\{\Delta_i\})$, by Lemma \ref{Existence of delta C}, there exists a map of $\mathbb{S}$-modules 

\[
\Delta_\C: \C \longrightarrow \overline{\mathscr{T}}^{\wedge}(\C)~.
\]

If $\C$ is a conilpotent partial cooperad, this map factors through the inclusion $\overline{\mathscr{T}}(\C) \hookrightarrow \overline{\mathscr{T}}^{\wedge}(\C)~.$ One can check that this map endows $\C$ with a structure of coalgebra over the reduced tree comonad.
\end{proof}

\begin{Remark}
There are two main obstructions for a partial cooperad $(\C,\{\Delta_i\})$ to be conilpotent. 

\begin{enumerate}
\item The \textit{vertical} obstruction, which comes from $(\C(1), \Delta(1))$ and its non-unital coassociative coalgebra structure. If $(\C(1), \Delta(1))$ is not a conilpotent coassociative coalgebra, since each $\C(n)$ is a $\C(1)$-bicomodule, one could have an operation in $\C(n)$ which decomposes as another operation in $\C(n)$ and an unary operation in $\C(1)$, and this unary operation could be decomposed \textit{ad libitum}. Nevertheless, if $\C(1)$ is conilpotent, then there can be no such obstruction as the each $\C(n)$ is conilpotent as a $\C(1)$-bicomodule (admits an exhausitve filtration as coming from its $\C(1)$-bicomodule structure).

\medskip

\item The \textit{horizontal} obstruction, which comes from the existence of 0-arity elements in $\C$. For a operation in $\C(n)$, one could iterate the following decompositions $\C(n) \longrightarrow \C(n+1) \otimes \C(0) \longrightarrow \C(n+2) \otimes \C(0)^{\otimes 2} \longrightarrow \cdots$ which can go on \textit{ad libitum}. If we add a cofree counit to $\C$, and view it as a May cooperad, this is the obstruction for $\C^{u}$ to be a cooperad. 
\end{enumerate}
\end{Remark}

Let $(\C,\{\Delta_i\})$ be a partial cooperad and let $\mathscr{R}_{\omega}\C$ be its coradical fltration. There is a coradical filtration on $\C^{u}$ which starts at $0$ by declaring that $\mathscr{R}_{0}\C^{u} = 0$ and $\mathscr{R}_{\omega}\C^{u} = \I \oplus \mathscr{R}_{\omega}\C$ for all $\omega \geq 1$. We denote $\iota_\omega$ the inclusion morphism of partial cooperads $\mathscr{R}_{\omega}\C \hookrightarrow \C$ and $\pi_\omega$ the projection morphism of $\mathbb{S}$-modules $\C \twoheadrightarrow \C/\mathscr{R}_\omega \C~.$

\begin{Proposition}
Let $(\C,\{\Delta_i\})$ be a conilpotent partial cooperad. There is a coassociative morphism of $\mathbb{S}$-modules $\Delta: \C \longrightarrow \C \circ \C$ that induces a coaugmented cooperad structure on $(\C^{u}, \Delta, \epsilon, \nu)~.$ 
\end{Proposition}

\begin{proof}
The outline of this proof was done similarly for the coproperadic case using the graph comonad. See \cite[Proposition 2.18]{hoffbeck2019properadic} for more details. Let $(\C,\{\Delta_i\})$ be a conilpotent partial cooperad. Then it is a coalgebra over the reduced tree comonad and the structural map
\[
\Delta_\C: \C \longrightarrow \overline{\mathscr{T}}^c(\C)
\]
can be projected onto the sub-$\mathbb{S}$-module $\C \circ \C \hookrightarrow \overline{\mathscr{T}}^c(\C)$. The axioms of a partial cooperad ensure that this projection is coassociative. Finally, $(\C^{u}, \Delta, \epsilon, \nu)$ forms a coaugmented cooperad since the counit is added cofreely. 
\end{proof}

It defines a functor 
\[
\mathsf{Conil}: \mathsf{pCoop}^\mathsf{conil} \longrightarrow \mathsf{Coop}
\]
from the category of conilpotent partial cooperads to the category of cooperads defined as comonoids.

\begin{Definition}[Conilpotent cooperad]\label{def: conilpotent cooperad}
Let $(\C,\Delta,\epsilon)$ be a cooperad. It is said to be \textit{conilpotent} if it is in the essential image of the functor $\mathsf{Conil}$ defined above. 
\end{Definition}

\begin{Remark}
Since $(\C, \Delta, \epsilon, \nu)$ is a cooperad, the two leveled decompositions of its associated partial cooperad are always finite. Hence there cannot be any \textit{horizontal obstruction} to its conilpotency. 
\end{Remark}

This notion of conilpotency for coaugmented cooperads is not equivalent to the notion of a conilpotent cooperad given in \cite{LodayVallette12}, but it is equivalent to the notion of "locally conilpotent cooperad" given in \cite[Section 2.3]{grignoulejay18}. Their coradical filtration for coaugmented cooperads is equal to the one defined here for partial cooperads.

\section{Coalgebras over operads and algebras over cooperads}\label{Section: coalgebras and algebras}
If one wants to encode non-necessarily conilpotent types of coalgebras, then one needs to use operads in order to encode them. This point of view can be seen as a \textit{properadic} point of view, where coalgebras are coded by a reversed operad. These are "Koszul dual" to the new notion of an \textit{algebra over a cooperad}. In this section, we review the main definitions from \cite{grignoulejay18}. 

\medskip

The \textit{dual Schur realization} functor $\widehat{\mathscr{S}}^c$ is given by

\[
\begin{tikzcd}[column sep=4pc,row sep=0.5pc]
\widehat{\mathscr{S}}^c : \smod^{\mathsf{op}} \arrow[r]
&\mathsf{End}(\kmod) \\
M \arrow[r,mapsto]
&\widehat{\mathscr{S}}^c(M)(-) \coloneqq \displaystyle \prod_{n \geq 0} \mathrm{Hom}_{\mathbb{S}_n}(M(n),(-)^{\otimes n})~,
\end{tikzcd}
\]
where $\widehat{\mathscr{S}}^c(M)$ is the dual Schur functor associated to $M$.

\begin{Proposition}[{\cite[Corollary~3.4]{grignoulejay18}}]\label{prop: dual Schur lax contravariant functor}
The contravariant functor $\widehat{\mathscr{S}}^c(-)$ is can be endowed with a lax monoidal functor structure, meaning there a morphism 

\[
\varphi_{M,N}: \widehat{\mathscr{S}}^c(M) \circ \widehat{\mathscr{S}}^c(N) \rightarrowtail \widehat{\mathscr{S}}^c(M \circ N)
\]

that is natural in $M,N$ and satisfies the associativity and unitary conditions. Furthermore, $\varphi_{M,N}$ is a monomorphism for all $M$ and $N$.
\end{Proposition}

\begin{proof}
Let $N$ be an $\mathbb{S}$-module and let $R$ be a $\mathbb{K}$-module, there is a map

\[
\begin{tikzcd}
\left( \displaystyle \prod_{n \geq 0} \mathrm{Hom}_{\mathbb{S}_n}(N(n), R^{\otimes n})\right)^{\otimes k} \arrow[d] \\
\displaystyle \prod_{n \geq 0} \left(\bigoplus_{i_1 + \cdots + i_k = n} \left( \mathrm{Hom}_{{\mathbb{S}_{i_1}}}(N(i_1), R^{\otimes i_1}) \otimes \cdots \otimes \mathrm{Hom}_{\mathbb{S}_{i_k}}(N(i_k), R^{\otimes i_k})\right) \right)
\end{tikzcd}
\]

which is $\mathbb{S}_k$-equivariant and natural in $N$ and $R$. One composes it with

\[
\begin{tikzcd}
\displaystyle \prod_{n \geq 0} \left(\bigoplus_{i_1 + \cdots + i_k = n} \left( \mathrm{Hom}_{{\mathbb{S}_{i_1}}}(N(i_1), R^{\otimes i_1}) \otimes \cdots \otimes \mathrm{Hom}_{\mathbb{S}_{i_k}}(N(i_k), R^{\otimes i_k})\right) \right) \arrow[d] \\
\displaystyle \prod_{n \geq 0} \left( \bigoplus_{i_1 + \cdots + i_k = n} \left( \mathrm{Hom}_{{\mathbb{S}_{i_1} \times \cdots \times \mathbb{S}_{i_k}}}(N(i_1) \otimes \cdots \otimes N(i_k), R^{\otimes n})\right)\right)
\end{tikzcd}
\]

which is also $\mathbb{S}_k$-equivariant and natural in $N$ and $R$. Hence we get a map

\[
\begin{tikzcd}
\displaystyle \prod_{k \geq 0} \mathrm{Hom}_{\mathbb{S}_k} \left( M(k), \left( \prod_{n \geq 0} \mathrm{Hom}_{\mathbb{S}_n}(N(n), R^{\otimes n}) \right)^{\otimes k} \right) \arrow[d] \\
\displaystyle \prod_{k \geq 0} \mathrm{Hom}_{\mathbb{S}_k} \left( M(k), \prod_{n \geq 0} \left( \bigoplus_{i_1 + \cdots + i_k = n} \mathrm{Hom}_{ \mathbb{S}_{i_1} \times \cdots \times \mathbb{S}_{i_k} }(N(i_1) \otimes \cdots \otimes N(i_k), R^{\otimes n}) \right) \right) \arrow[d,"\cong"] \\
\displaystyle \prod_{n \geq 0} \mathrm{Hom}_{\mathbb{S}_n} \left(\bigoplus_{k \geq 0} M(k) \otimes_{\mathbb{S}_k} \left(\bigoplus_{i_1 + \cdots + i_k = n} \mathsf{Ind}_{\mathbb{S}_{i_1} \times \cdots \times \mathbb{S}_{i_k}}^{\mathbb{S}_n}(N(i_1) \otimes \cdots \otimes N(i_k)\right), R^{\otimes n}\right)
\end{tikzcd}
\]

which is natural in $M,N$ and $R$. The last isomorphism is obtained by using the universal properties of the products and the coproducts, as well as the adjunctions between tensor-hom and between induction-restriction. This gives the desired natural transformation 
\[
\varphi_{M,N}: \widehat{\mathscr{S}}^c(M) \circ \widehat{\mathscr{S}}^c(N) \longrightarrow \widehat{\mathscr{S}}^c(M \circ N)~.
\]

One can check that $\varphi_{M,N}$ satisfies the compatibility conditions of a lax monoidal natural transformation.
\end{proof}

\begin{Corollary}
Let $(\C,\Delta,\epsilon)$ be a cooperad. The functor $\widehat{\mathscr{S}}^c(\C)$ admits a canonical monad structure, where the structural map is given by:
\[
\begin{tikzcd}[column sep=4pc,row sep=0.5pc]
\widehat{\mathscr{S}}^c(\C) \circ \widehat{\mathscr{S}}^c(\C) \arrow[r,"\varphi_{\C,\C}"]
&\widehat{\mathscr{S}}^c(\C \circ \C) \arrow[r,"\widehat{\mathscr{S}}^c(\Delta)"]
&\widehat{\mathscr{S}}^c(\C)~.
\end{tikzcd}
\]
\end{Corollary} 

\begin{proof}
Direct consequence of the previous proposition.
\end{proof}

\begin{Definition}[$\C$-algebra]\label{def: C algebra}
Let $(\C, \Delta, \epsilon)$ be a cooperad. A $\mathcal{C}$\textit{-algebra} $B$ amounts to the data $(B,\gamma_B)$ of an algebra over the monad $\widehat{\mathscr{S}}^c(\mathcal{C})$. That is, a structural map
\[
\gamma_B: \prod_{n \geq 0} \mathrm{Hom}_{\mathbb{S}_n}(\C(n),B^{\otimes n}) \longrightarrow B~,
\]

such that the following diagram commutes: 

\[
\begin{tikzcd}[column sep=4pc,row sep=3.5pc]
\widehat{\mathscr{S}}^c(\C) \circ \widehat{\mathscr{S}}^c(\C)(B) \arrow[d, "\widehat{\mathscr{S}}^c(\mathrm{id}_\C)(\gamma_B)",swap] \arrow[r,"\varphi_{\C,\C}(B)"]
&\widehat{\mathscr{S}}^c(\C \circ \C)(B) \arrow[r,"\widehat{\mathscr{S}}^c(\Delta)"]
&\widehat{\mathscr{S}}^c(\C)(B) \arrow[d,"\gamma_B "] \\
\widehat{\mathscr{S}}^c(\C)(B) \arrow[rr,"\gamma_B "]
&
&B~.
\end{tikzcd}
\]
\end{Definition}

\begin{Remark}
The notion of an algebra over a cooperad defines a new type of algebraic structures. The reason is that the structural map
\[
\gamma_B: \displaystyle \prod_{n \geq 0} \mathrm{Hom}_{\mathbb{S}_n}(\mathcal{C}(n),B^{\otimes n}) \longrightarrow B 
\]
associates a element in $B$ to any infinite series of operations. Thus algebras over a cooperad are endowed with a notion of infinite summation without presupposing any underlying topology. 
\end{Remark}

\begin{Proposition}
Let $(\C, \Delta, \epsilon)$ be a cooperad. The category of $\C$-algebras is presentable.
\end{Proposition}

\begin{proof}
The dual Schur functor can also be shown to be accessible, therefore the result follows.
\end{proof}

The notion of an algebra over a cooperad admits a further description in the case where the cooperad is \textit{conilpotent}. From now on, suppose that $(\C,\Delta,\epsilon, \nu)$ is a conilpotent cooperad as defined in Definition \ref{def: conilpotent cooperad}. Each $\mathscr{R}_\omega \C$ defines a sub-cooperad, and there is a short exact sequence of $\mathbb{S}$-modules 

\[
\begin{tikzcd}
0 \arrow[r]
&\mathscr{R}_\omega \C \arrow[r,"\iota_\omega",hook]
&\C \arrow[r,"\pi_\omega"]
&\C / \mathscr{R}_\omega \C \arrow[r]
&0
\end{tikzcd}
\]

for every $\omega \geq 0$. It induces a short exact sequence 

\[
\begin{tikzcd}[column sep=4pc,row sep=0.5pc]
0 \arrow[r]
&\widehat{\mathscr{S}}^c(\C / \mathscr{R}_\omega \C)(V) \arrow[r,"\widehat{\mathscr{S}}^c(\pi_\omega)",hook]
&\widehat{\mathscr{S}}^c(\C)(V) \arrow[r,"\widehat{\mathscr{S}}^c(\iota_\omega)"]
&\widehat{\mathscr{S}}^c(\mathscr{R}_\omega \C)(V) \arrow[r]
&0
\end{tikzcd}
\]

for all $\kk$-modules $V$.

\begin{Definition}[Canonical filtration on a $\mathcal{C}$-algebra]\label{def: canonical filtration C alg}
Let $(\C,\Delta,\epsilon, \nu)$ be a conilpotent cooperad and let $(B, \gamma_B)$ be a $\C$-algebra. The \textit{canonical filtration} of $B$ is the decreasing filtration given by 
\[ 
\mathrm{W}_\omega B \coloneqq \mathrm{Im}\left(\gamma_B \circ \widehat{\mathscr{S}}^c(\pi_\omega)(\mathrm{id}_B): \widehat{\mathscr{S}}^c(\C / \mathscr{R}_\omega \C)(B) \longrightarrow B \right)
\]
where $\mathscr{R}_\omega \C$ denotes the $\omega$-th term of the coradical filtration, for all $\omega \geq 0~.$ Notice that we have 
\[
B = \mathrm{W}_0 B \supseteq \mathrm{W}_1 B \supseteq \mathrm{W}_2 B \supseteq \cdots \supseteq \mathrm{W}_\omega B \supseteq \cdots.
\]
\end{Definition}

\begin{Definition}[Completion of a $\C$-algebra]
Let $(B, \gamma_B)$ be a $\C$-algebra. Its \textit{completion} is given by 
\[ 
\widehat{B} \coloneqq \lim_{\omega \geq 0} B/\mathrm{W}_\omega B~,
\]
where the limit is taken in the category of $\C$-algebras.
\end{Definition}

It comes equipped with a canonical morphism of $\C$-algebras $\varphi_B: B \longrightarrow \widehat{B}~.$ 

\begin{Definition}[Complete $\mathcal{C}$-algebra]
The $\C$-algebra $(B,\gamma_B)$ is said to be \textit{complete} if $\varphi_B$ is an isomorphism of $\C$-algebras. 
\end{Definition}

\begin{Proposition}[{\cite[Proposition~4.24]{grignoulejay18}}]\label{prop: varphi is an epi}
Let $(B, \gamma_B)$ be a $\C$-algebra. The canonical morphism 
\[
\varphi_B: B \longrightarrow \widehat{B}
\]
is an epimorphism.
\end{Proposition}
 
\begin{Remark}\label{Remark: varphi is an epimorphism}
Conceptually, this comes from the fact that any $\mathcal{C}$-algebras already carries a meaningful notion of infinite summation. Thus "nothing is added" when one applies the completion functor. On the other hand, the topology induced by the canonical filtration of a $\mathcal{C}$-algebra might not be Hausdorff. Meaning that the canonical morphism $\varphi_B$ might not be a monomorphism. This completion functor should be considered a \textit{sifting functor}.
\end{Remark}

\begin{Proposition}\label{prop: complete C algebras are a reflexiv subcat}
Let $(\C,\Delta,\epsilon, \nu)$ be a conilpotent cooperad. Any free $\C$-algebra is complete. Furthermore, the category of complete $\C$-algebras forms a reflective subcategory of the category of $\C$-algebras, where the reflector is given by the completion.
\end{Proposition}

\begin{Remark}
Contrary to $\mathcal{C}$-coalgebras, which are always conilpotent, there are examples of $\mathcal{C}$-algebras which are not complete. See for instance \cite[Section 4.5]{grignoulejay18}.
\end{Remark}

Given an operad $\PP$, the dual Schur functor $\widehat{\mathscr{S}}^c(\PP)$ does not inherit a comonad structure since the dual Schur realization functor is only lax monoidal. Nevertheless, one can still define $\PP$-coalgebras as coalgebras over the associated \textit{lax comonad}. This notion, in principle, fails to have all the nice properties that coalgebras over comonads have. 

\begin{Definition}[$\mathcal{P}$-coalgebra]\label{def: P coalgebra}
A $\mathcal{P}$-coalgebra $D$ amounts to the data $(D, \Delta_D)$ of a $\kk$-module $D$ endowed with a structural map
\[
\Delta_D: D \longrightarrow \displaystyle \prod_{n \geq 0} \mathrm{Hom}_{\mathbb{S}_n}(\mathcal{P}(n),D^{\otimes n})~,
\]
such that the following diagram commutes 
\[
\begin{tikzcd}[column sep=4.5pc,row sep=3pc]
D \arrow[r,"\Delta_D"] \arrow[d,"\Delta_D",swap] 
&\widehat{\mathscr{S}}^c(\PP)(D) \arrow[r,"\widehat{\mathscr{S}}^c(\mathrm{id})(\Delta_D)"]
&\widehat{\mathscr{S}}^c(\PP)(D) \circ \widehat{\mathscr{S}}^c(\PP)(D) \arrow[d,"\varphi_{\PP,\PP}(D)"] \\
\widehat{\mathscr{S}}^c(\PP)(D) \arrow[rr,"\widehat{\mathscr{S}}^c(\gamma)(\mathrm{id})"]
&
&\widehat{\mathscr{S}}^c(\PP \circ \PP)(D)~.
\end{tikzcd}
\]
\end{Definition}

\begin{Definition}[Coendomorphism operad]
Let $D$ be a $\mathbb{K}$-module. The \textit{coendomorphism operad} of $D$ is defined by the $\mathbb{S}$-module:
\[
\mathrm{Coend}_D(n) \coloneqq \text{Hom}(D,D^{\otimes n})~,
\]
where the composition maps are given by the composition of functions and the unit is given by $\mathrm{id}_D$. 
\end{Definition}

\begin{Proposition}\label{Equivalence coalg and coend}
Let $D$ be a $\mathbb{K}$-module. A $\PP$-coalgebra structure $\Delta_D$ on $D$ is equivalent to a morphism of operads $\delta_D: \PP \longrightarrow \mathrm{Coend}_D$.
\end{Proposition}

\begin{proof}
The structural morphism 
\[
\Delta_D: D \longrightarrow \prod_{n \geq 0} \mathrm{Hom}_{\mathbb{S}_n}(\PP(n),D^{\otimes n})
\]
can be decomposed into its components $\Delta_D^n: D \longrightarrow \mathrm{Hom}_{\mathbb{S}_n}(\PP(n),D^{\otimes n})$, which we can be transposed under the $\otimes_{\mathbb{S}_{n}}$ - $\mathrm{Hom}_{\mathbb{S}_{n}}$ adjunction into $(\delta_D^n)^t: D \otimes_{\mathbb{S}_{n}} \PP(n) \longrightarrow D^{\otimes n}$. Since the action of $\mathbb{S}_n$ is trivial on $D$, the maps $(\delta_D^n)^t$ are equivalent to $\mathbb{S}_n$-equivariant maps $\delta_D^n: \PP(n) \longrightarrow \mathrm{Hom}(D,D^{\otimes n})$. This forms morphism of $\mathbb{S}$-modules $\delta: \PP \longrightarrow \mathrm{Coend}_D$. It is straightforward to check that $\Delta_D$ endows $D$ with a $\PP$-coalgebra structure if and only if $\delta: \PP \longrightarrow \mathrm{Coend}_D$ is a morphism of operads. 
\end{proof}

\begin{Remark}
Since coalgebras over an operad $\PP$ are given by morphisms of operads to the coendomorphism operad, one can view $\PP$ as a reversed operad in the category of properads when encoding its coalgebras.
\end{Remark}

\begin{theorem}[{\cite[Theorem 2.7.11]{anelcofree2014}}]
Let $(\PP,\gamma,\eta)$ be an operad. The category of $\PP$-coalgebras is comonadic, i.e: there exists a comonad $(\mathscr{C}(\PP), \omega, \zeta)$ in the category of $\kk$-modules such that the category of $\mathscr{C}(\PP)$-coalgebras is equivalent to the category of $\PP$-coalgebras.
\end{theorem}

In particular, this entails the existence of a cofree $\PP$-coalgebra. While in the general setting of \cite{anelcofree2014}, the construction of the comonad $\mathscr{C}(\PP)$ is given by an infinite recursion, the construction of $\mathscr{C}(\PP)$ in the category of $\kk$-modules stops at the first step. 

\begin{theorem}[{\cite[Theorem 3.3.1]{anelcofree2014}}]\label{thm: existence of the cofree P cog}
The endofunctor $\mathscr{C}(\PP)$ is given by the pullback
\[
\begin{tikzcd}[column sep=3pc,row sep=3pc]
\mathscr{C}(\PP) \arrow[r,"p_2"] \arrow[d,"p_1",swap,rightarrowtail] \arrow[dr, phantom, "\ulcorner", very near start]
&\widehat{\mathscr{S}}^c(\PP) \circ \widehat{\mathscr{S}}^c(\PP) \arrow[d,"\varphi_{\PP,\PP}",rightarrowtail] \\
\widehat{\mathscr{S}}^c(\PP) \arrow[r,"\widehat{\mathscr{S}}^c(\gamma)"]
&\widehat{\mathscr{S}}^c(\PP \circ \PP)
\end{tikzcd}
\]
in the category $\mathsf{End}(\kmod)$. For $\kk$-module $V$ it gives a pullback
\[
\begin{tikzcd}[column sep=3pc,row sep=3pc]
\mathscr{C}(\PP)(V) \arrow[r,"p_2(V)"] \arrow[d,"p_1(V)",swap,rightarrowtail] \arrow[dr, phantom, "\ulcorner", very near start]
&\widehat{\mathscr{S}}^c(\PP)(V) \circ \widehat{\mathscr{S}}^c(\PP)(V) \arrow[d,"\varphi_{\PP,\PP}(V)",rightarrowtail] \\
\widehat{\mathscr{S}}^c(\PP)(V) \arrow[r,"\widehat{\mathscr{S}}^c(\gamma)(V)"]
&\widehat{\mathscr{S}}^c(\PP \circ \PP)(V)
\end{tikzcd}
\]
in the category of $\kk$-modules. Here $p_1$ is a degree-wise monomorphism since $\varphi_{\PP,\PP}(V)$ is a degree-wise monomorphism. The structural map of the comonad 
\[
\omega(V): \mathscr{C}(\PP)(V) \longrightarrow \mathscr{C}(\PP) \circ \mathscr{C}(\PP)(V)
\]
is given by the map $p_2$ in the previous pullback. The counit is given of the comonad $\mathscr{C}(\PP)$ is given by
\[
\begin{tikzcd}[column sep=4pc,row sep=0.5pc]
\xi(V): \mathscr{C}(\PP)(V) \arrow[r,"p_1(V)"] 
&\widehat{\mathscr{S}}^c(\PP)(V) \arrow[r,"\widehat{\mathscr{S}}^c(\eta)(\mathrm{id}_V)"]
&V~.
\end{tikzcd}
\]
\end{theorem}

\begin{Remark}
The subspace $\mathscr{C}(\PP)(V)$ of $\widehat{\mathscr{S}}^c(\PP)(V)$ admits an explicit description in terms of \textit{representative functions}. See \cite[Section 3.1]{anelcofree2014} or \cite{blockleroux} for the original reference about representative functions in the case of coassociative and cocommutative coalgebras.
\end{Remark}

\begin{Corollary}
Let $\PP$ be an operad, then for any module $V$, the cofree $\PP$-coalgebra on $V$ is given by $\mathscr{C}(\PP)(V)$ together with the map: 
\[
\begin{tikzcd}[column sep=4pc,row sep=0.5pc]
\omega_V: \mathscr{C}(\PP)(V) \arrow[r,"\omega(V)"] 
&\mathscr{C}(\PP) \circ \mathscr{C}(\PP)(V) \arrow[r,"p_1 \circ \hspace{1pt} \mathrm{id} "]
&\widehat{\mathscr{S}}^c(\PP) \circ \mathscr{C}(\PP)(V)~,
\end{tikzcd}
\]
where $\omega$ is the structural map of the comonad $\mathscr{C}(\PP)$.
\end{Corollary}

\section{Twisting morphisms and Bar-Cobar adjunctions}\label{Section: Complete Bar-Cobar}
In this section, we recall the definition of a twisting morphism $\alpha: \C \longrightarrow \PP$ between a conilpotent partial cooperad and a partial operad, which is a Maurer-Cartan element in the convolution pre-Lie algebra. From such a twisting morphism, one can construct a Bar-Cobar adjunction $\Omega_{\alpha} \dashv \text{B}_{\alpha}$ relative to $\alpha$ between the category of dg $\PP$-algebras and the category of dg $\C$-coalgebras. This adjunction has been a very powerful tool to study the homotopy theory of dg algebras encoded by an operad. Indeed, using this adjunction, one can transfer the model structure of dg $\PP$-algebras where weak-equivalences are given by quasi-isomorphisms onto the category of dg $\C$-coalgebras. If $\alpha$ is a Koszul twisting morphism, this adjunction is a Quillen equivalence. This approach, first developed by V. Hinich for dg Lie algebras in \cite{Hinich01} and by K. Lefèvre-Hasegawa in \cite{LefevreHasegawa03} for dg associative algebras, was extended to all types of dg $\mathcal{P}$-algebras in \cite{Vallette14}. 

\medskip

As we explained in the previous sections, one cannot encode non-conilpotent types of coalgebras using cooperads, one needs to use operads to encode them. But given a category of dg $\PP$-coalgebras, how can we study their homotopy theory? An answer to this question was provided by Brice Le Grignou and Damien Lejay in \cite{grignoulejay18}. Using again a twisting morphism $\alpha$, the authors of \textit{loc.cit} construct a \textit{complete Bar-Cobar adjunction relative to} $\alpha$ between the category of dg $\PP$-coalgebras and the category of dg $\C$-algebras. And they show that one can transfer (when it exists) the model category structure of dg $\PP$-coalgebras where weak-equivalences are given by quasi-isomorphisms to the category of dg $\C$-algebras. Furthermore, when we restrict to the case of the twisting morphism $\iota: \C \longrightarrow \Omega \C$, this adjunction is a Quillen equivalence. This approach allows us to study the homotopy theory of these non-necessarily conilpotent coalgebras using these new algebraic objects that are algebras over cooperads.

\begin{Remark}
The results concerning the complete Bar-Cobar adjunction can be found in \cite{grignoulejay18}, but their notations are drastically different. For isntance, the aforementioned adjunction is called "Cobar and its right adjoint" in \textit{loc.cit}.
\end{Remark}

From now on, the base category will be the category of differential graded (dg) $\kk$-modules. Everything stated so far in this chapter can be generalized \textit{mutatis mutandis} to the category of dg modules. 

\begin{Definition}[Convolution partial operad]
Let $(\C,\{\Delta_i\},d_\C)$ be a dg partial cooperad and let $(\mathcal{P},\{\circ_i\},d_\PP)$ be a dg partial operad. The \textit{convolution partial operad} of $\C$ and $\PP$ is given by the dg $\mathbb{S}$-module 
\[
\mathcal{H}om(\C,\PP)(n) \coloneqq \mathrm{Hom}_{\mathsf{gr}~\mathsf{mod}}(\C(n), \PP(n))~,
\]
with its natural $\mathbb{S}_n$-action, where the differential is given by $\partial(f) = d_\PP \circ f - (-1)^{|f|} f \circ d_\C$. It can be endowed with partial composition maps given by
\[
\begin{tikzcd}[column sep=3pc,row sep=0pc]
\alpha \circ_i \beta : \C(n+k-1) \arrow[r,"\Delta_{i}"]
&\C(n) \otimes \C(k) \arrow[r,"\alpha \hspace{1pt} \otimes \hspace{1pt} \beta"]
&\PP(n) \otimes \PP(k) \arrow[r,"\circ_i"]
&\PP(n+k-1)~.
\end{tikzcd}
\]
\end{Definition}

\begin{Definition}[Totalization of a partial operad]
Let $(\PP,\{\circ_i\},d_\PP)$ be a dg partial operad, the \textit{totalization} of $\PP$ given by 
\[
\prod_{n \geq 0} \PP(n)^{\mathbb{S}_n}~.
\]
It can be endowed with a dg pre-Lie algebra structure by setting
\[ 
\mu \star \nu \coloneqq \sum_{i=1}^{n} \sum_{\sigma \in \mathbb{S}_k} (\mu \circ_i \nu)^{\sigma^\flat}~,
\]
for $\mu$ in $\PP(n)$ and $\nu$ in $\PP(k)$, where $\sigma^\flat$ is the unique permutation in $\mathbb{S}_{n+k-1}$ that acts as $\sigma$ on $\{i,\cdots,i+n-\}$ and as the identity elsewhere. 
\end{Definition}

\begin{Remark}
One can check by direct computation that the axioms of a partial operad make the associator of $\star$ right symmetric. See 
\end{Remark}

\begin{Definition}[Twisting morphism]
Let $(\C,\{\Delta_i\},d_\C)$ be a dg partial cooperad and let $(\mathcal{P},\{\circ_i\},d_\PP)$ be a dg partial operad. Let 
\[
\mathfrak{g}_{\C,\PP} \coloneqq \prod_{n \geq 0} \mathrm{Hom}_{\mathbb{S}_n}(\C(n), \PP(n))
\]
be the pre-Lie algebra given by the totalization of the convolution partial operad of $\C$ and $\PP$. A \textit{twisting morphism} is a Maurer-Cartan element $\alpha$ of $\mathfrak{g}_{\C,\PP}$, that is, a morphism $\alpha: \C \longrightarrow \PP$ graded $\mathbb{S}$-modules of degree $-1 $ satisfying :
\[
\partial(\alpha) + \alpha \star \alpha = 0~.
\]
The set of twisting morphisms between $\C$ and $\PP$ is denoted by $\mathrm{Tw}(\C,\PP)~.$
\end{Definition}

For the rest of this section, let us fix a dg partial operad $(\mathcal{P},\{\circ_i\},d_\PP)$, a conilpotent dg partial cooperad $(\C,\{\Delta_i\},d_\C)$, and a twisting morphism $\alpha: \C \longrightarrow \PP$.

\input cyracc.def
\font\tencyr=wncysc10
\def\cyr{\tencyr\cyracc}
\def\diracComb{\mbox{\cyr SH}}

\begin{Notation}\label{not: le dirac}
Let $f: X \longrightarrow Y$ be a morphism of graded modules and $g: X \longrightarrow Y$ be a map of degree $p$. We denote
\[
\diracComb_n(f,g) \coloneqq \sum_{i=0}^n f^{\otimes i-1} \otimes g \otimes f^{\otimes n-i} : X^{\otimes n} \longrightarrow Y^{\otimes n}
\]
the resulting $\mathbb{S}_n$-equivariant map of degree $p$. Let $M$ be a graded $\mathbb{S}$-module. It induces first a map of degree $p$
\[
\bigoplus_{n \geq 0} M(n) \otimes_{\mathbb{S}_n} X^{\otimes n} \longrightarrow \bigoplus_{n \geq 0} M(n) \otimes_{\mathbb{S}_n} Y^{\otimes n}
\]
by applying $\mathrm{id}_M \otimes \diracComb_n(f,g)$ at each arity. By a slight abuse of notation, this map will be denoted by $\mathscr{S}(\mathrm{id}_M)(\diracComb(f,g))~.$ Likewise, it induces a second map of degree $p$
\[
\prod_{n \geq 0} \mathrm{Hom}_{\mathbb{S}_n}(M(n),X^{\otimes n}) \longrightarrow \prod_{n \geq 0} \mathrm{Hom}_{\mathbb{S}_n}(M(n),Y^{\otimes n})
\]
by applying $\mathrm{Hom}(\mathrm{id}_M,\diracComb_n(f,g))$ at each arity. By a slight abuse of notation, this map will be denoted by $\widehat{\mathscr{S}}^c(\mathrm{id}_M)(\diracComb(f,g))~.$
\end{Notation}

Let $\PP^{u}$ the dg operad given by the unital extension of $\PP$ and let $\C^{u}$ be the conilpotent dg cooperad given by the counital extension of $\C$. From this data we can construct two adjunctions as follows. 

\medskip

\subsection{First Bar-Cobar adjunction.} The classical Bar-Cobar adjunction relative to $\alpha$ is given by a pair of adjoint functors

\[
\begin{tikzcd}[column sep=5pc,row sep=3pc]
          \mathsf{dg}~\mathcal{C}^{u}\text{-}\mathsf{coalg} \arrow[r, shift left=1.1ex, "\Omega_{\alpha}"{name=F}] & \mathsf{dg}~\mathcal{P}^{u}\text{-}\mathsf{alg}, \arrow[l, shift left=.75ex, "\mathrm{B}_{\alpha}"{name=U}]
            \arrow[phantom, from=F, to=U, , "\dashv" rotate=-90]
\end{tikzcd}
\]

between the category of dg $\PP^{u}$-algebras and the category of dg $\C^{u}$-coalgebras.

\begin{Definition}[Bar construction relative to $\alpha$]
Let $(A,\gamma_A,d_A)$ be a dg $\PP^{u}$-algebra. The \textit{Bar construction relative to} $\alpha$ of $A$ is given by
\[
\mathrm{B}_{\alpha}A \coloneqq \left(\mathscr{S}(\C^{u})(A), d_{\text{bar}} \coloneqq d_1 + d_2 \right)~,
\]
where $\mathscr{S}(\C^{u})(A)$ denotes the cofree $\C^{u}$-coalgebra generated by $A$. The differential $d_{\text{bar}}$ is given as the sum of two terms $d_1$ and $d_2$. The term $d_1$ is given by 
\[
d_1 \coloneqq \mathscr{S}(d_\C)(\mathrm{id}_A) + \mathscr{S}(\mathrm{id}_\C)(\diracComb(\mathrm{id}_A,d_A))~.
\]
The term $d_2$ is the unique coderivation extending
\[
\begin{tikzcd}[column sep=5pc,row sep=1pc]
\mathscr{S}(\C^{u})(A) \arrow[r,"\mathscr{S}(\alpha)(\mathrm{id}_A)"]
&\mathscr{S}(\PP^{u})(A) \arrow[r,"\gamma_A"]
&A~.
\end{tikzcd}
\]
\end{Definition}

\begin{Proposition}
Let $(A,\gamma_A,d_A)$ be a dg $\PP^{u}$-algebra. The Bar construction $\mathrm{B}_{\alpha}A$ forms a dg $\C^{u}$-coalgebra, and it defines a functor

\[
\mathrm{B}_{\alpha}: \mathsf{dg}~\PP^{u}\text{-}\mathsf{alg} \longrightarrow \mathsf{dg}~\C^{u}\text{-}\mathsf{coalg}~.
\]
\end{Proposition}

\begin{Definition}[Cobar construction relative to $\alpha$]
Let $(C,\Delta_C,d_C)$ be a dg $\C^{u}$-coalgebra. The \textit{Cobar construction relative to} $\alpha$ of $C$ is given by
\[
\Omega_\alpha C \coloneqq \left(\mathscr{S}(\PP^{u})(C), d_{\text{cobar}} \coloneqq d_1 - d_2 \right)~,
\]
where $\mathscr{S}(\PP^{u})(C)$ denotes the free $\PP^{u}$-algebra generated by $C$. The differential $d_{\text{cobar}}$ is given as the difference of two terms $d_1$ and $d_2$. The term $d_1$ is given by 
\[
d_1 \coloneqq - \mathscr{S}(d_\PP)(\mathrm{id}_C) + \mathscr{S}(\mathrm{id}_\PP)(\diracComb(\mathrm{id}_C,d_C))~.
\]
The term $d_2$ is the unique coderivation extending
\[
\begin{tikzcd}[column sep=5pc,row sep=1pc]
C \arrow[r," \Delta_C"]
&\mathscr{S}(\C^{u})(C) \arrow[r,"\mathscr{S}(\alpha)(\mathrm{id}_C)"]
&\mathscr{S}(\PP^{u})(C)~.
\end{tikzcd}
\]
\end{Definition}

\begin{Proposition}
Let $(C,\Delta_C,d_C)$ be a dg $\C^{u}$-coalgebra. The Cobar construction $\Omega_\alpha C$ forms a dg $\PP^{u}$-algebra, and it defines a functor

\[
\Omega_\alpha: \mathsf{dg}~\C^{u}\text{-}\mathsf{coalg} \longrightarrow \mathsf{dg}~\PP^{u}\text{-}\mathsf{alg}~.
\]
\end{Proposition}

\begin{Definition}[Twisting morphism relative to $\alpha$]
Let $(A,\gamma_A,d_A)$ be a dg $\PP^{u}$-algebra and let $(C,\Delta_C,d_C)$ be a dg $\C^{u}$-coalgebra. A graded morphism 
\[
\nu: C \longrightarrow A
\]

is said to be a twisting morphism relative to $\alpha$ if it satisfies the following equation 

\[
\gamma_A \cdot \mathscr{S}(\alpha)(\nu) \cdot \Delta_C + \partial(\nu) = 0~.
\]

The set of curved twisting morphisms relative to $\alpha$ are denoted by $\mathrm{Tw}^{\alpha}(D,B)$.
\end{Definition}

\begin{Proposition}\label{prop: Bar-Cobar adjunction relative to alpha}
There are bijections

\[
\mathrm{Hom}_{\mathsf{dg}~\PP\text{-}\mathsf{alg}}\left(\Omega_{\alpha}C, A\right) \cong  \mathrm{Tw}^{\alpha}(C,A)  \cong \mathrm{Hom}_{\mathsf{dg}~\C\text{-}\mathsf{coalg}} \left(C, \mathrm{B}_{\alpha}A \right)~,
\]
\vspace{0.1pc}

which are natural in $C$ and $A$. 
\end{Proposition}

\begin{proof}
There natural isomorphism:

\[
\mathrm{Hom}_{\mathsf{gr}~\PP^{u}\text{-}\mathsf{alg}}(\mathscr{S}(\PP^{u})(C),A) \cong \mathrm{Hom}_{\mathsf{gr}~\mathsf{mod}}(C,A) \cong \mathrm{Hom}_{\mathsf{gr}~\C^{u}\text{-}\mathsf{coalg}^{\mathsf{conil}}}(C,\mathscr{S}(\C^{u})(A))~,
\]
\vspace{0.1pc}

given by the fact that $\mathscr{S}(\PP^{u})(C)$ is the free $\PP^{u}$-algebra generated by $C$ and that $\mathscr{S}(\C^{u})$ is the cofree conilpotent $\C^{u}$-coalgebra generated by $A$. One can check that a graded morphism $f: C \longrightarrow A$ induces a morphism of \textit{dg} $\PP^{u}$-algebras on the left hand side if and only if it is a twisting morphism relative to $\alpha$, if and only if it induces a morphism of \textit{dg} conilpotent $\C^{u}$-coalgebras on the right hand side. See \cite[Section 11]{LodayVallette12} for more details.
\end{proof}

Over a field of characteristic zero, the category of dg $\PP^{u}$-algebras always admits a model category structure which is right-transferred along the free-forgetful adjunction, where weak equivalences are given by quasi-isomorphisms

\begin{theorem}[\cite{Hinich97}]\label{thm: model structure on P-alg}
There is a model structure on the category of dg $\mathcal{P}^{u}$-algebras right-transferred along the free-forgetful adjunction
\[
\begin{tikzcd}[column sep=7pc,row sep=3pc]
\mathsf{dg}\textsf{-}\mathsf{mod} \arrow[r, shift left=1.1ex, "\mathscr{S}(\mathcal{P}^{u})(-)"{name=F}]      
&\mathsf{dg}~\PP^{u}\textsf{-}\mathsf{alg}~, \arrow[l, shift left=.75ex, "\mathrm{U}"{name=U}]
\arrow[phantom, from=F, to=U, , "\dashv" rotate=-90]
\end{tikzcd}
\]
where 
\begin{enumerate}
\medskip 

\item the class of weak equivalences is given by quasi-isomorphisms,

\medskip 

\item the class of fibrations is given by degree-wise epimorphisms,

\medskip 

\item the class of cofibrations is given by left lifting property with respect to acyclic fibrations.
\end{enumerate}
\end{theorem}

\begin{Remark}
Any morphism of dg operads $f: \mathcal{P} \longrightarrow \mathcal{Q}$ induces a Quillen adjunction between the category of dg $\mathcal{P}$-algebras and the category of dg $\mathcal{Q}$-algebras. If $f$ is a quasi-isomorphism of dg operads, the induced adjunction automatically becomes a Quillen equivalence.
\end{Remark}

This model structure can always be transferred along the Bar-Cobar adjunction relative to $\alpha$ in order to obtain a model structure on the category of dg $\C^{u}$-coalgebras. 

\begin{theorem}[{\cite{Vallette14}}]
There is a model structure on the category of dg $\mathcal{C}^{u}$-coalgebras left-transferred along the complete Bar-Cobar adjunction

\[
\begin{tikzcd}[column sep=7pc,row sep=3pc]
            \mathsf{dg}~\mathcal{P}^{u}\textsf{-}\mathsf{alg} \arrow[r, shift left=1.1ex, "\mathrm{B}_{\alpha}"{name=F}] &\mathsf{dg}~\mathcal{C}^{u}\textsf{-}\mathsf{coalg}~, \arrow[l, shift left=.75ex, "\Omega_{\alpha}"{name=U}]
            \arrow[phantom, from=F, to=U, , "\dashv" rotate=90]
\end{tikzcd}
\]

where 
\begin{enumerate}

\medskip

\item the class of weak equivalences is given by morphisms $f$ such that $\Omega_\alpha(f)$ is a quasi-isomorphism,

\medskip

\item the class of cofibrations is given by degree-wise monomorphisms,

\medskip

\item  and the class of fibrations is given by right lifting property with respect to acyclic cofibrations.
\end{enumerate}
\end{theorem}

\begin{Remark}
Cofibrations should be defined, in principle, by morphisms $f$ such that $\Omega_\alpha(f)$ is a cofibration. Nevertheless, this class of morphisms turns out be equivalent to degree-wise monomorphisms.
\end{Remark}

\begin{theorem}[\cite{Vallette14}]
Let $\alpha$ be a Koszul twisting morphism, then the Bar-Cobar adjunction relative to $\alpha$
\[
\begin{tikzcd}[column sep=7pc,row sep=3pc]
            \mathsf{dg}~\mathcal{P}^{u}\textsf{-}\mathsf{alg} \arrow[r, shift left=1.1ex, "\mathrm{B}_{\alpha}"{name=F}] &\mathsf{dg}~\mathcal{C}^{u}\textsf{-}\mathsf{coalg}~, \arrow[l, shift left=.75ex, "\Omega_{\alpha}"{name=U}]
            \arrow[phantom, from=F, to=U, , "\dashv" rotate=90]
\end{tikzcd}
\]
is a Quillen equivalence. 
\end{theorem}

Another powerful tool to study the homotopy theory of dg $\PP^{u}$-algebra is given by $\infty$-morphisms. These are also in fact given by the Bar-Cobar adjunction relative to $\alpha$. For more details, we refer to \cite{RobertNicoudWierstra17}.

\begin{Definition}[$\infty_\alpha$-morphism of algebras]
Let $A$ and $B$ be two dg $\PP^{u}$-algebras. An $\infty_\alpha$\textit{-morphism} $f: A \rightsquigarrow B$ is the data of a morphism $f: \mathrm{B}_{\alpha}A \longrightarrow \mathrm{B}_{\alpha}B$ of dg $\C^{u}$-coalgebras. 
\end{Definition}

\begin{Example}
Consider the dg operad $\Omega \mathcal{C}om$ which encodes $\mathcal{L}_\infty$-algebras. Morphisms sometimes called "$\mathcal{L}_\infty$-morphisms" or $\infty$-morphisms are in fact $\infty_\iota$-morphisms.
\end{Example}

\begin{Remark}
This data is equivalent to a family of maps $f_n: \C^{u}(n) \otimes_{\mathbb{S}_n} A^{\otimes n} \longrightarrow B$ which satisfies certain relations. Since $\C^{u}(1) = \C(1) \oplus \I$ the morphism $f_1$ has a linear part $f_\I: \I \otimes A \cong A \longrightarrow B$ which is a morphism of dg modules. 
\end{Remark}

Note that to any morphism of dg $\PP$-algebras $f: A \longrightarrow B$, one can associate an $\infty_\alpha$-morphism by setting $f_\I = f$ and zero to all other components.

\begin{Definition}[$\infty_\alpha$-quasi-isomorphisms]
An $\infty_\alpha$-\textit{quasi-isomorphism} is a $\infty_\alpha$-morphism $f: A \rightsquigarrow B$ where $f_\I$ is a quasi-isomorphism of dg modules.
\end{Definition}

\begin{Proposition}
Let $f: A \rightsquigarrow B$ be an $\infty_\alpha$-quasi-isomorphism between two dg $\PP^{u}$-algebras. 

\medskip

\begin{enumerate}

\item The morphism $f: \mathrm{B}_{\alpha}A \longrightarrow \mathrm{B}_{\alpha}B$ is a weak-equivalence of dg $\C^{u}$-coalgebras.

\medskip

\item There exists an $\infty_\alpha$-quasi-isomorphism $g: B \rightsquigarrow A$ such that $g_\I$ is the inverse of $f_\I$ on the level of the homology groups. 
\end{enumerate}
\end{Proposition}

\begin{Corollary}
Let $\alpha$ be a Koszul twisting morphism and let $A$ and $B$ be two dg $\PP^{u}$-algebras. There exists a zig-zag of quasi-isomorphisms 
\[
A \lqi \cdot \qi \cdot \lqi \cdot \qi B 
\]
of dg $\mathcal{P}^{u}$-algebras if and only if there exists an $\infty_\alpha$-quasi-isomorphism $A \rightsquigarrow B$.
\end{Corollary}

These adjunctions are "functorial" in the following sense. See \cite{DrummondColeHirsh14} for a general exposition.

\begin{Proposition}
Let $\mathcal{P}$ and $\mathcal{Q}$ be dg operads, and let $\mathcal{C}$ and $\mathcal{D}$ be conilpotent dg cooperads. 

\medskip

\begin{enumerate}
\item Let $\alpha: \C \longrightarrow \PP$ and $\beta: \mathcal{D} \longrightarrow \mathcal{Q}$ be two twisting morphisms. 

\medskip

\item Let $f: \PP \longrightarrow \mathcal{Q}$ be a morphism of dg operads and let $g: \mathcal{C} \longrightarrow \mathcal{D}$ be a morphism of conilpotent dg cooperads. 
\end{enumerate}

\medskip

Such that the following diagram commutes 

\[
\begin{tikzcd}[column sep=4pc,row sep=4pc]
\C \arrow[r,"\alpha"] \arrow[d,"g",swap] 
&\mathcal{P} \arrow[d,"f"]\\
\mathcal{D} \arrow[r,"\beta"]
&\mathcal{Q}~.
\end{tikzcd}
\]

Then the following square 

\[
\begin{tikzcd}[column sep=5pc,row sep=5pc]
\mathsf{dg}~\mathcal{P}\text{-}\mathsf{alg} \arrow[r,"\mathrm{B}_\alpha"{name=B},shift left=1.1ex] \arrow[d,"\mathrm{Ind}_f "{name=SD},shift left=1.1ex ]
&\mathsf{dg}~\mathcal{C}\text{-}\mathsf{coalg} \arrow[d,"\mathrm{Res}_g"{name=LDC},shift left=1.1ex ] \arrow[l,"\Omega_\alpha"{name=C},,shift left=1.1ex]  \\
\mathsf{dg}~\mathcal{Q}\text{-}\mathsf{alg} \arrow[r,"\mathrm{B}_\beta "{name=CC},shift left=1.1ex]  \arrow[u,"\mathrm{Res}_f"{name=LD},shift left=1.1ex ]
&\mathsf{dg}~\mathcal{D}\text{-}\mathsf{coalg} \arrow[l,"\Omega_\beta"{name=CB},shift left=1.1ex] \arrow[u,"\mathrm{Coind}_g"{name=TD},shift left=1.1ex] \arrow[phantom, from=SD, to=LD, , "\dashv" rotate=180] \arrow[phantom, from=C, to=B, , "\dashv" rotate=90]\arrow[phantom, from=TD, to=LDC, , "\dashv" rotate=180] \arrow[phantom, from=CC, to=CB, , "\dashv" rotate=90]
\end{tikzcd}
\] 

of Quillen adjunctions commutes. 
\end{Proposition}

\begin{Remark}
This theory was generalized by Brice Le Grignou to the case of a curved twisting morphism $\alpha: \C \longrightarrow \PP$ between a non-augmented dg operad $\PP$ and a conilpotent curved cooperad $\C$ in \cite{grignou2019}. See Chapter 2 for more information on these notions. 
\end{Remark}

\subsection{Second Bar-Cobar adjunction.} There is a second adjunction that one can construct using a twisting morphism $\alpha: \C \longrightarrow \PP$. The \textit{complete Bar-Cobar adjunction relative to} $\alpha$ is given by

\[
\begin{tikzcd}[column sep=5pc,row sep=3pc]
            \mathsf{dg}~\PP^{u}\text{-}\mathsf{coalg} \arrow[r, shift left=1.1ex, "\widehat{\Omega}_{\alpha}"{name=F}] & \mathsf{dg}~\C^{u}\text{-}\mathsf{alg}^{\mathsf{comp}}~. \arrow[l, shift left=.75ex, "\widehat{\text{B}}_{\alpha}"{name=U}]
            \arrow[phantom, from=F, to=U, , "\dashv" rotate=-90]
\end{tikzcd}
\]

Recall that the category of dg $\PP^{u}$-coalgebras has cofree objects given by $\mathscr{C}(\PP^{u})(V)$. These objects behave well with respect to coderivations, which will allow us to define the differentials on the Bar constructions. 

\begin{Notation}
Let $X$ be a graded module, we denote $\pi: \mathscr{C}(\PP^{u})(X) \longrightarrow X$ the composition:
\[
\begin{tikzcd}[column sep=4pc,row sep=1pc]
\mathscr{C}(\PP^{u})(X) \arrow[r,"p_1",hook]
&\widehat{\mathscr{S}}^c(\PP^{u})(X) \arrow[r,"\widehat{\mathscr{S}}^c(\eta)(\mathrm{id})"]
&X~.
\end{tikzcd}
\]
\end{Notation}

\begin{lemma}[{\cite[Lemma 6.20]{grignoulejay18}}]
Let $f: X \longrightarrow Y$ be a morphism graded modules of degree $0$ and $g: X \longrightarrow Y$ be a morphism of degree $p$, then the degree $p$ map: 

\[
\begin{tikzcd}[column sep=4pc,row sep=1pc]
&\widehat{\mathscr{S}}^c(\mathrm{id}_{\PP^{u}})(\diracComb(f,g)): \widehat{\mathscr{S}}^c(\PP^{u})(X) \arrow[r]
&\widehat{\mathscr{S}}^c(\PP^{u})(Y)
\end{tikzcd}
\]
\vspace{0.1pc}

restricts to a degree $p$ morphism $\mathscr{C}(\mathrm{id})(\diracComb(f,g)): \mathscr{C}(\PP^{u})(X) \longrightarrow \mathscr{C}(\PP^{u})(Y)$.
\end{lemma}

\begin{proof}
This is a consequence of the universal property of pullbacks: it suffices to construct a map
\[
\mathscr{C}(\PP)(X) \longrightarrow \widehat{\mathscr{S}}^c(\PP) \circ \widehat{\mathscr{S}}^c(\PP)(Y)
\]

that coincides with the morphism $\widehat{\mathscr{S}}^c(\mathrm{id}_{\PP^{u}})(\diracComb(f,g))$ on $\widehat{\mathscr{S}}^c(\PP \circ \PP)(Y)$. One can check that this morphism is given by 

\[
\widehat{\mathscr{S}}^c(\mathrm{id})(\diracComb(\widehat{\mathscr{S}}^c(\mathrm{id})(f), \allowbreak \widehat{\mathscr{S}}^c(\mathrm{id})(\diracComb(f,g)))) \circ p_2~.
\]
\end{proof}

\begin{lemma}[{\cite[Lemmas 6.21 and 6.22]{grignoulejay18}}]\label{lemma: on coderivations of cofree 1}
Let $X$ be a graded module. The map

\[
\widehat{\mathscr{S}}^c(d_\PP)(\mathrm{id}_X): \widehat{\mathscr{S}}^c(\PP^{u})(X) \longrightarrow \widehat{\mathscr{S}}^c(\PP^{u})(X)
\]
\vspace{0.1pc}

restricts to a unique map $\mathscr{C}(d_\PP)(\mathrm{id}): \mathscr{C}(\PP^{u})(X) \longrightarrow \mathscr{C}(\PP^{u})(X)$ of degree $-1$ which is a coderivation with respect to the cofree $\PP^{u}$-coalgebra structure.
\end{lemma}

\begin{proof}
The first part of the proof is analogous to the previous lemma. The map 

\[
\mathscr{C}(\PP)(X) \longrightarrow \widehat{\mathscr{S}}^c(\PP) \circ \widehat{\mathscr{S}}^c(\PP)(X)
\]

is given by the morphism 

\[
\widehat{\mathscr{S}}^c(d_\PP) \circ \widehat{\mathscr{S}}^c(\mathrm{id})(\mathrm{id}) + \widehat{\mathscr{S}}^c(\mathrm{id})(\diracComb(\mathrm{id},\widehat{\mathscr{S}}^c(d_\PP)(\mathrm{id})))~.
\]

Since $d_\PP$ is a derivation with respect to the operad structure of $\PP^{u}$, it makes the resulting morphism $\mathscr{C}(d_\PP)(\mathrm{id})$ a coderivation of $\mathscr{C}(\PP^{u})(X)~.$
\end{proof}

\begin{Proposition}[Coderivations on the cofree construction]
Let $X$ be a graded module and let $\mathscr{C}(\PP^{u})(X)$ be the cofree graded $\PP^{u}$-coalgebra generated by $X$. There is a natural bijection between maps 

\[
\varphi: \mathscr{C}(\PP^{u})(X) \longrightarrow X
\]
of degree $p$ and coderivations 

\[
d_\varphi: \mathscr{C}(\PP^{u})(X) \longrightarrow \mathscr{C}(\PP^{u})(X)
\]
of degree $p$. This bijection sends $\varphi$ to the degree $p$ coderivation given by 

\[
d_\varphi \coloneqq \mathscr{C}(\mathrm{id})(\diracComb(\pi_X,\varphi)) \cdot \omega(X)~,
\]

where $\omega(X)$ is the comonad structure map of $\mathscr{C}(\PP^{u})$.
\end{Proposition}

\begin{proof}
See \cite[Proposition 6.23]{grignoulejay18}.
\end{proof}

\begin{Definition}[Complete Bar construction relative to $\alpha$]\label{def: complete Bar}
Let $(B,\gamma_B,d_B)$ be a dg $\C$-algebra. The \textit{complete Bar construction relative to} $\alpha$ of $B$ is given by
\[
\widehat{\mathrm{B}}_{\alpha} B \coloneqq (\mathscr{C}(\PP^{u})(B), d_{\mathrm{bar}} \coloneqq d_1 + d_2)~,
\]
where $\mathscr{C}(\PP)(B)$ denotes the cofree graded $\PP^{u}$-coalgebra generated by $B$. The differential $d_{\mathrm{bar}}$ is given by the sum of two terms $d_1$ and $d_2$. The term $d_1$ is given by
\[
d_1 = \mathscr{C}(d_\PP)(\mathrm{id}) +  \mathscr{C}(\mathrm{id})(\diracComb(\mathrm{id}_B,d_B))~.
\]
The term $d_2$ is given by the unique coderivation extending  
\[
\begin{tikzcd}[column sep=4pc,row sep=1pc]
\mathscr{C}(\PP^{u})(B) \arrow[r,"p_1 (B)",rightarrowtail]
&\widehat{\mathscr{S}}^c(\PP^{u})(B)\arrow[r,"\widehat{\mathscr{S}}^c(\alpha)(\mathrm{id}_B)"]
&\widehat{\mathscr{S}}^c(\C^{u})(B)  \arrow[r,"\gamma_B "]
&B~.
\end{tikzcd}
\]
\end{Definition}

\begin{Proposition}
For any dg $\C$-algebra $B$, the complete Bar construction $\widehat{\mathrm{B}}_{\alpha}B$ forms a dg $\PP$-coalgebra, and it defines a functor
\[
\widehat{\mathrm{B}}_{\alpha}: \mathsf{dg}~\C^{u}\text{-}\mathsf{alg} \longrightarrow \mathsf{dg}~\PP^{u}\text{-}\mathsf{coalg}~.
\]
\end{Proposition}

\begin{Definition}[Complete Cobar construction relation to $\alpha$]\label{def: complete Cobar}
Let $(D,\delta_D,d_D)$ be a dg $\PP$-coalgebra. The \textit{complete Cobar construction relative to} $\alpha$ of $D$ is given by
\[
\widehat{\Omega}_\alpha D \coloneqq (\widehat{\mathscr{S}}^c(\C^{u})(D), d_{\mathrm{cobar}} \coloneqq d_1 - d_2)~,
\]
where $\widehat{\mathscr{S}}^c(\C^{u})(D)$ denotes the free complete graded $\C$-coalgebra generated by $D$. The differential $d_{\mathrm{cobar}}$ is given by the difference of two terms $d_1$ and $d_2$. The term $d_1$ is given by
\[
d_1 = -\widehat{\mathscr{S}}^c(d_\C)(\mathrm{id}) + \widehat{\mathscr{S}}^c(\mathrm{id})(\diracComb(\mathrm{id}_D,d_D))~.
\]
The term $d_2$ is given by the unique derivation extending  
\[
\begin{tikzcd}[column sep=4pc,row sep=1pc]
D \arrow[r," \Delta_D"]
&\widehat{\mathscr{S}}^c(\PP^{u})(D) \arrow[r,"\widehat{\mathscr{S}}^c(\alpha)(\mathrm{id}_D )"]
&\widehat{\mathscr{S}}^c(\C^{u})(D)~.
\end{tikzcd} 
\]
\end{Definition}

\begin{Proposition}
For any dg $\PP$-coalgebra $D$, the complete Cobar construction $\widehat{\Omega}_{\alpha}D$ forms a complete dg $\C$-algebra, and it defines a functor
\[
\widehat{\Omega}_\alpha: \mathsf{dg}~\PP^{u}\text{-}\mathsf{coalg} \longrightarrow \mathsf{dg}~\C^{u}\text{-}\mathsf{alg}^{\mathsf{comp}}~.
\]
\end{Proposition}

\begin{Definition}[Twisting morphism relative to $\alpha$]
Let $(D,\delta_D,d_D)$ be a dg $\PP$-coalgebra and let $(B,\gamma_B,d_B)$ be a dg $\C$-algebra. A graded morphism 
\[
\nu: D \longrightarrow B
\]
is said to be a twisting morphism relative to $\alpha$ if it satisfies the following equation 

\[
\gamma_B \cdot \widehat{\mathscr{S}}^c(\alpha)(\nu) \cdot \Delta_D + \partial(\nu) = 0~.
\]

The set of twisting morphisms relative to $\alpha$ are denoted by $\mathrm{Tw}^{\alpha}(D,B)$.
\end{Definition}

\begin{Proposition}
There are bijections

\[
\mathrm{Hom}_{\mathsf{dg}~\C\text{-}\mathsf{alg}^{\mathsf{comp}}}\left(\widehat{\Omega}_{\alpha}D, A\right) \cong  \mathrm{Tw}^{\alpha}(D,B)  \cong \mathrm{Hom}_{\mathsf{dg}~\PP\text{-}\mathsf{coalg}} \left(D, \widehat{\mathrm{B}}_{\alpha}B \right)~,
\]

which are natural in $D$ and $B$. 
\end{Proposition}

\begin{proof}
The proof is very similar to the proof of Proposition \ref{prop: Bar-Cobar adjunction relative to alpha}.
\end{proof}

Therefore we get an adjunction 
\[
\begin{tikzcd}[column sep=5pc,row sep=3pc]
            \mathsf{dg}~\PP\text{-}\mathsf{coalg} \arrow[r, shift left=1.1ex, "\widehat{\Omega}_{\alpha}"{name=F}] & \mathsf{dg}~\C\text{-}\mathsf{alg}^{\mathsf{comp}} \arrow[l, shift left=.75ex, "\widehat{\mathrm{B}}_{\alpha}"{name=U}]
            \arrow[phantom, from=F, to=U, , "\dashv" rotate=-90]
\end{tikzcd}
\]
between the category of dg $\PP^{u}$-coalgebras and the category of complete dg $\C$-algebras. This adjunction is called the \textit{complete Bar-Cobar adjunction} relative to $\alpha$. In some case, one can promote this adjunction into a Quillen adjunction. But first, one needs a model structure on the category of dg $\PP$-coalgebras. 

\begin{theorem}[{\cite[Section 9]{grignoulejay18}}]\label{thm: model structure on P-cog}
Suppose that $\mathcal{P}$ is a cofibrant dg operad. Then there is a model structure on the category of dg $\mathcal{P}$-coalgebras left-transferred along the cofree-forgetful adjunction
\[
\begin{tikzcd}[column sep=7pc,row sep=3pc]
\mathsf{dg}\textsf{-}\mathsf{mod} \arrow[r, shift left=1.1ex, "\mathcal{L}(\mathcal{P}^{u})(-)"{name=F}]      
&\mathsf{dg}~\PP^{u}\textsf{-}\mathsf{coalg}~, \arrow[l, shift left=.75ex, "U"{name=U}]
\arrow[phantom, from=F, to=U, , "\dashv" rotate=90]
\end{tikzcd}
\]
where 

\begin{enumerate}

\medskip 

\item the class of weak equivalences is given by quasi-isomorphisms,

\medskip

\item the class of cofibrations is given by degree-wise monomorphisms,

\medskip 

\item the class of fibrations is given by right lifting property with respect to acyclic cofibrations.

\medskip
\end{enumerate}
\end{theorem}

\begin{Remark}
The assumption that $\PP$ is a cofibrant dg operad is necessary. Indeed, one can show that for $\PP = u\mathcal{C}om$, the category of dg $ u\mathcal{C}om$-coalgebras does not admit a model structure where weak-equivalences are given by quasi-isomorphisms and where cofibrations are given by degree-wise monomorphisms. Therefore, while all dg operads are \textit{admissible} in characteristic zero, not all dg operads are \textit{coadmissible}.

\medskip

Nevertheless, it might be the case $\PP$ is not cofibrant and this model category structure still exists. For example, this model structure does exist for $u\mathcal{A}ss$, which encodes non-necessarily conilpotent dg coassociative coalgebras. See the model structure constructed in \cite{GetzlerGoerss99}.
\end{Remark}

Model structures on coalgebras over dg operads behave well with respect to quasi-isomorphisms at the operadic level.

\begin{theorem}[{\cite[Section 9]{grignoulejay18}}]
Let $f: \mathcal{P} \qi \mathcal{Q}$ be a quasi-isomorphism of cofibrant dg operads. The induced adjunction 

\[
\begin{tikzcd}[column sep=7pc,row sep=3pc]
\mathsf{dg}~\mathcal{P}\textsf{-}\mathsf{coalg} \arrow[r, shift left=1.1ex, "\mathrm{Coind}_{f}"{name=F}]      
&\mathsf{dg}~\mathcal{Q}\textsf{-}\mathsf{coalg}~, \arrow[l, shift left=.75ex, "\mathrm{Res}_{f}"{name=U}]
\arrow[phantom, from=F, to=U, , "\dashv" rotate=90]
\end{tikzcd}
\]

is a Quillen equivalence.
\end{theorem}

When this model structure exists, it can be transferred along the complete Bar-Cobar adjunction relative to $\alpha$. 

\begin{theorem}[{\cite[Section 10]{grignoulejay18}}]
Suppose that $\mathcal{P}$ is a cofibrant dg operad. There is a model structure on the category of complete dg $\mathcal{C}$-algebras right-transferred along the complete Bar-Cobar adjunction
\[
\begin{tikzcd}[column sep=7pc,row sep=3pc]
            \mathsf{dg}~\mathcal{P}^{u}\textsf{-}\mathsf{coalg} \arrow[r, shift left=1.1ex, "\widehat{\Omega}_{\alpha}"{name=F}] &\mathsf{dg}~\mathcal{C}^{u}\textsf{-}\mathsf{alg}^{\mathsf{comp}}~, \arrow[l, shift left=.75ex, "\widehat{\mathrm{B}}_{\alpha}"{name=U}]
            \arrow[phantom, from=F, to=U, , "\dashv" rotate=-90]
\end{tikzcd}
\]
where 
\begin{enumerate}

\medskip

\item the class of weak equivalences is given by morphisms $f$ such that $\widehat{\mathrm{B}}_\alpha(f)$ is a quasi-isomorphism,

\medskip

\item the class of fibrations is given by degree-wise epimorphisms,

\medskip

\item  and the class of cofibrations is given by left lifting property with respect to acyclic fibrations.
\end{enumerate}
\end{theorem}

In the case where $\PP$ is the Cobar construction of $\mathcal{C}$, this adjunction can be promoted to a Quillen equivalence.

\begin{theorem}[{\cite[Section 11]{grignoulejay18}}] \label{thm: complete Bar-Cobar relative to iota Quillen equivalence}
The complete Bar-Cobar adjunction relative to $\iota: \C \longrightarrow \Omega \C$
\[
\begin{tikzcd}[column sep=7pc,row sep=3pc]
            \mathsf{dg}~\Omega\mathcal{C}\textsf{-}\mathsf{coalg} \arrow[r, shift left=1.1ex, "\widehat{\Omega}_{\iota}"{name=F}] 
            &\mathsf{dg}~\mathcal{C}\textsf{-}\mathsf{alg}^{\mathsf{comp}}~. \arrow[l, shift left=.75ex, "\widehat{\text{B}}_{\iota}"{name=U}]
            \arrow[phantom, from=F, to=U, , "\dashv" rotate=-90]
\end{tikzcd}
\]
is a Quillen equivalence. 
\end{theorem}

\begin{Remark}
The model structure on the category of complete dg $\mathcal{C}$-algebras transferred using the complete Bar-Cobar adjunction relative to the twisting morphism $\iota: \mathcal{C} \longrightarrow \Omega \mathcal{C}$ is called the \textit{canonical model structure} in \textit{loc.cit}.
\end{Remark}

This adjunction allows us to define $\infty_\alpha$-morphisms of dg $\PP^{u}$-coalgebras. The theory is developed in \cite[Section 12]{grignoulejay18}.

\begin{Definition}[$\infty_\alpha$-morphism of coalgebras]
Let $D$ and $E$ be two dg $\PP^{u}$-coalgebras. An $\infty_\alpha$\textit{-morphism} $f: D \rightsquigarrow E$ is the data of a morphism 
\[
f: \widehat{\Omega}_{\alpha}D \longrightarrow \widehat{\Omega}_{\alpha} E
\]
of dg $\C^{u}$-algebras. 
\end{Definition}

\begin{Remark}
This data is equivalent to a family of maps $f_n: D \longrightarrow \mathrm{Hom}_{\mathbb{S}_n}(\C^{u}(n),E^{\otimes n})$ which satisfy certain relations. Since $\C^{u}(1) = \C(1) \oplus \I$ the morphism $f_1$ has a linear part $f_\I: \cong D \longrightarrow \mathrm{Hom}(\I, B) \cong B$ which is a morphism of dg modules. 
\end{Remark}

Notice that any morphism of dg $\PP^{u}$-algebras $f: D \longrightarrow E$ is an $\infty_\alpha$-morphism by setting $f_\I = f$ and zero to all other components.

\begin{Definition}[$\infty_\alpha$-quasi-isomorphisms]
An $\infty_\alpha$-\textit{quasi-isomorphism} is an $\infty_\alpha$-morphism $f: D \rightsquigarrow E$ where $f_\I$ is a quasi-isomorphism of dg modules.
\end{Definition}

We now restrict to the case of the twisting morphism $\iota$, which endows dg $\C^{u}$-algebras with the canonical model structure.

\begin{Proposition}\label{prop: infinty mor of coalgebras equivalences and invertible}
Let $f: D \rightsquigarrow E$ be an $\infty_\iota$-quasi-isomorphism between two dg $\Omega \C^{u}$-algebras. 

\medskip

\begin{enumerate}

\item The morphism $f: \widehat{\Omega}_{\iota}D \longrightarrow \widehat{\Omega}_{\iota} E$ is a weak-equivalence of dg $\C^{u}$-algebras.

\medskip

\item There exists an $\infty_\iota$-quasi-isomorphism $g: D \rightsquigarrow E$ such that $g_\I$ is the inverse of $f_\I$ on the level of the homology groups. 
\end{enumerate}
\end{Proposition}

\begin{proof}
The first point is shown in \cite[Theorem 12.10]{grignoulejay18}. The second point can be deduced from \cite[Theorem 4.18]{hoffbeck2019properadic}, where they show that $\infty_\iota$-morphisms of dg properads admit an inverse. Since coalgebras over an dg operad can be seen as gebras over the properad obtained by reversing the dg operad, this allows us to conclude. 
\end{proof}

\begin{Corollary}
Let $D$ and $E$ be two dg $\Omega \C^{u}$-coalgebras. There exists a zig-zag of quasi-isomorphisms 
\[
D \lqi \cdot \qi \cdot \lqi \cdot \qi E
\]
of dg $\Omega \C^{u}$-algebras if and only if there exists an $\infty_\iota$-quasi-isomorphism $D \rightsquigarrow E$.
\end{Corollary}

\begin{proof}
Follows from Theorem \ref{thm: complete Bar-Cobar relative to iota Quillen equivalence} combined with Proposition \ref{prop: infinty mor of coalgebras equivalences and invertible}.
\end{proof}

Furthermore, these adjunctions are "functorial" in the following sense.

\begin{Proposition}[{\cite[Lemma 9.8]{grignoulejay18}}]\label{prop: compatibility of induced adjunctions}
Let $\mathcal{P}$ and $\mathcal{Q}$ be two cofibrant dg operads, and let $\mathcal{C}$ and $\mathcal{D}$ be conilpotent curved cooperads. 

\medskip

\begin{enumerate}
\item Let $\alpha: \C \longrightarrow \PP$ and $\beta: \mathcal{D} \longrightarrow \mathcal{Q}$ be two curved twisting morphisms. 

\medskip

\item Let $f: \PP \longrightarrow \mathcal{Q}$ be a morphism of dg operads and let $g: \mathcal{C} \longrightarrow \mathcal{D}$ be a morphism of conilpotent curved cooperads. 
\end{enumerate}

\medskip

Such that the following diagram commutes 

\[
\begin{tikzcd}[column sep=4pc,row sep=4pc]
\C \arrow[r,"\alpha"] \arrow[d,"g",swap] 
&\mathcal{P} \arrow[d,"f"]\\
\mathcal{D} \arrow[r,"\beta"]
&\mathcal{Q}~.
\end{tikzcd}
\]

Then the following square 

\[
\begin{tikzcd}[column sep=5pc,row sep=5pc]
\mathsf{dg}~\mathcal{P}\text{-}\mathsf{coalg} \arrow[r,"\widehat{\Omega}_\alpha"{name=B},shift left=1.1ex] \arrow[d,"\mathrm{Coind}_f "{name=SD},shift left=1.1ex ]
&\mathsf{curv}~\mathcal{C}\text{-}\mathsf{alg}^{\mathsf{comp}} \arrow[d,"\mathrm{Res}_g"{name=LDC},shift left=1.1ex ] \arrow[l,"\widehat{\mathrm{B}}_\alpha"{name=C},,shift left=1.1ex]  \\
\mathsf{dg}~\mathcal{Q}\text{-}\mathsf{coalg} \arrow[r,"\widehat{\Omega}_\beta "{name=CC},shift left=1.1ex]  \arrow[u,"\mathrm{Res}_f"{name=LD},shift left=1.1ex ]
&\mathsf{curv}~\mathcal{D}\text{-}\mathsf{alg}^{\mathsf{comp}} \arrow[l,"\widehat{\mathrm{B}}_\beta"{name=CB},shift left=1.1ex] \arrow[u,"\mathrm{Ind}_g"{name=TD},shift left=1.1ex] \arrow[phantom, from=SD, to=LD, , "\dashv" rotate=0] \arrow[phantom, from=C, to=B, , "\dashv" rotate=-90]\arrow[phantom, from=TD, to=LDC, , "\dashv" rotate=0] \arrow[phantom, from=CC, to=CB, , "\dashv" rotate=-90]
\end{tikzcd}
\] 

of Quillen adjunctions commutes. 
\end{Proposition}

\begin{Remark}
The cofibrancy assumption on $\mathcal{P}$ and $\mathcal{Q}$ is not needed for the commutativity of this square, only for the existence of model structures.
\end{Remark}

\begin{Remark}
Please note that in \cite{grignoulejay18}, the authors develop the general case where units/curvature are present. For pedagogical reasons, we have chosen to restrict to the more classical case their results. In Chapter 2, we will develop the general theory of curved operadic calculus. Their results will be later explained in full generality.
\end{Remark}

\section{Duality squares}\label{Section: Magic squares}
Let us fix a dg operad $\mathcal{P}$ and a conilpotent dg cooperad $\mathcal{C}$, together with a twisting morphism $\alpha: \mathcal{C} \longrightarrow \mathcal{P}$. The goal of this section is to construct two duality adjunctions that interrelate the classical Bar-Cobar constructions relative to $\alpha$ with the complete Bar-Cobar constructions relative to $\alpha$. These results admit a further generalization  given in Section \ref{Section: Curved duality square}.

\subsection{Sweedler functor}
The linear dual of a coalgebra over a given operad is naturally an algebra over the same operad. 

\begin{lemma}
The linear dual defines a functor
\[
\begin{tikzcd}[column sep=4pc,row sep=0pc]
\left(\mathsf{dg}~\mathcal{P}\text{-}\mathsf{coalg}\right)^{\mathsf{op}} \arrow[r,"(-)^*"] 
&\mathsf{dg}~\mathcal{P}\text{-}\mathsf{alg}~.
\end{tikzcd}
\]
\end{lemma}

\begin{proof}
Let $(C,\Delta_C,d_C)$ be a dg $\mathcal{P}$-coalgebra. Given an element $\mu$ in $\mathcal{P}(n)$, it induces a decomposition map 
\[
\Delta_{\mu}: C \longrightarrow C^{\otimes n}~.
\]
It induces a map
\[
\begin{tikzcd}[column sep=4pc,row sep=0pc]
\gamma_\mu^*: (C^*)^{\otimes n} \arrow[r,"i",rightarrowtail] 
&(C^{\otimes n})^* \arrow[r,"(\Delta_{\mu})^* "] 
&C^*~. 
\end{tikzcd}
\]
One can check that 
\[
\left\{
\begin{tikzcd}[column sep=4pc,row sep=0pc]
\mathcal{P} \arrow[r, "\Gamma_{C^* } "] 
& \mathrm{End}_{C^*} \\
\mu \arrow[r,mapsto]
&\gamma_\mu^*
\end{tikzcd}
\right.
\]
is indeed a morphism of dg operads. Therefore the linear dual of a dg $\mathcal{P}$-coalgebra has a canonical dg $\mathcal{P}$-algebra structure, and one checks that this association is indeed functorial with respect to morphisms of dg $\mathcal{P}$-coalgebras since the morphism $i: (C^*)^{\otimes n} \longrightarrow (C^{\otimes n})^*$ is a natural transformation.
\end{proof}

\begin{Proposition}\label{prop: adjoint à droite}
The linear duality functor 
\[
\begin{tikzcd}[column sep=4pc,row sep=0pc]
\left(\mathsf{dg}~\mathcal{P}\text{-}\mathsf{coalg}\right)^{\mathsf{op}} \arrow[r,"(-)^*"] 
&\mathsf{dg}~\mathcal{P}\text{-}\mathsf{alg}~.
\end{tikzcd}
\]
admits a left adjoint.
\end{Proposition}

\begin{proof}
Consider the following square of functors
\[
\begin{tikzcd}[column sep=4pc,row sep=4pc]
\left(\mathsf{dg}~\mathcal{P}\text{-}\mathsf{coalg}\right)^{\mathsf{op}} \arrow[r,"(-)^*"]  \arrow[d,"\mathrm{U}^\mathsf{op}"{name=SD},shift left=1.1ex ]
&\mathsf{dg}~\mathcal{P}\text{-}\mathsf{alg} \arrow[d,"\mathrm{U}"{name=LDC},shift left=1.1ex ] \\
\mathsf{dg}~\mathsf{mod}^{\mathsf{op}} \arrow[r,"(-)^*"{name=CC},,shift left=1.1ex] \arrow[u,"\left(\mathscr{C}(\mathcal{P})(-)\right)^\mathsf{op}"{name=LD},shift left=1.1ex ] \arrow[phantom, from=SD, to=LD, , "\dashv" rotate=0]
&\mathsf{dg}~\mathsf{mod}~,  \arrow[l,"(-)^*"{name=CB},shift left=1.1ex] \arrow[u,"\mathscr{S}(\mathcal{P})(-)"{name=TD},shift left=1.1ex] \arrow[phantom, from=TD, to=LDC, , "\dashv" rotate=0] \arrow[phantom, from=CC, to=CB, , "\dashv" rotate=90]
\end{tikzcd}
\] 
where $\mathscr{S}(\mathcal{P})(-)$ is the free dg $\mathcal{P}$-algebra functor and where $\mathscr{C}(\mathcal{P})(-)$ is the cofree dg $\mathcal{P}$-coalgebra functor given by Theorem \ref{thm: existence of the cofree P cog}. The left hand side adjunction is monadic, since it is the opposite of a comonadic adjunction. All categories involved are complete and cocomplete. We also have that $(-)^* \cdot \mathrm{U}^\mathsf{op} \cong \mathrm{U} \cdot (-)^*~.$ Thus we can apply the Adjoint Lifting Theorem \cite[Theorem 2]{AdjointLifting}, which concludes the proof. 
\end{proof}

\begin{Definition}[Sweedler dual]\label{def: Sweedler dual functor}
The \textit{Sweedler duality functor}
\[
\begin{tikzcd}[column sep=4pc,row sep=0pc]
\mathsf{dg}~\mathcal{P}\text{-}\mathsf{alg} \arrow[r,"(-)^\circ"] 
&\left(\mathsf{dg}~\mathcal{P}\text{-}\mathsf{coalg}\right)^{\mathsf{op}}
\end{tikzcd}
\]
is defined as the functor left adjoint of the linear dual functor. 
\end{Definition}

\begin{Remark}
The proof of the Adjoint Lifting Theorem \cite[Theorem 2]{AdjointLifting} gives an explicit construction of this left adjoint. Let $(A,\gamma_A,d_A)$ be a dg $\mathcal{P}$-algebra. The Sweelder dual dg $\mathcal{P}$-coalgebra $(A)^\circ$ is given by the following equalizer: 
\[
\begin{tikzcd}[column sep=4pc,row sep=4pc]
\mathrm{Eq}\Bigg(\mathscr{C}(\mathcal{P})(A^*) \arrow[r,"(\gamma_A)^*",shift right=1.1ex,swap]  \arrow[r,"\varrho"{name=SD},shift left=1.1ex ]
&\mathscr{C}(\mathcal{P})\left((\mathscr{S}(\mathcal{P})(A))^*\right) \Bigg)~,
\end{tikzcd}
\]
where $\varrho$ is an arrow constructed using the comonadic structure of $\mathscr{C}(\mathcal{P})$ and the canonical inclusion of a dg module into its double linear dual.
\end{Remark}

\begin{Proposition}\label{prop: natural mono for Sweedler dual}
There is a natural monomorphism
\[
\epsilon: \mathrm{U}^{\mathsf{op}} \cdot (-)^\circ \rightarrowtail (-)^* \cdot \mathrm{U}~,
\]
which implies that the Sweedler dual is a sub-dg module of the linear dual functor. 
\end{Proposition}

\begin{proof}
Let $V$ be a dg module, there is a monomorphism
\[
\begin{tikzcd}[column sep=4pc,row sep=4pc]
\epsilon_V: \left(\mathscr{S}(\mathcal{P})(V) \right)^\circ \cong \mathscr{C}(\mathcal{P})(V^*) \arrow[r,rightarrowtail,"(p_1)_V"]
&\widehat{\mathscr{S}}^c(\mathcal{P})(V^*) \arrow[r,rightarrowtail]
&\left(\mathscr{S}(\mathcal{P})(V) \right)^*~,
\end{tikzcd}
\]
where the monomorphism $p_1$ is given by Theorem \ref{thm: existence of the cofree P cog}, hence the proposition is true on free dg $\mathcal{P}$-algebras. Any dg $\mathcal{P}$-algebra can be written as a split coequalizer of free dg $\mathcal{P}$-algebras. Both $\mathrm{U}^{\mathsf{op}} \cdot (-)^\circ$ and $(-)^* \cdot \mathrm{U}$ send these split coequalizers to split equalizers of dg modules, therefore the monomorphism $\epsilon$ extends to all dg $\mathcal{P}$-algebras.
\end{proof}

\begin{Remark}[Beck--Chevalley condition]
Let $A$ be a dg $\mathcal{P}$-algebra which is degree-wise finite dimensional. Then the map
\[
\epsilon_A: \mathrm{U}^{\mathsf{op}} \cdot A^\circ \cong A^* \cdot \mathrm{U}~,
\]
is an isomorphism of dg modules. Furthermore, there is a canonical structure of dg $\mathcal{P}$-coalgebra on the linear dual $A^*$. Therefore the adjunction $(-)^\circ \dashv (-)^*$ restricts to an anti-equivalence of categories between the category of degree-wise finite dimensional dg $\mathcal{P}$-algebra and the category of degree-wise finite dimensional dg $\mathcal{P}$-coalgebras.
\end{Remark}

\begin{Example}
Consider the operad $\mathcal{A}ss$ which encodes dg associative algebras as its algebras and dg coassociative coalgebras as its coalgebras. Then the adjunction 

\[
\begin{tikzcd}[column sep=7pc,row sep=3pc]
\mathsf{dg}~\mathcal{A}ss\text{-}\mathsf{alg}  \arrow[r, shift left=1.1ex, "(-)^\circ"{name=F}] 
& \left(\mathsf{dg}~\mathcal{A}ss\text{-}\mathsf{coalg}\right)^{\mathsf{op}} ~. \arrow[l, shift left=.75ex, "(-)^*"{name=U}] \arrow[phantom, from=F, to=U, , "\dashv" rotate=-90]
\end{tikzcd}
\]

coincides with the original Sweedler adjunction constructed in \cite{Sweedler69}.
\end{Example}

\subsection{Topological dual functor}
Let's turn to the other side of the Koszul duality, where $\mathcal{C}$ is a conilpotent dg cooperad.

\begin{lemma}
The linear duality defines a functor
\[
\begin{tikzcd}[column sep=4pc,row sep=0pc]
\left(\mathsf{dg}~\mathcal{C}\text{-}\mathsf{coalg}\right)^{\mathsf{op}} \arrow[r,"(-)^*"] 
&\mathsf{dg}~\mathcal{C}\text{-}\mathsf{alg}^{\mathsf{comp}}
\end{tikzcd}
\]
from the category of dg $\mathcal{C}$-coalgebras to the category of complete dg $\mathcal{C}$-algebras.
\end{lemma}

\begin{proof}
Let $(C,\Delta_C, d_C)$ be a dg $\mathcal{C}$-coalgebra, where 
\[
\Delta_C: C \longrightarrow \bigoplus_{n \geq 0} \mathcal{C}(n) \otimes_{\mathbb{S}_n} C^{\otimes n} 
\]
is the structural morphism. By applying the linear duality, we get a map
\[
\begin{tikzcd}[column sep=3.5pc,row sep=0pc]
\gamma_{C^*}: \displaystyle \prod_{n \geq 0} \mathrm{Hom}_{\mathbb{S}_n}\left(\mathcal{C}(n), (C^*)^{\otimes n} \right) \arrow[r,rightarrowtail]
&\displaystyle \prod_{n \geq 0} \mathrm{Hom}_{\mathbb{S}_n}\left(\mathcal{C}(n), (C^{\otimes n})^* \right) \arrow[r,"(\Delta_C)^* "]
&C^*~.
\end{tikzcd}
\]
One can check that it defines a dg $\mathcal{C}$-algebra structure on $C^*$. Furthermore, let 
\[
\mathrm{F}_\omega C \coloneqq \mathrm{Ker}\left(\Delta_C^\omega: C \longrightarrow \bigoplus_{n \geq 0} \mathcal{C}/\mathscr{R}_\omega\mathcal{C}(n) \otimes_{\mathbb{S}_n} C^{\otimes n} \right)
\] 
be the canonical coradical filtration on $C$ induced by the coradical filtration on $\mathcal{C}$. Since $\mathcal{C}$ is conilpotent, the coradical filtration of any dg $\mathcal{C}$-coalgebra is exhaustive, therefore 
\[
C \cong \colim_{\omega}\mathrm{F}_\omega C~,
\]
which in turn implies that
\[
C^* \cong \lim_{\omega} ~ (\mathrm{F}_\omega C)^*~.
\]
One can check that $(\mathrm{F}_\omega C)^* \cong C^*/\mathrm{W}_\omega C^*$, therefore the image of the linear duality functor $(-)^*$ lies in the sub-category of complete dg $\mathcal{C}$-coalgebras.
\end{proof}

\begin{Proposition}\label{prop: adjoint à gauche dual topo}
The linear duality functor
\[
\begin{tikzcd}[column sep=4pc,row sep=0pc]
\left(\mathsf{dg}~\mathcal{C}\text{-}\mathsf{coalg}\right)^{\mathsf{op}} \arrow[r,"(-)^*"] 
&\mathsf{dg}~\mathcal{C}\text{-}\mathsf{alg}^{\mathsf{comp}}
\end{tikzcd}
\]
admits a left adjoint. 
\end{Proposition}

\begin{proof}
We consider the following square of functors
\[
\begin{tikzcd}[column sep=4pc,row sep=4pc]
\left(\mathsf{dg}~\mathcal{C}\text{-}\mathsf{coalg}\right)^{\mathsf{op}} \arrow[r,"(-)^*"]  \arrow[d,"\mathrm{U}^\mathsf{op}"{name=SD},shift left=1.1ex ]
&\mathsf{dg}~\mathcal{C}\text{-}\mathsf{alg}^{\mathsf{comp}} \arrow[d,"\mathrm{U}"{name=LDC},shift left=1.1ex ] \\
\mathsf{dg}~\mathsf{mod}^{\mathsf{op}} \arrow[r,"(-)^*"{name=CC},,shift left=1.1ex] \arrow[u,"\left(\mathscr{S}(\mathcal{C})(-)\right)^\mathsf{op}"{name=LD},shift left=1.1ex ] \arrow[phantom, from=SD, to=LD, , "\dashv" rotate=0]
&\mathsf{dg}~\mathsf{mod}~,  \arrow[l,"(-)^*"{name=CB},shift left=1.1ex] \arrow[u,"\widehat{\mathscr{S}}^c(\mathcal{C})(-)"{name=TD},shift left=1.1ex] \arrow[phantom, from=TD, to=LDC, , "\dashv" rotate=0] \arrow[phantom, from=CC, to=CB, , "\dashv" rotate=90]
\end{tikzcd}
\] 
where $\mathscr{S}(\mathcal{C})(-)$ is the cofree dg $\mathcal{C}$-algebra functor and where $\widehat{\mathscr{S}}^c(\mathcal{C})(-)$ is the free dg $\mathcal{C}$-algebra functor, which is always complete. Again, vertical adjunctions are monadic and it is clear that $(-)^* \cdot \mathrm{U}^\mathsf{op} \cong \mathrm{U} \cdot (-)^*~.$ Thus we can apply the Adjoint Lifting Theorem \cite[Theorem 2]{AdjointLifting}, which concludes the proof. 
\end{proof}

\begin{Definition}[Topological dual functor]\label{def: topological dual functor}
The \textit{topological dual functor} 
\[
\begin{tikzcd}[column sep=4pc,row sep=0pc]
\mathsf{dg}~\mathcal{C}\text{-}\mathsf{alg}^{\mathsf{comp}} \arrow[r,"(-)^\vee"]
&\left(\mathsf{dg}~\mathcal{C}\text{-}\mathsf{coalg}\right)^{\mathsf{op}}
\end{tikzcd}
\]
is defined as the functor left adjoint to the linear dual functor.
\end{Definition}

\begin{Remark}\label{Rmk: formule pour le dual topologique}
Given a complete dg $\mathcal{C}$-algebra $(B, \gamma_B, d_B)$, its topological dual $B^\vee$ is given by the following equalizer:
\[
\begin{tikzcd}[column sep=4pc,row sep=4pc]
\mathrm{Eq}\Bigg(\mathscr{S}(\mathcal{C})(B^*) \arrow[r,"(\gamma_B)^*",shift right=1.1ex,swap]  \arrow[r,"\varrho"{name=SD},shift left=1.1ex ]
&\mathscr{S}(\mathcal{C})\left((\widehat{\mathscr{S}}^c(\mathcal{C})(B))^*\right) \Bigg)~,
\end{tikzcd}
\]
where $\varrho$ is an arrow constructed using the comonadic structure of $\mathscr{S}(\mathcal{C})(-)$ and the canonical inclusion of a dg module into its double linear dual.
\end{Remark}

\begin{Proposition}\label{prop: natural mono for topo dual}
There is a natural monomorphism
\[
\epsilon: \mathrm{U}^{\mathsf{op}} \cdot (-)^\vee \rightarrowtail (-)^* \cdot \mathrm{U}~,
\]
which implies that the topological dual is a sub-dg module of the linear dual functor. 
\end{Proposition}

\begin{proof}
Let $V$ be a dg module, there is a monomorphism
\[
\begin{tikzcd}[column sep=4pc,row sep=0pc]
\epsilon_V: \left(\widehat{\mathscr{S}}^c(\mathcal{C})(V)\right)^\vee \cong \mathscr{S}(\mathcal{C})(V^*) \arrow[r,rightarrowtail]
&\left(\widehat{\mathscr{S}}^c(\mathcal{C})(V)\right)^*~,
\end{tikzcd}
\]
which is given by the following composite
\[
\begin{tikzcd}[column sep=0pc,row sep=2pc]
\mathscr{S}(\mathcal{C})(V^*) \arrow[d,rightarrowtail,"\eta_{\mathscr{S}(\mathcal{C})(V^*)} "] \\
\left(\mathrm{Hom}_{\mathsf{dg}~\mathsf{mod}}\left(\bigoplus_{n \geq 0} \mathcal{C}(n) \otimes_{\mathbb{S}_n} (V^*)^{\otimes n},\kk\right)\right)^* \arrow[d,"\cong"] \\
\left(\prod_{n \geq 0} \mathrm{Hom}_{\mathbb{S}_n}\left(\mathcal{C}(n),((V^*)^{\otimes n})^*\right)\right)^* \arrow[d,twoheadrightarrow] \\
\left(\prod_{n \geq 0} \mathrm{Hom}_{\mathbb{S}_n}\left(\mathcal{C}(n),(V^{**})^{\otimes n}\right)\right)^* \arrow[d,twoheadrightarrow] \\
\left(\prod_{n \geq 0} \mathrm{Hom}_{\mathbb{S}_n}\left(\mathcal{C}(n),V^{\otimes n}\right)\right)^* \arrow[d,"\cong"] \\
\left(\widehat{\mathscr{S}}^c(\mathcal{C})(V)\right)^*~,
\end{tikzcd}
\]
where $\eta$ is the inclusion of a dg module into its double dual. Let us show that $\epsilon_V$ is a monomorphism. Using the explicit formula that appears in the proof of \cite[Theorem 2]{AdjointLifting}, we compute that the morphism of dg modules $\rho$ that appears in the equalizer that defines $(-)^\vee$ is given by the following composite 

\[
\begin{tikzcd}[column sep=4pc,row sep=2pc]
\rho: \mathscr{S}(\mathcal{C})(V^*) \arrow[r,"\mathscr{S}(\Delta_\mathcal{C})"]
&\mathscr{S}(\mathcal{C}) \left(\mathscr{S}(\mathcal{C})(V^*)\right) \arrow[r,"\mathscr{S}(\mathrm{id})(\epsilon_V)"]
&\mathscr{S}(\mathcal{C}) \left(\widehat{\mathscr{S}}^c(\mathcal{C})(V)\right)^*~,
\end{tikzcd}
\]

see Remark \ref{Rmk: formule pour le dual topologique}. The aforementioned equalizer is $\mathrm{U}$-split, where the splitting is given by the morphism $\mathscr{S}(\iota_V^*)$, where $\iota_V$ is the unit $V \longrightarrow \widehat{\mathscr{S}}^c(\mathcal{C})(V)$. Therefore $\rho$ is a split monomorphism of dg modules. 

\medskip

This implies that $\mathscr{S}(\mathrm{id})(\epsilon_V)$ is also a monomorphism. Using the fact that the functor $\mathscr{S}(\mathrm{id})(-)$ is faithful, and therefore reflects monomorphisms, we conclude that $\epsilon_V$ must also be a monomorphism.

\medskip

Therefore the proposition is true on free complete dg $\mathcal{C}$-algebras. Any dg $\mathcal{C}$-algebra can be written as a split coequalizer of free dg $\mathcal{C}$-algebras. Both $\mathrm{U}^{\mathsf{op}} \cdot (-)^\vee$ and $(-)^* \cdot \mathrm{U}$ send these split coequalizers to split equalizers of dg modules, therefore the monomorphism $\epsilon$ extends to all dg $\mathcal{C}$-algebras.
\end{proof}

\begin{Remark}
Given a complete dg $\mathcal{C}$-algebra $B$ which is degree-wise finite dimensional, it is not clear to us whether its linear dual $B^*$ admits a dg $\mathcal{C}$-coalgebra structure or not. Furthermore, determining when one has that $B^\vee \cong B^*$ do not seem to be trivial. Should the interested reader find something in this direction, we would very much appreciate it.
\end{Remark}

\begin{Remark}
The results of this subsection are still valid if one does not suppose that $\mathcal{C}$ is conilpotent, with the exemption of results which concern the canonical filtration of $\C$-algebras and their "completeness". 
\end{Remark}

\subsection{The algebraic duality square} Let $\alpha: \mathcal{C} \longrightarrow \mathcal{P}$ be a twisting morphism. Using the two adjunctions constructed so far, one can interrelate the classical Bar-Cobar adjunction relative to $\alpha$ with the complete Bar-Cobar adjunction relative to $\alpha$ using a square of commuting adjunctions.

\begin{theorem}[Duality square]\label{thm: magical square}
The square of adjunctions 
\[
\begin{tikzcd}[column sep=5pc,row sep=5pc]
\left(\mathsf{dg}~\mathcal{P}\text{-}\mathsf{alg}\right)^{\mathsf{op}} \arrow[r,"\mathrm{B}_\alpha^{\mathsf{op}}"{name=B},shift left=1.1ex] \arrow[d,"(-)^\circ "{name=SD},shift left=1.1ex ]
&\left(\mathsf{dg}~\mathcal{C}\text{-}\mathsf{coalg}\right)^{\mathsf{op}} \arrow[d,"(-)^*"{name=LDC},shift left=1.1ex ] \arrow[l,"\Omega_\alpha^{\mathsf{op}}"{name=C},,shift left=1.1ex]  \\
\mathsf{dg}~\mathcal{P}\text{-}\mathsf{coalg} \arrow[r,"\widehat{\Omega}_\alpha "{name=CC},shift left=1.1ex]  \arrow[u,"(-)^*"{name=LD},shift left=1.1ex ]
&\mathsf{dg}~\mathcal{C}\text{-}\mathsf{alg}^{\mathsf{comp}}~, \arrow[l,"\widehat{\mathrm{B}}_\alpha"{name=CB},shift left=1.1ex] \arrow[u,"(-)^\vee"{name=TD},shift left=1.1ex] \arrow[phantom, from=SD, to=LD, , "\dashv" rotate=0] \arrow[phantom, from=C, to=B, , "\dashv" rotate=-90]\arrow[phantom, from=TD, to=LDC, , "\dashv" rotate=0] \arrow[phantom, from=CC, to=CB, , "\dashv" rotate=-90]
\end{tikzcd}
\] 
commutes in the following sense: right adjoints going from the top right to the bottom left are naturally isomorphic.
\end{theorem}

\begin{proof}
Let $V$ be a graded module. There is an isomorphism
\[
\left(\mathscr{S}(\mathcal{P})(V)\right)^\circ \cong \mathscr{C}(\PP)(V^*)
\]
of graded $\PP$-coalgebras. Let $(C,\Delta_C,d_C)$ be a dg $\mathcal{C}$-coalgebra, one can check by direct inspection that the above isomorphism commutes with the differential and therefore gives an isomorphism 
\[
\left(\Omega_\alpha(C)\right)^\circ \cong \widehat{\mathrm{B}}_\alpha(C^*)
\]
of dg $\PP$-coalgebras. This isomorphism is manifestly natural in $C$.
\end{proof}

\begin{Remark}
Let $(D,\Delta_D,d_D)$ be a dg $\PP$-coalgebra. Then we have an isomorphism
\[
\mathrm{B}_\alpha(D^*) \cong \left(\widehat{\Omega}_\alpha(D) \right)^\vee
\]
of dg $\mathcal{C}$-coalgebras which is natural in $D$, which is given by the mate of the isomorphism constructed in the proof of Theorem \ref{thm: magical square}. 
\end{Remark}

\begin{Proposition}
There is a natural monomorphism
\[
\zeta: \widehat{\Omega}_\alpha \cdot (-)^\circ \longrightarrow  (-)^* \cdot \mathrm{B}_\alpha~.
\]
of complete dg $\C$-algebras.
\end{Proposition}

\begin{proof}
This monomorphism is build using the monomorphisms constructed in Proposition \ref{prop: natural mono for Sweedler dual} and Proposition \ref{prop: natural mono for topo dual}.
\end{proof}

\begin{Proposition}\label{prop: finite dual commutes}
Let $(A,\gamma_A, d_A)$ be a dg $\PP$-algebra degree-wise finite dimensional. There is an isomorphism 
\[
\widehat{\Omega}_\alpha(A^*) \cong \left(\mathrm{B}_\alpha(A)\right)^*
\]
of complete dg $\mathcal{C}$-algebras. 
\end{Proposition}

\begin{proof}
Let $V$ be a graded module which is degree-wise finite dimensional. There is an isomorphism
\[
\left(\mathscr{S}(\mathcal{C})(V)\right)^* \cong \widehat{\mathscr{S}}^c(\C)(V^*)
\]
\vspace{0.3pc}

of complete graded $\mathcal{C}$-algebras. Let $(A,\gamma_A, d_A)$ be a dg $\PP$-algebra degree-wise finite dimensional, one can check that the above isomorphism commutes with the differentials.
\end{proof}

\begin{Remark}
The subcategory of dg $\PP$-algebras which satisfy the \textit{Beck-Chevalley condition} with respect to the duality square of Theorem \ref{thm: magical square} is exactly given by the subcategory of degree-wise finite dimensional dg $\PP$-algebras.
\end{Remark}

\subsection{Homotopical duality square} In this subsection, we show that the duality square of Theorem \ref{thm: magical square} behaves well with respect to model structures. We now restrict to the case where $\PP = \Omega \C$ and where the curved twisting morphism considered is the canonical twisting morphism $\iota: \C \longrightarrow \Omega \C$, in order to ensure the existence of a model category structure on the category of dg $\PP$-coalgebras. 

\begin{lemma}\label{lemma: Sweedler dual is a Quillen adjunction}
The adjunction 
\[
\begin{tikzcd}[column sep=7pc,row sep=3pc]
            \mathsf{dg}~\Omega \C\text{-}\mathsf{coalg} \arrow[r, shift left=1.1ex, "(-)^*"{name=F}] &\left(\mathsf{dg}~\Omega \C\text{-}\mathsf{alg}\right)^{\mathsf{op}} ~. \arrow[l, shift left=.75ex, "(-)^\circ"{name=U}]
            \arrow[phantom, from=F, to=U, , "\dashv" rotate=-90]
\end{tikzcd}
\]
is a Quillen adjunction. 
\end{lemma}

\begin{proof}
The left adjoint $(-)^*$ sends degree-wise monomorphisms to degree-wise epimorphisms. Thus it preserves cofibrations. It also preserves quasi-isomorphisms. Therefore we have a Quillen adjunction.
\end{proof}

\begin{theorem}[Homotopical properties of the duality square]\label{thm: homotopical magical square}
All the adjunctions in the square 
\[
\begin{tikzcd}[column sep=5pc,row sep=5pc]
\left(\mathsf{dg}~\Omega \C\text{-}\mathsf{alg}\right)^{\mathsf{op}} \arrow[r,"\mathrm{B}_\iota^{\mathsf{op}}"{name=B},shift left=1.1ex] \arrow[d,"(-)^\circ "{name=SD},shift left=1.1ex ]
&\left(\mathsf{dg}~\mathcal{C}\text{-}\mathsf{coalg}\right)^{\mathsf{op}} \arrow[d,"(-)^*"{name=LDC},shift left=1.1ex ] \arrow[l,"\Omega_\iota^{\mathsf{op}}"{name=C},,shift left=1.1ex]  \\
\mathsf{dg}~\Omega \C\text{-}\mathsf{coalg} \arrow[r,"\widehat{\Omega}_\iota "{name=CC},shift left=1.1ex]  \arrow[u,"(-)^*"{name=LD},shift left=1.1ex ]
&\mathsf{dg}~\mathcal{C}\text{-}\mathsf{alg}^{\mathsf{comp}}~, \arrow[l,"\widehat{\mathrm{B}}_\iota"{name=CB},shift left=1.1ex] \arrow[u,"(-)^\vee"{name=TD},shift left=1.1ex] \arrow[phantom, from=SD, to=LD, , "\dashv" rotate=0] \arrow[phantom, from=C, to=B, , "\dashv" rotate=-90]\arrow[phantom, from=TD, to=LDC, , "\dashv" rotate=0] \arrow[phantom, from=CC, to=CB, , "\dashv" rotate=-90]
\end{tikzcd}
\] 
are Quillen adjunctions.
\end{theorem}

\begin{proof}
The only thing left to check is that the adjunction 
\[
\begin{tikzcd}[column sep=7pc,row sep=3pc]
           \mathsf{dg}~\mathcal{C}\text{-}\mathsf{alg}^{\mathsf{comp}} \arrow[r,"(-)^\vee "{name=F}, shift left=1.1ex] 
           &\left(\mathsf{dg}~\mathcal{C}\text{-}\mathsf{coalg}\right)^{\mathsf{op}}~. \arrow[l, shift left=.75ex, "(-)^*"{name=U}]
            \arrow[phantom, from=F, to=U, , "\dashv" rotate=-90]
\end{tikzcd}
\]
is indeed a Quillen adjunction, where the model structure considered on dg $\C$-coalgebras is obtained by transfer along the adjunction $\Omega_\iota \dashv \mathrm{B}_\iota$ and where the model structure considered on the category of complete dg $\C$-algebras is obtained by transfer along the adjunction $\widehat{\Omega}_\iota \dashv \widehat{\mathrm{B}}_\iota$.

\medskip

Let us check that $(-)^*$ is a right Quillen functor. It sends monomorphisms to epimorphisms, thus preserves fibrations. Since every dg $\mathcal{C}$-coalgebra is cofibrant (fibrant in the opposite category), we are left to show that $(-)^*$ preserves weak equivalences of dg $\C$-coalgebras. Let 
\[
f: C_1 \qi C_2
\]
be a weak equivalence of curved $\C$-coalgebras, that is,
\[
\Omega_\iota(f): \Omega_\iota C_1 \qi \Omega_\iota C_2 
\]
is a quasi-isomorphism of dg $\Omega \C$-algebras, by Lemma \ref{lemma: Sweedler dual is a Quillen adjunction}, we know that the Sweedler dual functor $(-)^\circ$ is a right Quillen functor. Therefore it preserves weak-equivalences between fibrant objects (i.e: quasi-isomorphisms between cofibrant dg $\Omega \C$-algebras). Thus 
\[
(\Omega_\iota(f))^\circ: (\Omega_\iota C_1)^\circ \qi (\Omega_\iota C_2)^\circ 
\]
is a quasi-isomorphism of dg $\Omega \C$-coalgebras. Using the commutativity of the square, we get that
\[
\widehat{\mathrm{B}}_\iota(f^*): \widehat{\mathrm{B}}_\iota C_1^* \longrightarrow \widehat{\mathrm{B}}_\iota C_2^*
\]
is also a quasi-isomorphism of dg $\Omega \C$-coalgebras. Therefore $f^*: C_1^* \longrightarrow C_2^*$ is a weak equivalence in the model category of complete dg $\C$-algebras. 
\end{proof}

\begin{Proposition}
Let $(A,\gamma_A, d_A)$ be a dg $\Omega \C$-algebra whose homology is degree-wise finite dimensional. Then the derived unit of adjunction

\[
\mathbb{R}(\eta_A): A \qi \left(\mathbb{R}(A^\circ)\right)^*
\]
\vspace{0.2pc}

is a quasi-isomorphism. Therefore the functor $(-)^\circ$ is homotopically fully faithful on the full sub-$\infty$-category of dg $\Omega \C$-algebras with degree-wise finite dimensional homology.
\end{Proposition}

\begin{proof}
Let $A$ be a dg $\Omega \C$-algebra such that its homology is degree-wise finite dimensional. Since $\kk$ is a field of characteristic zero, using the Homotopy Transfer Theorem, one can replace $A$ by its homology $\mathrm{H}_*(A)$ in order to compute this derived unit of adjunction. Indeed, there is a direct quasi-isomorphism 

\[
\Omega(i_A): \Omega_\iota \mathrm{B}_\iota \mathrm{H}_*(A) \qi \Omega_\iota \mathrm{B}_\iota A
\]

of dg $\Omega \C$-algebras between their cofibrant resolutions given by the Bar-Cobar adjunction relative to $\iota$. Here $i_A: \mathrm{B}_\iota \mathrm{H}_*(A) \qi \mathrm{B}_\iota A$ is the $\infty$-quasi-isomorphism given by the choice of a homotopy contraction between $A$ and $\mathrm{H}_*(A)$. Therefore we obtain the following commutative square

\[
\begin{tikzcd}[column sep=4pc,row sep=4pc]
\Omega_\iota \mathrm{B}_\iota \mathrm{H}_*(A) \arrow[d,"\mathbb{R}(\eta_{\mathrm{H}_*(A)})"] \arrow[r,"\Omega(i_A)"]
&\Omega_\iota \mathrm{B}_\iota A    \arrow[d,"\mathbb{R}(\eta_A)"]        \\
((\Omega_\iota \mathrm{B}_\iota \mathrm{H}_*(A))^\circ)^* \arrow[r,"((\Omega(i_A))^\circ)^*"]
&((\Omega_\iota \mathrm{B}_\iota A)^\circ)^*~,
\end{tikzcd}
\]

Using 2 out of 3, we conclude that $\mathbb{R}(\eta_A)$ is a quasi-isomorphism of dg $\Omega \C$-algebras if and only if $\mathbb{R}(\eta_{\mathrm{H}_*(A)})$ is. We compute that 

\[
\left(\Omega_\iota \mathrm{B}_\iota \mathrm{H}_*(A)\right)^\circ \cong \widehat{\mathrm{B}}_\iota \left(\mathrm{B}_\iota \mathrm{H}_*(A)\right)^* \cong \widehat{\mathrm{B}}_\iota \widehat{\Omega}_\iota (\mathrm{H}_*(A))^*~,
\]

first using the commutativity of the square of Theorem \ref{thm: homotopical magical square} and then by applying Proposition \ref{prop: finite dual commutes} since $\mathrm{H}_*(A)$ is degree-wise finite dimensional. And the unit of the complete Bar-Cobar adjunction $\widehat{\Omega}_\iota \dashv \widehat{\mathrm{B}}_\iota$ 

\[
\eta_{(\mathrm{H}_*(A))^*}: (\mathrm{H}_*(A))^* \qi \widehat{\mathrm{B}}_\iota \widehat{\Omega}_\iota (\mathrm{H}_*(A))^*
\]

is a quasi-isomorphism of dg $\Omega \C$-coalgebras since this adjunction is a Quillen equivalence. We can apply the linear dual functor to this quasi-isomorphism, resulting in

\[
\mathrm{H}_*(A) \cong (\mathrm{H}_*(A))^{**} \qi \left(\widehat{\mathrm{B}}_\iota \widehat{\Omega}_\iota (\mathrm{H}_*(A))^*\right)^* =\left(\mathbb{R}(\mathrm{H}_*(A)^\circ)\right)^*~,
\]

using again that $\mathrm{H}_*(A)$ is degree-wise finite dimensional. Using 2 out of 3, we conclude that 

\[
\mathbb{R}(\eta_{\mathrm{H}_*(A)}): \mathrm{H}_*(A) \qi \left(\mathbb{R}(\mathrm{H}_*(A)^\circ)\right)^*
\]

must also be a quasi-isomorphism of dg $\Omega \C$-algebras.
\end{proof}

\section{Examples: contramodules and absolute associative algebras}\label{Section: Absolute algebras and contramodules}

\subsection{Contramodules}
The notion of a \textit{contramodule} over a dg counital coassociative coalgebra appears naturally as the dual definition of a comodule. This notion was introduced in the seminal work of S. Eilenberg and J. C. Moore in \cite{eilenbergmoore65}. After being almost forgotten for a number of years, this notion reemerged in the work of L. Positseski, see \cite{positselski2021contramodules} for a extensive account. In this subsection, we show that contramodules appear as a particular case of the new, more general notion of algebras over cooperads. This allows us to subsume the theory of contramodules, obtaining many examples of interest for this new notion. To the best of our knowledge, the previous duality square also provides a new result for the theory of contramodules.

\begin{Definition}[Contramodule]
Let $(C, \Delta, \epsilon, d_C)$ be a dg counital coassociative coalgebra. A \textit{dg} $C$-\textit{contramodule} $M$ is the data $(M,\gamma_M,d_M)$ of a  dg module $(M,d_M)$ together with a morphism of dg modules 

\[
\gamma_M: \mathrm{Hom}(C,M) \longrightarrow M~,
\]
\vspace{0.3pc}

such that the following diagrams commute

\[
\begin{tikzcd}[column sep=3pc,row sep=3pc]
\mathrm{Hom}(C \otimes C, M) \arrow[r,"\Delta^* "] \arrow[d,"\varsigma",swap]
& \mathrm{Hom}(C, M) \arrow[dd,"\gamma_M"] \\
\mathrm{Hom}(C, \mathrm{Hom}(C, M)) \arrow[d,"(\gamma_M)_*",swap]
& \\
\mathrm{Hom}(C, M) \arrow[r,"\gamma_M"]
&M ~,
\end{tikzcd} \quad \quad
\begin{tikzcd}[column sep=3pc,row sep=3pc]
M \cong \mathrm{Hom}(\kk, M)  \arrow[r,"\epsilon^*"] \arrow[rd,"\cong",swap]
&\mathrm{Hom}(C, M) \arrow[d,"\gamma_M"]\\
&M~.
\end{tikzcd}
\]
\end{Definition}

\begin{Remark}
If $(C, \Delta, \epsilon, d_C)$ is a dg counital coassociative coalgebra, one can show that the endofunctor
\[
\mathrm{Hom}(C,-): \mathsf{dg}~\mathsf{mod} \longrightarrow \mathsf{dg}~\mathsf{mod}
\]
admits a monad structure such that dg $C$-contramodules are algebras over this monad. In fact, it admits \textit{two} monad structures, depending on the choice of $\varsigma$ in the above definition. Indeed, there are two isomorphisms
\[
\mathrm{Hom}(C \otimes C, M) \cong \mathrm{Hom}(C, \mathrm{Hom}(C, M)) ~,
\]
depending on the choice of either the left or the right factor in $C \otimes C$. One should speak of \textit{left} or \textit{right} $C$-contramodules in each of these cases. We omit these subtleties when possible.
\end{Remark}

\begin{Example}
Let $(C, \Delta, \epsilon, d_C)$ be a dg counital coassociative coalgebra and $(D,d_D)$ be a dg module. The \textit{free} dg $C$-contramodule on $D$ is given by $\mathrm{Hom}(C,D)$ with the obvious structure maps. 
\end{Example}

\begin{Example}
Let $(C, \Delta, \epsilon, d_C)$ be a dg counital coassociative coalgebra and let $(V,\Delta_V,d_V)$ be a right dg $C$-comodule. For any dg module $(D,d_D)$, the dg module $\mathrm{Hom}(V,D)$ inherits a canonical left dg $C$-contramodule structure. 
\end{Example}

\begin{Proposition} \leavevmode

\medskip 

\begin{enumerate}
\item The data of a dg counital coassociative coalgebra $(C, \Delta, \epsilon, d_C)$ is equivalent to the data of a dg cooperad $(\mathcal{C},\Delta_{\mathcal{C}},\epsilon, d_{\C})$ whose dg $\mathbb{S}$-module $\mathcal{C}$ is concentrated in arity one.

\medskip

\item In the case above, the definition of a dg $C$-contramodule is equivalent to the definition of a dg $\C$-algebra.
\end{enumerate}
\end{Proposition}

\begin{proof}
It is a straightforward computation from the definitions.
\end{proof}

This equivalence is a great source of examples for algebras over cooperads. Moreover, it is also a great source of counter-examples.

\begin{Example}[Contramodules over formal power series]
Consider the cofree conilpotent coalgebra $\kk[t]^c$ cogenerated by a single element $t$ of degree zero. In this case, one can see that the data of $\kk[t]^c$-contramodule structure on $M$ is equivalent to the data of a structural map
\[
\gamma_M: \kk[[t]] ~\widehat{\otimes}~ M \cong \mathrm{Hom}(\kk[t]^c,M) \longrightarrow M~.
\]
which satisfies unital and associativity conditions with respect to the algebra of formal power series in one variable $\kk[[t]]$. Therefore, for any formal power series
\[
\sum_{n \geq 0} m_n \otimes t^n ~~\quad~~\text{in} ~~\quad~~ \kk[[t]] ~\widehat{\otimes}~ M~,
\]
the structural morphism $\gamma_M$ gives a well-defined image $\gamma_M\left(\sum_{n \geq 0} m_n \otimes t^n \right)$ in $M$, without \textit{presupposing} any topology on $M$.
\end{Example}

Let $(C, \Delta, \epsilon, \eta, d_C)$ be a dg \textit{conilpotent} counital coassociative coalgebra and let $\overline{C}$ be its coaugmentation ideal. Any dg $C$-contramodule $M$ admits a \textit{canonical decreasing filtration} given by 
\[
\mathrm{W}_\omega M \coloneqq \mathrm{Im}\left(\gamma_M^{ \geq \omega}: \mathrm{Hom}\left(\overline{C}/\mathscr{R}_\omega \overline{C}, M\right) \longrightarrow M \right)~,
\]
where $\mathscr{R}_\omega \overline{C}$ denotes the $\omega$-stage of the coradical filtration. One says that $M$ is \textit{complete} if the canonical morphism of dg $C$-contramodules 
\[
\varphi_M: M \twoheadrightarrow \lim_{\omega} M/\mathrm{W}_\omega M
\]
is an isomorphism. Again, this is a particular case of the more general notion of canonical filtration for algebras over conilpotent cooperads given in Definition \ref{def: canonical filtration C alg}.

\begin{Counterexample}
Consider again the cofree conilpotent coalgebra $\kk[t]^c$ cogenerated by a single element $t$ of degree zero. There exists a $\kk[t]^c$-contramodule $M$ and a family of elements $\{m_n\}_{n \geq 0}$ in $M$ such that 

\[
\gamma_M\left(m_n \otimes t^n \right) = 0 \quad \text{for all} \quad n\geq 0 \quad \text{and} \quad \gamma_M\left(\sum_{n \geq 0} m_n \otimes t^n \right) \neq 0~.
\]
\vspace{0.2pc}

In particular, the topology induced by the canonical filtration on $M$ is not complete. See \cite[Section 1.5]{positselski2021contramodules} for this construction.
\end{Counterexample}

Let $(A, \mu_A, \eta, d_A)$ be a dg unital associative algebra (viewed as a dg operad concentrated in arity one), and let $(C, \Delta, \epsilon, \eta, d_C)$ be a dg conilpotent counital coassociative algebra (viewed as a conilpotent dg cooperad concentrated in arity one). Let $\alpha: C \longrightarrow A$ be a twisting morphism between them. It induces a first adjunction 

\[
\begin{tikzcd}[column sep=7pc,row sep=3pc]
           \mathsf{dg}~C\text{-}\mathsf{comod}^{l} \arrow[r,"A \otimes_{\alpha} - "{name=F}, shift left=1.1ex] 
           &\mathsf{dg}~A\text{-}\mathsf{mod}^{l} ~. \arrow[l, shift left=.75ex, "C \otimes_{\alpha} - "{name=U}]
            \arrow[phantom, from=F, to=U, , "\dashv" rotate=-90]
\end{tikzcd}
\]
\vspace{0.3pc}

between the category of left dg $C$-comodules and the category of left dg $A$-module, where $\otimes_\alpha$ denotes the twisted tensor product, see \cite[Section 2.1]{LodayVallette12} for an account of this twisted tensor product. One can check that this adjunction is nothing but the classical Bar-Cobar adjunction relative to $\alpha$ when we view $A$ as a dg operad concentrated in arity one and $C$ as a dg conilpotent cooperad concentrated in arity one. When we endow the category of left dg $A$-module with the \textit{projective model structure} and the category of left dg $C$-comodules with the transferred model structure, this becomes a Quillen adjunction. 

\begin{lemma}
Let $(A, \mu_A, \eta, d_A)$ be a dg unital associative algebra viewed as a dg operad concentrated in arity one. The category of dg $A$-coalgebras is isomorphism to the category of right dg $A$-modules. Furthermore, the cofree right dg $A$-module functor is given by $\mathrm{Hom}(A,-)$.
\end{lemma}

\begin{proof}
It is a straightforward computation.
\end{proof}

\begin{Remark}
When $A$ is the group algebra $\kk[G]$ for some group $G$, the functor $\mathrm{Hom}(\kk[G],-)$ is usually called the \textit{coinduced representation}. See Section \ref{appendice: induction et coinduction}.
\end{Remark}

The twisting morphism $\alpha: C \longrightarrow A$ induces a second adjunction 

\[
\begin{tikzcd}[column sep=7pc,row sep=3pc]
           \mathsf{dg}~A\text{-}\mathsf{mod}^{r}\arrow[r," \mathrm{Hom}^{\alpha}(C \text{,}-) "{name=F}, shift left=1.1ex] 
           &\mathsf{dg}~C\text{-}\mathsf{contra}^{r} ~. \arrow[l, shift left=.75ex, "\mathrm{Hom}^{\alpha}(A \text{,}-)"{name=U}]
            \arrow[phantom, from=F, to=U, , "\dashv" rotate=-90]
\end{tikzcd}
\]
\vspace{0.3pc}

between the category of right dg $C$-contramodules and the category of right dg $A$-modules, where $\mathrm{Hom}^\alpha(-,-)$ denotes the twisted hom space, see \cite[Section 2.1]{LodayVallette12} for an account of this twisted hom space. Once again, one can check that this adjunction corresponds to the complete Bar-Cobar adjunction relative to $\alpha$. When we endow the category of right dg $A$-module with the \textit{injective model structure} and the category of right dg $C$-contramodules with the transferred model structure, this becomes a Quillen adjunction. 

\medskip

The two structures defined above are compatible under the following duality square.

\begin{Proposition}
The following square diagram 

\[
\begin{tikzcd}[column sep=5pc,row sep=5pc]
\left(\mathsf{dg}~A\text{-}\mathsf{mod}^{l}\right)^{\mathsf{op}} \arrow[r,"C \otimes_{\alpha} -"{name=B},shift left=1.1ex] \arrow[d,"(-)^\circ "{name=SD},shift left=1.1ex ]
&\left(\mathsf{dg}~C\text{-}\mathsf{comod}^{l}\right)^{\mathsf{op}} \arrow[d,"(-)^*"{name=LDC},shift left=1.1ex ] \arrow[l,"A \otimes_{\alpha} -"{name=C},,shift left=1.1ex]  \\
\mathsf{dg}~A\text{-}\mathsf{mod}^{r} \arrow[r,"\mathrm{Hom}^{\alpha}(C \text{,}-)"{name=CC},shift left=1.1ex]  \arrow[u,"(-)^*"{name=LD},shift left=1.1ex ]
&\mathsf{dg}~C\text{-}\mathsf{contra}^{r} \arrow[l,"\mathrm{Hom}^{\alpha}(A \text{,}-)"{name=CB},shift left=1.1ex] \arrow[u,"(-)^\vee"{name=TD},shift left=1.1ex] \arrow[phantom, from=SD, to=LD, , "\dashv" rotate=0] \arrow[phantom, from=C, to=B, , "\dashv" rotate=-90]\arrow[phantom, from=TD, to=LDC, , "\dashv" rotate=0] \arrow[phantom, from=CC, to=CB, , "\dashv" rotate=-90]
\end{tikzcd}
\] 

is commutative and is made of Quillen adjunctions.
\end{Proposition}

\begin{proof}
This is a direct consequence of Theorem \ref{thm: homotopical magical square}.
\end{proof}

One can view the adjunction 

\[
\begin{tikzcd}[column sep=7pc,row sep=3pc]
           \left(\mathsf{dg}~C\text{-}\mathsf{comod}^{l}\right)^{\mathsf{op}} \arrow[r,"(-)^*"{name=F}, shift left=1.1ex] 
           &\mathsf{dg}~C\text{-}\mathsf{contra}^{r} ~. \arrow[l, shift left=.75ex, "(-)^\vee"{name=U}]
            \arrow[phantom, from=F, to=U, , "\dashv" rotate=90]
\end{tikzcd}
\]

as a contravariant version of the \textit{co-contra correspondence} constructed of \cite{PositselskiTwoKinds}. Indeed, this adjunction identifies finitely generated cofree $C$-comodules and finitely generated $C$-contramodules: let $M$ be a dg module which is finite dimensional degree-wise, we have that 

\[
(C \otimes M)^* \cong \mathrm{Hom}(C,M^*) \quad \text{and} \quad \left(\mathrm{Hom}(C,M) \right)^\vee \cong C \otimes M~.
\]
\vspace{0.2pc}

Furthermore, when both of these categories are endowed with the canonical model structure transferred using the twisting morphism $\iota: C \longrightarrow \Omega C$, then this correspondences becomes a correspondence between the finitely generated fibrant-cofibrant objects of each of these categories. Therefore it gives a contravariant version of the \textit{derived co-contra correspondence}, see again \cite{PositselskiTwoKinds} for more on this.

\begin{Remark}
The data of a curved cooperad concentrated in arity one is equivalent to the data of a curved counital coassociative coalgebra. This subsection can be generalized \textit{mutatis mutandis} to the curved setting.
\end{Remark}

\subsection{Absolute associative algebras and absolute Lie algebras}
In this subsection, we explore the notion of algebras over the conilpotent cooperad $\mathcal{A}ss^*$, which we call \textit{absolute associative algebras}. We compare this category with the standard notion of non-unital associative algebras, encoded by the operad $\mathcal{A}ss$. Furthermore, we also introduce \textit{absolute Lie algebras}, encoded by the cooperad $\mathcal{L}ie^*$, and construct the universal enveloping absolute algebra functor. These two examples provide supplementary intuition on the notion of an algebra over a cooperad, in one of its simplest cases. Many of the subsequent results can be generalizes, \textit{mutatis mutandis}, to any well-behaved binary cooperad.

\begin{Definition}[dg absolute associative algebra]
A \textit{dg absolute associative algebra} $(A,\gamma_A,d_A)$ is the data of a dg $\mathcal{A}ss^*$-algebra. 
\end{Definition}

Let us unravel this definition. Recall that a dg $\mathcal{A}ss^*$-algebra $(A,\gamma_A,d_A)$ is the data of a dg module $(A,d_A)$ together with a structural morphism of dg modules

\[
\gamma_A: \displaystyle \prod_{n \geq 0} \mathrm{Hom}_{\mathbb{S}_n}(\mathcal{A}ss^*(n), A^{\otimes n}) \longrightarrow A~,
\]

which satisfies the axioms of Definition \ref{def: C algebra}. 

\begin{lemma}
There is an isomorphism of dg modules

\[
\displaystyle \prod_{n \geq 0} \mathrm{Hom}_{\mathbb{S}_n}(\mathcal{A}ss^*(n), A^{\otimes n}) \cong \displaystyle \prod_{n \geq 1} A^{\otimes n}~.
\]
\end{lemma}

\begin{proof}
The $\mathbb{S}_n$-module $\mathcal{A}ss^*(n)$ is given by the regular representation of $\mathbb{S}_n$, for $n \geq 1$, and by $0$, when $n=0$.
\end{proof}

\begin{Remark}
This implies that there is a monad structure on the endofunctor 

\[
\prod_{n \geq 1} (-)^{\otimes n}: \mathsf{dg}~\mathsf{mod} \longrightarrow \mathsf{dg}~\mathsf{mod}~,
\]

such that a dg absolute associative algebra is the data of an algebra over this monad.
\end{Remark}

The structural map $\gamma_A$ of an dg absolute associative algebra $A$ is equivalently given by a morphism of dg modules 

\[
\gamma_A: \displaystyle \prod_{n \geq 1} A^{\otimes n} \longrightarrow A~.
\]

It associates to any series
\[
\displaystyle \sum_{n \geq 1} \sum_{i \in \mathrm{I}_n } a_1^{(i)} \otimes \cdots \otimes a_n^{(i)}~,
\]
where the $a_j$ are elements in $A$ and where $\mathrm{I}_n$ is a finite set, a well-defined element
\[
\gamma_A \left(\displaystyle \sum_{n \geq 1} \sum_{i \in \mathrm{I}_n } a_1^{(i)} \otimes \cdots \otimes a_n^{(i)} \right) \quad \text{in} \quad A~,
\]
without \textbf{presupposing any topology} on the dg module $A$. 

\begin{Remark}
For a general dg absolute associative algebra $(A,\gamma_A,d_A)$, notice that 
\[
\gamma_A \left(\displaystyle \sum_{n \geq 1} \sum_{i \in \mathrm{I}_n } a_1^{(i)} \otimes \cdots \otimes a_n^{(i)} \right) \neq \displaystyle \sum_{n \geq 1} \sum_{i \in \mathrm{I}_n } \gamma_A\left(a_1^{(i)} \otimes \cdots \otimes a_n^{(i)}\right)~, 
\]
as the latter sum is not even well-defined in $A$.
\end{Remark}

\textbf{Differential condition:}
The condition that $\gamma_A$ is a morphism of dg modules can be rewritten as 

\begin{equation*}\label{condition: dg}
\begin{small}
\gamma_A \left(\displaystyle \sum_{n \geq 1} \sum_{i \in \mathrm{I}_n } \sum_{j = 0}^n (-1)^j  a_1^{(i)} \otimes \cdots \otimes d_A\left(a_j^{(i)}\right)\otimes \cdots \otimes a_n^{(i)} \right) = d_A \left(\gamma_A \left(\displaystyle \sum_{n \geq 1} \sum_{i \in \mathrm{I}_n } a_1^{(i)} \otimes \cdots \otimes a_n^{(i)} \right) \right)~,
\end{small}
\end{equation*}

for any series in  $\prod_{n \geq 1} A^{\otimes n}$.

\medskip

\textbf{Associativity condition:}
The condition that the structural morphism $\gamma_A$ defines an algebra over the monad $\prod_{n \geq 1} (-)^{\otimes n}$ can be rewritten as 

\begin{align*}\label{condition: associativité}
\gamma_A \left(\displaystyle \sum_{k \geq 1} \gamma_A \left(\displaystyle \sum_{i_1 \geq 1} \sum_{j_1 \in \mathrm{I}_{i_1} } a_1^{(j_1)} \otimes \cdots \otimes a_{i_1}^{(j_1)} \right) \otimes \cdots \otimes \gamma_A \left(\displaystyle \sum_{i_k \geq 1} \sum_{j_k \in \mathrm{I}_{i_k} } a_1^{(j_k)} \otimes \cdots \otimes a_{i_k}^{(j_k)} \right) \right)\\
=\gamma_A \left(\displaystyle \sum_{n \geq 1} \sum_{k \geq 1} \sum_{i_1 + \cdots + i_k = n} \sum_{j_1 \in  \mathrm{I}_{i_1}~, \cdots~, j_k \in \mathrm{I}_{i_k}} a_1^{(j_1)} \otimes \cdots \otimes a_{i_k}^{(j_1)} \otimes \cdots \otimes a_{n - i_k}^{(j_k)} \otimes \cdots \otimes a_{n}^{(j_k)} \right)~.
\end{align*}

\begin{Example}
Let $(V,d_V)$ be a dg module. The \textit{free} dg absolute associative algebra on $V$ is given by 
\[
\overline{\mathcal{T}}^{~\wedge}(V) \coloneqq \prod_{n \geq 1} V^{\otimes n}~,
\]

that is, the completed non-unital tensor algebra on $V$. If $V$ is a $\kk$-vector space of dimension $n$, this algebra is the algebra of non-commutative formal power series in $n$ variables without constant terms.
\end{Example}

\begin{lemma}\label{lemma: morphism of monads}
The inclusion 

\[
\iota: \bigoplus_{n \geq 1} (-)^{\otimes n} \hookrightarrow \prod_{n \geq 1} (-)^{\otimes n}
\]

defines a morphism of monads in the category of dg modules.
\end{lemma}

\begin{proof}
It is straightforward to check. 
\end{proof}

\begin{Proposition}\label{prop: adjunction absolute envelope}
There is an adjunction 

\[
\begin{tikzcd}[column sep=7pc,row sep=3pc]
            \mathsf{dg}~\mathsf{abs}~\mathsf{assoc}\text{-}\mathsf{alg}\arrow[r,"\mathrm{Res}"{name=F}, shift left=1.1ex] 
           &\mathsf{dg}~\mathsf{assoc}\text{-}\mathsf{alg}~, \arrow[l, shift left=.75ex, "\mathrm{Abs}"{name=U}]
            \arrow[phantom, from=F, to=U, , "\dashv" rotate=90]
\end{tikzcd}
\]

between the category of dg absolute associative algebras and the category of dg associative algebras, where the left adjoint $\mathrm{Abs}$ is called the \textit{absolute envelope functor}.
\end{Proposition}

\begin{proof}
This follows from Lemma \ref{lemma: morphism of monads}.
\end{proof}

Let us describe these functors. Let $(A,\gamma_A,d_A)$ be a dg absolute associative algebra. Its restriction algebra $(\mathrm{Res}(A), \mathrm{Res}(\gamma_A), d_A)$ is given by the dg module $(A,d_A)$ together with the structural map

\[
\begin{tikzcd}
\displaystyle \mathrm{Res}(\gamma_A): \bigoplus_{n \geq 1} A^{\otimes n} \arrow[r,"\iota_A"]
&\displaystyle \prod_{n \geq 1} A^{\otimes n} \arrow[r,"\gamma_A"]
&A~,
\end{tikzcd}
\]
which endows $\mathrm{Res}(A)$ with a dg associative algebra structure. Let $(B,\gamma_B,d_B)$ be a dg associative algebra, its absolute envelope $\mathrm{Abs}(B)$ is given by the following coequalizer 

\[
\begin{tikzcd}[column sep=4pc,row sep=4pc]
\mathrm{Coeq}\Bigg(\displaystyle \prod_{n \geq 1} \Bigg(\bigoplus_{k \geq 1} B^{\otimes k}\Bigg)^{\otimes n} \arrow[r,"\prod_{n \geq 0} (\gamma_B)^{\otimes n}",shift right=1.1ex,swap]  \arrow[r,"\psi_B"{name=SD},shift left=1.1ex ]
&\prod_{n \geq 1} B^{\otimes n}\Bigg)~,
\end{tikzcd}
\]

in the category of dg modules, where $\psi_B$ is given by
\[
\begin{tikzcd}[column sep=4pc,row sep=4pc]
\displaystyle \prod_{n \geq 1} \left(\bigoplus_{k \geq 1} B^{\otimes k}\right)^{\otimes n} \arrow[r,"\prod_{n \geq 0} (\iota_B)^{\otimes n}"]
&\displaystyle\prod_{n \geq 1} \left(\prod_{k \geq 1} B^{\otimes k}\right)^{\otimes n} \arrow[r, "\gamma_{\prod_{k \geq 1} B^{\otimes k}}"]
&\displaystyle \prod_{n \geq 1} B^{\otimes n}~.
\end{tikzcd}
\]

\textbf{Completeness:}
Any dg absolute associative algebra $(A,\gamma_A,d_A)$ comes equipped with a canonical decreasing filtration given by 

\[
\mathrm{W}_\omega A \coloneqq \mathrm{Im}\left(\gamma_A^{\geq \omega}: \prod_{n \geq \omega + 1} A^{\otimes n} \longrightarrow A \right)~.
\]

This filtration is the same filtration defined in Definition \ref{def: canonical filtration C alg} since the cooperad $\mathcal{A}ss^*$ is a binary cooperad. A dg absolute associative algebra $(A,\gamma_A,d_A)$ is therefore said to be \textit{complete} if the canonical epimorphism

\[
\varphi_A: A \twoheadrightarrow \lim_{\omega} A/\mathrm{W}_\omega A
\]

is an isomorphism of dg absolute associative algebras. 

\begin{Remark}
Let $(A,\gamma_A,d_A)$ be a complete dg absolute associative algebra, then  
\[
\gamma_A \left(\displaystyle \sum_{n \geq 1} \sum_{i \in \mathrm{I}_n } a_1^{(i)} \otimes \cdots \otimes a_n^{(i)} \right) = \displaystyle \sum_{n \geq 1} \sum_{i \in \mathrm{I}_n } \gamma_A\left(a_1^{(i)} \otimes \cdots \otimes a_n^{(i)}\right)~, 
\]
using the completeness of $A$. 
\end{Remark}

Similarly, when $(B,\gamma_B,d_B)$ is a dg associative algebra, one has a canonical decreasing filtration as well given by

\[
\mathrm{F}_\omega B \coloneqq \mathrm{Im}\left(\gamma_B^{\geq \omega}: \bigoplus_{n \geq \omega + 1} B^{\otimes n} \longrightarrow B \right)~.
\]

A dg associative algebra $(B,\gamma_B,d_B)$ is said to be \textit{complete} if the canonical morphism 
\[
\lambda_B: B \longrightarrow \lim_{\omega} B/ \mathrm{F}_\omega B 
\]
is an isomorphism of dg associative algebras. Let us try to compare these two filtrations. 

\begin{Definition}[Nilpotent dg absolute associative algebra]
A dg absolute associative algebra $(A,\gamma_A,d_A)$ is said to be \textit{nilpotent} if there exists $\omega_{0} \geq 1$ such that $\mathrm{W}_{\omega_0} A = 0$. The integer $\omega_0$ is called the \textit{nilpotency degree}.
\end{Definition}

\begin{Remark}
A nilpotent dg absolute associative algebra is in particular complete.
\end{Remark}

\begin{Definition}[Nilpotent dg associative algebra]
A dg associative algebra $(B,\gamma_B,d_B)$ is said to be \textit{nilpotent} if there exists $\omega_{0} \geq 1$ such that $\mathrm{F}_{\omega_0} B = 0$. The integer $\omega_0$ is called the nilpotency degree.
\end{Definition}

\begin{Remark}
This definition coincides with the standard notion of nilpotency for associative algebras present in the literature.
\end{Remark}

\begin{Proposition}\label{prop: les nilpotentes sont des absolues}
The data of a nilpotent dg absolute associative algebra is equivalent to the data of a nilpotent dg associative algebra with the same nilpotency degree.
\end{Proposition}

\begin{proof}
Let $\omega_0 \geq 1$. A nilpotent dg absolute associative algebra $(A,\gamma_A,d_A)$ with nilpotency degree $\omega_0$ amounts to the data of an algebra over the monad $\prod_{n = 1}^{\omega_0 + 1} (-)^{\otimes n}$. On the other hand, a nilpotent dg associative algebra $(B,\gamma_B,d_B)$ amounts to the data of an algebra over the monad $\bigoplus_{n = 1}^{\omega_0 + 1} (-)^{\otimes n}$. There is an isomorphism of monads 
\[
\prod_{n = 1}^{\omega_0 + 1} (-)^{\otimes n} \cong \bigoplus_{n = 1}^{\omega_0 + 1} (-)^{\otimes n}~,
\]
since the product only involves a finite amount of terms.
\end{proof}

\begin{Remark}\label{Remark: nilpotent Yoneda}
The image of a nilpotent dg associative algebra $(B,\gamma_B,d_B)$ via the absolute envelope functor $\mathrm{Abs}$ is simply given by $(B,\gamma_B,d_B)$ considered as a dg absolute associative algebra. Indeed, we have a natural isomorphism

\[
\mathrm{Hom}_{\mathsf{dg}~\mathsf{abs}~\mathsf{assoc}\text{-}\mathsf{alg}}(B,-) \cong \mathrm{Hom}_{\mathsf{dg}~\mathsf{assoc}\text{-}\mathsf{alg}}(B, \mathrm{Res}(-))~,
\]
\vspace{0.2pc}

hence, by Yoneda lemma, $\mathrm{Abs}(B) \cong B$. 
\end{Remark}

\begin{Proposition}
Let $(A,\gamma_A,d_A)$ be a dg absolute associative algebra. There is a monomorphism of dg modules 

\[
\mathrm{F}_\omega \mathrm{Res}(A) \rightarrowtail \mathrm{W}_\omega A~.
\]

Therefore, if $A$ is a complete dg absolute associative algebra, the topology induced by the canonical filtration on $(\mathrm{Res}(A), \mathrm{Res}(\gamma_A), d_A)$ is separated. 
\end{Proposition}

\begin{proof}
It is straightforward to notice that $\mathrm{F}_\omega \mathrm{Res}(A)$ consists of finite sums of products of arity greater than $\omega +1$, hence it is included in $\mathrm{W}_\omega A$. For the second point, notice that if $(A,\gamma_A,d_A)$ is complete, then 

\[
\bigcap_{\omega \geq 0}\mathrm{F}_\omega \mathrm{Res}(A) \hookrightarrow \bigcap_{\omega \geq 0}\mathrm{W}_\omega A = \{0\}~,
\]
hence the canonical morphism of dg associative algebras 
\[
\lambda_{\mathrm{Res}(A)}: \mathrm{Res}(A) \longrightarrow \lim_{\omega} \left(\mathrm{Res}(A)/\mathrm{F}_\omega \mathrm{Res}(A)\right)
\]
is a monomorphism and the topology induced by $\mathrm{F}_\omega \mathrm{Res}(A)$ is separated.
\end{proof}

\begin{Counterexample}
Suppose $(A,\gamma_A,d_A)$ is a complete dg absolute associative algebra, then $\mathrm{Res}(A)$ might not be complete as a dg associative algebra. This is for instance the case of the free absolute associative algebra on a dg module $(V,d_V)$: its restriction algebra $\mathrm{Res}(\prod_{n \geq  1} V^{\otimes n}))$ is not complete as a dg associative algebra. Indeed, one can check that the canonical morphism $\lambda_{\prod_{n \geq 1}V^{\otimes n}}$ is not an epimorphism in this case.
\end{Counterexample}

\begin{Remark}
The above results can be generalized \textit{mutatis mutandis} to algebras over a binary cooperad $\mathcal{C}$ with homogeneous corelations. Particular examples of these include: $\mathcal{C}om^*$, which encodes dg absolute commutative algebras, and $\mathcal{L}ie^*$, which encodes dg absolute Lie algebras.
\end{Remark}

\textbf{Universal enveloping absolute algebra functor.} As with operads, a morphism of cooperads induces an adjunction between their respective categories of algebras. This allows us to construct the universal enveloping absolute algebra of a dg absolute Lie algebra. 

\begin{Definition}[dg absolute Lie algebra]
A \textit{dg absolute Lie algebra} $(\mathfrak{g},\gamma_{\mathfrak{g}},d_\mathfrak{g})$ is the data of a dg $\mathcal{L}ie^*$-algebra.
\end{Definition}

There is a morphism of operads $\mathcal{S}kew: \mathcal{L}ie \longrightarrow \mathcal{A}ss$ given by the skew-symmetrization of the associative product. By taking the linear dual $\mathcal{S}kew^*: \mathcal{A}ss^* \longrightarrow \mathcal{L}ie^*$, we obtain a morphism of conilpotent cooperads. 

\begin{Proposition}\label{prop: adjunction universal enveloping absolute algebra}
There is an adjunction

\[
\begin{tikzcd}[column sep=7pc,row sep=3pc]
            \mathsf{dg}~\mathsf{abs}~\mathsf{Lie}\text{-}\mathsf{alg} \arrow[r,"\widehat{\mathfrak{U}}"{name=F}, shift left=1.1ex] 
           &\mathsf{dg}~\mathsf{abs}~\mathsf{assoc}\text{-}\mathsf{alg}~, \arrow[l, shift left=.75ex, "\mathrm{Res}_{\mathcal{S}kew^*}"{name=U}]
            \arrow[phantom, from=F, to=U, , "\dashv" rotate=-90]
\end{tikzcd}
\]

between the category of dg absolute Lie algebras and the category of dg absolute associative algebras. The left adjoint functor $\widehat{\mathfrak{U}}$ is called the universal enveloping absolute algebra functor. 
\end{Proposition}

\begin{proof}
The morphism $\mathcal{S}kew^*: \mathcal{A}ss^* \longrightarrow \mathcal{L}ie^*$ induces a natural transformation 

\[
\widehat{\mathscr{S}}^c(\mathcal{S}kew^*): \widehat{\mathscr{S}}^c(\mathcal{L}ie^*) \longrightarrow \widehat{\mathscr{S}}^c(\mathcal{A}ss^*)
\]

between the corresponding monads encoding dg absolute associative algebras and dg absolute Lie algebras.
\end{proof}

\begin{Remark}
Let $(\mathfrak{g},\gamma_{\mathfrak{g}},d_\mathfrak{g})$ be a dg absolute Lie algebra. The universal enveloping absolute algebra $\widehat{\mathfrak{U}}(\mathfrak{g})$ of $\mathfrak{g}$ is given by 

\[
\begin{tikzcd}[column sep=4pc,row sep=4pc]
\mathrm{Coeq}\Bigg(\displaystyle  \overline{\mathcal{T}}^{~\wedge}\bigg(\widehat{\mathcal{L}ie}(\mathfrak{g})\bigg) \arrow[r,"\overline{\mathcal{T}}^{~\wedge}(\gamma_{\mathfrak{g}})",shift right=1.1ex,swap]  \arrow[r,"\psi_\mathfrak{g}"{name=SD},shift left=1.1ex ]
&\overline{\mathcal{T}}^{~\wedge}(\mathfrak{g})\Bigg)~,
\end{tikzcd}
\]

where $\widehat{\mathcal{L}ie}(\mathfrak{g})$ denotes the free completed Lie algebra on $\mathfrak{g}$ and where $\psi_\mathfrak{g}$ is given by the composition

\[
\begin{tikzcd}[column sep=4.5pc,row sep=4pc]
\psi_\mathfrak{g}: \overline{\mathcal{T}}^{~\wedge}(\widehat{\mathcal{L}ie}(\mathfrak{g})) \arrow[r,"\overline{\mathcal{T}}^{~\wedge}(\mathcal{S}kew^*)"]
&\overline{\mathcal{T}}^{~\wedge}(\overline{\mathcal{T}}^{~\wedge}(\mathfrak{g})) \arrow[r, "\gamma_{\overline{\mathcal{T}}^{~\wedge}(\mathfrak{g})}"]
&\overline{\mathcal{T}}^{~\wedge}(\mathfrak{g})~.
\end{tikzcd}
\]

Hence $\widehat{\mathfrak{U}}(\mathfrak{g})$ is quotient of the completed tensor algebra $\overline{\mathcal{T}}^{~\wedge}(\mathfrak{g})$ on $\mathfrak{g}$ where one not only identifies $x \otimes y - (-1)^{|x|} y \otimes x$ with the Lie bracket $[x,y]$, but this identification has also to be done for all formal power series of brackets of elements of $\mathfrak{g}$.
\end{Remark}

\begin{Proposition}
Let $(\mathfrak{g},\gamma_{\mathfrak{g}},d_\mathfrak{g})$ be a nilpotent dg absolute Lie algebra. There is an isomorphism 

\[
\widehat{\mathfrak{U}}(\mathfrak{g}) \cong \frac{\overline{\mathcal{T}}^{~\wedge}(\mathfrak{g})}{\left(x \otimes y - (-1)^{|x|} y \otimes x - [x,y] \right)}
\]
\vspace{0.2pc}

of dg absolute associative algebras.
\end{Proposition}

\begin{proof}
The following square of adjunctions 

\[
\begin{tikzcd}[column sep=5pc,row sep=5pc]
\mathsf{dg}~\mathsf{assoc}\text{-}\mathsf{alg} \arrow[r,"\mathrm{Abs}"{name=B},shift left=1.1ex] \arrow[d,"\mathrm{Skew}"{name=SD},shift left=1.1ex ]
&\mathsf{dg}~\mathsf{abs}~\mathsf{assoc}\text{-}\mathsf{alg} \arrow[d,"\mathrm{Skew}"{name=LDC},shift left=1.1ex ] \arrow[l,"\mathrm{Res}"{name=C},,shift left=1.1ex]  \\
\mathsf{dg}~\mathsf{Lie}\text{-}\mathsf{alg} \arrow[r,"\mathrm{Abs}"{name=CC},shift left=1.1ex]  \arrow[u,"\mathfrak{U}"{name=LD},shift left=1.1ex ]
&\mathsf{dg}~\mathsf{abs}~\mathsf{Lie}\text{-}\mathsf{alg} \arrow[l,"\mathrm{Res}"{name=CB},shift left=1.1ex] \arrow[u,"\widehat{\mathfrak{U}}"{name=TD},shift left=1.1ex] \arrow[phantom, from=SD, to=LD, , "\dashv" rotate=0] \arrow[phantom, from=C, to=B, , "\dashv" rotate=-90]\arrow[phantom, from=TD, to=LDC, , "\dashv" rotate=0] \arrow[phantom, from=CC, to=CB, , "\dashv" rotate=-90]
\end{tikzcd}
\] 

commutes, where $\mathfrak{U}$ is the classical universal enveloping algebra functor. Indeed, one easily checks that 

\[
\mathrm{Res} \cdot \mathrm{Skew} \cong \mathrm{Skew} \cdot \mathrm{Res}~.
\]

This implies that $\widehat{\mathfrak{U}} \cdot \mathrm{Ab} \cong \mathrm{Ab} \cdot \mathfrak{U}$. Therefore, by Remark \ref{Remark: nilpotent Yoneda}, we have that 

\[
\widehat{\mathfrak{U}}(\mathfrak{g}) \cong \mathrm{Ab}\left(\frac{\overline{\mathcal{T}}(\mathfrak{g})}{\left(x \otimes y - (-1)^{|x|} y \otimes x - [x,y] \right)} \right) \cong \frac{\overline{\mathcal{T}}^{~\wedge}(\mathfrak{g})}{\left(x \otimes y - (-1)^{|x|} y \otimes x - [x,y] \right)}~.
\]
\end{proof}

\section{Applications to absolute Lie theory}
In this section, we give a generalization of one the main results obtained by R. Campos, D. Petersen, D. Robert-Nicoud and F. Wierstra in \cite{campos2020lie}, by rephrasing and reinterpreting their constructions with the language of algebras over cooperads, and by using the formalism developed so far. We now consider the dg operads $\Omega \mathcal{L}ie^*$ and $\Omega \mathcal{A}ss^*$ which encode, respectively, \textit{shifted} $\mathcal{C}_\infty$ and \textit{shifted} $\mathcal{A}_\infty$ algebras (and coalgebras). From now on, the adjective \textit{shifted} will be implicit.

\begin{lemma}
There is a morphism of dg operads
\[
\varphi: \Omega \mathcal{A}ss^* \longrightarrow \Omega \mathcal{L}ie^*~,
\]
given by $\varphi \coloneqq \Omega(\mathcal{S}kew^*)$.
\end{lemma}

\begin{proof}
It is straightforward from the definition.
\end{proof}

\begin{Proposition}
The morphism of dg operads $\varphi: \Omega \mathcal{A}ss^* \longrightarrow \Omega \mathcal{L}ie^*$ induces a Quillen adjunction

\[
\begin{tikzcd}[column sep=7pc,row sep=3pc]
           \mathcal{C}_\infty\text{-}\mathsf{coalg} \arrow[r,"\mathrm{Res}_{\varphi}"{name=F}, shift left=1.1ex] 
           &\mathcal{A}_\infty\text{-}\mathsf{coalg}~. \arrow[l, shift left=.75ex, "\mathrm{Coind}_{\varphi}"{name=U}]
            \arrow[phantom, from=F, to=U, , "\dashv" rotate=-90]
\end{tikzcd}
\]

between the model category of $\mathcal{C}_\infty$-coalgebras and the model category of $\mathcal{A}_\infty$-coalgebras.
\end{Proposition}

\begin{proof}
Since both dg operads $\Omega \mathcal{A}ss^*$ and $\Omega \mathcal{L}ie^*$ are cofibrant, their respective categories of dg coalgebras admit a left-transferred model structure where weak equivalences are given by quasi-isomorphisms and where cofibrations are given by degree-wise monomorphisms, by Theorem \ref{thm: model structure on P-cog}. Both categories are comonadic by Theorem \ref{thm: existence of the cofree P cog}, and the morphism $\varphi$ induces a morphism between their respective comonads. Thus it induces an adjunction. Finally, it is straightforward to check that $\mathrm{Res}_{\varphi}$ preserves cofibrations and quasi-isomorphisms, since it does not change the underlying chain complex of the $\mathcal{C}_\infty$-coalgebra.
\end{proof}

\begin{theorem}[{\cite[Theorem 4.27]{campos2020lie}}]\label{thm: cited thm about C infinity}
Let $C_1$ and $C_2$ be two $\mathcal{C}_\infty$-coalgebras. There exists a zig-zag of quasi-isomorphisms of $\mathcal{C}_\infty$-coalgebras 
\[
C_1 \lqi \cdot \qi C_2
\]
if and only if there exists a zig-zag of $\mathcal{A}_\infty$-coalgebras 
\[
\mathrm{Res}_{\varphi}(C_1) \lqi \cdot \qi \mathrm{Res}_{\varphi}(C_2)~.
\]
\end{theorem}

We now use the compatibility of the complete Bar-Cobar adjunctions of Proposition \ref{prop: compatibility of induced adjunctions} in order to state a stronger version of \cite[Theorem B]{campos2020lie}.

\begin{Proposition}\label{Cor: commutative of the square of absolute Lie}
There is a commuting square

\[
\begin{tikzcd}[column sep=5pc,row sep=5pc]
\mathcal{A}_\infty\text{-}\mathsf{coalg} \arrow[r,"\widehat{\Omega}_\iota"{name=B},shift left=1.1ex] \arrow[d,"\mathrm{Coind}_{\varphi} "{name=SD},shift left=1.1ex ]
&\mathsf{dg}~\mathsf{abs}~\mathsf{assoc}\text{-}\mathsf{alg}^{\mathsf{comp}} \arrow[d,"\mathrm{Skew}"{name=LDC},shift left=1.1ex ] \arrow[l,"\widehat{\mathrm{B}}_\iota"{name=C},,shift left=1.1ex]  \\
\mathcal{C}_\infty\text{-}\mathsf{coalg} \arrow[r,"\widehat{\Omega}_\iota "{name=CC},shift left=1.1ex]  \arrow[u,"\mathrm{Res}_{\varphi}"{name=LD},shift left=1.1ex ]
&\mathsf{dg}~\mathsf{abs}~\mathsf{Lie}\text{-}\mathsf{alg}^{\mathsf{comp}} \arrow[l,"\widehat{\mathrm{B}}_\iota"{name=CB},shift left=1.1ex] \arrow[u,"\widehat{\mathfrak{U}}"{name=TD},shift left=1.1ex] \arrow[phantom, from=SD, to=LD, , "\dashv" rotate=0] \arrow[phantom, from=C, to=B, , "\dashv" rotate=-90]\arrow[phantom, from=TD, to=LDC, , "\dashv" rotate=0] \arrow[phantom, from=CC, to=CB, , "\dashv" rotate=-90]
\end{tikzcd}
\] 

of Quillen adjunctions.
\end{Proposition}

\begin{proof}
It is immediate to check that the square 
\[
\begin{tikzcd}[column sep=3pc,row sep=3pc]
\mathcal{A}ss^* \arrow[r,"\iota "] \arrow[d,"\mathcal{S}kew^*",swap] 
&\Omega \mathcal{A}ss^* \arrow[d,"\varphi"]\\
\mathcal{L}ie^* \arrow[r,"\iota"]
&\Omega \mathcal{L}ie^*~.
\end{tikzcd}
\]

is commutative. The result follows from Proposition \ref{Cor: commutative of the square of absolute Lie}.
\end{proof}

\begin{theorem}\label{thm: inclusion infini cat des Lie absolues}
Let $\mathfrak{g}$ and $\mathfrak{h}$ be two dg absolute Lie algebras. There exists a zig-zag of weak equivalences of dg absolute Lie algebras
\[
\mathfrak{g} \lqi \cdot \qi \mathfrak{h}
\]
if and only if there exists a zig-zag of weak equivalences of dg absolute associative algebras
\[
\widehat{\mathfrak{U}}(\mathfrak{g}) \lqi \cdot \qi \widehat{\mathfrak{U}}(\mathfrak{h})  ~.
\]
\end{theorem}

\begin{proof}
It is a direct consequence of Theorem \ref{thm: cited thm about C infinity}, using the fact that the horizontal Quillen adjunctions of Proposition \ref{Cor: commutative of the square of absolute Lie} are Quillen equivalences.
\end{proof}

Let us make more explicit what these weak equivalences of dg absolute associative algebras or dg absolute Lie algebras look like. We state the analogue of \cite[Proposition 2.5]{Vallette14} for the model structure on algebras over a conilpotent dg cooperad transferred along the complete Bar-Cobar adjunction.

\begin{Proposition}\label{thm: weak equiv inclues dans les quasi-isos}
Let $f: \mathfrak{g} \qi \mathfrak{h}$ be a weak-equivalence of complete dg absolute Lie algebras. It is in particular a quasi-isomorphism.
\end{Proposition}

\begin{proof}
First notice that the dg operad $\Omega \mathcal{L}ie^*$ is augmented, that is, there is a morphism of dg operads $\varepsilon: \Omega \mathcal{L}ie^* \longrightarrow \I$. Now we consider the following commutative square 

\[
\begin{tikzcd}[column sep=3pc,row sep=3pc]
\mathcal{L}ie^* \arrow[r,"\iota "] \arrow[d,"\mathrm{id}_{\mathcal{L}ie^*}",swap] 
&\Omega \mathcal{L}ie^* \arrow[d,"\varepsilon"]\\
\mathcal{L}ie^* \arrow[r,"\epsilon_{*}\iota"]
&\I~,
\end{tikzcd}
\]
where $\varepsilon_{*}\iota$ is the push forward of the twisting morphism $\iota$ by $\varepsilon$. It induces a commutative square of Quillen adjunctions 

\[
\begin{tikzcd}[column sep=5pc,row sep=5pc]
\mathcal{C}_\infty\text{-}\mathsf{coalg} \arrow[r,"\widehat{\Omega}_\iota"{name=B},shift left=1.1ex] \arrow[d,"\mathrm{Coind}_{\varepsilon} "{name=SD},shift left=1.1ex ]
&\mathsf{dg}~\mathsf{abs}~\mathsf{Lie}\text{-}\mathsf{alg}^{\mathsf{comp}} \arrow[d,"\mathrm{Id}"{name=LDC},shift left=1.1ex ] \arrow[l,"\widehat{\mathrm{B}}_\iota"{name=C},,shift left=1.1ex]  \\
\mathsf{dg}\text{-}\mathsf{mod} \arrow[r,"\widehat{\Omega}_{\varepsilon_{*}\iota}"{name=CC},shift left=1.1ex]  \arrow[u,"\mathrm{Res}_{\varepsilon}"{name=LD},shift left=1.1ex ]
&\mathsf{dg}~\mathsf{abs}~\mathsf{Lie}\text{-}\mathsf{alg}^{\mathsf{comp}}~. \arrow[l,"\widehat{\mathrm{B}}_{\varepsilon_{*}\iota}"{name=CB},shift left=1.1ex] \arrow[u,"\mathrm{Id}"{name=TD},shift left=1.1ex] \arrow[phantom, from=SD, to=LD, , "\dashv" rotate=0] \arrow[phantom, from=C, to=B, , "\dashv" rotate=-90]\arrow[phantom, from=TD, to=LDC, , "\dashv" rotate=0] \arrow[phantom, from=CC, to=CB, , "\dashv" rotate=-90]
\end{tikzcd}
\] 
 
One can easily see that the bottom adjunction is the free-forgetful adjunction. The category of dg modules is equivalent to the category of coalgebras over the trivial operad $I$, thus we can transfer its model category structure along the free-forgetful adjunction using the results of \cite[Section 10]{grignoulejay18} to obtain a Quillen adjunction, where the model category structure on dg absolute Lie algebras has quasi-isomorphisms as weak-equivalences. Since the identity functor 

\[
\mathrm{Id}: \left(\mathsf{dg}~\mathsf{abs}~\mathsf{Lie}\text{-}\mathsf{alg}^{\mathsf{comp}}, \mathcal{W} \right) \longrightarrow \left(\mathsf{dg}~\mathsf{abs}~\mathsf{Lie}\text{-}\mathsf{alg}^{\mathsf{comp}}, \mathcal{Q}uasi\text{-}isos \right) 
\]
\vspace{0.2pc}

is a right Quillen functor and since every object is fibrant, it sends weak-equivalences of dg absolute Lie algebras to quasi-isomorphisms. Therefore any weak-equivalence of dg absolute Lie algebras is a quasi-isomorphisms.
\end{proof}

\begin{Remark}
The above proposition can be generalized \textit{mutatis mutandis} to the case of any twisting morphism $\alpha: \mathcal{P} \longrightarrow \mathcal{C}$ between a reduced cofibrant dg operad $\mathcal{P}$ and a conilpotent dg cooperad $\mathcal{C}$. In particular, this also holds for weak-equivalences of dg absolute associative algebras.
\end{Remark}

\begin{Remark}
Model category structures on algebras over a conilpotent dg cooperads behave in an analogous way to what happens with model category structures on coalgebras over conilpotent dg cooperads, as described in \cite{DrummondColeHirsh14}.
\end{Remark}

\begin{Corollary}
Let $\mathfrak{g}$ be a complete dg absolute Lie algebra. Then the canonical morphism 
\[
\epsilon_{\mathfrak{g}}: \widehat{\Omega}_\iota \widehat{\mathrm{B}}_\iota \mathfrak{g} \qi \mathfrak{g}
\]
is in fact a quasi-isomorphism of dg absolute Lie algebras. 
\end{Corollary}

\begin{proof}
This is automatic from Proposition \ref{thm: weak equiv inclues dans les quasi-isos}. 
\end{proof}

\begin{theorem}\label{thm: iso envelopantes Lie absolues}
Let $\mathfrak{g}$ and $\mathfrak{h}$ be two complete graded absolute Lie algebras. They are isomorphic as complete graded absolute Lie algebras if and only if their universal enveloping absolute algebras are isomorphic.
\end{theorem}

\begin{proof}
This follows from Theorems \ref{thm: inclusion infini cat des Lie absolues} and Proposition \ref{thm: weak equiv inclues dans les quasi-isos}, considering the fact that any weak-equivalence is in particular a quasi-isomorphism, and that any quasi-isomorphism between graded modules with zero differential is an isomorphism.
\end{proof}

\begin{Remark}
Any nilpotent graded Lie algebra is a complete graded absolute Lie algebra by the analogue of Proposition \ref{prop: les nilpotentes sont des absolues} for absolute Lie algebras. Therefore this result is a generalization of \cite[Corollary 0.12]{campos2020lie}.
\end{Remark}

This approach also allows us to generalize the above theorems to the universal enveloping absolute $\mathcal{A}_\infty$-algebra of an absolute $\mathcal{L}_\infty$-algebra. We change our \textit{shifting conventions}, considering this time unshifted $\mathcal{C}_\infty$ or $\mathcal{A}_\infty$-coalgebras, and therefore shifting absolute $\mathcal{L}_\infty$-algebras and absolute $\mathcal{A}_\infty$-algebras.

\begin{theorem}\label{thm: inclusion infini cat des L infinies}
Let $\mathfrak{g}$ and $\mathfrak{h}$ be two absolute $\mathcal{L}_\infty$-algebras. There exists a zig-zag of weak equivalences of absolute $\mathcal{L}_\infty$-algebras
\[
\mathfrak{g} \lqi \cdot \qi \mathfrak{h}
\]
if and only if there exists a zig-zag of weak equivalences of absolute $\mathcal{A}_\infty$-algebras
\[
\widehat{\mathfrak{U}}_\infty(\mathfrak{g}) \lqi \cdot \qi \widehat{\mathfrak{U}}_\infty(\mathfrak{h})  ~.
\]
\end{theorem}

\begin{proof}
We consider the dg operad $\Omega \mathrm{B} \mathcal{C}om$. Since it is a cofibrant resolution for the operad $\mathcal{C}om$, there exists a quasi-isomorphism of dg operads $f: \Omega \mathrm{B} \mathcal{C}om \qi \Omega s\mathcal{L}ie^*$. Therefore there is a Quillen equivalence between dg $\Omega \mathrm{B} \mathcal{C}om$-coalgebras and $\mathcal{C}_\infty$-coalgebras. Likewise, there is a Quillen equivalence between dg $\Omega \mathrm{B} \mathcal{A}ss$-coalgebras and $\mathcal{A}_\infty$-coalgebras. We consider the adjunction 

\[
\begin{tikzcd}[column sep=7pc,row sep=3pc]
           \mathsf{dg}~\Omega \mathrm{B} \mathcal{C}om\text{-}\mathsf{coalg} \arrow[r,"\mathrm{Res}_{\rho}"{name=F}, shift left=1.1ex] 
           &\mathsf{dg}~\Omega \mathrm{B} \mathcal{A}ss\text{-}\mathsf{coalg}~. \arrow[l, shift left=.75ex, "\mathrm{Coind}_{\rho}"{name=U}]
            \arrow[phantom, from=F, to=U, , "\dashv" rotate=-90]
\end{tikzcd}
\]

induced by the morphism of dg operads $\rho: \Omega \mathrm{B} \mathcal{A}ss \longrightarrow \Omega \mathrm{B} \mathcal{C}om$. Two dg $\Omega \mathrm{B} \mathcal{C}om$-coalgebras $C_1$ and $C_2$ are linked by a zig-zag of quasi-isomorphisms if and only if $\mathrm{Res}_{\rho}(C_1)$ and $\mathrm{Res}_{\rho}(C_2)$ are linked by a zig-zag of quasi-isomorphisms of dg $\Omega \mathrm{B} \mathcal{A}ss$-coalgebras. Using the commutative square 

\[
\begin{tikzcd}[column sep=5pc,row sep=5pc]
\mathsf{dg}~\Omega \mathrm{B} \mathcal{A}ss\text{-}\mathsf{coalg} \arrow[r,"\widehat{\Omega}_\iota"{name=B},shift left=1.1ex] \arrow[d,"\mathrm{Coind}_{\rho} "{name=SD},shift left=1.1ex ]
&\mathsf{abs}~\mathcal{A}_\infty\text{-}\mathsf{alg}^{\mathsf{comp}} \arrow[d,"\mathrm{Skew}"{name=LDC},shift left=1.1ex ] \arrow[l,"\widehat{\mathrm{B}}_\iota"{name=C},,shift left=1.1ex]  \\
\mathsf{dg}~\Omega \mathrm{B} \mathcal{C}om\text{-}\mathsf{coalg} \arrow[r,"\widehat{\Omega}_\iota "{name=CC},shift left=1.1ex]  \arrow[u,"\mathrm{Res}_{\rho}"{name=LD},shift left=1.1ex ]
&\mathsf{abs}~\mathcal{L}_\infty\text{-}\mathsf{alg}^{\mathsf{comp}} \arrow[l,"\widehat{\mathrm{B}}_\iota"{name=CB},shift left=1.1ex] \arrow[u,"\widehat{\mathfrak{U}}"{name=TD},shift left=1.1ex] \arrow[phantom, from=SD, to=LD, , "\dashv" rotate=0] \arrow[phantom, from=C, to=B, , "\dashv" rotate=-90]\arrow[phantom, from=TD, to=LDC, , "\dashv" rotate=0] \arrow[phantom, from=CC, to=CB, , "\dashv" rotate=-90]
\end{tikzcd}
\] 

and the fact that the horizontal adjunction are Quillen equivalences concludes the proof.
\end{proof}

\begin{Proposition}\label{prop: weak equiv de L infinie absolues sont des quasi-isos}
Let $f: \mathfrak{g} \qi \mathfrak{h}$ be a weak equivalence of complete absolute $\mathcal{L}_\infty$-algebras. It is in particular a quasi-isomorphism.
\end{Proposition}

\begin{proof}
The arguments are the same as in the proof of Theorem \ref{thm: weak equiv inclues dans les quasi-isos} apply in this situation.
\end{proof}

\begin{Definition}[Minimal absolute $\mathcal{L}_\infty$-algebra]
Let $(\mathfrak{g},\gamma_\mathfrak{g},d_\mathfrak{g})$ be an absolute $\mathcal{L}_\infty$-algebra. It is \textit{minimal} if the differential $d_\mathfrak{g}$ is equal to zero.
\end{Definition}

\begin{theorem}\label{thm: isos envelopantes absolues de L infinies}
Let $\mathfrak{g}$ and $\mathfrak{h}$ be two complete minimal absolute $\mathcal{L}_\infty$-algebras. They are isomorphic as complete minimal absolute $\mathcal{L}_\infty$-algebras if and only if their universal enveloping absolute $\mathcal{A}_\infty$-algebras are isomorphic.
\end{theorem}

\begin{proof}
This is a direct corollary of Proposition \ref{prop: weak equiv de L infinie absolues sont des quasi-isos} and Theorem \ref{thm: inclusion infini cat des L infinies}, using the same arguments as in the proof of Theorem \ref{thm: iso envelopantes Lie absolues}.
\end{proof}

\begin{Example}[Arity-wise nilpotent $\mathcal{L}_\infty$-algebras]
Nilpotent $\mathcal{L}_\infty$-algebras in the sense of \cite{Getzler09} are particular examples of absolute $\mathcal{L}_\infty$-algebras. Therefore the above theorem applies to \textit{minimal} nilpotent $\mathcal{L}_\infty$-algebras without any degree restriction.
\end{Example}

\begin{Remark}
The analogues of Theorems \ref{thm: inclusion infini cat des Lie absolues} and \ref{thm: inclusion infini cat des L infinies} should also hold when we replace the categories of absolute Lie/$\mathcal{L}_\infty$-algebras by their \textit{curved} counterparts. Indeed, using Remark \ref{Rmk: PBW curved}, one should get analogue statements as in \cite{campos2020lie} concerning the deformation complexes of unital $\mathcal{C}_\infty$-coalgebras and unital $\mathcal{A}_\infty$-coalgebras. Then applying the same formalism is straightforward.
\end{Remark}

\section*{Appendix: Induction and coinduction}\label{appendice: induction et coinduction}
Let $f: H \longrightarrow G$ be a morphism of groups. It admits an unique extension into a morphism of $\kk$-algebras $f: \kk[H] \longrightarrow \kk[G]$. Let $V$ be a left $\kk[H]$-module, the \textit{induced representation} $\mathrm{Ind}_f(V)$ is the left $\kk[G]$-module given by:

\[
\mathsf{Ind}_f(V) \coloneqq \kk[G] \otimes_{H} V~,
\]

where $\kk[G]$ is considered as a right $\kk[H]$-module via the morphism $f ~.$ The \textit{coinduced representation} $\mathrm{Coind}_f(V)$ is the left $\kk[G]$-module given by:

\[
\mathsf{Coind}_f(V) \coloneqq \mathrm{Hom}_{H}(\kk[G], V)~,
\]

where $\kk[G]$ is considered as a left $\kk[H]$-module via the morphism $f ~.$ Let $L$ be a left $\kk[G]$-module, one can endow it with a left $\kk[H]$-module structure by restricting the action of $\kk[G]$ along the morphism $f ~.$ These constructions are functorial and assemble into the following adjunction diagram:

\[
\begin{tikzcd}[column sep=8pc,row sep=2pc]
\kk[H]\text{-}\mathsf{mod} \arrow[r,"\mathsf{Res}_f",""{name=B, below}] & \kk[G]\text{-}\mathsf{mod} \arrow[l,bend right=60,"\mathsf{Ind}_f",swap,""{name=A,above}] \arrow[l,bend left=60,"\mathsf{Coind}_f",""{name=C,above}] \arrow[phantom, from=B, to=A, "\dashv" rotate=-90] \arrow[phantom, from=C, to=B, "\dashv" rotate=-90]
\end{tikzcd}~.
\]

When the morphism $f: H \longrightarrow G$ is the inclusion morphism of a subgroup $H$, the induced representation is denoted $\mathsf{Ind}_H^G(V)$ and the coinduced representation $\mathsf{Coind}_H^G(V)~.$

\begin{lemma}\label{lemmaindandcoind}
Let $H$ be a subgroup of $G$. Let $V$ be a left $\kk[H]$-module. Then there is a natural isomorphism of right $\kk[G]$-modules:

\[
\mathrm{Hom}(\mathsf{Ind}_H^G(V),\kk) \cong \mathsf{Coind}_H^G(V^*)~.
\]
\end{lemma} 

\begin{proof}
The adjunction tensor-hom gives the isomorphism:
\[
\mathrm{Hom}(\kk[G]\otimes_H V, \kk) \cong \mathrm{Hom}(\kk[G],\mathrm{Hom}_H(V,\kk))~.
\]
The set of $H$-equivariant morphism $\mathrm{Hom}_H(V,\kk)$ is given by the invariants $\mathrm{Hom}(V,\kk)^H$, hence:
\[
\mathrm{Hom}(\kk[G],\mathrm{Hom}(V,\kk)^H) \cong \mathrm{Hom}(\kk[G],\mathrm{Hom}(V,\kk))^H \cong \mathrm{Hom}_H(\kk[G],V^*)~,
\]
since invariants are given by a limit diagram.
\end{proof}

\begin{lemma}\label{lemma invariants and coinduction}
Let $H$ be a subgroup of $G$. Let $V$ be a left $\kk[H]$-module. There is a natural isomorphism of $\kk$-modules: 

\[
\left(\mathsf{Coind}_H^G(V)\right)^G \cong V^H ~.
\]
\end{lemma}

\begin{proof}
We have that:
\[
\mathrm{Hom}_H(\kk[G],V)^G \cong \mathrm{Hom}(\kk[G],V)^{G\times H} \cong \mathrm{Hom}_G(\kk[G],V^H) \cong V^H~.
\]
\end{proof}

\chapter{Curved operadic calculus.}\label{Chapter 2}

\epigraph{"Je le vis, je rougis, je pâlis à sa vue. \\ Un trouble s'éleva dans mon âme éperdue. \\ Mes yeux ne voyaient plus, je ne pouvais parler. \\ Je sentis tous mon corps et transir et brûler. \\ Je reconnus ce formalisme admirable, \\ où toute algèbre et cogèbre a son semblable."}{Phèdre, dans la pièce homonyme de Jean Racine.}

\section*{Introduction}
\textbf{Global picture.} The theory of operads shifts the point of view of universal algebra: instead of working by hand with specific types of algebras, one works with the operads that encode them. Examples of algebraic structures encoded by operads include associative algebras, Lie algebras, Poisson algebras, Batalin--Vilkovisky algebras, etc. One calls \textit{operadic calculus} the set of techniques that allow us to work with operads themselves. There is a purely algebraic side to the operadic calculus. For instance, one might consider morphisms between operads. A morphism between operads induces a structure of an operadic bimodule. Via the theory of operadic bimodules, one can construct universal functors between the categories of algebras over operads. Thus this theory gives, for any morphism of operads, a universal adjunction between their respective categories of algebras. The universal enveloping algebra of a Lie algebra can be recovered in this way. This approach allows vast generalizations, and new universal enveloping algebras can be constructed for other types of algebraic structures in this way. It can also generalize well-known theorems such as the Poincaré--Birkoff--Witt theorem in a functorial way, see \cite{Tamaroff2020}. 

\medskip

The operadic calculus has also an homotopical side. When operads themselves live in the category of differential graded $\mathbb{S}$-modules, they admit a notion of weak equivalence given by arity-wise quasi-isomorphisms. Understanding the homotopy theory of operads gives results on the homotopy theory in the category of algebras it encodes. Two weakly equivalent operads encode the same homotopy category of algebras. For any operad $\mathcal{P}$, algebras over a cofibrant resolution of $\mathcal{P}$ provide us with a suitable notion of a $\mathcal{P}$-\textit{algebra up to homotopy}. This recovers the seminal notions of $\mathcal{A}_\infty$-algebras, $\mathcal{L}_\infty$-algebras, $\mathcal{C}_\infty$-algebras as particular examples. These resolutions also allow one to construct universal André-Quillen cohomology theories for algebras over an operad, see \cite{JoanQuillen}. The main tools for studying the homotopy theory of operads and computing cofibrant resolutions are the Bar-Cobar adjunction and the Koszul duality theory. Combining both the algebraic and the homotopical aspects of the operadic calculus provides us with new tools to solve problems, see \cite{LodayVallette12}. For example, in \cite{campos2020lie}, R. Campos, D. Petersen, D. Robert-Nicoud, and F. Wierstra show that a nilpotent Lie algebra is completely characterized up to isomorphism by its universal enveloping algebra. The proof of this purely algebraic result requires a combination of both of the aforementioned methods. Another example is given by \cite{robertnicoud2020higher}, where D. Robert-Nicoud and B. Vallette develop the integration theory of $\mathcal{L}_\infty$-algebras using methods coming from the operadic calculus. One of the main motivations of this chapter was to lay down the operadic tools required to generalize the results of \textit{loc.cit} to the case of \textit{curved} $\mathcal{L}_\infty$-algebras, which will be done in Chapter \ref{Chapter 3}.

\medskip

In the mist of algebraic structures, there are the so-called \textit{curved algebras}. The prototype of curved algebras are curved associative algebras. These are graded associative algebras $(A, \mu_A)~,$ endowed with a derivation $d$ of degree $-1$ and a distinguished element $\theta$ of degree $-2$ called the \textit{curvature}, such that 
\[
d^2(-) = \mu_A(\theta,-) - \mu_A(-,\theta)~. \quad \quad (*)
\]
Here, the element $\theta$ is the obstruction for the derivation $d$ to square to zero. When $\theta$ is non-trivial, this type of algebras do not have underlying homology groups, thus there is no notion of quasi-isomorphism for them in general. In these examples, the distinction between homological algebra and the more general concept of homotopical algebra becomes apparent. Other examples include curved Lie algebras, curved $\mathcal{A}_\infty$-algebras, curved $\mathcal{L}_\infty$-algebras and so on. These types of algebras are playing an increasingly important role in various areas of mathematics: curved $\mathcal{A}_\infty$-algebras in Floer cohomology in symplectic geometry \cite{FOOO7}, curved $\mathcal{L}_\infty$-algebras in derived differential geometry \cite{behrend2021derived}, derived deformation theory \cite{calaque2021lie} and $\mathcal{L}_\infty$-spaces \cite{costello2011geometric}. 

\medskip

\textbf{Main results.} The goal of this foundation chapter, based on \cite{lucio2022curved}, is to settle the operadic calculus to \textit{curved operads}. Curved operads appear naturally when one tries to encode types of curved algebras with operad-like structures: to encode the curvature relation $(*)$ at the algebra level, one needs to add a curvature on the operad level. Working with curved operads forces the underlying category to be the category of pre-differential graded (pdg) modules, which are given by graded modules with a degree $-1$ endomorphism. So we start by settling the basic properties related to this new framework. In order to develop the algebraic side of the curved operad calculus, we introduce the notion of a curved operadic bimodule and we develop the subsequent theory. There are obstructions to this generalization: the "free algebra" over a curved operad does not, in general, satisfy the curvature relation. This implies that a curved operad is not a "naive" curved left module over itself. Nevertheless, we develop the theory of curved bimodules for curved operads. This allows us to construct universal functors between the categories of curved algebras. For instance, we construct for the first time the universal curved enveloping algebra for curved Lie algebras and the universal curved enveloping $\mathcal{A}_\infty$-algebra for curved $\mathcal{L}_\infty$-algebras. Along the way, we explain again the theory of \cite{grignoulejay18} in Section \ref{Section: Complete Bar-Cobar} in the curved setting and generalize the duality square of adjunctions constructed in Section \ref{Section: Magic squares} to the curved setting.

\medskip

Once this algebraic framework is established, the rest of the chapter is devoted to generalizing the homotopical side of the operadic calculus to the curved setting. We start \textit{ab initio} from a conceptual point of view. We use the groupoid-colored formalism developed by B. Ward in \cite{ward19}. We introduce the unital groupoid-colored operad $u\mathcal{O}$ that encodes unital partial operads as its algebras and counital partial cooperads as its coalgebras. By partial operad, we mean the definition of operads introduced by M. Markl in \cite{Markl96} in terms of partial composition maps $\{\circ_i\}$. The notion of partial cooperads is the dual notion, defined in terms of partial decomposition maps $\{\Delta_i\}$. We generalize the main point of the inhomogeneous Koszul duality of Hirsh--Millès in \cite{HirshMilles12} to extend them to the groupoid-colored framework. This allows us to compute the conilpotent curved groupoid-colored dual cooperad $u\mathcal{O}^{\ac} \cong c\mathcal{O}^\vee$ (up to suspension). We show that coalgebras over $c\mathcal{O}^\vee$ correspond to conilpotent curved partial cooperads. 

\begin{theoremintro}[Theorem \ref{Koszulity of uO}]
The unital groupoid-colored operad $u\mathcal{O}$ is a Koszul operad and its Koszul dual curved groupoid-colored cooperad is given by $c\mathcal{O}^\vee$ up to suspension.
\end{theoremintro}

In particular, this result implies that there is a Koszul curved twisting morphism 
\[
\kappa: c\mathcal{O}^\vee \longrightarrow u\mathcal{O}~.
\]
From this curved twisting morphism, one obtains a first Bar-Cobar adjunction using the classical methods of \cite{LodayVallette12} and \cite{grignou2019}. This adjunction recovers the Bar-Cobar adjunction between unital partial operads and conilpotent curved partial cooperads constructed by B. Le Grignou in \cite{grignou2021}. In \textit{loc.cit.}, the author endows the category of unital operads with a model structure where weak-equivalences are given by arity-wise quasi-isomorphisms, and then transfers it along this adjunction to conilpotent curved coaugmented cooperads. This endows conilpotent curved coaugmented cooperads with a meaningful model structure, and the Bar-Cobar adjunction is shown to be a Quillen equivalence. The above theorem provides a new proof that the Bar-Cobar adjunction recovered is indeed a Quillen equivalence.

\medskip

Our goal here is to construct another adjunction using the curved twisting morphism $\kappa$. In order to do this, we generalize the results of D. Lejay and B. Le Grignou in \cite{grignoulejay18} to the groupoid-colored case. In \textit{loc.cit}, from a curved twisting morphism the authors construct a "complete" Bar-Cobar adjunction between the category of coalgebras over an operad and the category of curved algebras over a curved cooperad. In our case, we get an adjunction between the category of coalgebras over the groupoid-colored operad $u\mathcal{O}$ and the category of curved algebras over the curved groupoid-colored cooperad $c\mathcal{O}^\vee$. We call this adjunction the \textit{complete Bar-Cobar adjunction}. Coalgebras over $u\mathcal{O}$ are simply counital partial cooperads. Algebras over cooperads are a new notion that provide us with algebraic objects where infinite sums of compositions of operations are well-defined by definition. There is still a lot of work to be done in this area in order to understand these new algebraic structures that emerge naturally from this operadic calculus. In this particular case, curved algebras over $c\mathcal{O}^\vee$ give rise to new objects in the operadic calculus, which we call \textit{curved absolute partial operads}. There is also a more basic notion of \textit{absolute partial operad}. Both of these structures appear naturally when operations of arity $0$ and $1$ are taken into account. We provide a characterization of them in Appendix \ref{Appendix B}. Summarizing, we obtain a complete Bar-Cobar adjunction

\[
\begin{tikzcd}[column sep=7pc,row sep=3pc]
            \mathsf{dg}~\mathsf{uCoop} \arrow[r, shift left=1.1ex, "\widehat{\Omega}"{name=F}] &\mathsf{curv}~\mathsf{abs}~\mathsf{Op}~. \arrow[l, shift left=.75ex, "\widehat{\mathrm{B}}"{name=U}]
            \arrow[phantom, from=F, to=U, , "\dashv" rotate=-90]
\end{tikzcd}
\]

Notice first that this adjunction interrelates "counits" with "curvature", which are known to be Koszul dual. See \cite{PolischukPositselski05}. But its main novelty lies in the fact that it also lifts the usual assumption of \textit{conilpotency} on the cooperad side. Lifting conilpotency on one side of the Koszul duality is what makes infinite sums appear in the other side. This is the conceptual explanation for the existence of these curved \textit{absolute} partial operads. Absolute partial operads admit a canonical filtration, dual to the coradical filtration for partial cooperads. And like the coradical filtration, this filtration comes from the structure of the object. By definition, infinite sums of partial compositions are well-defined. But the topology induced by this filtration might not be Hausdorff. We restrict ourselves to the sub-category of those which are complete. Notice that unlike many other approaches to curved objects, where one \textit{changes the base category} from graded modules to filtered/complete graded modules in order to deal with the infinite sums that appear, here we introduce new types of algebraic structures which admit these infinite sums \textit{without enriching further the underlying category of graded modules}. And the complete filtrations that appear, come from the structure of these new objects, and are therefore \textit{canonical}.

\medskip

Our next step is to endow the category of complete curved absolute partial operads with a satisfactory model structure. Since the notion of a quasi-isomorphism does not exist in the curved context, our approach is to obtain a notion of weak equivalences via a transfer theorem from a Koszul dual category. For this purpose, the category of counital partial cooperads is too narrow. We define the notion of counital partial cooperads up to homotopy. Then, we endow the category of counital partial cooperad up to homotopy with strict morphisms with a model category structure where weak equivalences are given by arity-wise quasi-isomorphisms. Finally, we construct another complete Bar-Cobar adjunction and transfer this model structure to the category of complete curved absolute partial operads. 

\begin{theoremintro}[Theorem \ref{thm: model structure on curved operads}]
The category of complete curved absolute partial operads admits a model structure transferred along the adjunction
\[
\begin{tikzcd}[column sep=7pc,row sep=3pc]
            \mathsf{uCoop}_\infty \arrow[r, shift left=1.1ex, "\widehat{\Omega}_\iota"{name=F}] &\mathsf{curv}~\mathsf{abs}~\mathsf{Op}^{\mathsf{comp}}~, \arrow[l, shift left=.75ex, "\widehat{\mathrm{B}}_\iota"{name=U}]
            \arrow[phantom, from=F, to=U, , "\dashv" rotate=-90]
\end{tikzcd}
\]
where the model category structure considered on the left-hand side has arity-wise quasi-isomorphisms as weak-equivalences and monomorphisms as cofibrations.
\end{theoremintro}

We then relate the "classical" Bar-Cobar adjunction with the complete Bar-Cobar adjunction constructed here. For this purpose, we construct a pair of duality adjunctions that relate these Bar-Cobar adjunctions, analogue to Section \ref{Section: Curved duality square}. First we do this in the algebraic context, where we prove the following result.

\begin{theoremintro}[Theorem \ref{thm: carré magique}]
The following square of adjunctions 
\[
\begin{tikzcd}[column sep=5pc,row sep=5pc]
\mathsf{dg}~\mathsf{uOp}^{\mathsf{op}} \arrow[r,"\mathrm{B}^{\mathsf{op}}"{name=B},shift left=1.1ex] \arrow[d,"(-)^\circ "{name=SD},shift left=1.1ex ]
&\left(\mathsf{curv}~\mathsf{Coop}^{\mathsf{conil}}\right)^{\mathsf{op}}  \arrow[d,"(-)^*"{name=LDC},shift left=1.1ex ] \arrow[l,"\Omega^{\mathsf{op}}"{name=C},,shift left=1.1ex]  \\
\mathsf{dg}~\mathsf{uCoop} \arrow[r,"\widehat{\Omega}"{name=CC},shift left=1.1ex]  \arrow[u,"(-)^*"{name=LD},shift left=1.1ex ]
&\mathsf{curv}~\mathsf{abs}~\mathsf{Op}^{\mathsf{comp}}~, \arrow[l,"\widehat{\mathrm{B}}"{name=CB},shift left=1.1ex] \arrow[u,"(-)^\vee"{name=TD},shift left=1.1ex] \arrow[phantom, from=SD, to=LD, , "\dashv" rotate=0] \arrow[phantom, from=C, to=B, , "\dashv" rotate=-90]\arrow[phantom, from=TD, to=LDC, , "\dashv" rotate=0] \arrow[phantom, from=CC, to=CB, , "\dashv" rotate=-90]
\end{tikzcd}
\] 
commutes in the following sense: right adjoints going from top right corner to bottom left corner are naturally isomorphic. 
\end{theoremintro}

Here the functor $(-)^\circ$ is a direct generalization of the Sweedler duality functor defined in \cite{Sweedler69} for associative algebras, and the functor $(-)^\vee$ can be thought of as a topological dual. Then, by replacing (co)unital partial (co)operads on the left-hand side of the square with their "up to homotopy" counterparts, we construct another commuting square of duality adjunctions which are, in this case, all Quillen adjunctions. Using this enhanced duality square, we construct explicit cofibrant resolutions for a vast class of complete curved absolute partial operads. For instance, we construct a minimal cofibrant resolution complete of the curved absolute partial operads encoding curved Lie algebras, which provides us with a suitable notion of curved Lie algebra up to homotopy. Applications of these results will be found in Chapter \ref{Chapter 3}, where we develop the integration theory of curved absolute $\mathcal{L}_\infty$-algebras, and apply it to deformation theory and rational homotopy theory. 

\medskip

Finally, using our complete Bar construction, we extend the classical Homotopy Transfer Theorem to the case of algebras over a general class of cofibrant complete curved absolute operads. Let $V$ and $H$ be two complete pdg modules. A \textit{homotopy contraction} amounts to the data of 
\[
\begin{tikzcd}[column sep=5pc,row sep=3pc]
V \arrow[r, shift left=1.1ex, "p"{name=F}] \arrow[loop left]{l}{h}
&H~, \arrow[l, shift left=.75ex, "i"{name=U}]
\end{tikzcd}
\]
where $p$ and $i$ are two morphisms of filtered pre-differential graded (pdg) modules and $h$ is a morphism of filtered graded modules of degree $-1$, which satisfies standard conditions. This data allows us to construct a \textit{Van der Laan morphism} between the complete Bar constructions of the curved endomorphisms operads of $V$ and $H$. See \cite[Chapter 10]{LodayVallette12} for the classical Van der Laan morphism. Using this new morphism, we obtain the following result.

\begin{theoremintro}[Theorem \ref{thm: curved HTT}]
Let $\mathcal{C}$ be a dg counital partial cooperad. Let $H$ be a homotopy retract of $V$ and let 
\[
\varphi: \widehat{\Omega}\mathcal{C} \longrightarrow \mathrm{end}_V
\]
be a curved algebra structure on $V$. There is a transferred curved $\widehat{\Omega}\mathcal{C}$-algebra on $H$ constructed from the Van der Laan morphism.
\end{theoremintro}

In the particular case of complete curved $\mathcal{L}_\infty$-algebras, the transferred structure coincides with that constructed by K. Fukaya in \cite{FukayaHTT}. 

\section{Curved operads and their curved algebras}\label{Section: Curved Operads}
In the section, we introduce a new type of operad-like structure, called curved operads. They naturally encode curved algebras. The proper framework for this theory will be the underlying category of pre-differential graded modules. This framework, and the existence of the curved endomorphisms operad, was discovered independently of \cite{JoanCurved}.

\begin{Definition}[Pre-differential graded module]
A \textit{pre-differential graded module} (pdg module for short) $(V,d_V)$ is the data of a graded module $V$ together with a linear map $d_V: V \longrightarrow V$ of degree $-1$. A \textit{morphism} $f: (V,d_V) \longrightarrow (W,d_W)$ is a morphism of graded modules $f:V \longrightarrow W$ that commutes with the pre-differentials.
\end{Definition}

Pre-differential graded modules form a symmetric monoidal category $(\mathsf{pdg}~\mathsf{mod}, \otimes, \mathbb{K})~,$ where the pre-differential on the graded tensor product $A \otimes B$ is given by 
\[
d_{A \otimes B} (a \otimes b) \coloneqq d_A(a) \otimes b + (-1)^{|a|}\hspace{1pt} a \otimes d_B(b)~.
\]
The symmetric monoidal category $\mathsf{pdg}~\mathsf{mod}$ will be our \textit{underlying} category. Notice that we \textbf{do not} ask the condition $d_V^2 =0$. The category of dg modules is a full subcategory of the category of pdg modules.

\medskip

\subsection{First definitions} The data of a pdg operad is the data $(\PP,\gamma,\eta,d_\PP)$ of a graded operad $(\PP,\gamma,\eta)$ together with a derivation $d_\PP$ of degree $-1$. Similarly, a pdg partial operad amounts to a graded partial operad $(\PP,\{\circ_i\})$ equipped with a degree $-1$ derivation $d_\PP$. Again, a pdg cooperad $(\C,\Delta, \epsilon, d_\C)$ amounts to a graded cooperad $(\C,\Delta, \epsilon)$ with a coderivation $d_\C$ of degree $-1$ and a pdg partial cooperad $(\C,\{\Delta_i\},d_\C)$ amounts to a graded partial cooperad $(\C,\{\Delta_i\})$ with a coderivation $d_\C$ of degree $-1$. The notions of conilpotentcy for partial cooperads or completeness for partial operads are the same in this framework.

\begin{Definition}[Curved operad]
A \textit{curved operad} $(\mathcal{P},\gamma,\eta,d_\PP, \Theta_\mathcal{P})$ amounts to the data of a pdg operad $(\PP,\gamma,\eta,d_\PP)$ together with a morphism of pdg $\mathbb{S}$-modules $\Theta_\PP: (\I,0) \longrightarrow (\PP,d_\PP)$ of degree $-2$, such that the following diagram commutes: 
\[
\begin{tikzcd}[column sep=7pc,row sep=3pc]
\PP \arrow[r,"\mathrm{diag}"] \arrow[rrd,"(d_\PP)^2", bend right =10]
&\PP \oplus \PP \cong (\I \circ \PP) \oplus (\PP \circ \I) \arrow[r,"(\Theta_\PP ~ \circ~ \mathrm{id})~ -~(\mathrm{id}~ \circ'~ \Theta_\PP)"] 
& \PP \circ_{(1)} \PP \arrow[d,"\gamma_{(1)}"]\\
&
&\PP~,
\end{tikzcd}
\]
where $\mathrm{diag}$ is given by $\mathrm{diag}(\mu) \coloneqq (\mu,\mu)$. 
\end{Definition}

The data of the morphism of pdg $\mathbb{S}$-modules $\Theta_\PP$ is equivalent to a element $\Theta_\PP(\mathrm{id}) \coloneqq \theta$ in $\PP(1)_{-2}$ such that $d_\PP(\theta)=0$. The commutativity of the diagram amounts to the following condition: for every $\mu$ in $\PP(n)$, we ought to have
\[
\vcenter{\hbox{
\begin{tikzpicture}
\node at (3.25,1.5) {$d_\PP^2(\mu)$};

\node at (4.5,1.5) {$=$};

\node at (10,1.5) {$\displaystyle - \sum_{i=1}^n$};

\draw(13,-0.5)--(13,0.5);
\draw(13,0.5)--(12.5,1.5);
\draw(13,0.5)--(13.5,1.5);
\draw(13,0.5)--(13,2.5);
\draw(13,0.5)--(11.5,1.5);
\draw(13,0.5)--(14.5,1.5);
\node at (13,2) {$\bullet$};
\node at (13.3,0.4) {$\mu$};
\node at (13.3,2) {$\theta$};
\node at (13,2.75) {$i$};

\node at (12.5,1.7) {$\cdots$};
\node at (13.5,1.7) {$\cdots$};
\node at (11.5,1.7) {$1$};
\node at (14.5,1.7) {$n$};

\draw(7,-0.5)--(7,0.5);
\draw(7,0.5)--(6.5,1.5);
\draw(7,0.5)--(8.5,1.5);
\draw(7,0.5)--(7,1.5);
\draw(7,0.5)--(7.5,1.5);
\draw(7,0.5)--(5.5,1.5);
\draw(7,0.5)--(8.5,1.5);
\node at (7,1.7) {$\cdots$};
\node at (5.5,1.7) {$1$};
\node at (8.5,1.7) {$n$};
\node at (7.3,0.4) {$\mu$};
\node at (7.3,-0.2) {$\theta$};
\node at (7,-0.2) {$\bullet$};
\end{tikzpicture}.}}
\]

This states that the square of the pre-differential $d_\PP^2$ is equal to the operadic Lie bracket $[\theta,-]$ on the totalization of $\PP$. By a slight abuse of notation, we denote this equality by 
\[
d_\PP^2 = \gamma_{(1)} \cdot (\Theta_\PP \circ \mathrm{id} - \mathrm{id} \circ' \Theta_\PP)~,
\]
forgetting $\mathrm{diag}$ and the identifications made.

\begin{Notation}
A curved operad $(\PP,\gamma,\eta,d_\PP,\Theta_\PP)$ will be denoted by $(\PP,d_\PP,\Theta_\PP)$ for short, making the composition map and the unit implicit. When regarded as an element of $\mathcal{P}(1)$, the curvature is denoted $\theta$. Likewise for pdg operads, which will be sometimes denoted by $(\PP,d_\PP)$ when the context is clear.
\end{Notation}

\begin{Definition}
A \textit{morphism} of curved operads $f: (\PP,d_\PP,\Theta_\PP) \longrightarrow \allowbreak (\D,d_\D,\Theta_\D)$ amounts to the data of a pdg operad morphism $f: (\PP,d_\PP) \longrightarrow (\D,d_\D)$ that preserves the curvatures $f \circ \Theta_\PP = \Theta_\D$. 
\end{Definition}

Let $(V,d_V)$ be a pdg module, we consider the endomorphism pdg $\mathbb{S}$-module given by: 
\[
\mathrm{End}_V(n) \coloneqq \mathrm{Hom}(V^{\otimes n},V)~.
\]
Endowed with the composition of functions, and the pre-differential $\partial \coloneqq [d_V,-]~,$ it forms a pdg operad, called the \textit{endomorphism operad} of $V$.

\begin{lemma}\label{lemma: endo curved is curved}
Equipped with the curvature $\Theta_V: \I \longrightarrow \text{End}_V$ given by $\Theta_V(\mathrm{id}) = d_V^2$, the data $(\text{End}_V,\partial, \Theta_V)$ forms a curved operad, called the \textit{curved endomorphism operad} associated to $(V,d_V)$. 
\end{lemma}

\begin{proof}
It is straightforward to check that $\partial_V^2=[d_V^2,-]$ and that $[d_V,d_V^2]=0$.
\end{proof}

\begin{Definition}[Pdg $\PP$-algebra]
Let $(\PP,d_\PP)$ be a pdg operad. A \textit{pdg} $\PP$-\textit{algebra} $(A,\gamma_A,d_A)$ amounts to the data of a graded $\PP$-algebra structure $\gamma_A: \mathscr{S}(\PP)(A)\longrightarrow A$ on a pdg module $(A,d_A)$ such that 
\begin{equation}\label{derivation de P alg}
\gamma_A \cdot \left( \mathscr{S}(d_\PP)(\mathrm{id}) + \mathscr{S}(\mathrm{id})(\diracComb(\mathrm{id},d_A))\right) = d_A \cdot \gamma_A~. 
\end{equation}
A \textit{morphism} of pdg $\PP$-algebras $g: (A,\gamma_A,d_A) \longrightarrow (B,\gamma_B,d_B)$ amounts to a morphism of pdg modules $f:(A,d_A) \longrightarrow (B,d_B)$ that commutes with the $\PP$-algebra structures. The category of pdg $\PP$-algebras is denoted by $\mathsf{pdg}~\palg$.
\end{Definition}

\begin{Remark}
The notation $\mathscr{S}(\mathrm{id})(\diracComb(\mathrm{id},d_A))$ is explained in \ref{not: le dirac}.
\end{Remark}

\begin{Definition}[Curved $\PP$-algebra]
Let $(\PP,d_\PP,\Theta_\mathcal{P})$ be a curved operad. A \textit{curved} $\PP$-\textit{algebra} $(A,\gamma_A,d_A)$ amounts to the data of a pdg $\PP$-algebra structure $(A,\gamma_A,d_A)$ such that the following diagram commutes:
\begin{equation}\label{curved}
\begin{tikzcd}[column sep=5pc,row sep=2pc]
\mathscr{S}(\I)(A) \arrow[r,"\mathscr{S}(\Theta_\PP)(\mathrm{id})"] \arrow[rd, "d_A^2",swap]
&\mathscr{S}(\PP)(A) \arrow[d,"\gamma_A"] \\
&A~.
\end{tikzcd}
\end{equation}
Otherwise stated, we have that $\gamma_A(\theta,a) = d_A^2(a)$ for any $a$ in $A$. Morphisms of curved $\PP$-algebras are just morphisms of pdg $\PP$-algebras. The category of curved $\PP$-algebras will be denoted $\mathsf{curv}~\PP\text{-}\mathsf{alg}~.$
\end{Definition}

\begin{Remark}
Notice that a curved operad $(\PP,d_\PP,\Theta_\mathcal{P})$ is the data of a pdg operad $(\PP,d_\PP)$ together with extra \textit{structure} given by the curvature $\Theta_\mathcal{P}~.$ On the other hand, the data of a curved algebra structure over a curved operad is the same as the data of a pdg algebra structure over its underlying pdg operad, but this structure has to satisfy an extra \textit{property} given by diagram \ref{curved}.
\end{Remark}

There is an obvious inclusion functor $\mathsf{Inc}: \mathsf{curv}~\PP\text{-}\mathsf{alg} \hookrightarrow \mathsf{pdg}~\palg$ which is fully faithful.

\begin{Proposition}[{\cite[Proposition C.31]{JoanCurved}}]\label{reflexive}
Let $(\PP,d_\PP,\Theta_\mathcal{P})$ be a curved operad. The inclusion functor has a left adjoint 
\[
\mathsf{Curv}: \mathsf{pdg}~\palg \longrightarrow \mathsf{curv}~\PP\text{-}\mathsf{alg}~.
\]
Hence $\mathsf{curv}~\PP\text{-}\mathsf{alg}$ is a reflexive subcategory of $\mathsf{pdg}~\palg$. For a pdg $\PP$-algebra $(A,\gamma_A,d_A)$, its image under this functor is given by the following quotient: 
\[
\mathsf{Curv}(A) \coloneqq \frac{A}{(\mathsf{Im}(\gamma_A(\theta,-) - d_A^2(-)))}~.
\]
where $(\mathsf{Im}(\gamma_A(\theta,-) - d_A^2(-)))$ is the ideal generated by $\mathsf{Im}(\gamma_A(\theta,-) - d_A^2(-))$. Its pdg $\PP$-algebra structure is induced by $\gamma_A$ and $d_A$.
\end{Proposition}

\begin{proof}
First, lets show that $\mathsf{Im}(\gamma_A(\theta,-) - d_A^2(-))$ is stable under the pre-differential $d_A$. We have that $d_A(\gamma_A(\theta,a) - d_A^2(a)) = \gamma_A(d_\PP(\theta),a) + \gamma_A(\theta,d_A(a)) - d_A^3(a)~,$ which is again in $\mathsf{Im}(\gamma_A(\theta,-) - d_A^2(-))$ since $\gamma_A(d_\PP(\theta),a)=0~.$ Hence the quotient forms a pdg $\PP$-algebra, which is curved by definition. 

\medskip

Let $A$ be a pdg $\PP$-algebra and $B$ be a curved $\PP$-algebra. Since the pdg $\PP$-algebra structure on $B$ satisfies diagram \ref{curved}, any pdg $\PP$-algebra morphism $f: A \longrightarrow B$ factors through $\mathsf{Curv}(A)$. This gives the adjunction isomorphism by the universal property of the quotient. 
\end{proof}

\begin{Corollary}\label{bicomplete}
Let $(\PP,d_\PP,\Theta_\mathcal{P})$ be a curved operad. The category of curved $\PP$-algebras is complete and cocomplete. 
\end{Corollary}

\begin{proof}
The category $\mathsf{pdg}~\palg$ is both complete and cocomplete, since $\mathsf{curv}~\PP\text{-}\mathsf{alg}$ is a reflexive subcategory of $\mathsf{pdg}~\palg$ by Proposition \ref{reflexive}, it is also bicomplete.
\end{proof}

\begin{lemma}
Let $(\PP,d_\PP,\Theta_\mathcal{P})$ be a curved operad and let $(A,d_A)$ be a pdg module. The data of a curved $\PP$-algebra structure $\gamma_A$ on $(A,d_A)$ is equivalent to the data of a morphism of curved operads $\Gamma_A :(\PP,d_\PP,\Theta_\PP) \longrightarrow (\mathrm{End}_A, \partial, \Theta_{A})~.$
\end{lemma}

\begin{proof}
The data of a graded $\PP$-algebra structure $\gamma_A: \mathcal{S}(\PP)(A) \longrightarrow A$ is equivalent to the data of a morphism of graded operads $\Gamma_A : \PP \longrightarrow \mathrm{End}_A~.$ One can check that condition of equation \ref{derivation de P alg} is equivalent to $\Gamma_A$ being a morphism of pdg operads. The commutativity of the diagram \ref{curved} means that $\gamma_A (\mathscr{S}(\Theta_\PP)(\mathrm{id})) = d_A^2~,$ which is equivalent to $\Gamma_A (\Theta_\PP) = \Theta_{A}~.$
\end{proof}

Given a pdg module $(V,d_V)$, we can endow the \textit{coendomorphism operad} $\mathrm{Coend}_V$ with a curved operad structure using the same curvature $\Theta_V(\mathrm{id}) \coloneqq d_V^2$. 

\begin{Definition}[Curved $\mathcal{P}$-coalgebra]
A \textit{curved} $\PP$-\textit{coalgebra} structure on $(V,d_V)$ is the data of a morphism of curved operads $\Delta_V: (\PP,d_\PP,\Theta_\mathcal{P}) \longrightarrow (\mathrm{Coend}_V,\partial,\Theta_V)~.$ 
\end{Definition}

\begin{Remark}
The notion of a curved $\PP$-algebra can encode types of curved coalgebras which are not necessarily conilpotent.
\end{Remark}

\begin{Remark}\label{extends to curved case}
We introduce the analogue notion of a curved partial operad, since the condition imposed on the pre-differential involves the partial compositions of the operad and not the total compositions. The category of curved partial operads is denoted by $\mathsf{curv}~\mathsf{pOp}$. Furthermore, the category of curved unital partial operads is equivalent to the category of curved operads. See Proposition \ref{prop: unital partial operads are operads} for the non-curved statement.
\end{Remark}

\subsection{First examples}Here are some classical examples of curved algebras that one can encode with curved operads. See for instance \cite{PositselskiTwoKinds} for curved associative algebras. See for instance \cite{LazarevCocom} for curved Lie algebras.

\begin{Definition}[Curved Lie algebra]\label{clie}
A \textit{curved Lie algebra} $(\mathfrak{g},[-,-],d_{\mathfrak{g}},\vartheta)$  is the data of a graded Lie algebra $(\mathfrak{g},[-,-])$ together with a pre-differential $d_{\mathfrak{g}}$ of degree $-1$ which is a derivation with respect to the bracket $[-,-]$ , and a morphism of graded modules $\vartheta: \mathbb{K} \longrightarrow \mathfrak{g}$ of degree $-2$ such that:
\[
d_{\mathfrak{g}}^2=[\vartheta(1), -]~, \quad \text{and} \quad d_{\mathfrak{g}}(\vartheta(1))=0~.
\]
A \textit{morphism} $f: (\mathfrak{g}, d_\mathfrak{g},\vartheta_\mathfrak{g}) \longrightarrow (\mathfrak{h},d_\mathfrak{h},\vartheta_\mathfrak{h})$ is the data of a graded Lie algebra morphism $f: \mathfrak{g} \longrightarrow \mathfrak{h}$ that commutes with the pre-differentials and such that $f(\vartheta_\mathfrak{g}) = \vartheta_\mathfrak{h}~.$
\end{Definition}

Recall that the classical partial operad $\mathcal{L}ie$, encoding Lie algebras, can be defined as the free partial operad generated by one binary skew-symmetric operation, modulo the operadic ideal generated by the Jacobi relation. See \cite[Section 13.2]{LodayVallette12} for a complete account.

\medskip 

Let $M$ be the pdg $\mathbb{S}$-module given by $(\mathbb{K}.\zeta,0, \mathbb{K}.\beta,0,\cdots)$ with zero pre-differential, where $\zeta$ is an arity $0$ operation of degree $-2$, and $\beta$ is an arity $2$ operation of degree $0$, basis of the signature representation of $\mathbb{S}_2$. 

\begin{Definition}[$\Liec$ operad]
The \textit{curved partial operad} $\Liec$ is given by the free pdg partial operad generated by $M$ modulo the operadic ideal generated by the Jacobi relation on the generator $\beta$. It is endowed with the curvature $\Theta_\mathcal{L}$ given by $\Theta_\mathcal{L}(\mathrm{id}) \coloneqq \beta \circ_1 \zeta~.$  
\end{Definition}

\begin{lemma}\label{lemmalie}
The data $(\Liec, 0, \Theta_\mathcal{L})$ forms a curved partial operad. Furthermore, the category of curved $\Liec$-algebras is isomorphic to the category of curved Lie algebras.
\end{lemma}

\begin{proof}
First, let us check that $(\Liec, 0, \Theta_\mathcal{L})$ is indeed a curved partial operad. The pre-differential is $0$, hence $d(\theta)=0$. Now we need to check that $[\theta,-] = 0$ . Since $[\theta,-]$ is a derivation, it is enough to check it on the two generators of $\Liec$. The result $[\theta, \vartheta] = 0$ follows from the anti-symmetry of $\beta$. A straightforward computation gives that $[\theta, \beta] = 0$, using the Jacobi relation.  

\medskip

A pdg algebra over the curved partial operad $\Liec$ amounts to a pdg-module $(A,d_A)$ endowed with a graded Lie bracket $[-,-] \coloneqq \gamma_A(\beta)$ and a morphism $\vartheta \coloneqq \gamma_A(\zeta): \kk \longrightarrow \mathfrak{g}$ such that $d_A(\vartheta(1))=0$.
They form a curved $\Liec$-algebra if and only if Diagram \ref{curved} commutes, which equivalent to  $[\vartheta,-] = d_A^2(-)$. Morphisms of curved $\Liec$-algebras must commute with the structure, i.e: the bracket and curvature $\vartheta$; and with the underlying pre-differentials. 
\end{proof}

\begin{Definition}[Curved associative algebra]\label{cass}
A \textit{curved associative algebra} $(A,\mu_A,d_A,\vartheta)$ amounts to the data a non-unital graded associative algebra $(A,\mu_A)$ together with a pre-differential $d_A$ of degree $-1$ which is a derivation with respect to the associative product $\mu_A$, and a morphism of graded modules $\vartheta: \mathbb{K} \longrightarrow A$ of degree $-2$ such that:
\[
d_A^2=\mu_A(\vartheta(1),-) - \mu_A(-, \vartheta(1))~, \quad \text{and} \quad d_A(\vartheta(1))=0~.
\]
A \textit{morphism} $f: (A,\mu_A, d_A,\vartheta_A) \longrightarrow (B,\mu_B, d_B, \vartheta_B)$ is the data of a morphism of graded non-unital associative algebras $f:(A,\mu_A) \longrightarrow (B,\mu_B)$ that commutes with the pre-differentials and preserves the curvatures $f(\vartheta_A) = \vartheta_B~.$
\end{Definition}

Let $N$ be the pdg $\mathbb{S}$-module given by $(\mathbb{K}.\phi,0,\mathbb{K}[\mathbb{S}_2].\mu,0,\cdots)$ with zero pre-differential, where $\phi$ is an arity $0$ operation of degree $-2$, and $\mu$ is an binary operation of degree $0$, basis of the regular representation of $\mathbb{S}_2$.

\begin{Definition}[$\Assc$ operad]
The \textit{curved partial operad} $\Assc$ is given by the free pdg partial operad generated by $N$ modulo the operadic ideal generated by the associativity relation on the generator $\mu$. It is endowed with the curvature $\Theta_\mathcal{A}$ given by $\Theta_\mathcal{A}(\mathrm{id}) \coloneqq \mu \circ_1 \phi - \mu \circ_2 \phi.$ 
\end{Definition}

\begin{lemma}
The data $(\Assc, 0, \Theta_\mathcal{A})$ forms a curved partial operad. Furthermore, the category of curved $\Assc$-algebras is isomorphic to the category of curved associative algebras.
\end{lemma}

\begin{proof}
The proof of this lemma follows the same steps as the proof of Lemma \ref{lemmalie}. 
\end{proof}

\begin{Remark}
See Appendix \ref{Appendix B} for the "absolute analogues" of these curved operads.
\end{Remark}

\begin{Proposition}\label{assliem}
The classical morphism of partial operads $\mathcal{L}ie \longrightarrow \mathcal{A}ss$ induced by $\beta \mapsto \mu - \mu^{(12)}~,$ extends to a morphism of curved partial operads $\Liec \longrightarrow \Assc$ by sending $\zeta$ to $\phi$.
\end{Proposition}

\begin{proof}
It is straightforward to check that the latter morphism commutes with the respective curvatures. 
\end{proof}

\begin{Corollary}\label{Antisymmetrization}
The morphism $\Liec \longrightarrow \Assc$ induces a functor
\[
\begin{tikzcd}[column sep=3pc,row sep=0pc]
\mathsf{Skew}: \mathsf{curv}~\Assc\text{-}\mathsf{alg} \arrow[r]
&\mathsf{curv}~\Liec\text{-}\mathsf{alg} \\
(A,\mu_A,d_A,\vartheta) \arrow[r,mapsto]
&(A,[-,-] \coloneqq \mu_A - \mu_A^{(12)},d_A,\vartheta)~.
\end{tikzcd}
\]
given by the skew-symmetrization of the associative product. 
\end{Corollary} 

\begin{proof}
Given a morphism of curved partial operads $\Gamma_A: \Assc \longrightarrow \mathrm{End}_A$, one pulls back along $\Liec \longrightarrow \Assc$.
\end{proof}

\begin{Remark}
In the next section we will construct the left adjoint of this functor using curved bimodules.
\end{Remark}

\section{Curved bimodules and universal functors}\label{Section: Curved operadic bimodules}
The goal of this section is to define the "curved" generalization of modules over operads. Their appeal comes from the fact that operadic modules encode universal functors relating categories of algebras over operads. The generalization is not immediate since in a curved operad, the curvature condition intertwines the underlying category with the operadic structure. Nevertheless, there exists a good notion of bimodule that encodes universal functors between the categories of curved algebras over curved operads. Note that in order for a bimodule to define a functor between categories of algebras over operads, we need to use the monoidal definition of operads, like in \cite{Fresse09}. If one is working with curved partial operads, one can simply add an augmented unit to fit into this framework.

\medskip

\subsection{Curved bimodules} One major difference between curved operads and classical operads is the following fact. Let $\cP$ be a curved operad and $(V,d_V)$ be a pdg module. Then $\mathscr{S}(\PP)(V)$ endowed with its canonical $\PP$-algebra structure does not form a curved $\PP$-algebra in general. The pre-differential on $\mathscr{S}(\PP)(V)$ is given by: 
\[
d_{\mathscr{S}(\PP)(V)} = \mathscr{S}(d_\PP)(\mathrm{id}_V) + \mathscr{S}(id_\PP)(\diracComb(\mathrm{id}_V,d_V))~.
\]
Hence we have that 
\[
d_{\mathscr{S}(\PP)(V)}^2 = \mathscr{S}(d_\PP^2)(\mathrm{id}_V) + \mathscr{S}(\mathrm{id}_\PP)(\diracComb(\mathrm{id}_V,d_V^2))~,
\]
because of the Koszul sign convention. Since $\PP$ is a curved operad:
\[
d_\PP^2 = \gamma_{(1)} \cdot (\Theta_\PP \circ \mathrm{id}_\PP - \mathrm{id}_\PP \circ' \Theta_\PP)~.
\]
Therefore
\begin{align*}
d_{\mathscr{S}(\PP)(V)}^2 = &~\mathscr{S}(\gamma)(\mathrm{id}_V) \cdot \Big( \mathscr{S}(\Theta_\PP) \circ \mathscr{S}(\mathrm{id}_\PP)(\mathrm{id}_V) - \mathscr{S}(\mathrm{id}_\PP) \circ \mathscr{S}(\diracComb(\mathrm{id}_\PP,\Theta_\PP))(\mathrm{id}_V) \Big) \\ 
&+ \mathscr{S}(\mathrm{id}_\PP)(\diracComb(\mathrm{id}_V,d_V^2))
\end{align*}
and does not satisfy the condition imposed by the diagram \ref{curved}. In general
\[
d_{\mathscr{S}(\PP)(V)}^2 \neq \mathscr{S}(\gamma)(\mathrm{id}_V) \cdot \Big(\mathscr{S}(\Theta_\PP) \circ \mathscr{S}(\mathrm{id}_\PP)(\mathrm{id}_V)\Big)~,
\]
even if we impose the extra condition that $d_V^2=0~.$

\begin{Proposition}[{\cite[Proposition C.31]{JoanCurved}}]\label{freecurvedalg}
Let $\cP$ be a curved operad and $(V,d_V)$ be a pdg module. The free curved $\PP$-algebra is given by:
\[
\mathsf{F}(\PP)(V) \coloneqq \frac{\mathscr{S}(\PP)(A)}{\mathsf{Im}\Big(\mathscr{S}(\mathrm{id}_\PP)(d_V^2) - \mathscr{S}(\Theta_\PP)(\mathrm{id}_V)\Big)}~.
\]
\end{Proposition}

\begin{proof}
One can check that $\mathsf{F}(\PP)(V)$ is equal to the composition of the free pdg $\PP$-algebra functor $\mathscr{S}(\PP)(V)$ followed by the reflector $\mathsf{Curv}$ constructed in Proposition \ref{reflexive}.
\end{proof}

This type of construction is encoded by a left $\PP$-module in the classical case. The above Proposition shows that in the curved case, modules over operads are more intricate, since it is not obvious that this quotient can be interpreted as a "curved" left $\PP$-module. We begin by recalling standard definitions.

\begin{Definition}[Left and right $\PP$-modules]
Let $(\PP, \gamma,\eta, d_\PP)$ be a pdg operad. 
\begin{enumerate}
\item A \textit{left pdg} $\PP$-\textit{module} $(N,\lambda,d_N)$ is the data of a pdg $\mathbb{S}$-module $(N,d_N)$ together with a morphism of pdg $\mathbb{S}$-modules $\lambda: \PP \circ N \longrightarrow N$, such that the following diagram
\[
\begin{tikzcd}[column sep=3pc,row sep=3pc]
\PP \circ \PP \circ N \arrow[r, "\gamma_\PP ~ \circ ~ \mathrm{id}"] \arrow[d,swap,"\mathrm{id} ~ \circ ~ \lambda"]
&\PP \circ N \arrow[d,"\lambda"]\\
\PP \circ N \arrow[r,"\lambda"]
&N 
\end{tikzcd}
\]
commutes and such that $\lambda \cdot (\eta \circ \mathrm{id}_N) = \mathrm{id}_N~.$ 

\item A \textit{right pdg} $\PP$-\textit{module} $(M,\rho,d_M)$ is the data of a pdg $\mathbb{S}$-module $(M,d_M)$ together with a morphism of pdg $\mathbb{S}$-modules $\rho: M \circ \PP \longrightarrow \PP$, such that the following diagram
\[
\begin{tikzcd}[column sep=3pc,row sep=3pc]
M \circ \PP \circ \PP \arrow[r, "\mathrm{id} ~ \circ ~ \gamma"] \arrow[d,swap,"\rho ~ \circ ~ \mathrm{id}"]
&M \circ \PP \arrow[d,"\rho"]\\
M \circ \PP  \arrow[r,"\rho"]
&M
\end{tikzcd}
\]
commutes and such that $\rho \cdot (\mathrm{id}_M \circ \rho) = \mathrm{id}_M~.$ 
\end{enumerate}
\end{Definition}

\begin{Remark}
One can identify left pdg $\PP$-modules concentrated in arity $0$ with the category of pdg $\PP$-algebras.
\end{Remark} 

\begin{Definition}[Relative composition product]
Let $(\PP, d_\PP)$ and $(\Q, d_\Q)$ be two pdg operads. Let $(N,\lambda_N,d_N)$ be a left pdg $\Q$-module and let $(M,\rho_M,d_M)$ be a right pdg $\PP$-module, the \textit{relative composition product} $M \circ_\Q N$, given by the following coequalizer
\[
M \circ_\Q N \coloneqq
\begin{tikzcd}[column sep=4.5pc,row sep=0pc]
\mathsf{Coeq}~ \Big( M \circ \Q \circ N \arrow[r,shift left=.75ex,"\rho_M ~ \circ ~ \mathrm{id}"] \arrow[r,shift right=.75ex,swap,"\mathrm{id} ~ \circ ~ \lambda_N"]
&M \circ N \Big) ~,
\end{tikzcd}
\]
in the category of pdg $\mathbb{S}$-modules. 
\end{Definition} 

\begin{Definition}[Operadic bimodule]
Let $(\PP, d_\PP)$ and $(\Q, d_\Q)$ be two pdg operads. A \textit{pdg} $(\PP, \Q)$-\textit{bimodule} is the data $(M,\lambda,\rho,d_M)$ of a pdg $\mathbb{S}$-module $(M,d_M)$ together with two morphisms of pdg $\mathbb{S}$-modules $\lambda: \PP \circ M \longrightarrow M$ and $\rho: M \circ \Q \longrightarrow M$ which endow $M$ with a left $\PP$-module structure and a right $\PP$-module structure. Those structures are compatible with each other in the following sense: 
\[
\begin{tikzcd}[column sep=3pc,row sep=3pc]
\PP \circ M \circ \Q \arrow[d,swap, "\mathrm{id} ~ \circ ~ \rho "] \arrow[r,"\lambda ~ \circ ~ \mathrm{id}"]
&M \circ \Q \arrow[d,"\rho"] \\
\PP \circ M  \arrow[r,"\lambda"]
&M~.
\end{tikzcd}
\]
\end{Definition}

Operadic bimodules encode functors between the categories of algebras over operads, see \cite[Chapter 9]{Fresse09} for a detailed account. 

\begin{Definition}[Relative Schur functor of a bimodule]
Let $(M,\lambda,\rho,d_M)$ be a pdg $(\PP,\Q)$-bimodule. The \textit{relative Schur functor} $\mathscr{S}_\Q(M)(-)$ associated to $M$ is given, for $(A,\gamma_A,d_A)$ a pdg $\Q$-algebra, by the following coequalizer:
\[
\mathscr{S}_\Q(M)(A) \coloneqq
\begin{tikzcd}[column sep=4.5pc,row sep=0pc]
\mathsf{Coeq}~ \Big( \mathscr{S}(M) \circ \mathscr{S}(\Q)(A) \arrow[r,shift left=.75ex,"\mathscr{S}(\rho)(\mathrm{id}_A)"] \arrow[r,shift right=.75ex,swap,"\mathscr{S}(\mathrm{id}_M)(\gamma_A)"]
&\mathscr{S}(M)(A) \Big) ~.
\end{tikzcd}
\]
It defines a functor: 
\[
\mathscr{S}_\Q(M)(-): \mathsf{pdg}~\Q\text{-}\mathsf{alg} \longrightarrow \mathsf{pdg}~\PP\text{-}\mathsf{alg}~.
\]
\end{Definition}

Under the identification of left pdg $\Q$-modules concentrated in arity $0$ with pdg $\Q$-algebras, $\mathscr{S}_\Q(M)(A)$ is in fact given by the relative composition product $M \circ_\Q A$ of the left pdg $\Q$-module $(A,\gamma_A,d_A)$ with the $(\PP, \Q)$-bimodule $(M,\lambda,\rho,d_M)$.

\begin{Definition}[Curved Bimodules]\label{curvedbimodule}
Let $(\PP, d_\PP,\Theta_\PP)$ and $(\Q, d_\Q,\Theta_\Q)$ be two curved operads. A \textit{curved} $(\PP, \Q)$-\textit{bimodule} is the data of a pdg $(\PP, \Q)$-bimodule $(M,\lambda_M,\rho_M,d_M)$ such that the following diagram commutes:
\[
\begin{tikzcd}[column sep=7pc,row sep=3pc]
M \arrow[r,"\mathrm{diag}"] \arrow[rrd,swap,"d_M^2", bend right = 10]
&(\I \circ M) \oplus (M \circ \I) \arrow[r,"(\Theta_\PP~ \circ~ \mathrm{id}) - (\mathrm{id}~ \circ' ~\Theta_\Q)"] 
&(\PP \circ M) \oplus (M \circ_{(1)} \Q) \arrow[d,"\lambda + \rho_{(1)}"]\\
&
&M~.
\end{tikzcd}
\]
This equality is denoted by $d_M^2 = \lambda \cdot (\Theta_\PP \circ \mathrm{id}) - \rho_{(1)} \cdot (\mathrm{id} \circ' \Theta_\Q)~.$ 
\end{Definition}

\begin{Remark}
A curved generalization of infinitesimal bimodules for curved partial operads is immediate since the curvature is an infinitesimal notion. For instance, curved infinitesimal bimodules are the coefficients for the André-Quillen cohomology of a curved operad, which can be defined using the same methods as for operads and standard bimodules. See \cite{MerkulovVallette09I} for more details.
\end{Remark}

\begin{Proposition} \label{functorbimod}
Any curved $(\PP, \Q)$-bimodule $(M,\lambda,\rho,d_M)$ induces a functor 
\[
\mathscr{S}_\Q(M)(-): \mathsf{curv}~\Q\text{-}\mathsf{alg} \longrightarrow \mathsf{curv}~\PP\text{-}\mathsf{alg}~,
\]
given by the relative Schur functor associated to $M$ as a pdg $(\PP, \Q)$-bimodule.
\end{Proposition}

\begin{proof}
Let $(A,\gamma_A,d_A)$ be a curved $\Q$-algebra. Let $\pi_A: \mathscr{S}(M)(A) \twoheadrightarrow \mathscr{S}_\Q(M)(A)$ be the projection map. The pre-differential of $\mathscr{S}_\Q(M)(A)$ is given by the image in the quotient of 
\[
\mathscr{S}(d_M)(\mathrm{id}_A) + \mathscr{S}(\mathrm{id}_M)(\diracComb(\mathrm{id}_A,d_A))~,
\]
which we will denote $d$ for simplicity. The pdg $\PP$-algebra structure on $\mathscr{S}_\Q(M)(A)$ is given by the map $\pi_A \cdot (\mathscr{S}(\lambda)(\mathrm{id}))~.$ We need to check that $\mathscr{S}_\Q(M)(A)$ is indeed a curved $\PP$-algebra, i.e: that the following diagram commutes:
\[
\begin{tikzcd}[column sep=5.5pc,row sep=3pc]
\mathscr{S}(\I)(\mathscr{S}_\Q(M)(A)) \arrow[r,"\mathscr{S}(\Theta_\PP)(\mathrm{id})"] \arrow[rd,"d^2",swap]
&\mathscr{S}(\PP)(\mathscr{S}_\Q(M)(A)) \arrow[d,"\pi_A \cdot (\mathscr{S}(\lambda_M)(\mathrm{id}))"]\\
&\mathscr{S}_\Q(M)(A)~.
\end{tikzcd}
\]
We know that $d^2$ will be induced by $\mathscr{S}(d_M^2)(\mathrm{id}_A) + \mathscr{S}(\mathrm{id}_M)(\diracComb(\mathrm{id}_A,d_A^2))$. On one hand, since $M$ is a curved $(\PP, \Q)$-bimodule, $d_M^2 = \lambda \cdot (\Theta_\PP \circ \mathrm{id}) - \rho_{(1)} \cdot (\mathrm{id} \circ' \Theta_\Q)~.$ On the other hand, since $A$ is a curved $\Q$-algebra, $d_A^2 = \gamma_A \cdot (\mathscr{S}(\Theta_\Q)(\mathrm{id}))~.$ Therefore: 
\begin{align*}
d^2 = &~\pi_A \cdot (\mathscr{S}(\lambda)(\mathrm{id}_A)) \cdot (\mathscr{S}(\Theta_\PP) \circ \mathscr{S}(\mathrm{id}_M)(\mathrm{id}_A)) - \pi_A \cdot (\mathscr{S}(\rho)(\mathrm{id}_A)) \cdot (\mathscr{S}(\diracComb(\mathrm{id}_M,\Theta_\Q)(A)) \\
	  &+\pi_A \cdot \mathscr{S}(\mathrm{id}_M)\Big(\diracComb(\mathrm{id}_A,\gamma_A \cdot (\mathscr{S}(\Theta_\Q)(\mathrm{id})))\Big)~.	 
\end{align*}
By definition of $\mathscr{S}_\Q(M)(A)~,$ we have that the right action of $\Q$ on $M$ given by $\rho$ is equal to the action of $\Q$ on $A$ given by $\gamma_A$, hence the last two terms are equal and therefore cancel each other. As a result:
\[
d^2 = \pi_A \cdot (\mathscr{S}(\lambda)(\mathrm{id}_A)) \cdot (\mathscr{S}(\Theta_\PP) \circ \mathscr{S}(\mathrm{id}_M)(\mathrm{id}_A))~,
\]
which is precisely the condition imposed by the above diagram. 
\end{proof}

\begin{Example}
Given a curved operad $(\PP, d_\PP,\Theta_\mathcal{P})$, its underlying pdg $\mathbb{S}$-module $(\PP, d_\PP)$ can be endowed with a canonical curved $(\PP,\PP)$-bimodule structure given by the composition of $\PP$. This curved bimodule encodes the identity endofunctor of curved $\PP$-algebras. 
\end{Example}

\begin{theorem}\label{Ind and Res adjunction}
Let $(\PP, d_\PP,\Theta_\PP)$ and $(\Q, d_\Q,\Theta_\Q)$ be two curved operads, and let $f: \PP \longrightarrow \Q$ be a morphism of curved operads. The morphism $f$ induces an adjunction at the level of curved algebras:
\[
\begin{tikzcd}[column sep=7pc,row sep=3pc]
            \mathsf{Ind}_f: \mathsf{curv}~\PP\text{-}\mathsf{alg} \arrow[r, shift left=1.1ex, ""{name=F}] & \mathsf{curv}~\Q\text{-}\mathsf{alg} : \mathsf{Res}_f \arrow[l, shift left=.75ex, ""{name=U}]
            \arrow[phantom, from=F, to=U, , "\dashv" rotate=-90]
\end{tikzcd}
\]
given by the restriction and the induction functors. 
\end{theorem}

\begin{proof}
The right adjoint is given by the restriction along $f$ as follows. Given a curved $\D$-algebra structure on a pdg module $(A,d_A)$, that is, a morphism of curved operads $\Gamma_A: \Q \longrightarrow \mathrm{End}_A~,$ one can always pre-compose $\Gamma_A$ with $f$. This curved operadic morphism $\Gamma_A \cdot f: \PP \longrightarrow \mathrm{End}_A$ endows $(A,d_A)$ with a curved $\PP$-algebra structure, denoted by $\mathsf{Res}_f(A)~.$

\medskip

On the other hand, the morphism $f$ allows us to endow $\Q$ with a curved $(\Q,\PP)$-bimodule structure, where right action of $\PP$ is given by: 
\[
\begin{tikzcd}[column sep = 4pc,row sep=2pc]
\rho :\Q \circ \PP \arrow[r,"\mathrm{id}_\Q ~\circ ~f"]
&\Q \circ \Q  \arrow[r,"\gamma_\Q"]
&\Q~,
\end{tikzcd}
\]
and the left action of $\Q$ is simply given by $\gamma_\Q~.$ Let us check that this endows $\Q$ with a curved bimodule structure. We have that:
\[
d_\Q^2 = \gamma_\Q \cdot (\Theta_\Q \circ \mathrm{id}) - (\gamma_\Q)_{(1)} \cdot (\mathrm{id} \circ' \Theta_\Q) = \gamma_\Q \cdot (\Theta_\Q \circ \mathrm{id}) - \rho_{(1)} \cdot (\mathrm{id} \circ' \Theta_\PP)
\]
since $f \cdot \Theta_\PP  = \Theta_\Q~.$ Thus, the curved $(\Q,\PP)$-bimodule $\Q$ induces a functor 
\[
\mathsf{Ind}_f(-) \coloneqq \mathscr{S}_\PP(\Q)(-): \mathsf{curv}~\PP\text{-}\mathsf{alg} \longrightarrow \mathsf{curv}~\Q\text{-}\mathsf{alg}~.
\]
Since morphism of curved algebras are morphisms of pdg algebras, this general operadic construction is still an adjunction. 
\end{proof}

One important application of this new result is the construction of the curved universal enveloping algebra of a curved Lie algebra. To the best of our knowledge, this construction is new. The morphism $\Liec \longrightarrow \Assc$, defined in Proposition \ref{assliem}, induces the following adjunction.

\begin{Corollary}[Curved universal enveloping algebra]
There an adjunction between the category of curved Lie algebras and the category of curved associative algebras
\[
\begin{tikzcd}[column sep=7pc,row sep=3pc]
            \mathfrak{U} : \mathsf{curv}~\Liec\text{-}\mathsf{alg} \arrow[r, shift left=1.1ex, ""{name=F}] & \mathsf{curv}~\Assc \text{-}\mathsf{alg} : \mathsf{Skew}~, \arrow[l, shift left=.75ex, ""{name=U}]
            \arrow[phantom, from=F, to=U, , "\dashv" rotate=-90]
\end{tikzcd}
\]
where $\mathfrak{U}$ denotes the curved universal enveloping algebra and $\mathsf{Skew}$ denotes the functor obtained by the skew-symmetrization of the associative product of Corollary \ref{Antisymmetrization}.
\end{Corollary}

\begin{proof}
Simply apply Theorem \ref{Ind and Res adjunction} to the morphism $f: \Liec \longrightarrow \Assc$, and set $\mathfrak{U} \coloneqq \mathsf{Ind}_f$ and $\mathsf{Anti} \coloneqq \mathsf{Res}_f~.$
\end{proof}

\begin{Proposition}
Let $\mathfrak{g}$ be a curved Lie algebra. Its universal enveloping algebra is isomorphic to
\[
\mathfrak{U}(\mathfrak{g}) \cong \frac{\overline{\mathcal{T}}(\mathfrak{g})}{\left(x \otimes y - (-1)^{|x|} y \otimes x - [x,y]\right)}~,
\]
where $\overline{\mathcal{T}}(-)$ denotes the non-unital tensor algebra, endowed with the curvature $\vartheta: \kk \longrightarrow \mathfrak{g} \hookrightarrow \mathfrak{U}(\mathfrak{g})$. 
\end{Proposition}

\begin{proof}
One can check this adjunction by hand. Indeed, the adjunction holds between between pdg $c\mathcal{L}ie$-algebras and pdg $c\mathcal{A}ss$-algebras. It is straightforward to check that $\mathfrak{U}(\mathfrak{g})$ endowed with the curvature $\vartheta$ forms a curved $c\mathcal{L}ie$-algebra. 
\end{proof}

\begin{Remark}\label{Rmk: PBW curved}
One has that 
\[
c\mathcal{A}ss \cong \mathcal{A}ss ~ \sqcup ~ \phi.\kk~,
\]
as a pdg right $c\mathcal{L}ie$-module, where $\phi.\kk$ is an arity one degree $-2$ operation. This implies that 
\[
c\mathcal{A}ss \circ_{c\mathcal{L}ie} - ~~ \cong \mathcal{C}om \circ -
\]
since $\mathcal{A}ss \cong \mathcal{C}om \circ \mathcal{L}ie$ as a right $\mathcal{L}ie$-module. Thus the universal enveloping algebra of a curved Lie algebra also satisfies the Poincaré-Birkoff-Witt property. See \cite{Tamaroff2020}. 
\end{Remark}

This machinery can also be applied to curved mixed $\mathcal{L}_\infty$-algebras. These are curved $\mathcal{L}_\infty$-algebras with two pre-differentials, one coming from the underlying pdg module and another one from the structure, such that their difference satisfies the axioms of a classical curved $\mathcal{L}_\infty$-algebra. See Appendix \ref{Section: Appendix B} of Chapter \ref{Chapter 3}. These algebras are exactly curved algebras over $\widehat{\Omega}(u\mathcal{C}om^*)$, where $\widehat{\Omega}$ is the complete Cobar construction of Section \ref{Section: Constructions Bar-Cobar operadiques}.

\begin{Corollary}[Curved $\mathcal{A}_\infty$ universal enveloping algebra]
There is the following adjunction between the category of curved mixed $\mathcal{L}_\infty$-algebras and the category of curved mixed $\mathcal{A}_\infty$ algebras given by:
\[
\begin{tikzcd}[column sep=7pc,row sep=3pc]
            \mathfrak{U} : \mathsf{curv.mix}~\mathcal{L}_\infty\text{-}\mathsf{alg} \arrow[r, shift left=1.1ex, ""{name=F}] & \mathsf{curv.mix}~\mathcal{A}_\infty \text{-}\mathsf{alg} : \mathsf{Anti}~, \arrow[l, shift left=.75ex, ""{name=U}]
            \arrow[phantom, from=F, to=U, , "\dashv" rotate=-90]
\end{tikzcd}
\]
where $\mathfrak{U}$ denotes the curved $\mathcal{A}_\infty$ universal enveloping algebra.
\end{Corollary}

\begin{proof}
There is a morphism of counital partial cooperads $\iota: u\mathcal{C}om^* \longrightarrow u\mathcal{A}ss^*$ given by the linear dual of the classical morphism of unital partial operads $u\mathcal{A}ss \twoheadrightarrow u\mathcal{C}om$. Therefore, using the complete Cobar construction of Section \ref{Section: Constructions Bar-Cobar operadiques}, we get a morphism of curved partial operads $\widehat{\Omega}(\iota): \widehat{\Omega}(u\mathcal{C}om^*) \longrightarrow \widehat{\Omega}(u\mathcal{A}ss^*)$. By Theorem \ref{Ind and Res adjunction} it induces an adjunction.
\end{proof}

\subsection{The analogues of left and right modules}One may ask whether there exists a natural notion of a left (resp. right) curved $\PP$-module. This is in fact not obvious. Since every algebra over an operad is a particular case of left module in the classical setting, a naive generalization could be: a left pdg $\PP$-module $(M,\lambda,d_M)$ is a  left curved $\PP$-module if the following diagram commutes
\[
\begin{tikzcd}[column sep=5pc,row sep=2pc]\label{curvedmod}
\I \circ M \arrow[r,"\Theta_\PP \hspace{1pt} \circ \hspace{1pt} \mathrm{id}"] \arrow[rd, "d_M^2",swap]
&\PP \circ M \arrow[d,"\lambda"]\\
&M~.
\end{tikzcd}
\]
But in this case $\mathscr{S}(M)(-)$ does not produce a functor $\mathsf{pdg}~\mathsf{mod} \longrightarrow \mathsf{curv}~\PP\text{-}\mathsf{alg}~,$ since $\mathscr{S}(M)(A)$ is not, in general, a curved $\PP$-algebra. One can notice  by doing a straightforward computation that for this notion of curved left $\PP$-module, $\mathscr{S}(M)(A)$ is a curved $\PP$-algebra if and only if $d_A^2=0~.$

\medskip

The conceptual explanation is that the data of a dg module is equivalent to the data of a curved algebra over the curved operad $(\I, 0, 0)$, where $\I$ is the trivial operad with the zero pre-differential and the zero curvature. The above-mentioned naive definition of a left curved $\PP$-module is in fact a curved $(\PP, \I)$-bimodule, which is coherent with the fact that this notion encodes functors $\mathsf{dg}~\mathsf{mod} \longrightarrow \mathsf{curv}~\PP\text{-}\mathsf{alg}~.$ Likewise, a symmetric naive definition of a right curved $\PP$-module would in fact be a $(\I, \PP)$-bimodule that would encode functors $\mathsf{curv}~\PP\text{-}\mathsf{alg} \longrightarrow \mathsf{dg}~\mathsf{mod}$ .

\begin{Remark}
Let $(\PP, \gamma,\eta, d_\PP, \Theta_\PP)$ be a curved operad. The unit $\eta: \I \longrightarrow \PP$ is not a morphism of curved operads. Otherwise, we could endow $\PP$ with a curved $(\PP, \I)$-bimodule structure and $\mathscr{S}(\PP)(A)$ would be a curved $\PP$-algebra for $(A,d_A)$ a dg module, which is not the case.
\end{Remark}

The natural question is then to find a curved operad that plays the same role as the unit operad $\I$ plays for classical operads in the theory of bimodules. Otherwise stated, a curved operad that encodes pdg modules as its curved algebras.

\begin{Definition}[$\mathcal{I}\mathcal{C}$ operad]
The pdg operad $\mathcal{IC}$ is the free pdg operad generated by an operation of arity $1$ and degree $-2$. It is given by the pdg $\mathbb{S}$-module $(0,\mathbb{K}.\mathrm{id} ~\oplus~ \bar{S}(\theta),0,\cdots)$ with zero pre-differential, where $\bar{S}(\theta)$ denotes the free non-unital commutative algebra generated by $\theta~.$
\end{Definition}

\begin{lemma}
The pdg operad $\mathcal{IC}$, endowed the curvature $\Theta_{\mathcal{IC}}(\mathrm{id}) \coloneqq \theta$ forms a curved operad. The category of curved $\mathcal{IC}$-algebras is isomorphic to the category of pdg modules. 
\end{lemma}

\begin{proof}
The commutator with $\theta$ is zero, hence $\mathcal{IC}$ is a curved operad. Let $(V,d_V)$ be a pdg module and let $\Gamma: \mathcal{IC} \longrightarrow \mathrm{End}_V$ be a morphism of curved operads. The image of $\theta$ determines $\Gamma$, since it is a morphism of curved operads then $\Gamma(\theta) = d_V^2~.$ 
\end{proof}

For any curved operad $(\PP, d_\PP, \Theta_\PP)$, there is an unique morphism of curved operads $\varphi_\PP: \mathcal{IC} \longrightarrow \PP$ given by $\varphi_\PP(\theta) = \theta_\PP$, where $\theta_\PP = \Theta_\PP(\mathrm{id})$. Therefore the pdg $\mathbb{S}$-module $\PP$ can be endowed canonically with a curved $(\mathcal{IC},\PP)$-bimodule structure and with a curved $(\PP,\mathcal{IC})$-bimodule structure.

\begin{Proposition}
Let $(\PP, d_\PP, \Theta_\PP)$ be a curved operad and let 
\[
\begin{tikzcd}[column sep=7pc,row sep=3pc]
            \mathsf{F}(\PP) : \mathsf{pdg}~\mathsf{mod} \arrow[r, shift left=1.1ex, ""{name=F}] & \mathsf{curv}~\PP \text{-}\mathsf{alg} : \mathsf{U}~, \arrow[l, shift left=.75ex, ""{name=U}]
            \arrow[phantom, from=F, to=U, , "\dashv" rotate=-90]
\end{tikzcd}
\]
be the free-forgetful adjunction of Proposition \ref{freecurvedalg}. The free functor $\mathsf{F}(\PP)$ is naturally isomorphic to the functor $\mathscr{S}_{\mathcal{IC}}(\PP)$ given by the canonical curved $(\PP,\mathcal{IC})$-bimodule structure on $\PP$. The forgetful functor $\mathsf{U}$ is naturally isomorphic to $\mathscr{S}_{\PP}(\PP)$ given by the canonical curved $(\mathcal{IC},\PP)$-bimodule structure on $\PP$.
\end{Proposition}

\begin{proof}
Let $(V,d_V)$ be a pdg module, $\mathscr{S}_{\mathcal{IC}}(\PP)(V)$ is defined by the same quotient as $\mathsf{F}(\PP)(V)$: the right action of $\mathcal{IC}$ on $\PP$ is given by $\theta_\PP$ and the left action of $\mathcal{IC}$ on $(V,d_V)$ is given by $d_V^2$, which are identified in $\mathscr{S}_{\mathcal{IC}}(\PP)(V)~.$ It is straightforward to check that $\mathscr{S}_{\PP}(\PP)(V)$ amounts to $(V,d_V)$ endowed with its pdg module structure.
\end{proof}

This indicates that curved $(\PP,\mathcal{IC})$-bimodules are a good analogue to left $\PP$-modules over operads in the case of curved operads; likewise, curved $(\mathcal{IC},\PP)$-bimodules are a good analogue to right $\PP$-modules. 

\begin{Remark}
To the best of our knowledge, these definitions also provide new notions of left and right modules for curved associative algebras when we consider them as curved operads concentrated in arity one. 
\end{Remark}

\section{Curved cooperads and curved partial cooperads}\label{Section: Curved cooperads}
In this section, we briefly recall the notion of a curved cooperad which first appeared in \cite{HirshMilles12}. Conilpotent curved partial cooperads will play the role of the Koszul dual of unital partial operads. The notions of coalgebras and algebras over cooperads also extend to the case of curved cooperads.

\begin{Definition}[Curved cooperad]\label{def curved cooperad}
A \textit{curved cooperad} $(\C,\Delta,\epsilon,d_\C,\Theta_\C)$ amounts to the data of a pdg cooperad $(\C,\Delta,\epsilon,d_\C)$ and a morphism of pdg $\mathbb{S}$-modules $\Theta_\C: (\C,d_\C) \longrightarrow (\I,0)$ of degree $-2$, such that the following diagram commutes: 
\[
\begin{tikzcd}[column sep=7.5pc,row sep=3pc]
\C \arrow[r,"\Delta_{(1)}"] \arrow[rrd,"d_\C^2", bend right =10]
&\C \circ_{(1)} \C \arrow[r,"(\mathrm{id}~ \circ ~ \Theta_\C)~-~(\Theta_\C~ \circ_{(1)}~ \mathrm{id})~"] 
&(\C \circ \I) \oplus (\I \circ \C) \cong  \C \oplus \C  \arrow[d,"\mathrm{proj}"]\\
&
&\C~,
\end{tikzcd}
\]
where $\mathrm{proj}$ is given by $\mathrm{proj}(\mu,\nu) \coloneqq \mu + \nu$. A \textit{morphism} of curved cooperads $f: (\C,,d_\C,\Theta_\C) \longrightarrow (\D,d_\D,\Theta_\D)$ is the data of a morphism of pdg cooperads $f: (\C,,d_\C) \longrightarrow (\D,d_\D)$ such that $\Theta_\D \circ f = \Theta_\C~.$
\end{Definition}

By a slight abuse of notation, we denote this equality by 
\[
d_\C^2 = (\mathrm{id} \circ_{(1)} \Theta_\C - \Theta_\C \circ \mathrm{id}) \cdot \Delta_{(1)}~,
\]
forgetting $\mathrm{proj}$ and the identifications made.

\begin{Remark}
One also defines \textit{curved partial cooperads} and \textit{curved counital partial cooperads} as one would expect, since the condition on the curvature only involves partial decompositions. In this context, the comparison results of Section \ref{Section: partial (co)operads} between different types of cooperads extend without any difficulty to the curved case, \textit{mutatis mutandis}.
\end{Remark}

Any curved cooperad $(\C,d_\C,\Theta_\C)$ induces a comonad structure on its Schur functor $\mathscr{S}(\C)$, and a monad structure on its dual Schur functor $\widehat{\mathscr{S}}^c(\C)~.$ Hence we can define curved coalgebras and curved algebras over any given curved cooperad.

\begin{Definition}[Curved $\C$-coalgebra]\label{def curved conil coalg}
Let $(\C,d_\C,\Theta_\C)$ be a curved cooperad and let $(C,\Delta_C,d_C)$ be a pdg $\C$-coalgebra. It is a \textit{curved} $\C$-\textit{coalgebra} if the following diagram commutes:
\[
\begin{tikzcd}[column sep=3pc,row sep=3pc]
C  \arrow[r,"\Delta_C "] \arrow[rd,"d_C^2",swap]
&\mathscr{S}(\C)(C) \arrow[d,"\mathscr{S}(\Theta_\C)(\mathrm{id})"]\\
&C \cong \mathscr{S}(\I)(C)~.
\end{tikzcd}
\]
A \textit{morphism} of curved conilpotent $\C$-coalgebras $f: (C,\Delta_C,d_C) \longrightarrow (C',\Delta_{C'},d_{C'})$ is the data of a morphism of pdg $\C$-coalgebras. 
\end{Definition}

Thus the category of curved $\C$-coalgebras is a full-subcategory of the category of pdg $\mathcal{C}$-coalgebras. It is not obvious to us whether the category of curved $\C$-coalgebras admits a cofree object or not. It would be interesting to have a comonadicity result like in the case of curved algebras over a curved operad. Nevertheless, the following result guarantees the existence of limits and colimits.

\begin{theorem}[{\cite{grignou2019}}]
The category curved $\C$-coalgebras is a presentable category. In particular, it is both complete and cocomplete.
\end{theorem}

\begin{Definition}[Curved algebra over a cooperad]\label{def curved alg over a coop}
Let $(\C,d_\C,\Theta_\C)$ be a curved cooperad and let $(B,\gamma_B,d_B)$ be a pdg $\C$-algebra. It is a \textit{curved} $\C$\textit{-algebra} if the following diagram commutes: 
\[
\begin{tikzcd}[column sep=4pc,row sep=3pc]
B \cong \widehat{\mathscr{S}}^c(\I)(B) \arrow[r,"\widehat{\mathscr{S}}^c(\Theta_\C)(\mathrm{id}) "] \arrow[rd,"- d_B^2",swap]
&\widehat{\mathscr{S}}^c(\C)(B) \arrow[d,"\gamma_B"]\\
&B ~.
\end{tikzcd}
\]
A \textit{morphism} of curved $\C$-algebras $f: (B,\gamma_B,d_B) \longrightarrow (B',\gamma_{B'},d_{B'})$ is the data of a morphism of pdg $\C$-algebras. 
\end{Definition}

Thus the category of curved $\mathcal{C}$-algebras is also a full sub-category of the category of pdg $\mathcal{C}$-algebras. The following proposition gives a reflector.

\begin{Proposition}[{\cite[Theorem 7.5]{grignoulejay18}}]
Let $(\C,d_\C,\Theta_\C)$ be a curved cooperad. The inclusion functor 
\[
\mathsf{Inc}: \mathsf{curv}~\C\text{-}\mathsf{alg} \hookrightarrow \mathsf{pdg}~\C\text{-}\mathsf{alg}
\]
has a left adjoint 
\[
\mathsf{Curv}: \mathsf{pdg}~\C\text{-}\mathsf{alg} \longrightarrow \mathsf{curv}~\C\text{-}\mathsf{alg}~.
\]
Hence $\mathsf{curv}~\C\text{-}\mathsf{alg}$ is a reflexive subcategory of $\mathsf{pdg}~\C\text{-}\mathsf{alg}$. For a pdg $\C$-algebra $(B,\gamma_B,d_B)$, its image under this functor is given by the following quotient: 

\[
\mathsf{Curv}(B) \coloneqq \frac{B}{\left(\gamma_B \cdot \widehat{\mathscr{S}}^c(\Theta_\C)(\mathrm{id}) + d_B^2(-)\right)}~.
\]
where $\left(\gamma_B \cdot \widehat{\mathscr{S}}^c(\Theta_\C)(\mathrm{id}) + d_B^2(-)\right)$ denoted the ideal generated by $\mathsf{Im}\left(\gamma_B \cdot \widehat{\mathscr{S}}^c(\Theta_\C)(\mathrm{id}) + d_B^2(-)\right)~.$ Its pdg $\C$-algebra structure is induced by $\gamma_B$ and $d_B$.
\end{Proposition}

\begin{Corollary}
The category of curved $\C$-algebras is presentable. In particular it is complete and cocomplete.
\end{Corollary}

The notion of conilpotentcy also generalizes to the case of curved cooperads.

\begin{Definition}[Conilpotent curved partial cooperad]
Let $(\C,\{\Delta_i\},d_\C,\Theta_\C)$ be a curved partial cooperad. It is \textit{conilpotent} if its underlying partial pdg cooperad $(\C,\{\Delta_i\},d_\C)$ is a conilpotent partial cooperad. 
\end{Definition}

There is a functor 
\[
\mathsf{Conil}: \mathsf{curv}~\mathsf{pCoop}^\mathsf{conil} \longrightarrow \mathsf{curv}~\mathsf{Coop}~,
\]
from the category of conilpotent curved partial cooperads to the category of curved cooperads defined as comonoids given by adding a counit, see Definition \ref{def: conilpotent cooperad}.

\begin{Definition}[Conilpotent curved cooperad]
A curved cooperad $(\C,d_\C,\Theta_\C)$ is said to be \textit{conilpotent} if it is in the essential image of the functor $\mathsf{Conil}~.$
\end{Definition}

In the case of a conilpotent curved cooperad, one can induce a canonical filtration on the category of curved $\mathcal{C}$-algebras.

\begin{Definition}[Complete curved algebra over a conilpotent cooperad]
Let $(\C^u,d_\C,\Theta_\C)$ be a conilpotent curved cooperad. A \textit{complete} curved $\C$-algebra $B$ amounts to the data of a curved $\C$-algebra $(B,\gamma_B,d_B)$ such that the canonical morphism of curved $\C$-algebras
\[
\varphi_B: B \longrightarrow \lim_{\omega} B/\mathrm{W}_\omega B
\]
is an isomorphism, where $\mathrm{W}_\omega B$ is the canonical filtration of the $\C$-algebra $B$.
\end{Definition}

\begin{Remark}
These notions where defined in Section \ref{Section: coalgebras and algebras}.
\end{Remark}

\begin{Proposition}\label{bicomplete C alg}
Let $(\C,d_\C,\Theta_\C)$ be a conilpotent curved cooperad. The category of complete curved $\C$-algebras is a reflexive subcategory of pdg $\C$-algebras. It is thus presentable, and bicomplete.
\end{Proposition}

\begin{proof}
Its reflector is given by the composition of the previous two reflectors: the functor $\mathsf{Curv}$ and the completion functor with respect to the canonical filtration.
\end{proof}

\section{Complete Bar-Cobar adjunction in the curved case}\label{Section: Complete Bar-Cobar curved}
Let us fix a dg operad $\mathcal{P}$ and a conilpotent curved cooperad $\mathcal{C}$, together with a curved twisting morphism $\alpha: \mathcal{C} \longrightarrow \mathcal{P}$ in the sense of Definition \ref{def: curved twisting entre P et curved C}. These are the curved generalization of Section \ref{Section: Complete Bar-Cobar}.

\begin{Definition}[Complete Bar construction relative to $\alpha$]\label{def: complete Bar}
Let $(B,\gamma_B,d_B)$ be a curved $\C$-algebra. The \textit{complete Bar construction relative to} $\alpha$ of $B$ is given by
\[
\widehat{\mathrm{B}}_{\alpha} B \coloneqq (\mathscr{C}(\PP)(B), d_{\mathrm{bar}} \coloneqq d_1 + d_2)~,
\]
where $\mathscr{C}(\PP)(B)$ denotes the cofree graded $\PP$-coalgebra generated by $B$. The differential $d_{\mathrm{bar}}$ is given by the sum of two terms $d_1$ and $d_2$. The term $d_1$ is given by
\[
d_1 = \mathscr{C}(d_\PP)(\mathrm{id}) +  \mathscr{C}(\mathrm{id})(\diracComb(id_B,d_B))~.
\]
The term $d_2$ is given by the unique coderivation extending  
\[
\begin{tikzcd}[column sep=4pc,row sep=1pc]
\mathscr{C}(\PP)(B) \arrow[r,"p_1 (B)",rightarrowtail]
&\widehat{\mathscr{S}}^c(\PP)(B)\arrow[r,"\widehat{\mathscr{S}}^c(\alpha)(\mathrm{id}_B)"]
&\widehat{\mathscr{S}}^c(\C)(B)  \arrow[r,"\gamma_B "]
&B~.
\end{tikzcd}
\]
\end{Definition}

\begin{Proposition}
For any curved $\C$-algebra $B$, the complete Bar construction $\widehat{\mathrm{B}}_{\alpha}B$ forms a dg $\PP$-coalgebra, and it defines a functor
\[
\widehat{\mathrm{B}}_{\alpha}: \mathsf{curv}~\C\text{-}\mathsf{alg} \longrightarrow \mathsf{dg}~\PP\text{-}\mathsf{coalg}~.
\]
\end{Proposition}

\begin{Definition}[Complete Cobar construction relation to $\alpha$]\label{def: complete Cobar}
Let $(D,\delta_D,d_D)$ be a dg $\PP$-coalgebra. The \textit{complete Cobar construction relative to} $\alpha$ of $D$ is given by
\[
\widehat{\Omega}_\alpha D \coloneqq (\widehat{\mathscr{S}}^c(\C)(D), d_{\mathrm{cobar}} \coloneqq d_1 - d_2)~,
\]
where $\widehat{\mathscr{S}}^c(\C)(D)$ denotes the free complete pdg $\C$-coalgebra generated by $D$. The differential $d_{\mathrm{cobar}}$ is given by the difference of two terms $d_1$ and $d_2$. The term $d_1$ is given by
\[
d_1 = -\widehat{\mathscr{S}}^c(d_\C)(\mathrm{id}) + \widehat{\mathscr{S}}^c(\mathrm{id})(\diracComb(\mathrm{id}_D,d_D))~.
\]
The term $d_2$ is given by the unique derivation extending  
\[
\begin{tikzcd}[column sep=4pc,row sep=1pc]
D \arrow[r," \Delta_D"]
&\widehat{\mathscr{S}}^c(\PP)(D) \arrow[r,"\widehat{\mathscr{S}}^c(\alpha)(\mathrm{id}_D )"]
&\widehat{\mathscr{S}}^c(\C)(D)~.
\end{tikzcd} 
\]
\end{Definition}

\begin{Proposition}
For any dg $\PP$-coalgebra $D$, the complete Cobar construction $\widehat{\Omega}_{\alpha}D$ forms a complete dg curved $\C$-algebra, and it defines a functor
\[
\widehat{\Omega}_\alpha: \mathsf{dg}~\PP\text{-}\mathsf{coalg} \longrightarrow \mathsf{curv}~\C\text{-}\mathsf{alg}^{\mathsf{comp}}~.
\]
\end{Proposition}

\begin{Definition}[Curved twisting morphism relative to $\alpha$]
Let $(D,\delta_D,d_D)$ be a dg $\PP$-coalgebra and let $(B,\gamma_B,d_B)$ be a curved $\C$-algebra. A graded morphism 
\[
\nu: D \longrightarrow B
\]
is said to be a curved twisting morphism relative to $\alpha$ if it satisfies the following equation 

\[
\gamma_B \cdot \widehat{\mathscr{S}}^c(\alpha)(\nu) \cdot \Delta_D + \partial(\nu) = 0~.
\]

The set of curved twisting morphisms relative to $\alpha$ are denoted by $\mathrm{Tw}^{\alpha}(D,B)$.
\end{Definition}

\begin{Proposition}
There are bijections

\[
\mathrm{Hom}_{\mathsf{curv}~\C\text{-}\mathsf{alg}^{\mathsf{comp}}}\left(\widehat{\Omega}_{\alpha}D, A\right) \cong  \mathrm{Tw}^{\alpha}(D,B)  \cong \mathrm{Hom}_{\mathsf{dg}~\PP\text{-}\mathsf{coalg}} \left(D, \widehat{\mathrm{B}}_{\alpha}B \right)~,
\]

which are natural in $D$ and $B$. 
\end{Proposition}

Therefore we get an adjunction 
\[
\begin{tikzcd}[column sep=5pc,row sep=3pc]
            \mathsf{dg}~\PP\text{-}\mathsf{coalg} \arrow[r, shift left=1.1ex, "\widehat{\Omega}_{\alpha}"{name=F}] & \mathsf{curv}~\C\text{-}\mathsf{alg}^{\mathsf{comp}} \arrow[l, shift left=.75ex, "\widehat{\mathrm{B}}_{\alpha}"{name=U}]
            \arrow[phantom, from=F, to=U, , "\dashv" rotate=-90]
\end{tikzcd}
\]
between the category of dg $\PP$-coalgebras and the category of complete curved $\C$-algebras. This adjunction is called the \textit{complete Bar-Cobar adjunction} relative to $\alpha$. See \cite[Section 8]{grignoulejay18} for more details.
 
\medskip

In some case, one can promote this adjunction into a Quillen adjunction. But first, one needs a model structure on the category of dg $\PP$-coalgebras. 

\begin{theorem}[{\cite[Section 9]{grignoulejay18}}]\label{thm: model structure on P-cog}
Let $\mathcal{P}$ be a cofibrant dg operad. There is a model structure on the category of dg $\mathcal{P}$-coalgebras left-transferred along the cofree-forgetful adjunction
\[
\begin{tikzcd}[column sep=7pc,row sep=3pc]
\mathsf{dg}\textsf{-}\mathsf{mod} \arrow[r, shift left=1.1ex, "\mathcal{L}(\mathcal{P})(-)"{name=F}]      
&\mathsf{dg}~\PP\textsf{-}\mathsf{coalg}~, \arrow[l, shift left=.75ex, "U"{name=U}]
\arrow[phantom, from=F, to=U, , "\dashv" rotate=90]
\end{tikzcd}
\]
where 

\medskip
\begin{enumerate}
\item the class of weak equivalences is given by quasi-isomorphisms,

\medskip

\item the class of cofibrations is given by degree-wise monomorphisms,

\medskip

\item the class of fibrations is given by right lifting property with respect to acyclic cofibrations.
\end{enumerate}
\end{theorem}

\begin{Remark}
The assumption that $\PP$ is a cofibrant dg operad is necessary. Indeed, one can show that for $\PP = u\mathcal{C}om$, the category of dg $ u\mathcal{C}om$-coalgebras does not admit a model structure where weak-equivalences are given by quasi-isomorphisms and where cofibrations are given by degree-wise monomorphisms. Therefore, while all dg operads are \textit{admissible} in characteristic zero, not all dg operad are \textit{coadmissible}.
\end{Remark}

Model structures on coalgebras over dg operads behave well with respect to quasi-isomorphisms at the operadic level.

\begin{theorem}[{\cite[Section 9]{grignoulejay18}}]
Let $f: \mathcal{P} \qi \mathcal{Q}$ be a quasi-isomorphism of cofibrant dg operads. The induced adjunction 

\[
\begin{tikzcd}[column sep=7pc,row sep=3pc]
\mathsf{dg}~\mathcal{P}\textsf{-}\mathsf{coalg} \arrow[r, shift left=1.1ex, "\mathrm{Coind}_{f}"{name=F}]      
&\mathsf{dg}~\mathcal{Q}\textsf{-}\mathsf{coalg}~, \arrow[l, shift left=.75ex, "\mathrm{Res}_{f}"{name=U}]
\arrow[phantom, from=F, to=U, , "\dashv" rotate=90]
\end{tikzcd}
\]

is a Quillen equivalence.
\end{theorem}

When this model structure exists, it can be transferred along the complete Bar-Cobar adjunction relative to $\alpha$. 

\begin{theorem}[{\cite[Section 10]{grignoulejay18}}]
Let $\mathcal{P}$ be a cofibrant dg operad. There is a model structure on the category of complete $\mathcal{C}$-algebras right-transferred along the complete Bar-Cobar adjunction
\[
\begin{tikzcd}[column sep=7pc,row sep=3pc]
            \mathsf{dg}~\mathcal{P}\textsf{-}\mathsf{coalg} \arrow[r, shift left=1.1ex, "\widehat{\Omega}_{\alpha}"{name=F}] &\mathsf{curv}~\mathcal{C}\textsf{-}\mathsf{alg}^{\mathsf{comp}}~, \arrow[l, shift left=.75ex, "\widehat{\text{B}}_{\alpha}"{name=U}]
            \arrow[phantom, from=F, to=U, , "\dashv" rotate=-90]
\end{tikzcd}
\]
where 

\medskip 

\begin{enumerate}
\item the class of weak equivalences is given by morphisms $f$ such that $\widehat{\text{B}}_\alpha(f)$ is a quasi-isomorphism,

\medskip

\item the class of fibrations is given by degree-wise epimorphisms,

\medskip

\item  and the class of cofibrations is given by left lifting property with respect to acyclic fibrations.
\end{enumerate}
\end{theorem}

In the case where $\PP$ is the Cobar construction of $\mathcal{C}$ in the sense of \cite{grignou2021}, this adjunction can be promoted to a Quillen equivalence.

\begin{theorem}[{\cite[Section 11]{grignoulejay18}}]
The complete Bar-Cobar adjunction relative to $\iota: \C \longrightarrow \Omega \C$
\[
\begin{tikzcd}[column sep=7pc,row sep=3pc]
            \mathsf{dg}~\Omega\mathcal{C}\textsf{-}\mathsf{coalg} \arrow[r, shift left=1.1ex, "\widehat{\Omega}_{\iota}"{name=F}] 
            &\mathsf{curv}~\mathcal{C}\textsf{-}\mathsf{alg}^{\mathsf{comp}}~. \arrow[l, shift left=.75ex, "\widehat{\text{B}}_{\iota}"{name=U}]
            \arrow[phantom, from=F, to=U, , "\dashv" rotate=-90]
\end{tikzcd}
\]
is a Quillen equivalence. 
\end{theorem}

\begin{Remark}
The model structure on the category of complete curved $\mathcal{C}$-algebras transferred using the complete Bar-Cobar adjunction relative to the curved twisting morphism $\iota: \mathcal{C} \longrightarrow \Omega \mathcal{C}$ is called the \textit{canonical model structure}.
\end{Remark}

\section{Curved duality squares}\label{Section: Curved duality square}
Let us fix a dg operad $\mathcal{P}$ and a conilpotent curved cooperad $\mathcal{C}$, together with a curved twisting morphism $\alpha: \mathcal{C} \longrightarrow \mathcal{P}$. The goal of this section is to construct two duality adjunctions that interrelate the classical Bar-Cobar constructions relative to $\alpha$ with the complete Bar-Cobar constructions relative to $\alpha$. This a generalization of the results of Section \ref{Section: Magic squares}.

\subsection{Sweedler functor and topological dual functors}
Recall that there was an adjunction

\[
\begin{tikzcd}[column sep=7pc,row sep=3pc]
            \mathsf{dg}~\PP\text{-}\mathsf{coalg} \arrow[r, shift left=1.1ex, "(-)^*"{name=F}] &\left(\mathsf{dg}~\PP \text{-}\mathsf{alg}\right)^{\mathsf{op}} ~. \arrow[l, shift left=.75ex, "(-)^\circ"{name=U}]
            \arrow[phantom, from=F, to=U, , "\dashv" rotate=-90]
\end{tikzcd}
\]

between the category of dg $\PP$-algebras and the category of dg $\PP$-coalgebras, given by the linear dual functor and the Sweedler dual functor. Let's turn to the other side of the Koszul duality, when $\mathcal{C}$ is a conilpotent \textit{curved} cooperad.

\begin{lemma}
The linear dual defines a functor
\[
\begin{tikzcd}[column sep=4pc,row sep=0pc]
\left(\mathsf{pdg}~\mathcal{C}\text{-}\mathsf{coalg}\right)^{\mathsf{op}} \arrow[r,"(-)^*"] 
&\mathsf{pdg}~\mathcal{C}\text{-}\mathsf{alg}^{\mathsf{comp}}
\end{tikzcd}
\]
from the category of pdg $\mathcal{C}$-coalgebras to the category of complete pdg $\mathcal{C}$-algebras.
\end{lemma}

\begin{proof}
Analogue to the non-curved case.
\end{proof}

\begin{Proposition}\label{prop: adjoint à gauche dual topo curved}
The linear dual defines a functor
\[
\begin{tikzcd}[column sep=4pc,row sep=0pc]
\left(\mathsf{pdg}~\mathcal{C}\text{-}\mathsf{coalg}\right)^{\mathsf{op}} \arrow[r,"(-)^*"] 
&\mathsf{pdg}~\mathcal{C}\text{-}\mathsf{alg}^{\mathsf{comp}}
\end{tikzcd}
\]
admits a left adjoint. 
\end{Proposition}

\begin{proof}
Analogue to the non-curved case. 
\end{proof}

\begin{Definition}[Topological dual functor]\label{def: topological dual functor curved}
The \textit{topological dual functor} 
\[
\begin{tikzcd}[column sep=4pc,row sep=0pc]
\mathsf{pdg}~\mathcal{C}\text{-}\mathsf{alg}^{\mathsf{comp}} \arrow[r,"(-)^\vee"]
&\left(\mathsf{pdg}~\mathcal{C}\text{-}\mathsf{coalg}\right)^{\mathsf{op}}
\end{tikzcd}
\]
is defined as the left adjoint of the linear dual functor.
\end{Definition}

\begin{Remark}
Given a complete pdg $\mathcal{C}$-algebra $(B, \gamma_B, d_B)$, we have that $B^\vee$ is given by the following equalizer:
\[
\begin{tikzcd}[column sep=4pc,row sep=4pc]
\mathrm{Eq}\Bigg(\mathscr{S}(\mathcal{C})(B^*) \arrow[r,"(\gamma_B)^*",shift right=1.1ex,swap]  \arrow[r,"\varrho"{name=SD},shift left=1.1ex ]
&\mathscr{S}(\mathcal{C})\left((\widehat{\mathscr{S}}^c(\mathcal{C})(B))^*\right) \Bigg)~,
\end{tikzcd}
\]
where $\varrho$ is an arrow constructed using the comonadic structure of $\mathscr{S}(\mathcal{C})(-)$ and the canonical inclusion of a dg module into its double linear dual.
\end{Remark}

\begin{Proposition}\label{prop: natural mono for topo dual curved}
There is a natural monomorphism
\[
\epsilon: \mathrm{U}^{\mathsf{op}} \cdot (-)^\vee \rightarrowtail (-)^* \cdot \mathrm{U}~,
\]
which implies that the topological dual is a sub-pdg module of the linear dual functor. 
\end{Proposition}

\begin{proof}
Analogue to the non-curved case.
\end{proof}

\begin{Proposition}\label{prop: restriction aux curved}
The adjunction 
\[
\begin{tikzcd}[column sep=7pc,row sep=3pc]
           \mathsf{pdg}~\mathcal{C}\text{-}\mathsf{alg}^{\mathsf{comp}} \arrow[r,"(-)^\vee "{name=F}, shift left=1.1ex] 
           &\left(\mathsf{pdg}~\mathcal{C}\text{-}\mathsf{coalg}\right)^{\mathsf{op}}~, \arrow[l, shift left=.75ex, "(-)^*"{name=U}]
            \arrow[phantom, from=F, to=U, , "\dashv" rotate=-90]
\end{tikzcd}
\]
restricts to an adjunction 
\[
\begin{tikzcd}[column sep=7pc,row sep=3pc]
           \mathsf{curv}~\mathcal{C}\text{-}\mathsf{alg}^{\mathsf{comp}} \arrow[r,"(-)^\vee "{name=F}, shift left=1.1ex] 
           &\left(\mathsf{curv}~\mathcal{C}\text{-}\mathsf{coalg}\right)^{\mathsf{op}}~. \arrow[l, shift left=.75ex, "(-)^*"{name=U}]
            \arrow[phantom, from=F, to=U, , "\dashv" rotate=-90]
\end{tikzcd}
\]
\end{Proposition}

\begin{proof}
Denote $\Theta_\mathcal{C}: \C \longrightarrow \I$ the curvature of $\C$. Let $(C,\Delta_C,d_C)$ be a pdg $\mathcal{C}$-coalgebra. Recall that it is curved if the following diagram commutes
\[
\begin{tikzcd}[column sep=3pc,row sep=3pc]
C \arrow[r,"\Delta_C"] \arrow[rd,"d_C^2",swap] 
&\mathscr{S}(\mathcal{C})(C) \arrow[d,"\mathscr{S}(\Theta_\mathcal{C})(\mathrm{id}_C)"] \\
&C \cong \mathscr{S}(\I)(C)~.
\end{tikzcd}
\]
Therefore we have 
\[
\begin{tikzcd}[column sep=3pc,row sep=3pc]
\widehat{\mathscr{S}}^c(\mathcal{C})(C^*) \arrow[r,rightarrowtail] 
&\left(\mathscr{S}(\mathcal{C})(C)\right)^* \arrow[r,"(\Delta_C)^* "] 
&C^* \\
\widehat{\mathscr{S}}^c(\I)(C^*) \arrow[r,"\cong"] \arrow[u,"\widehat{\mathscr{S}}^c(\Theta_\mathcal{C})(\mathrm{id}_C)"]
&\left(\mathscr{S}(\I)(C)\right)^* \cong C^*~. \arrow[ru,"(d_C^2)^*",swap] \arrow[u,"(\mathscr{S}(\Theta_\mathcal{C})(\mathrm{id}_C))^*"]
\end{tikzcd}
\]
The right square commutes by naturality of the inclusion. The left triangle commutes since its the image of a commutative triangle by the functor $(-)^*$. Thus the big diagram commutes. Finally, one can observe that $(d_C^2)^* = - d_{C^*}^2$. Therefore $C^*$ is indeed a complete curved $\mathcal{C}$-algebra. 

\medskip 

Let $(B, \gamma_B, d_B)$ be a complete curved $\mathcal{C}$-algebra. Consider the following diagram
\[
\begin{tikzcd}[column sep=3pc,row sep=3pc]
\left(\widehat{\mathscr{S}}^c(\mathcal{C})(B)\right)^*
&
&B^* \arrow[ll,"(\gamma_B)^* ",swap] \\
&\left(\widehat{\mathscr{S}}^c(\I)(B)\right)^* \cong B^* \arrow[lu,"\left(\widehat{\mathscr{S}}^c(\Theta_\mathcal{C})(\mathrm{id})\right)^* " {xshift= 0.35pc}]\arrow[ru,"(-d_B^2)^* "]
& \\
\mathscr{S}(\C)(B^\vee) \arrow[uu,rightarrowtail] 
&
&B^\vee  \arrow[ll, "\Delta_{B^\vee}" {xshift= -0.7pc},swap] \arrow[uu,rightarrowtail]  \\
&\mathscr{S}(\I)(B^\vee) \cong B^\vee~, \arrow[ru,"d_{B^\vee}^2 "] \arrow[lu,"\mathscr{S}(\Theta_\mathcal{C})(\mathrm{id}_{B^\vee})"] \arrow[uu,rightarrowtail, crossing over]
\end{tikzcd}
\]
in the category of pdg modules, where the vertical arrows are given by Proposition \ref{prop: natural mono for topo dual}. The top triangle commutes since it is the image of a commuting triangle via the functor $(-)^*$. Each of the vertical faces commutes as well, where $d_{B^\vee}^2$ is simply given by $(d_B)^\vee \circ (d_B)^\vee$. Since every vertical map is a monomorphism, the bottom triangle also commutes and thus $B^\vee$ is indeed a curved $\mathcal{C}$-coalgebra.
\end{proof}

\begin{Remark}
As far as we know, the category of curved $\mathcal{C}$-coalgebras might not be comonadic. Thus in order to prove Proposition \ref{prop: restriction aux curved}, one cannot use the same arguments as in the proof of Proposition \ref{prop: adjoint à gauche dual topo}.
\end{Remark}

\begin{Remark}
Given a complete curved $\mathcal{C}$-algebra $B$ which is degree-wise finite dimensional, it is not clear to us whether its linear dual $B^*$ admits a curved $\mathcal{C}$-coalgebra structure or not. Furthermore, determining when one has that $B^\vee \cong B^*$ do not seem to be trivial. Should the interested reader find something in this direction, we would very much appreciate it.
\end{Remark}

\begin{Remark}
The results of this subsection are still valid if one does not suppose that $\mathcal{C}$ is conilpotent, and consequently forgets about the "completeness" of $\mathcal{C}$-algebras.
\end{Remark}

\subsection{The algebraic duality square} Let $\alpha: \mathcal{C} \longrightarrow \mathcal{P}$ be a curved twisting morphism. Using the two adjunctions constructed so far, one can interrelate the classical Bar-Cobar adjunction relative to $\alpha$ with the complete Bar-Cobar adjunction relative to $\alpha$ is a square of commuting adjunctions.

\begin{theorem}[Duality square]\label{thm: magical square curved}
The square of adjunctions 
\[
\begin{tikzcd}[column sep=5pc,row sep=5pc]
\left(\mathsf{dg}~\mathcal{P}\text{-}\mathsf{alg}\right)^{\mathsf{op}} \arrow[r,"\mathrm{B}_\alpha^{\mathsf{op}}"{name=B},shift left=1.1ex] \arrow[d,"(-)^\circ "{name=SD},shift left=1.1ex ]
&\left(\mathsf{curv}~\mathcal{C}\text{-}\mathsf{coalg}\right)^{\mathsf{op}} \arrow[d,"(-)^*"{name=LDC},shift left=1.1ex ] \arrow[l,"\Omega_\alpha^{\mathsf{op}}"{name=C},,shift left=1.1ex]  \\
\mathsf{dg}~\mathcal{P}\text{-}\mathsf{coalg} \arrow[r,"\widehat{\Omega}_\alpha "{name=CC},shift left=1.1ex]  \arrow[u,"(-)^*"{name=LD},shift left=1.1ex ]
&\mathsf{curv}~\mathcal{C}\text{-}\mathsf{alg}^{\mathsf{comp}}~, \arrow[l,"\widehat{\mathrm{B}}_\alpha"{name=CB},shift left=1.1ex] \arrow[u,"(-)^\vee"{name=TD},shift left=1.1ex] \arrow[phantom, from=SD, to=LD, , "\dashv" rotate=0] \arrow[phantom, from=C, to=B, , "\dashv" rotate=-90]\arrow[phantom, from=TD, to=LDC, , "\dashv" rotate=0] \arrow[phantom, from=CC, to=CB, , "\dashv" rotate=-90]
\end{tikzcd}
\] 
commutes in the following sense: right adjoints going from the top right to the bottom left are naturally isomorphic.
\end{theorem}

\begin{proof}
Analogue to the non-curved case.
\end{proof}

\begin{Remark}
Let $(D,\Delta_D,d_D)$ be a dg $\PP$-coalgebra. Then we have an isomorphism
\[
\mathrm{B}_\alpha(D^*) \cong \left(\widehat{\Omega}_\alpha(D) \right)^\vee
\]
of curved $\mathcal{C}$-coalgebras which is natural in $D$, which is given by the mate of the isomorphism constructed in the proof of Theorem \ref{thm: magical square}. 
\end{Remark}

\begin{Proposition}
There is a natural monomorphism
\[
\zeta: \widehat{\Omega}_\alpha \cdot (-)^\circ \longrightarrow  (-)^* \cdot \mathrm{B}_\alpha~.
\]
of complete curved $\C$-algebras.
\end{Proposition}

\begin{proof}
This monomorphism is build using the monomorphisms constructed in Proposition \ref{prop: natural mono for Sweedler dual} and Proposition \ref{prop: natural mono for topo dual}.
\end{proof}

\begin{Proposition}\label{prop: finite dual commutes curved}
Let $(A,\gamma_A, d_A)$ be a dg $\PP$-algebra degree-wise finite dimensional. There is an isomorphism 
\[
\widehat{\Omega}_\alpha(A^*) \cong \left(\mathrm{B}_\alpha(A)\right)^*
\]
of complete curved $\mathcal{C}$-algebras. 
\end{Proposition}

\begin{proof}
Analogue to the non-curved case.
\end{proof}

\begin{Remark}
The subcategory of dg $\PP$-algebras which satisfy the \textit{Beck-Chevalley condition} with respect to the duality square of Theorem \ref{thm: magical square} is exactly given by the subcategory of degree-wise finite dimensional dg $\PP$-algebras.
\end{Remark}

\subsection{Homotopical properties of the duality square} In this subsection, we show that the duality square of Theorem \ref{thm: magical square} behaves well with respect to model structures. We now restrict to the case where $\PP = \Omega \C$ and where the curved twisting morphism considered is the canonical twisting morphism $\iota: \C \longrightarrow \Omega \C$, in order to ensure the existence of a model category structure on the category of dg $\PP$-coalgebras. 

\begin{theorem}[Homotopical properties of the duality square]\label{thm: homotopical magical square curved}
All the adjunctions in the square 
\[
\begin{tikzcd}[column sep=5pc,row sep=5pc]
\left(\mathsf{dg}~\Omega \C\text{-}\mathsf{alg}\right)^{\mathsf{op}} \arrow[r,"\mathrm{B}_\iota^{\mathsf{op}}"{name=B},shift left=1.1ex] \arrow[d,"(-)^\circ "{name=SD},shift left=1.1ex ]
&\left(\mathsf{curv}~\mathcal{C}\text{-}\mathsf{coalg}\right)^{\mathsf{op}} \arrow[d,"(-)^*"{name=LDC},shift left=1.1ex ] \arrow[l,"\Omega_\iota^{\mathsf{op}}"{name=C},,shift left=1.1ex]  \\
\mathsf{dg}~\Omega \C\text{-}\mathsf{coalg} \arrow[r,"\widehat{\Omega}_\iota "{name=CC},shift left=1.1ex]  \arrow[u,"(-)^*"{name=LD},shift left=1.1ex ]
&\mathsf{curv}~\mathcal{C}\text{-}\mathsf{alg}^{\mathsf{comp}}~, \arrow[l,"\widehat{\mathrm{B}}_\iota"{name=CB},shift left=1.1ex] \arrow[u,"(-)^\vee"{name=TD},shift left=1.1ex] \arrow[phantom, from=SD, to=LD, , "\dashv" rotate=0] \arrow[phantom, from=C, to=B, , "\dashv" rotate=-90]\arrow[phantom, from=TD, to=LDC, , "\dashv" rotate=0] \arrow[phantom, from=CC, to=CB, , "\dashv" rotate=-90]
\end{tikzcd}
\] 
are Quillen adjunctions.
\end{theorem}

\begin{proof}
Analogue to the non-curved case.
\end{proof}

\begin{Proposition}\label{Prop: plongement pleinement fidele de Sweedler}
Let $(A,\gamma_A, d_A)$ be a dg $\Omega \C$-algebra whose homology is degree-wise finite dimensional. Then the derived unit of adjunction

\[
\mathbb{R}(\eta_A): A \qi \left(\mathbb{R}(A^\circ)\right)^*
\]
\vspace{0.2pc}

is a quasi-isomorphism. Therefore the functor $(-)^\circ$ is homotopically fully faithful on the full sub-$\infty$-category of dg $\mathcal{P}$-algebras with finite dimensional homology.
\end{Proposition}

\begin{proof}
Analogue to the non-curved case.
\end{proof}

\section{The groupoid-colored level}
In this section, we develop the formalism of groupoid-colored (co)operads introduced in \cite{ward19} in order to encode the different (co)operadic structures that we have encountered so far. More precisely, we construct a unital groupoid-colored operad, denoted by $u\mathcal{O}$, that encodes unital partial operads as its algebras and counital partial cooperads as its coalgebras. The goal of this section is to compute its Koszul dual conilpotent curved groupoid-colored cooperad.

\medskip

In order to do this, we need to extend to groupoid-colored operads the inhomogeneous Koszul duality for operads introduced in \cite{HirshMilles12} (by inhomogeneous we mean relations that involve constant-linear-quadratic terms in this case). Extending the quadratic Koszul duality for groupoid-colored operads of \cite{ward19} to the inhomogeneous case is conceptually quite straightforward. Since this theory will only be applied to the case that interests us, we give an overview of the main results, without entering in full generality.
\[
\begin{tikzpicture} 
\tikzstyle{block} = [draw, rectangle, text width = 12em, 
  text centered, minimum height = 7mm, node distance = 7em];
\tikzstyle{line} = [draw, -stealth, thick]
\node [block] at (0, 0) (one) {\scriptsize Koszul duality \\ homogeneous operads \cite{GinzburgKapranov95},\cite{GetzlerJones94}};
\node [block] at (0, 2) (two) {\scriptsize Koszul duality \\ inhomogeneous operads};
\node [block] at (9, 0) (three) {\scriptsize Koszul duality \\ groupoid-colored homogeneous operads};
\node [block] at (9, 2) (four) {\scriptsize Koszul duality \\ groupoid-colored inhomogeneous operads};
\path [line] (one) -- (two);
\path [line] (two) -- (four);
\path [line] (three) -- (four);
\path [line] (one) -- (three);
\node at (-0.75, 1) [align = left]{\scriptsize \cite{HirshMilles12}};
\node at (4.5, 0.3) [align = center]{\scriptsize \cite{ward19}};
\end{tikzpicture}
\]
This new Koszul duality gives us a conilpotent curved groupoid-colored cooperad $(u\mathcal{O})^{\ac}$, which encodes conilpotent curved partial cooperads as its curved coalgebras. Its category of curved algebras provides us with a new notion of operads which we call \textit{curved absolute partial operads}. See the Appendix \ref{Appendix B} for a detailed description. Using the Koszul curved twisting morphism given by this theory, we will construct two Bar-Cobar adjunctions that interrelate these objects in the next section.

\medskip

\subsection{$\mathbb{S}$-colored (co)operads}
We begin by applying the formalism of \cite{ward19} to the specific case that interests us. We consider the groupoid 
\[
\mathbb{S} \coloneqq \coprod_{n \geq 0} \mathrm{B}\mathbb{S}_n
\]
where $\mathrm{B}\mathbb{S}_n$ is the classifying space of $\mathbb{S}_n$ (the category with one object $\{*\}$ and $\mathrm{Aut}(*) = \mathbb{S}_n$). We denote by $n$ the single object in $\mathrm{B}\mathbb{S}_n~.$ In this context, an $\mathbb{S}$-color scheme as defined in \cite[Definition 2.2.1]{ward19} amounts to the following definition:

\begin{Definition}[pdg $\mathbb{S}$-color scheme]
An \textit{pdg} $\mathbb{S}$-\textit{color scheme} $E$ amounts to the data of a family $\{E(n_1,\cdots,n_r;n)\}$ of pdg modules for all $r$-tuples of natural numbers $(n_1,\cdots,n_r)$ in $\mathbb{N}^{r}$ and all $n$ in $\mathbb{N}$. Each pdg module $E(n_1,\cdots,n_r;n)$ comes equipped with an action of $\mathbb{S}_{n_1} \times \cdots \times \mathbb{S}_{n_r}$ on the right and an action of $\mathbb{S}_n$ on the left as part of the structure. There is a supplementary action of the symmetric groups on this family given by permuting the entries. For any $r$-tuples of natural numbers $(n_1,\cdots,n_r)$, there is an isomorphism
\[
\varphi_\sigma : E(n_1,\cdots,n_r;n) \longrightarrow  E(n_{\sigma^{-1}(1)},\cdots, n_{\sigma^{-1}(r)} ; n) 
\]
of right pdg $\kk[\mathbb{S}_{n_1} \times \cdots \times \mathbb{S}_{n_r}] \otimes \kk[\mathbb{S}_n]^{\text{op}}$-module for all $\sigma$ in $\mathbb{S}_r$ which defines a group action. This endows each $E(n_1,\cdots,n_r;n)$ with an left action of $\left(\mathbb{S}_{n_1} \times \cdots \times \mathbb{S}_{n_r}\right)~\wr~ \mathbb{S}_r$, where $\wr$ denotes the wreath product of groups.

\medskip

A \textit{morphism} $f: E \longrightarrow F$ is the data of a morphism of pdg $\kk[\mathbb{S}_{n_1} \times \cdots \times \mathbb{S}_{n_r}] \otimes \kk[\mathbb{S}_n]^{\text{op}}$-modules 
\[
f(n_1,\cdots,n_r;n): E(n_1,\cdots,n_r;n) \longrightarrow F(n_1,\cdots,n_r;n)~,
\]
for all $r$-tuples $(n_1,\cdots,n_r)$ in $\mathbb{N}^{r}$ which commutes with the permutation maps $\varphi_\sigma$. Pdg $\mathbb{S}$-color schemes form a category denoted by $\mathsf{pdg}~\mathsf{Col}_\mathbb{S}~.$
\end{Definition}

Let $E$ and $F$ be two pdg $\mathbb{S}$-color schemes, their \textit{composition product} $E \circ F$ is given by: 
\begin{align*}
&E \circ F (n_1, \cdots, n_r; n) \coloneqq \bigoplus_{j \geq 1}~ \bigoplus_{(a_1,\cdots ,a_j) \in \mathbb{N}^{j}} \Bigg( E(a_1, \cdots, a_j ;n) \otimes_{(\mathbb{S}_{a_1} \times \cdots \times \mathbb{S}_{a_j})~\wr ~\mathbb{S}_j } \\
&  \Bigg( \bigoplus_{i_1 + \cdots + i_j = r} \mathsf{Ind}_{\mathbb{S}_{i_1} \times \cdots \times \mathbb{S}_{i_j}}^{\mathbb{S}_r} \Big( F(n_1, \cdots ,n_{i_1} ; a_1) \otimes \cdots \otimes F(n_{r-i_j}, \cdots , n_r ; a_j) \Big) \Bigg)\Bigg)_{\mathbb{S}_j} ~.
\end{align*}
where the second sum runs over all $j$-tuples $(a_1,\cdots ,a_j)$ and the third sum runs over all $j$-tuples $(i_1,\cdots ,i_j)$ such that $i_1 + \cdots + i_j = r$. It is the analogue of the composition product $\circ$ of $\mathbb{S}$-modules. In this case, leaves are colored by elements of the symmetric groups, hence the extra requirement that the colors have to match. The unit of this product is given by the $\mathbb{S}$-color scheme $\I_{\mathbb{S}}$, defined by $\I_{\mathbb{S}}(n;n) \coloneqq \mathbb{K}[\mathbb{S}_n]$ as a bimodule over itself and $0$ elsewhere, endowed with the trivial pre-differential. 

\begin{lemma}
The data $(\mathsf{pdg}~\mathsf{Col}_\mathbb{S},\circ, \I_{\mathbb{S}})$ forms a monoidal category. 
\end{lemma}

\begin{proof}
A straightforward computation analogous to the standard case.
\end{proof}

\begin{Definition}[$\mathbb{S}$-colored operad]
A \textit{pdg} $\mathbb{S}$-\textit{colored operad} $\mathcal{G}$ is the data of a monoid $(\mathcal{G}, \gamma_\mathcal{G}, \eta, d_\mathcal{G})$ in the monoidal category $(\mathsf{pdg}~\mathsf{Col}_\mathbb{S},\circ, \I_{\mathbb{S}})$. 
\end{Definition}

\begin{Example}[$\mathbb{S}$-colored endomorphism operad] 
Let $M$ be a pdg $\mathbb{S}$-module. The collection
\[
\mathrm{End}_M(n_1,\cdots,n_r;n) \coloneqq \mathrm{Hom}(M(n_1) \otimes \cdots \otimes M(n_r), M(n)) \quad \text{for} \quad (n_1,\cdots,n_r) \in \mathbb{N}^{r}
\]
has a natural structure of pdg $\mathbb{S}$-color scheme where the right pdg $\kk[\mathbb{S}_{n_1} \times \cdots \times \mathbb{S}_{n_r}] \otimes \kk[\mathbb{S}_n]^{\text{op}}$-module structure comes from the $\mathbb{S}$-module structure of $M$ and where the permutation morphism $\varphi_\sigma$ are given by the natural action of the symmetric group $\mathbb{S}_r$ on the tensor product $M(n_1) \otimes \cdots \otimes M(n_r)$. The composition of morphisms and the identity $\mathrm{id}_n$ in $\mathrm{End}_M(n;n)$ endow this pdg $\mathbb{S}$-color scheme with an $\mathbb{S}$-colored operad structure, called the $\mathbb{S}$-colored \textit{endomorphism operad} of $M$.
\end{Example}

\medskip

\begin{Example}[$\mathbb{S}$-colored coendomorphism operad] 
Let $M$ be a pdg $\mathbb{S}$-module. The collection
\[
\mathrm{Coend}_M(n_1,\cdots,n_r;n) \coloneqq \mathrm{Hom}(M(n), M(n_1) \otimes \cdots \otimes M(n_r)) \quad \text{for} \quad (n_1,\cdots,n_r) \in \mathbb{N}^{r}
\]
has a natural structure of pdg $\mathbb{S}$-color scheme as in the above example. The composition of morphisms and the identity $\mathrm{id}_n$ in $\mathrm{End}_M(n;n)$ endow this pdg $\mathbb{S}$-color scheme with an $\mathbb{S}$-colored operad structure, called the $\mathbb{S}$-colored \textit{coendomorphism operad} of $M$.
\end{Example}

\begin{Definition}[$\mathbb{S}$-colored cooperad]
A \textit{pdg} $\mathbb{S}$-\textit{colored cooperad} $\mathcal{V}$ is the data of a comonoid $(\mathcal{V}, \Delta_\mathcal{V}, \epsilon, d_\mathcal{V})$ in the monoidal category $(\mathsf{Col}_\mathbb{S},\circ, \I_{\mathbb{S}})$. 
\end{Definition}

To any pdg $\mathbb{S}$-color scheme $E$ one can associate an endofunctor in the category of pdg $\mathbb{S}$-modules via the \textit{Schur realization functor}:
\[
\begin{tikzcd}[column sep=4pc,row sep=0.5pc]
\mathscr{S}_{\mathbb{S}} : \mathsf{pdg}~\mathsf{Col}_\mathbb{S} \arrow[r]
&\mathsf{End}(\mathsf{pdg}~\smod) \\
E \arrow[r,mapsto]
&\mathscr{S}_{\mathbb{S}}(E)(-)~.
\end{tikzcd}
\]
The Schur endofunctor $\mathscr{S}_{\mathbb{S}}(E)$ is given, for $M$ a pdg $\mathbb{S}$-module, by:
\[
\mathscr{S}_{\mathbb{S}}(E)(M)(n) \coloneqq \bigoplus_{(n_1,\cdots,n_r) \in \mathbb{N}^{r}~,~r\geq1} E(n_1,\cdots,n_r;n) \otimes_{\left(\mathbb{S}_{n_1} \times \cdots \times \mathbb{S}_{n_r}\right)~ \wr ~ \mathbb{S}_r} \left(M(n_1) \otimes \cdots \otimes M(n_r) \right)~.
\]

\begin{lemma}
The colored Schur realization functor $\mathscr{S}_{\mathbb{S}}$ is a strong monoidal functor between the monoidal categories $(\mathsf{pdg}~\mathsf{Col}_\mathbb{S},\circ, \I_{\mathbb{S}})$ and $(\mathsf{End}(\mathsf{pdg}~\smod),\circ, \mathrm{Id})~.$
\end{lemma}

\begin{proof}
The computation is analogous to the classical case.
\end{proof}

\begin{Corollary}
Let $E$ be an $\mathbb{S}$-color scheme. Any $\mathbb{S}$-colored operad structure on $E$ induces a monad structure on $\mathscr{S}_{\mathbb{S}}(E)$ and any $\mathbb{S}$-colored cooperad structure on $E$ induces a comonad structure on $\mathscr{S}_{\mathbb{S}}(E)$.
\end{Corollary}

\begin{Definition}[Algebra over an $\mathbb{S}$-colored operad]
Let $\mathcal{G}$ be a pdg $\mathbb{S}$-colored operad. A $\mathcal{G}$-algebra $M$ amounts to the data $(M,\gamma_M,d_M)$ of an algebra over the monad $\mathscr{S}_{\mathbb{S}}(\mathcal{G})$.
\end{Definition}

\begin{lemma}
Let $\mathcal{G}$ be a pdg $\mathbb{S}$-colored operad. The data of a pdg $\mathcal{G}$-algebra structure $\gamma_M$ on a pdg $\mathbb{S}$-module $(M,d_M)$ is equivalent to a morphism of pdg $\mathbb{S}$-colored operads $\Gamma_M: \mathcal{G} \longrightarrow \mathrm{End}_M$. 
\end{lemma}

\begin{proof}
Straightforward generalization of the standard case.
\end{proof}

\begin{Definition}[Coalgebra over an $\mathbb{S}$-colored cooperad]
Let $\mathcal{V}$ be a pdg $\mathbb{S}$-colored cooperad. A $\mathcal{V}$-coalgebra $K$ amounts to the data $(K,\Delta_K,d_K)$ of a coalgebra over the comonad $\mathscr{S}_{\mathbb{S}}(\mathcal{V})$.
\end{Definition}

The \textit{dual Schur realization functor} $\widehat{\mathscr{S}}_{\mathbb{S}}^c$ also associates to any pdg $\mathbb{S}$-color scheme $E$ an endofunctor in the category of pdg $\mathbb{S}$-modules 
\[
\begin{tikzcd}[column sep=4pc,row sep=0.5pc]
\widehat{\mathscr{S}}_{\mathbb{S}}^c : \mathsf{pdg}~\mathsf{Col}_\mathbb{S}^{\text{op}} \arrow[r]
&\mathsf{End}(\mathsf{pdg}~\smod) \\
E \arrow[r,mapsto]
&\widehat{\mathscr{S}}_{\mathbb{S}}^c(E)(-)~.
\end{tikzcd}
\]
The dual Schur endofunctor $\widehat{\mathscr{S}}_{\mathbb{S}}^c(E)$ is given, for $M$ a pdg $\mathbb{S}$-module, by:
\[
\widehat{\mathscr{S}}_{\mathbb{S}}^c (E)(M)(n) \coloneqq \prod_{(n_1,\cdots,n_r) \in \mathbb{N}^{r}~,~r\geq1} \mathrm{Hom}_{\left(\mathbb{S}_{n_1} \times \cdots \times \mathbb{S}_{n_r}\right)~ \wr ~ \mathbb{S}_r} \left( E(n_1,\cdots,n_r;n), M(n_1) \otimes \cdots \otimes M(n_r) \right)~.
\]

\begin{lemma}
The functor $\widehat{\mathscr{S}}_{\mathbb{S}}^c: (\mathsf{pdg}~\mathsf{Col}_\mathbb{S},\circ, \I_{\mathbb{S}})^\mathsf{op} \longrightarrow (\mathsf{End}(\mathsf{pdg}~\smod),\circ, \mathrm{Id})$ can be endowed with a lax monoidal structure. That is, there exists a natural transformation
\[
\varphi_{E,F}: \widehat{\mathscr{S}}_{\mathbb{S}}^c(E) \circ \widehat{\mathscr{S}}_{\mathbb{S}}^c(F) \longrightarrow \widehat{\mathscr{S}}_{\mathbb{S}}^c(E \circ F)~,
\]
which satisfies associativity and unitality compatibility conditions with respect to the monoidal structures. Furthermore, this natural transformation is a monomorphism for all pdg $\mathbb{S}$-color schemes $E,F$. 
\end{lemma}

\begin{proof}
The construction of $\varphi_{E,F}$ and the proof are completely analogue to \cite[Corollary 3.4]{grignoulejay18}.
\end{proof}

\begin{Corollary}
Let $E$ be an $\mathbb{S}$-color scheme. Any $\mathbb{S}$-colored cooperad structure on $E$ induces a monad structure on $\widehat{\mathscr{S}}_{\mathbb{S}}^c(E)$.
\end{Corollary}

\begin{Definition}[Algebra over an $\mathbb{S}$-colored cooperad]
Let $\mathcal{V}$ be an $\mathbb{S}$-colored cooperad. A $\mathcal{V}$-algebra $L$ amounts to the data $(L,\gamma_L,d_L)$ of an algebra over the monad $\widehat{\mathscr{S}}_{\mathbb{S}}^c(\mathcal{V})$.
\end{Definition}

\begin{Definition}[Coalgebra over an $\mathbb{S}$-colored operad]
Let $\mathcal{G}$ be a pdg $\mathbb{S}$-colored operad. A $\mathcal{G}$-coalgebra $N$ amount to the $(N,\Delta_N,d_N)$ of a morphism of pdg $\mathbb{S}$-modules
\[
\Delta_N: N \longrightarrow \displaystyle \prod_{(n_1,\cdots,n_r) \in \mathbb{N}^{r}~,~r\geq1} \mathrm{Hom}_{\left(\mathbb{S}_{n_1} \times \cdots \times \mathbb{S}_{n_r}\right)~\wr~ \mathbb{S}_r } \left( \mathcal{G}(n_1,\cdots,n_r;n), N(n_1) \otimes \cdots \otimes N(n_r) \right)~,
\]
such that the following diagram commutes 
\[
\begin{tikzcd}[column sep=4.5pc,row sep=3pc]
N \arrow[r,"\Delta_N"] \arrow[d,"\Delta_N",swap] 
&\widehat{\mathscr{S}}_{\mathbb{S}}^c(\mathcal{G})(N) \arrow[r,"\widehat{\mathscr{S}}_{\mathbb{S}}^c(\mathrm{id})(\Delta_N)"]
&\widehat{\mathscr{S}}_{\mathbb{S}}^c(\mathcal{G})(N) \circ \widehat{\mathscr{S}}_{\mathbb{S}}^c(\PP)(D) \arrow[d,"\varphi_{\mathcal{G},\mathcal{G}}(N)"] \\
\widehat{\mathscr{S}}_{\mathbb{S}}^c(\mathcal{G})(N) \arrow[rr,"\widehat{\mathscr{S}}_{\mathbb{S}}^c(\gamma_\mathcal{G})(\mathrm{id})"]
&
&\widehat{\mathscr{S}}_{\mathbb{S}}^c(\mathcal{G} \circ \mathcal{G})(N)~.
\end{tikzcd}
\]
\end{Definition}

\begin{Remark}
Let $\mathcal{G}$ be a pdg $\mathbb{S}$-colored operad. The data of a pdg $\mathcal{G}$-coalgebra structure $\Delta_N$ on a pdg $\mathbb{S}$-module $(N,d_N)$ is equivalent to a morphism of pdg $\mathbb{S}$-colored operads $\delta_N: \mathcal{G} \longrightarrow \mathrm{Coend}_N$. 
\end{Remark}

\medskip

\subsection{Partial $\mathbb{S}$-colored (co)operads}\label{subsection: comparison for partial colored (co)operads}
As the careful reader might expect by now, all the other definitions of Section \ref{Section: partial (co)operads} generalize to the $\mathbb{S}$-colored setting \textit{mutatis mutandis}.

\begin{Definition}[Partial $\mathbb{S}$-colored operad]
An \textit{pdg partial} $\mathbb{S}$-\textit{colored operad} amounts to the data of $(\mathcal{G},\{\circ_i\},d_\mathcal{G})$ a pdg $\mathbb{S}$-color scheme $\mathcal{G}$ endowed with partial composition operations 
\[
\circ_i : \mathcal{G}(n_1,\cdots, n_i, \cdots, n_r;n) \otimes_{\mathbb{S}_{n_i}} \mathcal{G}(p_1,\cdots,p_l;n_i) \longrightarrow \mathcal{G}(n_1,\cdots,n_{i-1},p_1,\cdots,p_l,n_{i+1},\cdots,n_r;n)~.
\]
This family of partial compositions maps $\{\circ_i\}$ satisfies sequential and parallel axioms analogue to those of a partial operad. They satisfy an equivariance condition with respect to the permutations of the entries $\varphi_\sigma$ which is also completely analogue to Definition \ref{def: partialoperad}. 
\end{Definition}

\begin{Definition}[Unital partial $\mathbb{S}$-colored operad]
A \textit{pdg unital partial} $\mathbb{S}$-colored operad amounts to the data of $(\mathcal{G},\{\circ_i\},\eta,d_\mathcal{G})$ a partial pdg $\mathbb{S}$-colored operad $(\mathcal{G},\{\circ_i\},d_\mathcal{G})$ together with a morphism of pdg $\mathbb{S}$-color schemes $\eta: \I_{\mathbb{S}} \longrightarrow \mathcal{G}$ which acts as a unit for the partial compositions maps. 
\end{Definition}

The data of a unit $\eta: \I_{\mathbb{S}} \longrightarrow \mathcal{G}$ amounts to a family of elements $\{id_n \}$ in $\mathcal{G}(n;n)$. 

\begin{lemma}\label{lemma: add a unit colored case}
The category of pdg unital partial $\mathbb{S}$-colored operads is equivalent to the category of pdg $\mathbb{S}$-colored operads defined as monoids.
\end{lemma}

\begin{proof}
Straightforward generalization of Proposition \ref{prop: unital partial operads are operads}.
\end{proof}

Thus, given a pdg partial $\mathbb{S}$-colored operad $(\mathcal{G},\{\circ_i\},\eta,d_\mathcal{G})$ one can obtain a pdg $\mathbb{S}$-colored operad defined as a monoid by freely adding a unit to it $\mathcal{G}^{u} \coloneqq \mathcal{G} \oplus \I_\mathbb{S}$.

\begin{theorem}[{\cite[Theorems 2.10, 2.11 and 2.24]{ward19}}]
Let $\mathbb{V}$ be a groupoid. 
\begin{enumerate}
\item The exists a monad $\mathscr{T}_\mathbb{V}$ in the category of pdg $\mathbb{V}$-color schemes, called the $\mathbb{V}$-colored tree monad, such that the category of $\mathscr{T}_\mathbb{V}$-algebras is equivalent to the category of pdg unital partial $\mathbb{V}$-colored operads.

\item  The exists a monad $\overline{\mathscr{T}}_\mathbb{V}$ in the category of pdg $\mathbb{V}$-color schemes, called the reduced $\mathbb{V}$-colored tree monad, such that the category of $\overline{\mathscr{T}}_\mathbb{V}$-algebras is equivalent to the category of pdg partial $\mathbb{V}$-colored operads.
\end{enumerate}
\end{theorem}

\begin{Definition}[Partial $\mathbb{S}$-colored cooperad]
A \textit{pdg partial} $\mathbb{S}$-\textit{colored cooperad} amounts to the data of $(\mathcal{V},\{\Delta_i\},d_\mathcal{V})$ a pdg $\mathbb{S}$-color scheme $\mathcal{V}$ endowed with partial decomposition operations 
\[
\Delta_i :\mathcal{V}(n_1,\cdots,n_{i-1},p_1,\cdots,p_l,n_{i+1},\cdots,n_r;n) \longrightarrow \mathcal{V}(n_1,\cdots, n_i, \cdots, n_r;n) \otimes_{\mathbb{S}_{n_i}} \mathcal{V}(p_1,\cdots,p_l;n_i)~.
\]
This family of partial decompositions maps $\{\Delta_i\}$ satisfies cosequential and coparallel axioms analogue to those of a partial cooperad. They satisfy an equivariance condition with respect to the permutations of the entries $\varphi_\sigma$ which is also completely analogue to Definition \ref{def: partialcoop}. 
\end{Definition}

\begin{Definition}[Counital partial $\mathbb{S}$-colored cooperad]
A \textit{pdg counital partial} $\mathbb{S}$-colored cooperad  amounts to the data of $(\mathcal{V},\{\Delta_i\},\epsilon, d_\mathcal{V})$ a pdg partial $\mathbb{S}$-colored cooperad $(\mathcal{V},\{\Delta_i\},d_\mathcal{V})$ together with a morphism of pdg $\mathbb{S}$-color schemes $\epsilon: \mathcal{V} \longrightarrow \I_{\mathbb{S}}$ which acts as a counit for the partial decompositions maps. 
\end{Definition}

Any pdg partial $\mathbb{S}$-colored cooperad structure induces a morphism of pdg $\mathbb{S}$-color schemes
\[
\Delta_\mathcal{V}: \mathcal{V} \longrightarrow \overline{\mathscr{T}}^\wedge_\mathbb{S}(\mathcal{V})~,
\]
where $\overline{\mathscr{T}}^\wedge_\mathbb{S}$ is the endofunctor given by the completion of the reduced $\mathbb{V}$-colored tree monad with respect to its canonical weight filtration given by the number of internal edges of the rooted trees. Using this morphism one defines an analogue version of the coradical filtration of a partial cooperad and an analogue definition of a conilpotent pdg partial $\mathbb{S}$-colored cooperad. \textit{Mutatis mutandis} the same characterization of this type of partial cooperads still holds. See Section \ref{Section: Filtration on partial (co)operads}.

\begin{theorem}[{\cite[Section 2.4.1]{ward19}}]
Let $\mathbb{V}$ be a groupoid such that for any $v$ in $\mathbb{V}$, the set $\mathrm{Aut}(v)$ is finite.
\begin{enumerate}
\item The exists a comonad structure on the underlying endofunctor of the reduced tree monad $\overline{\mathscr{T}}_\mathbb{V}$ on the category of $\mathbb{V}$-colored schemes.

\item The category of $\overline{\mathscr{T}}_\mathbb{V}$-coalgebras is equivalent to the category of conilpotent pdg partial $\mathbb{S}$-colored cooperads.
\end{enumerate}
\end{theorem}

Let $(\mathcal{V},\{\Delta_i\},d_\mathcal{V})$ be a conilpotent pdg partial cooperad. Then one can cofreely add a counit to it by setting $\mathcal{V}^{u} \coloneqq \mathcal{V} \oplus \I_{\mathbb{V}}$. The pdg $\mathbb{S}$-color scheme has an unique cooperad structure induced by the partial decompositions maps $\{\Delta_i\}$. It defines a functor 
\[
\mathsf{Conil}: \mathsf{pdg}~\mathsf{pCoop}_\mathbb{S}^\mathsf{conil} \longrightarrow \mathsf{pdg}~\mathsf{Coop}_\mathbb{S}~.
\]
\begin{Definition}[Conilpotent $\mathbb{S}$-colored cooperad]
A pdg $\mathbb{S}$-colored cooperad $\mathcal{V}$ is said to be \textit{conilpotent} if its in the essential image of the functor defined above.
\end{Definition}

\begin{Remark}
For any conilpotent pdg $\mathbb{S}$-colored cooperad $\mathcal{V}$, there is a canonical filtration on $\mathcal{V}$-algebras and a notion of \textit{complete} $\mathcal{V}$-algebra analogous to the case studied in Section \ref{Section: coalgebras and algebras}.
\end{Remark}

\subsection{Koszul duality for quadratic $\mathbb{S}$-colored operads}\label{subsection: Koszul duality for quadratic S-colored operads}
We now state the quadratic Koszul duality of \cite{ward19} in the context of $\mathbb{S}$-colored operads. In order to do this, we briefly restrict from the underlying category of pdg modules to the underlying category of dg modules. 

\begin{Proposition}\label{prop S colored totalization}
Let $(\mathcal{G},\{\circ_i\},d_\mathcal{G})$ be a dg partial $\mathbb{S}$-colored operad. The \textit{totalization} of $\mathcal{G}$ is given by 
\[
\prod_{(n_1,\cdots,n_r;n) \in \mathbb{N}^{r+1}} \left(\mathcal{G}(n_1,\cdots,n_r;n)^{\left(\mathbb{S}_{n_1} \times \cdots \times \mathbb{S}_{n_r}\right)~\wr~ \mathbb{S}_r}\right)^{\mathbb{S}_n}~.
\]
It can be endowed with a dg pre-Lie algebra structure. Let $\mu$ be in $\mathcal{G}(n_1,\cdots,n_r;n)$ and $\nu$ be in $\mathcal{G}(p_1,\cdots,p_l;p)$, where $\mu$ is colored by $(\sigma_1, \cdots, \sigma_r;\sigma)$ in $\mathbb{S}_{n_1} \times \cdots \times \mathbb{S}_{n_r} \times \mathbb{S}_n$ and where $\nu$ is colored by $(\tau_1, \cdots, \tau_l;\tau)$ in $\mathbb{S}_{p_1} \times \cdots \times \mathbb{S}_{p_l} \times \mathbb{S}_p$. The pre-Lie bracket of $\mu$ and $\nu$ is given by
\[
\mu \star \nu \coloneqq \sum_{\{(n_i,\sigma_i) = (p,\tau)\}} \sum_{\lambda \in \mathbb{S}_l} (\mu \circ_i \nu)^{\lambda^\flat}~,
\]
where the sum runs over all pairs $(n_i,\sigma_i)$ which are equal to $(p,\tau)$, and where $\lambda^{\flat}$ is the unique permutation of $\mathbb{S}_{r+l-1}$ which acts as $\lambda$ on $\{i, \cdots, i + r-1\}$ and as the identity elsewhere. If there are none, then $\mu \star \nu \coloneqq 0$.
\end{Proposition}

\begin{proof}
The associator of the product $\star$ is right symmetric like in the case of partial operads, and is compatible with the differential by definition.
\end{proof}

\begin{Definition}[Convolution partial $\mathbb{S}$-colored operad]
Let $(\mathcal{V},\{\Delta_i\}, d_\mathcal{V})$ be a conilpotent dg partial $\mathbb{S}$-colored cooperad and let $(\mathcal{G},\{\circ_i\},d_\mathcal{G})$ be a dg partial $\mathbb{S}$-colored operad. The dg $\mathbb{S}$-color scheme
\[
\mathcal{H}om(\mathcal{V},\mathcal{G})(n_1,\cdots,n_r;n) \coloneqq \mathrm{Hom}_{\mathsf{dg}}(\mathcal{V}(n_1,\cdots,n_r;n), \mathcal{G}(n_1,\cdots,n_r;n))
\]
endowed with the differential 
\[
\partial(\alpha) \coloneqq d_{\mathcal{G}} \hspace{2pt} \circ \hspace{2pt} \alpha - (-1)^{|\alpha|} \alpha \hspace{2pt} \circ \hspace{2pt} d_{\mathcal{V}}~.
\]
forms a dg partial $\mathbb{S}$-colored operad structure where the partial compositions maps are given by 
\[
\alpha \circ_i \beta \coloneqq \circ_i \cdot(\alpha \otimes \beta) \cdot \Delta_i~.
\]
\end{Definition}

\begin{Definition}[Twisting morphism]
Let 
\[
\mathfrak{g}_{\mathcal{V},\mathcal{G}} \coloneqq \prod_{(n_1,\cdots,n_r;n) \in \mathbb{N}^{r+1}} \mathrm{Hom}_{\mathbb{S}\text{-}\mathsf{col}}(\mathcal{V}(n_1,\cdots,n_r;n), \mathcal{G}(n_1,\cdots,n_r;n))
\]
be the dg pre-Lie algebra given by the totalization of the convolution operad of $\mathcal{V}$ and $\mathcal{G}$. A \textit{twisting morphism} is a Maurer-Cartan element $\alpha$ of $\mathfrak{g}_{\mathcal{V},\mathcal{G}}$, that is, a morphism $\alpha: \mathcal{V} \longrightarrow \mathcal{G}$ of dg $\mathbb{S}$-color schemes of degree $-1$  satisfying :
\[
\partial(\alpha) + \alpha \star \alpha = 0~.
\]
The set of twisting morphism between $\mathcal{V}$ and $\mathcal{G}$ is denoted by $\mathrm{Tw}(\mathcal{V},\mathcal{G})~.$
\end{Definition}

Twisting morphisms induce Bar-Cobar adjunctions relative to them.

\begin{Proposition}[Bar-Cobar adjunction relative to $\alpha$]\label{Bar-Cobar colored adjunction, classical case}
Let $\alpha: \mathcal{V} \longrightarrow \mathcal{G}$ be a twisting morphism between a conilpotent partial $\mathbb{S}$-colored cooperad and a partial $\mathbb{S}$-colored operad. It induces a Bar-Cobar adjunction
\[
\begin{tikzcd}[column sep=5pc,row sep=3pc]
          \mathsf{dg}~\mathcal{G}^{u}\text{-}\mathsf{alg} \arrow[r, shift left=1.1ex, "\Omega_\alpha"{name=F}] &\mathsf{dg}~\mathcal{V}^{u}\text{-}\mathsf{coalg}^{\mathsf{conil}} \arrow[l, shift left=.75ex, "\text{B}_\alpha"{name=U}]
            \arrow[phantom, from=F, to=U, , "\dashv" rotate=-90]
\end{tikzcd}
\]
relative to the twisting morphism $\alpha$. 
\end{Proposition}

\begin{proof}
This adjunction is constructed using the free $\mathcal{G}^{u}$-algebra given by $\mathscr{S}_{\mathbb{S}}(\mathcal{G}^{u})$ and the cofree conilpotent $\mathcal{V}^{u}$-coalgebra given by $\mathscr{S}_{\mathbb{S}}(\mathcal{V}^{u})$. They are endowed with differentials which are analogous to those constructed in the standard case \cite[Section 11.2]{LodayVallette12}. See \cite[Section 2.7.1]{ward19} for a detailed exposition of this construction.
\end{proof}

\begin{Definition}[Koszul twisting morphism]
Let $\alpha: \mathcal{V} \longrightarrow \mathcal{G}$ be a twisting morphism between a conilpotent partial $\mathbb{S}$-colored cooperad and a   partial $\mathbb{S}$-colored operad. It is a \textit{Koszul twisting morphism} if the twisted complex $\mathcal{V}^{u} \circ_\alpha \mathcal{G}^{u}$ is acyclic or if the twisted complex $\mathcal{G}^{u} \circ_\alpha \mathcal{V}^{u}$ is acyclic. 
\end{Definition}

See \cite[Section 2.7]{ward19} or \cite[Section 6.4.5]{LodayVallette12} for a precise definition of the twisted complexes $\mathcal{G}^{u} \circ_\alpha \mathcal{V}^{u}$ and $\mathcal{V}^{u} \circ_\alpha \mathcal{G}^{u}$.

\begin{Remark}
There is a natural isomorphism $\Omega_\alpha \circ \mathrm{B}_\alpha \cong \mathscr{S}_\mathbb{S}(\mathcal{G}^{u} \circ_\alpha \mathcal{V}^{u})$. This implies that $\alpha$ is Koszul, then the counit of adjunction 
\[
\Omega_\alpha \circ \mathrm{B}_\alpha \longrightarrow \mathrm{Id}_{\mathcal{G}\text{-}\mathsf{alg}}
\]
is a natural quasi-isomorphism. The converse is true under some assumptions, see \cite[Theorem 11.3.3]{LodayVallette12}.
\end{Remark}

\begin{lemma}[{\cite[Section 2.48]{ward19}}]\label{Koszul duality Ward}
Let $\alpha: \mathcal{V} \longrightarrow \mathcal{G}$ be a twisting morphism between a conilpotent partial $\mathbb{S}$-colored cooperad and a partial $\mathbb{S}$-colored operad. The following are equivalent:
\begin{enumerate}
\item The twisting morphism $\alpha$ is a Koszul twisting morphism.
\item The morphism $f_\alpha: \Omega_{\mathbb{S}}\mathcal{V} \longrightarrow \mathcal{G}$ induced by $\alpha$ is a quasi-isomorphism.
\item The morphism $g_\alpha: \mathcal{V} \longrightarrow \mathrm{B}_{\mathbb{S}}\mathcal{G}$ induced by $\alpha$ is a quasi-isomorphism.
\end{enumerate}
Here $\mathrm{B}_{\mathbb{S}}$ and $\Omega_{\mathbb{S}}$ denote respectively the Bar and Cobar constructions for partial groupoid-colored (co)operads constructed in \cite{ward19}.
\end{lemma}

\begin{Definition}[{\cite[Section 2.5.3]{ward19}}]
Let $(\mathcal{G}, \{\circ_i\}, d_\mathcal{G})$ be a partial $\mathbb{S}$-colored operad. A \textit{quadratic presentation} of $\mathcal{G}$ is a pair $(E,R)$ such that there is an isomorphism of partial $\mathbb{S}$-colored operads
\[
\mathcal{G} \cong \overline{\mathscr{T}}_{\mathbb{S}}(E)/(R)~,
\]
where $\overline{\mathscr{T}}_{\mathbb{S}}(E)$ denotes the free partial $\mathbb{S}$-colored operad on the $\mathbb{S}$-color scheme $E$ modulo the quadratic relations $R \subset \overline{\mathscr{T}}_{\mathbb{S}}(E)^{(2)}~.$
\end{Definition}

\begin{Definition}[{\cite[Section 2.5.4]{ward19}}]
Let $(\mathcal{G}, \{\circ_i\}, d_\mathcal{G})$ be a partial $\mathbb{S}$-colored operad and $(E,R)$ be a quadratic presentation of $\mathcal{G}$. The \textit{Koszul dual conilpotent partial} $\mathbb{S}$-\textit{colored cooperad} $\mathcal{G}^{\ac}$ of $\mathcal{G}$ is given by the conilpotent partial $\mathbb{S}$-colored cooperad 
\[
\mathcal{G}^{\ac} \coloneqq \overline{\mathscr{T}}_{\mathbb{S}}^c(sE,s^2R)~,
\]
where $\overline{\mathscr{T}}_{\mathbb{S}}^c(sE,s^2R)$ denotes the cofree conilpotent partial $\mathbb{S}$-colored cooperad  cogenerated by the $\mathbb{S}$-color scheme $sE$, with $s^2R$ as the corelations. (Smallest sub-cooperad of $\overline{\mathscr{T}}_{\mathbb{S}}^c(sE)$ containing $s^2R$). Here $s$ denotes the suspension of graded $\mathbb{S}$-color schemes.
\end{Definition}

It comes equipped with a canonical twisting morphism $\kappa: \mathcal{G}^{\ac} \longrightarrow \mathcal{G}$ given by
\[
\mathcal{G}^{\ac} \twoheadrightarrow sE \cong E \hookrightarrow \mathcal{G}~. 
\]

\begin{Definition}[Koszul partial $\mathbb{S}$-colored operad]\label{def: Koszul quadratic S-colored operad}
Let $(\mathcal{G}, \{\circ_i\}, d_\mathcal{G})$ be a partial $\mathbb{S}$-colored operad. The partial $\mathbb{S}$-colored operad $\mathcal{G}$ is said to be a Koszul quadratic $\mathbb{S}$-colored operad if there exists a quadratic presentation $(E,R)$ of $\mathcal{G}$ such that the canonical twisting morphism 
\[
\kappa: \mathcal{G}^{\ac} \twoheadrightarrow sE \cong E \hookrightarrow \mathcal{G}~. 
\]
is a Koszul twisting morphism.
\end{Definition}

\subsection{The partial $\mathbb{S}$-colored operad encoding partial (co)operads and its Koszul dual}
We now define the partial $\mathbb{S}$-colored operad $\mathcal{O}$ that encodes partial (co)operads as its (co)algebras. This is a direct generalization of the set-colored operad encoding non-symmetric partial (co)operads defined in \cite[Definition 4.1]{VanderLaan03}. The ground category is still dg modules in this section. 

\begin{Definition}[The partial $\mathbb{S}$-colored operad encoding partial (co)operads]\label{The colored operad O}
The partial $\mathbb{S}$-colored operad $\mathcal{O}$ is given by the following quadratic presentation: 
\[
\mathcal{O} \coloneqq \overline{\mathscr{T}}_{\mathbb{S}}(E)/(R)~,
\]
where $\overline{\mathscr{T}}_{\mathbb{S}}(E)$ denotes the free partial $\mathbb{S}$-colored operad on the $\mathbb{S}$-color scheme $E$ modulo the quadratic relations $R \subset \overline{\mathscr{T}}_{\mathbb{S}}(E)^{(2)}$. The $\mathbb{S}$-color scheme $E$ is the generated by:
\[
\gamma_i^{n,k} \in E(n,k;n+k-1) \hspace{1pc} \textrm{for} \hspace{1pc} 1 \leq i \leq n~,
\]
as a right $\mathbb{K}[\mathbb{S}_n \times \mathbb{S}_k] \otimes_\mathbb{K} \mathbb{K}[\mathbb{S}_{n+k-1}]^{\text{op}}$-module, endowed with the zero differential.

\medskip

The right action of the symmetric groups $\mathbb{S}_n$ and $\mathbb{S}_k$ can be written as in both cases as the left action of an element $\mathbb{S}_{n+k-1}$. For $\sigma$ in $\mathbb{S}_k$, its action is equal to 
\[
\sigma \cdot (\gamma_i^{n,k}) = (\gamma_i^{n,k}) \cdot \sigma' \hspace{1pc} \textrm{for} \hspace{1pc} 1 \leq i \leq n~,
\]
where $\sigma' \in \mathbb{S}_{n+k-1}$ is the unique permutation that acts as the identity everywhere except for $\{i,\cdots,i+k-1\}$, where it acts as $\sigma$. For $\tau$ in $\mathbb{S}_n$, its action is equal to
\[
\tau \cdot (\gamma_i^{n,k}) = (\gamma_{\tau(i)}^{n,k}) \cdot \tau' \hspace{1pc} \textrm{for} \hspace{1pc} 1 \leq i \leq n~,
\]
where $\tau'$ is the unique permutation of $\mathbb{S}_{n+k-1}$ that acts as $\tau$ on the block $\{1,\cdots,n+n-1\} - \{i,\cdots, i+k-1\}$ with values in $\{1,\cdots,n+k-1\} - \{ \tau(i),\cdots, \tau(i)+k-1\}$ and as the identity on $\{i,\cdots,i+k-1\}$ with values in $\{\tau(i), \cdots, \tau(i) + k -1 \}~.$ The $\mathbb{S}_2$ action that permutes the entries is simply given by $\varphi_{(12)} (\gamma_i^{n,k}) \coloneqq \gamma_i^{k,n}~.$

\medskip

The relations $R$ are the given by: 
\[
\begin{tikzcd}
&\left\{ 
    \arraycolsep=1.4pt\def\arraystretch{1.8}\begin{array}{lllll}
         \gamma_{i+j-1}^{n+k-1,l} \circ_1 \gamma_i^{n,k}= \gamma_i^{l+n-1,k} \circ_1 \gamma_j^{n,l} \hspace{1pc}\mathrm{if} \hspace{1pc} 1 \leq i \leq n, 1 \leq j \leq k \hspace{1pc}~, (1)\\
         \gamma_{j+k-1}^{n+k-1,l} \circ_1 \gamma_i^{n,k} = \gamma_j^{n+k-1,l}\circ_1 \gamma_i^{n,k} \hspace{1pc}\mathrm{if} \hspace{1pc} 1\leq i < j \leq n \hspace{1pc} (2)~.
    \end{array}
\right.
\end{tikzcd}
\]
The relation $(1)$ is the \textit{parallel} axiom and $(2)$ is the \textit{sequential} axiom. 
\end{Definition}

\begin{Remark}
The equivariance relations of the partial (de)compositions of partial (co)operads are not coded as relations in the partial $\mathbb{S}$-colored operad $\mathcal{O}$. They are encoded by the action of the symmetric groups on the generators $\{\gamma_i^{n,k}\}$.
\end{Remark}

By adding freely a unit to $\mathcal{O}$ one obtains an $\mathbb{S}$-colored operad $\mathcal{O}^u$ defined as a monoid.

\begin{Proposition}\label{Prop: O-algebres}
The category of dg $\mathcal{O}^u$-algebras is equivalent to the category of dg partial operads.
\end{Proposition}

\begin{proof}
Let $(P,d_P)$ be a dg-$\mathbb{S}$-module and let $\Gamma_P: \mathcal{O}^{u} \longrightarrow \mathrm{End}_P$ be a morphism of dg $\mathbb{S}$-colored operads. It is determined by the image of the generators $\gamma_i^{n,k}$ in $\mathcal{O}(n,k;n+k-1)$. Let's denote $\circ_i^{n,k} \coloneqq \Gamma_P(\gamma_i^{n,k}): P(n) \otimes P(k) \longrightarrow P(n+k-1)$. The family $\{\circ_i^{n,k}\}$ satisfies the equivariance axiom of a partial operad since $\Gamma_P$ is $\kk[\mathbb{S}_n \times \mathbb{S}_k] \otimes \kk[\mathbb{S}_{n+k-1}]$-equivariant; it satisfies the parallel and sequential axioms by definition of $\mathcal{O}$. Since $\Gamma_P$ commutes with the pre-differentials, $d_P$ is a derivation with respect to the partial composition maps. Hence it endows $(P,d_P)$ with a dg partial operad structure. 
\end{proof}

\begin{Proposition}\label{Prop: O-coalgebras}
The category of dg $\mathcal{O}^u$-coalgebras is equivalent to the category of dg partial cooperads.
\end{Proposition}

\begin{proof}
The proof is analogue to the proof of Proposition \ref{Prop: O-algebres}. One replaces the endomorphism operad with the coendomorphism operad.
\end{proof}

\begin{lemma}
The Koszul dual conilpotent partial $\mathbb{S}$-colored cooperad $\mathcal{O}^{\ac}$ is isomorphic to the operadic suspension of the linear dual of $\mathcal{O}$. The operad $\mathcal{O}$ is Koszul autodual. 
\end{lemma}

\begin{proof}
First, since each $\mathcal{O}(n_1,\cdots,n_r;n)$ is finite dimensional over $\kk$ and since it is \textit{reduced} in the sense of \cite[Definition 2.35]{ward19}, the arity-wise linear dual of $\mathcal{O}$ is a conilpotent partial $\mathbb{S}$-colored operad, and we have that $(\mathcal{O}^*)^* \cong \mathcal{O}$. Therefore it is only necessary to compute its Koszul dual partial $\mathbb{S}$-colored operad, which is given by the operadic suspension of $(\mathcal{O}^{\ac})^*$.

\medskip

We adapt the proof of \cite[Theorem 4.3]{VanderLaan03} to the groupoid-colored. The Koszul dual operad $\mathcal{O}^!$ is given by $\overline{\mathscr{T}}_{\mathbb{S}}(E^*)/(R^{\bot})$, where $\overline{\mathscr{T}}_{\mathbb{S}}(E)$ is the free partial $\mathbb{S}$-colored partial operad generated by $E^*$, and $(R^{\bot})$ denotes the operadic ideal generated by the orthogonal of $R$ inside $\overline{\mathscr{T}}_{\mathbb{S}}(E^*)^{(2)}$. Using the same arguments as in the set-colored case, one can show that the dimension of the relations $R$ as a $\mathbb{S}$-color scheme is exactly half of the dimension over $\kk$ of $\overline{\mathscr{T}}_{\mathbb{S}}(E)^{(2)}$. Since, the relations $R$ must be contained in the ideal generated by the orthogonal, concludes by a dimension argument that the relations in $R$ form a basis of $R^\bot$. 
\end{proof}

\begin{lemma}\label{lemma with the trees bijection}
The following statements hold: 
\begin{enumerate}
\item There is an isomorphism of monads in the category of $\mathbb{S}$-modules between the monad $\mathscr{S}_{\mathbb{S}}(\mathcal{O}^{u})$ and the reduced tree monad $\overline{\mathscr{T}}$ encoding partial operads.

\item The is an isomorphism of comonads in the category of $\mathbb{S}$-modules between the comonad $\mathscr{S}_{\mathbb{S}}((\mathcal{O}^*)^{u})$ and the reduced tree comonad $\overline{\mathscr{T}}^c$ encoding conilpotent partial cooperads.
\end{enumerate}
\end{lemma}

\begin{proof}
The first results follows immediately from the fact that the category of $\mathcal{O}$-algebras is equivalent to the category of partial operads. Nevertheless, for any $\mathbb{S}$-module $M$, there is an explicit bijection between the $\mathbb{S}$-modules $\mathscr{S}_{\mathbb{S}}(\mathcal{O}^{u})(M)$ and $\overline{\mathscr{T}}(M)$ that identifies the partial operad structures. For any element 
\[
(\psi;m_1, \cdots m_r) \quad \text{in} \quad \mathcal{O}(n_1,\cdots,n_r;n) \allowbreak \otimes_{\left(\mathbb{S}_{n_1} \times \cdots \times \mathbb{S}_{n_r}\right)~\wr~ \mathbb{S}_r} M(n_1) \otimes \cdots \otimes M(n_r)~,
\]
one can represent $\psi$ as a equivalence class of binary trees with vertices labeled by the generators $\gamma_i^{n,k}$ and with edges labeled by elements of the symmetric groups that match the coloring of the edge. To any $(\psi;m_1, \cdots m_r)$, one associates the rooted tree $\tau_\psi$ obtained by composing the $n_i$-corollas labeled by $m_i$ in the way the binary tree $\psi$ indicates, applying at each step the permutations that label the edges of $\psi$. Pictorially, the bijection is given by:

\begin{center}
\includegraphics[width=130mm,scale=1.7]{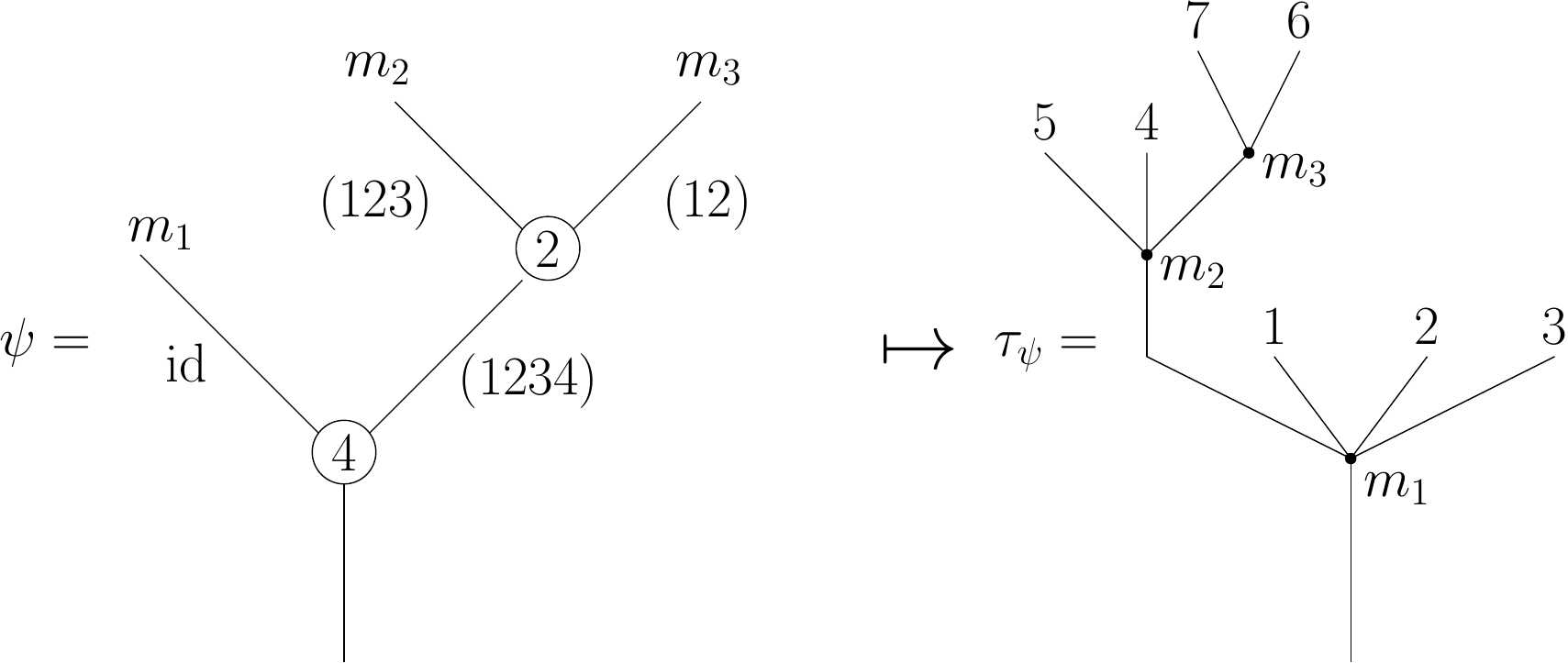}
\end{center}

This type of bijections can be found in \cite[Section 1.3]{dehlingvallette15}). It is straightforward to check that is bijection identifies the partial operad structures of $\mathscr{S}_{\mathbb{S}}(\mathcal{O}^{u})(M)$ and $\overline{\mathscr{T}}(M)$. Moreover, this bijection on the underlying $\mathbb{S}$-modules of $\mathscr{S}_{\mathbb{S}}((\mathcal{O}^*)^{u})(M)$ and $\overline{\mathscr{T}}^c(M)$ identifies the two conilpotent partial cooperad structures as well.
\end{proof}

\begin{Proposition}\label{Prop les alge/cog conil sur la coop O*}
Let $\mathcal{O}^{\ac}$ be the Koszul dual conilpotent partial $\mathbb{S}$-colored cooperad of $\mathcal{O}$.

\begin{enumerate}
\item The category of dg $(\mathcal{O}^{\ac})^{u}$-coalgebras is equivalent to the category of shifted conilpotent dg partial cooperads.

\item The category of complete dg $(\mathcal{O}^{\ac})^{u}$-algebras is equivalent to the category of shifted complete dg absolute partial operads.
\end{enumerate}
\end{Proposition}

\begin{proof}
The first statement follows from the previous Lemma. For a definition of \textit{absolute partial operads}, as well as their characterization, we refer to the Appendix \ref{Appendix B}.
\end{proof}

\begin{Proposition}
Let $\kappa: \mathcal{O}^{\ac} \longrightarrow \mathcal{O}$ be the canonical twisting morphism. The adjunction induced by the twisting morphism $\kappa$
\[
\begin{tikzcd}[column sep=5pc,row sep=3pc]
          \mathsf{dg}~\partialop \arrow[r, shift left=1.1ex, "\Omega_\kappa"{name=F}] &\mathsf{dg}~\partialcoop^{\mathsf{conil}} \arrow[l, shift left=.75ex, "\mathrm{B}_\kappa"{name=U}]
            \arrow[phantom, from=F, to=U, , "\dashv" rotate=-90]
\end{tikzcd}
\]
between dg partial operads and conilpotent dg partial cooperads is naturally isomorphic to the classical Bar-Cobar adjunction of \cite[Section 6.5]{LodayVallette12}.
\end{Proposition}

\begin{proof}
In order to check that two adjunctions are isomorphic, it is only necessary to check that there exists a natural isomorphism between left adjoints, since the mate of this natural isomorphism will also induce an natural isomorphism of right adjoints. The isomorphism between the reduced tree monad $\overline{\mathscr{T}}$ and the monad $\mathscr{S}_\mathbb{S}(\mathcal{O})$ induces a natural isomorphism of partial operads between $\Cobar$ and $\Cobar_\kappa$. It is straightforward to check that this bijection commutes with the respective differentials of the two Cobar constructions.
\end{proof}

\begin{theorem}\label{the Koszulity of O}
The canonical twisting morphism $\kappa: \mathcal{O}^{\ac} \longrightarrow \mathcal{O}$ is a Koszul twisting morphism. Hence the partial $\mathbb{S}$-colored operad $\mathcal{O}$ is Koszul.
\end{theorem}

\begin{proof}
We extend the arguments of \cite[Theorem 11.3.3]{LodayVallette12} to this specific case. Indeed, we know that for any dg partial operad $(\PP,\{\circ_i\},d_\PP)$, the counit $\epsilon_\kappa: \Cobar_\kappa \text{B}_\kappa \PP \longrightarrow \PP$ of the adjunction induced by $\kappa$ is a quasi-isomorphism of dg partial operads. In particular, any dg $\mathbb{S}$-module can be endowed with the trivial dg partial operad structure given by the family of partial composition maps $\{ 0 \}$. Therefore we conclude that $\mathcal{O} \circ_\kappa \mathcal{O}^{\ac}$ and $\mathcal{O}^{\ac} \circ_\kappa \mathcal{O}$ are both acyclic. 
\end{proof}

\begin{Remark}
Generalizing the formalism of \cite{grignoulejay18}, one can construct a complete Bar-Cobar adjunction relative to the twisting morphism $\kappa: \mathcal{O}^{\ac} \longrightarrow \mathcal{O}$ between complete dg absolute partial operads and non-necessarily conilpotent dg partial cooperads. 
\end{Remark}

\subsection{Inhomogeneous $\mathbb{S}$-colored Koszul duality}
We introduce the unital partial $\mathbb{S}$-colored operad $u\mathcal{O}$ which encodes (co)unital partial (co)operads as its (co)algebras. In order to compute its Koszul dual conilpotent curved $\mathbb{S}$-colored cooperad, we generalize the inhomogeneous Koszul duality of \cite{HirshMilles12} to the $\mathbb{S}$-colored case. This provides us with a Koszul dual conilpotent curved partial $\mathbb{S}$-colored cooperad $c\mathcal{O}^\vee$. It encodes conilpotent curved partial cooperads as its curved coalgebras and \textit{complete curved absolute partial operads} as its complete curved algebras. Absolute partial operads are a new type of operad-like structure, for which we provide useful descriptions in the Appendix \ref{Appendix B}. The ground category considered is once again the category of pdg modules. 

\medskip

Let $U$ be the $\mathbb{S}$-colored scheme given by $U(0;1) \coloneqq \kk.u$, where $u$ is an element of degree $0$, and zero elsewhere. 

\begin{Definition}[The unital partial $\mathbb{S}$-colored operad encoding (co)unital partial (co)operads]\label{def colored operad uO}
Let $(E,R)$ be the quadratic presentation of the partial $\mathbb{S}$-colored operad $\mathcal{O}$. The unital partial $\mathbb{S}$-colored operad $u\mathcal{O}$ is given by
\[
u\mathcal{O} \coloneqq \mathscr{T}_{\mathbb{S}}(E \oplus U)/(R')~,
\]
where $\mathscr{T}_{\mathbb{S}}(E \oplus U)$ is the free unital partial $\mathbb{S}$-colored operad generated by the $\mathbb{S}$-colored scheme $E \oplus U$, and where $(R')$ is the operadic ideal generated by $R'$. Here $R'$ is the sub-$\mathbb{S}$-color scheme of $\mathscr{T}_{\mathbb{S}}(E \oplus U)^{(\leq 2)}$ given by $R$, together with the additional relations 
\[
\begin{tikzcd}
\left\{ 
    \arraycolsep=1.4pt\def\arraystretch{1.8}\begin{array}{lllll}
         \gamma_{1}^{1,n} \circ_1 u = |_n \hspace{1pc}\mathrm{for} \hspace{1pc} n \in \mathbb{N}~,\\
         \gamma_{i}^{n,1} \circ_2 u = |_n \hspace{1pc}\mathrm{for}  \hspace{1pc} 1 \leq i \leq n \hspace{1pc} \mathrm{and} \hspace{1pc} n \in \mathbb{N}~,
    \end{array}
\right.
\end{tikzcd}
\]
where $|_n \in \mathscr{T}_{\mathbb{S}}(E \oplus U)^{(0)}(n;n)$ is the trivial tree. 
\end{Definition}

\begin{lemma}
We have that:
\begin{enumerate}
\item The category of dg $u\mathcal{O}$-algebras is equivalent to the category of dg unital partial operads.

\item The category of dg $u\mathcal{O}$-coalgebras is equivalent to the category of dg counital partial cooperad.
\end{enumerate}
\end{lemma}

\begin{proof}
It follows immediately from Proposition \ref{Prop: O-algebres} and Proposition \ref{Prop: O-coalgebras}.
\end{proof}

\begin{Definition}[Inhomogeneous quadratic presentation]
Let $(\mathcal{G},\{\circ_i\},\eta)$ be a unital partial $\mathcal{S}$-colored operad and let 
\[
\mathcal{G} \cong \mathscr{T}_{\mathbb{S}}(V)/(S)
\]
be a presentation of $\mathcal{G}$, where $S \subset \{ \I_\mathbb{S} \oplus V \oplus \mathscr{T}_{\mathbb{S}}(V)^{(2)} \}$. The presentation is a \textit{inhomogeneous quadratic presentation} if the $\mathbb{S}$-color scheme of relations $S$ satisfies the following conditions.

\begin{enumerate}
\item That the space of generators is \textit{minimal}, that is $S \cap \{\I_{\mathbb{S}} \oplus V\} = \{0\}$.

\item That the space of relations is \textit{maximal}, that is $(S) \cap \{\I_{\mathbb{S}} \oplus V \oplus \mathscr{T}_{\mathbb{S}}(V)^{(2)} \} = R'$. 
\end{enumerate}
\end{Definition}

\begin{lemma}
The presentation of the unital partial $\mathbb{S}$-colored operad $u\mathcal{O}$ given in its definition is a \textit{inhomogeneous quadratic presentation}.
\end{lemma}

\begin{proof}
It is straightforward to check given the aforementioned presentation.
\end{proof}

\begin{Definition}[Curved partial $\mathbb{S}$-colored cooperad]\label{def: curved colored partial cooperad}
A \textit{curved cooperad} $(\mathcal{V},\{\Delta_i\},\epsilon,d_\mathcal{V},\Theta_\mathcal{V})$ amounts to the data of a partial pdg $\mathbb{S}$-colored operad $(\mathcal{V},\{\Delta_i\},\epsilon,d_\mathcal{V})$ and a morphism of pdg $\mathbb{S}$-color schemes $\Theta_\mathcal{V}: (\mathcal{V},d_\mathcal{V}) \longrightarrow (\I_\mathbb{S},0)$ of degree $-2$, such that the following diagram commutes: 
\[
\begin{tikzcd}[column sep=7.5pc,row sep=3pc]
\mathcal{V} \arrow[r,"\Delta_{(1)}"] \arrow[rrd,"d_\mathcal{V}^2", bend right =10]
&\mathcal{V} \circ_{(1)} \mathcal{V} \arrow[r,"(\mathrm{id}~ \circ ~ \Theta_\mathcal{V})~-~(\Theta_\mathcal{V}~ \circ_{(1)}~ \mathrm{id})~"] 
&(\mathcal{V} \circ \I_\mathbb{S}) \oplus (\I_\mathbb{S} \circ \mathcal{V}) \cong  \mathcal{V} \oplus \mathcal{V}  \arrow[d,"\mathrm{proj}"]\\
&
&\mathcal{V}~,
\end{tikzcd}
\]
where $\mathrm{proj}$ is given by $\mathrm{proj}(\mu,\nu) \coloneqq \mu + \nu$.
\end{Definition}

\begin{Remark}
The definitions developed in Section \ref{Section: Curved cooperads} generalize \textit{mutatis mutandis} to the $\mathbb{S}$-colored case. 
\end{Remark}

We extend the formalism of semi-augmented operads of \cite{HirshMilles12} to the $\mathbb{S}$-colored case in order to compute the Koszul dual conilpotent curved partial $\mathbb{S}$-colored cooperad of $u\mathcal{O}$. 

\begin{Definition}[Semi-augmented unital partial $\mathbb{S}$-colored operad]
A \textit{semi-augmented} unital partial $\mathbb{S}$-colored operad $(\mathcal{G},\{\circ_i\},\eta,\iota)$ is the data of a unital partial $\mathbb{S}$-colored operad $(\mathcal{G},\{\circ_i\}, \allowbreak \eta,d_\mathcal{G})$ together with a morphism of $\mathbb{S}$-color schemes $\iota: \mathcal{G} \longrightarrow \I_{\mathbb{S}}$ of degree $0$ such that $\iota \cdot \eta = \mathrm{id}_{\I_\mathbb{S}}~.$
\end{Definition}

\begin{Remark}
The unital partial $\mathbb{S}$-colored operad $u\mathcal{O}$ is canonically semi-augmented by the identity morphism of $\I_\mathbb{S}$, we denote this semi-augmentation by $\iota_{u\mathcal{O}}$.
\end{Remark}

\begin{Proposition}\label{Proposition: S-colored inhomogeneous Koszul duality}
Let $(\mathcal{G},\{\circ_i\},\eta,\iota)$ be a semi-augmented unital partial $\mathcal{S}$-colored operad that admits a inhomogeneous quadratic presentation $(V,S)$. Let $qS \coloneqq S \cap \mathscr{T}_{\mathbb{S}}(V)^{(2)}$. Let $\varphi: qS \longrightarrow \I_\mathbb{S} \oplus V$ be the linear map that gives $S$ as its graph. The Koszul dual conilpotent partial $\mathbb{S}$-colored cooperad of $\mathcal{G}$ is given by 
\[
\mathcal{G}^{\ac} \coloneqq \overline{\mathscr{T}}_{\mathbb{S}}^c(sV, s^2qR)~.
\]
It is endowed with a coderivation of degree $-1$ $d_{\mathcal{G}^{\ac}}$ given by the unique extension of:
\[
\begin{tikzcd}[column sep=4pc,row sep=1pc]
\mathcal{G}^{\ac} \arrow[r,twoheadrightarrow]
&s^2qR \arrow[r,"s^{-1} \varphi_1 "]
&sE~,
\end{tikzcd}
\]
and with a curvature $\Theta_{\mathcal{G}^{\ac}}$ given by the degree $-2$ map: 
\[
\begin{tikzcd}[column sep=4pc,row sep=1pc]
\mathcal{G}^{\ac} \arrow[r,twoheadrightarrow]
&s^2qR \arrow[r,"s^{-2} \varphi_0 "]
&\I_{\mathbb{S}}~.
\end{tikzcd}
\]
The data of $(\mathcal{G}^{\ac}, d_{\mathcal{G}^{\ac}}, \Theta_{\mathcal{G}^{\ac}})$ forms a conilpotent curved partial $\mathbb{S}$-colored cooperad. 
\end{Proposition}

\begin{proof}
The proof is completely analogous to the non-$\mathbb{S}$-colored case developed in \cite[Section 4]{HirshMilles12}.
\end{proof}

\begin{Definition}[Koszul unital partial $\mathbb{S}$-colored operad]
Let $(\mathcal{G},\{\circ_i\},\eta,\iota)$ be a semi-augmented unital partial $\mathcal{S}$-colored operad that admits a inhomogeneous quadratic presentation $(V,S)$. Let $qS = S \cap \mathscr{T}_{\mathbb{S}}(V)^{(2)}$ and let 
\[
q\mathcal{G} \coloneqq \overline{\mathscr{T}}(V)/(qS)
\]
be the quadratic partial $\mathbb{S}$-colored operad associated to $\mathcal{G}$. The semi-augmented unital partial $\mathcal{S}$-colored operad $\mathcal{G}$ is said to be \textit{Koszul} if the quadratic operad $q\mathcal{G}$, endowed with the quadratic presentation $(V,qS)$, is a Koszul quadratic partial $\mathbb{S}$-colored in the sense of Definition \ref{def: Koszul quadratic S-colored operad}.
\end{Definition}

Let $J$ be the graded $\mathbb{S}$-color scheme given by $J(0;1) \coloneqq \kk.\theta$, where $\theta$ is an element of degree $-2$, and zero elsewhere.

\begin{Definition}[Curved partial $\mathbb{S}$-colored $c\mathcal{O}^\vee$]\label{def: cO vee}
Let $(E,R)$ be the quadratic presentation of the partial $\mathbb{S}$-colored operad $\mathcal{O}$. The \textit{curved partial} $\mathbb{S}$-\textit{colored cooperad} $c\mathcal{O}^\vee$ is given by the presentation 
\[
c\mathcal{O}^\vee \coloneqq \overline{\mathscr{T}}^c_{\mathbb{S}}(E \oplus J,R)~,
\]
where $E \oplus J$ are the cogenerators and $R$ are the corelations. It is endowed with the curvature $\Theta_{c\mathcal{O}^\vee}: c\mathcal{O}^\vee \longrightarrow \I_{\mathbb{S}}$ defined on $c\mathcal{O}^\vee(n;n)$ by the following map
\[
\left\{ \begin{tikzcd}[column sep=1.5pc,row sep=0pc]
\gamma_1^{1,n} \circ_1 \theta - \sum_{i=0}^n \gamma_i^{n,1} \circ_2 \theta \arrow[r,mapsto]
&\mathrm{id}_n ~,\\
\mu \arrow[r,mapsto]
&0~, \\
\end{tikzcd}
\right.
\]
if $\mu$ is not contained in the sub-$\mathbb{S}_n$-module generated by $\gamma_1^{1,n} \circ_1 \theta - \sum_{i=0}^n \gamma_i^{n,1} \circ_2 \theta$.
\end{Definition}

\begin{theorem}\label{Koszulity of uO}
Let $u\mathcal{O}$ be the unital partial $\mathbb{S}$-colored operad encoding (co)unital partial (co)operads. 

\begin{enumerate}
\item The Koszul dual conilpotent curved partial $\mathbb{S}$-colored cooperad $(u\mathcal{O})^{\ac}$ is isomorphic to the suspension of $c\mathcal{O}^\vee$.

\vspace{0.5pc}

\item The unital partial $\mathbb{S}$-colored operad $u\mathcal{O}$ is a Koszul, meaning that $q(u\mathcal{O})$ is a Koszul quadratic partial $\mathbb{S}$-colored operad.
\end{enumerate}
\end{theorem}

\begin{proof}
The Koszul dual cooperad $(u\mathcal{O})^{\ac}$ is given by $\overline{\mathscr{T}}^c_{\mathbb{S}}(s(E \oplus U), s^2 qR')$. It is cogenerated by a degree $1$ operation $u \in sU(0;1)$ and the suspension of the family $\{\gamma_i^{n,k}\}$. The corelations $s^2 qR'$ are given by $s^{2}R$. 

\medskip

Let us compute the rest of the structure. The projection of $R'$ into $E \oplus U$ is zero, hence $(u\mathcal{O})^{\ac}$ has a zero pre-differential. The projection of $R'$ into $\I_{\mathbb{S}}$ will be non zero only on the relations that involve the element $u \in U(0;1)$. It is straightforward to compute that the pre-image of $\mathrm{id}_n \in \I_{\mathbb{S}}$ by the projection $s^{-2} \varphi_0$ is given by 
\[
\sum_{i=0}^n s\gamma_{i}^{n,1} \circ_2 su - s\gamma_{1}^{1,n} \circ_1 su~.
\]
This completely characterizes the curvature of $(u\mathcal{O})^{\ac}$. A direct inspection identifies $(u\mathcal{O})^{\ac}$ with the operadic suspension of $c\mathcal{O}^\vee$.

\medskip

Let us prove the second assertion. The space $qR'$ is in fact $R$, and therefore $q(u\mathcal{O})$ is given by the coproduct $\mathcal{O} \coprod \{u\}$. By Theorem \ref{the Koszulity of O}, we know that $\mathcal{O}$ is a quadratic Koszul $\mathbb{S}$-colored operad. Therefore the Koszul complex of $q(u\mathcal{O})$ is also acyclic by an argument analogue to that of \cite[Proposition 6.16]{HirshMilles12}.
\end{proof}

The last point of this section is to understand what curved coalgebras and complete curved algebras over the conilpotent curved partial $\mathbb{S}$-colored cooperad $c\mathcal{O}^\vee$ are. 

\begin{Proposition}\label{Prop: cO-cogebres.}
The category of curved $(c\mathcal{O}^\vee)^u$-coalgebras is isomorphic to the category of conilpotent curved partial cooperads.
\end{Proposition}

\begin{proof}
First, notice that, for any pdg $\mathbb{S}$-module $M$, there is an isomorphism of pdg $\mathbb{S}$-modules 
\[
\mathscr{S}_{\mathbb{S}}(c\mathcal{O}^\vee)(M) \cong \overline{\mathscr{T}}^c(M \oplus \nu), 
\]
where $\nu$ is an arity $1$ and degree $-2$ generator. Indeed, this can by shown by extending the bijection given in the proof of Lemma \ref{lemma with the trees bijection}. This extension is defined by sending $\theta$ in $c\mathcal{O}^\vee(0;1)$ to the generator $\nu$. Pictorially it is given by
\begin{center}
\includegraphics[width=110mm,scale=1.3]{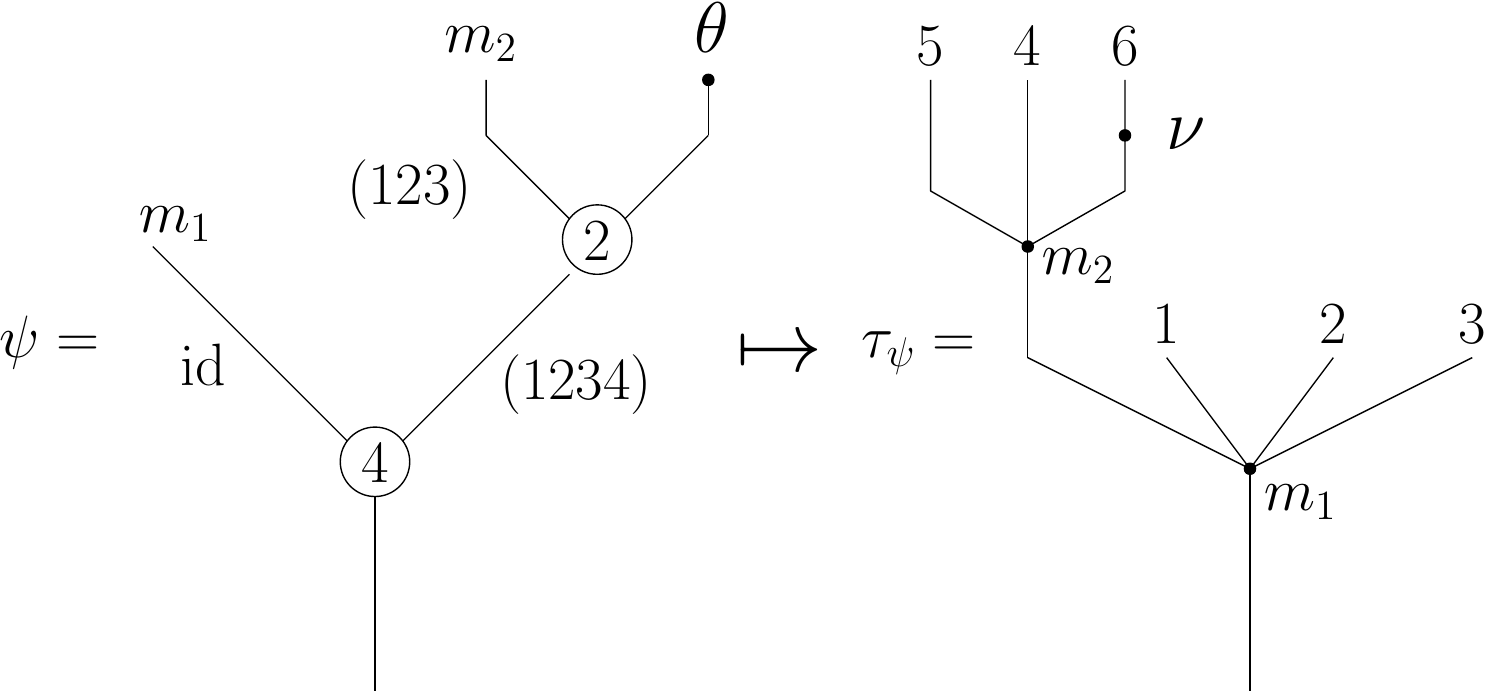}.
\end{center}
This isomorphism is natural in $M$. Therefore the endofunctor $\overline{\mathscr{T}}^c(- \oplus \nu)$ in the category of pdg $\mathbb{S}$-modules as a comonad structure. A direct computation shows that this comonad structure coincides with the reduced tree comonad structure. 

\medskip

Let $(\C,\{\Delta_i\},d_\C,\Theta_\C)$ be a conilpotent curved partial cooperad. Let 
\[
\Delta_\C: \C \longrightarrow \overline{\mathscr{T}}^c(\C)
\]
be its structure map as a coalgebra over the reduced tree comonad. We construct an extension
\[
\Delta_\C^+: \C \longrightarrow \overline{\mathscr{T}}^c(\C \oplus \nu)
\]
using $\Theta_\C$. The data $(\C,\Delta_\C^+,d_\C)$ forms a pdg $(c\mathcal{O}^\vee)^u$-coalgebra. It is in fact a curved $(c\mathcal{O}^\vee)^u$-coalgebra. Indeed, the diagram 
\[
\begin{tikzcd}[column sep=3pc,row sep=3pc]
\C \arrow[r,"\Delta_\C^+ "] \arrow[rd,"d_\C^2",swap]
&\mathscr{S}_\mathbb{S}(c\mathcal{O}^\vee)(\C) \arrow[d,"\mathscr{S}_\mathbb{S}(\Theta_{c\mathcal{O}^\vee})(\mathrm{id})"]\\
&\C \cong \mathscr{S}_\mathbb{S}(\I_\mathbb{S})(\C)
\end{tikzcd}
\]
commutes since $(\C,\{\Delta_i\},d_\C,\Theta_\C)$ forms a \textit{curved} partial cooperad. The other way around, given a curved $(c\mathcal{O}^\vee)^u$-coalgebra $(\C, \Delta_C^+, d_\C)$, one can compose the structural map
\[
\Delta_\C^+: \C \longrightarrow \overline{\mathscr{T}}^c(\C \oplus \nu)
\]
with the projection $\overline{\mathscr{T}}^c(\C \oplus \nu) \twoheadrightarrow \overline{\mathscr{T}}^c(\C)$, which endows $\C$ with a conilpotent partial cooperad structure $\{\Delta_i\}$. By composing $\Delta_\C^+$ with the projection $\overline{\mathscr{T}}^c(\C \oplus \nu) \twoheadrightarrow \I.\nu$, one obtains a map 
\[
\Theta_\C: \C \longrightarrow \I
\]
of pdg $\mathbb{S}$-modules of degree $-2$. The data $(\C,\{\Delta_i\},d_\C,\Theta_\C)$ forms a conilpotent curved partial cooperad, since $(\C, \Delta_C^+, d_\C)$ is a \textit{curved} $(c\mathcal{O}^\vee)^u$-coalgebra.
\end{proof}

\begin{Proposition}\label{Prop: cO-algebres.}
The category of complete curved $(c\mathcal{O}^\vee)^u$-algebras is isomorphic to the category of complete curved absolute partial operads.
\end{Proposition}

\begin{proof}
For the proof of this statement, see Proposition \ref{Prop: cO algebres vrai}. For an explicit description of these objects, we refer to the Appendix \ref{Appendix B}.
\end{proof}

\medskip

\section{Curved twisting morphisms and Bar-Cobar adjunctions at the operadic level}\label{Section: Constructions Bar-Cobar operadiques}
The curved Koszul duality established in the previous section at the groupoid-colored level gives a \textit{curved twisting morphism} $\kappa$ between the groupoid-colored cooperad encoding curved partial (co)operads $c\mathcal{O}^\vee$ and the groupoid-colored operad encoding unital partial (co)operads $u\mathcal{O}$. This curved twisting morphism induces two different Bar-Cobar adjunctions.

\medskip

The first Bar-Cobar adjunction induced by $\kappa$ is between dg unital partial operads and conilpotent curved partial cooperads. It will be shown to be isomorphic to the Bar-Cobar adjunction defined in \cite[Section 4.1]{grignou2019}: 

\hspace{9.3pc}
\begin{tikzcd}[column sep=5pc,row sep=3pc]
            \mathsf{curv}~\mathsf{pCoop}^{\mathsf{conil}} \arrow[r, shift left=1.1ex, "\Omega"{name=F}] &\mathsf{dg}~\mathsf{upOp}~.  \arrow[l, shift left=.75ex, "\text{B}"{name=U}]
            \arrow[phantom, from=F, to=U, , "\dashv" rotate=-90]
\end{tikzcd}

On the other hand, using the generalization of \cite{grignoulejay18} to the groupoid-colored setting, we obtain a new complete Bar-Cobar adjunction:

\hspace{10.5pc}
\begin{tikzcd}[column sep=6pc,row sep=3pc]
            \mathsf{dg}~\mathsf{upCoop} \arrow[r, shift left=1.1ex, "\widehat{\Omega}"{name=F}] &\mathsf{curv}~\mathsf{abs.pOp}^{\mathsf{comp}}~, \arrow[l, shift left=.75ex, "\widehat{\text{B}}"{name=U}]
            \arrow[phantom, from=F, to=U, , "\dashv" rotate=-90]
\end{tikzcd}

between complete curved \textit{absolute} partial operads and counital partial cooperads. Constructing this second adjunction is the main goal of this section.

\medskip

\subsection{Curved twisting morphisms and curved pre-Lie algebras}
In order to defined curved twisting morphisms in the first place, we introduce a new type of structure, called \textit{curved pre-Lie algebras}. In this type of algebras, the curvature has to satisfy a special condition, called the \textit{left-nucleus condition}. This condition comes from the deformation theory. More precisely, it appears in \cite[Proposition 1.1, Chapter 4]{DSV18} for the following reason: one can show that a dg pre-Lie algebra is \textit{twistable} by a Maurer-Cartan element if and only if this element satisfies the extra condition of being left-nucleus.   

\begin{Definition}[Curved pre-Lie algebra]\label{curvedprelie}
A \textit{curved pre-Lie algebra} $(\mathfrak{g}, \{-,-\},d_{\mathfrak{g}},\vartheta)$ amounts to the data of a pre-Lie algebra $(\mathfrak{g},\{-,-\})$, a derivation $d_{\mathfrak{g}}$ with respect to $\{-,-\}$ of degree $-1$, and a morphism of pdg modules of degree $-2$ $\vartheta: \mathbb{K} \longrightarrow \mathfrak{g}$. The data of this morphism is equivalent to the data of an element $\vartheta(1) \coloneqq \vartheta$ in $\mathfrak{g}_{-2}$. They are subject to the following conditions.
\begin{enumerate}
    \item The element $\vartheta$ has to be \textit{left-nucleus}, that is, for all $\mu,\nu$ in $\mathfrak{g}$: 
    \[
    \{\vartheta,\{\mu,\nu\}\} = \{\{\vartheta,\mu\},\nu\}~.
    \]
    \item Moreover, for all $\mu$ in $\mathfrak{g}$:
    \[
    d_{\mathfrak{g}}(\mu) = \{\vartheta,\mu\} - \{\mu,\vartheta\}~.
    \]
    \item And finally, $d_{\mathfrak{g}}(\vartheta)=0$.
\end{enumerate}
\end{Definition}

In the same spirit as for curved Lie algebras and curved associative algebras, one can define a curved partial operad, $cp\mathcal{L}ie$, that encodes curved pre-Lie algebras. Let $H$ be the pdg $\mathbb{S}$-module $(\mathbb{K}.\vartheta, 0, \mathbb{K}[\mathbb{S}_2].\nu, 0, \cdots)$ endowed with the zero pre-differential.

\begin{Definition}[Curved operad encoding curved pre-Lie algebras]
The curved partial operad $cp\mathcal{L}ie$ is given by the presentation: 
\[
cp\mathcal{L}ie \coloneqq \overline{\mathscr{T}}(H)/(D)~,
\]
where $(D)$ is the operadic ideal generated by the following relations: 
\begin{enumerate}
\item The right pre-Lie relation, already present in the classical pre-Lie operad, given by:
\[
\includegraphics[width=120mm,scale=1]{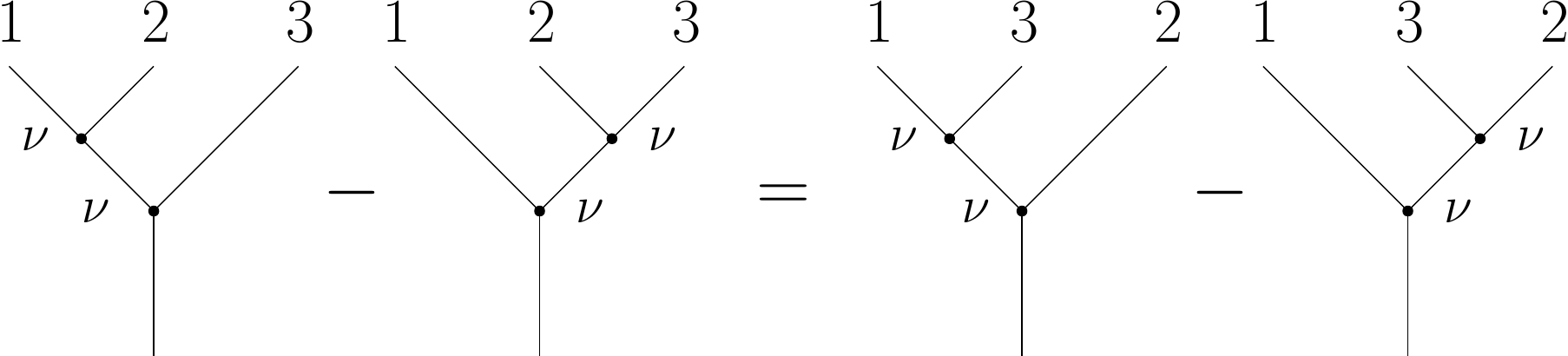}~.
\]
\item The left-nucleus relation stating for the curvature $\vartheta$, given by:
\[
\includegraphics[width=60mm,scale=1]{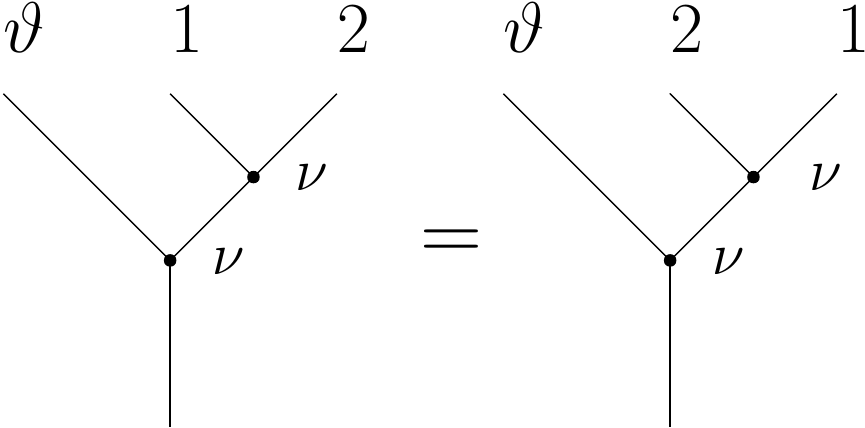}~.
\]
\end{enumerate}

It is endowed with the curvature $\Theta_{cp\mathcal{L}ie}(\mathrm{id}) \coloneqq \nu \circ_{1} \vartheta - \nu \circ_2 \vartheta~.$
\end{Definition}

\begin{Proposition}
The data $(cp\mathcal{L}ie, 0, \Theta_{cp\mathcal{L}ie})$ forms a curved partial operad. The category of curved $cp\mathcal{L}ie$-algebras is equivalent to the category of curved pre-Lie algebras. 

\medskip

Furthermore, the morphism of curved partial operads $ \Liec \longrightarrow \Assc$ given by the skew-symmetrization of the Lie bracket factors through the curved partial operad $cp\mathcal{L}ie$.
\end{Proposition}

\begin{proof}
Let us show that the data forms a curved partial operad. The proof is quite similar to that of Lemma \ref{lemmalie}. In order to check that $[\nu \circ_{1} \vartheta - \nu \circ_2 \vartheta,-] = 0$, since it is a derivation with respect partial compositions, it is enough to test this equality on the generators. First, $[\nu \circ_{1} \vartheta - \nu \circ_2 \vartheta, \vartheta] = 0$ is evident. Secondly, expanding $[\nu \circ_{1} \vartheta - \nu \circ_2 \vartheta, \nu]$ makes six terms that appear: it is straightforward to check that four of them cancel because of the right pre-Lie relation and the last two because of the left-nucleus relation. Proving that this curved partial operad encodes curved pre-Lie algebras is completely analogous to Lemma \ref{lemmalie}. One checks that the map $\Liec \longrightarrow cp\mathcal{L}ie$ given on generators by $\beta \mapsto \nu - \nu^{(12)}$ and $\zeta \mapsto \vartheta$ is indeed a morphism of curved partial operads. Furthermore, there is a morphism of curved partial operads $cp\mathcal{L}ie \longrightarrow \Assc$ simply given by $\nu \mapsto \mu$ and $\vartheta \mapsto \phi$. Their composition gives back the morphism constructed in Lemma \ref{assliem}.
\end{proof}

\begin{Remark}
Contrary to the case of Lie algebras and associative algebras, we had to introduce a new relation into the curved partial operad that encodes curved pre-Lie algebras. In the previous cases, this meant that the curved homotopy version of them was given simply by adding a curvature to the classical homotopy version. For instance, a curved $\mathcal{A}_{\infty}$-algebra is just an $\mathcal{A}_{\infty}$-algebra with an added structure of a curvature; the relations of a curved $\mathcal{A}_{\infty}$-algebra are clearly analogous to those of an $\mathcal{A}_{\infty}$-algebra. In the curved pre-Lie case, more structure will appear when one resolves the left-nucleus relation up to homotopy. Nevertheless, since the left-nucleus relation is not even quadratic, one would first need an expanded Koszul duality in order to treat this case.
\end{Remark}

We leave it to the reader to generalize the definition of curved partial operads to the $\mathbb{S}$-colored setting. (See Definition \ref{def: curved colored partial cooperad} for a similar definition).

\begin{lemma}[Totalization of a curved partial $\mathbb{S}$-colored operad]\label{lemma S colored totalization}
Let $(\mathcal{G}, \{\circ_i\}, d_{\mathcal{G}}, \Theta_\mathcal{G})$ be a curved partial $\mathbb{S}$-colored operad. The \textit{totalization} of $\mathcal{G}$ given by 
\[
\prod_{(n_1,\cdots,n_r;n) \in \mathbb{N}^{r+1}} \left(\mathcal{G}(n_1,\cdots,n_r;n)^{\left(\mathbb{S}_{n_1} \times \cdots \times \mathbb{S}_{n_r}\right)~\wr~ \mathbb{S}_r}\right)^{\mathbb{S}_n}
\]
forms a curved pre-Lie algebra, where the bracket is defined like in Proposition \ref{prop S colored totalization} and where the curvature is given by 
\[
\vartheta(1) = \sum_{n\geq 0} \theta_{n}~.
\]
\end{lemma}

\begin{proof}
Let $g$ be in $\mathcal{G}(n_1,\cdot,n_r;n)$, we have that $d_\mathcal{G}^2(g) = \theta_n \circ_1 g - \sum_{i=0}^r g \circ_i \theta_{n_i}$ since $\mathcal{G}$ is a curved partial $\mathbb{S}$-colored operad. The bracket $\vartheta(1) \star g = \theta_n \circ_1 g$ since only $\theta_n$ matches the appropriate color. Similarly, $g \star \vartheta(1) = \sum_{i=0}^r g \circ_i \theta_{n_i}$. Since each $\theta_{n}$ is an arity one operation, by the sequential and parallel axioms of a partial $\mathbb{S}$-colored operad, the curvature satisfies the \textit{left-nucleus} relation of a curved pre-Lie algebra.
\end{proof}

\begin{lemma}[Curved convolution partial $\mathbb{S}$-colored operad]\label{lemma curved S colored convolution}
Let $(\mathcal{G},\{\circ_i\},\eta,d_\mathcal{G})$ be a dg unital partial $\mathbb{S}$-colored operad and let $(\mathcal{V},\{\Delta_i\},d_\mathcal{V},\Theta_\mathcal{V})$ be a curved partial $\mathbb{S}$-colored cooperad. The convolution partial $\mathbb{S}$-colored operad $\mathcal{H}om(\mathcal{V},\mathcal{G})$ forms a curved partial $\mathbb{S}$-colored operad endowed with the curvature:
\[
\begin{tikzcd}
\Theta_{\mathcal{H}om}: \mathcal{V} \arrow[r,"\Theta_\mathcal{V}"]
&\I_{\mathbb{S}} \arrow[r,"\eta"]
&\mathcal{G}~.
\end{tikzcd}
\]
\end{lemma}

\begin{proof}
Let $\alpha$ be in $\mathcal{H}om(\mathcal{V},\mathcal{G})(n_1,\cdots,n_r;n)$, we have that:
\[
\partial^2(\alpha) = - (-1)^{2.|\alpha|} \alpha \circ d_\mathcal{V}^2 = \Theta_{\mathcal{H}om} \circ_1 \alpha - \sum_{i=0}^r \alpha \circ_i \Theta_{\mathcal{H}om}~,
\]
since $d_\mathcal{G}^2 = 0$ and $d_\mathcal{V}^2 = (\mathrm{id} \circ_{(1)} \Theta_\mathcal{V} - \Theta_\mathcal{V} \circ \mathrm{id} ) \cdot \Delta_{(1)}~.$
\end{proof}

\begin{Definition}[Maurer-Cartan of a curved pre-Lie]
Let $(\mathfrak{g}, \{-,-\},d_{\mathfrak{g}},\vartheta)$ be a curved pre-Lie algebra. A \textit{Maurer-Cartan element} $\alpha$ is a degree $-1$ element that of $\mathfrak{g}$ that satisfies the following equation:
\[
d_{\mathfrak{g}}(\alpha) + \{\alpha,\alpha\} = \vartheta~ .
\]
\end{Definition}

\begin{Remark}
The set of Maurer-Cartan elements of a curved pre-Lie (or similarly a curved Lie) algebra can be empty. A curvature term $\vartheta \neq 0$ stops $0$ in $\mathfrak{g}_{-1}$ from being a canonical Maurer-Cartan element in $\mathfrak{g}$.
\end{Remark}

\begin{Definition}[Curved $\mathbb{S}$-colored twisting morphism]
Let $(\mathcal{G},\{\circ_i\},\eta,d_\mathcal{G})$ be a unital partial dg $\mathbb{S}$-colored operad and let $(\mathcal{V},\{\Delta_i\},d_\mathcal{V},\Theta)$ be a curved partial $\mathcal{S}$-colored cooperad. A \textit{curved twisting morphism} $\alpha$ between $\mathcal{V}$ and $\mathcal{G}$ is a morphism of pdg $\mathbb{S}$-color schemes $\alpha: \mathcal{V} \longrightarrow \mathcal{G}$ of degree $-1$ that satisfies the Maurer-Cartan equation in the totalization of the curved convolution $\mathbb{S}$-colored operad $\mathcal{H}om(\mathcal{V},\mathcal{G})$. Otherwise stated, $\alpha$ satisfies:
\[
\partial(\alpha) + \alpha \star \alpha = \Theta_{\mathcal{H}om}~.
\]
\end{Definition}

\begin{lemma}\label{lemma: kappa is a curved twisting morphism}
The $\mathbb{S}$-color scheme morphism $\kappa: \mathcal{S}\otimes(c\mathcal{O})^\vee \longrightarrow u\mathcal{O}$ given by
\[
\kappa: \mathcal{S}\otimes(c\mathcal{O}^\vee) \twoheadrightarrow sE \oplus sJ \cong E \oplus U \hookrightarrow u\mathcal{O}
\]
is a curved twisting morphism. Here $\mathcal{S}$ denotes the operadic suspension of $c\mathcal{O}^\vee$.
\end{lemma}

\begin{proof}
Notice that $\partial(\kappa) = 0$ since the pre-differentials are null. Let us show that
\[
\kappa \star \kappa = \Theta_{\mathcal{H}om}~.
\]
The morphism of $\mathbb{S}$-color schemes $\kappa \star \kappa$ is non-zero only on elements of weight two in $c\mathcal{O}^\vee$. These are elements which can informally be written as  $\gamma_i^{n+m-1,k} \circ_1 \gamma_j^{n,m}$ and $\gamma_i^{n,k+m-1} \circ_2 \gamma_j^{k,m}$, or as $\gamma_i^{1,n} \circ_1 \theta$ and $\gamma_i^{n,1} \circ_2 \theta$. One computes that $\kappa \star \kappa$ of the first kind of weight two elements is zero because of the Koszul signs, by using the sequential and parallel relations. On the second kind of weight two elements, $\kappa \star \kappa$ is equal to $\gamma_i^{1,n} \circ_1 u$ and $\gamma_i^{n,1} \circ_2 u$. Using the unital relation in $u\mathcal{O}$, these are both equal to $\mathrm{id}_n$, the operadic unit of the $\mathbb{S}$-colored operad $u\mathcal{O}$. Hence $\kappa \star \kappa$ is equal to $\Theta_{\mathcal{H}om}$. 
\end{proof}

\begin{Notation}\label{not: le dirac S-colore}
Let $f: M \longrightarrow N$ be a morphism of graded $\mathbb{S}$-modules of degree $0$ and $g: M \longrightarrow N$ be a morphism of graded $\mathbb{S}$-modules degree $p$. We denote $\diracComb_{\mathbb{S}}(n_1,\cdots,n_r)(f,g)$ the map
\[
\sum_{i=1}^r f(n_1) \otimes \cdots \otimes f(n_{i-1}) \otimes g(n_i) \otimes f(n_{i+1}) \otimes \cdots \otimes f(n_r) : M(n_1) \otimes \cdots \otimes M(n_r) \longrightarrow N(n_1) \otimes \cdots \otimes N(n_r)
\]
which is an $\mathbb{S}_{n_1} \times \cdots \times \mathbb{S}_{n_r} \wr \mathbb{S}_r$-equivariant morphism of degree $p$. Let $E$ be a graded $\mathbb{S}$-color scheme. The family of maps $\{\diracComb_{\mathbb{S}}(n_1,\cdots,n_r)(f,g)\}$ induces a morphism of graded $\mathbb{S}$-modules of degree $p$:
\[
\mathscr{S}_{\mathbb{S}}(E)(M) \longrightarrow \mathscr{S}_{\mathbb{S}}(E)(N)
\]
by applying $\mathrm{id}_E \otimes \diracComb_{\mathbb{S}}(n_1,\cdots,n_r)(f,g)$ to each component. By a slight abuse of notation, this morphism will be denoted by $\mathscr{S}_{\mathbb{S}}(\mathrm{id}_E)(\diracComb_{\mathbb{S}}(f,g))~.$ Likewise, it induces a morphism of graded $\mathbb{S}$-modules of degree $p$:

\[
\widehat{\mathscr{S}}_{\mathbb{S}}^c (E)(M) \longrightarrow \widehat{\mathscr{S}}_{\mathbb{S}}^c (E)(N)
\]
by applying $\mathrm{Hom}(id_E,\diracComb_{\mathbb{S}}(n_1,\cdots,n_r)(f,g))$ to each component. By a slight abuse of notation, this morphism will be denoted by $\widehat{\mathscr{S}}_{\mathbb{S}}^c(\mathrm{id}_E)(\diracComb_{\mathbb{S}}(f,g))~.$
\end{Notation}

\subsection{Classical Bar-Cobar adjunction relative to $\kappa$}\label{subsection: classical bar-cobar relative to kappa}
The first adjunction induced by $\kappa$ will be an adjunction between dg $u\mathcal{O}$-algebras and curved $\mathcal{S}\otimes(c\mathcal{O}^\vee)^{u}$-coalgebras. That is, between dg unital partial operads and shifted conilpotent curved partial cooperads.

\begin{Definition}[Bar-Cobar constructions relative to $\kappa$]\label{Bar Cobar constructions relative to kappa}
Using $\kappa$, one can define two functor:
\begin{enumerate}
\item Let $(\PP,\gamma_\PP,d_\PP)$ a dg $u\mathcal{O}$-algebra. Its \textit{Bar construction relative to} $\kappa$, denoted by $\mathrm{B}_\kappa\PP$, is given by the cofree $\mathcal{S}\otimes (c\mathcal{O}^\vee)^u$-coalgebra $\mathscr{S}_{\mathbb{S}}(\mathcal{S}\otimes (c\mathcal{O}^\vee)^u)(\PP)$. Its pre-differential $d_{\mathrm{bar}}$ is given by the sum of two terms $d_1$ and $d_2$. The first term is given by 
\[
d_1 \coloneqq \mathscr{S}_{\mathbb{S}}(\mathrm{id})(\diracComb_{\mathbb{S}}(\mathrm{id},d_\PP))~. 
\]
The second term $d_2$ is the unique coderivation extending:
\[
\begin{tikzcd}[column sep=4pc,row sep=3pc]
\mathscr{S}_{\mathbb{S}}(\mathcal{S}\otimes c\mathcal{O}^\vee)(\PP) \arrow[r,"\mathscr{S}_{\mathbb{S}}(\kappa)(\mathrm{id})"]
&\mathscr{S}_{\mathbb{S}}(u\mathcal{O})(\PP) \arrow[r,"\gamma_\PP"]
&\PP~.
\end{tikzcd}
\]
\item Let $(\C,\Delta_\C,d_\C)$ be a curved $\mathcal{S}\otimes c\mathcal{O}^\vee$-coalgebra. Its \textit{Cobar construction relative to} $\kappa$, denoted by $\Omega_\kappa \C$, is given by the free $u\mathcal{O}$-algebra $\mathscr{S}_{\mathbb{S}}(u\mathcal{O})(\C)$. Its differential $d_{\mathrm{cobar}}$ is the sum of two terms $d_1$ and $d_1$. The first term is given by 
\[
d_1 \coloneqq -\mathscr{S}_{\mathbb{S}}(\mathrm{id})(\diracComb_{\mathbb{S}}(\mathrm{id},d_\C))~.
\]
The second term $d_2$ is the unique derivation extending:
\[
\begin{tikzcd}[column sep=4pc,row sep=3pc]
\C \arrow[r,"\Delta_\C"]
&\mathscr{S}_{\mathbb{S}}(\mathcal{S}\otimes c\mathcal{O}^\vee)(\C) \arrow[r,"\mathscr{S}_{\mathbb{S}}(\kappa)(\mathrm{id})"]
&\mathscr{S}_{\mathbb{S}}(u\mathcal{O})(\C)~.
\end{tikzcd}
\]
\end{enumerate}
\end{Definition}

\begin{lemma}
There is an adjunction 

\hspace{7pc}
\begin{tikzcd}[column sep=5pc,row sep=3pc]
            \mathsf{curv}~\mathcal{S}\otimes (c\mathcal{O}^\vee)^u\text{-}\mathsf{coalg} \arrow[r, shift left=1.1ex, "\Omega_\kappa"{name=F}] &\mathsf{dg}~u\mathcal{O}\text{-}\mathsf{alg}~.  \arrow[l, shift left=.75ex, "\text{B}_\kappa"{name=U}]
            \arrow[phantom, from=F, to=U, , "\dashv" rotate=-90]
\end{tikzcd}

\end{lemma}

\begin{proof}
The proof is a minor generalization of a the adjunction induced by a curved twisting morphism at the level of algebras in \cite[Section 5]{HirshMilles12}.
\end{proof}

This adjunction relative to $\kappa$ is in fact isomorphic to the adjunction between conilpotent curved coaugmented cooperads and operads introduced in \cite[Section 4.1]{grignou2019}. Before proving this result, we briefly recall the definition of this adjunction. 

\begin{Definition}[{\cite[Definition 58]{grignou2019}}]
Let $(\PP,\{\circ_i\},\eta, d_\PP)$ be a dg unital partial operad. The \textit{Bar construction} $\mathrm{B}\PP$ of $\PP$ is given by:
\[
\mathrm{B}\PP \coloneqq \left( \overline{\mathscr{T}}^c(s\PP \oplus \nu), d_{\mathrm{bar}} = d_1 + d_2, \Theta_{\mathrm{bar}} \right)~,
\]
where $\overline{\mathscr{T}}^c(s\PP \oplus \nu)$ is the cofree conilpotent partial pdg cooperad generated by the dg $\mathbb{S}$-module $s\PP \oplus \vartheta$. Here $\nu$ is an arity $1$ and degree $-2$ generator. It is endowed the pre-differential $d_{\mathrm{bar}}$, given by the sum of $d_1$ and $d_2$. The term $d_1$ is the unique coderivation extending
\[
\begin{tikzcd}[column sep=4pc,row sep=3pc]
\overline{\mathscr{T}}^c(s\PP \oplus \nu)  \arrow[r, twoheadrightarrow]
&s\PP \arrow[r,"s d_{\PP}"]
&s\PP~.
\end{tikzcd} 
\]
The term $d_2$ comes from the \textit{structure} of a unital partial operad on $\PP$, it is given by the unique coderivation extending 
\[
\begin{tikzcd}[column sep=4pc,row sep=3pc]
\overline{\mathscr{T}}^c(s\PP \oplus \nu) \arrow[r, twoheadrightarrow]
&\I.\nu \oplus (s\PP \circ_{(1)} s\PP) \arrow[r,"s \eta + s^2\gamma_{(1)}"]
&s\PP~.
\end{tikzcd}
\]
It is also endowed with the following curvature:
\[
\begin{tikzcd}[column sep=4pc,row sep=3pc]
\Theta_{\mathrm{bar}}: \overline{\mathscr{T}}^c(s\PP \oplus \nu) \arrow[r,twoheadrightarrow] 
&\I.\nu \arrow[r,"s^{-2}"]
&\I~.
\end{tikzcd}
\]
The resulting Bar construction of $\PP$ forms a conilpotent curved partial cooperad. 
\end{Definition}

\begin{Definition}[{\cite[Definition 60]{grignou2019}}]\label{Def: Cobar de Brice}
Let $(\C,\{\Delta_i\},d_\C,\Theta_\C)$ be a curved partial cooperad. The \textit{Cobar construction} $\Omega \C$ of $\C$ is given by:
\[
\Omega \C \coloneqq \left(\mathscr{T}(s^{-1}\C), d_{\mathrm{cobar}} = d_1 - d_2 \right)~,
\]
where $\mathscr{T}(s^{-1}\C)$ is the free unital partial operad generated by the pdg $\mathbb{S}$-module $s^{-1}\C$. It is endowed with the differential $d_{\mathrm{cobar}}$ given by the difference of $d_1$ and $d_2$. The term $d_1$ is the unique derivation extending
\[
\begin{tikzcd}[column sep=5pc,row sep=3pc]
s^{-1}\C \arrow[r,"s^{-1}d_\C"]
&s^{-1}\C \arrow[r,rightarrowtail]
&\mathscr{T}(s^{-1}\C)~.
\end{tikzcd}
\]
The term $d_2$ comes from the \textit{structure} of a conilpotent curved partial cooperad on $\C$, it is given by the unique derivation extending
\[
\begin{tikzcd}[column sep=5pc,row sep=3pc]
s^{-1}\C \arrow[r,"s^{-2}\Delta_{(1)} + s^{-1}\Theta_\C"]
&(s^{-1}\C \circ_{(1)} s^{-1}\C) \oplus \I \arrow[r,rightarrowtail]
&\mathscr{T}(s^{-1}\C)~.
\end{tikzcd}
\]
The resulting Cobar construction of $\C$ forms a dg unital partial operad.
\end{Definition}

\begin{Remark}
The Cobar constructions $\Omega \C$ is not augmented because the canonical morphism $\mathscr{T}(s^{-1}\C) \twoheadrightarrow \I$ does not commute with the differentials in general. Indeed, $d_{\mathrm{cobar}}(\Theta_\C(\mathrm{id}))$ is the trivial tree $|$. It is therefore augmented if and only if the curvature $\Theta_\C (\mathrm{id})$ is zero.
\end{Remark}

The Bar-Cobar constructions described above also form an adjunction 
\[
\begin{tikzcd}[column sep=5pc,row sep=3pc]
            \mathsf{curv}~\mathsf{pCoop}^{\mathsf{conil}} \arrow[r, shift left=1.1ex, "\Omega"{name=F}] &\mathsf{dg}~\mathsf{upOp}~.  \arrow[l, shift left=.75ex, "\mathrm{B}"{name=U}]
            \arrow[phantom, from=F, to=U, , "\dashv" rotate=-90]
\end{tikzcd}
\]
\begin{Proposition}[Mise en abîme]\label{Mise en abime}
The Bar-Cobar adjunction relative to $\kappa$ given by $\Omega_\kappa \dashv \mathrm{B}_\kappa$ is naturally isomorphic to the Bar-Cobar adjunction $\Omega \dashv \mathrm{B}$ constructed in \cite{grignou2021}.
\end{Proposition}

\begin{proof}
Let $(\PP,\{\circ_i\},\eta, d_\PP)$ be a dg unital partial operad. One uses the bijection given in the proof of Proposition \ref{Prop: cO-cogebres.} to construct a natural isomorphism $\mathrm{B} \PP \cong \mathrm{B}_\kappa \PP~.$
\end{proof}

\begin{Corollary}
This gives another proof that the Bar-Cobar adjunction $\Omega \dashv \mathrm{B}$ constructed in \cite{grignou2019} is a Quillen equivalence.
\end{Corollary}

The notion of curved twisting morphism between conilpotent curved partial cooperad and dg unital partial operad is again encoded by a curved pre-Lie algebra.

\begin{lemma}
Let $(\PP,\{\circ_i\},d_\PP,\Theta_\PP)$ be a curved partial operad. The totalization of $\PP$ given by 
\[
\prod_{n \geq 0} \PP(n)^{\mathbb{S}_n}~,
\]
together with its pre-Lie bracket and endowed with the curvature $\vartheta(1) \coloneqq \Theta_\PP(\mathrm{id})$ forms a curved pre-Lie algebra.
\end{lemma}

\begin{proof}
The proof is completely analogous to Lemma \ref{lemma S colored totalization}.
\end{proof}

\begin{lemma}[Curved convolution operad]\label{def: convolution curved partial operad}
Let $(\PP,\{\circ_i\},\eta,d_\PP)$ be a dg unital partial operad and let $(\C,\{\Delta_i\},d_\C,\Theta)$ be a curved partial cooperad. The convolution partial pdg operad $(\mathcal{H}om(\C,\PP),\{\circ_i\},\partial)$ forms a curved partial operad endowed with the curvature given by
\[
\begin{tikzcd}
\Theta_{\mathcal{H}om}(\mathrm{id}): \mathcal{C} \arrow[r,"\Theta"]
&\I \arrow[r,"\eta"]
&\mathcal{P}~.
\end{tikzcd}
\]
\end{lemma}

\begin{proof}
The proof is completely analogous to Lemma \ref{lemma curved S colored convolution}.
\end{proof}

\begin{Definition}[Curved twisting morphism]\label{def: curved twisting entre P et curved C}
Let $(\PP,\{\circ_i\},\eta,d_\PP)$ be a dg unital partial operad and let $(\C,\{\Delta_i\},d_\C,\Theta)$ be a curved partial cooperad. A \textit{curved twisting morphism} $\alpha$ is a Maurer-Cartan element curved pre-Lie algebra given by the curved convolution operad:
\[
\mathfrak{g}_{\C,\PP} \coloneqq \prod_{n\geq 0}\mathrm{Hom}_{\mathbb{S}}(\C(n),\PP(n))~.
\]
This is the data of a  morphism of $\mathbb{S}$-modules $\alpha: \C \longrightarrow \PP$ of a degree $-1$ such that:
\[
\partial(\alpha) + \alpha \star \alpha = \Theta_{\mathcal{H}om}(\mathrm{id})~.
\]
The set of twisting morphism between $\C$ and $\PP$ will be denoted $\mathrm{Tw}(\C,\PP)$. 
\end{Definition}

The set of curved twisting morphism between conilpotent curved partial cooperads and unital partial dg operads defines a bifunctor
\[
\begin{tikzcd}
\mathrm{Tw}(- , -) : \left(\mathsf{curv}~\mathsf{pCoop}^{\mathsf{conil}}\right)^{\text{op}} \times \mathsf{dg}~\mathsf{upOp}  \arrow[r]
&\mathsf{Set}~,
\end{tikzcd}
\]
which is represented on both sides by the Bar-Cobar construction defined before.

\begin{Proposition}[{\cite[Proposition 63]{grignou2019}}]
Let $(\PP,\{\circ_i\},\eta,d_\PP)$ be a dg unital partial operad and let $(\C,\{\Delta_i\},d_\C,\Theta)$ be a conilpotent curved partial cooperad. There are isomorphism:
\[
\mathrm{Hom}_{\mathsf{dg}~\mathsf{upOp}}(\Omega\C,\PP) \cong \mathrm{Tw}(\C,\PP) \cong \mathrm{Hom}_{\mathsf{curv}~\mathsf{pCoop}^{\mathsf{conil}}}(\C,\mathrm{B}\PP)~,
\]
which are natural in $\C$ and $\PP$.
\end{Proposition}

\subsection{Complete Bar-Cobar adjunction relative to $\kappa$}
We introduce a new adjunction, which we call the \textit{complete Bar-Cobar adjunction}, between complete curved absolute partial operads and dg counital partial cooperads. This adjunction is again induced by the curved twisting morphism $\kappa$, using the techniques developed in \cite{grignoulejay18}. In order to construct it, we need a mild generalization of the results in \cite{anelcofree2014}, in which the author only considers set-colored operads. 

\begin{theorem}\label{thm: existence de la cogèbre colibre}
Let $(\mathcal{G},\gamma,\eta,d_\G)$ be a dg $\mathbb{S}$-colored operad. The category of dg $\mathcal{G}$-coalgebras is comonadic, that is, there exists a comonad $(\mathscr{C}(\G),\omega,\xi)$ in the category of dg $\mathbb{S}$-modules such that the category of dg $\mathscr{C}(\G)$-coalgebras is equivalent to the category of dg $\mathcal{G}$-coalgebras.

\medskip

The endofunctor $\mathscr{C}(\G)$ is given by the following pullback:
\[
\begin{tikzcd}[column sep=3pc,row sep=3pc]
\mathscr{C}(\G) \arrow[r,"p_2"] \arrow[d,"p_1",swap] \arrow[dr, phantom, "\ulcorner", very near start]
&\widehat{\mathscr{S}}_\mathbb{S}^c(\G) \circ \widehat{\mathscr{S}}_\mathbb{S}^c(\G) \arrow[d,"\varphi_{\G,\G}"] \\
\widehat{\mathscr{S}}_\mathbb{S}^c(\G) \arrow[r,"\widehat{\mathscr{S}}_\mathbb{S}^c(\gamma)"]
&\widehat{\mathscr{S}}_\mathbb{S}^c(\G \circ \G)
\end{tikzcd}
\]

in the category $\mathsf{End}(\mathsf{dg}~\smod)$. It gives a pullback 
\[
\begin{tikzcd}[column sep=3pc,row sep=3pc]
\mathscr{C}(\G)(M) \arrow[r,"p_2(M)"] \arrow[d,"p_1(M)",swap] \arrow[dr, phantom, "\ulcorner", very near start]
&\widehat{\mathscr{S}}_\mathbb{S}^c(\G)(M) \circ \widehat{\mathscr{S}}_\mathbb{S}^c(\G)(M) \arrow[d,"\varphi_{\G,\G}(M)"] \\
\widehat{\mathscr{S}}_\mathbb{S}^c(\G)(M) \arrow[r,"\widehat{\mathscr{S}}_\mathbb{S}^c(\gamma)(M)"]
&\widehat{\mathscr{S}}_\mathbb{S}^c(\G \circ \G)(M)
\end{tikzcd}
\]
for any dg $\mathbb{S}$-module $M$. Here $p_1$ is a monomorphism since $\varphi_{\G,\G}(M)$ is also a monomorphism. The structural map of the comonad $\omega(M): \mathscr{C}(\G)(M) \longrightarrow \mathscr{C}(\G) \circ \mathscr{C}(\G)(M)$ is given by the map $p_2$ in the previous pullback, as its image lies in $\mathscr{C}(\G) \circ \mathscr{C}(\G)(M)$ instead of $\widehat{\mathscr{S}}_\mathbb{S}^c(\G) \circ \widehat{\mathscr{S}}_\mathbb{S}^c(\G)(M)$. The counit is given of the comonad $\mathscr{C}(\G)$ is given by
\[
\begin{tikzcd}[column sep=4pc,row sep=0.5pc]
\xi(M): \mathscr{C}(\G)(M) \arrow[r,"p_1(M)"] 
&\widehat{\mathscr{S}}_\mathbb{S}^c(\G)(M) \arrow[r,"\widehat{\mathscr{S}}_\mathbb{S}^c(\eta)(id_M)"]
&M~.
\end{tikzcd}
\]
\end{theorem}

\begin{proof}
In order to prove the first result, one has to generalize the framework of \cite{anelcofree2014} from set-colored operads to groupoid-colored operads. This is purely formal: the underlying theory upon which the result is obtained is the formalism of analytic functors and operads developed in \cite{Gambino_2017}. In \textit{loc.cit}, the formalism is already developed in the case where the colors form a category. \textit{Mutatis mutandis}, one can perform the same proofs using this general formalism. Thus the category of dg $\mathcal{G}$-coalgebras is indeed comonadic, by \cite[Theorem 2.7.11]{anelcofree2014}.

\medskip

The second statement is equivalent to claiming that the infinite recursion that defines the comonad $\mathscr{C}(\G)$ stops at the first step. We know this is the case for operads in the category of dg modules. We have to check that the lemmas proven for the category of dg modules also hold for the category of dg $\mathbb{S}$-modules. This comes from the fact that the base field $\kk$ is of characteristic $0$, and hence for any $n$, the $\kk$-algebra $\kk[\mathbb{S}_n]$ is a semi-simple finite dimensional $\kk$-algebra where every module is injective and projective. One can check that the different lemmas of \cite[Section 3]{anelcofree2014} that make the proof work are also true in this context.
\end{proof}

Before defining the complete Bar construction relative to $\kappa$ using the cofree dg-$u\mathcal{O}$-coalgebra $\mathscr{C}(u\mathcal{O})$ construction, one needs to understand coderivations on this object. We generalize the results of \cite[Section 6]{grignoulejay18} to our framework. Let $M$ be a graded $\mathbb{S}$-module, we denote by $\pi_M$ the map given by
\[
\begin{tikzcd}[column sep=5pc,row sep=3pc]
\pi_M: \mathscr{C}(u\mathcal{O})(M)  \arrow[r,"p_1"]
&\widehat{\mathscr{S}}_\mathbb{S}^c(u\mathcal{O})(M) \arrow[r,"\widehat{\mathscr{S}}_\mathbb{S}^c(\eta_{u\mathcal{O}})(\mathrm{id})"]
&M~,
\end{tikzcd} 
\]
where $\eta_{u\mathcal{O}}: \I_\mathbb{S} \longrightarrow u\mathcal{O}$ is the unit of the unital partial $\mathbb{S}$-colored operad $u\mathcal{O}$.

\begin{lemma}\label{lemma: coderivations on the cofree cooperad}
Let $f: M \longrightarrow N$ be a morphism graded $\mathbb{S}$-modules of degree $0$ and $g:M \longrightarrow N$ be a morphism of graded $\mathbb{S}$-modules of degree $p$, then the degree $p$ map: 
\[
\begin{tikzcd}[column sep=4pc,row sep=1pc]
&\widehat{\mathscr{S}}_{\mathbb{S}}^c(\mathrm{id}_{u\mathcal{O}})(\diracComb_{\mathbb{S}}(f,g)): \widehat{\mathscr{S}}_{\mathbb{S}}^c(u\mathcal{O})(M) \arrow[r]
&\widehat{\mathscr{S}}_{\mathbb{S}}^c(u\mathcal{O})(N)
\end{tikzcd}
\]
restricts to a degree $p$ morphism $\mathscr{C}(\mathrm{id})(\diracComb_{\mathbb{S}}(f,g)): \mathscr{C}(u\mathcal{O})(M) \longrightarrow \mathscr{C}(u\mathcal{O})(N)$.
\end{lemma}

\begin{proof}
The proof is essentially the same as \cite[Lemma 6.20]{grignoulejay18}. It is a consequence of the universal property of pullbacks: it suffices to construct a morphism of graded $\mathbb{S}$-modules 
\[
\mathscr{C}(u\mathcal{O})(M) \longrightarrow \widehat{\mathscr{S}}_{\mathbb{S}}^c(u\mathcal{O}) \circ \widehat{\mathscr{S}}_{\mathbb{S}}^c(u\mathcal{O})(N)
\]
that coincides with the morphism $\widehat{\mathscr{S}}_{\mathbb{S}}^c(u\mathcal{O})(\diracComb_{\mathbb{S}}(f,g))$ on $\widehat{\mathscr{S}}_{\mathbb{S}}^c(u\mathcal{O} \circ u\mathcal{O})(N)$. One can check that this morphism is given by 
\[
\widehat{\mathscr{S}}_{\mathbb{S}}^c(\mathrm{id}) \left(\diracComb\left(\widehat{\mathscr{S}}_{\mathbb{S}}^c(\mathrm{id})(f), \allowbreak \widehat{\mathscr{S}}_{\mathbb{S}}^c(\mathrm{id})(\diracComb(f,g))\right)\right) \circ p_2~.\]
\end{proof}

\begin{Proposition}[Coderivations on the cofree construction]
Let $M$ be a graded $\mathbb{S}$-module and let $\mathscr{C}(u\mathcal{O})(M)$ be the cofree graded $u\mathcal{O}$-coalgebra generated by $M$. There is a natural bijection between maps 
\[
\varphi: \mathscr{C}(u\mathcal{O})(M) \longrightarrow M
\]
of degree $p$ and coderivations 
\[
d_\varphi: \mathscr{C}(u\mathcal{O})(M) \longrightarrow \mathscr{C}(u\mathcal{O})(M)
\]
of degree $p$. This bijection sends $\varphi$ to the degree $p$ coderivation given by 
\[
d_\varphi \coloneqq \mathscr{C}(\mathrm{id})(\diracComb(\pi_M,\varphi)) \cdot \omega(M)~,
\]
where $\omega(M)$ is the comonad structure map of $\mathscr{C}(u\mathcal{O})$.
\end{Proposition}

\begin{proof}
For any $\varphi$, the morphism $d_\varphi$ is an endomorphism of $\mathscr{C}(u\mathcal{O})(M)$ by Lemma \ref{lemma: coderivations on the cofree cooperad}. One has to check that it is a coderivation. The main point is to use the fact that $p_1$ is a monomorphism to be able to check it directly on the dual Schur functors: the same diagrams as in \cite[Proposition 6.23]{grignoulejay18} commute, thus $d_\varphi$ is a coderivation.
\end{proof}

We first present the formal constructions before giving a more explicit description of this adjunction.

\begin{Definition}[Complete Bar construction relative to $\kappa$]
Let $(\PP,\gamma_\PP, d_\PP)$ be a complete curved $\mathcal{S}\otimes c\mathcal{O}^\vee$-algebra. The \textit{complete Bar construction relative to} $\kappa$, denoted by $\widehat{\text{B}}_\kappa\PP$, is given by:
\[
\widehat{\mathrm{B}}_\kappa\PP \coloneqq \left(\mathscr{C}(u\mathcal{O})(\PP), d_{\mathrm{bar}} = d_1 + d_2\right)
\]
where $\mathscr{C}(u\mathcal{O})$ is the cofree dg $u\mathcal{O}$-coalgebra generated by the pdg $\mathbb{S}$-module $\PP$. It is endowed with the pre-differential $d_{\mathrm{bar}}$ given by the sum of two terms $d_1$ and $d_2$. The term $d_1$ is of the pre-differential is the unique coderivation extending the following map
\[
\begin{tikzcd}[column sep=5pc,row sep=3pc]
\phi_1: \mathscr{C}(u\mathcal{O})(\PP) \arrow[r,"\mathscr{C}(u\mathcal{O})(d_\PP)"]
&\mathscr{C}(u\mathcal{O})(\PP) \arrow[r,"\pi_\PP"]
&\PP~.
\end{tikzcd}
\]
It is given by 
\[
d_1 = \mathscr{C}(\mathrm{id})(\diracComb_{\mathbb{S}}(\mathrm{id},d_\mathcal{P}))~.
\]
The term $d_2$ is given by the unique coderivation that extends the map: 
\[
\begin{tikzcd}[column sep=4pc,row sep=3pc]
\phi_2: \mathscr{C}(u\mathcal{O})(\PP)  \arrow[r,"p_1"]
&\widehat{\mathscr{S}}_\mathbb{S}^c(u\mathcal{O})(\PP) \arrow[r,"\widehat{\mathscr{S}}_\mathbb{S}^c(\kappa)(\mathrm{id})"]
&\widehat{\mathscr{S}}_\mathbb{S}^c(\mathcal{S} \otimes c\mathcal{O}^\vee)(\PP) \arrow[r,"\gamma_\PP "]
&\PP~.
\end{tikzcd}
\]
\end{Definition}

\begin{Proposition}
Let $(\PP,\gamma_\PP, d_\PP)$ be a complete curved $\mathcal{S}\otimes c\mathcal{O}^\vee$-algebra. The complete Bar construction $\widehat{\mathrm{B}}_\kappa\PP$ forms a dg counital partial cooperad. Meaning that
\[
d_{\mathrm{bar}}^2=0~.
\]
\end{Proposition}

\begin{proof}
The proof is very similar to the proof of \cite[Proposition 9.2]{grignoulejay18}. In order to check that $d_{\mathrm{bar}}^2=0$, by a straightforward generalization of \cite[Proposition 6.25]{grignoulejay18} to the $\mathbb{S}$-colored case, it is enough to check that $(\phi_1 + \phi_2) \cdot d_{\mathrm{bar}} = 0$. We have that:
\[
d_{\mathrm{bar}} = \widehat{\mathscr{S}}_\mathbb{S}^c(\mathrm{id})(\diracComb(\mathrm{id},d_\PP)) + \widehat{\mathscr{S}}_\mathbb{S}^c(\mathrm{id}_{u\mathcal{O}}) \circ \widehat{\mathscr{S}}_\mathbb{S}^c(\diracComb(\pi_{\PP},\gamma_\PP \cdot \widehat{\mathscr{S}}_\mathbb{S}^c(\kappa)(\mathrm{id}_\PP))(\mathrm{id}_\PP) \cdot p_2(\PP)~,
\]
where the first term is $d_1$ and the second is $d_2$. One computes that: 
\[
(\phi_1 + \phi_2) \cdot d_{\mathrm{bar}} = \gamma_\PP \cdot \widehat{\mathscr{S}}_\mathbb{S}^c(\kappa \star \kappa - \Theta_{\mathcal{S} \otimes c\mathcal{O}^\vee} \cdot \eta_{u\mathcal{O}})(\mathrm{id}_\PP) = 0~. 
\]
since $\kappa$ is a curved twisting morphism. Therefore $d_{\mathrm{bar}}^2=0$ and $\widehat{\mathrm{B}}_\kappa \PP$ is a dg counital partial cooperad.
\end{proof}

Derivations on $\widehat{\mathscr{S}}_\mathbb{S}^c(\mathcal{S}\otimes c\mathcal{O}^\vee)(M)$, the free pdg $\mathcal{S}\otimes c\mathcal{O}^\vee$-algebra on an pdg $\mathbb{S}$-module $M$, are completely characterized by their restrictions to the generators $M$. See \cite[Proposition 6.15]{grignoulejay18} for an analogue statement, which holds \textit{mutatis mutandis} in our case. Let $M$ be an pdg $\mathbb{S}$-module, we denote $\iota_M$ be the following map:
\[
\begin{tikzcd}[column sep=5pc,row sep=3pc]
\iota_M: M \arrow[r,"\mathscr{S}_\mathbb{S}^c(\epsilon)(id_M)"]
&\mathscr{S}_\mathbb{S}^c(\mathcal{S}\otimes c\mathcal{O}^\vee)(M)~.
\end{tikzcd}
\]
\begin{Definition}[Complete Cobar construction relative to $\kappa$]
Let $(\C, \{\Delta_C\}, d_\C)$ be a dg $u\mathcal{O}$-coalgebra. The \textit{complete Cobar construction relative to} $\kappa$, denoted by $\widehat{\Omega}_\kappa \C$, is given by:
\[
\widehat{\Omega}\C \coloneqq \left(\widehat{\mathscr{S}}_\mathbb{S}^c(\mathcal{S}\otimes c\mathcal{O}^\vee)(\C), d_{\mathrm{cobar}} = d_1 - d_2\right)
\]

where $\widehat{\mathscr{S}}_\mathbb{S}^c(\mathcal{S}\otimes c\mathcal{O}^\vee)(\C)$ is the free pdg $\mathcal{S}\otimes c\mathcal{O}^\vee$-algebra. It is endowed with the pre-differential $d_{\mathrm{cobar}}$ which is the difference of two terms $d_1$ and $d_2$. The term $d_1$ is the unique derivation that extends the map:
\[
\begin{tikzcd}[column sep=5pc,row sep=3pc]
\psi_1: \C \arrow[r,"\iota_\C"]
&\mathscr{S}_\mathbb{S}^c(\mathcal{S}\otimes c\mathcal{O}^\vee)(\C) \arrow[r,"\mathscr{S}_\mathbb{S}^c(\mathrm{id})(d_\C)" ]
&\mathscr{S}_\mathbb{S}^c(\mathcal{S}\otimes c\mathcal{O}^\vee)(\C)~. 
\end{tikzcd}
\]
It is given by
\[
d_1 = \widehat{\mathscr{S}}_{\mathbb{S}}^c(\mathrm{id})(\diracComb(\mathrm{id},d_\C))~.
\]
The term $d_2$ is given by the unique derivation that extends the map: 
\[
\begin{tikzcd}[column sep=5pc,row sep=3pc]
\psi_2: \C \arrow[r,"\Delta_\C"]
&\mathscr{S}_\mathbb{S}^c(u\mathcal{O})(\C) \arrow[r,"\mathscr{S}_\mathbb{S}^c(\kappa)(\mathrm{id}_\C)" ]
&\mathscr{S}_\mathbb{S}^c(\mathcal{S}\otimes c\mathcal{O}^\vee)(\C)~. 
\end{tikzcd}
\]
\end{Definition}

\begin{Proposition}
Let $(\C, \{\Delta_C\}, d_\C)$ be a dg $u\mathcal{O}$-coalgebra. The complete Cobar construction relative to $\kappa$ $\widehat{\Omega}_\kappa\C$ forms a complete curved $\mathcal{S}\otimes c\mathcal{O}^\vee$-algebra. In other words, the following diagram commutes: 
\[
\begin{tikzcd}[column sep=9pc,row sep=3pc]
\widehat{\mathscr{S}}_\mathbb{S}^c(\I_\mathbb{S}) \circ \widehat{\mathscr{S}}_\mathbb{S}^c(\mathcal{S}\otimes c\mathcal{O}^\vee)(\C) \arrow[r,"\widehat{\mathscr{S}}_\mathbb{S}^c(\Theta_{\mathcal{S}\otimes c\mathcal{O}^\vee}) \circ \widehat{\mathscr{S}}_\mathbb{S}^c(\mathrm{id})(\mathrm{id}) "] \arrow[rd,"- d_{\mathrm{cobar}}^2",swap]
&\widehat{\mathscr{S}}_\mathbb{S}^c(\mathcal{S}\otimes c\mathcal{O}^\vee) \circ \widehat{\mathscr{S}}_\mathbb{S}^c(\mathcal{S}\otimes c\mathcal{O}^\vee)(\C)\arrow[d,"\widehat{\mathscr{S}}_\mathbb{S}^c(\Delta_{\mathcal{S}\otimes c\mathcal{O}^\vee}) \cdot \varphi_{\mathcal{S}\otimes c\mathcal{O}^\vee}"]\\
&\widehat{\mathscr{S}}_\mathbb{S}^c(\mathcal{S}\otimes c\mathcal{O}^\vee)(\C)~.
\end{tikzcd}
\]
\end{Proposition}

\begin{proof}
The proof is again very similar to the proof of \cite[Proposition 9.1]{grignoulejay18}. In order to check that 
\[
- d_{\mathrm{cobar}}^2 = \widehat{\mathscr{S}}_\mathbb{S}^c(\Delta_{\mathcal{S}\otimes c\mathcal{O}^\vee}) \cdot \varphi_{\mathcal{S}\otimes c\mathcal{O}^\vee} \cdot \widehat{\mathscr{S}}_\mathbb{S}^c(\Theta_{\mathcal{S}\otimes c\mathcal{O}^\vee}) \circ \widehat{\mathscr{S}}_\mathbb{S}^c(\mathrm{id})(\mathrm{id}) ~,
\]
it is in fact enough to check that 
\[
- d_{\mathrm{cobar}} \cdot (\psi_1 - \psi_2) = \widehat{\mathscr{S}}_\mathbb{S}^c(\Theta_{\mathcal{S}\otimes c\mathcal{O}^\vee})(\mathrm{id}_\C)~.
\]
This comes from a straightforward generalization of \cite[Proposition 7.4]{grignoulejay18}. The pre-differential $d_{\mathrm{cobar}}$ is given by
\[
d_{\mathrm{cobar}} = \widehat{\mathscr{S}}_\mathbb{S}^c(\mathrm{id})(\diracComb(\iota_{\C},d_\C)) - \widehat{\mathscr{S}}_\mathbb{S}^c(\Delta_{\mathcal{S}\otimes c\mathcal{O}^\vee})(\mathrm{id}_\C) \cdot \widehat{\mathscr{S}}_\mathbb{S}^c(\diracComb(\iota_{\C},\widehat{\mathscr{S}}_\mathbb{S}^c(\kappa)(\mathrm{id}_\C) \cdot \Delta_\C))(\mathrm{id}_\C)~,
\]
where the first term is $d_1$ and the second is $d_2$. One computes that:
\[
- d_{\mathrm{cobar}} \cdot (\psi_1 - \psi_2) = \widehat{\mathscr{S}}_\mathbb{S}^c(\kappa \star \kappa)(\mathrm{id}_\C)~.
\]
Since $\kappa$ is a curved twisting morphism, this concludes the proof. 
\end{proof}

Let us give explicit descriptions of these functors.

\begin{Definition}[Complete Cobar construction]
Let $(\mathcal{C},\{\Delta_i\},\epsilon_\C,d_\C)$ be a dg counital partial cooperad. The \textit{complete Cobar construction} $\widehat{\Omega}\C$ of $\mathcal{C}$ is given by:
\[
\widehat{\Omega}\C \coloneqq \left( \overline{\mathscr{T}}^\wedge(s^{-1}\C \oplus \nu), d_{\mathrm{cobar}} = d_1 - d_2, \Theta_{\mathrm{cobar}} \right)~,
\]
where $\overline{\mathscr{T}}^\wedge(s^{-1}\C \oplus \nu)$ is the completed reduced tree monad applied to the dg $\mathbb{S}$-module $s^{-1}\C \oplus \nu$. Here $\nu$ is an arity $1$ degree $-2$ generator. It is endowed with the differential $d_{\mathrm{cobar}}$ given by the difference of $d_1$ and $d_2$. The term $d_1$ is the unique derivation extending
\[
\begin{tikzcd}[column sep=5pc,row sep=3pc]
s^{-1}\C \arrow[r,"s^{-1}d_\C"]
&s^{-1}\C \arrow[r,rightarrowtail]
&\overline{\mathscr{T}}^\wedge(s^{-1}\C \oplus \nu)~.
\end{tikzcd}
\]
The term $d_2$ comes from the \textit{structure} of a dg counital partial cooperad on $\C$, it is given by the unique derivation extending
\[
\begin{tikzcd}[column sep=5pc,row sep=3pc]
s^{-1}\C \arrow[r,"s^{-2}\Delta_{(1)} + s^{-1}\epsilon_\C"]
&(s^{-1}\C \circ_{(1)} s^{-1}\C) \oplus \I.\nu \arrow[r,rightarrowtail]
&\overline{\mathscr{T}}^\wedge(s^{-1}\C \oplus \I.\nu)~.
\end{tikzcd}
\]
Its curvature $\Theta_{\mathrm{cobar}}$ is given by the following map
\[
\begin{tikzcd}[column sep=4pc,row sep=3pc]
\Theta_{\mathrm{cobar}}: \I \arrow[r,"s^{-2}"]
&\I.\nu \arrow[r,hookrightarrow] 
&\overline{\mathscr{T}}^\wedge(s^{-1}\C \oplus \I.\nu)~.
\end{tikzcd}
\]
The resulting complete Cobar construction of $\C$ forms a complete curved absolute partial operad.
\end{Definition}

\begin{Notation}
We denote by $\mathscr{T}^\vee(-)$ the cofree counital partial cooperad endofunctor in the category of pdg $\mathbb{S}$-modules. It is given by the cofree graded $u\mathcal{O}$-coalgebra $\mathscr{C}(u\mathcal{O})(-)$.
\end{Notation}

\begin{Definition}[Complete Bar construction]
Let $(\PP, \gamma_\PP, d_\PP,\Theta_\PP)$ a complete curved absolute partial operad. The \textit{complete Bar construction} $\widehat{\mathrm{B}}\PP$ of $\PP$ is given by:
\[
\widehat{\mathrm{B}\PP} \coloneqq \left(\mathscr{T}^\vee(s\PP), d_{\mathrm{bar}} = d_1 + d_2 \right)~,
\]
where $\mathscr{T}^\vee(s\PP)$ is the cofree counital partial cooperad generated by the pdg $\mathbb{S}$-module $s\PP$. It is endowed the pre-differential $d_{\mathrm{bar}}$, given by the sum of $d_1$ and $d_2$. The term $d_1$ is the unique coderivation extending
\[
\begin{tikzcd}[column sep=4pc,row sep=3pc]
\mathscr{T}^\vee(s\PP)  \arrow[r, twoheadrightarrow]
&s\PP \arrow[r,"sd_{\PP}"]
&s\PP~.
\end{tikzcd} 
\]
The term $d_2$ comes from the \textit{structure} of $\PP$, it is given by the unique coderivation extending 
\[
\begin{tikzcd}[column sep=5pc,row sep=3pc]
\mathscr{T}^\vee(s\PP) \arrow[r, twoheadrightarrow]
&\I \oplus (s\PP \circ_{(1)} s\PP) \arrow[r,"s\Theta_\PP + s^2\gamma_{(1)}"]
&s\PP~.
\end{tikzcd}
\]
The resulting complete Bar construction of $\PP$ forms a dg counital partial cooperad. 
\end{Definition}

\begin{Remark}
The complete Bar construction of a complete curved absolute partial operad $\PP$ is, in general, not coaugmented. Indeed, one has that
\[
d_{\mathrm{bar}}(|) = d_2(|) = \Theta_\PP(\mathrm{id})~,
\]
where $|$ denotes the coaugmented counit of $\mathscr{T}^\vee(s\PP)$.
\end{Remark}

\begin{Proposition}
The following holds:
\begin{enumerate}
\item Let $(\mathcal{C},\{\Delta_i\},\epsilon_\C,d_\C)$ be a dg counital partial cooperad. There is a natural isomorphism
\[
\widehat{\Omega}_\kappa \C \cong \widehat{\Omega}\C~.
\]
\item Let $(\PP, \gamma_\PP, \Theta_\PP, d_\PP)$ a complete curved absolute partial operad. There is a natural isomorphism 
\[
\widehat{\mathrm{B}}_\kappa \PP \cong \widehat{\mathrm{B}} \PP~.
\]
\end{enumerate}
\end{Proposition}

\begin{proof}
This can be shown by direct inspection, extending the isomorphism of Lemma \ref{lemma: dual Schur of cOvee}.
\end{proof}

\begin{lemma}[Curved convolution operad]
Let $(\PP,\gamma_\PP,d_\PP,\Theta_\PP)$ be a complete curved absolute partial operad and let $(\C,\{\Delta_i\},\epsilon_\C,d_\C)$ be a dg counital partial cooperad. The convolution pdg partial operad $(\mathcal{H}om(\C,\PP),\{\circ_i\},\partial)$ forms a curved partial operad endowed with the curvature given by
\[
\begin{tikzcd}
\Theta_{\mathcal{H}om}(\mathrm{id}): \mathcal{C} \arrow[r,"\epsilon_\C"]
&\I \arrow[r,"\Theta_\PP"]
&\mathcal{P}~.
\end{tikzcd}
\]
\end{lemma}

\begin{proof}
The proof is completely analogous to Lemma \ref{lemma curved S colored convolution}.
\end{proof}

\begin{Remark}
We only used the operations $\{\circ_i\}$ on $\PP$ in order to define this convolution operad. These operations are obtained by restricting the structural morphism $\gamma_\PP$ to finite sums inside the completed reduced tree monad. See Appendix B \ref{Appendix B} for more details. In fact, using the morphism $\gamma_\PP$, one can show that this convolution curved partial operad is in fact an curved \textit{absolute} partial operad. 
\end{Remark}

Once we have the curved convolution partial operad, we can define the notion of curved twisting morphism between complete curved absolute partial operads and dg counital partial cooperads.

\begin{Definition}[Curved twisting morphism]
Let $(\PP,\gamma_\PP,d_\PP,\Theta_\PP)$ be a complete curved absolute partial operad and let $(\C,\{\Delta_i\},\epsilon,d_\C)$ be a dg counital partial cooperad. A \textit{curved twisting morphism} $\alpha$ is a Maurer-Cartan element curved pre-Lie algebra given by the curved convolution operad:
\[
\mathfrak{g}_{\C,\PP} \coloneqq \prod_{n\geq 0}\mathrm{Hom}_{\mathbb{S}}(\C(n),\PP(n))~.
\]
This is the data of a morphism of graded $\mathbb{S}$-modules $\alpha: \C \longrightarrow \PP$ of a degree $-1$ such that:
\[
\partial(\alpha) + \alpha \star \alpha = \Theta_{\mathcal{H}om}(\mathrm{id})~.
\]
The set of twisting morphism between $\C$ and $\PP$ will be denoted $\mathrm{Tw}(\C,\PP)$. 
\end{Definition}

The set of curved twisting morphism between between complete curved partial operads and counital partial dg cooperads
\[
\begin{tikzcd}
\mathrm{Tw}(- , -) : (\mathsf{dg}~\mathsf{upCoop})^{\text{op}} \times \mathsf{curv}~\mathsf{abs.pOp}^{\mathsf{comp}}  \arrow[r]
&\mathsf{Set}~,
\end{tikzcd}
\]
which is represented on both sides by the complete Bar-Cobar constructions.

\begin{Proposition}[Complete Bar-Cobar adjunction]
Let $(\PP,\gamma_\PP,\Theta_\PP,d_\PP)$ be a complete curved partial operad and let $(\C,\{\Delta_i\},\epsilon,d_\C)$ be a dg counital partial  cooperad. There are isomorphisms:
\[
\mathrm{Hom}_{\mathsf{curv}~\mathsf{abs.pOp}^{\mathsf{comp}}}(\widehat{\Omega}\C,\PP) \cong \mathrm{Tw}(\C,\PP) \cong \mathrm{Hom}_{\mathsf{dg}~\mathsf{upCoop}}(\C,\widehat{\mathrm{B}}\PP)~,
\]
which are natural in $\C$ and $\PP$. 
\end{Proposition}

\begin{proof}
This can be shown by a straightforward computation.
\end{proof}

\begin{Corollary}\label{Complete Bar-Cobar adjunction}
There is a pair of adjoint functors:
\[
\begin{tikzcd}[column sep=7pc,row sep=3pc]
            \mathsf{dg}~\mathsf{upCoop} \arrow[r, shift left=1.1ex, "\widehat{\Omega}"{name=F}] &\mathsf{curv}~\mathsf{abs.pOp}^{\mathsf{comp}}~. \arrow[l, shift left=.75ex, "\widehat{\mathrm{B}}"{name=U}]
            \arrow[phantom, from=F, to=U, , "\dashv" rotate=-90]
\end{tikzcd}
\]
\end{Corollary}

\section{Counital partial cooperads up to homotopy and transfer of model structures}
In this section, we introduce the notion of a counital partial cooperad up to homotopy. This notion is dual to that of unital partial operads up to homotopy developed in \cite{grignou2021}. In fact, there is a $\mathbb{S}$-colored operad which encodes unital partial operads up to homotopy as its algebras and counital partial cooperads up to homotopy as its coalgebras. The advantage of counital partial cooperads up to homotopy with respect to counital partial cooperads is that the admit a canonical cylinder object. Thus, using a left transfer theorem, we endow them with a model structure where weak equivalences are arity-wise quasi-isomorphisms. Notice that here, the category of counital partial cooperads up to homotopy with \textit{strict morphisms} is the one being considered. Afterwards, we transfer this model structure to the category of complete curved absolute partial operads via a complete Bar-Cobar adjunction.

\begin{Notation}
We denote by $\mathrm{RT}_\omega^n$ the set of rooted trees of arity $n$ with $\omega$ internal edges. There is no restriction on the number of incoming edges of the vertices considered.
\end{Notation}

\begin{Definition}[Counital partial cooperad up to homotopy]
Let $(\C,d_\C)$ be a dg $\mathbb{S}$-module. A \textit{counital partial cooperad up to homotopy structure} on $\C$ amounts to the data of a derivation of degree $-1$
\[
d: \overline{\mathscr{T}}^\wedge(s^{-1}\C \oplus \nu) \longrightarrow \overline{\mathscr{T}}^\wedge(s^{-1}\C \oplus \nu)
\]
such that the restriction of $d$ to $s^{-1}\C$ is given by $s^{-1}d_\C$ and such that for any series of rooted trees of arity $n$ labeled by elements of $s^{-1}\C$:
\[
d^2\left(\sum_{\omega \geq 1} \sum_{\tau \in \mathrm{RT}_\omega^n} \tau \right) = \sum_{\omega \geq 1} \sum_{\tau \in \mathrm{RT}_\omega^n} \left( v \circ_1 \tau - \sum_{i=0}^n \tau \circ_i v \right)~.
\]
In order words, the data of $\left(\overline{\mathscr{T}}^\wedge(s^{-1}\C \oplus \nu),d \right)$ forms a complete curved absolute partial operad.
\end{Definition}

\begin{Remark}
The morphism 
\[
\begin{tikzcd}[column sep=3pc,row sep=3pc]
\C  \arrow[r,"\cong"]
&s^{-1}\C \arrow[r,"d"]
&\overline{\mathscr{T}}^\wedge(s^{-1}\C \oplus \nu) \arrow[r,twoheadrightarrow]
&\I.\nu
\end{tikzcd} 
\]
endows $\C$ with a counit $\epsilon_\C: \C \longrightarrow \I$ which satisfies the counital axiom of counital partial cooperad only up to higher homotopies. The morphism 
\[
\begin{tikzcd}[column sep=3pc,row sep=3pc]
\C  \arrow[r,"\cong"]
&s^{-1}\C \arrow[r,"d"]
&\overline{\mathscr{T}}^\wedge(s^{-1}\C \oplus \nu) \arrow[r,twoheadrightarrow]
&\overline{\mathscr{T}}^{\wedge~(2)}(s^{-1}\C \oplus \nu)
\end{tikzcd} 
\]
endows $\C$ with a family of partial decompositions maps $\{\Delta_i\}$ which satisfy the coparallel and cosequential axioms of a counital partial cooperad only up to higher homotopies. All the higher homotopies are given by the morphism
\[
\begin{tikzcd}[column sep=3pc,row sep=3pc]
\C  \arrow[r,"\cong"]
&s^{-1}\C \arrow[r,"d"]
&\overline{\mathscr{T}}^\wedge(s^{-1}\C \oplus \nu) \arrow[r,twoheadrightarrow]
&\overline{\mathscr{T}}^{\wedge~(\geq 3)}(s^{-1}\C \oplus \nu)~.
\end{tikzcd} 
\]
\end{Remark}

\begin{Example}
Any dg counital partial cooperad $(\C,\{\Delta_i\},\epsilon,d_\C)$ is an example of counital partial cooperad up to homotopy via its complete Cobar construction $\widehat{\Omega}\C$. 
\end{Example}

\begin{Definition}[Strict morphisms]
Let $(\C,d_1,d_\C)$ and $(\D,d_2,d_\D)$ be two counital partial cooperads up to homotopy. A \textit{morphism} $f: \C \longrightarrow \D$ amounts to the data of a morphism of graded $\mathbb{S}$-modules $f: \C \longrightarrow \D$ such that the following diagram commute
\[
\begin{tikzcd}[column sep=3pc,row sep=2pc]
\C  \arrow[r,"\cong"] \arrow[d,"f",swap]
&s^{-1}\C \arrow[r,"d_1 "]
&\overline{\mathscr{T}}^\wedge(s^{-1}\C \oplus \nu) \arrow[d,"\overline{\mathscr{T}}^\wedge(s^{-1}f \oplus \nu)"]\\
\D  \arrow[r,"\cong"]
&s^{-1}\D \arrow[r,"d_2 "]
&\overline{\mathscr{T}}^\wedge(s^{-1}\D \oplus \nu)~. 
\end{tikzcd} 
\]
\end{Definition}

\begin{Remark}[$\infty$-morphisms]
Let $(\C,d_1,d_\C)$ and $(\D,d_2,d_\D)$ be two counital partial cooperads up to homotopy. An $\infty$\textit{-morphism} $f: \C \rightsquigarrow \D$ amounts to the data of a morphism of complete curved absolute partial operads 
\[
f: \left(\overline{\mathscr{T}}^\wedge(s^{-1}\C \oplus \nu),d_1 \right) \longrightarrow \left(\overline{\mathscr{T}}^\wedge(s^{-1}\D \oplus \nu),d_2 \right)~.
\]
These $\infty$-morphisms of counital partial cooperads up to homotopy should have all the expected properties of $\infty$-morphism between unital partial operads described in \cite{grignou2021}. In particular, they are invertible and they describe the hom-sets of the homotopy category. For the sake of brevity, we do not explore this path here.
\end{Remark}

Counital partial cooperads are coalgebras over a dg $\mathbb{S}$-colored operad, which is given by the Cobar construction on the conilpotent curved $\mathbb{S}$-colored partial cooperad $\mathcal{S} \otimes c\mathcal{O}^\vee$. This Cobar construction is the $\mathbb{S}$-colored analogue of the Cobar construction of Definition \ref{Def: Cobar de Brice}. (We swear, dear reader, this is the last Cobar construction of the chapter !)

\begin{Definition}[The dg $\mathbb{S}$-colored operad encoding (co)unital partial (co)operads up to homotopy]
The dg $\mathbb{S}$-colored operad $\Omega_\mathbb{S}(\mathcal{S} \otimes c\mathcal{O}^\vee)$ is given by
\[
\Omega_\mathbb{S}\left(\mathcal{S} \otimes c\mathcal{O}^\vee \right) \coloneqq \left(\mathscr{T}_\mathbb{S}\left(s^{-1}(\mathcal{S} \otimes c\mathcal{O}^\vee\right), d_{\mathrm{cobar}} = -d_2 \right)~,
\]
where $\mathscr{T}_\mathbb{S}(s^{-1}(\mathcal{S} \otimes c\mathcal{O}^\vee))$ is the free unital partial $\mathbb{S}$-colored operad generated by the desuspension of the pdg $\mathbb{S}$-colored scheme $\mathcal{S} \otimes c\mathcal{O}^\vee$. It is endowed with a differential $d_2$, which is given by the unique derivation extending the map:
\[
\begin{tikzcd}[column sep=6.5pc,row sep=3pc]
s^{-1}(\mathcal{S} \otimes c\mathcal{O}^\vee) \arrow[r, "s^{-2}\Delta_{(1)} \oplus s^{-1}\Theta_{\mathcal{S} \otimes c\mathcal{O}^\vee}"]
&\left(s^{-1}(\mathcal{S} \otimes c\mathcal{O}^\vee) \circ_{(1)} s^{-1}(\mathcal{S} \otimes c\mathcal{O}^\vee)\right) \oplus \I_\mathbb{S} \hookrightarrow \mathscr{T}_\mathbb{S}\left(s^{-1}(\mathcal{S} \otimes c\mathcal{O}^\vee)\right)~.
\end{tikzcd} 
\]
\end{Definition}

\begin{Proposition}
The following holds:
\begin{enumerate}
\item The category of unital partial operads up to homotopy with strict morphisms is equivalent to the category of dg algebras over $\Omega_\mathbb{S}\left(\mathcal{S} \otimes c\mathcal{O}^\vee \right)$.

\item The category of counital partial cooperads up to homotopy with strict morphism is equivalent to the category of dg coalgebras over $\Omega_\mathbb{S}\left(\mathcal{S} \otimes c\mathcal{O}^\vee \right)$.
\end{enumerate}
\end{Proposition}

\begin{proof}
For the first statement, the proof is \textit{mutatis mutandis} the same as the proof of \cite[Proposition 36]{grignou2021}. For the second statement, the proof is \textit{mutatis mutandis} the same as the proof of \cite[Theorem 12.1]{grignoulejay18}.
\end{proof}

\begin{Corollary}
There forgetful functor from counital partial cooperads up to homotopy to the category of dg $\mathbb{S}$-modules admits a right adjoint
\[
\begin{tikzcd}[column sep=7pc,row sep=3pc]
             \mathsf{upCoop}_\infty \arrow[r, shift left=1.3ex, "\mathrm{U}"{name=F}] &\mathsf{dg}~\mathbb{S}\text{-}\mathsf{mod}~,  \arrow[l, shift left=1.2ex, "\mathscr{C}\left(\Omega_\mathbb{S}\left(\mathcal{S} \otimes c\mathcal{O}^\vee \right)\right)"{name=U}]
            \arrow[phantom, from=F, to=U, , "\dashv" rotate=-90]
\end{tikzcd}
\]
which is the cofree counital partial cooperad up to homotopy functor.
\end{Corollary}

\begin{proof}
This is an immediate consequence of Theorem \ref{thm: existence de la cogèbre colibre}.
\end{proof}

Using this adjunction, we want to transfer the standard model structure on dg $\mathbb{S}$-modules via a \textit{left-transfer theorem}. The key ingredient to apply the left-transfer theorem in this situation is the construction of a \textit{functorial cylinder object} in the category of counital partial cooperads up to homotopy.

\medskip

Let $\mathrm{I}$ be the commutative Hopf monoid given by the cellular chains on the interval $[0,1]$. Its commutative algebra structure is given by the cup-product and the coassociative coalgebra structure is given by the choice of a cellular diagonal on the topological interval. It is an interval object in the category of dg modules, meaning that there is a factorization
\[
\kk \oplus \kk \rightarrowtail \mathrm{I} \qi \kk
\]
of the codiagonal map $\mathrm{id} \oplus \mathrm{id}: \kk \oplus \kk \longrightarrow \kk$.

\begin{Definition}[Interval dg $\mathbb{S}$-modules]
The \textit{interval} dg $\mathbb{S}$-module $\mathrm{Int}$ is the dg $\mathbb{S}$-module given by 
\[
\mathrm{Int}(n) \coloneqq \mathrm{I}
\]
for all $n \geq 0$. 
\end{Definition}

\begin{Remark}
If we endow $\mathrm{Int}$ with the partial decomposition maps $\Delta_i: \mathrm{Int}(n+k-1) \longrightarrow \mathrm{Int}(n) \otimes \mathrm{Int}(k)$ given by the coalgebra structure of $\mathrm{I}$
\[
\Delta: \mathrm{I} \longrightarrow \mathrm{I} \otimes \mathrm{I}~,
\]
then $\mathrm{Int}$ \textit{almost} forms a dg counital partial cooperad. It satisfies the cosequential axiom but, since $\Delta$ is not \textit{cocommutative}, it does not satisfy the coparallel axiom. If there existed a \textit{cocommutative} interval object in the category of dg module, we could directly transfer the standard model structure of dg $\mathbb{S}$-module to the category of counital partial cooperads, in a dual version of \cite[Theorem 3.2]{bergermoerdieck03}. This is one reason to work with counital partial cooperads up to homotopy.
\end{Remark}

\begin{Proposition}
Let $(\C,d_1,d_\C)$ be a counital partial cooperad up to homotopy. The Hadamard product $\C \otimes \mathrm{Int}$ has a canonical structure of counital partial cooperad up to homotopy.
\end{Proposition}

\begin{proof}
The proof is completely dual to the proof of \cite[Proposition 31]{grignou2021}, where given a unital partial operad up to homotopy, the author builds a unital partial operad up to homotopy structure on $\PP \otimes \mathrm{Int}$ by induction. 
\end{proof}

\begin{Corollary}[Functorial cylinder object]\label{lemma: functorial cylinder}
Let $(\C,d_1,d_\C)$ be a counital partial cooperad up to homotopy. Then the codiagonal morphism $\mathrm{id} \oplus \mathrm{id}: \mathcal{C} \oplus \mathcal{C} \longrightarrow \mathcal{C}$ factors as follows
\[
\begin{tikzcd}[column sep=4pc,row sep=0pc]
\mathcal{C} \oplus \mathcal{C} \arrow[r,rightarrowtail]
&\mathcal{C} \otimes \mathrm{Int} \arrow[r,"\sim"]
&\mathcal{C} ~,
\end{tikzcd}
\]
where the first arrow is a degree-wise monomorphism and the second arrow is an arity-wise quasi-isomorphism. Thus $\C \otimes \mathrm{Int}$ is a cylinder object which is functorial in $\mathcal{C}$. 
\end{Corollary}

\begin{theorem}[Transferred model structure]\label{thm: model structure on counital cooperads}
There is a cofibrantly generated model structure on the category of counital partial cooperads up to homotopy defined by the following classes of morphisms:
\begin{enumerate}
\item The class of weak-equivalences $\mathrm{W}$ is given by strict morphisms of counital partial cooperads up to homotopy $f: \mathcal{C} \longrightarrow \mathcal{D}$ such that $f(n): \mathcal{C}(n) \longrightarrow \mathcal{D}(n)$ is a quasi-isomorphism for all $n \geq 0$. 
\vspace{0.5pc}
\item The class of cofibrations is $\mathrm{Cof}$ is given by strict morphisms of counital partial cooperads up to homotopy $f: \mathcal{C} \longrightarrow \mathcal{D}$ such that $f(n): \mathcal{C}(n) \rightarrowtail \mathcal{D}(n)$ is a degree-wise monomorphism for all $n \geq 0$. 
\vspace{0.5pc}
\item The class of fibrations $\mathrm{Fib}$ is given by strict morphisms of counital partial cooperads up to homotopy which have the right lifting property against all morphism in $\mathrm{W} \cap \mathrm{Cof}$. 
\end{enumerate}
\end{theorem}

\begin{proof}
In order to apply the left-transfer theorem such as \cite{lefttransfer}, essentially need a functorial cofibrant replacement functor and a functorial cylinder object. Since every object is cofibrant in this model structure, we only need a functorial cylinder object. This follows from Corollary \ref{lemma: functorial cylinder}.
\end{proof}

Now our aim is to transfer this structure to category of complete curved absolute partial operads, replicating the methods of Section 10 in \cite{grignoulejay18}. For this purpose, we induce a complete Bar-Cobar adjunction using a curved twisting morphism between the $\mathbb{S}$-colored objects encoding these categories. (Yes, dear reader, we lied).

\begin{lemma}[Curved twisting morphism $\iota$]
The $\mathbb{S}$-color scheme map given by the inclusion 
\[
\iota: \mathcal{S} \otimes c\mathcal{O}^\vee \hookrightarrow \Omega_\mathbb{S}\left(\mathcal{S} \otimes c\mathcal{O}^\vee \right)
\]
is a curved twisting morphism.
\end{lemma}

\begin{proof}
This is immediate to check from the definition. 
\end{proof}

\begin{Proposition}[Complete Bar-Cobar adjunction relative to $\iota$]\label{Prop: complete Bar-Cobar adjunction relative to iota}
There is a pair of adjoint functors 
\[
\begin{tikzcd}[column sep=7pc,row sep=3pc]
            \mathsf{upCoop}_\infty \arrow[r, shift left=1.1ex, "\widehat{\Omega}_\iota"{name=F}] &\mathsf{curv}~\mathsf{abs.pOp}^{\mathsf{comp}}~. \arrow[l, shift left=.75ex, "\widehat{\mathrm{B}}_\iota"{name=U}]
            \arrow[phantom, from=F, to=U, , "\dashv" rotate=-90]
\end{tikzcd}
\]
\end{Proposition}

\begin{proof}
The constructions of the complete Bar and the complete Cobar functors are \textit{mutatis mutandis} the same as those described in the previous section, where we constructed the complete Bar-Cobar relative to $\kappa$. The proof that these functors form an adjunction is also \textit{mutatis mutandis} the same.
\end{proof}

This allows us in turn to transfer the above model structure from counital partial cooperads up to homotopy to complete curve absolute partial operads using a \textit{right transfer theorem}.

\begin{theorem}\label{thm: model structure on curved operads}
There is a cofibrantly generated model structure on the category of complete curved absolute partial operads defined by the following classes of morphisms:
\begin{enumerate}
\item The class of weak-equivalences $\mathrm{W}$ is given by morphisms of complete curved absolute partial operads $f: \mathcal{P} \longrightarrow \mathcal{Q}$ such that $\widehat{\mathrm{B}}_\iota(f)(n): \widehat{\mathrm{B}}_\iota(\mathcal{P})(n) \longrightarrow \widehat{\mathrm{B}}_\iota(\mathcal{Q})(n)$ is a quasi-isomorphism for all $n \geq 0$. 
\vspace{0.5pc}
\item The class of fibrations is $\mathrm{Fib}$ is given by morphisms of complete curved absolute partial operads $f: \mathcal{P} \longrightarrow \mathcal{Q}$ such that $\widehat{\mathrm{B}}_\iota (f): \widehat{\mathrm{B}}_\iota(\mathcal{P}) \longrightarrow \widehat{\mathrm{B}}_\iota(\mathcal{Q})$ is a fibration of counital partial cooperads up to homotopy. 
\vspace{0.5pc}
\item The class of cofibrations $\mathrm{Cof}$ is given by morphisms of complete curved absolute partial operads which have the left lifting property against all morphism in $\mathrm{W} \cap \mathrm{Fib}$. 
\end{enumerate}
\end{theorem}

\begin{proof}
In order to apply the right transfer theorem, we need to check that the acyclicity conditions of \cite[Section 2.5]{bergermoerdieck03} are satisfied. To the best of our knowledge, one can not do this by constructing a path-object in the category of complete curved absolute partial operads. Thus these conditions need to be checked by direct computation, as it is done in \cite[Section 10.5]{grignoulejay18}. In is a straightforward but tedious exercise to check that the arguments used in \textit{loc.cit} generalize to our framework.
\end{proof}

\begin{Proposition}
Let $f: (\mathcal{P},\gamma_\PP,d_\PP,\Theta_\PP) \longrightarrow (\mathcal{Q},\gamma_\Q,d_\Q,\Theta_\Q)$ be a morphism of complete curved absolute partial operads. If $f(n): \mathcal{P}(n) \longrightarrow \mathcal{Q}(n)$ is a degree-wise epimorphism, then it is a fibration. Thus all complete curved absolute partial operads are fibrant. 
\end{Proposition}

\begin{proof}
The arguments of Section 10.2 in \cite{grignoulejay18} generalize to the operadic setting \textit{mutatis mutandis}. 
\end{proof}

\begin{Definition}[Graded quasi-isomorphism]
Let $f: (\mathcal{P},\gamma_\PP,d_\PP,\Theta_\PP) \longrightarrow (\mathcal{Q},\gamma_\Q,d_\Q,\Theta_\Q)$ be a morphism of complete curved absolute partial operads. It is a \textit{graded quasi-isomorphism} if the induced morphism of dg $\mathbb{S}$-modules
\[
\mathrm{gr}(f): \mathrm{gr}(\mathcal{P}) \cong \bigoplus_{\omega \geq 1} \overline{\mathscr{F}}_{\omega} \PP/ \overline{\mathscr{F}}_{\omega +1} \PP \longrightarrow  \mathrm{gr}(\mathcal{Q}) \cong \bigoplus_{\omega \geq 1} \overline{\mathscr{F}}_{\omega} \mathcal{Q}/ \overline{\mathscr{F}}_{\omega +1} \Q
\]
is an arity-wise quasi-isomorphism. Here $\overline{\mathscr{F}}_{\omega} \PP$ denotes the $\omega$-term of the canonical filtration of an absolute partial operad defined in Appendix \ref{Appendix B}.
\end{Definition}

\begin{Remark}
One checks easily that since $d_\mathcal{P}^2$ raises the weigh of an operation in $\mathcal{P}$ by one, it is equal to zero in the associated graded. Therefore $\mathrm{gr}(\mathcal{P})$ forms a dg $\mathbb{S}$-module.
\end{Remark}

\begin{Proposition}
Let $f: (\mathcal{P},\gamma_\PP,d_\PP,\Theta_\PP) \longrightarrow (\mathcal{Q},\gamma_\Q,d_\Q,\Theta_\Q)$  be a graded quasi-isomorphism between two complete curved absolute partial operads. Then it is a weak-equivalence, meaning that $\widehat{\mathrm{B}}_\iota(f): \widehat{\mathrm{B}}(\mathcal{P}) \longrightarrow \widehat{\mathrm{B}}_\iota(\mathcal{Q})$ is an arity-wise quasi-isomorphism.
\end{Proposition}

\begin{proof}
The proof is completely analogous to that of \cite[Theorem 10.23]{grignoulejay18}. It is essentially done by induction. The cofree counital cooperad functor preserves quasi-isomorphisms. Thus it sends graded quasi-isomorphisms between curved partial operads endowed with the trivial operad structure to quasi-isomorphisms. Then one shows a weak version of the five-lemma \cite[Lemma 10.20]{grignoulejay18} which allows us to go from weight $\omega$ to weight $\omega +1$. See Section 10.3 for more details on that matter. 
\end{proof}

\begin{Remark}
The above proposition implies that the weak-equivalences of complete curved operads as defined in \cite{JoanCurved} by graded quasi-isomorphisms should be strictly included inside the weak-equivalences of our model structure. Nevertheless, the objects and the underlying categories considered are quite different in nature, so we do not attempt a precise comparison result.
\end{Remark}

There is a morphism of dg unital partial $\mathbb{S}$-colored operads $f_\kappa: \Omega_\mathbb{S}\left(\mathcal{S} \otimes c\mathcal{O}^\vee \right) \longrightarrow u\mathcal{O}$, which in turn induces a morphism between the associated comonads in the category of dg $\mathbb{S}$-modules. Thus one has an adjunction
\[
\begin{tikzcd}[column sep=7pc,row sep=3pc]
            \mathsf{dg}~\mathsf{upCoop} \arrow[r, shift left=1.1ex, "\mathrm{Res}_{f_\kappa}"{name=F}] &\mathsf{upCoop}_\infty~, \arrow[l, shift left=.75ex, "\mathrm{Coind}_{f_\kappa}"{name=U}]
            \arrow[phantom, from=F, to=U, , "\dashv" rotate=-90]
\end{tikzcd}
\]
where the forgetful functor $\mathrm{Res}_{f_\kappa}$ is fully faithful and preserves quasi-isomorphisms.

\begin{Proposition}\label{prop: commuting triangle of adjunctions}
There following triangle of adjunctions 
\[
\begin{tikzcd}[column sep=5pc,row sep=2.5pc]
&\hspace{1pc}\mathsf{upCoop}_\infty \arrow[dd, shift left=1.1ex, "\widehat{\Omega}_{\iota}"{name=F}] \arrow[ld, shift left=.75ex, "\mathrm{Coind}_{f_\kappa}"{name=C}]\\
\mathsf{dg}~\mathsf{upCoop}  \arrow[ru, shift left=1.5ex, "\mathrm{Res}_{f_\kappa}"{name=A}]  \arrow[rd, shift left=1ex, "\widehat{\Omega}_{\kappa}"{name=B}] \arrow[phantom, from=A, to=C, , "\dashv" rotate=-70]
& \\
&\hspace{2.5pc}\mathsf{curv}~\mathsf{abs.pOp}^{\mathsf{comp}}~, \arrow[uu, shift left=.75ex, "\widehat{\text{B}}_{\iota}"{name=U}] \arrow[lu, shift left=.75ex, "\widehat{\text{B}}_{\kappa}"{name=D}] \arrow[phantom, from=B, to=D, , "\dashv" rotate=-110] \arrow[phantom, from=F, to=U, , "\dashv" rotate=-180]
\end{tikzcd}
\]
commutes up to natural isomorphism.
\end{Proposition}

\begin{proof}
It is immediate to check that the left-adjoints are naturally isomorphic. Thus by mates of adjunction, the right adjoints are also naturally isomorphic.
\end{proof}

\begin{Remark}
We believe that the complete Bar-Cobar adjunction relative to $\iota$ is a Quillen equivalence. One could generalize the arguments in Section 11 of \cite{grignoulejay18} to the $\mathbb{S}$-colored framework in order to prove this. 

\medskip

On the other hand, one could try to induce a model structure on counital partial cooperads using the forgetful-cofree adjunction and transfer it to the category of complete curved absolute partial operads. It is not clear to us that these model structures should coincide. Indeed, this would amount to proving that the adjunction
\[
\begin{tikzcd}[column sep=7pc,row sep=3pc]
            \mathsf{dg}~\mathsf{upCoop} \arrow[r, shift left=1.1ex, "\mathrm{Res}_{f_\kappa}"{name=F}] &\mathsf{upCoop}_\infty~, \arrow[l, shift left=.75ex, "\mathrm{Coind}_{f_\kappa}"{name=U}]
            \arrow[phantom, from=F, to=U, , "\dashv" rotate=-90]
\end{tikzcd}
\]
is a Quillen equivalence. But this adjunction restricts to the adjunction
\[
\begin{tikzcd}[column sep=7pc,row sep=3pc]
            \mathsf{dg}~u\mathcal{A}ss\text{-}\mathsf{coalg} \arrow[r, shift left=1.1ex, "\mathrm{Res}_{f_\kappa}"{name=F}] &u\mathcal{A}_\infty\text{-}\mathsf{coalg}~, \arrow[l, shift left=.75ex, "\mathrm{Coind}_{f_\kappa}"{name=U}]
            \arrow[phantom, from=F, to=U, , "\dashv" rotate=-90]
\end{tikzcd}
\]
when one restricts to dg $\mathbb{S}$-modules concentrated in arity one. It is not clear at all that this adjunction is a Quillen equivalence, see \cite[Conjecture 8.10]{grignoulejay18}.
\end{Remark}

\section{Curved operadic duality squares}
In this section, we build duality functors that intertwine the two Bar-Cobar adjunctions defined in Section \ref{Section: Constructions Bar-Cobar operadiques}, forming an algebraic duality square of commuting adjunctions. On the homotopical side of things, we build a second duality square of commuting Quillen adjunctions when (co)unital partial (co)operads are replaced by their "up to homotopy" counterparts. Using these squares, we compute minimal cofibrant resolutions for complete curved absolute partial operads in certain cases of interest. 

\subsection{Algebraic duality square square}
Our goal is to construct left adjoint functors to the linear dual functor that sends "coalgebraic objects" into "algebraic objects".

\begin{lemma}\label{lemma: adjoint à droite}
The linear duality functor 
\[
\begin{tikzcd}[column sep=4pc,row sep=0pc]
\mathsf{dg}~\mathsf{upCoop}^{\mathsf{op}} \arrow[r,"(-)^*"] 
&\mathsf{dg}~\mathsf{upOp}
\end{tikzcd}
\]
admits a left adjoint.
\end{lemma}

\begin{proof}
Consider the following square of functors
\[
\begin{tikzcd}[column sep=4pc,row sep=4pc]
\mathsf{dg}~\mathsf{upCoop}^{\mathsf{op}} \arrow[r,"(-)^*"]  \arrow[d,"\mathrm{U}^\mathsf{op}"{name=SD},shift left=1.1ex ]
&\mathsf{dg}~\mathsf{upOp} \arrow[d,"\mathrm{U}"{name=LDC},shift left=1.1ex ] \\
\mathsf{dg}~\mathbb{S}\text{-}\mathsf{mod}^{\mathsf{op}} \arrow[r,"(-)^*"{name=CC},,shift left=1.1ex] \arrow[u,"\left(\mathscr{T}^\vee(-)\right)^\mathsf{op}"{name=LD},shift left=1.1ex ] \arrow[phantom, from=SD, to=LD, , "\dashv" rotate=0]
&\mathsf{dg}~\mathbb{S}\text{-}\mathsf{mod}~.  \arrow[l,"(-)^*"{name=CB},shift left=1.1ex] \arrow[u,"\mathscr{T}(-)"{name=TD},shift left=1.1ex] \arrow[phantom, from=TD, to=LDC, , "\dashv" rotate=0] \arrow[phantom, from=CC, to=CB, , "\dashv" rotate=90]
\end{tikzcd}
\] 
First notice that the adjunction on the left hand side is monadic, since we consider the opposite of a comonadic adjunction. Furthermore, all categories involved are complete and cocomplete. It is absolutely clear that $(-)^* \cdot \mathrm{U}^\mathsf{op} \cong \mathrm{U} \cdot (-)^*~.$ Thus we can apply the Adjoint Lifting Theorem \cite[Theorem 2]{AdjointLifting} to this situation, which concludes the proof. 
\end{proof}

\begin{Definition}[Sweedler dual]\label{def: Sweedler dual functor}
The \textit{Sweedler duality functor}
\[
\begin{tikzcd}[column sep=4pc,row sep=0pc]
\mathsf{dg}~\mathsf{upOp} \arrow[r,"(-)^\circ"] 
&\mathsf{dg}~\mathsf{upCoop}^{\mathsf{op}}
\end{tikzcd}
\]
is defined as the left adjoint of the linear dual functor. 
\end{Definition}

\begin{Remark}
The proof of the Adjoint Lifting Theorem \cite[Theorem 2]{AdjointLifting} gives an explicit construction of this left adjoint. Let $(\mathcal{P},\{\circ_i\},\eta,d_\PP)$ be a dg unital partial operad, and let 
\[
\gamma_\PP: \mathscr{T}(\mathcal{P}) \longrightarrow \PP
\]
be its structural morphism as an algebra over the tree monad. The Sweelder dual dg counital partial cooperad $\mathcal{P}^\circ$ is given by the following equalizer: 
\[
\begin{tikzcd}[column sep=4pc,row sep=4pc]
\mathrm{Eq}\Bigg(\mathscr{T}^\vee(\mathcal{P}^*) \arrow[r,"(\gamma_\PP)^*",shift right=1.1ex,swap]  \arrow[r,"\varrho"{name=SD},shift left=1.1ex ]
&\mathscr{T}^\vee\left((\mathscr{T}(\mathcal{P}))^*\right) \Bigg)~,
\end{tikzcd}
\]
where $\varrho$ is an arrow constructed using the comonadic structure of $\mathscr{T}^\vee$ and the canonical inclusio of a dg $\mathbb{S}$-module into its double linear dual.
\end{Remark}

\begin{Remark}
If we restrict to unital partial operads concentrated in arity $1$, that is, unital associative algebras, the Sweedler dual functor we have constructed is naturally isomorphic with the original Sweedler dual of \cite{Sweedler69}.
\end{Remark}

\begin{Remark}
Let $(\mathcal{P},\{\circ_i\},\eta,d_\PP)$ be a dg unital partial operad such that $\mathcal{P}(n)$ is degree-wise finite dimensional over $\kk$. Its Sweedler dual $\mathcal{P}^\circ$ is simply given by $(\mathcal{P}^*,\{ \circ_i^*\},\eta^*,d_\PP^*)$. The adjunction constructed restricts to an anti-equivalence of categories between dg counital partial cooperads which are arity and degree-wise finite dimensional and dg unital partial operads which are arity and degree-wise finite dimensional.
\end{Remark}

Let's turn to the other side of the Koszul duality. We postpone the following proofs and constructions to the Appendix \ref{Appendix B}, where there is a detailed discussion of absolute partial operads and their properties. 

\begin{lemma}
Let $(\mathcal{C},\{\Delta_i\},d_\C,\Theta_\C)$ be a conilpotent curved partial cooperad. Then its linear dual $\C^*$ has inherits a structure of a complete curved absolute partial operad. This defines a functor 
\[
\begin{tikzcd}[column sep=4pc,row sep=0pc]
\left(\mathsf{curv}~\mathsf{pCoop}^{\mathsf{conil}}\right)^{\mathsf{op}} \arrow[r,"(-)^*"]
&\mathsf{curv}~\mathsf{abs.pOp}^{\mathsf{comp}}~.
\end{tikzcd}
\]
\end{lemma}

\begin{proof}
See Lemma \ref{lemma: dual lin d'une curved conil coop}.
\end{proof}

\begin{Proposition}
The linear duality functor admits a left adjoint. There is an adjunction 
\[
\begin{tikzcd}[column sep=7pc,row sep=3pc]
\mathsf{curv}~\mathsf{pOp}^{\mathsf{comp}}  \arrow[r, shift left=1.1ex, "(-)^\vee"{name=F}] 
& \left(\mathsf{curv}~\mathsf{pCoop}^{\mathsf{conil}}\right)^{\mathsf{op}} ~. \arrow[l, shift left=.75ex, "(-)^*"{name=U}]
            \arrow[phantom, from=F, to=U, , "\dashv" rotate=-90]
\end{tikzcd}
\]
\end{Proposition}

\begin{proof}
See Proposition \ref{prop: adjonction topo et dual lin en curved}.
\end{proof}

\begin{Definition}[Topological dual functor]\label{def: topological dual functor}
The \textit{topological dual functor} 
\[
\begin{tikzcd}[column sep=4pc,row sep=0pc]
\mathsf{curv}~\mathsf{abs.pOp}^{\mathsf{comp}} \arrow[r,"(-)^\vee"]
&\left(\mathsf{curv}~\mathsf{pCoop}^{\mathsf{conil}}\right)^{\mathsf{op}}
\end{tikzcd}
\]
is defined as the left adjoint of the linear dual functor.
\end{Definition}

This allows us to construct the first duality square of commuting functors.

\begin{theorem}[Duality square]\label{thm: carré magique}
The following square of adjunction 
\[
\begin{tikzcd}[column sep=5pc,row sep=5pc]
\left(\mathsf{dg}~\mathsf{upOp}\right)^{\mathsf{op}} \arrow[r,"\mathrm{B}^{\mathsf{op}}"{name=B},shift left=1.1ex] \arrow[d,"(-)^\circ "{name=SD},shift left=1.1ex ]
&\left(\mathsf{curv}~\mathsf{pCoop}^{\mathsf{conil}}\right)^{\mathsf{op}}  \arrow[d,"(-)^*"{name=LDC},shift left=1.1ex ] \arrow[l,"\Omega^{\mathsf{op}}"{name=C},,shift left=1.1ex]  \\
\mathsf{dg}~\mathsf{upCoop} \arrow[r,"\widehat{\Omega}"{name=CC},shift left=1.1ex]  \arrow[u,"(-)^*"{name=LD},shift left=1.1ex ]
&\mathsf{curv}~\mathsf{pOp}^{\mathsf{comp}}~, \arrow[l,"\widehat{\mathrm{B}}"{name=CB},shift left=1.1ex] \arrow[u,"(-)^\vee"{name=TD},shift left=1.1ex] \arrow[phantom, from=SD, to=LD, , "\dashv" rotate=0] \arrow[phantom, from=C, to=B, , "\dashv" rotate=-90]\arrow[phantom, from=TD, to=LDC, , "\dashv" rotate=0] \arrow[phantom, from=CC, to=CB, , "\dashv" rotate=-90]
\end{tikzcd}
\] 
commutes in the following sense: right adjoints going from the top right to the bottom left are naturally isomorphic.
\end{theorem}

\begin{proof}
Let us show the commutativity of this square. We have, for any graded $\mathbb{S}$-module $M$, a natural isomorphism of graded counital partial cooperads 
\[
\mathscr{T}^\vee(sM^*) \cong \left(\mathscr{T}(s^{-1}M) \right)^\circ~.
\]
This isomorphism is obtained as the mate of the obvious isomorphism $(-)^* \cdot \mathrm{U} \cong \mathrm{U}^\mathrm{op} \cdot (-)^*~.$ Let $(\C, \{\Delta_i\}, d_\C, \Theta_\C)$ be a conilpotent curved partial cooperad. One can show by direct inspection that the isomorphism of graded counital partial cooperads 
\[
\mathscr{T}^\vee(s\C^*) \cong \left(\mathscr{T}(s^{-1}\C) \right)^\circ
\]
extends to an isomorphism of dg counital partial cooperads
\[
\widehat{\mathrm{B}}(\C^*) \cong \left(\Omega(\C) \right)^\circ~.
\]
\end{proof}

\begin{Remark}
Let $(\C,\{\Delta_i\},\epsilon,d_\C)$ be a dg counital partial cooperad. Then, by the above theorem, there is an isomorphism
\[
\mathrm{B}(\C^*) \cong \left(\widehat{\Omega}(\C)\right)^\vee~,
\]
which is natural in $\C$.
\end{Remark}

\begin{Proposition}\label{prop: finite dual commutes}
Let $(\mathcal{P},\{\circ_i\},\eta,d_\PP)$ be a dg unital partial operad which is arity-wise and degree-wise finite dimensional. There is an isomorphism
\[
\widehat{\Omega}(\mathcal{P}^*) \cong \left(\mathrm{B}(\mathcal{P})\right)^*
\]
of complete curved absolute partial operads.
\end{Proposition}

\begin{proof}
Let $M$ be a pdg $\mathbb{S}$-module which is arity-wise and degree-wise finite dimensional. There is an isomorphism of complete pdg absolute partial operads:
\[
\overline{\mathscr{T}}^\wedge(s^{-1}M^* \oplus \nu) \cong \left(\overline{\mathscr{T}}(sM \oplus \nu)\right)^*~. 
\]
One can show that this isomorphism extends to an isomorphism of complete curved absolute partial operads
\[
\widehat{\Omega}(\mathcal{P}^*) \cong \left(\mathrm{B}(\mathcal{P})\right)^*
\]
by direct inspection, looking at the morphisms that induce the pre-differentials on each of those constructions.
\end{proof}

\begin{Remark}[Beck-Chevalley condition]
In fact, one can show that there is a monomorphism of complete curved absolute partial operads
\[
\lambda_\PP: \widehat{\Omega}(\mathcal{P}^\circ) \hookrightarrow \left(\mathrm{B}(\mathcal{P})\right)^*
\]
which is natural in $\mathcal{P}$. And $\lambda_\PP$ is an isomorphism if and only if $\PP$ is arity-wise and degree-wise finite dimensional. So the sub-category of dg unital partial operad which satisfy the \textit{Beck-Chevalley condition} is exactly the sub-category of arity-wise and degree-wise finite dimensional dg unital partial operads. 
\end{Remark}

\begin{Example}
Let $\ucom$ be the unital partial operad encoding dg unital commutative algebras. There is an isomorphism of complete curved absolute partial operads: 
\[
(\text{B}\ucom)^* \cong \widehat{\Omega}\ucomd~.
\]
The complete curved absolute partial operad $\widehat{\Omega}\ucomd$ encodes the notion of mixed curved $\mathcal{L}_\infty$-algebras. See Section \ref{Section: curved HTT} or Appendix \ref{Section: Appendix B} for more details on this. 
\end{Example}

\subsection{Homotopical duality square} We now construct the homotopical version of the duality square by replacing (co)unital partial (co)operads by their "up to homotopy" counterparts. Before making those constructions, we state results that are a direct consequence of the results stated in \cite{grignou2021}. They give Quillen equivalence between unital partial operads up to homotopy and conilpotent curved partial cooperads.

\begin{theorem}
There is a cofibrantly generated model structure on the category of unital partial operads up to homotopy with strict morphisms defined by the following classes of morphisms:
\begin{enumerate}
\item The class of weak-equivalences $\mathrm{W}$ is given by strict morphisms of unital partial operads up to homotopy $f: \mathcal{P} \longrightarrow \mathcal{Q}$ such that $f(n): \mathcal{P}(n) \longrightarrow \mathcal{Q}(n)$ is a quasi-isomorphism for all $n \geq 0$. 
\vspace{0.5pc}
\item The class of fibrations is $\mathrm{Fib}$ is given by strict morphisms of unital partial operads up to homotopy $f: \mathcal{P} \longrightarrow \mathcal{Q}$ such that $f(n): \mathcal{P}(n) \twoheadrightarrow \mathcal{Q}(n)$ is a degree-wise epimorphism for all $n \geq 0$. 
\item The class of cofibrations $\mathrm{Cof}$ is given by strict morphisms of unital partial operads up to homotopy which have the left lifting property against all morphism in $\mathrm{W} \cap \mathrm{Fib}$. 
\end{enumerate}
\end{theorem}

\begin{proof}
This is a direct consequence of \cite[Section 4.5]{grignou2021}, which allows for a right transfer theorem using the free-forgetful adjunction.
\end{proof}

The author of \cite{grignou2021} endows the category of conilpotent curved partial cooperads with a transferred structure from unital partial operads, using the Bar-Cobar construction of Subsection \ref{subsection: classical bar-cobar relative to kappa}.

\begin{theorem}[{\cite[Theorem 5]{grignou2021}}]\label{thm: structure de modèles conil curved part coop}
There is a cofibrantly generated model structure on the category of conilpotent curved partial cooperads defined by the following classes of morphisms:
\begin{enumerate}
\item The class of weak-equivalences $\mathrm{W}$ is given by morphisms of conilpotent curved partial cooperads $f: \mathcal{C} \longrightarrow \mathcal{D}$ such that $\Omega(f)(n): \Omega(\mathcal{C})(n) \longrightarrow \Omega(\mathcal{D})(n)$ is a quasi-isomorphism for all $n \geq 0$. \vspace{0.00001pc}
\item The class of cofibrations is $\mathrm{Cof}$ is given by morphisms of conilpotent curved partial cooperads $f: \mathcal{C} \longrightarrow \mathcal{D}$ such that $\Omega(f): \Omega(\mathcal{C}) \longrightarrow \Omega(\mathcal{D})$ is a cofibration of unital partial operads up to homotopy.
\vspace{0.5pc}
\item The class of fibrations $\mathrm{Fib}$ is given by morphisms of conilpotent curved partial cooperads which have the right lifting property against all morphism in $\mathrm{W} \cap \mathrm{Cof}$. 
\end{enumerate}
\end{theorem}

\begin{Proposition}[{\cite[Proposition 13]{grignou2021}}]
In the model structure of conilpotent curved partial cooperads, the class of cofibrations is given by morphisms of conilpotent curved partial cooperads $f: \mathcal{C} \longrightarrow \mathcal{D}$ such that $f(n):\mathcal{C}(n) \longrightarrow \mathcal{D}(n)$ is a monomorphism for all $n \geq 0$. 
\end{Proposition}

The following result is a straightforward consequence of the above.

\begin{Proposition}
There is a Bar-Cobar adjunction
\[
\begin{tikzcd}[column sep=7pc,row sep=3pc]
           \mathsf{upOp}_\infty \arrow[r, shift left=1.1ex, "\Omega_\iota"{name=F}] &\mathsf{curv}~\mathsf{pCoop}^{\mathsf{conil}}~. \arrow[l, shift left=.75ex, "\mathrm{B}_\iota"{name=U}]
            \arrow[phantom, from=F, to=U, , "\dashv" rotate=-90]
\end{tikzcd}
\]
which is a Quillen equivalence.
\end{Proposition}

\begin{proof}
This adjunction is induced by the curved twisting morphism $\iota:\mathcal{S} \otimes c\mathcal{O}^\vee \longrightarrow \Omega_\mathbb{S}\left(\mathcal{S} \otimes c\mathcal{O}^\vee \right)$ using the same methods that in Subsection \ref{subsection: classical bar-cobar relative to kappa}. Consider the canonical morphism of dg unital partial $\mathbb{S}$-colored operads
\[
f_\kappa: \Omega_\mathbb{S}\left(\mathcal{S} \otimes c\mathcal{O}^\vee \right) \longrightarrow u\mathcal{O}~.
\]
Since $\kappa: \mathcal{S} \otimes c\mathcal{O}^\vee \longrightarrow u\mathcal{O}$ is a Koszul curved twisting morphism by Theorem \ref{Koszulity of uO}, then $f_\kappa$ is an arity-wise quasi-isomorphism. Thus the adjunction induced by $f_\kappa$
\[
\begin{tikzcd}[column sep=7pc,row sep=3pc]
            \mathsf{upOp}_\infty \arrow[r, shift left=1.1ex, "\mathrm{Ind}_{f_\kappa}"{name=F}] &\mathsf{dg}~\mathsf{upOp}~. \arrow[l, shift left=.75ex, "\mathrm{Res}_{f_\kappa}"{name=U}]
            \arrow[phantom, from=F, to=U, , "\dashv" rotate=-90]
\end{tikzcd}
\]
is a Quillen equivalence between the category of dg unital partial operads and unital partial operads up to homotopy. Using the fact that the following triangle of adjunctions commutes
\[
\begin{tikzcd}[column sep=5pc,row sep=2.5pc]
&\hspace{1pc}\mathsf{upOp}_\infty \arrow[dd, shift left=1.1ex, "\mathrm{B}_{\iota}"{name=F}] \arrow[ld, shift left=.75ex, "\mathrm{Ind}_{f_\kappa}"{name=C}]\\
\mathsf{dg}~\mathsf{upOp}  \arrow[ru, shift left=1.5ex, "\mathrm{Res}_{f_\kappa}"{name=A}]  \arrow[rd, shift left=1ex, "\mathrm{B}"{name=B}] \arrow[phantom, from=A, to=C, , "\dashv" rotate=110]
& \\
&\hspace{1.5pc}\mathsf{curv}~\mathsf{pCoop}^{\mathsf{conil}}~, \arrow[uu, shift left=.75ex, "\Omega_{\iota}"{name=U}] \arrow[lu, shift left=.75ex, "\Omega"{name=D}] \arrow[phantom, from=B, to=D, , "\dashv" rotate=70] \arrow[phantom, from=F, to=U, , "\dashv" rotate=0]
\end{tikzcd}
\]
we conclude that the Bar-Cobar adjunction induced by $\iota:\mathcal{S} \otimes c\mathcal{O}^\vee \longrightarrow \Omega_\mathbb{S}\left(\mathcal{S} \otimes c\mathcal{O}^\vee \right)$ is also a Quillen equivalence.
\end{proof}

\begin{lemma}\label{lemma: adjoint à droite Sweedler à homotopies près}
The linear duality functor 
\[
\begin{tikzcd}[column sep=4pc,row sep=0pc]
\left(\mathsf{upCoop}_\infty\right)^{\mathsf{op}} \arrow[r,"(-)^*"] 
&\mathsf{upOp}_\infty
\end{tikzcd}
\]
admits a left adjoint.
\end{lemma}

\begin{proof}
Recall that there exists a dg unital partial $\mathbb{S}$-colored operad $\Omega_\mathbb{S}\left(\mathcal{S} \otimes c\mathcal{O}^\vee \right)$ which encodes unital partial operads up to homotopy as its algebras and counital partial cooperads up to homotopy as its coalgebras. Thus the first category is monadic and the second is comonadic. The rest of the proof is \textit{mutatis mutandis} the same as the proof of Lemma \ref{lemma: adjoint à droite}.
\end{proof}

\begin{Definition}[Homotopical Sweedler dual]\label{def: Sweedler dual functor homotopy}
The \textit{homotopical Sweedler duality functor}
\[
\begin{tikzcd}[column sep=4pc,row sep=0pc]
\mathsf{upOp}_\infty \arrow[r,"(-)_h^\circ"] 
&\left(\mathsf{upCoop}_\infty\right)^{\mathsf{op}}
\end{tikzcd}
\]
is defined as the left adjoint of the linear dual functor. 
\end{Definition}

\begin{lemma}
The adjunction 
\[
\begin{tikzcd}[column sep=7pc,row sep=3pc]
            \mathsf{upCoop}_\infty \arrow[r, shift left=1.1ex, "(-)^*"{name=F}] &\left(\mathsf{upOp}_\infty\right)^{\mathsf{op}} ~. \arrow[l, shift left=.75ex, "(-)_h^\circ"{name=U}]
            \arrow[phantom, from=F, to=U, , "\dashv" rotate=-90]
\end{tikzcd}
\]
is a Quillen adjunction.
\end{lemma}

\begin{proof}
The left adjoint $(-)^*$ sends degree-wise monomorphisms to degree-wise epimorphisms. Thus it preserves cofibrations. It also preserves quasi-isomorphisms. Therefore we have a Quillen adjunction.
\end{proof}

\begin{theorem}[Homotopical duality square]\label{thm: homotopical carré magique}
The following square of adjunction 
\[
\begin{tikzcd}[column sep=5pc,row sep=5pc]
\left(\mathsf{upOp}_\infty\right)^{\mathsf{op}} \arrow[r,"\mathrm{B}_\iota^{\mathsf{op}}"{name=B},shift left=1.1ex] \arrow[d,"(-)_h^\circ "{name=SD},shift left=1.1ex ]
&\left(\mathsf{curv}~\mathsf{pCoop}^{\mathsf{conil}}\right)^{\mathsf{op}}  \arrow[d,"(-)^*"{name=LDC},shift left=1.1ex ] \arrow[l,"\Omega_\iota^{\mathsf{op}}"{name=C},,shift left=1.1ex]  \\
\mathsf{upCoop}_\infty \arrow[r,"\widehat{\Omega}_\iota"{name=CC},shift left=1.1ex]  \arrow[u,"(-)^*"{name=LD},shift left=1.1ex ]
&\mathsf{curv}~\mathsf{pOp}^{\mathsf{comp}}~, \arrow[l,"\widehat{\mathrm{B}}_\iota"{name=CB},shift left=1.1ex] \arrow[u,"(-)^\vee"{name=TD},shift left=1.1ex] \arrow[phantom, from=SD, to=LD, , "\dashv" rotate=0] \arrow[phantom, from=C, to=B, , "\dashv" rotate=-90]\arrow[phantom, from=TD, to=LDC, , "\dashv" rotate=0] \arrow[phantom, from=CC, to=CB, , "\dashv" rotate=-90]
\end{tikzcd}
\] 
commutes in the following sense: right adjoints going from the top right to the bottom left are naturally isomorphic. Furthermore, this square is a square of Quillen adjunctions. 
\end{theorem}

\begin{proof}
The commutativity of this square of functors can be shown using the same arguments as in Theorem \ref{thm: carré magique}. In order to check that it is a square of Quillen adjunctions, we need to check that the adjunction 
\[
\begin{tikzcd}[column sep=7pc,row sep=3pc]
\mathsf{curv}~\mathsf{pOp}^{\mathsf{comp}}  \arrow[r, shift left=1.1ex, "(-)^\vee"{name=F}] 
& \left(\mathsf{curv}~\mathsf{pCoop}^{\mathsf{conil}}\right)^{\mathsf{op}} ~. \arrow[l, shift left=.75ex, "(-)^*"{name=U}]
            \arrow[phantom, from=F, to=U, , "\dashv" rotate=-90]
\end{tikzcd}
\]
is a Quillen adjunction. The right adjoint $(-)^*$ sends monomorphisms to epimorphisms, thus preserves fibrations. Let $f: \mathcal{C} \longrightarrow \mathcal{D}$ be a weak equivalence between two conilpotent curved partial cooperads. This is equivalent to $\Omega_\iota(f): \Omega_\iota\mathcal{C} \longrightarrow \Omega_\iota \mathcal{D}$ being an arity-wise quasi-isomorphism. Since $(-)_h^\circ$ is right Quillen, it preserves quasi-isomorphisms between cofibrant unital partial operads up to homotopy (fibrant in the opposite category). Thus $\Omega_\iota(f)^\circ: \Omega_\iota(\mathcal{C})^\circ \longrightarrow \Omega_\iota(\mathcal{D})^\circ$ is an arity-wise quasi-isomorphism. Using the commutativity of this adjunction, we get that
\[
\widehat{\mathrm{B}}_\iota(f^*): \widehat{\mathrm{B}}_\iota(\mathcal{C}^*) \longrightarrow \widehat{\mathrm{B}}_\iota(\mathcal{D}^*)~,
\]
is an arity-wise quasi-isomorphism. Which is equivalent to $f^*: \mathcal{C}^* \longrightarrow \mathcal{D}^*$ being a weak equivalence in the transferred model structure. Thus $(-)^*$ is a right Quillen functor. This concludes the proof.
\end{proof}

The above homotopical square allows us, using the curved Koszul duality established in \cite{HirshMilles12}, to compute explicit cofibrant resolutions for a vast class of complete curved absolute partial operads.

\begin{Proposition}\label{Prop: cofibrant resolutions of absolutes}
Let $(\mathcal{P},\{\circ_i\},\eta)$ be a unital partial operad which is arity-wise degree-wise finite dimensional and let $\mathcal{P}^{\mathrm{\ac}}$ be its Koszul dual conilpotent curved partial cooperad. If $\mathcal{P}$ is Koszul, then 
\[
\widehat{\Omega}(\mathcal{P}^*) \longrightarrow (\mathcal{P}^{\mathrm{\ac}})^*
\]
is a cofibrant resolution of $(\mathcal{P}^{\mathrm{\ac}})^*$ in the model category of complete curved absolute partial operads.
\end{Proposition}

\begin{proof}
Recall that $\kappa: \mathcal{P}^{\mathrm{\ac}} \longrightarrow \mathcal{P}$ is Koszul if and only if $f_\kappa: \Omega(\mathcal{P}^{\ac}) \longrightarrow \mathcal{P}$ is an arity-wise quasi-isomorphism. This is equivalent to
\[
g_\kappa: \mathcal{P}^{\ac} \longrightarrow \mathrm{B}\mathcal{P}
\]
being a weak-equivalence in the category of conilpotent curved partial cooperads, since the two model categories are Quillen equivalent. Every conilpotent curved partial cooperad is cofibrant in the model structure of Theorem \ref{thm: structure de modèles conil curved part coop} (thus fibrant in the opposite category), therefore the linear dual functor $(-)^*$ preserves all weak equivalences. This implies that
\[
(g_\kappa)^*: (\mathrm{B}\mathcal{P})^* \longrightarrow (\mathcal{P}^{\ac})^*
\]
is a weak equivalence of complete curved absolute partial operads. Using Proposition \ref{prop: finite dual commutes}, we know that there is an isomorphism
\[
\widehat{\Omega}(\mathcal{P}^*) \cong (\mathrm{B}\mathcal{P})^*
\]
of complete curved absolute partial operads, which concludes the proof.
\end{proof}

\begin{Example}[Cofibrant resolution of $c\mathcal{L}ie^\wedge$]\label{Example: Omega(ucomd)}
Let $c\mathcal{L}ie^\wedge$ be the complete curved absolute partial operad encoding curved Lie algebras of Definition \ref{def: Lie curved absolute}. There is a weak equivalence
\[
\widehat{\Omega}(\ucomd) \qi c\mathcal{L}ie^\wedge
\]
of complete curved absolute partial operads and $\widehat{\Omega}\ucomd$ is the minimal cofibrant resolution of $c\mathcal{L}ie^\wedge$ in this model structure. 
\end{Example}

\begin{Example}[Cofibrant resolution of $c\mathcal{A}ss^\wedge$]
Let $c\mathcal{A}ss^\wedge$ be the complete curved absolute partial operad encoding curved associative algebras of Definition \ref{def: Ass curved absolute}. There is a weak equivalence
\[
\widehat{\Omega}(u\mathcal{A}ss^*) \qi c\mathcal{A}ss^\wedge
\]
of complete curved absolute partial operads and $\widehat{\Omega}(u\mathcal{A}ss^*)$ is the minimal cofibrant resolution of $c\mathcal{A}ss^\wedge$ in this model structure. 
\end{Example}

\section{Homotopy transfer theorem for curved algebras}\label{Section: curved HTT}
The Homotopy Transfer Theorem is a fundamental tool in homological algebra, as it allows to transfer algebraic structures up to homotopy using contractions of dg modules. The operadic reformulation of this theorem was given in \cite[Section 10.3]{LodayVallette12} in terms of a morphism between the Bar constructions of the endomorphisms operads associated to the contraction. This morphism is called the Van der Laan morphism in \textit{loc.cit}. The notion of a contraction is \textit{homotopical} not \textit{homological}, and extends easily to the setting of pdg modules. In this section, we prove a version of the Homotopy Transfer Theorem for curved algebraic structures up to homotopy, extending the Van der Laan morphism to the complete Bar construction of Section \ref{Section: Constructions Bar-Cobar operadiques}. We recover, in the case of "pro-nilpotent" curved $\mathcal{L}_\infty$-algebras in the sense of \cite{getzler2018maurercartan}, the Homotopy Transfer Theorem constructed by Fukaya using fixed-point equations in \cite{FukayaHTT}.

\subsection{Complete filtrations on pdg modules and complete curved absolute endomorphisms operad}
A complete filtration on a pdg modules allows to endow its curved endomorphisms operad with a structure of complete curved absolute partial operad. 

\begin{Definition}[Filtered pdg module]
A \textit{filtered pdg module} $V$ amounts to the data $(V,\mathrm{F}_\bullet V ,d_V)$ of a pdg module $(V,d_V)$ together with a degree-wise decreasing filtration:
\[
V = \mathrm{F}_0V \supset \mathrm{F}_1V \supset \mathrm{F}_2V \supset \cdots \supset \mathrm{F}_n V \supset \cdots~,
\]
such that $d_V(\mathrm{F}_nV) \subset \mathrm{F}_{n+1}V$ for all $n$ in $\mathbb{N}$. Morphisms $f: V \longrightarrow W$ of filtered pdg modules are morphisms of pdg modules which are compatible with the filtrations, that is, $f(\mathrm{F}_nV) \subset \mathrm{F}_n W$ for all $n$.
\end{Definition}

\begin{Remark}
We ask that $d_V$ raises the filtration degree by one for simplicity. These are not the optimal assumptions. In fact, we only need to impose that $d_V^2(\mathrm{F}_nV) \subset \mathrm{F}_{n+1}V$. These would correspond to "gr-dg filtered modules" in \cite{JoanCurved}. In any case, notice that the associated graded of these filtrations are dg modules. 
\end{Remark}

The category of filtered pdg modules can be endowed with a closed symmetric monoidal structure by considering the following filtration on the tensor product. For $V$ and $W$ two filtered pdg modules, their tensor product $V \otimes W$ is endowed with 
\[
\mathrm{F}_n \left(V \otimes W \right) \coloneqq \sum_{p+q = n} \mathrm{Im}\left(F_pV \otimes F_qW \rightarrow V \otimes W\right)~.
\]
The internal hom functor, denoted $\mathrm{hom}(A,B)$, is given by the internal hom of pdg modules endowed with the filtration 
\[
\mathrm{F}_n \left(\mathrm{hom}(V,W)\right) \coloneqq \Big\{ ~f ~\text{in}~ \mathrm{hom}(V,W)~~|~~ f(\mathrm{F}_kV) \subset \mathrm{F}_{k+n}W\Big\} ~.
\]
\begin{Definition}[Complete pdg module]\label{def: complete pdg module}
Let $V$ be a filtered pdg module. It is a \textit{complete} pdg module if the canonical morphism 
\[
\pi_V: V \longrightarrow \lim_{n} V/\mathrm{F}_n V 
\]
is an isomorphism of filtered pdg modules. 
\end{Definition}

Complete pdg modules form a reflexive subcategory of filtered pdg modules. The reflector is simply given by the completion functor, which sends a filtered pdg module $V$ to its completion
\[
\widehat{V} \coloneqq \lim_{n} A/\mathrm{F}_nA~.
\]
By setting $V \widehat{\otimes} W \coloneqq \widehat{V \otimes W}~,$ the category of complete pdg modules is endowed with a closed symmetric monoidal structure, the internal hom being the same as the one defined above. One can define operads in this context, \cite[Section 2]{DSV18} for more details. Here, we will only use the following construction. 

\begin{Definition}[Complete curved endomorphisms partial operad]
Let $V$ be a complete pdg module. Its \textit{complete curved endomorphisms partial operad} is given by the complete pdg $\mathbb{S}$-module 
\[
\mathrm{end}_V(n) \coloneqq \mathrm{hom}^{(\geq 1)}(V^{\otimes n}, V)~,
\]
where we only consider here morphisms of pdg modules which raise the filtration degree by at least one. The operad structure is given by the partial composition of functions. Its pre-differential is given by $\partial \coloneqq [d_V,-]$ and its curvature determined by $\Theta_V(\mathrm{id}) \coloneqq d_V^2~.$
\end{Definition}

\begin{lemma}
Let $V$ be a complete pdg module. Its complete curved endomorphisms partial operad is a complete curved \textit{absolute} partial operad. 
\end{lemma}

\begin{proof}
We postpone this proof to Lemma \ref{lemma: curved endomorphisms is absolute} in the Appendix \ref{Appendix B}, where there is a detailed discussion of complete curved absolute partial operads. 
\end{proof}

\subsection{Homotopy contractions and Van der Laan morphisms}
The data of a homotopy contraction between two complete pdg modules which is compatible with their underlying filtrations induces a "Van der Laan morphism" between the complete Bar construction of their complete curved endomorphisms operads. From this, one recovers a Homotopy Transfer theorem for curved algebras over a vast class of cofibrant complete curved absolute partial operads. 

\begin{Definition}[Homotopy contraction]\label{def: contraction of pdg modules}
Let $V$ and $H$ be two complete pdg modules. A \textit{homotopy contraction} amounts to the data of 
\[
\begin{tikzcd}[column sep=5pc,row sep=3pc]
V \arrow[r, shift left=1.1ex, "p"{name=F}] \arrow[loop left]{l}{h}
&H~, \arrow[l, shift left=.75ex, "i"{name=U}]
\end{tikzcd}
\]
where $p$ and $i$ are two morphisms of filtered pdg modules and $h$ is a morphism of filtered graded modules of degree $-1$. This data satisfies the following conditions:
\[
p i = \mathrm{id}_H~, \quad i p - \mathrm{id}_V = d_V h + h d_V~.
\]
\end{Definition}

\begin{Notation}
We say that $H$ is a \textit{homotopy retract} of $V$ if there exists a homotopy contraction as in the definition above.
\end{Notation}

\begin{Proposition}[Van der Laan morphism]
Let $V$ and $H$ be two complete pdg modules such $H$ is a homotopy retract of $V$. The data of this homotopy contraction induces a morphism of dg counital partial cooperads
\[
\mathrm{VdL}: \widehat{\mathrm{B}}(\mathrm{end}_V) \longrightarrow \widehat{\mathrm{B}}(\mathrm{end}_H)~. 
\]
This morphism is called the Van der Laan morphism associated to the contraction.
\end{Proposition}

\begin{proof}
Let 
\[
\mathrm{vdl}: \mathscr{T}^\wedge(\mathrm{end}_V) \longrightarrow \mathrm{end}_H
\]
be the morphism of pdg $\mathbb{S}$-modules which sends a series of rooted trees labeled by elements of $\mathrm{end}_V$ to the converging series in $\mathrm{end}_H$ given by applying the Van der Laan map to each rooted tree of the series. For a recollection on the standard Van der Laan map, see \cite[Section 10.3.2]{LodayVallette12}. One can restrict $\mathrm{vdl}$ along the inclusion of the cofree counital partial cooperad inside the completed tree endofunctor. By the universal property of the cofree counital partial cooperad, this induces a morphism of counital partial cooperads:
\[
\mathrm{VdL}: \mathscr{T}^\vee(\mathrm{end}_V) \longrightarrow \mathscr{T}^\vee(\mathrm{end}_H)~.
\]
Let us show that $\mathrm{VdL}$ commutes with the differentials. It equivalent to the following diagram commuting 
\[
\begin{tikzcd}[column sep=3pc,row sep=3pc]
\mathscr{T}^\vee(\mathrm{end}_V) \arrow[r,"\mathrm{VdL}"] \arrow[d,"d_{\mathrm{bar}} ",swap] 
&\mathscr{T}^\vee(\mathrm{end}_H) \arrow[d,"\psi_1 + \psi_2"] \\
\mathscr{T}^\vee(\mathrm{end}_V) \arrow[r,"\mathrm{vdl}"] 
&\mathrm{end}_H~,
\end{tikzcd}
\]
where $\psi_1 + \psi_2$ is the map that induces the differential of the complete Bar construction as its unique coderivation extending it. One can show that this diagram commutes for elements in $\mathscr{T}^\vee(\mathrm{end}_V)$ which do not start with the trivial tree $|~$ extending the same computations as in \cite[Section 10.3.2]{LodayVallette12} to infinite series of rooted trees. For the trivial tree $|~$, a small computation shows that:
\[
\mathrm{vdl} \cdot d_{\mathrm{bar}}(|) = p \cdot d_V^2 \cdot i = d_H \cdot p \cdot i \cdot d_H = d_H^2 = (\psi_1 + \psi_2) \cdot \mathrm{VdL}(|)~, 
\]
using the fact that $p$ and $i$ are morphisms of pdg modules. 
\end{proof}

\begin{theorem}[Curved Homotopy Transfer Theorem]\label{thm: curved HTT}
Let $\mathcal{C}$ be a dg counital partial cooperad. Let $V$ and $H$ be two complete pdg modules such $H$ is a homotopy retract of $V$. Then any curved $\widehat{\Omega}(\mathcal{C})$-algebra structure on $V$ can be transferred along the homotopy contraction to a curved $\widehat{\Omega}(\mathcal{C})$-algebra structure on $H$. 
\end{theorem}

\begin{proof}
One has that
\[
\mathrm{Hom}_{\mathsf{curv}~\mathsf{abs.pOp}^{\mathsf{comp}}}(\widehat{\Omega}\C,\mathrm{end}_V) \cong \mathrm{Tw}(\C,\mathrm{end}_V) \cong \mathrm{Hom}_{\mathsf{dg}~\mathsf{upCoop}}(\C,\widehat{\mathrm{B}}(\mathrm{end}_V))~,
\]
thus a curved $\widehat{\Omega}(\mathcal{C})$-algebra structure on $V$ is equivalent to a morphism of dg counital partial cooperads 
\[
\phi: \mathcal{C} \longrightarrow \widehat{\mathrm{B}}(\mathrm{end}_V)~.
\]
Any such morphism can be pushed forward by the Van der Laan morphism
\[
\mathrm{VdL} \cdot \phi: \mathcal{C} \longrightarrow \widehat{\mathrm{B}}(\mathrm{end}_H)~,
\]
thus inducing a curved $\widehat{\Omega}(\mathcal{C})$-algebra structure on $H$. 
\end{proof}

\begin{Remark}
All cofibrant complete curved absolute partial operads can be written as $\widehat{\Omega}(\mathcal{C})$, where $\mathcal{C}$ is a counital partial cooperad \textit{up to homotopy}. By restricting to the case where $\C$ is a dg counital partial cooperad, we are restricting to the case of Koszul resolutions. See Proposition \ref{Prop: cofibrant resolutions of absolutes}. In order to prove the statement for all cofibrant complete curved absolute partial operad, one should prove that a homotopy contraction induces a morphism between the complete Bar constructions \textit{relative to $\iota$} of Proposition \ref{Prop: complete Bar-Cobar adjunction relative to iota}.
\end{Remark}

Let us conclude with some examples of how this framework can be applied. Let $\widehat{\Omega}(\ucomd)$ be the complete curved absolute partial operad of Example \ref{Example: Omega(ucomd)}.

\begin{Example}
Let $V$ be a complete dg modules. Then a curved $\widehat{\Omega}(\ucomd)$-algebra structure on $V$ corresponds to a \textit{pro-nilpotent curved} $\mathcal{L}_\infty$-algebra structure on $V$ in the sense of \cite{getzler2018maurercartan}.
\end{Example}

\begin{Example}[Fukaya's Homotopy Transfer Theorem]
Let $V$ and $H$ be two complete dg modules and let 
\[
\begin{tikzcd}[column sep=5pc,row sep=3pc]
V \arrow[r, shift left=1.1ex, "p"{name=F}] \arrow[loop left]{l}{h}
&H~, \arrow[l, shift left=.75ex, "i"{name=U}]
\end{tikzcd}
\]
be homotopy contraction. Given a pro-nilpotent curved $\mathcal{L}_\infty$-algebra structure on $V$, one can transfer it to $H$ using a version of the Homotopy Transfer Theorem given in \cite{FukayaHTT}. The transferred structure is given as the solution of a fixed-point equation. If one computes the solution of this fixed point equation using the methods of \cite[Appendix A]{robertnicoud2020higher}, one obtains the same formulae for the transferred structure as the ones obtained from Theorem \ref{thm: curved HTT}.
\end{Example}

\newpage

\section*{Appendix: What is an absolute partial operad ?}\label{Appendix B}
The notion of an \textit{absolute partial operad} is an example of a vast class of algebraic structures that emerge as algebras over a cooperad. See Chapter \ref{Chapter 1} for more examples of these absolute structures. In our particular case, absolute partial operads appear as algebras over the conilpotent partial $\mathbb{S}$-colored cooperad $\mathcal{O}^*$. Here $\mathcal{O}^*$ denotes the linear dual of the partial $\mathbb{S}$-colored operad $\mathcal{O}$ which encodes partial operads as its algebras. The goal of this appendix is to give a somewhat explicit definition of what absolute partial operads are and how to characterize them. Then to compare them with standard partial operads. Afterwards, we will generalize these results to the curved case. 

\subsection{Absolute partial operads}
We work in the underlying category of $\mathbb{S}$-modules for simplicity. This subsection admits a straightforward generalization to dg $\mathbb{S}$-modules or pdg $\mathbb{S}$-modules.

\begin{lemma}\label{lemma: iso avec complete tree monad}
There is an isomorphism of endofunctors in the category of $\mathbb{S}$-modules:
\[
\widehat{\mathscr{S}}_{\mathbb{S}}^c (\mathcal{O}^*)(-) \cong \overline{\mathscr{T}}^\wedge(-)~,
\]
where $\overline{\mathscr{T}}^\wedge(-)$ is the reduced completed tree endofunctor.
\end{lemma}

\begin{proof}
We have that:
\begin{align*}
\widehat{\mathscr{S}}_{\mathbb{S}}^c (\mathcal{O}^*)(M)(n) \coloneqq \prod_{(n_1,\cdots,n_r) \in \mathbb{N}^{r}} \mathrm{Hom}_{\left(\mathbb{S}_{n_1} \times \cdots \times \mathbb{S}_{n_r}\right)~\wr~ \mathbb{S}_r } \left( \mathcal{O}^*(n_1,\cdots,n_r;n), M(n_1) \otimes \cdots \otimes M(n_r) \right) \cong \\
\cong \prod_{(n_1,\cdots,n_r) \in \mathbb{N}^{r}} \mathcal{O}(n_1,\cdots,n_r;n) \otimes_{\left(\mathbb{S}_{n_1} \times \cdots \times \mathbb{S}_{n_r}\right)~\wr~ \mathbb{S}_r }  M(n_1) \otimes \cdots \otimes M(n_r)~,
\end{align*}
since each $\mathcal{O}^*(n_1,\cdots,n_r;n)$ is finite dimensional over $\kk$ and that we are working over a field of characteristic, thus invariants and coinvariants turn out to be canonically isomorphic. Using the same bijection as in the proof of Lemma \ref{lemma with the trees bijection}, we identify this last term with the reduced completed tree endofunctor.
\end{proof}

\begin{Corollary}
There is a monad structure on the reduced completed tree endofunctor $\overline{\mathscr{T}}^\wedge(-)$. Its structural morphism 
\[
\mathsf{Sub}^\wedge: \overline{\mathscr{T}}^\wedge \circ \overline{\mathscr{T}}^\wedge (-) \longrightarrow \overline{\mathscr{T}}^\wedge(-)
\]
is given by the substitution of infinite series of rooted trees. 
\end{Corollary}

\begin{proof}
Let $M$ be an $\mathbb{S}$-module. An element of $\overline{\mathscr{T}}^\wedge \circ \overline{\mathscr{T}}^\wedge(M)$ amounts to the data of a series of rooted trees which are labeled in the following way: if $v$ is a vertex of with $k$ inputs of a rooted tree $\tau$, then it is labeled by a series of rooted trees of arity $k$ which are themselves labeled by elements of $M$ in the usual way. Let us describe the monad structure transported via the isomorphism of Lemma \ref{lemma: iso avec complete tree monad}. This monad structure is given by first distributing all the labeling series and then by substituting each of the node by the rooted tree it is labeled with.
Pictorially, the substitution is given by 

\begin{center}
\includegraphics[width=130mm,scale=1]{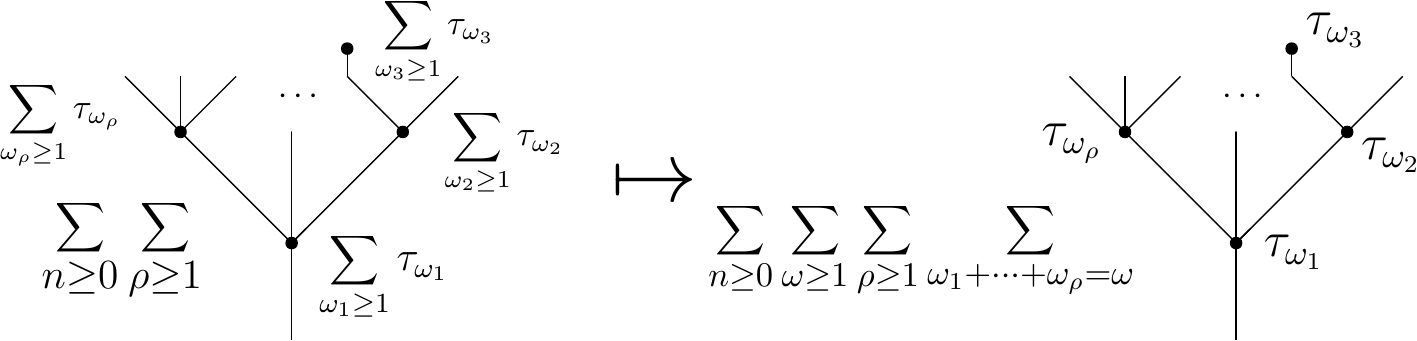}
\end{center}

where one should substitute each node by the corresponding labeling rooted tree $\tau_{\omega_i}$. 
\end{proof}

Thus we can rewrite the definition of an $\mathcal{O}^*$-algebra in terms of the reduced completed tree monad. 

\begin{Definition}[Absolute partial operad]\label{def: absolute partial operads}
An \textit{absolute partial operad} $\mathcal{Q}$ amounts to the data $(\mathcal{Q},\gamma_\Q)$ of an algebra structure over the reduced completed tree monad 
\[
\gamma_\Q: \overline{\mathscr{T}}^\wedge (\Q) \longrightarrow \Q~.
\]
\end{Definition}

The goal of this section is to make sense of this definition.

\begin{Remark}
The structural morphism $\gamma_\Q$ of an absolute partial operad gives a well-defined composition of any \textit{infinite sums} for rooted trees labeled by elements of $\mathcal{Q}$, \textit{without presupposing any underlying topology} on $\mathcal{Q}$. Let 
\[
\sum_{n \geq 0} \sum_{\omega \geq 1} \sum_{\tau \in \mathrm{RT}_n^\omega} \tau
\]
be a infinite sum of rooted trees with vertices labeled by elements of $\Q$. Notice that, in general:
\[
\gamma_\Q \left(\sum_{n \geq 0} \sum_{\omega \geq 1} \sum_{\tau \in \mathrm{RT}_n^\omega} \tau \right) = \sum_{n \geq 0} \gamma_\Q \left( \sum_{\omega \geq 1} \sum_{\tau \in \mathrm{RT}_n^\omega} \tau \right) \neq \sum_{n \geq 0} \sum_{\omega \geq 1} \sum_{\tau \in \mathrm{RT}_n^\omega} \gamma_\Q(\tau)~,
\]
as the later sum is not even well-defined in $\Q$ in general.
\end{Remark}

\begin{Proposition}
Let $M$ be an $\mathbb{S}$-module such that $M(0) = M(1) = 0$, then the data of a partial operad structure on $M$ is equivalent to the data of an absolute partial operad structure on $M$. 
\end{Proposition}

\begin{proof}
In this case, there are only a finite amount of rooted tree labeled by elements of $M$ for each possible arity, thus the reduced completed tree monad coincides with the reduced tree monad. 
\end{proof}

\begin{lemma}
Let $(\mathcal{Q},\gamma_\Q)$ be an absolute partial operad. The structure map $\gamma_\Q$ induces a partial operad structure on $\Q$. This defines a faithful functor
\[
\begin{tikzcd}[column sep=4pc,row sep=0pc]
\mathsf{Res}: \mathsf{abs}~\mathsf{pOp} \arrow[r]
&\mathsf{pOp}~.
\end{tikzcd}
\]
\end{lemma}

\begin{proof}
It is straightforward to check that the inclusion of endofunctors
\[
\overline{\mathscr{T}}(-) \hookrightarrow \overline{\mathscr{T}}^\wedge(-)
\]
is in fact an inclusion of monads. Thus, by pulling back along this inclusion, any absolute partial operad structure gives a partial operad structure. This pullback amounts to restrict $\gamma_\Q$ to finite sums of rooted trees inside the reduced completed tree monad.
\end{proof}

\begin{Remark}
Let $(\mathcal{Q},\gamma_\Q)$ be an absolute partial operad. It is in particular endowed with partial composition maps $\{\circ_i\}$ which satisfy the usual axioms of a partial operad. These maps are \textit{part of the structure} of an absolute partial operad.
\end{Remark}

\begin{Proposition}[Absolute envelope of a partial operad]
There is an adjunction
\[
\begin{tikzcd}[column sep=7pc,row sep=3pc]
\mathsf{pOp}  \arrow[r, shift left=1.1ex, "\mathsf{Abs}"{name=F}] 
&\mathsf{abs.pOp}~. \arrow[l, shift left=.75ex, "\mathsf{Res}"{name=U}] \arrow[phantom, from=F, to=U, , "\dashv" rotate=-90]
\end{tikzcd}
\]
We call the left adjoint \textit{the absolute envelope} of a partial operad, and we denote it by $\mathsf{Abs}$.
\end{Proposition}

\begin{proof}
Notice that both categories are categories of algebras over accessible monads, thus are presentable. Since $\mathsf{Res}$ is accessible and preserves all limits, it has a left adjoint, by \cite[Theorem 1.66]{AdamekRosicky}.
\end{proof}

There is a canonical topology on absolute partial operads induced by the structural morphism.

\begin{Definition}[Canonical filtration on an absolute partial operad]\label{def: canonical topology on absolutes}
Let $(\mathcal{Q},\gamma_\Q)$ be an absolute partial operad. The \textit{canonical filtration} is the decreasing filtration given by
\[ 
\overline{\mathscr{F}}_{\omega} \Q \coloneqq \mathrm{Im}\left(\gamma_\Q^{(\geq \omega)}: \overline{\mathscr{T}}^{\wedge~(\geq \omega)}(\Q) \longrightarrow \Q \right)
\]
for all $\omega \geq 0$, where $\overline{\mathscr{T}}^{\wedge~(\geq \omega)}(\Q)$ denotes the possibly infinite sums of rooted trees that have terms of at least $\omega$ internal edges. Each $\overline{\mathscr{F}}_{\omega} \Q$ defines an ideal of $\Q$. Notice that: 
\[
\Q = \overline{\mathscr{F}}_{0} \Q \supseteq \overline{\mathscr{F}}_{1} \Q \supseteq \cdots \supseteq \overline{\mathscr{F}}_{\omega} \Q \supseteq \cdots.
\]
\end{Definition}

\begin{Definition}[Nilpotent absolute partial operad]
Let $(\mathcal{Q},\gamma_\Q)$ be an absolute partial operad. It is said to be \textit{nilpotent} if there exists an $\omega \geq 1$ such that 
\[
\Q/\overline{\mathscr{F}}_{\omega} \Q \cong \Q~.
\]
The absolute partial operad is said to be $\omega_0$\textit{-nilpotent} if $\omega_0$ is the smallest integer such that the above isomorphism exits.
\end{Definition}

\begin{Definition}[Completion of an absolute partial operad]\label{def: completion functor for absolute operads}
Let $(\mathcal{Q},\gamma_\Q)$ be an absolute partial operad, its \textit{completion} $\widehat{\Q}$ is given by the following limit
\[
\widehat{\Q} \coloneqq \lim_{\omega} \Q/\overline{\mathscr{F}}_{\omega} \Q
\]
taken in the category of absolute partial operads.
\end{Definition}

The completion is functorial. There is a canonical morphism of absolute partial operads
\[
\psi: \Q \twoheadrightarrow \widehat{\Q}~,
\]
which is always an epimorphism. See Remark \ref{Remark: varphi is an epimorphism} as to why this is the case in the general setting.

\begin{Definition}[Complete absolute partial operad]\label{def: complete absolute operads}
Let $(\mathcal{Q},\gamma_\Q)$ be an absolute partial operad. It is \textit{complete} if the canonical morphism $\psi: \Q \twoheadrightarrow \widehat{\Q}$ is an isomorphism of absolute partial operads. 
\end{Definition} 

\begin{Example}
Any nilpotent absolute partial operad is also a complete absolute partial operad. Any complete absolute partial operad is the limit of a tower of nilpotent absolute partial operads. Any free absolute partial operad is also complete.
\end{Example}

\begin{Remark}
An absolute partial operad is complete if and only if the topology induced by its canonical filtration is Hausdorff.
\end{Remark}

\begin{Remark}
Let $(\Q, \gamma_\Q)$ be a \textit{complete} absolute partial operad and let 
\[
\sum_{n \geq 0} \sum_{\omega \geq 1} \sum_{\tau \in \mathrm{RT}_n^\omega} \tau
\]
be a infinite sum of rooted trees with vertices labeled by elements of $\Q$. In this case we have that
\[
\gamma_\Q \left(\sum_{n \geq 0} \sum_{\omega \geq 1} \sum_{\tau \in \mathrm{RT}_n^\omega} \tau \right) = \sum_{n \geq 0} \sum_{\omega \geq 1} \sum_{\tau \in \mathrm{RT}_n^\omega} \gamma_\Q(\tau)~,
\]
since $\gamma_\Q$ commutes with the sum over the weight of rooted trees and since there are only a finite amount of rooted trees of arity $n$ and of weight $\omega$.
\end{Remark}

\begin{Proposition}
There is an isomorphism of categories between the category of $\omega_0$-nilpotent absolute partial operads and the category of $\omega_0$-nilpotent partial operads, for all $\omega_0 \geq 1$. 
\end{Proposition}

\begin{proof}
Let $\overline{\mathscr{T}}^{(\leq \omega_0)}$ denote the reduce tree monad truncated at rooted trees of with more than $\omega_0$ internal edges. The data of an $\omega_0$-nilpotent partial operad amounts to the data of an algebra over the reduced $\omega_0$-truncated tree monad. The data of an $\omega_0$-nilpotent absolute partial operad also amounts to the data of an algebra over the reduced $\omega_0$-tree monad. Indeed, since there are only a finite number of rooted trees of weight less than $\omega_0$ at any given arity, the direct sum and the product coincide.
\end{proof}

The notion of a complete absolute partial operad appears naturally when one takes the linear dual of a conilpotent partial cooperad.

\begin{Proposition}\label{prop: linear dual of a conil coop}
Let $(\C, \{\Delta_i\})$ be a conilpotent partial cooperad. Then its linear dual $\C^*$ has an induced absolute partial operad structure. Furthermore, since $\C$ is conilpotent, then $\C^*$ is a complete absolute partial operad. It defines a functor
\[
\begin{tikzcd}[column sep=7pc,row sep=3pc]
\left(\mathsf{pCoop}^{\mathsf{conil}}\right)^{\mathsf{op}}  \arrow[r, "(-)^*"{name=F}] 
&\mathsf{abs.pOp}^{\mathsf{comp}}~. 
\end{tikzcd}
\]
\end{Proposition}

\begin{proof}
Let 
\[
\Delta_\C: \C \longrightarrow \overline{\mathscr{T}}^c(\C)
\]
denote the structural map of $\C$ as a coalgebra over the reduced tree comonad, given by Proposition \ref{Prop: conil coop = coalgebra over tree comonad}. By taking the linear dual we obtain a map
\[
\Delta_\C^*: \left(\overline{\mathscr{T}}^c(\C)\right)^* \longrightarrow \C^*~.
\]
There is a canonical inclusion
\[
\iota_\C:  \overline{\mathscr{T}}^{\wedge}(\C^*) \hookrightarrow \left(\overline{\mathscr{T}}^c(\C)\right)^*
\]
which is an isomorphism if and only if $\C$ is arity-wise finite dimensional. Thus by pulling back along $\iota_\C$, we obtain a map $\gamma_\C \coloneqq \Delta_\C^* \cdot \iota_\C$. One can check that this endows $\C^*$ with the structure of an algebra over the reduced complete tree monad. Since $\C$ is conilpotent, it can be written as 
\[
\C \cong \colim_{\omega} \mathscr{R}_\omega \C~,
\]
where $\mathscr{R}_\omega \C$ is the sub-cooperad given by operations which admit non-trivial $\omega$-iterated partial decompositions. Thus 
\[
\C^* \cong \lim_{\omega} \left(\mathscr{R}_\omega \C\right)^*~,
\]
where one can check that
\[
\left(\mathscr{R}_\omega \C\right)^* \cong \C^* /\overline{\mathscr{F}}_{\omega} \C^*~.
\]
Therefore the resulting absolute partial operad $\C^*$ is indeed complete.
\end{proof}

\begin{Proposition}\label{prop: adjonction topo dual lin}
The linear functor admits a left adjoint
\[
\begin{tikzcd}[column sep=7pc,row sep=3pc]
\left(\mathsf{pCoop}^{\mathsf{conil}}\right)^{\mathsf{op}}  \arrow[r,shift left=1.1 ex, "(-)^*"{name=F}] 
&\mathsf{abs.pOp}^{\mathsf{comp}}~. \arrow[l, shift left=.75ex, "(-)^\vee"{name=U}] \arrow[phantom, from=F, to=U, , "\dashv" rotate=90]
\end{tikzcd}
\]
We denote its left adjoint $(-)^\vee$ and call it the \textit{topological dual functor}.
\end{Proposition}

\begin{proof}
Consider the following square of functors
\[
\begin{tikzcd}[column sep=4pc,row sep=4pc]
\left(\mathsf{pCoop}^{\mathsf{conil}}\right)^{\mathsf{op}} \arrow[r,"(-)^*"]  \arrow[d,"\mathrm{U}^\mathsf{op}"{name=SD},shift left=1.1ex ]
&\mathsf{abs.pOp}^{\mathsf{comp}} \arrow[d,"\mathrm{U}"{name=LDC},shift left=1.1ex ] \\
\mathbb{S}\text{-}\mathsf{mod}^{\mathsf{op}} \arrow[r,"(-)^*"{name=CC},,shift left=1.1ex] \arrow[u,"\left(\overline{\mathscr{T}}^c(-)\right)^\mathsf{op}"{name=LD},shift left=1.1ex ] \arrow[phantom, from=SD, to=LD, , "\dashv" rotate=0]
&\mathbb{S}\text{-}\mathsf{mod}~.  \arrow[l,"(-)^*"{name=CB},shift left=1.1ex] \arrow[u,"\overline{\mathscr{T}}^\wedge(-)"{name=TD},shift left=1.1ex] \arrow[phantom, from=TD, to=LDC, , "\dashv" rotate=0] \arrow[phantom, from=CC, to=CB, , "\dashv" rotate=90]
\end{tikzcd}
\] 
The adjunction on the left hand side is monadic. All these categories are bicomplete. It is clear that $(-)^* \cdot \mathrm{U}^\mathsf{op} \cong \mathrm{U} \cdot (-)^*~.$ Thus we can apply the Adjoint Lifting Theorem \cite[Theorem 2]{AdjointLifting}, which concludes the proof. 
\end{proof}

\begin{Remark}
There is an explicit description of this left adjoint. Let $(\Q,\gamma_Q)$ be a complete absolute partial operad, its topological dual $\Q^\vee$ is given by the equalizer
\[
\begin{tikzcd}[column sep=4pc,row sep=4pc]
\mathrm{Eq}\Bigg(\overline{\mathscr{T}}^c(\mathcal{Q}^*) \arrow[r,"(\gamma_\Q)^*",shift right=1.1ex,swap]  \arrow[r,"\varrho"{name=SD},shift left=1.1ex ]
&\overline{\mathscr{T}}^c\left((\overline{\mathscr{T}}^\wedge(\mathcal{P}))^*\right) \Bigg)~,
\end{tikzcd}
\]
where $\varrho$ is a map constructed using the comonad structure of $\overline{\mathscr{T}}^c$ and the canonical inclusion into the double linear dual.
\end{Remark}

Lastly, one can compare the canonical filtrations induced by the structure of an absolute partial operad and by the structure of a partial operad. In general, they do not agree.

\begin{lemma}
Let $(\mathcal{Q},\gamma_\Q)$ be an absolute partial operad. Let $(\mathsf{Res}(\Q), \mathsf{Res}(\gamma_\Q))$ the partial operad given by its restriction. Then there is an inclusion of $\mathbb{S}$-modules:
\[
\mathscr{F}_{\omega}\mathsf{Res}(\Q) \subseteq \overline{\mathscr{F}}_{\omega} \Q
\]
for all $\omega \geq 0$, between its canonical filtration as a partial operad and its canonical filtration as an absolute partial operad. 
\end{lemma}

\begin{proof}
Elements in $\overline{\mathscr{F}}_{\omega} \Q$ are given by the images of the structural morphism $\gamma_\Q$ of series of rooted trees of at least weight $\omega$. Since any finite sum of at least weight $\omega$ is a particular example of this, any element in $\mathscr{F}_{\omega}\mathsf{Res}(\Q)$ is also included in $\overline{\mathscr{F}}_{\omega} \Q$.
\end{proof}

\begin{Corollary}
Let $(\mathcal{Q},\gamma_\Q)$ be an absolute partial operad. Then the canonical topology of $(\mathsf{Res}(\Q),\allowbreak \mathsf{Res}(\gamma_\Q))$ as a partial operad is Hausdorff.
\end{Corollary}

\begin{proof}
We have that
\[
\bigcap_{\omega \in \mathbb{N}} \mathscr{F}_{\omega}\mathsf{Res}(\Q) \subset \bigcap_{\omega \in \mathbb{N}} \overline{\mathscr{F}}_{\omega} \Q = 0~,
\]
hence the canonical topology on $\mathsf{Res}(\Q)$ is indeed Hausdorff.
\end{proof}

\begin{Counterexample}
Let $\overline{\mathscr{T}}^\wedge(M)$ be the free complete absolute partial operad generated by an $\mathbb{S}$-module $M$. Then 
\[
\mathsf{Res}\left(\overline{\mathscr{T}}^\wedge(M)\right)
\]
is not complete as a partial operad. 
\end{Counterexample}

\subsection{Complete curved absolute partial operads}
Recall Definition \ref{def: cO vee} of the conilpotent curved $\mathbb{S}$-colored partial cooperad $c\mathcal{O}^\vee$ which encodes conilpotent curved partial cooperads as its coalgebras. The goal of this subsection is to understand what algebras over this $\mathbb{S}$-colored cooperad are. For this purpose, we now consider pdg $\mathbb{S}$-modules as our ground category.

\begin{lemma}\label{lemma: dual Schur of cOvee}
Let $M$ be an pdg $\mathbb{S}$-module. There is an isomorphism of pdg $\mathbb{S}$-modules
\[
\widehat{\mathscr{S}}_\mathbb{S}^c(c\mathcal{O}^\vee)(M) \cong \overline{\mathscr{T}}^\wedge(M \oplus \nu)~,
\]
where $\nu$ is an arity $1$ generator of degree $-2$. This isomorphism is natural in $M$. 
\end{lemma}

\begin{proof}
This is a straightforward extension of Lemma \ref{lemma: iso avec complete tree monad}, using a similar bijection to the one defined in the proof of Proposition \ref{Prop: cO-cogebres.}.
\end{proof}

\begin{Definition}[Complete curved absolute partial operad]
A \textit{complete curved absolute partial operad} $(\Q,\gamma_\Q , d_\Q , \Theta_\Q)$ amounts to the data of $(\Q,\gamma_\Q ,d_\Q )$ a complete pdg absolute partial operad, as defined in the preceding subsection, endowed with a morphism of pdg $\mathbb{S}$-modules $\Theta_\Q: \I \longrightarrow \overline{\mathscr{F}}_{1}\Q$ of degree $-2$ such that the following diagram commutes
\[
\begin{tikzcd}[column sep=7.5pc,row sep=3pc]
\Q \arrow[r,"\mathrm{diag}"] \arrow[rrd,"(d_\Q)^2", bend right =10]
&\Q \oplus \Q \cong (\I \circ \Q) \oplus (\PP \circ \I) \arrow[r,"(\Theta_\Q~ \circ~ \mathrm{id})~ -~(\mathrm{id}~ \circ'~ \Theta_\Q)"] 
& \Q \circ_{(1)} \Q \arrow[d,"\gamma_{(1)}"]\\
&
&\Q~,
\end{tikzcd}
\]
where $\mathrm{diag}$ is given by $\mathrm{diag}(\mu) \coloneqq (\mu,\mu)$.
\end{Definition}

\begin{Proposition}\label{Prop: cO algebres vrai}
The category of complete curved $c\mathcal{O}^\vee$-algebras is equivalent to the category of complete curved absolute partial operads.
\end{Proposition}

\begin{proof}
Let $(\Q,\gamma_\Q , d_\Q , \Theta_\Q)$ be a complete curved absolute partial operad. One can extend $\gamma_\Q$ into a morphism of pdg $\mathbb{S}$-modules 
\[
\gamma_\Q^+: \overline{\mathscr{T}}^\wedge(\Q \oplus \nu) \longrightarrow \Q
\]
as follows 
\begin{enumerate}
\item It sends $\nu$ to $\Theta_\Q(\mathrm{id})$ in $\overline{\mathscr{F}}_{1}\Q~.$

\item It sends a rooted tree $\tau$ with vertices labeled by elements of $\Q$ and possibly to containing some unary vertices labeled by $\nu$ to the corresponding compositions of operations in $\Q$ where the unary vertices labeled by $\nu$ are replaced by $\Theta_\Q$. 

\item It defined on infinite sums of rooted trees in $\overline{\mathscr{T}}^\wedge(\Q \oplus \nu)$ as follows
\[
\sum_{n \geq 0} \sum_{\omega \geq 1} \sum_{\tau \in \mathrm{RT}_n^\omega} \tau \mapsto \sum_{n \geq 0} \sum_{\omega \geq 1} \sum_{\tau \in \mathrm{RT}_n^\omega} \gamma_\Q^+(\tau)~,
\]
which is well-defined since $\Theta_\Q(\mathrm{id})$ is in $\overline{\mathscr{F}}_{1}\Q$ and since $\Q$ is a complete absolute partial operad.
\end{enumerate}
This endows $\Q$ with the structure of a pdg $c\mathcal{O}^\vee$-algebra. Furthermore, one checks that it is complete for its canonical filtration as a $c\mathcal{O}^\vee$-algebra, since both filtrations are the same. It is straightforward to check that 
\[
\begin{tikzcd}[column sep=4.5pc,row sep=3pc]
\Q \cong \widehat{\mathscr{S}}_{\mathbb{S}}^c(\I_\mathbb{S})(\Q) \arrow[r,"\widehat{\mathscr{S}}_{\mathbb{S}}^c(\Theta_{c\mathcal{O}^\vee})(\mathrm{id}) "] \arrow[rd,"- d_\Q^2",swap]
&\widehat{\mathscr{S}}_{\mathbb{S}}^c(c\mathcal{O}^\vee)(\Q) \arrow[d,"\gamma_\Q"]\\
&\Q 
\end{tikzcd}
\]
commutes, making it a complete curved $c\mathcal{O}^\vee$-algebra. The other way around, let $(\Q,\gamma_\Q^+,d_\Q)$ be a complete curved $c\mathcal{O}^\vee$-algebra. Restricting $\gamma_\Q^+$ along the obvious inclusion
\[
\overline{\mathscr{T}}^\wedge(\Q) \hookrightarrow \overline{\mathscr{T}}^\wedge(\Q \oplus \nu)~,
\]
endows $\Q$ with an absolute partial operad structure. Furthermore, since its canonical filtration as an absolute partial operad is contained in its canonical filtration as a $c\mathcal{O}^\vee$-algebra, this implies that $\Q$ is a complete absolute partial operad. Define $\Theta_\Q(\mathrm{id}) \coloneqq \gamma_\Q^+(\nu)$, and the commutative of the above triangle implies that $\Q$ forms indeed a complete curved absolute partial operad.
\end{proof}

\begin{lemma}\label{lemma: dual lin d'une curved conil coop}
Let $(\mathcal{C},\{\Delta_i\},d_\C,\Theta_\C)$ be a conilpotent curved partial cooperad. Then its linear dual $\C^*$ has inherits a structure of a complete curved absolute partial operad. This defines a functor 
\[
\begin{tikzcd}[column sep=4pc,row sep=0pc]
\left(\mathsf{curv}~\mathsf{pCoop}^{\mathsf{conil}}\right)^{\mathsf{op}} \arrow[r,"(-)^*"]
&\mathsf{curv}~\mathsf{abs.pOp}^{\mathsf{comp}}~.
\end{tikzcd}
\]
\end{lemma}

\begin{proof}
We know that $\C^*$ is a complete absolute partial operad from Proposition \ref{prop: linear dual of a conil coop}. Therefore $\Theta_\C^*$ endows $\C^*$ with a complete curved absolute partial operad structure.
\end{proof}

\begin{Proposition}\label{prop: adjonction topo et dual lin en curved}
There is an adjunction 
\[
\begin{tikzcd}[column sep=7pc,row sep=3pc]
\mathsf{curv}~\mathsf{pOp}^{\mathsf{comp}}  \arrow[r, shift left=1.1ex, "(-)^\vee"{name=F}] 
& \left(\mathsf{curv}~\mathsf{pCoop}^{\mathsf{conil}}\right)^{\mathsf{op}} ~. \arrow[l, shift left=.75ex, "(-)^*"{name=U}]
            \arrow[phantom, from=F, to=U, , "\dashv" rotate=-90]
\end{tikzcd}
\]
\end{Proposition}

\begin{proof}
This is a particular case of Proposition \ref{prop: adjonction topo dual lin}. The only thing to prove is that if $(\Q,\gamma_\Q , d_\Q , \allowbreak \Theta_\Q)$ is a complete curved absolute partial operad, its topological dual $\Q^\vee$ conilpotent partial cooperad is indeed a conilpotent curved partial operad. Notice that by analogues arguments to Proposition \ref{prop: natural mono for topo dual}, we have that there is a natural monomorphism of pdg $\mathbb{S}$-modules 
\[
(-)^\vee \hookrightarrow (-)^*~,
\]
thus $\Q^\vee$ is a sub-object of $\Q^*$ and its pre-differential is induced by restriction. This in turn implies that $\Theta_\Q^\vee$ endows $\Q^\vee$ with a conilpotent curved partial cooperad structure. See Proposition \ref{prop: restriction aux curved} for an analogue statement.
\end{proof}

\subsection{Examples}
Here are some examples of complete curved absolute partial operads. As first examples, let us construct the "absolute analogues" of the curved operads $\Liec$ and $\Assc$ constructed in Section \ref{Section: Curved Operads}.

\medskip

Let $M$ be the pdg $\mathbb{S}$-module given by $(\mathbb{K}.\zeta,0, \mathbb{K}.\beta,0,\cdots)$ with zero pre-differential, where $\zeta$ is an arity $0$ operation of degree $-2$, and $\beta$ is an arity $2$ operation of degree $0$, basis of the signature representation of $\mathbb{S}_2$. 

\begin{Definition}[$\Liec^\wedge$ absolute operad]\label{def: Lie curved absolute}
The \textit{complete curved absolute partial operad} $\Liec^\wedge$ is given by the free pdg absolute partial operad generated by $M$ modulo the ideal generated by the Jacobi relation on the generator $\beta$. It is endowed with the curvature $\Theta$ given by $\Theta(\mathrm{id}) \coloneqq \beta \circ_1 \zeta~.$  
\end{Definition}

\begin{lemma}
The data $(\Liec^\wedge, 0, \Theta)$ forms a complete curved absolute partial operad. 
\end{lemma}

\begin{proof}
One can show by direct computation that this absolute partial operad is complete. The rest of the proof is analogous to Lemma \ref{lemmalie}.
\end{proof}

\begin{Proposition}
Let $u\mathcal{C}om$ be the unital partial operad encoding unital commutative algebras and let $u\mathcal{C}om^{\ac}$ be its Koszul dual conilpotent curved partial cooperad. There is an isomorphism of complete curved partial operads
\[
\left(u\mathcal{C}om^{\ac}\right)^* \cong \Liec^\wedge.
\]
\end{Proposition}

\begin{proof}
By direct inspection.
\end{proof}

Let $N$ be the pdg $\mathbb{S}$-module given by $(\mathbb{K}.\phi,0,\mathbb{K}[\mathbb{S}_2].\mu,0,\cdots)$ with zero pre-differential, where $\phi$ is an arity $0$ operation of degree $-2$, and $\mu$ is an binary operation of degree $0$, basis of the regular representation of $\mathbb{S}_2$.

\begin{Definition}[$\Assc^\wedge$ absolute operad]\label{def: Ass curved absolute}
The \textit{complete curved absolute partial operad} $\Assc^\wedge$ is given by the free pdg absolute partial operad generated by $N$ modulo the ideal generated by the associativity relation on the generator $\mu$. It is endowed with the curvature $\Theta$ given by $\Theta(\mathrm{id}) \coloneqq \mu \circ_1 \phi - \mu \circ_2 \phi.$ 
\end{Definition}

\begin{lemma}
The data $(\Assc^\wedge, 0, \Theta)$ forms a complete curved absolute partial operad.
\end{lemma}

\begin{proof}
Analogous to the proof of the previous lemma.
\end{proof}

\begin{Proposition}
Let $u\mathcal{A}ss$ be the unital partial operad encoding unital associative algebras and let $u\mathcal{A}ss^{\ac}$ be its Koszul dual conilpotent curved partial cooperad. There is an isomorphism of complete curved partial operads
\[
\left(u\mathcal{A}ss^{\ac}\right)^* \cong \Assc^\wedge.
\]
\end{Proposition}

\begin{proof}
By direct inspection.
\end{proof}

\begin{Proposition}\label{absolute assliem}
There is a morphism of complete curved absolute partial operads $\Liec^\wedge \longrightarrow \Assc^\wedge$ which is determined by sending $\beta \mapsto \mu - \mu^{(12)}~,$ and $\zeta$ to $\phi$. 
\end{Proposition}

\begin{proof}
The proof is identical to the proof of Proposition \ref{assliem}.
\end{proof}

\begin{Remark}
As one can guess by the above examples, one can construct "absolute analogues" of well-known operads when they are given by generators and relations. Furthermore, they show that if $\mathcal{P}$ is a unital partial operad, then its Koszul dual \textit{operad} $\mathcal{P}^!$ is in fact an \textit{absolute operad}.
\end{Remark}

Another class of examples is given by object which already carry an underlying filtration that makes them "complete". Let $V$ be a complete pdg module as defined in Definition \ref{def: complete pdg module}, its curved endomorphism operad carries a natural structure of a complete curved absolute partial operad.

\begin{Definition}[Complete curved endomorphisms partial operad]
Let $V$ be a complete pdg module. Its \textit{complete curved endomorphisms partial operad} is given by the complete pdg $\mathbb{S}$-module 
\[
\mathrm{end}_V(n) \coloneqq \mathrm{hom}^{(\geq 1)}(V^{\otimes n}, V)~,
\]
where we only consider here morphisms of pdg modules which raise the filtration degree by at least one. The operad structure is given by the partial composition of functions. Its pre-differential is given by $\partial \coloneqq [d_V,-]$ and its curvature determined by $\Theta_V(\mathrm{id}) \coloneqq d_V^2~.$
\end{Definition}

\begin{lemma}\label{lemma: curved endomorphisms is absolute}
Let $V$ be a complete pdg module. Its complete curved endomorphisms partial operad is a complete curved \textit{absolute} partial operad. 
\end{lemma}

\begin{proof}
Forgetting the filtration, $\mathrm{end}_V$ is in particular a partial operad, thus an algebra over the reduced tree monad. Let 
\[
\gamma_V: \overline{\mathscr{T}}(\mathrm{end}_V) \longrightarrow \mathrm{end}_V
\]
be its structural morphism. The morphism $\gamma_V$ can be extended to the completed reduced tree monad using the completeness of the underlying filtration of $\mathrm{end}_V$. Notice that it is crucial to restrict ourselves to operations in the endomorphisms operad which raise the degree by at least one in order for this to be true. One can check that the extension of $\gamma_V$ to the complete reduced tree monad satisfies the axioms of an algebra over this monad. It is thus an absolute partial operad ; moreover it forms a curved absolute partial operad endowed with its curvature, see Lemma \ref{lemma: endo curved is curved}. 

\medskip

Let us show it is a complete as a curved absolute partial operad. Let $\mathrm{F}_\bullet \mathrm{end}_V$ be its underlying filtration and let $\overline{\mathscr{F}}_{\bullet} \mathrm{end}_V$ the canonical filtration induced by its curved absolute partial operad structure. See Definition \ref{def: canonical topology on absolutes} for more details. There is an obvious inclusion $\overline{\mathscr{F}}_{n}\mathrm{end}_V \subset \mathrm{F}_n \mathrm{end}_V$, since any operation which can be written as partial compositions iterated $n$ times also must raise the filtration degree by at least $n$. Thus 
\[
\bigcap_{n \in \mathbb{N}} \overline{\mathscr{F}}_{n} \mathrm{end}_V \subset \bigcap_{n \in \mathbb{N}} \mathrm{F}_n \mathrm{end}_V = 0~,
\]
which implies that $\mathrm{end}_V$ is a complete curved absolute partial operad.
\end{proof}

\begin{Remark}
This provides for a natural setting that allows the canonical filtration of a complete curved absolute partial operad to be reflected upon its algebras. See Section \ref{Section: curved HTT} for an application in this direction.
\end{Remark}

\chapter{The integration theory of curved absolute homotopy Lie algebras}\label{Chapter 3}

\epigraph{"Annie aime les algèbres, \\
 les algèbres de Lie. \\
 Elles ont la courbure en degré deux, \\
 la courbure des jours heureux."}{France Gall, peut-être.}

\section*{Introduction}
\textbf{Global picture.} In deformation theory, one studies how to "deform" structures (algebraic, geometric, etc) on an object. Here "deform" can have multiple meanings. Let $X$ be some type of object and $\mathcal{P}$ some type of structure. Ideally, one has a \textit{moduli space}  $\mathrm{Def}_X(\mathcal{P})$ where the points are the $\mathcal{P}$-structures that we can endow $X$ with. If one considers these structures up to some equivalences, and remembers the equivalences between the equivalences and so on and so forth, one has an $\infty$-groupoid (or \textit{moduli stack}). In this context, a first meaning of "deform" is, given a $\mathcal{P}$-structure on $X$, a point $x$ in $\mathrm{Def}_X(\mathcal{P})$, to look at the formal neighborhood of $x$ inside $\mathrm{Def}_X(\mathcal{P})$. This formal neighborhood encodes the \textit{infinitesimal deformations}, structures "infinitely close" to $x$.

\medskip

Over the years, it was noticed by D. Quillen, P. Deligne, V. Drinfel'd and many others that in all practical examples of these deformations, the formal neighborhoods were "encoded" by a differential graded (dg) Lie algebras. Examples include the work of K. Kodaira and D. Spencer in \cite{KodairaSpencer58}, encoding deformations of complex structures on manifolds or the work of M. Gerstenhaber in \cite{Gerstenhaber64}, encoding deformations of associative algebra structures on a vector space. This point of view suggests that for every point $x$ in $\mathrm{Def}_X(\mathcal{P})$, there should be a dg Lie algebra $\mathfrak{g}_{x,X}$ such that infinitesimal deformations correspond to elements in $\mathfrak{g}_{x,X}$ which satisfy the \textit{Maurer-Cartan equation}. This principle was stated by Drinfel'd in one of his letter to Schechtman:
\[
\text{"Every infinitesimal deformation problem in characteristic zero is encoded by a dg Lie algebra."}
\]
This became recently a theorem by J. Lurie and J. P. Pridham, see \cite{Lurie11,Pridham10}. The method was to formalize what an "infinitesimal deformation" is using the notion of $\kk$-\textit{pointed formal moduli problem} and then show that their homotopy category is equivalent to the homotopy category of dg Lie algebras. But given a dg Lie algebra, how can one recover the $\infty$-groupoid of infinitesimal deformation it encodes?

\medskip

Another point of view on the deformation theory of algebraic structures is given by operadic deformation complexes. In this context, given an operad $\mathcal{P}$ and a dg module $M$, one can construct an explicit dg Lie algebra where Maurer-Cartan elements are exactly the $\mathcal{P}$-algebra structures on $M$ and where the equivalences (called $\infty$-isotopies) between these structures are given by the action of a group called the \textit{gauge group}. Here "deforming" structures can also mean understanding when two $\mathcal{P}$-algebra structures are gauge equivalent to each other or not. Again, this amounts to studying the $\infty$-groupoid that is encoded by the operadic deformation problem. 

\medskip

\textit{Integration theory} amounts to constructing a way to recover these $\infty$-groupoids from their dg Lie algebras. Its most basic example can be traced back to Lie theory and the Baker--Campbell--Hausdorff formula which produces a group out of a nilpotent Lie algebra. This allows us to recover the gauge group in operadic deformation complexes. A first general approach to the integration of dg Lie algebras is given by the seminal work of V. Hinich in \cite{Hinich01}, using methods from D. Sullivan \cite{Sullivan77}. A refined version of the integration procedure was constructed by E. Getzler in \cite{Getzler09}, where he extended the integration procedure to nilpotent $\mathcal{L}_\infty$-algebras (homotopy Lie algebras). Inspired by the ideas of \cite{Buijandco}, generalized using operadic calculus, D. Robert-Nicoud and B. Vallette were able in \cite{robertnicoud2020higher} to give a new characterizations of Getzler's functor and obtained an adjunction 

\[
\begin{tikzcd}[column sep=7pc,row sep=3pc]
\mathsf{sSet} \arrow[r, shift left=1.1ex, "\mathcal{L}"{name=F}]      
&\mathcal{L}_\infty\text{-}\mathsf{alg}^{\mathsf{comp}}~. \arrow[l, shift left=.75ex, "\mathcal{R}"{name=U}]
\arrow[phantom, from=F, to=U, , "\dashv" rotate=-90]
\end{tikzcd}
\]
\vspace{0.1pc}

between simplicial sets and complete $\mathcal{L}_\infty$-algebras. This adjunction lies at the crossroad of three domains: Lie theory, deformation theory, and rational homotopy. Here the right adjoint functor $\mathcal{R}$ produces an \textit{integration functor} from complete $\mathcal{L}_\infty$-algebras to $\infty$-groupoids (Kan complexes). But the horn-fillers of this $\infty$-groupoid are a \textit{structure}, not a \textit{property}. They are given by explicit formulas which extend with the classical Baker--Campbell--Hausdorff formula in the specific case of Lie algebras. Finally, the functor $\mathcal{L}$ is shown to produce faithful rational models for pointed connected finite type nilpotent spaces, and greatly simplifies the original approach of \cite{Quillen69} to rational homotopy. Similar rational models using complete dg Lie algebras where also constructed in \cite{Buijandco}. 

\medskip

\textbf{Motivations.} The first goal of this chapter is to generalize the integration procedure to \textit{curved} $\mathcal{L}_\infty$-algebras. Curved Lie or curved $\mathcal{L}_\infty$-algebras are a notion more general than their dg counterparts: these objects are endowed with a distinguished element, the \textit{curvature}, which perturbs the relationships satisfied by the bracket and higher operations in a dg Lie or $\mathcal{L}_\infty$-algebra. In particular, the "differential" no longer squares to zero, thus the notion of quasi-isomorphism so useful in homotopical algebra disappears. 

\medskip

The main reason to consider curved $\mathcal{L}_\infty$-algebras is that they are the "non-pointed analogue" of $\mathcal{L}_\infty$-algebras. In a classical $\mathcal{L}_\infty$-algebra, the element $0$ is always a Maurer-Cartan element. This gives a canonical base point in the corresponding $\infty$-groupoid. It is no longer the case with curved $\mathcal{L}_\infty$-algebras. From the point of view of Lie theory, these objects behave like Lie groupoids instead of Lie groups. From the point of view of rational homotopy theory, these objects provide us with rational models for non necessarily pointed and non necessarily connected spaces.

\medskip

Developing the integration theory of curved $\mathcal{L}_\infty$-algebras has also many applications to deformation theory. In his PhD. thesis \cite{Nuiten19}, J. Nuiten showed that, if $A$ is a dg unital commutative algebra, then the homotopy category of $A$-\textit{pointed formal moduli problems} is equivalent to the homotopy category of dg Lie algebroids over $A$. Later, it was showed in \cite{calaque2021lie} that the homotopy category of dg Lie algebroids over $A$ is equivalent to the homotopy category of certain curved $\mathcal{L}_\infty$-algebras over the de Rham algebra of $A$, under some assumptions on $A$. Therefore deformation problems "parametrized by $\mathrm{Spec}(A)$" can be encoded with curved $\mathcal{L}_\infty$-algebras. "Parametrized by $\mathrm{Spec}(A)$" means that instead of choosing a $\kk$-point in a moduli space $X$, one chooses a morphism $f: \mathrm{Spec}(A) \longrightarrow X$, and then looks at the formal neighborhood of $\mathrm{Spec}(A)$ inside of $X$. Another approach to parametrized deformation problems is also given by $\mathcal{L}_\infty$-spaces introduced in \cite{Costello}. These are families of curved $\mathcal{L}_\infty$-algebras parametrized by smooth manifolds. They were constructed in order to treat parametrized deformation problems arising in fundamental physics. On the other side, operadic deformation complexes of \textit{unital} algebraic structures form \textit{curved} Lie algebras as shown in Chapter \ref{Chapter 2}. Furthermore, the space of $\infty$-morphisms between types of unital algebras or the space of $\infty$-morphisms between types of counital coalgebras are encoded, in both cases, by convolution curved $\mathcal{L}_\infty$-algebras. In all the above cases, having a good \textit{integration theory} for curved $\mathcal{L}_\infty$-algebras is of primordial importance. 

\medskip
 
\textbf{Framework.} The methods of \cite{robertnicoud2020higher} cannot be easily generalized to the curved setting. First, we needed to develop the theory of curved operadic calculus in Chapter \ref{Chapter 2} in order to be able to correctly deal with curved algebras. Nevertheless, two main difficulties persist. The first one is that, in the curved setting, the notion of quasi-isomorphism is gone, thus one needs to find another notion of weak-equivalences. The second one is intrinsic to Koszul duality, as the dual notion of curved algebras are counital coalgebras. The operadic methods in \textit{loc.cit} rely on transferring a homotopy cocommutative coalgebra structure on the cellular chain functor. But in the curved case, this structure is a homotopy counital cocommutative coalgebra structure. Hence it is \textit{not conilpotent}, and therefore can not be encoded as coalgebras over a cooperad.

\medskip

The solution to these problems is to introduce new methods, principally the general framework of \cite{grignoulejay18} and the developments of it that we introduced in Chapter \ref{Chapter 1} and \ref{Chapter 2}. The only way to encode non-conilpotent coalgebras is with operads. Therefore our idea is to encode curved $\mathcal{L}_\infty$-algebras with \textit{a curved cooperad}. This gives rise to the new notion of \textit{curved absolute} $\mathcal{L}_\infty$\textit{-algebras}. A curved absolute $\mathcal{L}_\infty$-algebra can be thought as curved $\mathcal{L}_\infty$-algebra where all infinite sums of operations have a well defined image \textit{without supposing any underlying topology}. Infinite sums appear naturally already in the theory of $\mathcal{L}_\infty$-algebras. There has been two standard ways to deal with them either restrict to nilpotent $\mathcal{L}_\infty$-algebras, or use the somewhat \textit{ad hoc} approach of changing the base category in order to consider objects with an underlying complete topology, so that these infinite sums \textit{converge}. Our approach solves these problems altogether. Another good feature of curved absolute $\mathcal{L}_\infty$-algebras is that they admit a \textit{canonical filtration} from their algebraic structure which is the analogue of the coradical filtration for coalgebras. This enables us to define \textit{complete} curved absolute $\mathcal{L}_\infty$-algebras without changing the base category.

\medskip

We introduce $u\mathcal{CC}_\infty$-coalgebras, a version of homotopy counital cocommutative coalgebras encoded by a dg operad. Using the methods of \cite{grignoulejay18}, we endow the category of $u\mathcal{CC}_\infty$-coalgebras with a model structure where weak-equivalences are given by quasi-isomorphisms and cofibrations by monomorphisms. Then we transfer it onto the category of curved absolute $\mathcal{L}_\infty$-algebras, and thus obtain a Quillen equivalence. This solves our starting difficulties, and furthermore allows us to implement a new model category framework in the theory of integration. 

\medskip

\textbf{Main results.} The main results of this chapter are based on \cite{lucio2022integration}. We construct a functorial $u\mathcal{CC}_\infty$-coalgebra structure on the cellular chain functor $C_*^c(-)$. It automatically admits a right adjoint functor $\overline{R}$.

\begin{theoremintro}[Theorem \ref{thm: triangle of adjunctions}]
The following triangle of Quillen adjunctions commutes

\[
\begin{tikzcd}[column sep=5pc,row sep=2.5pc]
&\hspace{1pc}u\mathcal{CC}_\infty \textsf{-}\mathsf{coalg} \arrow[dd, shift left=1.1ex, "\widehat{\Omega}_{\iota}"{name=F}] \arrow[ld, shift left=.75ex, "\overline{\mathcal{R}}"{name=C}]\\
\mathsf{sSet}  \arrow[ru, shift left=1.5ex, "C^c_*(-)"{name=A}]  \arrow[rd, shift left=1ex, "\mathcal{L}"{name=B}] \arrow[phantom, from=A, to=C, , "\dashv" rotate=-70]
& \\
&\hspace{3pc}\mathsf{curv}~\mathsf{abs}~\mathcal{L}_\infty\textsf{-}\mathsf{alg}~. \arrow[uu, shift left=.75ex, "\widehat{\mathrm{B}}_{\iota}"{name=U}] \arrow[lu, shift left=.75ex, "\mathcal{R}"{name=D}] \arrow[phantom, from=B, to=D, , "\dashv" rotate=-110] \arrow[phantom, from=F, to=U, , "\dashv" rotate=-180]
\end{tikzcd}
\]
\end{theoremintro}

Here $\widehat{\Omega}_{\iota}$ is a completed Cobar construction, and $\widehat{\mathrm{B}}_{\iota}$ is a new complete Bar construction. This complete Bar construction is constructed using the \textit{real cofree} $u\mathcal{CC}_\infty$-coalgebra, not the usual cofree conilpotent coalgebra. Its homology can be thought as a \textit{higher Chevalley-Eilenberg} homology for curved absolute $\mathcal{L}_\infty$-algebras.

\medskip

The functor $\mathcal{R}$ is the \textit{integration functor} we where looking for. One of the main novelties of this is the fact that $\mathcal{R}$ is a right Quillen functor. This new methods give an intrinsic model category framework for integration theory, and embed classical Lie theory into the realm of homotopical algebra. This approach also allows us to obtain conceptually all the required properties that a well-behaved integration functor needs to satisfy. 

\begin{theoremintro}[Theorem \ref{thm: propriétés de l'intégration}]
\leavevmode

\medskip

\begin{enumerate}
\item For any curved absolute $\mathcal{L}_\infty$-algebra $\mathfrak{g}$, the simplicial set $\mathcal{R}(\mathfrak{g})$ is a Kan complex.

\medskip

\item For any degree-wise epimorphism $f: \mathfrak{g} \twoheadrightarrow \mathfrak{h}$ of curved absolute $\mathcal{L}_\infty$-algebras, the induced map
\[
\mathcal{R}(f): \mathcal{R}(\mathfrak{g}) \twoheadrightarrow \mathcal{R}(\mathfrak{h})
\]
is a fibrations of simplicial sets. 

\medskip

\item The functor $\mathcal{R}(-)$ preserves weak equivalences. In particular, it sends any graded quasi-isomorphism between complete curved absolute $\mathcal{L}_\infty$-algebras to a weak homotopy equivalence of simplicial sets.

\end{enumerate}
\end{theoremintro}

The fact that $\mathcal{R}(\mathfrak{g})$ is an $\infty$-groupoid is quintessential to integration theory as developed in \cite{Hinich01} and in \cite{Getzler09}. The second statement is a generalization of one of the main theorems of \cite{Getzler09}. The third point, the \textit{homotopy invariance} of the integration functor, is the generalization of the celebrated Goldman--Milson Theorem of  \cite{goldmanmillson} and its extension by Dolgushev--Rogers to $\mathcal{L}_\infty$-algebras in \cite{dolgushevrogers}. Moreover, the adjunction $\mathcal{L} \dashv \mathcal{R}$ is shown to be a non-abelian generalization of the Dold-Kan correspondence. 

\medskip

Having these constructions at hand, the rest of the second section is devoted to the study of higher absolute Lie theory. We introduce and characterize gauge equivalences for Maurer-Cartan elements in this setting. We construct higher Baker--Campbell--Hausdorff products for curved absolute $\mathcal{L}_\infty$-algebras and we show that these are given by the same explicit formulae as in \cite{robertnicoud2020higher}. We generalize Berglund's Theorem of \cite{Berglund15} to the case of curved absolute $\mathcal{L}_\infty$-algebras. If $\mathfrak{g}$ is curved absolute $\mathcal{L}_\infty$-algebra, it establishes an isomorphism between the homotopy groups of $\mathcal{R}(\mathfrak{g})$ and the homology groups of $\mathfrak{g}$ with a twisted differential. Finally, we establish comparison results with Getzler and Hinich integration functor using the restriction of a curved absolute $\mathcal{L}_\infty$-algebra to a curved $\mathcal{L}_\infty$-algebra. 

\medskip

The third section is devoted to constructing rational homotopy models using the new functor $\mathcal{L}$. The first main result is the following one. 

\begin{theoremintro}[Theorem \ref{thm: modèles d'homotopie rationnel type fini}]
Let $X$ be a finite type nilpotent simplicial set. The unit of the adjunction

\[
\eta_X: X \qi \mathcal{R}\mathcal{L}(X) 
\]

is a rational homotopy equivalence.
\end{theoremintro}

This allows us to obtain rational models for non necessarily pointed nor connected finite type nilpotent spaces using our functor $\mathcal{L}$. We then construct explicit models for each of the connected components of $X$ using absolute $\mathcal{L}_\infty$-algebras. This allows us to give a characterization of the essential homotopical image of the cellular chain functor $C_*^c(-)$ in the category of $u\mathcal{CC}_\infty$-coalgebras. Applying the theory of \cite{grignou2022mapping}, we get a convolution curved absolute $\mathcal{L}_\infty$-algebra on the space of graded morphisms between a $u\mathcal{CC}_\infty$-coalgebra and a curved absolute $\mathcal{L}_\infty$-algebra. Using this construction, we prove the following theorem. 

\begin{theoremintro}[Theorem \ref{thm: vrai thm mapping spaces}]
Let $\mathfrak{g}$ be a curved absolute $\mathcal{L}_\infty$-algebra and let $X$ be a simplicial set. There is a weak equivalence of Kan complexes

\[
\mathrm{Map}(X, \mathcal{R}(\mathfrak{g})) \qi \mathcal{R}\left(\mathrm{hom}(C^c_*(X),\mathfrak{g})\right)~,
\]
\vspace{0.25pc}

where $\mathrm{hom}(C^c_*(X),\mathfrak{g})$ denotes the convolution curved absolute $\mathcal{L}_\infty$-algebra. 
\end{theoremintro}

Notice that, to the best of our knowledge, there are no assumptions on $X$ nor on $\mathfrak{g}$ for the first time. If $Y$ is a finite type nilpotent simplicial set, the above theorem gives an explicit model for the mapping space of $X$ and the $\mathbb{Q}$-localization of $Y$. Futhermore, this model is constructed using the cellular chains on $X$, hence it is also relatively small in comparison to other rational models constructed before. See for instance \cite{Berglund15, BuijMapping, LazarevMapping}.

\medskip

The last section is devoted to applications in deformation theory. We first construct a convolution curved algebra $\mathcal{L}_\infty$-algebras associated to any curved twisting morphism $\alpha: \mathcal{C} \longrightarrow \mathcal{P}$ between a conilpotent curved cooperad and a dg operad. This allows to encode $\infty$-morphisms between unital algebras encoded by an dg operad as the Maurer-Cartan elements of a curved algebra $\mathcal{L}_\infty$-algebra. The same methods could be applied to $\infty$-morphisms of coalgebras over operads or, more generally, $\infty$-morphisms of gebras over properads as defined in \cite{hoffbeck2019properadic}.

\medskip

The last results concern the relationship between curved absolute $\mathcal{L}_\infty$-algebras and formal moduli problems. The first observation that one can make is that, since we are working with \textit{curved} absolute $\mathcal{L}_\infty$-algebras, the dg Artinian algebras that encode the formal neighborhoods of points need not be augmented. Starting from a derived affine stack $A$, that is, a dg unital commutative algebra, we construct an explicit curved absolute $\mathcal{L}_\infty$-algebra model $\mathfrak{g}_A$ that encodes the formal geometry of $A$ around points which are not specified. 

\begin{theoremintro}[Theorem \ref{thm: modèle géométrique sur les points Artiniens}]
Let $A$ be a dg $u\mathcal{C}om$-algebra. Let $B$ be a dg Artinian algebra. There is a weak equivalence of simplicial sets
\[
\mathrm{Spec}(A)(B) \simeq \mathcal{R}(\mathrm{hom}( (\mathrm{Res}_\varepsilon B)^\circ, \mathfrak{g}_A))~.
\]
\end{theoremintro}

More concretely, we can recover the formal neighborhood of \textit{any finite family} of points of $A$ living in any \textit{finite field extension} $\mathbb{L}$ of $\kk$ from its model $\mathfrak{g}_A$. Conversely, if we start with a curved absolute $\mathcal{L}_\infty$-algebra $\mathfrak{g}$, we can construct a deformation functor 
\[
\mathrm{Def}_{\mathfrak{g}}: \mathsf{dg}~\mathsf{Art}\text{-}\mathsf{alg}_{\geq 0} \longrightarrow \mathsf{sSet}
\]
from dg Artinian algebras to $\infty$-groupoids which preserves weak-equivalences and homotopy pullbacks. Furthermore, if two curved absolute $\mathcal{L}_\infty$-algebras are weakly equivalent, then their respective deformation functors are shown to be naturally weakly-equivalent. In the lack of a good definition of what a "non-pointed" formal moduli problem should correspond to, we limit ourselves to these constructions which make apparent the relationship between the new curved absolute $\mathcal{L}_\infty$-algebra we have introduced and formal geometry around points which are \textit{not specified in advance}. This will be the subject of further studies.

\section{Curved absolute $\mathcal{L}_\infty$-algebras}
In this section, we introduce the notion of curved \textit{absolute} $\mathcal{L}_\infty$-algebras. This is a new type of curved $\mathcal{L}_\infty$-algebras, encoded by a conilpotent curved cooperad. They posses a much richer algebraic structure than usual curved $\mathcal{L}_\infty$-algebras. Any infinite sum of operations has a well-defined image by definition. In particular, the Maurer-Cartan equation is always defined. The rest of this chapter is devoted to the study of this new notion and its applications. Notice that our degree conventions will correspond to \textit{shifted} curved $\mathcal{L}_\infty$-algebras. This \textit{shifted convention} will be implicit from now on. Curved absolute algebras appear as the Koszul dual of a specific model for \textit{non-necessarily conilpotent} $\mathcal{E}_{\infty}$-coalgebras, which we call $u\mathcal{CC}_{\infty}$-coalgebras. We construct a complete Bar-Cobar adjunction that relates these two types of algebraic objects and show that it is a Quillen equivalence.

\subsection{Curved absolute $\mathcal{L}_\infty$-algebras.}
Let $\ucom$ be the operad encoding unital commutative algebras and let $\epsilon: \ucom \longrightarrow \I$ be the canonical morphism of $\mathbb{S}$-modules given by the identity on $\ucom(1) \cong \kk~.$ We denote $\mathrm{B}^{\mathrm{s.a}}\ucom$ its semi-augmented Bar construction with respect to this canonical semi-augmentation. We refer to Appendix \ref{Section: Appendix B} for more details on this particular construction.

\begin{Definition}[Curved absolute $\mathcal{L}_\infty$-algebra]
A \textit{curved absolute} $\mathcal{L}_\infty$\textit{-algebra} $\mathfrak{g}$ amounts to the data $(\mathfrak{g},\gamma_\mathfrak{g},d_\mathfrak{g})$ of a curved $\mathrm{B}^{\mathrm{s.a}}\ucom$-algebra. 
\end{Definition}

\begin{Remark}
General recollections on the notion of curved algebras over cooperads are given in Section \ref{Section: Curved cooperads}.
\end{Remark}

Let us unravel this definition. The data of a curved absolute $\mathcal{L}_\infty$-algebra structure on a pdg module $(\mathfrak{g},d_\mathfrak{g})$ amounts to the data of a morphism of pdg modules
\[
\gamma_\mathfrak{g}: \prod_{n \geq 0} \mathrm{Hom}_{\mathbb{S}_n}\left(\mathrm{B}^{\mathrm{s.a}}\ucom(n), \mathfrak{g}^{\otimes n}\right) \longrightarrow \mathfrak{g}~,
\]
which satisfies the conditions of Definition \ref{def curved alg over a coop}. This map admits a simpler description.

\begin{lemma}\label{lemma: iso invariants avec les coinvariants}
There is an isomorphism of pdg modules
\[
\prod_{n \geq 0} \mathrm{Hom}_{\mathbb{S}_n}\left(\mathrm{B}^{\mathrm{s.a}}\ucom(n), \mathfrak{g}^{\otimes n}\right) \cong \prod_{n \geq 0} \widehat{\Omega}^{\mathrm{s.a}}\ucom^*(n) ~\widehat{\otimes}_{\mathbb{S}_n}~ \mathfrak{g}^{\otimes n}~,
\]
natural in $\mathfrak{g}$, where $\widehat{\otimes}$ denotes the completed tensor product with respect to the canonical filtration on the complete Cobar construction.
\end{lemma}

\begin{proof}
By definition, $\mathrm{B}^{\mathrm{s.a}}\ucom$ is a conilpotent curved cooperad. Let us denote by $\mathscr{R}_\omega \mathrm{B}^{\mathrm{s.a}}$ the $\omega$-term of its coradical filtration. There is an isomorphism of conilpotent curved cooperads
\[
\mathrm{B}^{\mathrm{s.a}}\ucom \cong \colim_{\omega}\mathscr{R}_\omega \mathrm{B}^{\mathrm{s.a}}\ucom~.
\]
Notice that for every $n \geq 0$ and every $\omega \geq 0$, $\mathscr{R}_\omega \mathrm{B}^{\mathrm{s.a}}\ucom(n)$ is degree-wise finite dimensional. Thus we have
\begin{align*}
\mathrm{Hom}_{\mathbb{S}_n}\left(\mathrm{B}^{\mathrm{s.a}}\ucom(n), \mathfrak{g}^{\otimes n}\right)
&\cong  \mathrm{Hom}_{\mathbb{S}_n}\left(\colim_{\omega}\mathscr{R}_\omega \mathrm{B}^{\mathrm{s.a}}\ucom(n), \mathfrak{g}^{\otimes n}\right) \\
&\cong \lim_{\omega} \mathrm{Hom}_{\mathbb{S}_n}\left(\mathscr{R}_\omega \mathrm{B}^{\mathrm{s.a}}\ucom(n), \mathfrak{g}^{\otimes n}\right)~.
\end{align*}
for all $n \geq 0$. By Lemma \ref{lemma: Bucom dual lin}, there is an isomorphism 
\[
\mathscr{R}_\omega \mathrm{B}^{\mathrm{s.a}}\ucom \cong \widehat{\Omega}^{\mathrm{s.c}}\ucom^*/\mathscr{F}_\omega \widehat{\Omega}^{\mathrm{s.c}}\ucom^*~,
\]
where $\mathscr{F}_\omega \widehat{\Omega}^{\mathrm{s.c}}\ucom^*$ denotes the $\omega$-term of its canonical filtration as an absolute partial operad. Therefore there are isomorphisms
\begin{align*}
\lim_{\omega} \mathrm{Hom}_{\mathbb{S}_n} \left(\mathscr{R}_\omega \mathrm{B}^{\mathrm{s.a}}\ucom(n), \mathfrak{g}^{\otimes n}\right)
&\cong \lim_{\omega} \left(\widehat{\Omega}^{\mathrm{s.c}}\ucom^*/\mathscr{F}_\omega \widehat{\Omega}^{\mathrm{s.c}}\ucom^*(n) \otimes \mathfrak{g}^{\otimes n}\right)^{\mathbb{S}_n} \\
&\cong \left(\widehat{\Omega}^{\mathrm{s.c}}\ucom^*(n) ~\widehat{\otimes}~ \mathfrak{g}^{\otimes n}\right)^{\mathbb{S}_n} \\
&\cong ~\widehat{\Omega}^{\mathrm{s.c}}\ucom^*(n) ~\widehat{\otimes}_{\mathbb{S}_n}~ \mathfrak{g}^{\otimes n}~.
\end{align*}
Notice that the last isomorphism identifies invariants with coinvariants, which is possible because of the characteristic zero assumption. Nevertheless, this identification carries non-trivial coefficients, see Remark \ref{Remark: renormalization}.
\end{proof}

\begin{Notation}
Let $\mathrm{CRT}_n^\omega$ denote the set of \textit{corked rooted trees} of arity $n$ and with $\omega$ internal edges. A corked rooted tree is a rooted tree where vertices either have at least two incoming edges or zero incoming edges, which are called \textit{corks}. The arity of a corked rooted tree is the number of non-corked leaves. The unique rooted tree of arity $n$ with one vertex is called the $n$-corolla, denoted by $c_n$. Notice that for each $n$, the set $\mathrm{CRT}_n$ of corked rooted trees of arity $n$ is infinite. Likewise, for each $\omega$, the set of corked rooted trees with $\omega$ internal edges is also infinite. Nevertheless, the set $\mathrm{CRT}_n^\omega$ is finite. The only corked rooted tree of weight $0$ is the trivial tree of arity one with zero vertices.
\end{Notation}

\begin{Proposition}
Let $(\mathfrak{g},d_\mathfrak{g})$ be a pdg module with a basis $\left\{ g_b ~|~b \in B \right\}~.$ The pdg module
\[
\prod_{n \geq 0} \widehat{\Omega}^{\mathrm{s.a}}\ucom^*(n) ~\widehat{\otimes}_{\mathbb{S}_n} ~ \mathfrak{g}^{\otimes n}
\]
admits a basis given by double series on $n$ and $\omega$ of corked rooted trees in $\mathrm{CRT}_n^\omega$ labeled by the basis elements of $\mathfrak{g}^{\otimes n}$. These basis elements can be written as

\[
\left\{ \sum_{n\geq 0} \sum_{\omega \geq 0} \sum_{\tau \in \mathrm{CRT}_n^\omega} \lambda_\tau \tau \left(g_{i_1}, \cdots, g_{i_n}\right) \right\}~,
\]

where $\lambda_\tau$ is a scalar in $\kk$ and $(i_1,\cdots,i_n)$ is in $B^n$. Here $\tau(g_{i_1}, \cdots, g_{i_n})$ is given by the rooted tree $\tau$ with input leaves decorated by the elements $\tau(g_{i_1}, \cdots, g_{i_n})$. The degree of $\tau(g_{i_1}, \cdots, g_{i_n})$ is $-\omega -1 + |g_{i_1}| + \cdots + |g_{i_n}|~.$ An example of such decoration is given by

\medskip

\begin{center}
\includegraphics[width=130mm,scale=1.3]{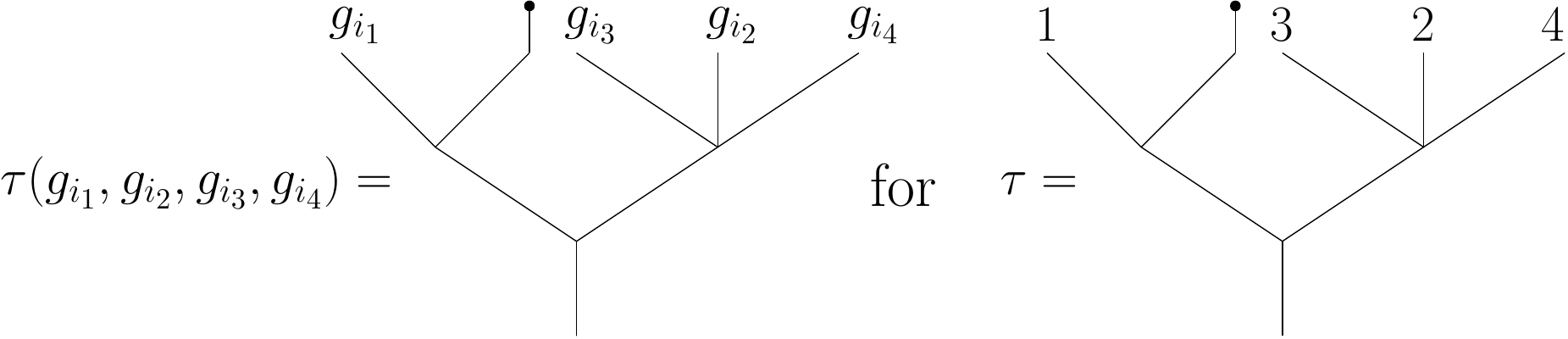}.
\end{center}

The pre differential of the pdg module is given by the sum of two terms: the first term is applies the pre-differential $d_\mathfrak{g}$ on each of the labels $g_{i_j}$ of a corked rooted tree $\tau$, the second splits each vertex into two vertices in all the possible ways including the splitting with a cork. 
\end{Proposition}

\begin{proof}
Corked rooted trees form a basis of the pdg $\mathbb{S}$-module $\widehat{\Omega}^{\mathrm{s.a}}\ucom^*$, hence the results follows by direct inspection.
\end{proof}

\begin{Remark}[Renormalization]\label{Remark: renormalization}
Let $\tau$ be a corked rooted tree. It can be written as $c_m \circ (\tau_1, \cdots, \tau_m)$. We define recursively the following coefficient for corked rooted trees:

\[
\mathcal{E}(c_n) \coloneqq n! \quad \text{and} \quad \mathcal{E}\left(c_m \circ (\tau_1, \cdots, \tau_m)\right) \coloneqq m! \prod_{i=1}^m\mathcal{E}(\tau_i )~.
\]

Then doubles series 

\[
\left\{ \sum_{n\geq 0} \sum_{\omega \geq 0} \sum_{\tau \in \mathrm{CRT}_n^\omega} \mathcal{E}(\tau)\lambda_\tau \tau \left(g_{i_1}, \cdots, g_{i_n}\right) \right\}~,
\]

form a basis of 

\[
\prod_{n \geq 0} \mathrm{Hom}_{\mathbb{S}_n}\left(\mathrm{B}^{\mathrm{s.a}}\ucom(n), \mathfrak{g}^{\otimes n}\right) \cong  \prod_{n \geq 0}   \left(\widehat{\Omega}^{\mathrm{s.a}}\ucom^*(n) ~\widehat{\otimes} ~ \mathfrak{g}^{\otimes n} \right)^{\mathbb{S}_n}~.
\]
\end{Remark}

\medskip

\textbf{Structural data.} Thus the data of a curved absolute $\mathcal{L}_\infty$-algebra structure on a pdg module $(\mathfrak{g},d_\mathfrak{g})$ amounts to the data of a morphism of pdg modules
\[
\gamma_\mathfrak{g}: \prod_{n \geq 0} \widehat{\Omega}^{\mathrm{s.a}}\ucom^*(n) ~\widehat{\otimes}_{\mathbb{S}_n} ~ \mathfrak{g}^{\otimes n} \longrightarrow \mathfrak{g}~,
\]
which satisfies the conditions of Definitions \ref{def: C algebra} and \ref{def curved alg over a coop}. In particular, for any infinite sum 
\[
\sum_{n\geq 0} \sum_{\omega \geq 0} \sum_{\tau \in \mathrm{CRT}_n^\omega} \lambda_\tau \tau(g_{i_1}, \cdots, g_{i_n})~,
\]
there is a well-defined image 
\[
\gamma_\mathfrak{g}\left(\sum_{n\geq 0} \sum_{\omega \geq 0} \sum_{\tau \in \mathrm{CRT}_n^\omega} \lambda_\tau \tau(g_{i_1}, \cdots, g_{i_n})\right)~,
\]
in the pdg module $\mathfrak{g}$. This \textbf{does not presuppose} an underlying topology on the pdg module $\mathfrak{g}$.

\medskip

\textbf{Pdg condition.} Let us make explicit each of the conditions satisfied by the structural morphism $\gamma_\mathfrak{g}$. It is a morphism of pdg modules if only if, we have that

\begin{small}

\begin{align}\label{pdg condition}
\gamma_\mathfrak{g}\left(\sum_{\substack{n \geq 0 \\ \omega \geq 0}} \sum_{\tau \in \mathrm{CRT}_n^\omega} \lambda_\tau d_{\widehat{\Omega}\ucom}(\tau)(g_{i_1}, \cdots, g_{i_n})\right) + \gamma_\mathfrak{g}\left(\sum_{\substack{n \geq 0 \\ \omega \geq 0}} \sum_{\tau \in \mathrm{CRT}_n^\omega}\sum_{j=0}^n (-1)^j \lambda_\tau\tau(g_{i_1}, \cdots, d_\mathfrak{g}(g_{i_j}), \cdots, g_{i_n})\right) = \nonumber \\
d_\mathfrak{g}\left(\gamma_\mathfrak{g}\left(\sum_{\substack{n \geq 0 \\ \omega \geq 0}} \sum_{\tau \in \mathrm{CRT}_n^\omega} \lambda_\tau\tau(g_{i_1}, \cdots, g_{i_n})\right) \right)~.
\end{align}

\end{small}

\textbf{Associativity condition.} The associativity condition on $\gamma_\mathfrak{g}$ imposed by the diagram of \ref{def: C algebra} is equivalent to the following equality

\begin{small}

\begin{align}\label{associativity condition}
\gamma_\mathfrak{g} \left(\sum_{\substack{k\geq 0 \\ \omega \geq 0}} \sum_{\tau \in \mathrm{CRT}_k^\omega} \lambda_\tau \tau \left(\gamma_\mathfrak{g}\left(\sum_{\substack{i_1 \geq 0 \\ \omega_1 \geq 0}} \sum_{\tau_1 \in \mathrm{CRT}_{i_1}^{\omega_1}} \lambda_{\tau_1}\tau_1(\bar{g}_{i_1} )\right) ,\cdots, \gamma_\mathfrak{g} \left(\sum_{\substack{i_k \geq 0 \\ \omega_k \geq 0}} \sum_{\tau_k \in \mathrm{CRT}_{i_k}^{\omega_k}} \lambda_{\tau_k} \tau_k(\bar{g}_{i_k} )\right)\right)\right) = \nonumber \\
\gamma_\mathfrak{g} \left(\sum_{n \geq 0} \sum_{k \geq 0} \sum_{i_1 + \cdots + i_k =n} \sum_{\substack{\omega^{\mathrm{tot}} \geq 0  \\ \omega +\omega_1 + \cdots +\omega_k = \omega^{\mathrm{tot}}}} \sum_{\substack{\tau \in \mathrm{CRT}_k^\omega \\ \tau_j \in \mathrm{CRT}_{i_j}^{\omega_j}}} \lambda_\tau(\lambda_{\tau_1}\cdots\lambda_{\tau_k})~ \tau \circ (\tau_1,\cdots, \tau_k)(\bar{g}_{i_1},\cdots,\bar{g}_{i_k}) \right)~, 
\end{align}
\end{small}

where $\bar{g}_{i_j}$ denotes an $i_j$-tuple of elements of $\mathfrak{g}^{\otimes i_j}$ and where $\tau \circ (\tau_1,\cdots, \tau_k)$ denotes the corked rooted tree obtained by grafting $(\tau_1,\cdots, \tau_k)$ onto the leaves of $\tau~.$ If $\gamma_\mathfrak{g}$ satisfies the above conditions, it endows $\mathfrak{g}$ with a pdg $\text{B}^{\mathrm{s.a}}\ucom$-algebra structure.

\begin{Definition}[Elementary operations of a curved absolute $\mathcal{L}_\infty$-algebra]
Let $(\mathfrak{g},\gamma_\mathfrak{g}, d_\mathfrak{g})$ be a curved absolute $\mathcal{L}_\infty$-algebra. The \textit{elementary operations} of $\mathfrak{g}$ are the family of symmetric operations of degree $-1$ given by
\[
\left\{l_n \coloneqq \gamma_\mathfrak{g}(c_n(-,\cdots,-)): \mathfrak{g}^{\odot n} \longrightarrow \mathfrak{g}~\right\}
\]
for all $n \neq 1$. 
\end{Definition}

\textbf{Curved condition.} The condition on $\gamma_\mathfrak{g}$ imposed by the diagram of Definition \ref{def curved alg over a coop} amounts in this case to the following equation on the elementary operations

\begin{equation}\label{curved condition}
d_\mathfrak{g}^2(g) = l_2(l_0,g)~.
\end{equation}

\textbf{Conclusion.} A curved absolute $\mathcal{L}_\infty$-algebra structure $\gamma_\mathfrak{g}$ on a pdg module $(\mathfrak{g},d_\mathfrak{g})$ amounts to the data of a degree $0$ map 
\[
\gamma_\mathfrak{g}: \prod_{n \geq 0} \widehat{\Omega}^{\mathrm{s.a}}\ucom^*(n) ~\widehat{\otimes}_{\mathbb{S}_n} ~ \mathfrak{g}^{\otimes n} \longrightarrow \mathfrak{g}~,
\]
satisfying conditions \ref{pdg condition}, \ref{associativity condition} and \ref{curved condition}.

\begin{Remark}[Warning]
In general, if $(\mathfrak{g},\gamma_\mathfrak{g}, d_\mathfrak{g})$ is a curved absolute $\mathcal{L}_\infty$-algebra, then 
\[
\gamma_\mathfrak{g}\left(\sum_{n\geq 0} \sum_{\omega \geq 0} \sum_{\tau \in \mathrm{CRT}_n^\omega} \lambda_\tau \tau(g_{i_1}, \cdots, g_{i_n})\right) \neq \sum_{n\geq 0} \sum_{\omega \geq 1} \sum_{\tau \in \mathrm{CRT}_n^\omega} \lambda_\tau \gamma_\mathfrak{g} \left(\tau(g_{i_1}, \cdots, g_{i_n})\right)~,
\]
as the latter expression \textit{is not well-defined} since infinite sums of elements \textit{in} $\mathfrak{g}$ are not well-defined in general.
\end{Remark}

\textbf{Morphisms.} Let $(\mathfrak{g},\gamma_\mathfrak{g}, d_\mathfrak{g})$ and $(\mathfrak{h},\gamma_\mathfrak{h}, d_\mathfrak{h})$ be two curved absolute $\mathcal{L}_\infty$-algebras and let $f: \mathfrak{g} \longrightarrow \mathfrak{h}$ be a morphism of pdg modules. The condition for $f$ to be a morphism of curved absolute $\mathcal{L}_\infty$-algebras can be written as 
\[
f\left(\gamma_\mathfrak{g}\left(\sum_{n\geq 0} \sum_{\omega \geq 0} \sum_{\tau \in \mathrm{CRT}_n^\omega} \lambda_\tau \tau(g_{i_1}, \cdots, g_{i_n})\right)\right) = \gamma_\mathfrak{h}\left(\sum_{n\geq 0} \sum_{\omega \geq 0} \sum_{\tau \in \mathrm{CRT}_n^\omega} \lambda_\tau \tau(f(g_{i_1}), \cdots, f(g_{i_n}))\right)~.
\]

Any curved absolute $\mathcal{L}_\infty$-algebra structure can be restricted along elementary operations to induce a curved $\mathcal{L}_\infty$-algebra structure in the standard sense.

\begin{Proposition}\label{prop: restriction functor}
There is a restriction functor 
\[
\mathrm{Res}: \mathsf{curv}~\mathsf{abs}~\mathcal{L}_\infty\text{-}\mathsf{alg} \longrightarrow \mathsf{curv}~\mathcal{L}_\infty\text{-}\mathsf{alg}~,
\]
from the category of curved absolute $\mathcal{L}_\infty$-algebras to the category of curved $\mathcal{L}_\infty$-algebras which is faithful.
\end{Proposition}

\begin{proof}
Let $(\mathfrak{g},\gamma_\mathfrak{g}, d_\mathfrak{g})$ be a curved absolute $\mathcal{L}_\infty$-algebra, we can restrict the structural map

\[
\begin{tikzcd}
\mathrm{Res}(\gamma_\mathfrak{g}): \displaystyle\bigoplus_{n \geq 0} \widehat{\Omega}^{\mathrm{s.c}}\ucom^*(n) \otimes_{\mathbb{S}_n} \mathfrak{g}^{\otimes n} \arrow[r,"\iota_\mathfrak{g}"]
&\displaystyle \prod_{n \geq 0} \widehat{\Omega}^{\mathrm{s.c}}\ucom^*(n) ~\widehat{\otimes}_{\mathbb{S}_n}~ \mathfrak{g}^{\otimes n} \arrow[r,"\gamma_\mathfrak{g}"]
&\mathfrak{g}
\end{tikzcd}
\]

along the natural inclusion $\iota_\mathfrak{g}$. It endows $\mathfrak{g}$ with a curved $\widehat{\Omega}^{\mathrm{s.c}}\ucom^*$-algebra structure. By Proposition \ref{Prop: Curved algebras over OmegauCom}, this is equivalent to a curved $\mathcal{L}_\infty$-algebra structure in the classical sense of Definition \ref{def: classical curved linfty alg}. Any morphism of curved absolute $\mathcal{L}_\infty$-algebras is in particular a morphism of curved $\mathcal{L}_\infty$-algebras, and it faithful.
\end{proof}

\begin{Proposition}
There restriction functor admits a left adjoint $\mathrm{Abs}$, which called the absolute envelope of a curved $\mathcal{L}_\infty$-algebra. Therefore there is an adjunction 

\[
\begin{tikzcd}[column sep=7pc,row sep=3pc]
\mathsf{curv}~\mathcal{L}_\infty\text{-}\mathsf{alg} \arrow[r, shift left=1.1ex, "\mathrm{Abs}"{name=F}]      
&\mathsf{curv}~\mathsf{abs}~\mathcal{L}_\infty\text{-}\mathsf{alg}~, \arrow[l, shift left=.75ex, "\mathrm{Res}"{name=U}]
\arrow[phantom, from=F, to=U, , "\dashv" rotate=-90]
\end{tikzcd}
\]
\end{Proposition}

\begin{proof}
Both categories are presentable, by \cite[Theorem 1.66]{AdamekRosicky} the functor $\mathrm{Res}$ admits as left-adjoint since it preserves all limits.
\end{proof}

\begin{Remark}
Let $(\mathfrak{g},\{l_n\},d_\mathfrak{g})$ be a curved $\mathcal{L}_\infty$-algebra. Its absolute envelope $\mathrm{Abs}(\mathfrak{g})$ is given by the following coequalizer 
\[
\begin{tikzcd}[column sep=4pc,row sep=4pc]
\mathrm{Coeq}\Bigg(\displaystyle \widehat{\mathscr{S}}^c(\mathrm{B}^{\mathrm{s.a}}\ucom) \circ \mathscr{S}(\widehat{\Omega}^{\mathrm{s.c}}\ucom^*)(\mathfrak{g}) \arrow[r,"\gamma_\mathfrak{g}",shift right=1.1ex,swap]  \arrow[r,"\psi_\mathfrak{g}"{name=SD},shift left=1.1ex ]
&\widehat{\mathscr{S}}^c(\mathrm{B}^{\mathrm{s.a}}\ucom)(\mathfrak{g})\Bigg)~,
\end{tikzcd}
\]
in the category of pdg modules, where $\psi_\mathfrak{g}$ is constructed using $\iota_\mathfrak{g}$ and the monad structure of the free pdg absolute $\mathcal{L}_\infty$-algebra functor. See Section \ref{Section: Absolute algebras and contramodules} for more examples of this type of constructions.
\end{Remark}

\textbf{Completeness.} The coradical filtration on the conilpotent curved cooperad $\text{B}^{\mathrm{s.a}}\ucom$ induces a canonical filtration on any curved absolute $\mathcal{L}_\infty$-algebra. See Definition \ref{def: canonical filtration C alg}. In this case, this \textit{canonical filtration} on $\mathfrak{g}$ is given by
\[
\mathrm{W}_\omega\mathfrak{g} \coloneqq \mathrm{Im}\left(\gamma_\mathfrak{g} \vert_{\mathscr{F}_\omega} : \prod_{n \geq 0} \mathscr{F}_\omega\widehat{\Omega}^{\mathrm{s.a}}\ucom^*(n) ~ \widehat{\otimes}_{\mathbb{S}_n}~  \mathfrak{g}^{\otimes n} \longrightarrow \mathfrak{g} \right)~,
\]
for $\omega \geq 0$, where $\mathscr{F}_\omega\widehat{\Omega}^{\mathrm{s.a}}\ucom^*(n)$ is spanned by corked rooted trees with a number of internal edges greater or equal to $\omega$. Therefore an element $g$ is in $\mathrm{W}_{\omega_0}\mathfrak{g}$ if and only if it can be written as
\[
g = \gamma_\mathfrak{g}\left(\sum_{n\geq 0} \sum_{\omega \geq \omega_0} \sum_{\tau \in \mathrm{CRT}_n^\omega} \lambda_\tau \tau(g_{i_1}, \cdots, g_{i_n})\right)~.
\]
One can check that $\mathfrak{g}/\mathrm{W}_\omega\mathfrak{g}$ has an unique curved absolute $\mathcal{L}_\infty$-algebra structure given by the quotient structural map. 

\begin{Remark}
Condition \ref{pdg condition} implies that $d_\mathfrak{g}(\mathrm{W}_\omega\mathfrak{g}) \subset \mathrm{W}_\omega\mathfrak{g}$ and therefore the pre-differential $d_\mathfrak{g}$ is \textit{continuous} with respect to the canonical filtration.
\end{Remark}

\begin{Definition}[Complete curved absolute $\mathcal{L}_\infty$-algebra] 
Let $(\mathfrak{g},\gamma_\mathfrak{g}, d_\mathfrak{g})$ be a curved absolute $\mathcal{L}_\infty$-algebra. It is \textit{complete} if the canonical epimorphism
\[
\varphi_\mathfrak{g}: \mathfrak{g} \twoheadrightarrow \lim_{\omega} \mathfrak{g}/\mathrm{W}_\omega\mathfrak{g}
\]
is an isomorphism of curved absolute $\mathcal{L}_\infty$-algebras.
\end{Definition}

\begin{Remark}
The map $\varphi_\mathfrak{g}$ is always an epimorphism by Proposition \ref{prop: varphi is an epi}, this phenomenon is explained in Remark \ref{Remark: varphi is an epimorphism}.
\end{Remark}

Let $(\mathfrak{g},\gamma_\mathfrak{g}, d_\mathfrak{g})$ be a complete curved absolute $\mathcal{L}_\infty$-algebra. One can write
\[
\gamma_\mathfrak{g}\left(\sum_{n\geq 0} \sum_{\omega \geq 0} \sum_{\tau \in \mathrm{CRT}_n^\omega} \lambda_\tau \tau(g_{i_1}, \cdots, g_{i_n})\right) =  
\sum_{\omega \geq 0} \gamma_\mathfrak{g}\left(\sum_{n\geq 0} \sum_{\tau \in \mathrm{CRT}_n^\omega} \lambda_\tau \tau(g_{i_1}, \cdots, g_{i_n})\right)~,
\]
using the fact that the canonical filtration on $\mathfrak{g}$ is complete.

\begin{Definition}[Maurer--Cartan element of a curved absolute $\mathcal{L}_\infty$-algebra]
Let $(\mathfrak{g},\gamma_\mathfrak{g}, d_\mathfrak{g})$ be a curved absolute $\mathcal{L}_\infty$-algebra. A \textit{Maurer--Cartan element} $\alpha$ is an element in $\mathfrak{g}$ of degree $0$ which satisfies the following equation
\[
\gamma_\mathfrak{g}\left(\sum_{n \geq 0,~n\neq 1} \frac{c_n(\alpha, \cdots, \alpha)}{n!} \right) + d_\mathfrak{g}(\alpha) = 0~.
\]
The set of Maurer--Cartan elements in $\mathfrak{g}$ is denoted by $\mathcal{MC}(\mathfrak{g})~.$
\end{Definition}

\begin{Remark}
The coefficients $1/n!$ appear because of the isomorphism between invariants and coinvariants in Lemma \ref{lemma: iso invariants avec les coinvariants}. When one considers the renormalization of this equation as in Remark \ref{Remark: renormalization}, these coefficients disappear.
\end{Remark}

\begin{Remark}
In a curved absolute $\mathcal{L}_\infty$-algebra, the element $0$ is not in general a Maurer--Cartan element since:
\[
\gamma_\mathfrak{g}\left(\sum_{n \geq 0,~n\neq 1} \frac{c_n(0, \cdots, 0)}{n!} \right) + d_\mathfrak{g}(0) = \gamma_\mathfrak{g}(c_0) \neq 0~.
\]
Therefore the set $\mathcal{MC}(\mathfrak{g})$ is not canonically pointed, neither non-necessarily non-empty.
\end{Remark}

\begin{Remark}
Notice that if one takes $\alpha$ to be a degree $0$ element of $\mathfrak{g}$, the value of 
\[
\gamma_\mathfrak{g}\left(\sum_{n \geq 0,~n\neq 1} \frac{c_n(\alpha, \cdots, \alpha)}{n!} \right)~.
\]
is not determined \textit{in general} by the values of the elementary operations $l_n(\alpha,\cdots, \alpha)$. It can happen that $l_n(\alpha,\cdots, \alpha) = 0$ for all $n \geq 0$ but $\alpha$ is \textbf{not} a Maurer--Cartan element.
\end{Remark}

\begin{lemma}\label{lemma: splitting of the Maurer-Cartan equation}
Let $(\mathfrak{g},\gamma_\mathfrak{g}, d_\mathfrak{g})$ be a complete curved absolute $\mathcal{L}_\infty$-algebra. Any degree $0$ element $\alpha$ in $\mathrm{W}_1 \mathfrak{g}$ satisfies
\[
\gamma_\mathfrak{g}\left(\sum_{n \geq 0,~n\neq 1} \frac{c_n(\alpha, \cdots, \alpha)}{n!} \right) = \sum_{n \geq 0,~n\neq 1} \frac{l_n(\alpha, \cdots, \alpha)}{n!} ~.
\]
\end{lemma}

\begin{proof}
Since $\alpha$ is in $\mathrm{W}_1 \mathfrak{g}$, it can be written as 
\[
\alpha = \gamma_\mathfrak{g}\left(\sum_{k\geq 0} \sum_{\omega \geq 1} \sum_{\tau \in \mathrm{CRT}_k^\omega} \lambda_\tau \tau(\bar{g}_k)\right)~.
\]
where $\bar{g}_k$ is an $k$-tuple in $\mathfrak{g}$. Using the associativity condition \ref{associativity condition}, one can write 
\[
\gamma_\mathfrak{g}\left(\sum_{n \geq 0,~n\neq 1} \frac{c_n(\alpha, \cdots, \alpha)}{n!} \right) = \gamma_\mathfrak{g}\left(\sum_{n \geq 0,~n\neq 1} \frac{1}{n!}
c_n \left( \sum_{\substack{k\geq 0 \\ \omega \geq 1}} \sum_{\tau \in \mathrm{CRT}_k^\omega} \lambda_\tau \tau(\bar{g}_k), \cdots, \sum_{\substack{k\geq 0 \\ \omega \geq 1}} \sum_{\tau \in \mathrm{CRT}_k^\omega} \lambda_\tau \tau(\bar{g}_k) \right) \right)~,
\]
where the term 
\[
c_n \left( \sum_{\substack{k\geq 0 \\ \omega \geq 1}} \sum_{\tau \in \mathrm{CRT}_k^\omega} \lambda_\tau \tau(\bar{g}_k), \cdots, \sum_{\substack{k\geq 0 \\ \omega \geq 1}} \sum_{\tau \in \mathrm{CRT}_k^\omega} \lambda_\tau \tau(\bar{g}_k) \right)
\]
can be rewritten as an infinite sum of terms of weight greater or equal to $k+1$. Thus, using the fact that $\mathfrak{g}$ is complete, we have that 
\[
\gamma_\mathfrak{g}\left(\sum_{n \geq 0,~n\neq 1} \frac{1}{n!}
c_n \left( \sum_{\substack{k\geq 0 \\ \omega \geq 1}} \sum_{\tau \in \mathrm{CRT}_k^\omega} \lambda_\tau \tau(\bar{g}_k), \cdots, \sum_{\substack{k\geq 0 \\ \omega \geq 1}} \sum_{\tau \in \mathrm{CRT}_k^\omega} \lambda_\tau \tau(\bar{g}_k) \right) \right) = \sum_{n \geq 0,~n\neq 1} \frac{\gamma_\mathfrak{g}(c_n(\alpha, \cdots, \alpha))}{n!}~.
\]
\end{proof}

\textbf{Graded homology groups.} Let $(\mathfrak{g},\gamma_\mathfrak{g}, d_\mathfrak{g})$ be a curved absolute $\mathcal{L}_\infty$-algebra. Recall that Condition \ref{curved condition} says that $d_\mathfrak{g}^2(-) = l_2(l_0,-)$. Thus 
\[
d_\mathfrak{g}^2(\mathrm{W}_\omega \mathfrak{g}) \subseteq \mathrm{W}_{\omega+1} \mathfrak{g}~.
\]
This implies that
\[
\mathrm{gr}_{\omega}(\mathfrak{g}) \coloneqq \mathrm{W}_\omega \mathfrak{g}/\mathrm{W}_{\omega+1} \mathfrak{g}~,
\]
forms a chain complex endowed with the differential induced by $d_\mathfrak{g}$, for all $\omega \geq 0$.

\begin{Remark}
Notice that any morphism $f: \mathfrak{g} \longrightarrow \mathfrak{h}$ of curved absolute $\mathcal{L}_\infty$-algebras is \textit{continuous} with respect to the canonical filtration, i.e: $f(\mathrm{W}_\omega \mathfrak{g}) \subset \mathrm{W}_\omega \mathfrak{h}$. Therefore the morphism of dg modules $\mathrm{gr}_{\omega}(f)$ is always well-defined.
\end{Remark}

\begin{Definition}[Graded quasi-isomorphisms] 
Let $f: \mathfrak{g} \longrightarrow \mathfrak{h}$ be a morphism between two curved absolute $\mathcal{L}_\infty$-algebras. It is a \textit{graded quasi-isomorphism} if 
\[
\begin{tikzcd}
\mathrm{gr}_{\omega}(f): \mathrm{gr}_{\omega}(\mathfrak{g}) \arrow[r]
&\mathrm{gr}_{\omega}(\mathfrak{h})
\end{tikzcd}
\]
is a quasi-isomorphism of dg modules, for all $\omega \geq 0$.
\end{Definition}

\subsection{On $u\mathcal{CC}_\infty$-coalgebras}
We study the Koszul dual notion of dg coalgebras over the dg operad $\Omega\mathrm{B}^{\mathrm{s.a}}\ucom$. We denote these coalgebras $u\mathcal{CC}_\infty$-coalgebras, since they correspond to counit cocommutative coalgebras relaxed up to homotopy. We endow the category of $u\mathcal{CC}_\infty$-coalgebras with a model structure where weak equivalences are given by quasi-isomorphisms and cofibrations by degree-wise monomorphisms. We use the canonical twisting morphism $\iota: \mathrm{B}^{\mathrm{s.a}}\ucom \longrightarrow \Omega\mathrm{B}^{\mathrm{s.a}}\ucom$ to construct a complete Bar-Cobar adjunction between $u\mathcal{CC}_\infty$-coalgebras and curved absolute $\mathcal{L}_\infty$-algebras. Using this adjunction, we transfer the aforementioned model category structure to curved absolute $\mathcal{L}_\infty$-algebras. 

\begin{Definition}[$u\mathcal{CC}_\infty$-coalgebra]
A $u\mathcal{CC}_\infty$\textit{-coalgebra} $C$ is the data $(C,\Delta_C,d_C)$ of a dg $\Omega\mathrm{B}^{\mathrm{s.a}}\ucom$-coalgebra.
\end{Definition}

\begin{lemma}
The data of a $u\mathcal{CC}_\infty$-coalgebra structure $\Delta_C$ on a dg module $(C,d_C)$ is equivalent to the data of elementary decomposition maps
\[
\Big\{ \Delta_\tau: C \longrightarrow C^{\otimes n} \Big\}
\]

of degree $\omega - 1$ for all corked rooted trees $\tau$ in $\mathrm{CRT}_n^\omega$ for all $n \geq 0$ and $\omega \geq 0$. These operations are subject to the following condition: let $E(\tau)$ denote the set of internal edges of $\tau$, then 
\[
\sum_{e \in E(\tau)} \Delta_{\tau_{(1)}} \circ_e \Delta_{\tau_{(2)}} - (-1)^{\omega -1} \sum_{e \in E(\tau)} \Delta_{\tau^{e}} + d_{C^{\otimes n}} \circ \Delta_\tau - (-1)^{\omega -1} \Delta_{\tau} \circ d_C = 0~,
\]
where in the first term the corked rooted tree is split along $e$ into $\tau_{(1)}$ and $\tau_{(2)}$ and where $\tau^{e}$ denotes the corked rooted tree obtained by contracting the internal edge $e$. 
\end{lemma}

\begin{proof}
The data of a $u\mathcal{CC}_\infty$-coalgebra structure $\Delta_C$ on a dg module $(C,d_C)$ is equivalent to the data of a morphism of dg operads
\[
\Delta_C: \Omega\mathrm{B}^{\mathrm{s.a}}\ucom \longrightarrow \mathrm{Coend}_C~,
\]
which in turn is equivalent to the data of a curved twisting morphism
\[
\delta_C: \mathrm{B}^{\mathrm{s.a}}\ucom \longrightarrow \mathrm{Coend}_C~.
\]
Thus the operations $\Delta_\tau$ are given by $\delta_C(\tau)$ and their relationships are given by the Maurer-Cartan equation that $\delta_C$ satisfies.
\end{proof}

\begin{Notation}
We denote by $\mathrm{PCRT}_n^{(\omega, \nu)}$ the set of partitioned corked rooted trees with $\omega$ vertices, $\nu$ parenthesis, of arity $n$. Such a tree amounts to the data of a corked rooted tree $\tau$ together with the data of a partition of the set of vertices of $\tau$.
\end{Notation}

\begin{Remark}
The dg operad $\Omega\mathrm{B}^{\mathrm{s.a}}\ucom$ admits a basis given by partitioned corked rooted trees. The full structure of a $u\mathcal{CC}_\infty$-coalgebra is given by decomposition maps 
\[
\Big\{ \Delta_{\tau}: C \longrightarrow C^{\otimes n} \Big\}
\]

of degree $\omega - \nu$ if $\tau$ is in $\mathrm{PCRT}_n^{(\omega, \nu)}$. Nevertheless, any such decomposition map is obtained as the composition of the elementary decomposition maps given by the sub-corked rooted trees contained inside each partition.
\end{Remark}

Let $(C,d_C)$ be a dg module endowed with a family of decomposition maps $\{\Delta_\tau\}$ for any partitioned corked rooted tree $\tau$. This data allows us to construct a morphism of dg modules 

\[
\begin{tikzcd}[column sep=4pc,row sep=-0.5pc]
C \arrow[r,"\Delta_{C}"]
&\displaystyle \prod_{n \geq 0} \mathrm{Hom}_{\mathbb{S}_n}\left(\Omega\mathrm{B}^{\mathrm{s.a}}\ucom(n), C^{\otimes n}\right)\\
c \arrow[r,mapsto]
&\Big[ \mathrm{ev}_c: \tau \mapsto \Delta_{\tau}(c) \Big]~.
\end{tikzcd}
\]

which satisfies the axioms of Definition \ref{def: P coalgebra}. This in particular gives a morphism
\[
\Delta_C^0: C \longrightarrow \mathrm{Hom}_{\mathsf{gr}} \left(\Omega\mathrm{B}^{\mathrm{s.a}}\ucom(0), \kk \right)
\]
which factors through the terminal coalgebra $\mathscr{C}(u\mathcal{CC}_\infty)(0) \rightarrowtail \mathrm{Hom}_{\mathsf{gr}} \left(\Omega\mathrm{B}^{\mathrm{s.a}}\ucom(0), \kk \right)$. This is the terminal morphism for any $u\mathcal{CC}_{\infty}$-coalgebra. Notice that $\kk$ forms a $u\mathcal{CC}_{\infty}$-coalgebra with its natural $u\mathcal{C}om$-coalgebra structure and that there is a map $k: \kk \longrightarrow \mathscr{C}(u\mathcal{CC}_\infty)(0)$.

\begin{lemma}\label{lemma: L(0) qi a k}
The morphism $k: \kk \qi \mathscr{C}(u\mathcal{CC}_\infty)(0)$ is a quasi-isomorphism of $u\mathcal{CC}_\infty$-coalgebras.
\end{lemma}

\begin{proof}
We know that $\Omega\mathrm{B}^{\mathrm{s.a}}\ucom(0) \qi \kk^*$ is a quasi-isomorphism of $u\mathcal{CC}_\infty$-algebras, therefore its adjoint map via the $(-)^\circ \dashv (-)^*$ adjunction must also be a quasi-isomorphism since the derived unit/counits are quasi-isomorphism on objects with degree-wise finite dimension homology by Proposition \ref{Prop: plongement pleinement fidele de Sweedler}. Thus
\[
\kk \qi (\Omega\mathrm{B}^{\mathrm{s.a}}\ucom(0))^\circ \cong \mathscr{C}(u\mathcal{CC}_\infty)(0)
\]
is a quasi-isomorphism of $u\mathcal{CC}_\infty$-coalgebras.
\end{proof}

\begin{Definition}[Strictly counital $u\mathcal{CC}_\infty$-coalgebra]
A $u\mathcal{CC}_\infty$-coalgebra $C$ is \textit{strictly counital} if the elementary decomposition maps 
\[
\Delta_{\tau}: C \longrightarrow C^{\otimes n}
\]
are the zero morphism for any $\tau$ in $\mathrm{CRT}_n^\omega$ with at least one cork, except for $\epsilon: C \longrightarrow \kk$, which is the decomposition associated to the single cork.
\end{Definition}

\begin{Remark}
Let $C$ be a strictly counital $u\mathcal{CC}_\infty$-coalgebra. Then the terminal morphism 
\[
\Delta_C^0: C \longrightarrow \mathscr{C}(u\mathcal{CC}_\infty)(0)
\]
factors through the morphism $k: \kk \longrightarrow \mathscr{C}(u\mathcal{CC}_\infty)(0)$. Therefore $\kk$ is the terminal object in the full sub-category of strictly counital $u\mathcal{CC}_\infty$-coalgebras.
\end{Remark}

\begin{Proposition}\label{prop: coreflective sub-category}
The full sub-category of strictly counital $u\mathcal{CC}_\infty$-coalgebras is a coreflective sub-category of the category of $u\mathcal{CC}_\infty$-coalgebras.
\end{Proposition}

\begin{proof}
One can encode strictly counital $u\mathcal{CC}_\infty$-coalgebras with a dg operad. This dg operad is given by $\Omega \mathrm{B}^{s.a} u\mathcal{C}om$ modulo all the elements corresponding to rooted trees with at least one cork. This dg operad is manifestly not cofibrant. Nevertheless, it shows that this sub-category is comonadic, therefore the inclusion functor admits a right adjoint.
\end{proof}

\begin{Definition}[Coaugmented $u\mathcal{CC}_\infty$-coalgebra]
A \textit{coaugmented} $u\mathcal{CC}_\infty$-coalgebra $(C,\Delta_C,d_C,\eta_C)$ is the data of a $u\mathcal{CC}_\infty$-coalgebra together with a $u\mathcal{CC}_\infty$-coalgebra morphism 
\[
\eta_C: \mathscr{C}(u\mathcal{CC}_\infty)(0) \longrightarrow C
\]
which is a section of the terminal morphism $\Delta_C^0$.
\end{Definition}

\begin{Remark}
A coaugmented strictly counital $u\mathcal{CC}_\infty$-coalgebra is also called a \textit{pointed} strictly counital $u\mathcal{CC}_\infty$-coalgebra.
\end{Remark}

\begin{Definition}[Group-like element]
Let $(C,\Delta_C,d_C)$ be a $u\mathcal{CC}_\infty$-coalgebra. A \textit{group-like element} is the data of a morphism of $u\mathcal{CC}_\infty$-coalgebras $\alpha: \kk \longrightarrow C$.
\end{Definition}

\begin{Remark}
If $C$ is strictly counital, the data of a coaugmentation is equivalent to the data of a group-like element.
\end{Remark}

\begin{Definition}[$\mathcal{CC}_\infty$-coalgebra]
A $\mathcal{CC}_\infty$-coalgebra $C$ is the data $(C,\Delta_C,d_C)$ of a dg $\Omega \mathrm{B} \mathcal{C}om$-coalgebra.
\end{Definition}

\begin{Proposition}\label{prop: equivalence entre les strictes pointes et les CC}
There is an adjunction 
\[
\begin{tikzcd}[column sep=7pc,row sep=3pc]
\mathsf{Strict}~u\mathcal{CC}_\infty\textsf{-}\mathsf{coalg}_{\bullet}    \arrow[r, shift left=.75ex, "\mathrm{Ker}(\epsilon)"{name=U}]  
&\mathcal{CC}_\infty\textsf{-}\mathsf{coalg}~, \arrow[l, shift left=1.1ex, "(-) \oplus \kk"{name=F}] 
\arrow[phantom, from=F, to=U, , "\dashv" rotate=-90]
\end{tikzcd}
\]
between the category of pointed (coaugmented) strictly counital $u\mathcal{CC}_\infty$-coalgebras and the category of $\mathcal{CC}_\infty$-coalgebras.

\begin{enumerate}
\medskip

\item The left adjoint $\mathrm{Ker}(\epsilon)$ is given by considering the $\mathcal{CC}_\infty$-coalgebra structure on the coaugmentation coideal of a coaugmented strictly counital $u\mathcal{CC}_\infty$-coalgebra.

\medskip 

\item The right adjoint $(-) \oplus \kk$ adds cofreely a strict counit to a $\mathcal{CC}_\infty$-coalgebra.

\medskip 

\end{enumerate}

Furthermore, this adjunction is an equivalence of categories. 
\end{Proposition}

\begin{proof}
Straightforward to check.
\end{proof}

\subsection{Model category structures} Since the dg operad encoding $u\mathcal{CC}_\infty$-coalgebras is cofibrant, the category of $u\mathcal{CC}_\infty$-coalgebras admits a canonical model category structure, left-transferred from the category of dg modules. See Section \ref{Section: Complete Bar-Cobar curved} for an overview of the general results used in this subsection.

\begin{Proposition}\label{prop: model structure on uCC coalgebras}
There is a model category structure on the category of $u\mathcal{CC}_\infty$-coalgebras left-transferred along the cofree-forgetful adjunction
\[
\begin{tikzcd}[column sep=7pc,row sep=3pc]
\mathsf{dg}\textsf{-}\mathsf{mod} \arrow[r, shift left=1.1ex, "\mathscr{C}(u\mathcal{CC}_\infty)(-)"{name=F}]      
&u\mathcal{CC}_\infty\textsf{-}\mathsf{coalg}~, \arrow[l, shift left=.75ex, "\mathrm{U}"{name=U}]
\arrow[phantom, from=F, to=U, , "\dashv" rotate=90]
\end{tikzcd}
\]
where 
\begin{enumerate}
\medskip
\item the class of weak equivalences is given by quasi-isomorphisms,
\medskip
\item the class of cofibrations is given by degree-wise monomorphisms,
\medskip 
\item the class of fibrations is given by right lifting property with respect to acyclic cofibrations.
\end{enumerate}
\end{Proposition}

\begin{proof}
Consequence of the results explained in Section \ref{Section: Complete Bar-Cobar curved} applied to this particular case, since $\Omega\mathrm{B}^{\mathrm{s.a}}\ucom$ is cofibrant in the model category of dg operads.
\end{proof}

Furthermore, using the canonical curved twisting morphism 
\[
\iota: \mathrm{B}^{\mathrm{s.a}}\ucom \longrightarrow \Omega\mathrm{B}^{\mathrm{s.a}}\ucom~,
\]
we can construct a complete Bar-Cobar adjunction relating $u\mathcal{CC}_\infty$-coalgebras and curved absolute $\mathcal{L}_\infty$-algebras.

\begin{Proposition}[Complete Bar-Cobar adjunction]
The curved twisting morphism $\iota$ induces a complete Bar-Cobar adjunction
\[
\begin{tikzcd}[column sep=5pc,row sep=3pc]
            u\mathcal{CC}_\infty\text{-}\mathsf{coalg} \arrow[r, shift left=1.1ex, "\widehat{\Omega}_{\iota}"{name=F}] & \mathsf{curv}~\mathsf{abs}~\mathcal{L}_\infty\text{-}\mathsf{alg}^{\mathsf{comp}} \arrow[l, shift left=.75ex, "\widehat{\mathrm{B}}_{\iota}"{name=U}]
            \arrow[phantom, from=F, to=U, , "\dashv" rotate=-90]
\end{tikzcd}
\]
between the category of $u\mathcal{CC}_\infty$-coalgebras and the category of complete curved absolute $\mathcal{L}_\infty$-algebras. 
\end{Proposition} 

\begin{proof}
Consequence of the results explained in Section \ref{Section: Complete Bar-Cobar curved} applied to this particular case.
\end{proof}

Using this adjunction, one can transfer the model category structure on $u\mathcal{CC}_\infty$-coalgebras in order to endow curved absolute $\mathcal{L}_\infty$-algebras with a model category structure.

\begin{theorem}\label{thm: Equivalence the Quillen CC infini et absolues}
There is a model category structure on the category of curved absolute $\mathcal{L}_\infty$-algebras right-transferred along the complete Bar-Cobar adjunction
\[
\begin{tikzcd}[column sep=5pc,row sep=3pc]
            u\mathcal{CC}_\infty\text{-}\mathsf{coalg} \arrow[r, shift left=1.1ex, "\widehat{\Omega}_{\iota}"{name=F}] & \mathsf{curv}~\mathsf{abs}~\mathcal{L}_\infty\text{-}\mathsf{alg}^{\mathsf{comp}} \arrow[l, shift left=.75ex, "\widehat{\mathrm{B}}_{\iota}"{name=U}]
            \arrow[phantom, from=F, to=U, , "\dashv" rotate=-90]
\end{tikzcd}
\]
where 
\begin{enumerate}
\medskip

\item the class of weak equivalences is given by morphisms $f$ such that $\widehat{\mathrm{B}}_\iota(f)$ is a quasi-isomorphism,

\medskip

\item the class of fibrations is given by morphisms $f$ such that $\widehat{\mathrm{B}}_\iota(f)$ is a fibration,

\medskip

\item  and the class of cofibrations is given by left lifting property with respect to acyclic fibrations.

\medskip
\end{enumerate}
Furthermore, this complete Bar-Cobar adjunction is a Quillen equivalence. 
\end{theorem}

\begin{proof}
Consequence of the results explained in Section \ref{Section: Complete Bar-Cobar curved} applied to this particular case.
\end{proof}

One can give a complete characterization of the fibrations in this model category structure as follows.

\begin{Proposition}
A morphism of curved absolute $\mathcal{L}_\infty$-algebras is a fibration if and only if it is a degree-wise epimorphism.
\end{Proposition}

\begin{proof}
Follows from \cite[Proposition 10.16]{grignoulejay18}.
\end{proof}

\begin{Remark}
In particular, all complete curved absolute $\mathcal{L}_\infty$-algebras are fibrant in this model structure.
\end{Remark}

There is a particular kind of weak-equivalences between curved absolute $\mathcal{L}_\infty$-algebras which admit an easy description.

\begin{Proposition}
Let $f: \mathfrak{g} \longrightarrow \mathfrak{h}$ be a graded quasi-isomorphism between two complete curved absolute $\mathcal{L}_\infty$-algebras. Then 
\[
\widehat{\mathrm{B}}_\iota(f): \widehat{\mathrm{B}}_\iota \mathfrak{g}  \qi \widehat{\mathrm{B}}_\iota \mathfrak{h} 
\]
is a quasi-isomorphism of $u\mathcal{CC}_\infty$-coalgebras, and $f$ is a weak-equivalence of complete curved absolute $\mathcal{L}_\infty$-algebras.
\end{Proposition}

\begin{proof}
Follows from \cite[Theorem 10.25]{grignoulejay18}.
\end{proof}

The same results hold for the categories of $\mathcal{CC}_\infty$-coalgebras and absolute $\mathcal{L}_\infty$-algebras.

\begin{Notation}
We denote $\widehat{\Omega}_{\iota}^{\flat} \dashv \widehat{\mathrm{B}}^\flat_{\iota}$ the complete Bar-Cobar adjunction relative to $\iota: \mathrm{B}\mathcal{C}om \longrightarrow \Omega\mathrm{B}\mathcal{C}om$ to distinguish it from the previous one. 
\end{Notation}

\begin{Proposition}\label{thm: structure de modeles sur les CC pas unitaires}
There is a Quillen equivalence
\[
\begin{tikzcd}[column sep=5pc,row sep=3pc]
            \mathcal{CC}_\infty\text{-}\mathsf{coalg} \arrow[r, shift left=1.1ex, "\widehat{\Omega}_{\iota}^\flat"{name=F}] & \mathsf{abs}~\mathcal{L}_\infty\text{-}\mathsf{alg}^{\mathsf{comp}} \arrow[l, shift left=.75ex, "\widehat{\mathrm{B}}^\flat_{\iota}"{name=U}]
            \arrow[phantom, from=F, to=U, , "\dashv" rotate=-90]
\end{tikzcd}
\]
where the model category structure on the left hand side is given by quasi-isomorphisms and degree-wise monomorphisms and the model category structure on the right hand side is the transferred structure.
\end{Proposition}

\begin{Remark}
By Proposition \ref{prop: equivalence entre les strictes pointes et les CC}, there is also a model structure on pointed strict $u\mathcal{CC}_\infty$-coalgebras where weak-equivalences are given by quasi-isomorphisms and where cofibrations are degree-wise monomorphisms. Furthermore, the inclusion functor into all $u\mathcal{CC}_\infty$-coalgebras is a left Quillen functor.
\end{Remark}

\subsection{Tensor product of $u\mathcal{CC}_\infty$-coalgebras}
The tensor product of two $u\mathcal{CC}_\infty$-coalgebras can naturally be endowed with a $u\mathcal{CC}_\infty$-coalgebra structure. This gives a closed monoidal structure on the category of $u\mathcal{CC}_\infty$-coalgebras which is compatible with its model category structure.

\begin{Proposition}
The dg operad $\Omega\mathrm{B}^{\mathrm{s.a}}u\mathcal{C}om$ is a dg Hopf operad, meaning there exists a diagonal morphism of operads
\[
\Delta_{\Omega\mathrm{B}\mathcal{C}om}: \Omega\mathrm{B}u\mathcal{C}om \longrightarrow \Omega\mathrm{B}u\mathcal{C}om \otimes_{H} \Omega\mathrm{B}u\mathcal{C}om
\]
which is coassociative.
\end{Proposition}

\begin{proof}
The operad $u\mathcal{C}om$ can be obtained as the cellular chains on the operad $\mathrm{uCom}$ in the category of sets given by $\{*\}$ in non-negative arities. Any set-theoretical operad is a cocommutative Hopf operad on the nose. By Theorem \ref{thm: Boardman-Vogt and cellular chains}, one has an isomorphism of dg operads

\[
C_*^c(\mathrm{W}(\mathrm{uCom})) \cong \Omega\mathrm{B}^{\mathrm{s.a}}\left(C_*^c(\mathrm{uCom})\right)~,
\]
\vspace{0.1pc}

where $\mathrm{W}$ denotes the Boardmann-Vogt construction of $\mathrm{uCom}$ seen as a discrete topological operad. Given a choice a cellular approximation of the diagonal on the topological interval $\mathrm{I}$, one has a coassociative Hopf operad structure on $\mathrm{W}(\mathrm{uCom})$ which might not be cocommutative. The cellular chain functor $C_*^c(-)$ is both lax and colax (Eilenberg-Ziber and Alexander-Whitney maps) in a compatible way by \cite[Proposition 73]{grignou2022mapping}, therefore by \cite[Proposition 70]{grignou2022mapping} it sends Hopf topological operads to dg Hopf operads.
\end{proof}

\begin{Corollary}\label{cor: tensor product of uCC coalgebras}
The category $u\mathcal{CC}_\infty$-coalgebras can be endowed with a monoidal structure given by the tensor product of the underlying dg modules. 
\end{Corollary}

\begin{proof}
Consider two $u\mathcal{CC}_\infty$-coalgebras $(C_1,\Delta_{C_1},d_{C_1})$ and $(C_2,\Delta_{C_2},d_{C_2})$. These $u\mathcal{CC}_\infty$-coalgebra structures amount to two morphisms of dg operads $f_1: \Omega\mathrm{B}\mathcal{C}om \longrightarrow \mathrm{Coend}(C_1)$ and $f_2: \Omega\mathrm{B}\mathcal{C}om \longrightarrow \mathrm{Coend}(C_2)$. Using the Hopf structure on gets 

\[
\begin{tikzcd}[column sep=2pc,row sep=0pc]
\Omega\mathrm{B}\mathcal{C}om \arrow[r,"\Delta_{\Omega\mathrm{B}\mathcal{C}om}"]
&\Omega\mathrm{B}\mathcal{C}om \otimes_{H} \Omega\mathrm{B}\mathcal{C}om \arrow[r, "f_1 ~ \otimes ~ f_2"]
&\mathrm{Coend}(C_1) \otimes_{H} \mathrm{Coend}(C_2) \arrow[r]
&\mathrm{Coend}(C_1 \otimes C_2)~,
\end{tikzcd}
\]
\vspace{0.1pc}

which endows the dg module $(C_1 \otimes C_2, d_{C_1} \otimes d_{C_2})$ with a $u\mathcal{CC}_\infty$-coalgebra structure. The counitality and coassociativity of $\Delta_{\Omega\mathrm{B}\mathcal{C}om}$ ensure that the category of $u\mathcal{CC}_\infty$-coalgebra together with the tensor product forms a monoidal category.
\end{proof}

\begin{Remark}
The dg Hopf operad structure on $\Omega\mathrm{B}^{\mathrm{s.a}}u\mathcal{C}om$ constructed here extends the dg Hopf structure of $u\mathcal{C}om$, therefore the tensor product of two counital cocommutative coalgebras seen as $u\mathcal{CC}_\infty$-coalgebras coincides with the usual structure on the tensor product of two counital cocommutative coalgebras.
\end{Remark}

\begin{Proposition}\label{prop: monoidal model category}
The category of $u\mathcal{CC}_\infty$-coalgebras together with their tensor product forms a monoidal model category.
\end{Proposition}

\begin{proof}
It is straightforward to check the pushout-product axiom, the unit axiom is automatic since all objects are cofibrant, see \cite{HoveyModel} for more details on this.
\end{proof}

In fact, the tensor product of two coalgebras over a dg Hopf operad always forms a \textit{biclosed} monoidal category. We specify the general constructions of the forthcoming paper \cite{grignou2022mapping} in the particular cases of interest for us. The construction of \textit{loc.cit} holds for a wider variety of cases.

\begin{Proposition}
The category of dg $u\mathcal{CC}_\infty$-coalgebras is a biclosed monoidal category, meaning that there exists an internal hom bifunctor 

\[
\{-,-\}: \left(u\mathcal{CC}_\infty\text{-}\mathsf{coalg}\right)^{\mathsf{op}} \times u\mathcal{CC}_\infty\text{-}\mathsf{coalg} \longrightarrow u\mathcal{CC}_\infty\text{-}\mathsf{coalg}
\]

and, for any triple of $u\mathcal{CC}_\infty)$-coalgebras $C,D,E$, there exists isomorphisms

\[
\mathrm{Hom}_{u\mathcal{CC}_\infty\text{-}\mathsf{coalg}}(C \otimes D, E) \cong \mathrm{Hom}_{u\mathcal{CC}_\infty\text{-}\mathsf{coalg}}(C, \{D,E\}) \cong \mathrm{Hom}_{u\mathcal{CC}_\infty\text{-}\mathsf{coalg}}(D, \{C,E\})~,
\]
which are natural in $C,D$ and $E$.
\end{Proposition}

\begin{proof}
Direct application of \cite[Corollary 20]{grignou2022mapping}.
\end{proof}

\begin{Remark}
Let $C_1$ and $C_2$ be two $u\mathcal{CC}_\infty$-coalgebras. The $u\mathcal{CC}_\infty$-coalgebra $\{C_1,C_2\}$ is given as the following equalizer 
\[
\begin{tikzcd}[column sep=4pc,row sep=4pc]
\mathrm{Eq}\Bigg(\mathscr{C}(u\mathcal{CC}_\infty)(\mathrm{hom}(C_1,C_2)) \arrow[r,"(\Delta_{C_2})_*",shift right=1.1ex,swap]  \arrow[r,"\varrho"{name=SD},shift left=1.1ex ]
&\mathscr{C}(u\mathcal{CC}_\infty)\left(\mathrm{hom}(C_1,\widehat{\mathscr{S}}^c(u\mathcal{CC}_\infty)(C_2) \right)\Bigg)~,
\end{tikzcd}
\]
where $\mathrm{hom}(C_1,C_2)$ denotes the dg module of graded morphisms, where $\varrho$ is a map constructed using the comonad structure of $\mathscr{C}(u\mathcal{CC}_\infty)$ and the Hopf structure of the dg operad $u\mathcal{CC}_\infty$. This kind of construction works for coalgebras over any Hopf dg operad. 
\end{Remark}
 
\begin{Definition}[Convolution curved absolute $\mathcal{L}_\infty$-algebra]
Let $C$ be a $u\mathcal{CC}_\infty$-coalgebra and let $\mathfrak{g}$ be a curved absolute $\mathcal{L}_\infty$-algebra. The pdg module of graded morphisms $(\mathrm{hom}(C,\mathfrak{g}),\partial)$ is endowed the following curved absolute $\mathcal{L}_\infty$-algebra structure.

\medskip

The structural map $\gamma_{\mathrm{hom}(C,\mathfrak{g})}$ is given by the following composition
\[
\begin{tikzcd}
\widehat{\mathscr{S}}^c(\mathrm{B}^{\mathsf{s.a}}u\mathcal{C}om)(\mathrm{hom}(C,\mathfrak{g})) \arrow[d, "\mathrm{coev}_C"]     \\
\mathrm{hom}\left(C, \widehat{\mathscr{S}}^c(\mathrm{B}^{\mathsf{s.a}}u\mathcal{C}om)(\mathrm{hom}(C,\mathfrak{g})) \otimes C \right) \arrow[d, "(\Delta_C)_*"] \\
\mathrm{hom} \left(C, \widehat{\mathscr{S}}^c(\mathrm{B}^{\mathsf{s.a}}u\mathcal{C}om)(\mathrm{hom}(C,\mathfrak{g})) \otimes \widehat{\mathscr{S}}^c(\Omega\mathrm{B}^{\mathsf{s.a}}u\mathcal{C}om)(C) \right) \arrow[d, "\xi"] \\
\mathrm{hom}\left(C, \widehat{\mathscr{S}}^c(\mathrm{B}^{\mathsf{s.a}}u\mathcal{C}om \otimes \Omega\mathrm{B}^{\mathsf{s.a}}u\mathcal{C}om)(\mathrm{hom}(C,\mathfrak{g}) \otimes C)\right) \arrow[d,"\widehat{\mathscr{S}}^c(\delta_{\mathrm{B}^{\mathrm{s.a}}u\mathcal{C}om})"] \\
\mathrm{hom}\left(C, \widehat{\mathscr{S}}^c(\mathrm{B}^{\mathsf{s.a}}u\mathcal{C}om)(\mathrm{hom}(C,\mathfrak{g}) \otimes C)\right) \arrow[d," \mathrm{ev}_C "] \\
\mathrm{hom}\left(C, \widehat{\mathscr{S}}^c(\mathrm{B}^{\mathsf{s.a}}u\mathcal{C}om)(\mathfrak{g})\right) \arrow[d,"(\gamma_\mathfrak{g})_* "] \\
\mathrm{hom}(C, \mathfrak{g})~, \\
\end{tikzcd}
\]
where $\mathrm{coev}_C$ and $\mathrm{ev}_C$ are respectively the unit and the counit of the tensor-hom adjunction, and where $\xi$ is the following natural inclusion 
\[
\begin{tikzcd}
\displaystyle \left(\prod_{n \geq 0} \mathrm{Hom}_{\mathbb{S}_n}(M(n), V^{\otimes n})\right) \otimes \left(\prod_{n \geq 0} \mathrm{Hom}_{\mathbb{S}_n}(N(n), W^{\otimes n})\right) \arrow[d,rightarrowtail] \\
\displaystyle \prod_{n \geq 0} \mathrm{Hom}_{\mathbb{S}_n}(M(n) \otimes N(n), (V \otimes W)^{\otimes n})~.
\end{tikzcd}
\]
Finally

\[
\delta_{\mathrm{B}^{\mathsf{s.a}}u\mathcal{C}om}: \mathrm{B}^{\mathsf{s.a}}u\mathcal{C}om \longrightarrow \mathrm{B}^{\mathsf{s.a}}u\mathcal{C}om \otimes \Omega\mathrm{B}^{\mathsf{s.a}}u\mathcal{C}om
\]
\vspace{0.1pc}

is a restriction of the diagonal $\Delta_{\Omega\mathrm{B}\mathcal{C}om}$.
\end{Definition}

Convolution curved absolute $\mathcal{L}_\infty$-algebras and the internal hom-set of $u\mathcal{CC}_\infty$-coalgebras are compatible in the following sense.

\begin{theorem}
Let $C$ be a $u\mathcal{CC}_\infty$-coalgebra and let $\mathfrak{g}$ be a curved absolute $\mathcal{L}_\infty$-algebra. There is an isomorphism of $u\mathcal{CC}_\infty$-coalgebras

\[
\left\{C, \widehat{\mathrm{B}}_\iota(\mathfrak{g}) \right\} \cong \widehat{\mathrm{B}}_\iota \left(\mathrm{hom}(C,\mathfrak{g}) \right)~,
\]
\vspace{0.1pc}

where $\mathrm{hom}(C,\mathfrak{g})$ denotes the convolution curved absolute $\mathcal{L}_\infty$-algebra of $C$ and $\mathfrak{g}$.
\end{theorem}

\begin{proof}
Direct application of \cite[Theorem 9]{grignou2022mapping}.
\end{proof}

\section{Higher absolute Lie theory}
In this section, we follow an analogue approach to \cite{robertnicoud2020higher} in order to integrate curved absolute $\mathcal{L}_\infty$-algebras. For the first time, the integration functor and its left adjoint form a Quillen adjunction between curved absolute $\mathcal{L}_\infty$-algebras and simplicial sets. The rest of this section is devoted to the study of the main properties of this adjunction.

\subsection{Dupont's contraction}\label{subsection: dupont contraction}
In \cite{DUPONT}, J-L Dupont proved that there is a homotopy contraction between the simplicial unital commutative algebra of polynomial differential forms on the geometrical simplicies and the simplicial sub-module of Whitney forms. Whitney forms on the simplicies are isomorphic to the cellular cochains on the simplicies and they are finite dimensional. This allows us to use the homotopy transfer theorem, and obtain a simplicial $u\mathcal{CC}_\infty$-coalgebra structure on the cellular chains of simplicies. This coalgebra is not conilpotent, hence we encode it with an operad. This enables us to use the complete Bar-Cobar construction in order to construct a commuting triangle of Quillen adjunctions. This approach allows us to extend and to refine automatically many of the standard results in the theory of integration. 

\medskip

The simplicial dg unital commutative algebra of piece-wise polynomial differential forms on the standard simplex $\Omega_\bullet$ is given by
\[
\Omega_n \coloneqq \frac{\kk[t_0,\cdots,t_n, dt_0,\cdots, dt_n]}{(t_0 + \cdots + t_n -1, dt_0 + \cdots + dt_n)}~,
\]
with the obvious simplicial structure.

\begin{theorem}[Dupont's contraction]\label{dupontcontraction}
There exists a simplicial contraction
\[
\begin{tikzcd}[column sep=5pc,row sep=3pc]
\Omega_{\bullet} \arrow[r, shift left=1.1ex, "p_{\bullet}"{name=F}] \arrow[loop left]{l}{h_{\bullet}}
& C_c^*(\Delta^{\bullet})~, \arrow[l, shift left=.75ex, "i_{\bullet}"{name=U}]
\end{tikzcd}
\]
where $C_c^*(\Delta^\bullet)$ denotes the simplicial dg module given by the cellular cochains of the standard simplex. 
\end{theorem}

A basis of $C^*_c(\Delta^n)$ is given by $\{\omega_I\}$ for $I = \{ i_0,\cdots,i_k \} \subset [n]$, where $|\omega_I| = -k$. Since we are working in the homological convention, both algebras are concentrated in degrees $\leq 0$, with $dt_i$ being of degree $-1$. 

\medskip

The data of a simplicial contraction amounts to the data of two morphisms of simplicial cochain complexes
\[
i_{\bullet}: \Omega_{\bullet} \longrightarrow C^*_c(\Delta^{\bullet})~, \quad \text{and} \quad p_{\bullet}: C^*_c(\Delta^{\bullet}) \longrightarrow \Omega_{\bullet}~;
\]
and a degree $1$ linear map $h_{\bullet}: \Omega_{\bullet}\longrightarrow \Omega_{\bullet}$, satisfying the following conditions:
\[
p_n i_n = \mathrm{id}_{C^*_c(\Delta^n)}~, \quad i_n p_n - \mathrm{id}_{\Omega_n} = d_n h_n +h_n d_n~, \quad h_n^2 = 0~, \hspace{1pc} p_n h_n = 0~, \hspace{1pc} h_n i_n = 0~,
\]
where $d_n$ stands for the differential of $\Omega_n$. 

\begin{Remark}
For $n \geq 1$, the dg unital commutative algebra $\Omega_n$ is not augmented, as we have $1 = t_0 + \cdots + t_n$. In order to keep track of the full structure on $\Omega^n$, one has to see it as a dg $\ucom$-algebra and not as a dg $\mathcal{C}om$-algebra.
\end{Remark}

\begin{lemma}\label{lemma: simplicial CC infinity}
There is a simplicial $u\mathcal{CC}_\infty$-algebra structure on $C^*_c(\Delta^{\bullet})$. 
\end{lemma}

\begin{proof}
The dg operad $\Omega\mathrm{B}^{\mathrm{s.a}}\ucom$ is a cofibrant resolution of $\ucom$. Therefore we can apply the homotopy transfer theorem given in \cite[Theorem 6.5.5]{HirshMilles12} using the Dupont contraction of Theorem \ref{dupontcontraction}. Since this contraction is compatible with the simplicial structure, we obtain a simplicial dg $\Omega\mathrm{B}^{\mathrm{s.a}}\ucom$-algebra.
\end{proof}

Let $C^c_*(\Delta^n)$ denote the cellular chains on the $n$-simplex. 

\begin{lemma}\label{lemma: cosimplicial CC infity}
There is a cosimplicial $u\mathcal{CC}_\infty$-coalgebra structure on $C^c_*(\Delta^\bullet)$.
\end{lemma}

\begin{proof}
The linear dual of a degree-wise finite dimensional algebra over an dg operad is naturally a dg coalgebra over the same operad. Since $C^*_c(\Delta^n)$ is finite dimensional degree-wise for all $n$, its linear dual $C^c_*(\Delta^n)$ carries a canonical a $u\mathcal{CC}_\infty$-coalgebra structure. Applying a contravariant functor to a simplicial object gives a cosimplicial object.
\end{proof}

\begin{Remark}
Notice that $C^c_*(\Delta^\bullet)$ is a non-conilpotent $u\mathcal{CC}_\infty$-coalgebra, therefore it cannot be encoded by a cooperad, like stated in \cite[Section 2]{robertnicoud2020higher}.
\end{Remark}

We use this cosimplicial $u\mathcal{CC}_\infty$-coalgebra to construct an adjunction between simplicial sets and $u\mathcal{CC}_\infty$-coalgebras, using the following seminal result.

\begin{theorem}[{\cite{Kan58}}]\label{thm: Kan seminal result}
Let $\mathsf{C}$ be a locally small cocomplete category. The data of an adjunction
\[
\begin{tikzcd}[column sep=7pc,row sep=3pc]
            \mathsf{sSet} \arrow[r, shift left=1.1ex, "L"{name=F}] &\mathscr{C} \arrow[l, shift left=.75ex, "R"{name=U}]
            \arrow[phantom, from=F, to=U, , "\dashv" rotate=-90]
\end{tikzcd}
\]
is equivalent to the data of a cosimplicial object $F: \Delta \longrightarrow \mathscr{C}$ in $\mathscr{C}$.
\end{theorem}

\begin{proof}
Given an adjunction $L \dashv R$, one can pullback the left adjoint $L$ along the Yoneda embedding $\mathrm{Yo}$
\[
\begin{tikzcd}
	\Delta \arrow[r," \mathrm{Yo} "] 
	&\mathsf{sSet} \arrow[r,"L"]
	&\mathsf{sSet}~,
\end{tikzcd}
\]
and obtain a cosimplicial object in $\mathscr{C}$. Given a cosimplicial object $F: \Delta \longrightarrow \mathscr{C}$ in $\mathscr{C}$, one can consider its left Kan extension $\mathrm{Lan}_{\mathrm{Yo}}(F)$ since $\mathscr{C}$ is cocomplete. Its right adjoint is given by 
\[
R(-)_\bullet \coloneqq \mathrm{Hom}_{\mathscr{C}}(F(\Delta^{\bullet}),-)~.
\]
\end{proof}

\begin{Proposition}
There is an adjunction 
\[
\begin{tikzcd}[column sep=7pc,row sep=3pc]
            \mathsf{sSet} \arrow[r, shift left=1.1ex, "C^c_*(-)"{name=F}] 
            &u\mathcal{CC}_\infty\text{-}\mathsf{coalg}~, \arrow[l, shift left=.75ex, "\overline{\mathcal{R}}"{name=U}]
            \arrow[phantom, from=F, to=U, , "\dashv" rotate=-90]
\end{tikzcd}
\]
where $C^c_*(-)$ denotes the cellular chain functor endowed with a canonical $u\mathcal{CC}_\infty$-coalgebra structure. 
\end{Proposition}

\begin{proof}
This is a straightforward application of Theorem \ref{thm: Kan seminal result} to the cosimplicial $u\mathcal{CC}_\infty$-coalgebra constructed in Lemma \ref{lemma: cosimplicial CC infity}. It is immediate to check that the left adjoint is given by the cellular chain functor $C^c_*(-)$ as a dg module, endowed with a $u\mathcal{CC}_\infty$-coalgebra structure. This structure is given by the following colimit

\[
C_*^c(X) \cong \colim_{\mathrm{E}(X)} C_*^c(\Delta^\bullet)~.
\]

on the category of elements $E(X)$ of $X$ in the category of $u\mathcal{CC}_\infty$-coalgebras.
\end{proof}

The right adjoint is therefore given, for a $u\mathcal{CC}_\infty$-coalgebra $C$, by the functor
\[
\overline{\mathcal{R}}(C)_\bullet \coloneqq \mathrm{Hom}_{u\mathcal{CC}_\infty\text{-}\mathsf{cog}}(C^c_*(\Delta^\bullet),C)~.
\]

By pushing forward along the complete Cobar construction $\widehat{\Omega}_\iota$ the cosimplicial $u\mathcal{CC}_\infty$-coalgebra $C^c_*(\Delta^\bullet)$, we also obtain a cosimplicial complete curved absolute $\mathcal{L}_\infty$-algebra.

\begin{Definition}[Maurer-Cartan cosimplicial algebra]
The \textit{Maurer-Cartan cosimplicial algebra} is given by the cosimplicial complete curved absolute $\mathcal{L}_\infty$-algebra
\[
\mathfrak{mc}^{\bullet} \coloneqq \widehat{\Omega}_{\iota}C^c_*(\Delta^{\bullet})~.
\]
\end{Definition}

\begin{Proposition}
There is an adjunction 
\[
\begin{tikzcd}[column sep=7pc,row sep=3pc]
            \mathsf{sSet} \arrow[r, shift left=1.1ex, "\mathcal{L}"{name=F}] 
            &\mathsf{curv}~\mathsf{abs}~\mathcal{L}_\infty\text{-}\mathsf{alg}^{\mathsf{comp}}~. \arrow[l, shift left=.75ex, "\mathcal{R}"{name=U}]
            \arrow[phantom, from=F, to=U, , "\dashv" rotate=-90]
\end{tikzcd}
\]
\end{Proposition}

\begin{proof}
This follows directly from Theorem \ref{thm: Kan seminal result}, applied to the cosimplicial object $\mathfrak{mc}^{\bullet}$.
\end{proof}

\begin{Remark}
The right adjoint is therefore given by, for a curved absolute $\mathcal{L}_\infty$-algebra $\mathfrak{g}$, by  
\[
\mathcal{R}(\mathfrak{g})_\bullet \coloneqq \mathrm{Hom}_{\mathsf{curv}~\mathsf{abs}~\mathcal{L}_\infty\text{-}\mathsf{alg}} \left(\widehat{\Omega}_\iota C^c_*(\Delta^\bullet), \mathfrak{g}\right)~.
\]
\end{Remark}

\begin{theorem}\label{thm: triangle of adjunctions}
The following triangle of adjunctions

\[
\begin{tikzcd}[column sep=5pc,row sep=2.5pc]
&\hspace{1pc}u\mathcal{CC}_\infty \textsf{-}\mathsf{coalg} \arrow[dd, shift left=1.1ex, "\widehat{\Omega}_{\iota}"{name=F}] \arrow[ld, shift left=.75ex, "\overline{\mathcal{R}}"{name=C}]\\
\mathsf{sSet}  \arrow[ru, shift left=1.5ex, "C^c_*(-)"{name=A}]  \arrow[rd, shift left=1ex, "\mathcal{L}"{name=B}] \arrow[phantom, from=A, to=C, , "\dashv" rotate=-70]
& \\
&\hspace{3pc}\mathsf{curv}~\mathsf{abs}\text{-}\mathcal{L}_\infty\textsf{-}\mathsf{alg}^{\mathsf{comp}}~, \arrow[uu, shift left=.75ex, "\widehat{\mathrm{B}}_{\iota}"{name=U}] \arrow[lu, shift left=.75ex, "\mathcal{R}"{name=D}] \arrow[phantom, from=B, to=D, , "\dashv" rotate=-110] \arrow[phantom, from=F, to=U, , "\dashv" rotate=-180]
\end{tikzcd}
\]

commutes. Furthermore, all the adjunctions are Quillen adjunctions when we consider the model category structures on $u\mathcal{CC}_\infty$-coalgebras and on complete curved absolute $\mathcal{L}_\infty$-algebras of Theorem \ref{thm: Equivalence the Quillen CC infini et absolues}; and the Kan-Quillen model category structure on simplicial sets.
\end{theorem}

\begin{proof}
Let $\mathfrak{g}$ be a curved absolute $\mathcal{L}_\infty$-algebra. We have that
\[
\mathcal{R}(\mathfrak{g})_\bullet = \mathrm{Hom}_{\mathsf{curv}~\mathsf{abs}~\mathcal{L}_\infty\text{-}\mathsf{alg}}\left(\widehat{\Omega}_\iota C^c_*(\Delta^\bullet), \mathfrak{g} \right) \cong \mathrm{Hom}_{u\mathcal{CC}_\infty\text{-}\mathsf{cog}}\left( C^c_*(\Delta^\bullet),\widehat{\mathrm{B}}_\iota \mathfrak{g} \right) \cong \overline{\mathcal{R}} \left(\widehat{\mathrm{B}}_\iota \mathfrak{g} \right)_\bullet~.
\]

The functor $C^c_*(-)$ sends monomorphisms of simplicial sets to degree-wise monomorphisms of $u\mathcal{CC}_\infty$-coalgebras and it sends weak homotopy equivalences to quasi-isomorphisms. Hence the triangle is made up of Quillen adjunctions.
\end{proof}

\subsection{The integration functor} The functor $\mathcal{R}$ in the above triangle is \textit{the integration functor} we were looking for. Before giving an explicit combinatorial description of it, let us first state some of its fundamental properties.

\begin{theorem}\label{thm: propriétés de l'intégration}\leavevmode

\begin{enumerate}
\item For any curved absolute $\mathcal{L}_\infty$-algebra $\mathfrak{g}$, the simplicial set $\mathcal{R}(\mathfrak{g})$ is a Kan complex.

\medskip

\item Let $f: \mathfrak{g} \twoheadrightarrow \mathfrak{h}$ be a degree-wise epimorphism of curved absolute $\mathcal{L}_\infty$-algebras. Then 
\[
\mathcal{R}(f): \mathcal{R}(\mathfrak{g}) \twoheadrightarrow \mathcal{R}(\mathfrak{h})
\]
is a fibrations of simplicial sets. 

\medskip

\item The functor $\mathcal{R}$ preserves weak equivalences. In particular, it sends any graded quasi-isomorphism $f: \mathfrak{g} \qi \mathfrak{h}$ between complete curved absolute $\mathcal{L}_\infty$-algebras to a weak homotopy equivalence of simplicial sets.

\end{enumerate}
\end{theorem}

\begin{proof}
The functor $\mathcal{R}$ is a right Quillen functor. 
\end{proof}

\begin{Remark}
Let us contextualize the results of Theorem \ref{thm: propriétés de l'intégration}.

\begin{enumerate}
\item The first result is quintessential to the classical integration theory of dg Lie algebras or nilpotent $\mathcal{L}_\infty$-algebras as developed in \cite{Hinich01} and in \cite{Getzler09}.

\vspace{0.5pc}

\item The second result generalizes one of the main theorems of \cite{Getzler09}, see Theorem 5.8 in \textit{loc.cit}.

\vspace{0.5pc}

\item The third point, the \textit{homotopy invariance} of the integration functor, is the generalization of Goldman--Milson's invariance theorem proved in \cite{goldmanmillson} and its generalization to $\mathcal{L}_\infty$-algebras given in \cite{dolgushevrogers}. 
\end{enumerate}
\end{Remark}

\begin{Proposition}
Let $\mathfrak{g}$ be a complete curved absolute $\mathcal{L}_\infty$-algebra. There is an isomorphism of simplicial sets
\[
\mathcal{R}(\mathfrak{g})_\bullet \cong \lim_{\omega} \mathcal{R}(\mathfrak{g}/\mathrm{W}_\omega \mathfrak{g})_\bullet~.
\]
\end{Proposition}

\begin{proof}
Since $\mathfrak{g}$ is complete, it can be written as 
\[
\mathfrak{g} \cong \lim_{\omega} \mathfrak{g}/\mathrm{W}_\omega \mathfrak{g}~,
\]
where the limit is taken in the category of curved absolute $\mathcal{L}_\infty$-algebras. Since $\mathcal{R}$ is right adjoint, it preserves all limits.
\end{proof}

\textbf{Dold-Kan correspondence.} The adjunction 
\[
\begin{tikzcd}[column sep=5pc,row sep=2.5pc]
\mathsf{sSet}  \arrow[r, shift left=1ex, "\mathcal{L}"{name=B}] 
&\mathsf{curv}~\mathsf{abs}\text{-}\mathcal{L}_\infty\textsf{-}\mathsf{alg}^{\mathsf{comp}} \arrow[l, shift left=.75ex, "\mathcal{R}"{name=D}] \arrow[phantom, from=B, to=D, , "\dashv" rotate=-90] 
\end{tikzcd}
\]
is a \textit{generalization of the Dold-Kan correspondence}. Notice that a chain complex $(V,d_V)$ is a particular example of a complete curved absolute $\mathcal{L}_\infty$-algebra where the structural morphism 
\[
\gamma_V: \displaystyle \prod_{n \geq 0} \widehat{\Omega}^{\mathrm{s.c}}\ucom^*(n) ~ \widehat{\otimes}_{\mathbb{S}_n} ~ V^{\otimes n} \longrightarrow V
\]
is the \textit{zero morphism}. We call them \textit{abelian} curved absolute $\mathcal{L}_\infty$-algebras.

\begin{Proposition}
Let $(V,d_V)$ be an abelian curved absolute $\mathcal{L}_\infty$-algebra. Then
\[
\mathcal{R}(V) \cong \Gamma(V)~,
\]
where $\Gamma(-)$ is the Dold-Kan functor. 
\end{Proposition}

\begin{proof}
Recall that 
\[
\mathcal{R}(V)_\bullet \cong \mathrm{Hom}_{u\mathcal{CC}_\infty\text{-}\mathsf{cog}}(C^c(\Delta^\bullet),\widehat{\mathrm{B}}_\iota V )~.
\]
Since $\gamma_V = 0$, the complete Bar construction $\widehat{\mathrm{B}}_\iota V $ is simply given by the cofree $u\mathcal{CC}_\infty$-coalgebra generated by the dg module $(V,d_V)$. Thus 

\[
\mathcal{R}(V)_\bullet \cong \mathrm{Hom}_{u\mathcal{CC}_\infty\text{-}\mathsf{cog}}(C^c_*(\Delta^\bullet),\mathscr{C}(u\mathcal{CC}_\infty)(V)) \cong \mathrm{Hom}_{\mathsf{dg}\text{-}\mathsf{mod}}(C^c_*(\Delta^\bullet),V) \cong \Gamma(V)~.
\]
\end{proof}

\begin{Remark}
Abelian curved absolute $\mathcal{L}_\infty$-algebras concentrated in positive degrees are models for simplicial abelian groups. One can think of general curved absolute $\mathcal{L}_\infty$-algebras concentrated in positive degrees as models for a non-commutative analogue of simplicial abelian groups. See Subsection \ref{subsection: higher BCH} for more on this.
\end{Remark}

\textbf{Derived adjunction.} The previous adjunction computes the \textit{derived functors} of the adjunction
\[
\begin{tikzcd}[column sep=5pc,row sep=2.5pc]
\mathsf{sSet}  \arrow[r, shift left=1.5ex, "C^c_{*}(-)"{name=A}]
&u\mathcal{CC}_\infty \textsf{-}\mathsf{coalg}~~. \arrow[l, shift left=.75ex, "\overline{\mathcal{R}}"{name=C}] \arrow[phantom, from=A, to=C, , "\dashv" rotate=-90]
\end{tikzcd}
\]
Let $C$ be a $u\mathcal{CC}_\infty$-coalgebra. It might not be fibrant, therefore one needs to take a fibrant resolution of $C$ in order to compute the right derived functor $\mathbb{R}\overline{\mathcal{R}}$ at $C$. The Theorem \ref{thm: Equivalence the Quillen CC infini et absolues} provides us with 
\[
\eta_C: C \qi \widehat{\mathrm{B}}_\iota \widehat{\Omega}_\iota C ~,
\]
a functorial fibrant resolution of $C$. Therefore we have that

\[
\mathbb{R}\overline{\mathcal{R}}(C)_\bullet = \mathrm{Hom}_{u\mathcal{CC}_\infty\text{-}\mathsf{cog}}(C^c(\Delta^\bullet),\widehat{\mathrm{B}}_\iota \widehat{\Omega}_\iota  C ) \cong \mathrm{Hom}_{\mathsf{curv}~\mathsf{abs}~\mathcal{L}_\infty\textsf{-}\mathsf{alg}}(\widehat{\Omega}_\iota C^c(\Delta^\bullet) , \widehat{\Omega}_\iota C ) \cong \mathcal{R}(\widehat{\Omega}_\iota  C )~.
\]

\begin{Remark}[$\infty$-morphisms]
Let $C$ and $D$ be two $u\mathcal{CC}_\infty$-coalgebras. The hom-set 

\[
\mathrm{Hom}_{\mathsf{curv}~\mathsf{abs}~\mathcal{L}_\infty\textsf{-}\mathsf{alg}}(\widehat{\Omega}_\iota C , \widehat{\Omega}_\iota  D ) 
\]

is the set of $\infty$-morphisms of $u\mathcal{CC}_\infty$-coalgebras between $C$ and $D$. See Sections \ref{Section: Complete Bar-Cobar} and \ref{Section: Complete Bar-Cobar curved}. 
\end{Remark}

\subsection{Explicit version of the integration functor} Let us now give a combinatorial description of the integration functor. Recall that the integration functor $\mathcal{R}$ is given by 
\[
\mathcal{R}(-)_\bullet \coloneqq \mathrm{Hom}_{\mathsf{curv}~\mathsf{abs}~\mathcal{L}_\infty\textsf{-}\mathsf{alg}}(\widehat{\Omega}_\iota(C^c_* \Delta^\bullet),-) \cong \mathrm{Hom}_{\mathsf{curv}~\mathsf{abs}~\mathcal{L}_\infty\textsf{-}\mathsf{alg}}(\mathfrak{mc}^\bullet,-)~.
\]
Thus it is primordial to describe the complete curved absolute $\mathcal{L}_\infty$-algebra structure on the cosimplicial Maurer-Cartan algebra $\mathfrak{mc}^\bullet$ which represents this functor. Since it is given by the complete Cobar construction of the cosimplicial $u\mathcal{CC}_\infty$-coalgebra $C^c_*(\Delta^\bullet)$, let us describe this structure.

\begin{Notation}
Consider $\Omega_n$ with its dg $u\mathcal{C}om$-algebra structure. Let
\[
\mu_k: \ucom(k) \otimes_{\mathbb{S}_k} \Omega_n^{\otimes k} \longrightarrow \Omega_n
\]
denote its dg $\ucom$-algebra structural morphisms for $k \geq 1$, and $u: \kk \longrightarrow \Omega_n$ its unit. 
\end{Notation}

\begin{Proposition}\label{prop: formulas for the HTT}
The transferred $u\mathcal{CC}_\infty$-algebra structure on $C_c^*(\Delta^n)$ via the Dupont contraction of Theorem \ref{dupontcontraction} is determined by morphisms
\[
\Big\{ \mu_\tau: (C_c^*(\Delta^n))^{\otimes m} \longrightarrow C_c^*(\Delta^n) \Big\}
\]
of degree $-\omega + 1$, where $\tau$ is a rooted tree in $\mathrm{RT}_m^\omega$ of arity $m$ and weight $\omega$ with no corks. 

\medskip

For $\tau$ in $\mathrm{RT}_m$, the operation $\mu_\tau$ is given by labeling all the vertices of $\tau$ with an operation $\mu_k$ of $\Omega_n$, where $k$ is the number of inputs of the vertex, by labeling all internal edges with the homotopy $h_n$, all the leaves with the map $i_n$ and the root with the map $p_n$, and composing the labeling operations along the rooted tree $\tau$. Pictorially, it is given by 

\begin{center}
\includegraphics[width=115mm,scale=1.15]{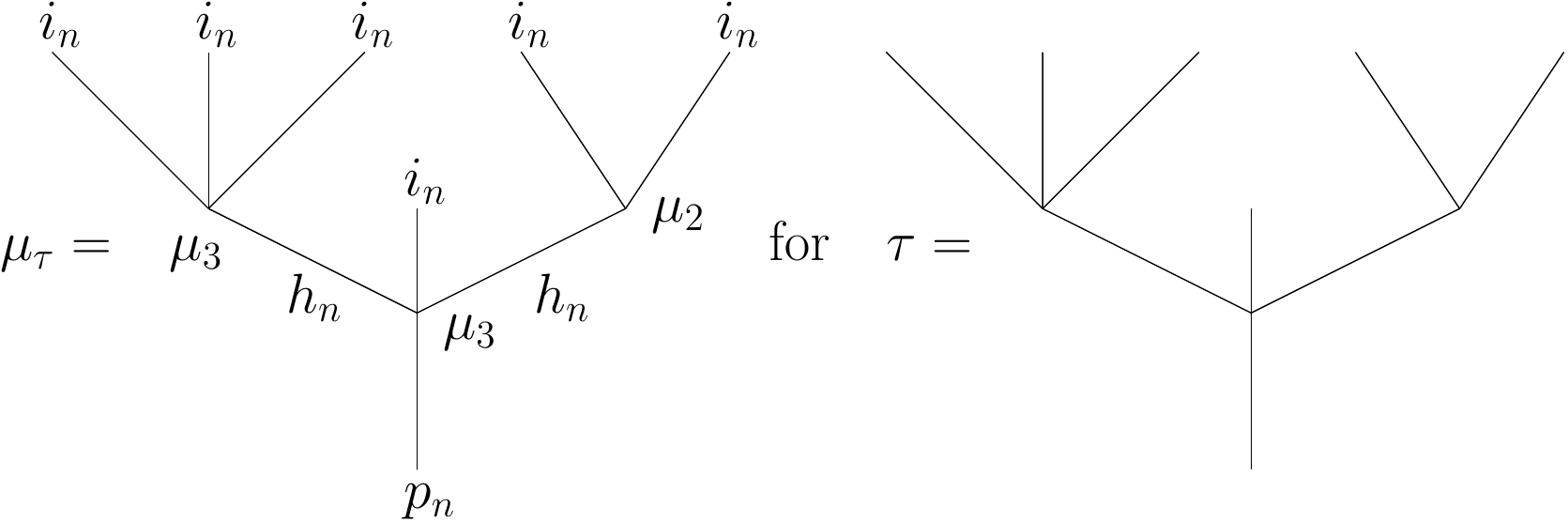}.
\end{center}

The only corked rooted tree acting non-trivially is the single cork morphism given by $p_n \circ u: \kk \longrightarrow C_c^*(\Delta^n) $.
\end{Proposition}

\begin{proof}
The standard transfer formula for the homotopy transfer theorem given in \cite[Theorem 6.5.5]{HirshMilles12} involves operations $\mu_\tau$ where $\tau$ runs over all corked rooted trees in $\mathrm{CRT}$. Two different types of corks can appear: the corks given by $u$ and the corks given by $u \circ h_n$. Pictorially, $\tau$ can be as follows

\begin{center}
\includegraphics[width=75mm,scale=1]{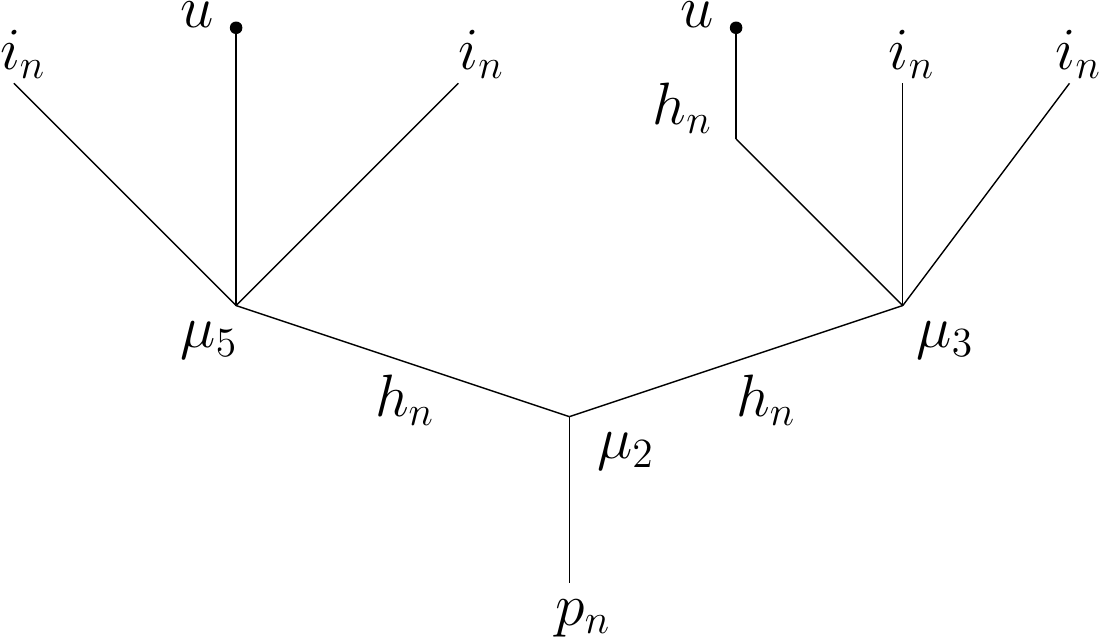},
\end{center}

where the first type of corks is represented on the left branch of the tree and the second type is represented on the right branch of the tree. Since we started with a strict $\ucom$-algebra structure on $\Omega_n$, the composition $\mu_k \circ_i u$ simplifies into $\mu_{k-1}$, hence the first type of corks disappears. Recall that $\Omega_n$ is concentrated in non-positive degrees and that the unit of $\Omega_n$ lies in degree $0$. Hence $h_n \circ u =0$ as $h_n$ raises the degree by one and the second type of corks disappears as well. We are left with operations $\mu_\tau$ where $\tau$ runs over all rooted trees without corks $\mathrm{RT}_m^\omega~,$ with a single operation involving the unit $u$ given by $p_n \circ u$.
\end{proof}

\begin{Remark}
We recover exactly the same structure on $C^*_c(\Delta^n)$ as the one given in \cite[Proposition 2.9]{robertnicoud2020higher} with the addition of the arity $0$ operation $p_n \circ u$.
\end{Remark}

Let us compute the $u\mathcal{CC}_\infty$-coalgebra structure on its linear dual $C^c_*(\Delta^{n})$. The elementary decompositions 
\[
\Big\{ \Delta_\tau: C^c_*(\Delta^{n}) \longrightarrow C^c_*(\Delta^{n})^{\otimes m} \Big\}
\]
which determine the $u\mathcal{CC}_\infty$-coalgebra structure of $C^c_*(\Delta^{n})$ are given by $\Delta_\tau \coloneqq (\mu_\tau)^*$. Let $\{a_I\} = \{\omega_I^*\}$ be the dual basis of $C^c_*(\Delta^{n})$. Let $\tau$ be a rooted tree in $\mathrm{RT}_m^\omega$. For any $m$-tuple $(I_1,\cdots,I_m)$ of non-trivial subsets of $[n]$, suppose we have that
\[
\mu_\tau(\omega_{I_1},\cdots,\omega_{I_m}) = \sum_{I \subseteq [n]} \lambda_I^{\tau,(I_1,\cdots,I_m)}\omega_I~.
\]
where $\{\omega_I\}$ is the basis of $C^*_c(\Delta^{n})$. Then
\[
\Delta_\tau(a_I) = \sum_{\substack{I_1,\cdots,I_m \subseteq [n],I_l \neq \emptyset \\ \lambda_I^{\tau,(I_1,\cdots,I_m)} \neq 0}} \frac{1}{\lambda_I^{\tau,(I_1,\cdots,I_m)}} a_{I_1} \otimes \cdots \otimes a_{I_m}~. 
\]
Notice that $a_{I_1} \otimes \cdots \otimes a_{I_m}$ appears in the decomposition of $a_I$ under $\Delta_\tau$ if and only if $\omega_I$ appears with a non-zero coefficient in $\mu_\tau(\omega_{I_1},\cdots,\omega_{I_m})$. 

\begin{Remark}
Computing the coefficients $\lambda_I^{\tau,(I_1,\cdots,I_m)}$ is a hard task, see \cite[Remark 5.11]{robertnicoud2020higher}.
\end{Remark}

\begin{lemma}\label{lemma: description de C(Delta) en uCC cogebre}
For all $n \geq 0$, the $u\mathcal{CC}_\infty$-coalgebra $C^c(\Delta^{n})$ is a strictly counital $u\mathcal{CC}_\infty$-coalgebra. 
\end{lemma}

\begin{proof}
This is immediate from Proposition \ref{prop: formulas for the HTT}.
\end{proof}

\begin{Proposition}\label{prop: C(X) is strictly counital}
Let $X$ be a simplicial set, then $C_*^c(X)$ is a strictly counital $u\mathcal{CC}_\infty$-coalgebra.
\end{Proposition}

\begin{proof}
By Lemma \ref{lemma: description de C(Delta) en uCC cogebre}, we know that this is true for the cellular chains of the standard $n$-simplex $C_*^c(\Delta^\bullet)$. Let $X$ be a simplicial set, it can be written as the following colimit
\[
X \cong \colim_{\mathrm{E}(X)} \Delta^\bullet~,
\]
indexed by $\mathrm{E}(X)$, the category of elements of $X$. Consequently, we have that 
\[
C_*^c(X) \cong \colim_{\mathrm{E}(X)} C_*^c(\Delta^\bullet)~.
\]
By Proposition \ref{prop: coreflective sub-category}, we know that strictly counital $u\mathcal{CC}_\infty$-coalgebras form a coreflective sub-category of $u\mathcal{CC}_\infty$-coalgebras, therefore they are closed under colimits. Hence $C_*^c(X)$ is a strictly counital $u\mathcal{CC}_\infty$-coalgebra.
\end{proof}

\textbf{Complete Cobar construction.} Let $(C,\Delta_C,d_C)$ be a $u\mathcal{CC}_\infty$-coalgebra, where 
\[
\Delta_C: C \longrightarrow \displaystyle \prod_{n \geq 0} \mathrm{Hom}_{\mathbb{S}_n}\left(\Omega\mathrm{B}^{\mathrm{s.a}}\ucom(n), C^{\otimes n}\right)
\]
denotes its structural morphism. The complete Cobar construction $\widehat{\Omega}_{\iota} C $ is given by the underlying graded module 
\[
\prod_{n \geq 0} \widehat{\Omega}^{\mathrm{s.c}}\ucom^*(n) ~ \widehat{\otimes}_{\mathbb{S}_n} ~ C^{\otimes n}~.
\]
It is endowed with a pre-differential $d_{\mathrm{cobar}}$ given by the difference of $d_1$ and $d_2$. The first term $d_1$ is given by 
\[
d_1= - \widehat{\mathscr{S}}^c(d_{\widehat{\Omega}^{\mathrm{s.c}}\ucom^*})(\mathrm{id}) + \widehat{\mathscr{S}}^c(\mathrm{id})(\diracComb(\mathrm{id},d_C))~.
\]

The second term $d_2$ is induced by the $u\mathcal{CC}_\infty$-coalgebra structure of $C$, it is the unique derivation extending the map $\varphi$, given by the following composites
\[
\begin{tikzcd}[column sep=4pc,row sep=3pc]
C \arrow[d,"\Delta_C",swap] \arrow[r,"\varphi"]
&\displaystyle \prod_{n \geq 0} \widehat{\Omega}^{\mathrm{s.c}}\ucom^*(n) ~ \widehat{\otimes}_{\mathbb{S}_n} ~ C^{\otimes n} \\
\displaystyle \prod_{n \geq 0} \mathrm{Hom}_{\mathbb{S}_n}\left(\Omega\mathrm{B}^{\mathrm{s.a}}\ucom(n), C^{\otimes n}\right)  \arrow[r,"\widehat{\mathscr{S}}^c(\iota)(\mathrm{id})"] 
&\displaystyle \prod_{n \geq 0} \mathrm{Hom}_{\mathbb{S}_n}\left(\mathrm{B}^{\mathrm{s.a}}\ucom(n), C^{\otimes n}\right)~. \arrow[u,"\cong"]
\end{tikzcd}
\]

\begin{theorem}\label{theom: formulas for the pre differentials of mc}
Let $I = \{i_0,\cdots,i_k\} \subseteq [n]$ with $k > 0$, and $a_I$ be the corresponding basis element in $C^c_*(\Delta^n)$. When $\mathrm{Card}(I) \geq 2$, the image of $a_I$ under the pre-differential $d_{\mathrm{cobar}}$ is given by 
\[
d_{\mathrm{cobar}}(a_I) = \sum_{l=0}^k (-1)^l a_{i_0\cdots\widehat{i_l}\cdots i_k} - \sum_{m \geq 2}\sum_{\substack{\tau \in \mathrm{RT}_m \\ I_1,\cdots,I_m \subseteq [n], I_l \neq \emptyset \\ \lambda_I^{\tau,(I_1,\cdots,I_m)} \neq 0}}\frac{1}{\mathcal{E}(\tau)\lambda_I^{\tau,(I_1,\cdots,I_m)}} \tau(a_{I_1},\cdots, a_{I_m})~,
\]

for $I = \{i_0,\cdots,i_k\}$, where $\mathcal{E}(\tau)$ is the renormalization coefficient of Remark \ref{Remark: renormalization}. When $\mathrm{Card}(I) = 1$, the image of $a_I$ under the pre-differential is given by 
\[
-\sum_{n \geq 0,~ n \neq 1} \frac{1}{n!} c_n(a_0,\cdots,a_0)~,
\]
where $c_n$ denotes the $n$-corolla. 
\end{theorem}

\begin{proof}
By Proposition \ref{prop: formulas for the HTT}, the only arity $0$ operation on $C^c_*(\Delta^n)$ is given by $p_n \circ u$. In order to know in which decompositions the lonely cork is going to appear, it is enough to compute the image of the unit $1$ in  $\Omega_n$ by the morphism $p_n$. Since
\[
p_n(1) = t_0 + \cdots + t_n \in C^*_c(\Delta^n)
\]
for all $n \geq 0$, the cork only appears in the decomposition of the elements of the form $a_{i}$ where $\{i\} \subseteq [n]$. Consequently, all the formulas for $d_{\mathrm{cobar}}(a_I)$ where $\mathrm{Card}(I) \geq 2$ are only indexed by rooted trees without corks. The renormalization coefficient $\mathcal{E}(\tau)$ appears because an identification between invariants and coinvariants was done in the definition of $d_2$.
\end{proof}

\begin{Corollary}\label{cor: same formulas as for RNV}
When $\mathrm{Card}(I) \geq 2$, the formula for $d_{\mathrm{cobar}}(a_I)$ coincides with the formula computed in \cite[Section 2.2]{robertnicoud2020higher} in the non-curved case.
\end{Corollary}

\begin{proof}
Follows directly from the above theorem.
\end{proof}

\begin{Corollary}
Let $(\mathfrak{g},\gamma_\mathfrak{g},d_\mathfrak{g})$ be a curved absolute $\mathcal{L}_\infty$-algebra. There is a bijection
\[
\mathcal{R}(\mathfrak{g})_0 \cong \mathcal{MC}(\mathfrak{g})~,
\]
between the $0$-simplicies of $\mathcal{R}(\mathfrak{g})$ and the Maurer-Cartan elements of $\mathfrak{g}$.
\end{Corollary}

\begin{proof}
At $n=0$, the Dupont contraction is trivial since $h_0 = 0$ and $i_0$ and $p_0$ are isomorphism. Thus the resulting $u\mathcal{CC}_\infty$-coalgebra structure on $\kk.a_0$ is simply given by 
\[
\begin{tikzcd}[column sep=4pc,row sep=-0.8pc]
\kk.a_0 \arrow[r]
&\displaystyle \prod_{n \geq 0} \mathrm{Hom}_{\mathbb{S}_n}\left(\Omega\mathrm{B}^{\mathrm{s.a}}\ucom(n), (\kk.a_0)^{\otimes n}\right) \\
a_0 \arrow[r,mapsto]
&\displaystyle \sum_{n \geq 0,~ n \neq 1} \left[\mathrm{ev}_{a_0}: c_n \mapsto a_0 \otimes \cdots \otimes a_0 \right]~.
\end{tikzcd}
\]
Thus 
\[
\mathfrak{mc}^0 \cong \left(\displaystyle \prod_{n \geq 0} \widehat{\Omega}^{\mathrm{s.c}}\ucom^*(n) ~ \widehat{\otimes}_{\mathbb{S}_n} ~ (\kk.a_0)^{\otimes n}, d_{\mathrm{cobar}}(a_0) = - \sum_{n \geq 0,~ n \neq 1} \frac{1}{n!} c_n(a_0,\cdots,a_0) \right)~.
\]

The data of a morphism of curved absolute $\mathcal{L}_\infty$-algebras $f: \mathfrak{mc}^0 \longrightarrow \mathfrak{g}$ is equivalent to the data of $\alpha$ in $\mathfrak{g}$ such that 
\[
d_\mathfrak{g}(\alpha) = - \gamma_\mathfrak{g}\left(\sum_{n \geq 0,~ n \neq 1} \frac{1}{n!} c_n(\alpha,\cdots,\alpha)\right)~,
\]
since $\mathfrak{mc}^0$ is freely generated by $a_0$ as a pdg $\mathrm{B}^{\mathrm{s.a}}\ucom$-algebra, and since the data of a morphism of pdg $\mathrm{B}^{\mathrm{s.a}}\ucom$-algebras is the same as the data of a morphism of curved $\mathrm{B}^{\mathrm{s.a}}\ucom$-algebras.
\end{proof}

\begin{Remark}
Let $\mathfrak{g}$ be a curved absolute $\mathcal{L}_\infty$-algebra, then
\[
\mathcal{MC}(\mathfrak{g}) \cong \mathrm{Hom}_{u\mathcal{CC}_\infty\text{-}\mathsf{cog}}(\kk,\widehat{\mathrm{B}}_\iota \mathfrak{g})~.
\]
Therefore the Maurer-Cartan elements of $\mathfrak{g}$ are in a one-to-one correspondence with \textit{group-like elements} of the complete Bar construction of $\mathfrak{g}$. 
\end{Remark}

\subsection{Gauge actions}
In this subsection, we define intrinsically the gauge action of an element $\lambda$ of degree $1$ on the set of Maurer-Cartan elements. For a complete curved absolute $\mathcal{L}_\infty$-algebra $\mathfrak{g}$, gauge actions correspond to paths in $\mathcal{R}(\mathfrak{g})$, thus we obtain an explicit description of $\pi_0(\mathcal{R}(\mathfrak{g}))$. 

\medskip

Let $\mathfrak{g}$ be a complete curved absolute $\mathcal{L}_\infty$-algebra. We can consider the following continuous function 
\[
\varphi_{\mathcal{MC}}(-) \coloneqq \left[d_\mathfrak{g}(-) + \gamma_\mathfrak{g}\left(\sum_{n \geq 0,~ n \neq 1} \frac{1}{n!} c_n(-,\cdots,-)\right)\right]: \mathfrak{g}_0 \longrightarrow \mathfrak{g}_{-1}
\]
between two complete vector spaces. The Maurer-Cartan set of $\mathfrak{g}$ is given by the $0$ locus in $\mathfrak{g}_{-1}$. Would this function be smooth or algebraic, the resulting set would be a differential or an algebraic variety. One can still formally consider its "tangent space" at a given point $\alpha$ in $\mathcal{MC}(\mathfrak{g})$, which in this case is given by the kernel of $d^\alpha_\mathfrak{g}$, where $d^\alpha_\mathfrak{g}$ is the \textit{twisted differential} by $\alpha$ given by 
\[
d^\alpha_\mathfrak{g}(-) \coloneqq d_\mathfrak{g}(-) + \gamma_\mathfrak{g}\left(\sum_{n \geq 2} \frac{1}{(n-1)!} c_n(\alpha,\cdots,\alpha,-)\right)~.
\]

\begin{lemma}
Let $\mathfrak{g}$ be a curved absolute $\mathcal{L}_\infty$-algebra and let $\alpha$ be a Maurer-Cartan element. Then $d^\alpha_\mathfrak{g}$ squares to zero.
\end{lemma}

\begin{proof}
The formal algebraic computations of \cite[Corollary 5.1]{DSV18} extend \textit{mutatis mutandis} to this case using the associativity Condition \ref{associativity condition}. 
\end{proof}

Therefore an easy way to produce elements in the tangent space $\mathrm{T}_\alpha \mathcal{MC}(\mathfrak{g}) = \mathrm{Ker}(d^\alpha_\mathfrak{g})$ is to take a degree one element $\lambda$ in $\mathfrak{g}$ and consider $d^\alpha_\mathfrak{g}(\lambda)$. Then, given this tangent vector at $\alpha$, one would like to "integrate it" inside $\mathcal{MC}(\mathfrak{g})$. As explained in \cite{DSV18}, these heuristics motivate the following definitions.

\begin{lemma}
Let $(\mathfrak{g},\gamma_\mathfrak{g},d_\mathfrak{g})$ be a complete curved absolute $\mathcal{L}_\infty$-algebra. The graded module  
\[
\kk[[t]] \widehat{\otimes} \mathfrak{g} \coloneqq \lim_\omega \left( \kk[t]/(t^\omega) \otimes \mathfrak{g}/\mathrm{W}_\omega \mathfrak{g} \right) \cong \Bigg\{ \sum_{\omega \geq 0} g_\omega \otimes t^\omega ~~| ~~ g_\omega \in \mathrm{W}_\omega \mathfrak{g} \Bigg\}
\]
has a complete curved absolute $\mathcal{L}_\infty$-algebra structure given by
\[
\gamma_{\kk[[t]] \widehat{\otimes} \mathfrak{g}} \left(\sum_{k \geq 0} \sum_{\omega \geq 0} \sum_{\tau \in \mathrm{CRT}_k^\omega} \lambda_\tau \tau\left(\sum_{i_1 \geq 0}g_{i_1} \otimes t^{i_1}, \cdots, \sum_{i_k \geq 0}g_{i_k} \otimes t^{i_k}\right)\right) \coloneqq 
\]
\[
 \sum_{n \geq 0} \sum_{k \geq 0} \sum_{\omega \geq 0}  \left( \sum_{i_1 + \cdots + i_k = n} \gamma_\mathfrak{g} \left(\sum_{\tau \in \mathrm{CRT}_n^\omega} \lambda_\tau \tau(g_{i_1}, \cdots, g_{i_k}) \right)\right) \otimes t^{n + \omega}~. 
\]
\end{lemma}

\begin{proof}
One directly and easily can check that this structure satisfies the axioms. 
\end{proof}

\begin{Remark}
This structure comes from the convolution complete curved absolute $\mathcal{L}_\infty$-algebra structure structure on 
\[
\mathrm{Hom}_{\mathsf{gr}\textsf{-}\mathsf{mod}}(\kk^c[t],\mathfrak{g})~,
\]
where $\kk^c[t]$ is the conilpotent cofree cocommutative coalgebra cogenerated by $t$ in degree $0$.
\end{Remark}

The algebraic way to integrate a degree one element $\lambda$ in $\mathfrak{g}$ amounts to solving the following ordinary differential equation

\begin{equation}\label{edo jauges}
\frac{\mathrm{d}\gamma(t)}{\mathrm{d}t} = d^{\gamma(t)}_\mathfrak{g}(\lambda)~,
\end{equation}

where $\gamma(t)$ is an element of $\kk[[t]]~ \widehat{\otimes}~ \mathfrak{g}$. Here $d^{\gamma(t)}_\mathfrak{g}$ is again given by 
\[
d^{\gamma(t)}_\mathfrak{g}(\lambda) \coloneqq  d_\mathfrak{g}(\lambda) + \gamma_{\kk[[t]] \widehat{\otimes} \mathfrak{g}}\left(\sum_{n \geq 2}  \frac{1}{(n-1)!} c_k \left(\gamma(t),\cdots,\gamma(t),\lambda\right) \right)~.
\]
where the term $\gamma(t)$ appears $n-1$ times on the right. The solution obviously depends on the initial value $\gamma(0)$. Let us denote it by $\alpha \coloneqq \gamma(0)$.

\begin{lemma}
Let $\mathfrak{g}$ be a curved absolute $\mathcal{L}_\infty$-algebra. Let $\lambda$ be in $\mathfrak{g}_{1}$. For every $\alpha$ in $\mathfrak{g}_0$, the above differential equation has an unique solution, given by: 
\[
\gamma(t) = \sum_{\omega \geq 0} \gamma_\mathfrak{g}\left(\sum_{n \geq 0} \sum_{\tau \in \mathrm{RT}^\omega_n} \frac{|\mathrm{Aut}(\tau)|}{C(\tau)} \tau^\lambda(\alpha, \cdots, \alpha) \right) \otimes t^\omega~.
\]
\end{lemma} 

\begin{proof}
We can apply \cite[Proposition A.11]{robertnicoud2020higher} to this context. This sum is here indexed by rooted trees instead of planar trees, and thus rooted trees automorphisms are taken into account, see \cite[Proposition 1.10]{robertnicoud2020higher} for the planar formula.
\end{proof}

\begin{Definition}[Gauge actions]
Let $\mathfrak{g}$ be a complete curved absolute $\mathcal{L}_\infty$-algebra. Let $\lambda$ be in $\mathfrak{g}_{1}$ and let $\alpha$ be a Maurer-Cartan element of $\mathfrak{g}$. The \textit{gauge action} of $\lambda$ on $\alpha$ is given by the element
\[
\lambda \bullet \alpha \coloneqq \gamma(1) = \sum_{\omega \geq 0} \gamma_\mathfrak{g}\left(\sum_{n \geq 0} \sum_{\tau \in \mathrm{RT}^\omega_n} \frac{|\mathrm{Aut}(\tau)|}{C(\tau)} \tau^\lambda(\alpha, \cdots, \alpha) \right)~,
\] 
in $\mathfrak{g}_0$, where $\gamma(t)$ is the unique solution to Equation \ref{edo jauges} having $\alpha$ as initial condition.
\end{Definition}

\begin{Definition}[Gauge equivalence]
Two Maurer-Cartan elements $\alpha$ and $\beta$ of $\mathfrak{g}$ are \textit{gauge equivalent} if there exists $\lambda$ in $\mathfrak{g}_1$ such that $\lambda \bullet \alpha = \beta$. 
\end{Definition}

\begin{Proposition}
The gauge equivalence relation on the set of Maurer-Cartan elements $\mathcal{MC}(\mathfrak{g})$ defines an equivalence relation. 
\end{Proposition}

\begin{proof}
Follows from Proposition \ref{prop: equivalence between gauge and paths} using the fact that paths equivalence in $\mathcal{R}(\mathfrak{g})$ is an equivalence relation on $\mathcal{MC}(\mathfrak{g})$ since it is a Kan complex.
\end{proof}

\begin{Proposition}\label{prop: equivalence between gauge and paths}
Let $\mathfrak{g}$ be a complete curved absolute $\mathcal{L}_\infty$-algebra and let $\alpha$ and $\beta$ be two Maurer-Cartan elements of $\mathfrak{g}$. The data of a gauge equivalence $\lambda$ between $\alpha$ and $\beta$ is equivalent to the data of a morphism of complete curved absolute $\mathcal{L}_\infty$-algebras 
\[
\varphi: \widehat{\Omega}_\iota(C_*^c(\Delta^1)) \longrightarrow \mathfrak{g}
\]
such that $\varphi(a_{01})= \lambda$, $\varphi(a_0) = \alpha$ and $\varphi(a_1) = \beta$. 
\end{Proposition}

\begin{proof}
A morphism of complete curved absolute $\mathcal{L}_\infty$-algebras
\[
\varphi: \widehat{\Omega}_\iota C_*^c(\Delta^1)  \longrightarrow \mathfrak{g}
\]
is equivalent to a twisting morphism 
\[
\zeta_\varphi: C_*^c(\Delta^1) \longrightarrow \mathfrak{g}~.
\]
By Corollary \ref{cor: same formulas as for RNV}, the decomposition of the element $a_{01}$ in $C_*^c(\Delta^1)$ under the coalgebra structure gives the same formulae as in \cite{robertnicoud2020higher}. The data of a twisting morphism $\zeta_\varphi$ is thus equivalent to the data of two Maurer-Cartan elements $\alpha$ and $\beta$ in $\mathfrak{g}$ and to the data of $\lambda$ in $\mathfrak{g}_1$ such that $\lambda \bullet \alpha = \beta$, see \cite[Proposition 4.3]{daniel19} for an analogue statement.
\end{proof}

\begin{theorem}\label{thm: characterisation pi zero et jauges}
Let $\mathfrak{g}$ be a complete curved absolute $\mathcal{L}_\infty$-algebra. There is a bijection 
\[
\pi_0(\mathcal{R}(\mathfrak{g})) \cong \mathcal{MC}(\mathfrak{g})/\sim_{\mathrm{gauge}}~,
\]
where the right-hand side denotes the set of Maurer-Cartan elements up to gauge equivalence. 
\end{theorem}

\begin{proof}
Follows directly from Proposition \ref{prop: equivalence between gauge and paths}.
\end{proof}

\begin{Remark}
Let $\mathfrak{g}$ be a complete curved absolute $\mathcal{L}_\infty$-algebra. Morphisms of the form
\[
\varphi: \widehat{\Omega} C_*^c(\Delta^1) \longrightarrow \mathfrak{g}
\]
are in a one-to-one correspondence with morphisms of $u\mathcal{CC}_\infty$-coalgebras
\[
\varphi^{\sharp}: C_*^c(\Delta^1) \longrightarrow \widehat{\mathrm{B}}_\iota \mathfrak{g}~.
\]
Since $C_*^c(\Delta^1)$ is the interval object in the model category of $u\mathcal{CC}_\infty$-coalgebras, the data of a gauge equivalence between two Maurer-Cartan elements is equivalent to the data of a left homotopy of morphisms of $u\mathcal{CC}_\infty$-coalgebras between the two associated group-like elements $\varphi^{\sharp}|_{a_0}$ and $\varphi^{\sharp}|_{a_1}$.
\end{Remark}

\subsection{Higher Baker--Campbell--Hausdorff formulae}\label{subsection: higher BCH}
In this subsection, we generalize the seminal results about the higher Baker--Campbell--Hausdorff formulae for complete $\mathcal{L}_\infty$-algebras developed in \cite[Section 5]{robertnicoud2020higher} to the case of curved absolute $\mathcal{L}_\infty$-algebras. The higher Baker--Campbell--Hausdorff provide us with \textit{explicit} formulae for the horn-fillers in $\mathcal{R}(\mathfrak{g})$, making it not only an $\infty$-groupoid, but an \textit{algebraic} $\infty$-groupoid with \textit{canonical} horn-fillers.

\medskip

Consider the graded pdg module $U(n)$ which is given by $\kk.u$ in degree $n$ and by $\kk.du$ in degree $n-1$, together with the pre-differential $d(u) \coloneqq du$. Notice that 

\[
\widehat{\mathscr{S}}^c(\mathrm{B}^{\mathrm{s.a}}\ucom)(U(n)) \cong  \prod_{m \geq 0} \widehat{\Omega}^{\mathrm{s.a}}\ucom^*(m) ~\widehat{\otimes}_{\mathbb{S}_m}~ U(n)^{\otimes m}
\]

is the free complete pdg $\mathrm{B}^{\mathrm{s.a}}\ucom$-algebra generated by $U(n)$, but it \textbf{does not} form a complete \textit{curved} $\mathrm{B}^{\mathrm{s.a}}\ucom$-algebra. For more explanations on this kind of phenomena, see Section \ref{Section: Curved operadic bimodules}.

\begin{lemma}\label{lemma: decomposition de L(Delta n)}
There is an isomorphism of complete pdg $\mathrm{B}^{\mathrm{s.a}}\ucom$-algebras 
\[
\mathcal{L}(\Delta^n) \cong \widehat{\mathscr{S}}^c(\mathrm{B}^{\mathrm{s.a}}\ucom)(U(n)) \amalg \mathcal{L}(\Lambda_k^n)
\]
for all $n \geq 2$ and for $1 \leq k \leq n$, where $\Lambda_k^n$ denotes the $k$-horn of dimension $n$. 

\end{lemma}

\begin{proof}
The proof is \textit{mutatis mutandis} the same as that of \cite[Lemma 5.1]{robertnicoud2020higher}, where this time we use the completeness of the underlying canonical filtration instead. The bijections constructed in \textit{loc.cit.} also work in this case, and the formula for the element $a_{\widehat{k}}$ is the same by Theorem \ref{theom: formulas for the pre differentials of mc}. 
\end{proof}

\begin{theorem}\label{thm: canonical horn fillers}
Let $\mathfrak{g}$ be a complete curved absolute $\mathcal{L}_\infty$-algebra. There is a bijection

\[
\mathrm{Hom}_{\mathsf{sSet}}(\Delta^n, \mathcal{R}(\mathfrak{g})) \cong  \mathfrak{g}_n \times \mathrm{Hom}_{\mathsf{sSet}}(\Lambda_k^n,\mathcal{R}(\mathfrak{g}))~,
\]

\vspace{0.5pc}
natural in $\mathfrak{g}$ with respect to morphisms of curved absolute $\mathcal{L}_\infty$-algebras.
\end{theorem}

\begin{proof}
Consider the following isomorphisms

\begin{align*}
\mathrm{Hom}_{\mathsf{sSet}}(\Delta^n, \mathcal{R}(\mathfrak{g})) &\cong \mathrm{Hom}_{\mathsf{curv}~\mathsf{abs}~\mathcal{L}_\infty\text{-}\mathsf{alg}}(\mathcal{L}(\Delta^n), \mathfrak{g}) \\
&\cong \mathrm{Hom}_{\mathsf{pdg}~\mathrm{B}^{\mathrm{s.a}}\ucom\text{-}\mathsf{alg}}(\mathcal{L}(\Delta^n), \mathfrak{g}) \\
&\cong \mathrm{Hom}_{\mathsf{pdg}~\mathrm{B}^{\mathrm{s.a}}\ucom\text{-}\mathsf{alg}}(\widehat{\mathscr{S}}^c(\mathrm{B}^{\mathrm{s.a}}\ucom)(U(n)) \amalg \mathcal{L}(\Lambda_k^n), \mathfrak{g}) \\
&\cong  \mathrm{Hom}_{\mathsf{pdg}~\mathrm{B}^{\mathrm{s.a}}\ucom\text{-}\mathsf{alg}}(\widehat{\mathscr{S}}^c(\mathrm{B}^{\mathrm{s.a}}\ucom)(U(n)), \mathfrak{g}) \times  \mathrm{Hom}_{\mathsf{pdg}~\mathrm{B}^{\mathrm{s.a}}\ucom\text{-}\mathsf{alg}}(\mathcal{L}(\Lambda_k^n), \mathfrak{g}) \\
&\cong \mathfrak{g}_n \times  \mathrm{Hom}_{\mathsf{curv}~\mathsf{abs}~\mathcal{L}_\infty\text{-}\mathsf{alg}}(\mathcal{L}(\Lambda_k^n), \mathfrak{g}) \\
&\cong \mathfrak{g}_n \times \mathrm{Hom}_{\mathsf{sSet}}(\Lambda_k^n,\mathcal{R}(\mathfrak{g}))~,
\end{align*}

where we used several times that curved $\mathrm{B}^{\mathrm{s.a}}\ucom$-algebras are a full subcategory of pdg $\mathrm{B}^{\mathrm{s.a}}\ucom$-algebras, and where we used Lemma \ref{lemma: decomposition de L(Delta n)} in line three.
\end{proof}

\begin{Corollary}
Let $\mathfrak{g}$ be a complete curved absolute $\mathcal{L}_\infty$-algebra. The $\infty$-groupoid $\mathcal{R}(\mathfrak{g})$ has a \textit{canonical structure} of $\infty$-groupoid. 
\end{Corollary}

\begin{proof}
For any lifting problem $\Lambda_k^n \longrightarrow \mathcal{R}(\mathfrak{g})$ there is a canonical lift given by the image of $0$ inside $\mathfrak{g}_n$ under the bijection of Theorem \ref{thm: canonical horn fillers}.
\end{proof}

\begin{Remark}
Thus the integration functor $\mathcal{R}$ lands on the category of \textit{algebraic} $\infty$-groupoids in the sense of \cite{algebraic}. This shifts the point of view in the following way: it not only satisfies a \textit{property} (existence of horn-filler) but instead these horn-fillers a given and thus are considered as a \textit{structure}.
\end{Remark}

An element in $\mathcal{R}(\mathfrak{g})_n$ amounts to the data of a curved twisting morphism
\[
\varphi: C_*^c(\Delta^n) \longrightarrow \mathfrak{g}~.
\]
It can be pictorially represented by geometric $n$-simplex where the $\mathrm{I}$-th face $\mathrm{I} \subseteq [n]$ of the $n$-simplex is labeled by the element $\varphi(a_\mathrm{I})$ in $\mathfrak{g}$. As an example, an element in $\mathcal{R}(\mathfrak{g})_2$ can be represented as 

\begin{center}
\includegraphics[width=60mm,scale=0.7]{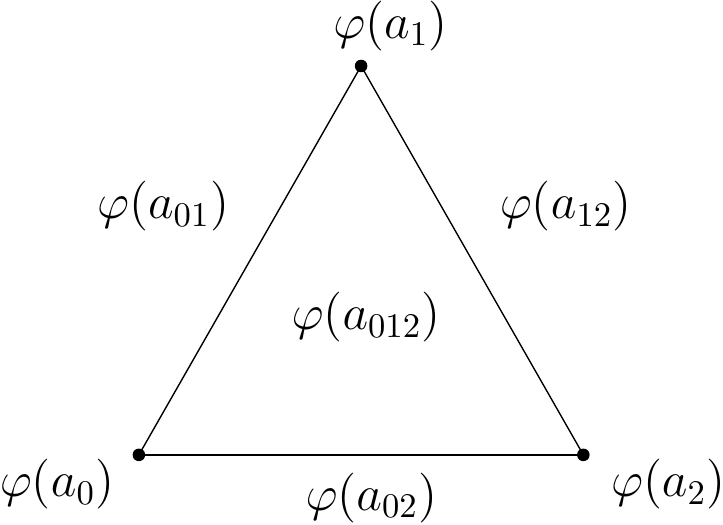},
\end{center}

where $\varphi(a_{i})$ are Maurer-Cartan elements in $\mathfrak{g}_0$ for $0 \leq i \leq 2$, where $\varphi(a_{ij})$ are degree $1$ elements in $\mathfrak{g}$ which induce gauge equivalences and where $\varphi(a_{012})$ is a degree $2$ element satisfying a higher compatibility condition imposed by the curved twisting morphism condition on $\varphi$. Similarly, the $\Lambda_k^n$-horns $\Lambda_k^n \longrightarrow \mathcal{R}(\mathfrak{g})$ admit an analogue pictorial description. As an example, a $\Lambda_1^2$-horn  in $\mathcal{R}(\mathfrak{g})$ can be depicted as

\begin{center}
\includegraphics[width=60mm,scale=0.7]{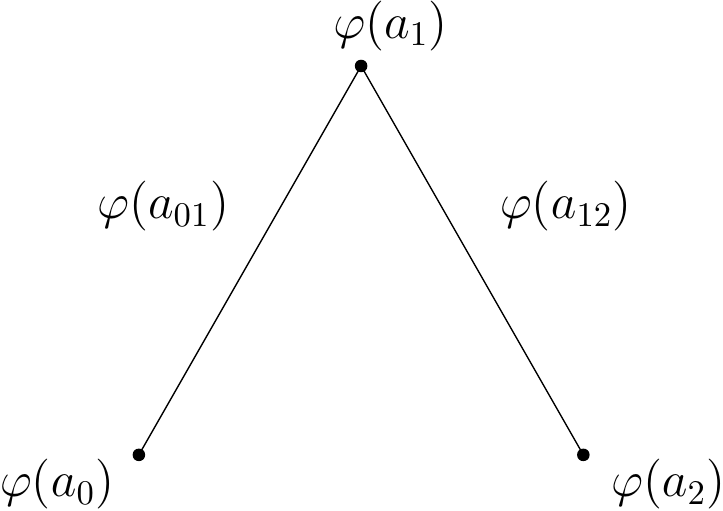},
\end{center}

\begin{Definition}[Higher Baker--Campbell--Hausdorff products]
Let $\mathfrak{g}$ be a complete curved absolute $\mathcal{L}_\infty$-algebra and let $y$ be an element in $\mathfrak{g}$ of degree $n$. The \textit{higher Baker--Campbell--Hausdorff product relative to} $y$ is given by 
\[
\begin{tikzcd}[column sep=4pc,row sep=-0pc]
\Gamma^y: \mathrm{Hom}_{\mathsf{sSet}}(\Lambda_k^n,\mathcal{R}(\mathfrak{g})) \arrow[r]
&\mathfrak{g}_{n-1} \\
x \arrow[r,mapsto]
&\Gamma^y(x) \coloneqq \varphi^{(x;y)}(a_{\widehat{k}})
\end{tikzcd}
\]
sends any $\Lambda_k^n$-horn $x: \Lambda_k^n \longrightarrow \mathcal{R}(\mathfrak{g})$ to the image $\varphi^{(x;y)}(a_{\widehat{k}})$ in $\mathfrak{g}_{n-1}$, where $\varphi^{(x;y)}$ is the curved twisting morphism associated to the element in $\mathcal{R}(\mathfrak{g})_n$ given by the bijection constructed in Theorem \ref{thm: canonical horn fillers}. 
\end{Definition}

These higher Baker--Campbell--Hausdorff products are given by explicit formulae. In our context, we recover the same formulas as the ones obtained in \cite[Proposition 5.10]{robertnicoud2020higher}.

\begin{Proposition}\label{prop: formules higher BCH}
Let $\mathfrak{g}$ be a complete curved absolute $\mathcal{L}_\infty$-algebra. Let $x: \Lambda_k^n \longrightarrow \mathcal{R}(\mathfrak{g})$ be a $\Lambda_k^n$-horn in $\mathcal{R}(\mathfrak{g})$. The higher Baker--Campbell--Hausdorff product relative to an element $y$ in $\mathfrak{g}_n$ is given by the following formula

\[
\Gamma^y(x)= \gamma_\mathfrak{g}\left(\sum_{\begin{subarray}{c}\tau\in\mathrm{PaPRT}\\
\chi\in\mathrm{Lab}^{[n],k}(\tau)\end{subarray}}\ \prod_{\begin{subarray}{c}%
\beta\text{ block of }\tau\\
\lambda^{\beta(\chi)}_{[n]}\neq 0\end{subarray}}\frac{(-1)^{k}}{\lambda^{\beta%
(\chi)}_{[n]}[\beta]!}\, \tau\left(x_{\chi(1)},\ldots,x_{\chi(p)};(-1)^{%
k}d_\mathfrak{g}(y)-\sum_{l\neq k}(-1)^{k+l}x_{\widehat{l}}\right) \right)~.
\]
See \cite[Proposition 5.10]{robertnicoud2020higher} for more details on the coefficients of this formula.
\end{Proposition}

\begin{proof}
The isomorphism of Lemma \ref{lemma: decomposition de L(Delta n)} is given by the same formulae that in \cite{robertnicoud2020higher}, and this isomorphism determines the higher Baker--Campbell-Hausdorff products.
\end{proof}

\begin{Remark}
Using the completeness of $\mathfrak{g}$, this last formula can be rewritten by splitting it along the weight filtration. 
\end{Remark}

These higher Baker--Campbell-Hausdorff products do recover the classical Baker--Campbell-Hausdorff product when one restricts to the case of nilpotent Lie algebras.

\begin{lemma}
Let $\mathfrak{g}$ be a pdg module concentrated in degree $1$. The data of a nilpotent curved $\mathcal{L}_\infty$-algebra structure on $\mathfrak{g}$ is equivalent to the data of a nilpotent Lie algebra on $\mathfrak{g}$.  
\end{lemma}

\begin{proof}
Let $\mathfrak{g}$ be a nilpotent curved $\mathcal{L}_\infty$-algebra concentrated in degree $1$ and let 
\[
\gamma_\mathfrak{g}: \prod_{n \geq 0} \widehat{\Omega}^{\mathrm{s.a}}\ucom^*(n) ~\widehat{\otimes}_{\mathbb{S}_n} ~ \mathfrak{g}^{\otimes n} \longrightarrow \mathfrak{g}~,
\]
be its structural morphism. For degree reasons, only series involving binary rooted trees can have a non-zero image via the structural morphism $\gamma_\mathfrak{g}$. Therefore it can be seen as an algebra over the conilpotent cooperad $\mathcal{L}ie^*$. Since the cooperad $\mathcal{L}ie^*$ is binary, nilpotent $\mathcal{L}ie^*$-algebras are equivalent to nilpotent $\mathcal{L}ie$-algebras, see Section \ref{Section: Absolute algebras and contramodules}.
\end{proof}

Let $\mathfrak{g}$ be a nilpotent curved $\mathcal{L}_\infty$-algebra concentrated in degree $1$ and let $x: \Lambda_1^2 \longrightarrow \mathcal{R}(\mathfrak{g})$ be a $\Lambda_1^2$-horn. Since $\mathfrak{g}$ is concentrated in degree $1$, the only Maurer-Cartan element is $0$ and the data of such a horn amount to the data of two elements $\alpha$ and $\beta$ in $\mathfrak{g}$. 
\vspace{0.4pc}

\begin{center}
\includegraphics[width=45mm,scale=0.5]{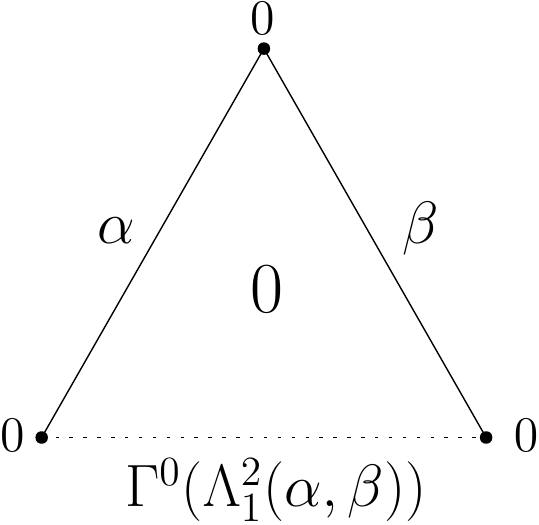},
\end{center}

Since $\mathfrak{g}_2 = \{0\}$, there is only one way to fill this horn. This filling corresponds to the classical Baker--Campbell--Hausdorff formula.

\begin{Proposition}\label{prop: horn 2,1 c'est BCH}
Let $\mathfrak{g}$ be a nilpotent curved $\mathcal{L}_\infty$-algebra concentrated in degree $1$. Let $\alpha$ and $\beta$ be two degree $1$ elements, then 
\[
\Gamma^{0}(\Lambda_1^2(\alpha,\beta)) = \mathrm{BCH}(\alpha,\beta)~,
\]
where $\Lambda_1^2(\alpha,\beta)$ stands for the horn labeled by $\alpha$ and $\beta$ and where $\mathrm{BCH}(\alpha,\beta)$ stands for the classical Baker--Campbell-Hausdorff formula.
\end{Proposition}

\begin{proof}
By Proposition \ref{prop: formules higher BCH}, we obtain the same formulas as in \cite{robertnicoud2020higher}. By Theorem 5.17 in \textit{loc.cit.}, the formula this higher BCH product coincides with the classical BCH formula. 
\end{proof}

\begin{Remark}
This above result was first proved in R. Bandiera's PhD Thesis, see \cite{bandiera}.
\end{Remark}

\begin{Remark}
It is known that the Baker--Campbell--Hausdorff formula lives in the free complete Lie algebra generated by two elements. This algebra is in fact the \textit{free absolute} Lie algebra, hence the Baker--Campbell--Hausdorff formula naturally belongs to the realm of absolute algebras. See Section \ref{Section: Absolute algebras and contramodules} for more details on absolute Lie algebras.
\end{Remark}

Let $\mathfrak{g}$ be a nilpotent Lie algebra. We denote 
\[
\mathrm{Exp}(\mathfrak{g}) \coloneqq \left(\mathfrak{g}, \mathrm{BCH}, 0 \right)
\]
the group obtained by considering $\mathfrak{g}$ equipped with the BCH formula. 

\begin{Corollary}\label{cor: integration of absolute Lie}
Let $\mathfrak{g}$ be a nilpotent curved $\mathcal{L}_\infty$-algebra concentrated in degree $1$. There is an isomorphism of simplicial sets 
\[
\mathcal{R}(\mathfrak{g}) \cong \mathcal{N}(\mathrm{Exp}(\mathfrak{g}))~,
\]
where $\mathcal{N}$ stands for the nerve of a category.
\end{Corollary}

\begin{proof}
It follows from Proposition \ref{prop: horn 2,1 c'est BCH}. 
\end{proof}

\begin{Corollary}
Let $\mathfrak{g}$ be a complete curved $\mathcal{L}_\infty$-algebra concentrated in degree $1$. There is an isomorphism of simplicial sets 
\[
\mathcal{R}(\mathfrak{g}) \cong \mathcal{N}\left(\lim_\omega \mathrm{Exp}(\mathfrak{g}/\mathrm{W}_\omega \mathfrak{g})\right)~.
\]
It is therefore the nerve of a pro-$\mathbb{K}$-unipotent group. 
\end{Corollary}

\begin{proof}
It follows from Corollary \ref{cor: integration of absolute Lie}, using the fact that the nerve functor commutes with limits.
\end{proof}

\subsection{Higher homotopy groups}
In this subsection, we compute the higher homotopy groups of the Kan complex $\mathcal{R}(\mathfrak{g})$. This is done by using the explicit $u\mathcal{CC}_\infty$-coalgebra structures on the cellular chains of the $n$-spheres. This new method allows us to generalize of the main result of \cite{Berglund15} to curved absolute $\mathcal{L}_\infty$-algebras. 

\medskip

Let $\mathbb{S}^n$ denote the simplicial set given by one non-degenerate $0$-simplex $[0]$ and one non-degenerate $n$-simplex $[n]$. This is a model for the simplicial set $\Delta^n/\partial \Delta^n$, which does not form a Kan complex. Nevertheless, since for any complete curved absolute $\mathcal{L}_\infty$-algebra $\mathfrak{g}$, $\mathcal{R}(\mathfrak{g})$ is a Kan complex, one has

\[
\pi_n(\mathcal{R}(\mathfrak{g}),\alpha) \cong \mathrm{Hom}_{\mathsf{sSets}_*}(\mathbb{S}^n, \mathcal{R}(\mathfrak{g}))/\sim_{\mathrm{hmt}}~,
\]
\vspace{0.1pc}

where the right-hand side is the set of morphisms of simplicial sets which send $[0]$ to $\alpha$ in $\mathcal{R}(\mathfrak{g})_0$ modulo the homotopy relation. Now notice that 

\[
\mathrm{Hom}_{\mathsf{sSets}}(\mathbb{S}^n, \mathcal{R}(\mathfrak{g}))/\sim_{\mathrm{hmt}} ~\cong \mathrm{Hom}_{u\mathcal{CC}_\infty\text{-}\mathsf{coalg}}\left(C_*^c(\mathbb{S}^n), \widehat{\mathrm{B}}_\iota(\mathfrak{g}) \right)/\sim_{\mathrm{hmt}}
\]
\vspace{0.1pc}

since $C_*^c(-) \dashv \overline{\mathcal{R}}$ is a Quillen adjunction and since both $\mathcal{R}(\mathfrak{g})$ and $\widehat{\mathrm{B}}_\iota \mathfrak{g}$ are fibrant objects. Hence by explicitly computing the $u\mathcal{CC}_\infty$-coalgebra structure on $C_*^c(\mathbb{S}^n)$ one can compute the homotopy groups of $\mathcal{R}(\mathfrak{g})$ for any complete curved absolute $\mathcal{L}_\infty$-algebra. 

\begin{lemma}\label{lemma: structure de cogèbre sur la sphere}
Let $n \geq 1$. The dg module $C_*^c(\mathbb{S}^n)$ is given by $\kk.a_0$ in degree $0$ and $\kk.a_{[n]}$ in degree $n$ with zero differential. The $u\mathcal{CC}_\infty$-coalgebra structure on $C_*^c(\mathbb{S}^n)$ is given by the elementary decomposition maps 

\[
\left\{
\begin{tikzcd}[column sep=4pc,row sep=0.5pc]
a_0 \arrow[r,mapsto]
&\Delta_{c_m}(a_0) = a_0 \otimes \cdots \otimes a_0~, \\
a_{[n]} \arrow[r,mapsto]
&\Delta_{c_m}(a_{[n]}) = a_0 \otimes \cdots \otimes a_0 \otimes a_{[n]}~, \\
\end{tikzcd}
\right.
\]

where $c_m$ denotes the $m$-corolla, for all $m \geq 2$, and where $\Delta_\tau$ is the zero morphism for any other corked rooted tree $\tau$ in $\mathrm{CRT}$. 
\end{lemma}

\begin{proof}
Consider the the following pushout in the category of simplicial sets 
\[
\begin{tikzcd}
\partial \Delta^n  \arrow[r,twoheadrightarrow] \arrow[d,rightarrowtail] \arrow[dr, phantom, "\lrcorner", very near end]
&\{*\} \arrow[d] \\
\Delta^n \arrow[r,twoheadrightarrow]
&\mathbb{S}^n~. 
\end{tikzcd}
\]
It induces a pullback in the category of $u\mathcal{CC}_\infty$-algebras 
\[
\begin{tikzcd}
C^*_c(\mathbb{S}^n) \arrow[r,twoheadrightarrow] \arrow[d,rightarrowtail] \arrow[dr, phantom, "\ulcorner", very near start]
&C^*_c(\{*\}) \arrow[d] \\
C^*_c(\Delta^n) \arrow[r,twoheadrightarrow]
&C^*_c(\partial \Delta^n)~.
\end{tikzcd}
\]
Therefore there is a epimorphism $C^*_c(\mathbb{S}^n) \twoheadrightarrow C^*_c(\Delta^n)$ of $u\mathcal{CC}_\infty$-algebras. Using Proposition \ref{prop: formulas for the HTT}, one can see that, for degree reasons, the only operations that survive on $C^*_c(\mathbb{S}^n)$ are the multiplications $\mu_{c_m}$. Furthermore, it is clear that $\mu_{c_m}(a_0,\cdots, a_0, a_{[n]}) = a_{[n]}$. Computing the $u\mathcal{CC}_\infty$-coalgebra structure of the linear dual $C_*^c(\mathbb{S}^n)$ from this is straightforward.
\end{proof}

\begin{Definition}[$\alpha$-homology groups]
Let $\mathfrak{g}$ be a curved absolute $\mathcal{L}_\infty$-algebra and let $\alpha$ be a Maurer--Cartan element. We define the $\alpha$-homology groups of $\mathfrak{g}$ to be
\[
\mathrm{H}^\alpha_*(\mathfrak{g}) \coloneqq \frac{\mathrm{Ker}(d^\alpha_\mathfrak{g})}{\mathrm{Im}(d^\alpha_\mathfrak{g})}~.
\]
\end{Definition}

\begin{Remark}
There should be a analogue notion of \textit{the twisting procedure} for curved absolute $\mathcal{L}_\infty$-algebras where given a Maurer-Cartan element $\alpha$, one can construct an \textit{absolute} $\mathcal{L}_\infty$-algebra $\mathfrak{g}^\alpha$. If this was the case, then we would have that 
\[
\mathrm{H}^\alpha_*(\mathfrak{g}) \cong \mathrm{H}_*(\mathfrak{g}^\alpha)~.
\]
Furthermore, given a curved absolute $\mathcal{L}_\infty$-algebra $\mathfrak{g}$, its image under the integration functor decomposes as 
\[
\mathcal{R}(\mathfrak{g}) \cong \underset{\alpha \in \pi_0(\mathfrak{R}(\mathfrak{g}))}{\coprod} \mathcal{R}(\mathfrak{g})^{\alpha}~,
\]
into a sum over its connected components. The twisted absolute $\mathcal{L}_\infty$-algebra $\mathfrak{g}^\alpha$ should be a model for the connected component $\mathcal{R}(\mathfrak{g})^{\alpha}$. 

\medskip

We could just set the twisted formulae in the \textit{absolute} setting and try to check by hand that they satisfy Condition \ref{pdg condition}, \ref{associativity condition} and \ref{curved condition}. We would nevertheless lack a conceptual explanation of these twisted algebras in terms of the action of a \textit{gauge group} as explained in \cite{DotsenkoShadrinVallette16}.
\end{Remark}

\begin{theorem}\label{thm: Berglund isomorphism}
Let $\mathfrak{g}$ be a complete curved absolute $\mathcal{L}_\infty$-algebra and let $\alpha$ be a Maurer--Cartan element. There is a bijection
\[
\pi_n(\mathcal{R}(\mathfrak{g}),\alpha) \cong \mathrm{H}^\alpha_n(\mathfrak{g})~,
\]
for all $n \geq 1$, which is natural in $\mathfrak{g}$.
\end{theorem}

\begin{proof}
Recall that the data of a morphism of $u\mathcal{CC}_\infty$-coalgebras $f_{\nu}: C_*^c(\mathbb{S}^n) \longrightarrow \widehat{\mathrm{B}}_\iota(\mathfrak{g})$ is equivalent to the data of a curved twisting morphism $\nu: C_*^c(\mathbb{S}^n) \longrightarrow \mathfrak{g}$ relative to $\iota$. In this case, by Lemma \ref{lemma: structure de cogèbre sur la sphere}, $\nu: C_*^c(\mathbb{S}^n) \longrightarrow \mathfrak{g}$ is a curved twisting morphism if and only if 
\[
\gamma_\mathfrak{g}\left(\sum_{n \geq 0,~ n \neq 1}\frac{1}{n!}c_n(\nu(a_0), \cdots , \nu(a_0)) \right) + d_\mathfrak{g}(\nu(a_0)) = 0~,
\]
that is, $\nu(a_0)$ is a Maurer-Cartan element, and if 

\[
\gamma_\mathfrak{g}\left(\sum_{n \geq 0,~ n \neq 1}\frac{1}{n!}c_n(\nu(a_0), \cdots , \nu(a_{[n]})) \right) + d_\mathfrak{g}(\nu(a_{[n]})) = 0~,
\]
that is, $d^{\nu(a_0)}_\mathfrak{g}(\nu(a_{[n]}) = 0$, where $d^{\nu(a_0)}$ is the twisted differential by $\nu(a_0)$. Therefore we have a bijection
\[
\mathrm{Hom}_{\mathsf{sSets}_*}((\mathbb{S}^n,a_0), (\mathcal{R}(\mathfrak{g}),\alpha)) \cong \mathrm{Z}_n^{\alpha}(\mathfrak{g})~.
\]
Let $f_\nu: C_*^c(\mathbb{S}^n) \longrightarrow \widehat{\mathrm{B}}_\iota(\mathfrak{g})$ and $g_{\rho}: C_*^c(\mathbb{S}^n) \longrightarrow \widehat{\mathrm{B}}_\iota(\mathfrak{g})$ be two morphisms of $u\mathcal{CC}_\infty$-coalgebras which send $a_0$ to $\alpha$. The data of an homotopy between them amounts to the data of a morphism
\[
h: C_*^c(\mathbb{S}^n) \otimes  C_*^c(\Delta^1) \longrightarrow \widehat{\mathrm{B}}_\iota(\mathfrak{g})
\]
such that $h(- \otimes b_0) = f_{\nu}$ and $h(- \otimes b_1) = g_{\rho}$, where $C_*^c(\Delta^1) = \kk.b_0 \oplus \kk.b_1 \oplus \kk.b_{01}$, since $C_*^c(\Delta^1)$ is the interval object in the category of $u\mathcal{CC}_\infty$-coalgebra. The $u\mathcal{CC}_\infty$-coalgebra structure on $C_*^c(\mathbb{S}^n) \otimes  C_*^c(\Delta^1)$ is given by the tensor product structure of Corollary \ref{cor: tensor product of uCC coalgebras}. One can check that the data of such of a morphism $h$ is equivalent to the data of an element $\lambda$ of degree $n+1$ such that 
\[
d^{\alpha}(\lambda) = \nu(a_{[n]}) - \rho(a_{[n]})~.
\]
Therefore there is a bijection 
\[
\mathrm{Hom}_{\mathsf{sSets}_*}((\mathbb{S}^n,a_0), (\mathcal{R}(\mathfrak{g}),\alpha))/\sim_{\mathrm{hmt}} ~\cong \mathrm{H}_n^{\alpha}(\mathfrak{g})~.
\]
\end{proof}

\begin{Corollary}
Let $\alpha$ and $\beta$ be two gauge equivalent Maurer-Cartan elements of $\mathfrak{g}$. There is an isomorphism
\[
\mathrm{H}^\alpha_n(\mathfrak{g}) \cong \mathrm{H}^\beta_n(\mathfrak{g})~,
\]
for all $n \geq 1$. 
\end{Corollary}

\begin{proof}
It follows from Theorem \ref{thm: Berglund isomorphism} and Theorem \ref{thm: characterisation pi zero et jauges}.
\end{proof}

\begin{Remark}[Berglund map]
One can also prove this result by following the same ideas of \cite{Berglund15}, that is, by constructing an explicit map
\[
\begin{tikzcd}[column sep=4pc,row sep=0pc]
\mathscr{B}_n^{\alpha}: \mathrm{H}_n^{\alpha}(\mathfrak{g}) \arrow[r]
&\pi_n(\mathcal{R}(\mathfrak{g}),\alpha) \\
\left[u\right] \arrow[r, mapsto]
& \left[ f_{[u]}: C^c_*(\Delta^n) \longrightarrow \mathfrak{g} \right]~,
\end{tikzcd}
\]
where $f_{[u]}$ is the curved twisting morphism defined by $f_{[u]}(a_{[n]}) = u$ and $f_{[u]}(a_{i}) = \alpha$, for all $0 \leq i \leq n$, and which is zero on any other element $a_I$ of $C^c_*(\Delta^n)$.
\end{Remark}

\begin{Remark}
The group structure on $\pi_n(\mathcal{R}(\mathfrak{g}),\alpha)$ is induced by the pinch map 
\[
\mathrm{pinch}: \mathbb{S}^n \longrightarrow \mathbb{S}^n \vee \mathbb{S}^n~.
\]
which on cellular chains is given by $C^c_*(\mathrm{pinch})(a_{[n]}) = a_{[n]}^{(1)} + a_{[n]}^{(2)}$. Using this, it can be shown that for $n \geq 2$, there is an isomorphism of abelian groups
\[
\pi_n(\mathcal{R}(\mathfrak{g}),\alpha) \cong \mathrm{H}_n^\alpha(\mathfrak{g})~,
\]
where we consider the sum of homology classes on the right-hand side. 
\end{Remark}

Let $h: \Lambda_1^2(\lambda_1,\lambda_2) \longrightarrow \mathfrak{g}$ be a $\Lambda_1^2$-horn in $\mathcal{R}(\mathfrak{g})$. Let $h(a_{01}) = \lambda_1$ and $h(a_{12}) = \lambda_2$. These are two degree one elements in $\mathfrak{g}$ which define a gauge equivalences $\lambda_1 \bullet \alpha = \alpha$ and $\lambda_2 \bullet \alpha = \alpha$. The \textit{Baker--Campbell--Hausdorff product relative to $\alpha$} of $\lambda_1$ and $\lambda_2$ is given by 

\[
\mathrm{BCH}^{\alpha}(\lambda_1,\lambda_2) \coloneqq \Gamma^0\left( \Lambda_1^2(\lambda_1,\lambda_2) \right)~.
\]

\begin{Remark}
The Baker--Campbell--Hausdorff product relative to $\alpha$ is not well-defined in general for any two degree one elements of $\mathfrak{g}$.
\end{Remark}

\begin{Corollary}
Let $\alpha$ be a Maurer--Cartan element of a curved absolute $\mathcal{L}_\infty$-algebra $\mathfrak{g}$. The Baker--Campbell--Hausdorff product relative to $\alpha$ defines a group structure with unit $[\alpha]$ on $\mathrm{H}_1^\alpha(\mathfrak{g})$, which is isomorphic to the group $\pi_1(\mathcal{R}(\mathfrak{g}),\alpha)$.
\end{Corollary}

\begin{proof}
This is straightforward from Theorem \ref{thm: Berglund isomorphism}.
\end{proof}

\subsection{Comparison with Getzler's functor}
In this subsection, we define the generalization of the integration functor as defined by E. Getzler in \cite{Getzler09} for complete curved $\mathcal{L}_\infty$-algebras in the sense of Definition \ref{def: complete curved L infinity algebra}. Then we state a comparison result between the integration functor $\mathcal{R}$ we constructed for complete curved absolute $\mathcal{L}_\infty$-algebras and the generalization of Getzler's functor for complete curved $\mathcal{L}_\infty$-algebras. 

\medskip 

Given a complete curved $\mathcal{L}_{\infty}$-algebra $\mathfrak{g}$, one can define the analogue of Hinich's functor as given in \cite{Hinich97Bis}: we consider the completed tensor product of $\mathfrak{g}$ with the simplicial unital commutative Sullivan algebra of polynomial differential forms $\Omega_\bullet$, which defines a complete curved $\mathcal{L}_{\infty}$-algebra $\mathfrak{g}$. In the end, we consider 
\[
\mathrm{MC}(\mathfrak{g})_{\bullet} \coloneqq \mathcal{MC}(\mathfrak{g} ~ \widehat{\otimes} ~ \Omega^\bullet)~. 
\]
This functor was already proven to be too big in the case of nilpotent $\mathcal{L}_\infty$-algebras, hence we impose a gauge condition on these sets of Maurer--Cartan element as it was done in \cite{Getzler09}.

\begin{Definition}[Generalized Getzler's functor]
Let $\mathfrak{g}$ be complete curved $\mathcal{L}_{\infty}$-algebra. The \textit{generalized Getzler's functor} is defined by

\[
\gamma(\mathfrak{g})_{\bullet} \coloneqq \Big\{ \alpha \in \mathcal{MC}(\mathfrak{g} ~ \widehat{\otimes} ~ \Omega^\bullet) ~| ~\hspace{1pt} (\mathrm{id} \otimes h_\bullet)(\alpha) = 0 \Big\}~. 
\]
\end{Definition}

\begin{Remark}
To the best of our knowledge, the tensor product of a unital commutative algebra with a curved \textit{absolute} $\mathcal{L}_\infty$-algebra does not seem to have a curved \textit{absolute} $\mathcal{L}_\infty$-algebra structure. Therefore, it is not clear to us whether the methods of \cite{Hinich97Bis} or \cite{Getzler09} could have been applied to the case of curved absolute $\mathcal{L}_\infty$-algebras. This is one justification for the approach chosen in the article.
\end{Remark}

Nevertheless, using the restriction functor constructed in Proposition \ref{prop: restriction functor}, we are able to compare the generalization of Getzler's functor with the integration functor $\mathcal{R}$.

\begin{Proposition}
Let $(\mathfrak{g},\gamma_\mathfrak{g},d_\mathfrak{g})$ be a complete curved absolute $\mathcal{L}_\infty$-algebra. Its image under the restriction functor  $(\mathrm{Res}(\mathfrak{g}),\mathrm{Res}(\gamma_\mathfrak{g}), \allowbreak \mathrm{W}_\omega \mathfrak{g},d_\mathfrak{g})$ forms a complete curved $\mathcal{L}_\infty$-algebra. This gives a functor 

\[
\mathrm{Res}: \mathsf{curv}~\mathsf{abs}~\mathcal{L}_\infty \textsf{-}\mathsf{alg}^{\mathsf{comp}} \longrightarrow \mathsf{curv}~\mathcal{L}_\infty \textsf{-}\mathsf{alg}^{\mathsf{comp}}
\]
\vspace{0.2pc}

from complete curved absolute $\mathcal{L}_\infty$-algebras to complete curved $\mathcal{L}_\infty$-algebras, which is faithful.
\end{Proposition}

\begin{proof}
The restriction $(\mathrm{Res}(\mathfrak{g}),\mathrm{Res}(\gamma_\mathfrak{g}),d_\mathfrak{g})$ forms a curved $\mathcal{L}_\infty$-algebra by Proposition \ref{prop: restriction functor}. The underlying filtration $\mathrm{W}_\omega \mathfrak{g}$ is complete and is indeed compatible with this curved $\mathcal{L}_\infty$-algebra structure 
\[
l_n(\mathrm{F}_{i_1}\mathfrak{g} \odot \cdots \odot \mathrm{F}_{i_n}\mathfrak{g}) \subseteq \mathrm{F}_{i_1 + \cdots i_k +1}\mathfrak{g}~,
\]
for all $i_1, \dots, i_k \geq 0$. Furthermore, since $\mathrm{Res}(\mathfrak{g})$ is a curved $\mathcal{L}_\infty$-algebra, $d_\mathfrak{g}^2(\mathrm{W}_\omega \mathfrak{g}) \subseteq \mathrm{W}_{\omega + 1} \mathfrak{g}~.$ Every morphism of complete curved absolute $\mathcal{L}_\infty$-algebras is continuous with respect to the canonical filtration, therefore it defines a functor which is clearly faithful.
\end{proof}

\begin{Proposition}\label{prop: comparaison avec le Getzler}
Let $\mathfrak{g}$ be a complete curved absolute $\mathcal{L}_\infty$-algebra. There is an isomorphism of simplicial sets 
\[
\gamma\left(\mathrm{Res}(\mathfrak{g})\right)_\bullet \cong \mathcal{MC}\left(\mathrm{Res}\left(\mathrm{hom}(C_*^c(\Delta^\bullet),\mathfrak{g})\right)\right)~,
\]
\vspace{0.1pc}

where $\mathrm{hom}(C_*^c(\Delta^\bullet),\mathrm{g})$ denotes the convolution curved absolute $\mathcal{L}_\infty$-algebra. 
\end{Proposition}

\begin{proof}
Notice that the Dupont contraction of Theorem \ref{dupontcontraction} induces the following contractions of complete pdg modules in the sense of Definition \ref{def: contraction of pdg modules}: 
\[
\begin{tikzcd}[column sep=5pc,row sep=3pc]
\mathfrak{g} ~ \widehat{\otimes} ~\Omega_{\bullet} \arrow[r, shift left=1.1ex, "\mathrm{id} ~ \otimes ~ p_{\bullet}"{name=F}] \arrow[loop left]{l}{\mathrm{id} ~ \otimes ~ h_{\bullet}}
&\mathfrak{g} ~ \widehat{\otimes} ~ C_c^*(\Delta^{\bullet})~, \arrow[l, shift left=.75ex, "\mathrm{id}_\mathfrak{g} ~ \otimes ~ i_{\bullet}"{name=U}]
\end{tikzcd}
\]
for all $\bullet \geq 0$, and these contractions are compatible with the simplicial structures. Thus one can apply Theorem \ref{thm: curved HTT}. Since the transferred structure coincides with the one of \cite{FukayaHTT}, we can use \cite[Theorem 1]{getzler2018maurercartan} to conclude that 
\[
\gamma(\mathfrak{g})_\bullet \cong \mathcal{MC}\left(\mathfrak{g} ~ \widehat{\otimes} ~ C_c^*(\Delta^{\bullet}) \right)~,
\]
where the transferred complete curved $\mathcal{L}_\infty$-algebra structure via the above simplicial contraction is considered on the right-hand side. The isomorphism of graded modules 
\[
\mathfrak{g} ~ \widehat{\otimes} ~ C_c^*(\Delta^{\bullet})  \cong \mathrm{Hom}_{\mathsf{gr}\text{-}\mathsf{mod}}(C_*^c(\Delta^{\bullet}), \mathfrak{g})~,
\]
allows us to rewrite the Maurer--Cartan equation in this complete curved $\mathcal{L}_\infty$-algebra. It is straightforward to check that there is an isomorphism of complete curved $\mathcal{L}_\infty$-algebras
\[
\mathrm{Hom}_{\mathsf{gr}\text{-}\mathsf{mod}}(C_*^c(\Delta^{\bullet}), \mathfrak{g}) \cong \mathrm{Res}(\mathrm{hom}(C_*^c(\Delta^\bullet),\mathfrak{g}))~.
\]
\end{proof}

\begin{Example}
It is difficult to say in general what Maurer--Cartan elements are well-defined in $\mathrm{Res}(\mathrm{hom}(C_*^c(\Delta^\bullet),\mathfrak{g}))$. For example, if $\mathfrak{g}$ is a nilpotent dg Lie algebra or an arity-wise nilpotent $\mathcal{L}_\infty$-algebra, then the Maurer--Cartan equation is always well-defined and we get an isomorphism 
\[
\gamma\left(\mathrm{Res}(\mathfrak{g})\right)_\bullet \cong \mathcal{R}(\mathfrak{g})~.
\]
As an example of this, consider Corollary \ref{cor: integration of absolute Lie}.
\end{Example}

\begin{Corollary}
Let $\mathfrak{g}$ be a complete curved absolute $\mathcal{L}_\infty$-algebra and let $\mathrm{W}_1\mathfrak{g}$ be the complete curved absolute $\mathcal{L}_\infty$-sub-algebra of weight $1$ elements. There is an isomorphism of simplicial sets 

\[
\gamma(\mathrm{W}_1\mathfrak{g})_\bullet \cong \mathcal{R}(\mathrm{W}_1\mathfrak{g})~.
\]
\end{Corollary}

\begin{proof}
Lemma \ref{lemma: splitting of the Maurer-Cartan equation} and Proposition \ref{prop: comparaison avec le Getzler} imply that

\[
\gamma(\mathrm{W}_1\mathfrak{g})_\bullet \cong \mathcal{MC}(\mathrm{Res}(\mathrm{hom}(C_*^c(\Delta^\bullet),\mathrm{W}_1\mathfrak{g}))) \cong \mathcal{MC}(\mathrm{hom}(C_*^c(\Delta^\bullet),\mathrm{W}_1\mathfrak{g})) \cong \mathcal{R}(\mathrm{W}_1\mathfrak{g})~,
\]
\vspace{0.1pc}

since the Maurer-Cartan infinite sum can be rewritten as an infinite sum of elements \textit{in} $\mathrm{W}_1\mathfrak{g}$ which \textit{converges} in its complete topology.
\end{proof}

\begin{Remark}
The hypothesis $\mathfrak{g} = \mathrm{W}_1 \mathfrak{g}$ is always assumed in \cite{robertnicoud2020higher}, so this result recovers their comparison with the Getzler functor in the non-curved case.
\end{Remark}

\begin{Remark}
The cosimplicial complete curved absolute $\mathcal{L}_\infty$-algebra $\mathfrak{mc}^\bullet$ induces a cosimplicial complete curved $\mathcal{L}_\infty$-algebra via the restriction functor $\mathrm{Res}$. This in turn produces an adjunction 

\[
\begin{tikzcd}[column sep=5pc,row sep=2.5pc]
\mathsf{sSet}  \arrow[r, shift left=1.5ex, "\mathcal{L}^{\mathrm{res}}"{name=A}]
&\mathsf{curv}~\mathcal{L}_\infty \textsf{-}\mathsf{alg}^{\mathsf{comp}}~, \arrow[l, shift left=.75ex, "\mathcal{R}^{\mathrm{res}}"{name=C}] \arrow[phantom, from=A, to=C, , "\dashv" rotate=-90]
\end{tikzcd}
\]

between the category of complete curved $\mathcal{L}_\infty$-algebras and the category of simplicial sets. The functor $\mathcal{R}^{\mathrm{res}}$ does \textbf{not} coincide with the generalization of Getzler's functor. Let $\mathfrak{g}$ be a complete curved absolute $\mathcal{L}_\infty$-algebra, we have inclusions 
\[
\gamma\left(\mathrm{Res}(\mathfrak{g})\right)_\bullet \subsetneq \mathcal{R}(\mathfrak{g})_\bullet \subset \mathcal{R}^{\mathrm{res}}(\mathfrak{g})_\bullet~,
\]
where the second inclusion is given by the fact that the restriction functor $\mathrm{Res}$ is faithful. Since the first inclusion is in general strict, the functor $\mathcal{R}^{\mathrm{res}}$ is not isomorphic to the generalization of Getzler's functor. 
\end{Remark}

\section{Rational homotopy models}
In this section, we show that curved absolute $\mathcal{L}_\infty$-algebras are models for finite type nilpotent rational spaces without any pointed or connectivity assumptions. We do this by relating the adjunction $\mathcal{L} \dashv \mathrm{R}$ constructed in the previous section to the derived adjunction constructed by Bousfield and Guggenheim in \cite{BousfieldGugenheim} using Sullivan's functor \cite{Sullivan77}. We then characterize the essential homotopical image of those spaces in the category of curved absolute $\mathcal{L}_\infty$-algebras. Finally, we construct explicit models for mapping spaces. 

\subsection{Comparing derived units of adjunctions}
From now on, we assume the base field $\kk$ to be the field of rational numbers $\mathbb{Q}$. Furthermore, we consider the category of simplicial sets endowed with the rational model structure constructed in \cite{Bousfield75}, where cofibrations are given by monomorphisms and where weak-equivalences are given by morphisms $f: X \longrightarrow Y$ such that $\mathrm{H}_*(f,\mathbb{Q}): \mathrm{H}_*(X,\mathbb{Q}) \longrightarrow \mathrm{H}_*(Y,\mathbb{Q})$ is an isomorphism.

\medskip

\textbf{Sullivan's adjunction.} One can use a dual version of Theorem \ref{thm: Kan seminal result} in order to induce the following contravariant adjunction 
\[
\begin{tikzcd}[column sep=5pc,row sep=2.5pc]
\mathsf{sSet} \arrow[r, shift left=1.5ex, "\mathrm{A}_{\mathrm{PL}}(-)"{name=A}]
&\mathsf{dg}~u\mathcal{C}om\textsf{-}\mathsf{alg}^{\mathsf{op}}~, \arrow[l, shift left=.75ex, "\langle - \rangle"{name=C}] \arrow[phantom, from=A, to=C, , "\dashv" rotate=-90]
\end{tikzcd}
\]
where $\mathrm{A}_{\mathrm{PL}}(-)$ is the piece-linear differential forms functor obtained by taking the left Kan extension of $\Omega^\bullet$. Its right adjoint functor, called the geometrical realization functor, is given by 
\[
\langle - \rangle_\bullet \coloneqq \mathrm{Hom}_{\mathsf{dg}~u\mathcal{C}om\textsf{-}\mathsf{alg}}(-, \Omega^\bullet)~.
\] 
This adjunction is in fact a Quillen adjunction when one considers the standard model structure on dg unital commutative algebras, where weak-equivalences are given by quasi-isomorphisms and fibrations by degree-wise epimorphisms.

\medskip

The quasi-isomorphism of dg operads $\varepsilon: \Omega \mathrm{B}^{\mathsf{s.a}}u\mathcal{C}om \qi u\mathcal{C}om$ induces the following Quillen equivalence 
\[
\begin{tikzcd}[column sep=5pc,row sep=2.5pc]
u\mathcal{CC}_\infty \textsf{-}\mathsf{alg} \arrow[r, shift left=1.5ex, "\mathrm{Ind}_\varepsilon"{name=A}]
&\mathsf{dg}~u\mathcal{C}om \textsf{-}\mathsf{alg} ~, \arrow[l, shift left=.75ex, "\mathrm{Res}_\varepsilon "{name=C}] \arrow[phantom, from=A, to=C, , "\dashv" rotate=-90]
\end{tikzcd}
\]
where $\mathrm{Res}_\varepsilon$ is fully faithful. Thus one obtains the following commutative triangle of Quillen adjunctions 
\[
\begin{tikzcd}[column sep=5pc,row sep=2.5pc]
&\hspace{2pc} \left(\mathsf{dg}~u\mathcal{C}om \textsf{-}\mathsf{alg}\right)^{\mathsf{op}} \arrow[dd, shift left=1.1ex, "\mathrm{Res}_\varepsilon^{\mathsf{op}}"{name=F}] \arrow[ld, shift left=.75ex, "\langle - \rangle"{name=C}]\\
\mathsf{sSet}  \arrow[ru, shift left=1.5ex, "\mathrm{A}_{\mathrm{PL}}(-)"{name=A}]  \arrow[rd, shift left=1ex, "\mathrm{C}_{\mathrm{PL}}(-)"{name=B}] \arrow[phantom, from=A, to=C, , "\dashv" rotate=-70]
& \\
&\hspace{1.5pc} \left(u\mathcal{CC}_\infty \textsf{-}\mathsf{alg}\right)^{\mathsf{op}} ~, \arrow[uu, shift left=.75ex, "\mathrm{Ind}_\varepsilon^{\mathsf{op}}"{name=U}] \arrow[lu, shift left=.75ex, "\langle - \rangle_\infty"{name=D}] \arrow[phantom, from=B, to=D, , "\dashv" rotate=-110] \arrow[phantom, from=F, to=U, , "\dashv" rotate=-180]
\end{tikzcd}
\]
by defining $\mathrm{C}_{\mathrm{PL}}(-) \coloneqq \mathrm{Res}_\varepsilon^{\mathsf{op}} \circ \mathrm{C}_{\mathrm{PL}}(-)$ and $\langle - \rangle_\infty \coloneqq \langle - \rangle \circ \mathrm{Ind}_\varepsilon^{\mathsf{op}}$.

\begin{lemma}\label{lemma: nieme equivalence}
Let $X$ be a simplicial set. There is a natural weak-equivalence of simplicial sets
\[
\mathbb{R}\langle \mathrm{C}_{\mathrm{PL}}(X) \rangle_\infty \simeq \mathbb{R}\langle \mathrm{A}_{\mathrm{PL}}(X)\rangle~.
\]
\end{lemma}

\begin{proof}
The adjunction induced by $\varepsilon: \Omega \mathrm{B}^{\mathsf{s.a}}u\mathcal{C}om \qi u\mathcal{C}om$ is a Quillen equivalence.
\end{proof}

\textbf{Cellular cochains functor.}
Consider again the simplicial $u\mathcal{CC}_\infty$-algebra $C_c^*(\Delta^\bullet)$ of Lemma \ref{lemma: simplicial CC infinity} given by the cellular cochains on the standard simplicies together with their transferred $u\mathcal{CC}_\infty$-algebra structure. It induces a Quillen adjunction

\[
\begin{tikzcd}[column sep=5pc,row sep=2.5pc]
\mathsf{sSet} \arrow[r, shift left=1.5ex, "C_c^*(-)"{name=A}]
&u\mathcal{CC}_\infty \textsf{-}\mathsf{alg}^{\mathsf{op}}~, \arrow[l, shift left=.75ex, "\mathcal{S}"{name=C}] \arrow[phantom, from=A, to=C, , "\dashv" rotate=-90]
\end{tikzcd}
\]

by applying again Theorem \ref{thm: Kan seminal result}, where $C_c^*(-)$ is the cellular cochain functor endowed with a $u\mathcal{CC}_\infty$-algebra structure, and where the right adjoint functor $\mathcal{S}$ is given by 

\[
\mathcal{S}(-) = \mathrm{Hom}_{u\mathcal{CC}_\infty \textsf{-}\mathsf{alg}}(-,C_c^*(\Delta^\bullet))~.
\]

Our goal is to compare the derived unit of this adjunction with the derived unit of the adjunction obtained by extending Sullivan's functor to the category of $u\mathcal{CC}_\infty$-algebras. In order, to do so, we need to be able to construct cofibrant resolutions. The canonical curved twisting morphism $\iota: \mathrm{B}^{\mathsf{s.a}}u\mathcal{C}om \longrightarrow \Omega \mathrm{B}^{\mathsf{s.a}}u\mathcal{C}om$ induces a Quillen equivalence
\[
\begin{tikzcd}[column sep=5pc,row sep=2.5pc]
\mathsf{curv}~ \mathrm{B}^{\mathsf{s.a}}u\mathcal{C}om \textsf{-}\mathsf{coalg} \arrow[r, shift left=1.5ex, "\Omega_\iota"{name=A}]
&u\mathcal{CC}_\infty \textsf{-}\mathsf{alg} ~, \arrow[l, shift left=.75ex, "\mathrm{B}_\iota"{name=C}] \arrow[phantom, from=A, to=C, , "\dashv" rotate=-90]
\end{tikzcd}
\]
where the model structure considered on the right hand side is the one obtained by transfer along this adjunction, see \cite{grignou2019} for more details. An $\infty$-morphism of $u\mathcal{CC}_\infty$-algebras $f: A \rightsquigarrow B$ is the data of a morphism 
\[
f: \mathrm{B}_\iota A \longrightarrow \mathrm{B}_\iota B
\]
of curved $\mathrm{B}^{\mathsf{s.a}}u\mathcal{C}om$-coalgebras. It is an $\infty$-quasi-isomorphism if the term $f_{\mathrm{id}}: A \longrightarrow B$ is a quasi-isomorphism. This is equivalent to $f$ being a weak-equivalence in the transferred model structure.

\begin{lemma}
Let $X$ be a simplicial set. There is a pair of inverse $\infty$-quasi-isomorphism of $u\mathcal{CC}_\infty$-algebras
\[
(p_\infty)_X: \mathrm{C}_{\mathrm{PL}}(X) \rightsquigarrow C_c^*(X)~,~(i_\infty)_X : C_c^*(X) \rightsquigarrow \mathrm{C}_{\mathrm{PL}}(X)~,
\]
which are natural in $X$.
\end{lemma}

\begin{proof}
Let $X$ be a simplicial set. Since Dupont's contraction is compatible with the simplicial structures, it induces a homotopy retract 

\[
\begin{tikzcd}[column sep=5pc,row sep=3pc]
\mathrm{C}_{\mathrm{PL}}(X) \arrow[r, shift left=1.1ex, "p_X"{name=F}] \arrow[loop left]{l}{h_X}
& C_c^*(X)~. \arrow[l, shift left=.75ex, "i_X"{name=U}]
\end{tikzcd}
\]
One can show, using analogue arguments to the proof of \cite[Proposition 7.12]{robertnicoud2020higher}, that the transferred $u\mathcal{CC}_\infty$-algebra structure from this contraction onto $C_c^*(X)$ is equal to the $u\mathcal{CC}_\infty$-algebra structure obtained by considering the left Kan extension of $C_c^*(\Delta^\bullet)$. Therefore the morphisms of dg modules $p_X$ and $i_X$ can be extended to two inverse $\infty$-quasi-isomorphism $(p_\infty)_X$ and $(i_\infty)_X$, which are both natural in $X$, see \cite[Remark 7.14]{robertnicoud2020higher}.
\end{proof}

\begin{Proposition}\label{prop: equivalence cochains CC with Apl}
Let $X$ be a simplicial set. There is a natural weak-equivalence of simplicial sets
\[
\mathbb{R}\langle \mathrm{C}_{\mathrm{PL}}(X) \rangle_\infty \simeq \mathbb{R}\mathcal{S}(C_c^*(X))~.
\]
\end{Proposition}

\begin{proof}
The above lemma implies that there is a natural weak-equivalence of curved $\mathrm{B}^{\mathsf{s.a}}u\mathcal{C}om$-coalgebras
\[
(p_\infty)_X: \mathrm{B}_\iota \mathrm{C}_{\mathrm{PL}}(X) \qi \mathrm{B}_\iota C_c^*(X)~.
\]
This in turn implies that there is a natural quasi-isomorphism 
\[
\Omega_\iota((p_\infty)_X) : \Omega_\iota\mathrm{B}_\iota \mathrm{C}_{\mathrm{PL}}(X) \qi \Omega_\iota\mathrm{B}_\iota C_c^*(X)
\]
of $u\mathcal{CC}_\infty$-algebras. Therefore we have 

\begin{align*}
\mathbb{R}\langle \mathrm{C}_{\mathrm{PL}}(X)\rangle_\infty &\simeq \mathrm{Hom}_{u\mathcal{C}om\text{-}\mathsf{alg}}(\mathrm{Ind}_\varepsilon \Omega_\iota\mathrm{B}_\iota \mathrm{C}_{\mathrm{PL}}(X), \Omega_\bullet)\\
& \cong \mathrm{Hom}_{u\mathcal{CC}_\infty\text{-}\mathsf{alg}}(\Omega_\iota\mathrm{B}_\iota \mathrm{C}_{\mathrm{PL}}(X), \mathrm{Res}_\varepsilon \Omega_\bullet) \\
&\simeq \mathrm{Hom}_{u\mathcal{CC}_\infty\text{-}\mathsf{alg}}(\Omega_\iota\mathrm{B}_\iota C_c^*(X), C_c^*(\Delta^\bullet)) \\
&\simeq \mathbb{R}\mathcal{S}(C_c^*(X))~,\\
\end{align*}

where we used the existence of an $\infty$-quasi-isomorphism of $u\mathcal{CC}_\infty$-algebras between $\mathrm{Res}_\varepsilon \Omega_\bullet$ and $C_c^*(\Delta^\bullet)$, given again by the Dupont contraction. Notice that the intermediate equivalences are obtained using the fact that these are hom-spaces between a cofibrant and a fibrant object, thus stable by weak-equivalences.
\end{proof}

\textbf{Finite type nilpotent spaces.} We recall what finite type nilpotent spaces are and we state the main theorem of this subsection.

\begin{Definition}[Finite type simplicial set]
Let $X$ be a simplicial set. It is said to be of \textit{finite type} if the homology groups $\mathrm{H}_n(X,\mathbb{Q})$ are finite dimensional for all $n \geq 0$.
\end{Definition}

\begin{Remark}
In particular, $X$ has a finitely many connected components since $\mathrm{H}_0(X,\mathbb{Q})$ is finite dimensional.
\end{Remark}

\begin{Definition}[Nilpotent simplicial set]
Let $X$ be a simplicial set. It is said to be \textit{nilpotent} if for every $0$-simplex $\alpha$ it satisfies the following conditions:

\begin{enumerate}
\item The group $\pi_1(X,\alpha)$ is nilpotent.

\medskip

\item The $\pi_1(X,\alpha)$-module $\pi_n(X,\alpha)$-module is a nilpotent $\pi_1(X,\alpha)$-module.
\end{enumerate}
\end{Definition}

\begin{Proposition}\label{prop: n-ieme equivalence}
Let $X$ be a finite type simplicial set. There is a weak-equivalence of simplicial sets

\[
\mathcal{R}\mathcal{L}(X) \simeq \mathbb{R}\mathcal{S}(C_c^*(X))~,
\]

which is natural on the subcategory of finite type simplicial sets.
\end{Proposition}

\begin{proof}
Recall from the proof above that 
\[
\mathbb{R}\mathcal{S}(C_c^*(X)) \simeq \mathrm{Hom}_{u\mathcal{CC}_\infty}(\Omega_\iota\mathrm{B}_\iota C_c^*(X), C_c^*(\Delta^\bullet))~.
\]

The following square of Quillen adjunctions
\[
\begin{tikzcd}[column sep=5pc,row sep=5pc]
u\mathcal{CC}_\infty\text{-}\mathsf{alg}^{\mathsf{op}} \arrow[r,"\mathrm{B}_\iota^{\mathsf{op}}"{name=B},shift left=1.1ex] \arrow[d,"(-)^\circ "{name=SD},shift left=1.1ex ]
&\mathsf{curv}~\mathrm{B}^{\mathsf{s.a}}u\mathcal{C}om\text{-}\mathsf{coalg}^{\mathsf{op}} \arrow[d,"(-)^*"{name=LDC},shift left=1.1ex ] \arrow[l,"\Omega_\iota^{\mathsf{op}}"{name=C},,shift left=1.1ex]  \\
u\mathcal{CC}_\infty\text{-}\mathsf{coalg} \arrow[r,"\widehat{\Omega}_\iota "{name=CC},shift left=1.1ex]  \arrow[u,"(-)^*"{name=LD},shift left=1.1ex ]
&\mathsf{curv}~\mathrm{B}^{\mathsf{s.a}}u\mathcal{C}om\text{-}\mathsf{alg}^{\mathsf{comp}}~, \arrow[l,"\widehat{\mathrm{B}}_\iota"{name=CB},shift left=1.1ex] \arrow[u,"(-)^\vee"{name=TD},shift left=1.1ex] \arrow[phantom, from=SD, to=LD, , "\dashv" rotate=0] \arrow[phantom, from=C, to=B, , "\dashv" rotate=-90]\arrow[phantom, from=TD, to=LDC, , "\dashv" rotate=0] \arrow[phantom, from=CC, to=CB, , "\dashv" rotate=-90]
\end{tikzcd}
\] 

commutes, by Theorem \ref{thm: homotopical magical square curved} applied to this situation. Since $X$ is a finite type simplicial set, the homology of $\Omega_\iota\mathrm{B}_\iota C_c^*(X)$ is degree-wise finite dimensional. Therefore the generalize Sweedler dual $(-)^\circ$ is homotopically fully-faithful by Proposition \ref{Prop: plongement pleinement fidele de Sweedler}, and we get 

\[
\mathrm{Hom}_{u\mathcal{CC}_\infty\text{-}\mathsf{alg}}(\Omega_\iota\mathrm{B}_\iota C_c^*(X), C_c^*(\Delta^\bullet)) \simeq \mathrm{Hom}_{u\mathcal{CC}_\infty\text{-}\mathsf{coalg}}(C_*^c(\Delta^\bullet), \widehat{\mathrm{B}}_\iota \widehat{\Omega}_\iota C_*^c(X) )~,
\]

which concludes the proof.
\end{proof}

\begin{Corollary}\label{cor: equivalences de foncteurs derives Apl et RL}
Let $X$ be a finite type simplicial set.  There is a weak-equivalence of simplicial sets

\[
\mathbb{R}\langle \mathrm{A}_{\mathrm{PL}}(X)\rangle \simeq \mathcal{R}\mathcal{L}(X)~,
\]

which is natural in the subcategory of finite type simplicial sets.
\end{Corollary}

\begin{proof}
This is a direct consequence of Lemma \ref{lemma: nieme equivalence} and Propositions \ref{prop: equivalence cochains CC with Apl} and \ref{prop: n-ieme equivalence}.
\end{proof}

Now we can transfer the known results about Sullivan's rational models along the above equivalence.

\begin{theorem}[{\cite[Theorem C]{MarklLazarev}}]\label{thm: Markl et Lazarev}
Let $X$ be a finite type nilpotent simplicial set. The unit 

\[
\eta_X: X \qi \mathbb{R}\langle \mathrm{A}_{\mathrm{PL}}(X)\rangle
\]

is a rational homotopy equivalence.
\end{theorem}

\begin{Remark}
The original construction of \cite{BousfieldGugenheim} only works for connected finite type nilpotent simplicial set and was subsequently extended to disconnected finite type nilpotent simplicial sets by Markl--Lazarev in \cite{MarklLazarev}.
\end{Remark}

\begin{theorem}\label{thm: modèles d'homotopie rationnel type fini}
Let $X$ be a finite type nilpotent simplicial set. The unit of adjunction

\[
\eta_X: X \qi \mathcal{R}\mathcal{L}(X) 
\]

is a rational homotopy equivalence.
\end{theorem}

\begin{proof}
Follows directly from Corollary \ref{cor: equivalences de foncteurs derives Apl et RL} and Theorem \ref{thm: Markl et Lazarev}.
\end{proof}

\begin{Corollary}
Let $X$ be a pointed connected finite type simplicial set. The unit of adjunction 
\[
\eta_X: X \longrightarrow \mathcal{RL}(X)
\]
is weakly equivalent to the $\mathbb{Q}$-completion of Bousfield-Kan. 
\end{Corollary}

\begin{proof}
This is implied by \cite[Theorem 12.2]{BousfieldGugenheim}, using the equivalence $\mathbb{R}\langle \mathrm{A}_{\mathrm{PL}}(X)\rangle \simeq \mathcal{RL}(X)$ of Corollary \ref{cor: equivalences de foncteurs derives Apl et RL}.
\end{proof}

\subsection{Minimal models}
In this subsection we show that curved absolute $\mathcal{L}_\infty$-algebras always admit a minimal resolutions, which are unique up to isomorphism. Furthermore, as a corollary of our constructions, we show that one can also recover the \textit{homology} of a space $X$ via the homology of the complete Bar construction of its models.

\begin{Definition}[Minimal model]
Let $\mathfrak{g}$ be a curved absolute $\mathcal{L}_\infty$-algebra. A \textit{minimal model} $(V,\varphi_{d_V},\psi_V)$ amounts to the data of 

\begin{enumerate}
\item A graded module $V$ together with a map

\[
\varphi_{d_V}: V \longrightarrow \prod_{n \geq 0} \widehat{\Omega}^{\mathrm{s.a}}\ucom^*(n)^{(\geq 1)} ~\widehat{\otimes}_{\mathbb{S}_n} ~ V^{\otimes n}~
\]

which lands on elements of weight greater or equal to $1$, such that the induced derivation $d_V$ satisfies

\[
d_V^2 = l_2 \circ_1 l_0~. 
\]

\item A weak equivalence of curved absolute $\mathcal{L}_\infty$-algebras 

\[
\psi_V: \left( \widehat{\Omega}^{\mathrm{s.a}}\ucom^*(n) ~\widehat{\otimes}_{\mathbb{S}_n} ~ V^{\otimes n},d_V \right) \qi \mathfrak{g}~.
\]
\end{enumerate}
\end{Definition}

\begin{Remark}
Given a graded module $V$, the data of the derivation $d_V$ amounts to a $u\mathcal{CC}_\infty$-coalgebra structure on $V$. Thus the data $(V,\varphi_{d_V})$ above amounts to the data of a \textit{minimal} $u\mathcal{CC}_\infty$-coalgebra, that is, a $u\mathcal{CC}_\infty$-coalgebra where the underlying differential is zero.
\end{Remark}

\begin{Proposition}
Let $\mathfrak{g}$ be a curved absolute $\mathcal{L}_\infty$-algebra and let $(V,\varphi_{d_V},\psi_V)$ and $(W,\varphi_{d_W},\psi_W)$ be two minimal models of $\mathfrak{g}$. Then there is an isomorphism of graded modules $V \cong W$. 
\end{Proposition}

\begin{proof}
If $(V,\varphi_{d_V},\psi_V)$ and $(W,\varphi_{d_W},\psi_W)$ are two minimal models of $\mathfrak{g}$, then, by definition, there exists a weak equivalence of curved absolute $\mathcal{L}_\infty$-algebras

\[
\left( \widehat{\Omega}^{\mathrm{s.a}}\ucom^*(n) ~\widehat{\otimes}_{\mathbb{S}_n} ~ V^{\otimes n},d_V \right) \qi \left( \widehat{\Omega}^{\mathrm{s.a}}\ucom^*(n) ~\widehat{\otimes}_{\mathbb{S}_n} ~ W^{\otimes n},d_W\right)~.
\]

Therefore there is a quasi-isomorphism between the graded modules $V$ and $W$, which implies they are isomorphic.
\end{proof}

\begin{Proposition}\label{prop; generateurs du model minimal, homologie de la Bar}
Let $\mathfrak{g}$ be a curved absolute $\mathcal{L}_\infty$-algebra. Then it admits a minimal model, where the graded module of generators is given by
\[
V \cong \mathrm{H}_*\left(\widehat{\mathrm{B}}_\iota \mathfrak{g} \right)~,
\]
that is, the homology of its complete Bar construction.
\end{Proposition}

\begin{proof}
Since $\kk$ is a field of characteristic $0$, one can always choose a contraction between $\widehat{\mathrm{B}}_\iota \mathfrak{g}$ and its homology in order to apply the Homotopy Transfer Theorem. The transferred $u\mathcal{CC}_\infty$-coalgebra structure provides us with the derivation by applying the complete Cobar construction. The image of the inclusion
\[
i:  \mathrm{H}_*\left(\widehat{\mathrm{B}}_\iota \mathfrak{g} \right) \twoheadrightarrow\widehat{\mathrm{B}}_\iota \mathfrak{g}~,
\]

under the complete Cobar construction allows us obtain the required weak-equivalence by post-composing it with $\eta_\mathfrak{g}: \mathfrak{g} \qi \widehat{\Omega}_\iota\widehat{\mathrm{B}}_\iota \mathfrak{g}~.$
\end{proof}

\begin{Remark}
The same minimal model constructions hold for absolute $\mathcal{L}_\infty$-algebras using  the complete Bar construction $\widehat{\mathrm{B}}_\iota^\flat$ relative to $\mathcal{CC}_\infty$-coalgebras.
\end{Remark}

\begin{Definition}[Rational curved absolute $\mathcal{L}_\infty$-algebras] 
Let $\mathfrak{g}$ be a curved absolute $\mathcal{L}_\infty$-algebra. It is \textit{rational} if there exists a simplicial set $X$ and a zig-zag of weak equivalences of curved absolute $\mathcal{L}_\infty$-algebras 
\[
\mathcal{L}(X) \lqi \cdot \qi \cdots \lqi \cdot \qi \mathfrak{g}~.
\]
\end{Definition}

\begin{Remark}
If $\mathfrak{g}$ is rational, then the homology of $\widehat{\mathrm{B}}_\iota(\mathfrak{g})$ is also concentrated in positive degrees. Indeed, $\mathcal{L}(X)$ is weakly equivalent to $\mathfrak{g}$ if and only if $C_*^c(X)$ is quasi-isomorphic to the complete Bar construction of $\mathfrak{g}$.
\end{Remark}

\begin{Definition}[Rational models] 
Let $\mathfrak{g}$ be a curved absolute $\mathcal{L}_\infty$-algebra. It is a \textit{rational model} if it is rational for some simplicial set $X$ and if furthermore $\mathcal{R}(\mathfrak{g})$ is weakly equivalent to $X$. 
\end{Definition}

\begin{Proposition}
Let $\mathfrak{g}$ be a curved absolute $\mathcal{L}_\infty$-algebra which is rational for a simplicial set $X$. Then its minimal model is generated by $\mathrm{H}_*(X)$.
\end{Proposition}

\begin{proof}
In this particular case, we have that 
\[
\mathrm{H}_*\left(\widehat{\mathrm{B}}_\iota \mathfrak{g} \right) \cong \mathrm{H}_*(X)~,
\]
therefore by Proposition \ref{prop; generateurs du model minimal, homologie de la Bar}, the graded module $\mathrm{H}_*(X)$ is the generator of the minimal model. 
\end{proof}

\begin{Remark}
Notice that the situation here is Koszul dual to Sullivan's minimal models, where the minimal Sullivan model for $\mathrm{A}_{\mathrm{PL}}(X)$ of a simply connected space $X$ is generated by the linear dual of the homotopy groups $\pi_*(X)^*$. 
\end{Remark}

\textbf{Recovering the homology groups.} Let $\mathfrak{g}$ be curved absolute $\mathcal{L}_\infty$-algebra, let us try to understand the homotopy type of $\mathcal{R}(\mathfrak{g})$. The \textit{homotopy groups} of $\mathcal{R}(\mathfrak{g})$ can be fully described by the homology groups $\mathrm{H}_*(\mathfrak{g}^\alpha)$, where $\alpha$ runs over the set of Maurer-Cartan elements of $\mathfrak{g}$. In the case where $\mathfrak{g}$ is a rational model, then one can also recover the \textit{homology groups} of $\mathcal{R}(\mathfrak{g})$ using only the complete Bar construction with respect to $\mathfrak{g}$. 

\begin{theorem}
Let $\mathfrak{g}$ be a curved absolute $\mathcal{L}_\infty$-algebra that is a rational model. The canonical morphism of $u\mathcal{CC}_\infty$-algebras 
\[
C_*^c(\mathcal{R}(\mathfrak{g})) \qi \widehat{\mathrm{B}}_\iota \mathfrak{g}
\]
is a quasi-isomorphism.
\end{theorem}

\begin{proof}
In this situation, $\mathcal{L}(\mathcal{R}(\mathfrak{g})) \qi \mathfrak{g}$ is a weak equivalence of curved absolute $\mathcal{L}_\infty$-algebras, which is equivalent to the canonical morphism $C_*^c(\mathcal{R}(\mathfrak{g})) \qi \widehat{\mathrm{B}}_\iota \mathfrak{g} $ being a quasi-isomorphism of $u\mathcal{CC}_\infty$-coalgebras.
\end{proof}

\begin{Remark}
One can think of $\widehat{\mathrm{B}}_\iota \mathfrak{g}$ as a \textit{higher Chevalley-Eilenberg} complex adapted to the setting of curved absolute $\mathcal{L}_\infty$-algebras. 
\end{Remark}

\subsection{Models for connected components of spaces}
Let $X$ be a pointed simplicial set. The base point $x: \{*\} \longrightarrow X$ induces a morphism $C_*^c(x): \kk \longrightarrow C_*^c(X)$ of $u\mathcal{CC}_\infty$-coalgebras. Since $C_*^c(X)$ is a strictly counital $u\mathcal{CC}_\infty$-coalgebra, the cellular chain functor then lands on the category of pointed strictly counital $u\mathcal{CC}_\infty$-coalgebras. 

\begin{Remark}
By a choice of a consistent base point, one can consider $C_*^c(\Delta^\bullet)$ as a cosimplicial pointed strictly counital $u\mathcal{CC}_\infty$-coalgebra.
\end{Remark}

\begin{Proposition}
There is a commuting square 
\[
\begin{tikzcd}[column sep=5pc,row sep=5pc]
\mathsf{sSet}_* \arrow[r,"C_*^c(-)"{name=B},shift left=1.1ex] \arrow[d,"\mathcal{L}_*"{name=SD},shift left=1.1ex ]
&\mathsf{Strict}~u\mathcal{CC}_\infty\text{-}\mathsf{coalg}_{\bullet} \arrow[d,"\mathrm{Ker}(\epsilon)"{name=LDC},shift left=1.1ex ] \arrow[l,"\overline{\mathcal{R}}"{name=C},,shift left=1.1ex]  \\
\mathsf{abs}~\mathcal{L}_\infty\text{-}\mathsf{alg}^{\mathsf{comp}} \arrow[r,"\widehat{\mathrm{B}}^\flat_\iota"{name=CC},shift left=1.1ex]  \arrow[u,"\mathcal{R}_*"{name=LD},shift left=1.1ex ]
&\mathcal{CC}_\infty\text{-}\mathsf{coalg}~. \arrow[l,"\widehat{\Omega}^\flat_\iota "{name=CB},shift left=1.1ex] \arrow[u,"(-)\oplus \kk"{name=TD},shift left=1.1ex] \arrow[phantom, from=SD, to=LD, , "\dashv" rotate=180] \arrow[phantom, from=C, to=B, , "\dashv" rotate=-90]\arrow[phantom, from=TD, to=LDC, , "\dashv" rotate=180] \arrow[phantom, from=CC, to=CB, , "\dashv" rotate=90]
\end{tikzcd}
\] 
of Quillen adjunctions.
\end{Proposition}

\begin{proof}
We set
\[
\mathcal{L}_*(X) \coloneqq \widehat{\Omega}_\iota^\flat\left( \widetilde{C}^c_*(X) \right)~,
\]
where the reduced cellular chain functor $\widetilde{C}^c_*(X)$ is given by the composition $\mathrm{Ker}(\epsilon) \circ C_*^c(-)$. We set
\[
\mathcal{R}_*(\mathfrak{g}) \coloneqq \mathrm{Hom}_{\mathsf{Strict}~u\mathcal{CC}_\infty\text{-}\mathsf{coalg}_{\bullet}}\left(C_*^c(\Delta^\bullet),\widehat{\mathrm{B}}_\iota^\flat(\mathfrak{g}) \oplus \kk \right)~.
\]
The square commutes by definition of the functors $\mathcal{L}_*$ and $\mathcal{R}_*$. Furthermore, they form a Quillen adjunction since the composition of Quillen adjunctions is still a Quillen adjunction.
\end{proof}

\begin{Proposition}\label{prop: modèle pour les composantes connexes}
Let $X$ be a pointed connected finite type nilpotent simplicial set. The unit
\[
\eta_X: X \qi \mathcal{R}_* \mathcal{L}_*(X)
\]
is a rational homotopy equivalence.
\end{Proposition}

\begin{proof}
There is a zig-zag of quasi-isomorphisms of $u\mathcal{CC}_\infty$-coalgebras:

\[
\widehat{\mathrm{B}}_\iota^\flat \widehat{\Omega}_\iota^\flat \widetilde{C}^c_*(X) \lqi C^c_*(X) \qi \widehat{\mathrm{B}}_\iota \widehat{\Omega}_\iota C^c_*(X)~.
\]

We get

\begin{align*}
\mathcal{R}_* \mathcal{L}_*(X) &= \mathrm{Hom}_{\mathsf{Strict}~u\mathcal{CC}_\infty\text{-}\mathsf{coalg}_{\bullet}}\left(C_*^c(\Delta^\bullet),\widehat{\mathrm{B}}_\iota^\flat \widehat{\Omega}_\iota^\flat \widetilde{C}^c_*(X) \oplus \kk \right) \\ 
&\simeq \mathrm{Hom}_{u\mathcal{CC}_\infty\text{-}\mathsf{coalg}_{\bullet}}\left(C_*^c(\Delta^\bullet),\widehat{\mathrm{B}}_\iota \widehat{\Omega}_\iota C^c_*(X) \right) \\
&\cong \mathcal{R}\mathcal{L}(X)~.
\end{align*}

We conclude by applying Theorem \ref{thm: modèles d'homotopie rationnel type fini}.
\end{proof}

\begin{Corollary}
Let $X$ be a finite type nilpotent simplicial set and let 
\[
X \cong \underset{\alpha \in \pi_0(X)}{\coprod} X^{\alpha}
\]
be its decomposition into connected components, where $\alpha: \{*\} \longrightarrow X$ are representatives of the equivalences classes in $\pi_0(X)$. We have 
\[
X \cong \underset{\alpha \in \pi_0(X)}{\coprod} X^{\alpha} \simeq \underset{\alpha \in \pi_0(X)}{\coprod} \mathcal{R}_* \mathcal{L}_*(X^{\alpha}) \simeq \mathcal{R}\mathcal{L}(X)~.
\]
\end{Corollary}

\begin{proof}
Follows directly from Proposition \ref{prop: modèle pour les composantes connexes}.
\end{proof}

\subsection{Characterization of rational models} In this subsection, we give a characterization of the essential homotopical image of finite type nilpotent spaces in the category of $u\mathcal{CC}_\infty$-coalgebras. 

\medskip

The essential homotopical image of the functor $\mathrm{A}_{\mathrm{PL}}$ of pointed connected finite type nilpotent spaces is given by augmented dg unital commutative algebras which admit a \textit{finite type Sullivan model}.

\begin{Definition}[Sullivan algebra]\label{def: Sullivan algebra}
A \textit{Sullivan algebra} $(\mathscr{S}(V),d)$ is the data of 
\begin{enumerate}
\item A graded vector space $V$ in negative degrees together with a filtration
\[
0 = V(-1) \subseteq V(0) \subseteq V(1) \subseteq \cdots \subseteq V~, 
\]
such that $\colim_{k} V(k) \cong V$.

\item A differential $d$ on the free unital commutative algebra $\mathscr{S}(V)$ such that 
\[
d(V(k)) \subseteq \mathscr{S}(V(k-1))~,
\]
for all $k \geq 0$. 
\end{enumerate}

It is of \textit{finite type} if $V$ is degree-wise finite dimensional.
\end{Definition}

\begin{theorem}[{\cite[Section 9]{BousfieldGugenheim}}]
There is a contravariant equivalence of $\infty$-categories

\[
\begin{tikzcd}[column sep=5pc,row sep=2.5pc]
\mathsf{sSet}_{*}^{\mathsf{f.t.c}~\mathsf{nilp}} \arrow[r, shift left=1.5ex, "\overline{\mathrm{A}_{\mathrm{PL}}}"{name=A}]
&\left(\mathsf{dg}~\mathcal{C}om\text{-}\mathsf{alg}_{\leq 1}^{\mathsf{f.t}~\mathsf{f.p}}\right)^{\mathsf{op}}~, \arrow[l, shift left=.75ex, "\left\langle - \right\rangle"{name=C}] \arrow[phantom, from=A, to=C, , "\dashv" rotate=-90]
\end{tikzcd}
\]

between the $\infty$-category of pointed connected finite type nilpotent simplicial sets and the $\infty$-category of dg unital commutative algebras with degree-wise finite dimensional homology concentrated in negative degrees which admit a finite type Sullivan algebra cofibrant replacement. 
\end{theorem}

\begin{Proposition}[{\cite[Theorem 2.3]{Berglund15}}]\label{prop: Sullivan algebra equal conilpotent L infinity}
The data of a Sullivan algebra $(\mathscr{S}(V),d)$ is equivalent to the data of a conilpotent $\mathcal{L}_\infty$-coalgebra $(V,\Delta_V,d_V)$ concentrated in negative degrees.
\end{Proposition}

\begin{proof}
We consider the following composition 
\[
\begin{tikzcd}
V \arrow[r,"d"]
&\displaystyle \bigoplus_{n \geq 1} \mathcal{C}om(n) \otimes V^{\otimes n} \arrow[r,rightarrowtail]
&\displaystyle \prod_{n \geq 1} \mathcal{C}om(n) \otimes V^{\otimes n}~,
\end{tikzcd}
\]
where we restrict the differential $d$ to the generators. It endows $V$ with a family of decomposition maps $l_n: V \longrightarrow V^{\odot n}$ for all $n \geq 1$. One can check that these maps satisfy the axioms of a $\mathcal{L}_\infty$-coalgebra. The filtration on $V$ makes it a conilpotent $\mathcal{L}_\infty$-coalgebra. On the other hand, we simply consider the Cobar construction
\[
\Omega_\pi: \mathsf{dg}~\mathcal{L}_\infty\text{-}\mathsf{coalg}^{\mathsf{conil}} \longrightarrow \mathsf{dg}~\mathcal{C}om\text{-}\mathsf{alg}
\]
applied to a conilpotent $\mathcal{L}_\infty$-coalgebra $\mathfrak{g}$ concentrated in negative degrees. The differential is build using the coalgebra structure and the coradical filtration gives the filtration on the generators. 
\end{proof}

\begin{Definition}[Finitely presented $\mathcal{CC}_\infty$-algebra]
Let $A$ be a $\mathcal{CC}_\infty$-algebra. It is \textit{finitely presented} if it admits a cofibrant replacement generated by a degree-wise finite dimensional conilpotent $\mathcal{L}_\infty$-coalgebra.
\end{Definition}

\begin{lemma}\label{lemma: uCC alg fin dim equiv uCom alg fin dim}
Let $\epsilon: \Omega \mathrm{B} \mathcal{C}om \qi \mathcal{C}om$ be the quasi-isomorphism of operads given by the counit of the operadic Bar-Cobar adjunction. The Quillen equivalence 

\[
\begin{tikzcd}[column sep=5pc,row sep=2.5pc]
\mathcal{CC}_\infty\text{-}\mathsf{alg} \arrow[r, shift left=1.5ex, "\mathrm{Ind}_\varepsilon"{name=A}]
&\mathsf{dg}~\mathcal{C}om\text{-}\mathsf{alg}~, \arrow[l, shift left=.75ex, "\mathrm{Res}_\varepsilon"{name=C}] \arrow[phantom, from=A, to=C, , "\dashv" rotate=-90]
\end{tikzcd}
\]

restricts to an equivalence of $\infty$-categories 

\[
\begin{tikzcd}[column sep=5pc,row sep=2.5pc]
\mathcal{CC}_\infty\text{-}\mathsf{alg}_{\leq 1}^{\mathsf{f.t}~\mathsf{f.p}} \arrow[r, shift left=1.5ex, "\mathrm{Ind}_\varepsilon"{name=A}]
&\mathsf{dg}~\mathcal{C}om\text{-}\mathsf{alg}_{\leq 1}^{\mathsf{f.t}~\mathsf{f.p}}~, \arrow[l, shift left=.75ex, "\mathrm{Res}_\varepsilon"{name=C}] \arrow[phantom, from=A, to=C, , "\dashv" rotate=-90]
\end{tikzcd}
\]

between the $\infty$-category of dg unital commutative degree-wise finite dimensional homology in negative degrees which admit a finite type Sullivan algebra cofibrant replacement and the $\infty$-category of $\mathcal{CC}_\infty$-algebras degree-wise finite dimensional homology in negative degrees which are finitely presented. 
\end{lemma}

\begin{proof}
One can check that the adjunction induces preserves the positive degree and degree-wise finite dimensional assumptions. A straightforward game playing with the Bar-Cobar adjunctions relative to $\pi$ and $\iota$ allows us to prove that this adjunction also preserves finitely presented objects, given Proposition \ref{prop: Sullivan algebra equal conilpotent L infinity}.
\end{proof}

\begin{Definition}[Finitely presented $\mathcal{CC}_\infty$-coalgebra]
Let $C$ be a $\mathcal{CC}_\infty$-coalgebra. It is \textit{finitely presented} if it admits a fibrant replacement generated by a degree-wise finite dimensional complete absolute $\mathcal{L}_\infty$-algebra.
\end{Definition}

\begin{Proposition}\label{prop: uCC coalg fin dim equiv uCC alg fin dim}
There is an equivalence of $\infty$-categories

\[
\begin{tikzcd}[column sep=5pc,row sep=2.5pc]
\left(\mathcal{CC}_\infty\text{-}\mathsf{alg}_{\leq 1}^{\mathsf{f.t}~\mathsf{f.p}}\right)^{\mathsf{op}} \arrow[r, shift left=1.5ex, "(-)^\circ"{name=A}]
&\mathcal{CC}_\infty\text{-}\mathsf{coalg}_{\geq 1}^{\mathsf{f.t}~\mathsf{f.p}}~. \arrow[l, shift left=.75ex, "(-)^*"{name=C}] \arrow[phantom, from=A, to=C, , "\dashv" rotate=-90]
\end{tikzcd}
\]
\end{Proposition}

\begin{proof}
By Proposition \ref{Prop: plongement pleinement fidele de Sweedler}, notice that this adjunction is indeed an equivalence between objects with degree-wise finite dimensional homology. It also preserves degrees since the (complete) Bar-Cobar resolution of a $\mathcal{CC}_\infty$-(co)algebra concentrated in positive degrees is still concentrated in positive degrees. By another adjunction game, one can show that it identifies finitely presented objects in both categories.
\end{proof}

\begin{Proposition}\label{prop: image essentielle dans le CC cog}
The adjunction 
\[
\begin{tikzcd}[column sep=5pc,row sep=2.5pc]
\mathsf{sSet}_{*}^{\mathsf{f.t.c}~\mathsf{nilp}} \arrow[r, shift left=1.5ex, "\widetilde{C}_*^c(-)"{name=A}]
&\mathcal{CC}_\infty\text{-}\mathsf{coalg}_{\geq 1}^{\mathsf{f.t}~\mathsf{f.p}}~, \arrow[l, shift left=.75ex, "\overline{\mathcal{R}}_*"{name=C}] \arrow[phantom, from=A, to=C, , "\dashv" rotate=-90]
\end{tikzcd}
\]
is an equivalence of $\infty$-categories.
\end{Proposition}

\begin{proof}
It follows directly from the propositions above. 
\end{proof}

\begin{theorem}
Let $C$ be a $u\mathcal{CC}_\infty$-coalgebra. There exists a finite type nilpotent simplicial set $X$ and a zig-zag of quasi-isomorphisms
\[
C_*^c(X) \lqi \cdot \qi \cdots \lqi \cdot \qi C
\]
of $u\mathcal{CC}_\infty$-coalgebras if and only if $C$ is quasi-isomorphic to a direct sum $\bigoplus_{i \in \mathrm{I}} C_i$ where $\mathrm{I}$ is a finite set and where $C_i$ are finitely presented pointed strict $u\mathcal{CC}_\infty$-coalgebras with degree-wise finite dimensional homology concentrated in non-negative degrees and such that $\mathrm{H}_0(C_i) = \kk$. 
\end{theorem}

\begin{proof}
The first implication is immediate by using the decomposition connected components
\[
C_*^c(X) \cong \bigoplus_{\alpha \in \pi_0(X)} C_*^c(X^\alpha)~.
\]
Suppose now that $C$ is quasi-isomorphic to some direct sum $\bigoplus_{i \in \mathrm{I}} C_i$. Then for each $\mathrm{H}_*(\widetilde{C}_i)$,  there is a pointed connected finite type nilpotent simplicial set $X^\alpha$ such that 
\[
\mathrm{H}_*(\widetilde{C}_i) \simeq \widetilde{C}_*^c(X^\alpha)~.
\]
Let 
\[
X \coloneqq \coprod_{\alpha} X^\alpha~,
\]
and since coproducts are also homotopy coproducts, we get that 
\[
C \simeq \bigoplus_{i} C_i \simeq \bigoplus_{\alpha} C_*^c(X^\alpha) \cong C_*^c(X)~.
\]
\end{proof}

\subsection{Models for mapping spaces}
In this section, we construct explicit rational models for mapping spaces, without any assumption on the source simplicial set. Furthermore, these models are relatively small, as they are constructed using the cellular chains on the source. 

\medskip

Recall that, if $X$ and $Y$ are simplicial sets, there is an explicit model for their mapping space given by 

\[
\mathrm{Map}(X,Y)_\bullet \coloneqq \mathrm{Hom}_{\mathsf{sSet}}(X \times \Delta^\bullet, Y)~,
\]
\vspace{0.25pc}

which forms a Kan complex when $Y$ is so.

\begin{lemma}\label{prop: les chaines sont lax monoidales en infini morphismes}
Let $X$ and $Y$ be two simplicial sets. There is an $\infty$-quasi-isomorphism 
\[
\psi_{X,Y}: C^c_*(X \times Y) \rightsquigarrow  C^c_*(X) \otimes C_*^c(Y)~,
\]
of $u\mathcal{CC}_\infty$-coalgebras which is natural in $X$ and $Y$.
\end{lemma}

\begin{proof}
For any simplicial sets $X,Y$, recall (see \cite{bookRHT}) that there is a natural quasi-isomorphism 
\[
\kappa: \mathrm{A}_{\mathrm{PL}}(X) \otimes \mathrm{A}_{\mathrm{PL}}(Y) \qi \mathrm{A}_{\mathrm{PL}}(X \times Y)
\]
of dg $u\mathcal{C}om$-algebras. This gives a quasi-isomorphism
\[
\mathrm{Res}_{\varepsilon}(\kappa): \mathrm{C}_{\mathrm{PL}}(X) \otimes \mathrm{C}_{\mathrm{PL}}(Y) \qi \mathrm{C}_{\mathrm{PL}}(X \times Y)
\]
of $u\mathcal{CC}_\infty$-algebras which is natural in $X,Y$, using the fact that the restriction functor $\mathrm{Res}_{\varepsilon}$ is strong monoidal and preserves all quasi-isomorphisms.

\medskip

Using the $\infty$-quasi-isomorphisms constructed in the proof of Proposition \ref{prop: equivalence cochains CC with Apl}, we construct 

\[
\begin{tikzcd}[column sep=4pc,row sep= 0pc]
C^*_c(X) \otimes C^*_c(Y) \arrow[r,"(i_\infty)_X \otimes (i_\infty)_Y ",rightsquigarrow]
&\mathrm{C}_{\mathrm{PL}}(X) \otimes \mathrm{C}_{\mathrm{PL}}(Y) \arrow[r,"\mathrm{Res}_{\varepsilon}(\kappa)"]
&\mathrm{C}_{\mathrm{PL}}(X \times Y) \arrow[r,"(p_\infty)_{X \times Y}",rightsquigarrow]
&C^*_c(X \times Y)~,
\end{tikzcd}
\]

where the tensor product of two $\infty$-morphisms is still an $\infty$-morphism, since the $2$-colored dg operad encoding $\infty$-morphisms of $u\mathcal{CC}_\infty$-algebras is cofibrant. This gives $\infty$-quasi-morphism of $u\mathcal{CC}_\infty$-algebras 

\[
C^*_c(X) \otimes C^*_c(Y) \rightsquigarrow C^*_c(X \times Y)~.
\]
\vspace{0.25pc}

which is natural in $X$ and $Y$. Equivalently, a weak-equivalence $\mathrm{B}_\iota(C^*_c(X) \otimes C^*_c(Y)) \qi \mathrm{B}_\iota(C^*_c(X \times Y))$. Now we take $X$ and $Y$ to be finite simplicial sets (finite total number of non-degenerate simplicies). By applying the linear dual functor of Theorem \ref{thm: homotopical magical square curved}, this gives a weak-equivalence

\[
\widehat{\Omega}_\iota(C^c_*(X) \otimes C_*^c(Y)) \qi \widehat{\Omega}_\iota(C_*^c(X \times Y))
\]
\vspace{0.25pc}

of complete curved absolute $\mathcal{L}_\infty$-algebras. Now suppose $X$ and $Y$ are arbitrary simplicial sets, one can write them as the filtered colimit of finite simplicial sets 
\[
X \cong \colim_\alpha X_\alpha \quad \text{and} \quad Y \cong \colim_\beta Y_\beta~.
\]
Now we have that 
\[
\widehat{\Omega}_\iota(C^c_*(X) \otimes C_*^c(Y)) \cong \colim_{\alpha,\beta} \widehat{\Omega}_\iota(C_*^c(X_{\alpha}) \otimes C_*^c(Y_{\beta}))~, 
\]
and 
\[
\widehat{\Omega}_\iota(C_*^c(X \times Y)) \cong \colim_{\alpha,\beta} \widehat{\Omega}_\iota(C_*^c(X_{\alpha} \times Y_{\beta}))~,
\]
since both bifunctors preserve filtered colimits in each variable. Both of this colimits are in fact homotopy colimits, hence there is a weak-equivalence 
\[
\widehat{\Omega}_\iota(C^c_*(X) \otimes C_*^c(Y)) \qi \widehat{\Omega}_\iota(C_*^c(X \times Y))~,
\]
which concludes the proof.
\end{proof}

\begin{theorem}\label{thm: vrai thm mapping spaces}
Let $\mathfrak{g}$ be a curved absolute $\mathcal{L}_\infty$-algebra and let $X$ be a simplicial set. There is a weak equivalence of Kan complexes

\[
\mathrm{Map}(X, \mathcal{R}(\mathfrak{g})) \simeq \mathcal{R}\left(\mathrm{hom}(C^c_*(X),\mathfrak{g})\right)~,
\]
\vspace{0.25pc}

which is natural in $X$ and in $\mathfrak{g}$, where $\mathrm{hom}(C^c_*(X),\mathfrak{g})$ denotes the convolution curved absolute $\mathcal{L}_\infty$-algebra. 
\end{theorem}

\begin{proof}
There is an isomorphism

\[
\mathrm{Map}(X,\mathcal{R}(\mathfrak{g}))_\bullet \coloneqq \mathrm{Hom}_{\mathsf{sSet}}(X \times \Delta^\bullet, \mathcal{R}(\mathfrak{g})) \cong \mathrm{Hom}_{\mathsf{sSet}}( C^c_*(X \times \Delta^\bullet), \widehat{\mathrm{B}}_\iota(\mathfrak{g}))~.
\]
\vspace{0.25pc}

We can pre-compose by the $\infty$-quasi-isomorphism of Lemma \ref{prop: les chaines sont lax monoidales en infini morphismes}, giving a weak-equivalence of simplicial sets

\[
\mathrm{Hom}_{u\mathcal{CC}_\infty}(C^c_*(X \times \Delta^\bullet), \widehat{\mathrm{B}}_\iota(\mathfrak{g})) \qi \mathrm{Hom}_{u\mathcal{CC}_\infty\text{-}\mathsf{cog}}( C^c_*(X)  \otimes C^c_*(\Delta^\bullet), \widehat{\mathrm{B}}_\iota(\mathfrak{g}))~,
\]
\vspace{0.25pc}

since both $C^c_*(X \times \Delta^\bullet)$ and $C^c_*(X)  \otimes C^c_*(\Delta^\bullet)$ are cofibrant. Let's compute this last simplicial set:

\begin{align*}
\mathrm{Hom}_{u\mathcal{CC}_\infty\text{-}\mathsf{cog}}( C^c_*(X)  \otimes C^c_*(\Delta^\bullet), \widehat{\mathrm{B}}_\iota(\mathfrak{g})) 
&\cong \mathrm{Hom}_{u\mathcal{CC}_\infty\text{-}\mathsf{cog}} \left( C^c_*(\Delta^\bullet), \left\{ C^c_*(X), \widehat{\mathrm{B}}_\iota(\mathfrak{g}) \right\} \right) \\
&\cong \mathrm{Hom}_{u\mathcal{CC}_\infty\text{-}\mathsf{cog}} \left( C^c_*(\Delta^\bullet), \widehat{\mathrm{B}}_\iota \left(\mathrm{hom}( C^c_*(X), \mathfrak{g}) \right) \right) \\
& \cong \mathrm{Hom}_{\mathsf{curv}~\mathsf{abs}~\mathcal{L}_\infty\text{-}\mathsf{alg}} \left( \widehat{\Omega}_\iota(C^c_*(\Delta^\bullet)), \mathrm{hom}( C^c_*(X), \mathfrak{g}) \right) \\
& \cong \mathcal{R}\left(\mathrm{hom}( C^c_*(X), \mathfrak{g}) \right)~.
\end{align*}
\end{proof}

\begin{Corollary}\label{Cor: True mapping spaces.}
Let $X$ be a simplicial set and let $Y$ be a finite type nilpotent simplicial set. There is a weak equivalence of simplicial sets 

\[
\mathrm{Map}(X, Y_\mathbb{Q}) \simeq \mathcal{R}\left(\mathrm{hom}(C^c_*(X), \mathcal{L}(Y))\right)~,
\]

where $Y_\mathbb{Q}$ denotes the $\mathbb{Q}$-localization of $Y$. Therefore, given a map $f: X \longrightarrow Y_\mathbb{Q}$, there is an isomorphism

\[
\pi_n\left(\mathrm{Map}(X, Y_\mathbb{Q}),f\right) \cong \mathrm{H}_n^{C_*^c(f)}\left(\mathrm{hom}(C^c_*(X), \mathcal{L}(Y))\right)~,
\]

of groups for for $n \geq 1$. 
\end{Corollary}

\begin{proof}
There is a weak equivalence of simplicial sets 
\[
Y_\mathbb{Q} \simeq \mathcal{RL}(Y)~.
\]
It induces equivalences
\[
\mathrm{Map}(X, Y_\mathbb{Q}) \simeq \mathrm{Map}(X, \mathcal{RL}(Y)) \simeq \mathcal{R}\left(\mathrm{hom}(C^c_*(X), \mathcal{L}(Y))\right)~.
\]
\vspace{0.1pc}

Finally, the second statement follows directly from Theorem \ref{thm: Berglund isomorphism}. 
\end{proof}

\begin{Remark}
Notice the following points about Corollary \ref{Cor: True mapping spaces.}.

\begin{enumerate}
\medskip

\item There is no finiteness hypothesis on the simplicial set $X$. 

\medskip

\item So far, the models for mapping spaces in the literature for mapping spaces are given in terms of complete tensor products
\[
A_X ~\widehat{\otimes}~ \mathfrak{g}_Y~,
\]
of a dg unital commutative algebra model $A_X$ of $X$ and a complete $\mathcal{L}_\infty$-algebra model $\mathfrak{g}_Y$. Our model in comparison is much smaller as it only uses the cellular chains on $X$, see for instance \cite{Berglund15, BuijMapping, LazarevMapping}.
\end{enumerate}
\end{Remark}

\begin{Corollary}
Let $X$ be a finite type nilpotent simplicial set. The space of rational endomorphisms of $X$ is weak-equivalent to 

\[
\mathrm{End}^h_\mathbb{Q}(X) \qi \mathcal{R}\left(\mathrm{hom}(C^c_*(X), \mathcal{L}(X))\right)~.
\]
\end{Corollary}

\begin{proof}
This is immediate from the previous corollary. 
\end{proof}

\section{Deformation theory}
In this section, we consider different applications of the theory developed so far. We define convolution curved absolute $\mathcal{L}_\infty$-algebras which encode $\infty$-morphisms between algebras over an operad as their Maurer-Cartan elements. We explain how this method can lead to simplicial enrichments for algebras over operads. Finally, we explore the way in which curved absolute $\mathcal{L}_\infty$-algebras encode non-pointed deformation problems.

\subsection{Convolution curved absolute $\mathcal{L}_\infty$-algebras and $\infty$-morphisms of algebras} Let $\mathcal{P}$ denote a dg operad and let $\mathcal{C}$ denote a conilpotent curved cooperad. 

\begin{Proposition}\label{prop: curved twisting morphisms}
There is a bijection 

\[
\mathrm{Tw}(\mathcal{C},\mathcal{P}) \cong \mathrm{Hom}_{\mathsf{curv}~\mathsf{ab}~\mathsf{pOp}}\left(\widehat{\Omega}^{\mathrm{s.c}} u\mathcal{C}om^*, \mathcal{H}om(\mathcal{C},\mathcal{P})\right)~,
\]

between morphisms of curved absolute partial operads $\phi_\alpha : \widehat{\Omega}^{\mathrm{s.c}} u\mathcal{C}om^* \longrightarrow \mathcal{H}om(\mathcal{C},\mathcal{P})$ and curved twisting morphisms $\alpha: \mathcal{C} \longrightarrow \mathcal{P}$ such that $\alpha(1)$ is the zero morphism.
\end{Proposition}

\begin{proof}
First, notice that the convolution curved partial operad of Definition \ref{def: convolution curved partial operad} is in fact a curved absolute partial operad in the sense of Definition \ref{def: absolute partial operads}. Indeed, one can check that any infinite sum of compositions in it has a well-defined image in it since any such sum is \textit{locally finite} because of the conilpotency of $\mathcal{C}$. There is an isomorphism 

\[
\mathcal{H}om\left(u\mathcal{C}om^*,\mathcal{H}om(\mathcal{C},\mathcal{P})\right) \cong \mathcal{H}om(\mathcal{C},\mathcal{P})
\]
\vspace{0.2pc}

of curved absolute partial operads since $u\mathcal{C}om^*(n) = \kk$ for all $n \geq 0$. Therefore,

\[
\mathrm{Hom}_{\mathsf{curv}~\mathsf{ab}~\mathsf{pOp}}\left(\widehat{\Omega} u\mathcal{C}om^*, \mathcal{H}om(\mathcal{C},\mathcal{P})\right) \cong \mathrm{Tw}\left(u\mathcal{C}om^*,\mathcal{H}om(\mathcal{C},\mathcal{P})\right) \cong \mathrm{Tw}(\mathcal{C},\mathcal{P})~,
\]
\vspace{0.2pc}

using the curved operadic Bar-Cobar adjunction of Section \ref{Section: Constructions Bar-Cobar operadiques}. Now one can check that $\widehat{\Omega}^{\mathrm{s.c}} u\mathcal{C}om^*$ represents exactly those curved twisting morphisms that have a trivial arity one component.
\end{proof}

\begin{Remark}
For an analogue statement concerning convolution $\mathcal{L}_\infty$-algebra structures, see the constructions in \cite{RobertNicoudWierstra17}.
\end{Remark}

Let $(A,\allowbreak \gamma_A,d_A)$ be a dg $\mathcal{P}$-algebra and let $(C,\Delta_C,d_C)$ be a curved $\mathcal{C}$-coalgebra. The graded module $\mathrm{hom}_{gr}(C,A)$ is naturally a pdg module endowed with the pre-differential $\partial(f) \coloneqq d_A \circ f - f \circ d_C$.

\begin{Proposition}
Let $(A,\gamma_A,d_A)$ be a dg $\mathcal{P}$-algebra and let $(C,\Delta_C,d_C)$ be a curved $\mathcal{C}$-coalgebra. The pdg module 
\[
\left(\mathrm{hom}_{\mathsf{gr}\text{-}\mathsf{mod}}(C,A),\partial \right)
\]
can be endowed with a curved $\mathcal{L}_\infty$-algebra structure is given by 
\[
\begin{tikzcd}[column sep=2.5pc,row sep=0.5pc]
\gamma_{\mathrm{hom}(C,A)}: \displaystyle \bigoplus_{n \geq 0} \widehat{\Omega} u\mathcal{C}om^*(n) \otimes_{\mathbb{S}_n} \mathrm{hom}(C,A)^{\otimes n} \arrow[r]
&\mathrm{hom}(C,A) \\
c_n (f_1 \otimes \cdots \otimes f_n) \arrow[r,mapsto]
&\displaystyle \gamma_A \circ \left[\phi_\alpha(c_n) \otimes (f_1 \otimes \cdots \otimes f_n) \right] \circ \Delta_C~,
\end{tikzcd}
\]
where $c_n$ denotes the $n$-corolla.
\end{Proposition}

\begin{proof}
This follows directly from Proposition \ref{prop: curved twisting morphisms} by considering the composition 
\[
\begin{tikzcd}[column sep=4pc,row sep=0.5pc]
\widehat{\Omega} u\mathcal{C}om^* \arrow[r," \phi_\alpha"] 
&\mathcal{H}om(\mathcal{C},\mathcal{P}) \arrow[r,"\lambda_{\mathrm{hom}(C,A)} "]
&\mathrm{End}_{\mathrm{hom}(C,A)}~,
\end{tikzcd}
\]
of morphisms of curved operads, where $\lambda_{\mathrm{hom}(C,A)}$ is the natural curved algebra structure of $\mathrm{hom}(C,A)$ over $\mathcal{H}om(\mathcal{C},\mathcal{P})$.
\end{proof}

\begin{Remark}[Local nilpotency]
One can notice that the Maurer-Cartan equation is always well-defined without the need of imposing any filtration on the space $\mathrm{hom}(C,A)$. Indeed, let $f: C \longrightarrow A$ be a degree zero map, then the sum
\[
\gamma_A \circ \left[\sum_{\substack{n \geq 0 \\ n \neq 1}} \phi_\alpha(c_n) \otimes f^{\otimes n} \right] \circ \Delta_C(c) + \partial(f)(c)
\]
is well defined for any element $c$ in $C$, since there are only a finite number of non-zero terms in $\Delta_C(c)$. Filtrations introduced to make this type of sum converge like in \cite{DolgushevHoffnungRogers14} are redundant. 
\end{Remark}

\begin{Proposition}
The extension of $\gamma_{\mathrm{hom}(C,A)}$ given by 
\[
\begin{tikzcd}[column sep=1pc,row sep=-0.5pc]
\displaystyle \prod_{n \geq 0} \widehat{\Omega} u\mathcal{C}om^*(n)~\widehat{\otimes}_{\mathbb{S}_n}~\mathrm{hom}(C,A)^{\otimes n} \arrow[r]
&\mathrm{hom}(C,A) \\
\displaystyle \sum_{\substack{n\geq 0 \\ \omega \geq 0}} \sum_{\tau \in \mathrm{CRT}_\omega^n } \lambda_\tau \tau(f_1 \otimes \cdots \otimes f_n) \arrow[r,mapsto]
&\displaystyle \gamma_A \circ \left[\sum_{\substack{n \geq 0 \\ \omega \geq 0}} \sum_{\tau \in \mathrm{CRT}_\omega^n } \phi_\alpha(\tau) \otimes (f_1 \otimes \cdots \otimes f_n) \right] \circ \Delta_C~,
\end{tikzcd}
\]
defines a structure of curved absolute $\mathcal{L}_\infty$-algebra, denoted by $\mathrm{hom}(C,A)$.
\end{Proposition}

\begin{proof}
One can check by hand that this formula satisfies conditions \ref{pdg condition}, \ref{associativity condition}, and \ref{curved condition}.
\end{proof}

\begin{Remark}
This convolution curved absolute $\mathcal{L}_\infty$-algebra structure coincides with the one constructed in \cite{grignou2022mapping} between a dg $\Omega \mathcal{C}$-algebra and a curved $\mathcal{C}$-coalgebra, using the Hopf comodule structure of the cofibrant dg operad $\Omega \mathcal{C}$.
\end{Remark}

\begin{Proposition}
Let $A$ and $B$ be two dg $\mathcal{P}$-algebras. The simplicial set $\mathcal{R}(\mathrm{hom}(\mathrm{B}_\alpha A, B))$ has as $0$-simplicies the set of $\infty_\alpha$-morphisms between $A$ and $B$.
\end{Proposition}

\begin{proof}
Maurer-Cartan elements correspond by definition to curved twisting morphisms $\mathrm{B}_\alpha A$ between and $B$, which are in bijection with morphisms of curved $\mathcal{C}$-coalgebras between $\mathrm{B}_\alpha A$ and $\mathrm{B}_\alpha B$.
\end{proof}

Therefore $1$-simplicies in $\mathcal{R}(\mathrm{hom}(\mathrm{B}_\alpha A, B))$ induce a notion of \textit{homotopies} between these $\infty_\alpha$-morphism.

\medskip

In the work of \cite{DolgushevHoffnungRogers14}, the authors construct convolution (complete) $\mathcal{L}_\infty$-algebras between a dg $\Omega \mathcal{C}$-algebras and dg $\mathcal{C}$-coalgebras. In \textit{loc.cit.} both $\mathcal{C}$ and $\Omega \mathcal{C}$ are \textit{reduced}, hence they only encode non-unital/non-curved types of algebraic objects. They show that the composition of $\infty_\iota$-morphisms

\[
\circ: \mathrm{hom}(\mathrm{B}_\iota A,  B) \oplus \mathrm{hom}(\mathrm{B}_\iota B,  C) \longrightarrow \mathrm{hom}(\mathrm{B}_\iota A,  C)
\]
\vspace{0.1pc}

can be extended into an $\infty_\iota$-morphism of $\mathcal{L}_\infty$-algebras. This allows them to enrich simplicially the category of dg $\Omega \mathcal{P}$-algebras with $\infty_\iota$-morphisms since Hinich's functor $\mathrm{MC}(- \otimes \Omega_\bullet)$ is functorial with respect to $\infty_\iota$-morphism of $\mathcal{L}_\infty$-algebras. 

\medskip

In order to construct a simplicial enrichment of dg $\Omega \mathcal{P}$-algebras with $\infty_\iota$-morphisms in full generality, one should be able to show that the composition defines a morphism

\[
\circ: \widehat{\mathrm{B}}_\iota(\mathrm{hom}(\mathrm{B}_\iota A,  B)) \times  \widehat{\mathrm{B}}_\iota(\mathrm{hom}(\mathrm{B}_\iota B,  C)) \longrightarrow \widehat{\mathrm{B}}_\iota(\mathrm{hom}(\mathrm{B}_\iota A,  C))
\]
\vspace{0.1pc}

of $u\mathcal{CC}_\infty$-coalgebras. This goes beyond the scope of this thesis and will be considered in a sequel work.

\begin{Remark}[The dual case of $\infty$-morphism of coalgebras]
In \cite{grignou2022mapping}, Brice Le Grignou constructs a convolution curved absolute $\mathcal{L}_\infty$-algebra from a dg $\Omega\mathcal{C}$-coalgebra and a curved $\mathcal{C}$-algebra. Integrating these convolution algebras provides us with the right set of $\infty_\iota$-morphisms between two coalgebras. This method also provides us with a way in which one can simplicially enrich the category of dg $\Omega\mathcal{C}$-coalgebras with $\infty_\iota$-morphisms.
\end{Remark}

\begin{Remark}[Extension to dg properads]
Let $\Omega \mathcal{C}$ be a properad. Consider the notion of $\infty$-morphism for gebras over $\Omega \mathcal{C}$ as defined in  \cite{hoffbeck2019properadic}. It was privately communicated to us that the convolution algebra that governs these $\infty$-morphisms does not appear to carry any sort of filtration that would make the Maurer-Cartant equation converge. As a consequence, the authors plant to apply the present integration theory of absolute $\mathcal{L}_\infty$-algebras in order to settle a suitable simplicial enrichment for the category of $\Omega \mathcal{C}$-gebras with $\infty$-morphisms.
\end{Remark}

\newcommand\htimes{\stackrel{\mathclap{\normalfont\mbox{\tiny{h}}}}{\times}}

\subsection{The formal geometry of Maurer-Cartan spaces}
Let $A$ be a dg $u\mathcal{C}om$-algebras, viewed equivalently as a \textit{derived affine stacks}. Its \textit{functor of points} is given by 

\[
\begin{tikzcd}[column sep=1pc,row sep=0.5pc]
\mathrm{Spec}(A)(-): \mathsf{dg}~u\mathcal{C}om\text{-}\mathsf{alg}_{\geq 0} \arrow[r]
&\mathsf{sSet} \\
B \arrow[r,mapsto]
&\mathrm{Spec}(A)(B)_\bullet \coloneqq \mathbb{R}\mathrm{Hom}_{\mathsf{dg}~u\mathcal{C}om\text{-}\mathsf{alg}}(A, B \otimes \Omega_\bullet)~,
\end{tikzcd}
\]
where $\Omega_\bullet$ is again the simplicial Sullivan algebra and where $B$ is a dg $u\mathcal{C}om$-algebra concentrated in non-negative homological degrees, that is, a \textit{derived affine scheme}. The simplicial set $\mathrm{Spec}(A)(B)_\bullet$ is called the $B$-points of $A$.

\medskip

The first goal of this subsection is to construct a functor from dg $u\mathcal{C}om$-algebras to curved absolute $\mathcal{L}_\infty$-algebras that recovers the "formal geometry" of $A$. Let $x: A \longrightarrow \kk$ be a $\kk$-point of $A$, it provides an $\kk$-\textit{augmentation} of $A$. The standard way of looking at \textit{infinitesimal thickenings} of $x$ inside $A$ has been to test $A$ against $\kk$-augmented dg Artinian algebras. 

\begin{Definition}[$\kk$-augmented dg Artinian algebra]
Let $R$ be a dg $u\mathcal{C}om$-algebra with homology concentrated in non-negative homological degrees. It is a $\kk$\textit{-augmented dg Artinian algebra} if it satisfies the following conditions.

\medskip

\begin{enumerate}
\item Its homology is degree-wise finite dimensional and concentrated in a finite number of degrees.

\medskip

\item There is an unique augmentation morphism $p: R \longrightarrow \kk$ (hence $R$ is local).

\medskip

\item Let $\overline{R}$ be the non-unital dg commutative algebra given by the kernel of $p$. Then $\overline{R}$ is nilpotent.
\end{enumerate}
\end{Definition}

This notion allows us to define the notion of a $\kk$-pointed formal moduli problem.

\begin{Definition}[$\kk$-pointed formal moduli problem]
Let 
\[
F: \mathsf{dg}~\mathsf{Art}\text{-}\mathsf{alg}^{\kk\text{-}\mathsf{aug}}_{\geq 0} \longrightarrow \mathsf{sSet}
\]

be a functor. It is a $\kk$-\textit{pointed formal moduli problem} if it satisfies the following conditions

\medskip

\begin{enumerate}
\item We have that $F(\kk) \simeq \{*\}$.

\medskip

\item The functor $F$ sends quasi-isomorphisms to weak-equivalences.

\medskip

\item The functor $F$ preserves homotopy pullbacks of the of dg Artinian algebras $X,Y$ and $Z$ 
\[
\begin{tikzcd}
X \htimes_Z Y \arrow[r] \arrow[d]
&Y \arrow[d,"\pi_1"] \\
X \arrow[r,"\pi_2 ",swap]
&Z \\
\end{tikzcd}
\]

such that $\pi_1$ and $\pi_2$ are surjections on the zeroth homology groups $\mathrm{H}_0$. 
\end{enumerate}
\end{Definition}

\begin{Example}
Any $\kk$-augmented derived affine stack $A$ defines such a pointed deformation problem by considering 

\[
\mathrm{Spec}^*(A)(-)_\bullet = \mathbb{R}\mathrm{Hom}_{\mathsf{dg}~u\mathcal{C}om\text{-}\mathsf{alg}^{\kk\text{-}\mathsf{aug}}}(A, - \otimes \Omega_\bullet)~,
\]

where morphisms of $\kk$-augmented dg commutative algebras are required to preserve the augmentation. This functor preserves homotopy limits and sends quasi-isomorphisms to weak homotopy equivalences, hence it defines a pointed formal moduli problem. 
\end{Example}

The $\infty$-category of $\kk$-pointed formal moduli problems is equivalent to the $\infty$-category of dg Lie algebras over $\kk$, or for that matter, to the $\infty$-category of $\mathcal{L}_\infty$-algebras, both of them constructed by consider quasi-isomorphisms as weak-equivalences. This was shown independently by J. Lurie in \cite{Lurie11} by working directly with $\infty$-categories and by J.P Pridham \cite{Pridham10} using model categories. 

\medskip

\textbf{Non-pointed version.} Heuristically speaking, curved absolute $\mathcal{L}_\infty$-algebras should correspond to a \textit{non-pointed} version of formal moduli problems. A first idea into what this notion might be is the following: in a pointed formal moduli problem, the point around which one considers the infinitesimal thickening is already specified in advance. In the non-pointed version, one should remove the $\kk$-augmented hypothesis. 

\medskip

In order to understand what these non-pointed deformation problems can correspond to, we build a curved absolute $\mathcal{L}_\infty$-algebra $\mathfrak{g}_A$ from a derived affine stack $A$. The goal is to understand how much "geometrical information" about $A$ one can recover from its geometrical model $\mathfrak{g}_A$. We consider the following commutative square given by Theorem \ref{thm: homotopical magical square curved}:

\[
\begin{tikzcd}[column sep=5pc,row sep=5pc]
\left(\mathsf{dg}~u\mathcal{C}om\text{-}\mathsf{alg}\right)^{\mathsf{op}} \arrow[r,"\mathrm{B}_\pi^{\mathsf{op}}"{name=B},shift left=1.1ex] \arrow[d,"(-)^\circ "{name=SD},shift left=1.1ex ]
&\left(\mathsf{curv}~\mathcal{L}_\infty\text{-}\mathsf{coalg}^{\mathsf{conil}}\right)^{\mathsf{op}} \arrow[d,"(-)^*"{name=LDC},shift left=1.1ex ] \arrow[l,"\Omega_\pi^{\mathsf{op}}"{name=C},,shift left=1.1ex]  \\
\mathsf{dg}~u\mathcal{C}om\text{-}\mathsf{coalg} \arrow[r,"\widehat{\Omega}_\pi "{name=CC},shift left=1.1ex]  \arrow[u,"(-)^*"{name=LD},shift left=1.1ex ]
&\mathsf{curv}~\mathsf{abs}~\mathcal{L}_\infty\text{-}\mathsf{alg}^{\mathsf{comp}}~, \arrow[l,"\widehat{\mathrm{B}}_\pi"{name=CB},shift left=1.1ex] \arrow[u,"(-)^\vee"{name=TD},shift left=1.1ex] \arrow[phantom, from=SD, to=LD, , "\dashv" rotate=0] \arrow[phantom, from=C, to=B, , "\dashv" rotate=-90]\arrow[phantom, from=TD, to=LDC, , "\dashv" rotate=0] \arrow[phantom, from=CC, to=CB, , "\dashv" rotate=-90]
\end{tikzcd}
\] 

which relates the category of dg $u\mathcal{C}om$-algebras to the category of curved absolute $\mathcal{L}_\infty$-algebras.

\begin{Remark}\label{Rmk: pas de top model}
The category of dg $u\mathcal{C}om$-coalgebras does not admit a model category structure transferred from dg modules, where weak-equivalences are quasi-isomorphisms and where cofibrations are degree-wise monomorphisms. Therefore the left adjunction on this square is not a Quillen adjunction in this case, see \cite[Section 9]{grignoulejay18} for more details.
\end{Remark}

\begin{Definition}[Geometrical model]
Let $A$ be a dg $u\mathcal{C}om$-algebra. Its \textit{geometrical model} $\mathfrak{g}_A$ is the complete curved absolute $\mathcal{L}_\infty$-algebra given by $(\mathrm{B}_\pi A)^*$. It defines a functor 
\[
\mathfrak{g}_{(-)} \coloneqq (\mathrm{B}_\pi(-))^*: \left(\mathsf{dg}~u\mathcal{C}om\text{-}\mathsf{alg}\right)^{\mathsf{op}} \longrightarrow \mathsf{curv}~\mathsf{abs}~\mathcal{L}_\infty\text{-}\mathsf{alg}^{\mathsf{comp}}~.
\]
\end{Definition}

\begin{Remark}
This approach is close to the one of \cite{calaqueformal}, where the authors consider the linear dual of the Bar construction relative to $\kappa$ of a dg $\mathcal{P}$-algebra. Nevertheless, they view it as a dg $\mathcal{P}^{\hspace{1pt}!}$-algebra instead of viewing this linear dual as dg algebra over the cooperad $\mathcal{P}^{\hspace{1pt}\ac}$. They use the linear dual of the Bar construction in order to show an equivalence between the $\infty$-category of pointed formal moduli problems defined using dg Artinian $\mathcal{P}$-algebras and the $\infty$-category of dg $\mathcal{P}^{\hspace{1pt}\ac}$-algebras, under certain hypothesis on the operad $\mathcal{P}$.
\end{Remark}

\begin{Proposition}
The geometrical model functor $\mathfrak{g}_{(-)}$ sends quasi-isomorphisms to weak-equivalences. 
\end{Proposition}

\begin{proof}
Let $f: A \qi B$ be a quasi-isomorphism of dg $u\mathcal{C}om$-algebras, then $\mathrm{B}_\pi(f): \mathrm{B}_\pi A \qi \mathrm{B}_\pi B$ is a weak-equivalence of conilpotent curved $\mathcal{L}_\infty$-coalgebras, since every dg $u\mathcal{C}om$-algebra is fibrant. The linear dual $(-)^*$ preserves all weak-equivalences since it is a right Quillen functor and since every conilpotent curved $\mathcal{L}_\infty$-coalgebra is cofibrant (fibrant in the opposite category). 
\end{proof}

\begin{Definition}[dg Artinian algebra]
Let $A$ be a dg $u\mathcal{C}om$-algebra in non-negative homological degrees. It is a \textit{dg Artinian algebra} if its homology is degree-wise finite dimensional.
\end{Definition}

\begin{Example}
Consider a dg Artinian algebra $A$ concentrated in degree zero. Since $A$ is a finite dimensional $\kk$-algebra, it can be written as a product 
\[
A \cong \prod_{i = 0}^n A_i~, 
\]
where $A_i$ are local Artinian algebras. Thus one can consider the following examples

\medskip

\begin{enumerate}
\item the algebra $\kk^n$, which gives a collection of $n$ points, 

\medskip

\item any finite field extension $\mathbb{L}$ of $\kk$,

\medskip

\item classical thickenings like $\kk[t]/(t^n)$. 
\end{enumerate}
\end{Example}

\begin{lemma}\label{lemma: equivalence dg Artinien avec les uCC}
There is a an equivalence of $\infty$-categories between the opposite $\infty$-category of $u\mathcal{CC}_\infty$-coalgebras with finite dimensional homology concentrated in non-positive degrees and the $\infty$-category of dg Artinian algebras.
\end{lemma}

\begin{proof}
There are equivalences 
\[
\begin{tikzcd}[column sep=5pc,row sep=2.5pc]
\left(u\mathcal{CC}_\infty\text{-}\mathsf{coalg}_{\leq 0}^{\mathsf{f.d}}\right)^{\mathsf{op}}  \arrow[r, shift left=1.5ex, "(-)^*"{name=B}]
&u\mathcal{CC}_\infty\text{-}\mathsf{alg}_{\geq 0}^{\mathsf{f.d}} \arrow[r, shift left=1.5ex, "\mathrm{Ind}_\varepsilon"{name=A}] \arrow[l, shift left=.75ex, "(-)^\circ"{name=D}]
&\mathsf{dg}~\mathcal{C}om\text{-}\mathsf{alg}_{\geq 0}^{\mathsf{f.d}}~, \arrow[l, shift left=.75ex, "\mathrm{Res}_\varepsilon"{name=C}] \arrow[phantom, from=A, to=C, , "\dashv" rotate=-90] \arrow[phantom, from=B, to=D, , "\dashv" rotate=90]
\end{tikzcd}
\]
where the first equivalence is given by Proposition \ref{Prop: plongement pleinement fidele de Sweedler} and the second equivalence is induced by the quasi-isomorphism of dg operads $\varepsilon: \Omega \mathrm{B}u\mathcal{C}om \qi u\mathcal{C}om$, which respects the degrees of the homology.
\end{proof}

\begin{theorem}\label{thm: modèle géométrique sur les points Artiniens}
Let $A$ be a dg $u\mathcal{C}om$-algebra. Let $B$ be a dg Artinian algebra. There is a weak homotopy equivalence of simplicial sets
\[
\mathrm{Spec}(A)(B) \simeq \mathcal{R}(\mathrm{hom}( (\mathrm{Res}_\varepsilon B)^\circ, \mathfrak{g}_A))~,
\]
which is natural in $B$ and in $A$. 
\end{theorem}

\begin{proof}
We start by computing the derived mapping space with the classical Bar-Cobar resolution

\[
\mathrm{Spec}(A)(B) \simeq \mathrm{Hom}_{\mathsf{dg}~u\mathcal{C}om\text{-}\mathsf{alg}}(\Omega_\pi \mathrm{B}_\pi A, B \otimes \Omega_\bullet)~,
\]
\vspace{0.1pc}

since $\Omega_\pi \mathrm{B}_\pi A$ is a cofibrant resolution of $A$. Using the same methods as in the proof of Proposition \ref{prop: equivalence cochains CC with Apl}, we show that

\[
\mathrm{Hom}_{\mathsf{dg}~u\mathcal{C}om\text{-}\mathsf{alg}}(\Omega_\pi \mathrm{B}_\pi A, B \otimes \Omega_\bullet) \simeq \mathrm{Hom}_{u\mathcal{CC}_\infty\text{-}\mathsf{alg}}(\Omega_\iota \mathrm{B}_\iota \mathrm{Res}_\epsilon A, \mathrm{Res}_\epsilon(B) \otimes C^*_c(\Delta^\bullet))~.
\]
\vspace{0.1pc}

We consider now the following quasi-isomorphism of $u\mathcal{CC}_\infty$-algebras 

\[
\begin{tikzcd}[column sep=4pc,row sep=0pc]
\mathrm{Res}_\epsilon(B) \otimes C^*_c(\Delta^\bullet) \arrow[r,"\eta_{\mathrm{Res}_\epsilon(B)} \otimes \mathrm{id}"]
&((\mathrm{Res}_\varepsilon B)^\circ)^* \otimes C^*_c(\Delta^\bullet) \arrow[r,"\Gamma"]
&\left((\mathrm{Res}_\varepsilon B)^\circ \otimes C_*^c(\Delta^\bullet))\right)^*~,
\end{tikzcd}
\]

where $\eta_{\mathrm{Res}_\epsilon(B)}$ is the unit of the $(-)^* \dashv (-)^\circ$ adjunction, which is a quasi-isomorphism since $B$ has degree-wise finite dimensional homology, and where $\Gamma$ is the lax monoidal structure of the linear dual functor $(-)^*$ with respect to the tensor product, which is a quasi-isomorphism since both objects have degree-wise finite dimensional homology. This quasi-isomorphism gives a direct natural weak-equivalence of simplicial sets 

\[
\begin{tikzcd}[column sep=0pc,row sep=2pc]
\mathrm{Hom}_{u\mathcal{CC}_\infty\text{-}\mathsf{alg}}(\Omega_\iota \mathrm{B}_\iota \mathrm{Res}_\epsilon A, \mathrm{Res}_\epsilon(B) \otimes C^*_c(\Delta^\bullet)) \arrow[d,"\simeq"] \\
\mathrm{Hom}_{u\mathcal{CC}_\infty\text{-}\mathsf{alg}}(\Omega_\iota \mathrm{B}_\iota \mathrm{Res}_\epsilon A, \left((\mathrm{Res}_\varepsilon B)^\circ \otimes C_*^c(\Delta^\bullet))\right)^*)
\end{tikzcd}
\]

There is an isomorphism of simplicial sets 

\[
\begin{tikzcd}[column sep=0pc,row sep=2pc]
\mathrm{Hom}_{u\mathcal{CC}_\infty\text{-}\mathsf{alg}}(\Omega_\iota \mathrm{B}_\iota \mathrm{Res}_\varepsilon A, \left((\mathrm{Res}_\varepsilon B)^\circ \otimes C_*^c(\Delta^\bullet)\right)^*) \arrow[d,"\cong"] \\
\mathrm{Hom}_{u\mathcal{CC}_\infty\text{-}\mathsf{coalg}}((\mathrm{Res}_\varepsilon B)^\circ \otimes C_*^c(\Delta^\bullet),\left(\Omega_\iota \mathrm{B}_\iota \mathrm{Res}_\varepsilon A\right)^\circ)
\end{tikzcd}
\]

induced by the adjunction $(-)^\circ \dashv (-)^*$. We compute that 

\[
\left(\Omega_\iota \mathrm{B}_\iota \mathrm{Res}_\varepsilon A\right)^\circ \cong \widehat{\mathrm{B}}_\iota (\mathrm{B}_\iota \mathrm{Res}_\varepsilon A)^* \cong \widehat{\mathrm{B}}_\iota (\mathrm{B}_\pi A)^*~.
\]
\vspace{0.1pc}

Finally we get

\begin{align*}
\mathrm{Hom}_{u\mathcal{CC}_\infty\text{-}\mathsf{coalg}}((\mathrm{Res}_\varepsilon B)^\circ \otimes C_*^c(\Delta^\bullet),\widehat{\mathrm{B}}_\iota (\mathrm{B}_\pi A)^*) \\
\cong \mathrm{Hom}_{u\mathcal{CC}_\infty\text{-}\mathsf{coalg}}(C_*^c(\Delta^\bullet),\{(\mathrm{Res}_\varepsilon B)^\circ,\widehat{\mathrm{B}}_\iota (\mathrm{B}_\pi A)^*\}) \\
\cong \mathrm{Hom}_{u\mathcal{CC}_\infty\text{-}\mathsf{coalg}}(C_*^c(\Delta^\bullet),\widehat{\mathrm{B}}_\iota(\mathrm{hom}((\mathrm{Res}_\varepsilon B)^\circ,(\mathrm{B}_\pi A)^*)) \\
\cong \mathrm{Hom}_{\mathsf{curv}~\mathsf{abs}~\mathcal{L}_\infty\text{-}\mathsf{alg}}(\widehat{\Omega}_\iota(C_*^c(\Delta^\bullet)),\mathrm{hom}((\mathrm{Res}_\varepsilon B)^\circ,(\mathrm{B}_\pi A)^*) \\
\cong \mathcal{R}(\mathrm{hom}((\mathrm{Res}_\varepsilon B)^\circ, \mathfrak{g}_A))~,
\end{align*}

which concludes the proof.
\end{proof}

\begin{Example}
Let $A$ be a derived affine scheme, that is, a dg $u\mathcal{C}om$-algebra concentrated in non-negative degrees. There is an isomorphism of constant simplicial sets
\[
\mathrm{Spec}(A)(\kk) \cong \mathcal{MC}(\mathfrak{g}_A) \cong \mathcal{R}(\mathfrak{g}_A)~,
\]
since $\mathfrak{g}_A$ is concentrated in non-positive degrees. For $A$ a derived affine stack, we have that
\[
\mathrm{Spec}(A)(\kk) \simeq \mathcal{R}(\mathfrak{g}_A)~.
\]
Geometrically speaking, these isomorphisms allow us to recover all the $\kk$-points of $A$ as well as the formal neighborhood \textit{of any of these points} from its geometrical model $\mathfrak{g}_A$. 

\medskip

Furthermore, given a finite field extension $\mathbb{L}$ of $\kk$, we have  
\[
\mathrm{Spec}(A)(\mathbb{L}) \cong \mathcal{MC}(\mathrm{hom}(\mathbb{L}^*,\mathfrak{g}_A) \cong \mathcal{MC}(\mathbb{L} \otimes \mathfrak{g}_A)~,
\]
so we can base-change $\mathfrak{g}_A$ to any finite extension $\mathbb{L}$, and recover the $\mathbb{L}$-points of $A$ as well as their formal neighborhoods. Finally, any these combinations can be done for a finite number $n$ of points in $A$. 
\end{Example}

On the other hand, given a curved absolute $\mathcal{L}_\infty$-algebra $\mathfrak{g}$, one can construct a functor from dg Artinian algebras to simplicial sets. 

\begin{Definition}[Deformation functor]
Let $\mathfrak{g}$ be a curved absolute $\mathcal{L}_\infty$-algebra. Its \textit{deformation functor}
\[
\mathrm{Def}_{\mathfrak{g}}: \mathsf{dg}~\mathsf{Art}\text{-}\mathsf{alg}_{\geq 0} \longrightarrow \mathsf{sSet}
\]
is given by 
\[
\mathrm{Def}_{\mathfrak{g}}(B) \coloneqq \mathcal{R}(\mathrm{hom}(\mathrm{Res}_\varepsilon (-)^\circ, \mathfrak{g}))~.
\]
\end{Definition}

\begin{Remark}
This construction is analogous to the following classical construction: given a dg Lie algebra $\mathfrak{g}$ and a dg Artinian algebra $A$, one can consider the simplicial set given by
\[
\mathrm{MC}(\mathfrak{g} \otimes \overline{A})_\bullet~,
\]
where $\overline{A}$ is the augmentation ideal of $A$ and $\mathrm{MC}(-)_\bullet$ is the integration functor constructed in \cite{Hinich01}. This construction has been the classical way to associate to a dg Lie/$\mathcal{L}_\infty$-algebra a "deformation problem". 
\end{Remark}

\begin{Proposition}
Let $f: \mathfrak{h} \qi \mathfrak{g}$ be a weak-equivalence of complete curved absolute $\mathcal{L}_\infty$-algebras. It induces a natural weak-equivalence of simplicial sets $\mathrm{Def}_{f} :  \mathrm{Def}_{\mathfrak{h}} \qi \mathrm{Def}_{\mathfrak{g}}$ between the associated deformation functors.
\end{Proposition}

\begin{proof}
Let $f: \mathfrak{h} \qi \mathfrak{g}$ is a weak-equivalence of complete curved absolute $\mathcal{L}_\infty$-algebras. It induces a natural weak-equivalence of curved absolute $\mathcal{L}_\infty$-algebras $f_* : \mathrm{hom}(-,\mathfrak{h}) \qi \mathrm{hom}(-,\mathfrak{g})$. Indeed, the map $f$ is a weak-equivalence if and only if
\[
\widehat{\mathrm{B}}_\iota(f): \widehat{\mathrm{B}}_\iota \mathfrak{h}  \qi \widehat{\mathrm{B}}_\iota \mathfrak{g}
\]
is a quasi-isomorphism of $u\mathcal{CC}_\infty$-coalgebras. The following composition is a quasi-isomorphism
\[
\widehat{\mathrm{B}}_\iota(\mathrm{hom}(-,\mathfrak{h})) \cong \{-, \widehat{\mathrm{B}}_\iota \mathfrak{h} \} \qi \{-, \widehat{\mathrm{B}}_\iota \mathfrak{g} \}  \cong  \widehat{\mathrm{B}}_\iota(\mathrm{hom}(-,\mathfrak{g}))~,
\]
since $\{-,-\}$ is a right Quillen functor by Proposition \ref{prop: monoidal model category} and thus preserves weak-equivalences between fibrant objects, guaranteeing that $\{-,\widehat{\mathrm{B}}_\iota(f)\}$ is a quasi-isomorphism. One then concludes using the fact that $\mathcal{R}$ sends weak-equivalences to weak homotopy equivalences of simplicial sets by Theorem \ref{thm: propriétés de l'intégration}. 
\end{proof}

\begin{Remark}
The above proposition does not hold in general for the classical case of dg Lie/$\mathcal{L}_\infty$-algebras and their associated deformation functor when one considers weak-equivalences as quasi-isomorphisms. Indeed, the functor 
\[
\mathrm{Def}_{(-)}: \mathfrak{g} \mapsto \left[\mathrm{Def}_\mathfrak{g} \coloneqq \mathrm{MC}(\mathfrak{g} \otimes \overline{(-)})_\bullet: \mathsf{dg}~\mathsf{Art}\text{-}\mathsf{alg}^{\kk\text{-}\mathsf{aug}}_{\geq 0} \longrightarrow \mathsf{sSet} \right]
\]

does not always send quasi-isomorphisms to natural weak-equivalence of simplicial sets. Some hypothesis are need for this, such as nilpotency. It just so happens that nilpotent dg Lie algebras are examples of dg absolute Lie algebras. 
\end{Remark}

Let us explore other properties of the deformation functor.

\begin{lemma}\label{lemma: Def g preserves}
Let $\mathfrak{g}$ be a curved absolute $\mathcal{L}_\infty$-algebra. The functor 
\[
\mathcal{R}(\mathrm{hom}(-,\mathfrak{g})): \left(u\mathcal{CC}_\infty\text{-}\mathsf{coalg}\right)^{\mathsf{op}} \longrightarrow \mathsf{sSet}
\]
sends any homotopy colimit of $u\mathcal{CC}_\infty$-coalgebras to a homotopy limit of simplicial sets. Furthermore, it sends quasi-isomorphisms to weak-equivalences.
\end{lemma}

\begin{proof}
Let $\mathrm{hocolim}~ C_{\alpha}$ be a homotopy colimit of $u\mathcal{CC}_\infty$-coalgebras, we consider the following equivalences

\begin{align*}
\mathcal{R}(\mathrm{hom}(\mathrm{hocolim}~ C_{\alpha},\mathfrak{g})) &\simeq \mathrm{Hom}_{u\mathcal{CC}_\infty\text{-}\mathsf{coalg}}(C_*^c(\Delta^\bullet),  \{ \mathrm{hocolim}~ C_{\alpha},\widehat{\mathrm{B}}_\iota \mathfrak{g} \}) \\
&\simeq \mathrm{holim}~ \mathrm{Hom}_{u\mathcal{CC}_\infty\text{-}\mathsf{coalg}}(C_*^c(\Delta^\bullet),  \{C_{\alpha},\widehat{\mathrm{B}}_\iota \mathfrak{g} \}) \\
&\simeq \mathrm{holim}~ \mathcal{R}(\mathrm{hom}(C_{\alpha},\mathfrak{g}))~.
\end{align*}

\end{proof}

\begin{Proposition}
Let $\mathfrak{g}$ be a curved absolute $\mathcal{L}_\infty$-algebra. 

\medskip

\begin{enumerate}
\item The deformation functor $\mathrm{Def}_{\mathfrak{g}}$ sends quasi-isomorphisms to weak-equivalences of simplicial sets.

\medskip

\item The deformation functor $\mathrm{Def}_{\mathfrak{g}}$ preserves homotopy limits.
\end{enumerate}
\end{Proposition}

\begin{proof}
This follows directly from Lemma \ref{lemma: equivalence dg Artinien avec les uCC} and Lemma \ref{lemma: Def g preserves}.
\end{proof}

\textbf{Conclusion.} It is beyond the scope of this chapter to give a meaningful definition of what non-pointed formal moduli problems are or to prove that these are equivalent to curved absolute $\mathcal{L}_\infty$-algebras. However, we can \textit{point} at some reasonable hypothesis on functors 
\[
F:  \mathsf{dg}~\mathsf{Art}\text{-}\mathsf{alg}_{\geq 0} \longrightarrow \mathsf{sSet}~.
\]

\begin{enumerate}
\item Asking that $F(\kk) \simeq \{*\}$ does not seem to make any sense in this context. The set of Maurer-Cartan elements of $\mathfrak{g}$ could be empty since it is curved.

\medskip

\item The functor $F$ should respect at least some homotopy pullbacks. This is the key property that allows to construct second, third, etc. order deformations from the previous ones. 
\end{enumerate}

\medskip

Here is another point of view: in \cite{Hinich01}, V. Hinich proposed the seminal point of view that dg counital commutative coalgebras should be regarded as \textit{formal stack}. These objects would encode these formal spaces as their distribution algebras. It is not a coincidence that (pro) $\kk$-augmented dg Artinian algebras are equivalent to pointed conilpotent coalgebras, nor that the extension of deformation problems J. Maunder considered in \cite{Maunder} using curved Lie algebra has "pseudo-compact algebras" as its source. This category is in fact equivalent to all counital cocommutative coalgebras.

\medskip

From our point of view, Remark \ref{Rmk: pas de top model} makes pseudo-compact algebras/counital cocommutative coalgebras not the right sources as they do not have an appropriate homotopy theory. Therefore it seems natural to consider $u\mathcal{CC}_\infty$-coalgebras as the right homotopy/$\infty$-category  for these formal stacks/formal spaces. Thus these non-pointed formal moduli problems could be pre-sheaves on this category of formal stack

\[
F: \left(u\mathcal{CC}_\infty\text{-}\mathsf{coalg}_{\leq 0} \right) ^{\mathsf{op}} \longrightarrow \mathsf{sSet}~, 
\]
\vspace{0.1pc}

which preserve some kind of "geometrical homotopy pushouts", such as the functor $\mathrm{Def}_\mathfrak{g}$ itself. This will be the subject of a future work.

\section*{Appendix A: Filtered framework}
In this appendix, we recall the notion of a curved $\mathcal{L}_{\infty}$-algebra as well as the definition of their Maurer--Cartan elements. Since these elements are given as solutions of an equation involving infinite sums, the standard approach in order to deal these infinite sums has been to consider curved $\mathcal{L}_{\infty}$-algebras \textit{on an underlying complete graded module}. We compare this approach with our new approach in terms of curved \textit{absolute} $\mathcal{L}_{\infty}$-algebras. Again, we only consider the \textit{shifted} version of these structures, as the appear more naturally and involve simpler signs. This adjective is therefore implicit from now on.

\begin{Definition}[Curved $\mathcal{L}_{\infty}$-algebra]\label{def: classical curved linfty alg}
Let $\mathfrak{g}$ be a graded module. A \textit{curved} $\mathcal{L}_{\infty}$-\textit{algebra} structure $\{l_n\}_{n \geq 0}$ on $\mathfrak{g}$ is the data of a family of symmetric morphisms $l_n: \mathfrak{g}^{\odot n} \longrightarrow \mathfrak{g}$ of degree $-1$ such that, for all $n \geq 0$, the following equation holds:
\[
\sum_{p+q= n+1} \sum_{\sigma \in \mathrm{Sh}^{-1}(p,q)}(l_{p} \circ_1 l_q)^{\sigma} = 0~,
\]
where $\mathrm{Sh}^{-1}(p,q)$ denotes the inverse of the $(p,q)$-shuffles. The morphism $l_0 : \mathbb{K} \longrightarrow \mathfrak{g}$ is equivalent to an element $\vartheta$ in $\mathfrak{g}_{-1}$ which is called the \textit{curvature} of $\mathfrak{g}$. 
\end{Definition}

\begin{Remark}
Let $(\mathfrak{g},\{l_n\})$ be a curved $\mathcal{L}_\infty$-algebra. By definition

\[
l_1^2 = l_2 \circ_1 l_0~,
\]

whereas in a classical $\mathcal{L}_\infty$-algebra $l_1^2 = 0$ and thus defines a differential on $\mathfrak{g}$. Therefore there is no obvious notion of \textit{homology} or \textit{quasi-isomorphism} for curved $\mathcal{L}_\infty$-algebras.
\end{Remark}

\begin{Definition}[Maurer-Cartan element]\label{def: maurercartan classique}
Let $(\mathfrak{g},\{l_n\}_{n \geq 0})$ be a curved $\mathcal{L}_{\infty}$-algebra. A \textit{Maurer-Cartan element} of $\mathfrak{g}$ is an element $\alpha$ of degree $0$ that satisfies the following equation:

\begin{equation}\label{Maurer-Cartan}
\sum_{n\geq 0}\frac{l_n(\alpha,\cdots,\alpha)}{n!} = 0 ~.
\end{equation}

We denote by $\mathcal{MC}(\mathfrak{g})$ the set of Maurer-Cartan elements of $\mathfrak{g}$. 
\end{Definition}

\begin{Remark}
Contrary to the case of classical $\mathcal{L}_{\infty}$-algebras, the element $0$ no longer satisfies the Maurer-Cartan equation when the curvature is non-trivial. The set of Maurer-Cartan elements might therefore be empty. This allows them to encode \textit{non-pointed} deformation problems and it allows them to be rational models for \textit{non-pointed} spaces.
\end{Remark}

The Maurer-Cartan equation is \textit{not even defined} in general for curved $\mathcal{L}_{\infty}$-algebras, as it involves an infinite sum of terms in $\mathfrak{g}$. The same problem appears for classical $\mathcal{L}_{\infty}$-algebras. A first solution considered by various authors has been to impose some "nilpotency" condition on $\mathfrak{g}$, such that this sum becomes finite. Another has been to switch to impose an "underlying filtration" on the graded module $\mathfrak{g}$, compatible with the operations $\{l_n\}_{n \geq 0}$, such that this sum converges, see for instance \cite{DSV18} or Section \ref{Section: curved HTT} for an more detailed review of this framework. 

\begin{Definition}[Complete curved $\mathcal{L}_{\infty}$-algebra]\label{def: complete curved L infinity algebra}
A \textit{complete} curved $\mathcal{L}_{\infty}$-algebra $(\mathfrak{g},\{l_n\}_{n \geq 0}, \mathrm{F}_\bullet \mathfrak{g})$ amounts to the data of a curved $\mathcal{L}_{\infty}$-algebra structure $\{l_n\}_{n \geq 0}$ on a complete graded module $(\mathfrak{g}, \mathrm{F}_\bullet \mathfrak{g})$, such that the structural operations are compatible with the filtration of $\mathfrak{g}$: 
\[
l_n(\mathrm{F}_{i_1}\mathfrak{g} \odot \cdots \odot \mathrm{F}_{i_n}\mathfrak{g}) \subseteq \mathrm{F}_{i_1 + \cdots i_k + 1}\mathfrak{g}~.
\]
This condition is equivalent to asking that $l_n \in \mathrm{F}_1 \mathrm{hom}(\mathfrak{g}^{\odot n } , \mathfrak{g})$ for all $n \geq 0$. 
\end{Definition}

\begin{Remark}\label{rmk: MC equation only well-defined for elements in degree one of filtration}
In our definition, the structural operations $l_n$ raise the degree of the filtration by one. Some authors ask instead that the structural operations $l_n$ preserve the degree of the filtration. Notice that this difference is not fundamental in terms of making the Maurer-Cartan equation converge. Indeed, in both cases, the Maurer-Cartan equation \ref{Maurer-Cartan} is well-defined \textbf{only} for elements in $\mathrm{F}_1\mathfrak{g}$. 
\end{Remark}

Let us now try to compare this \textit{filtered framework} with the framework of curved \textit{absolute} $\mathcal{L}_\infty$-algebras. Recall that, by definition, a curved absolute $\mathcal{L}_\infty$-algebra structure on a pdg module $(\mathfrak{g},d_\mathfrak{g})$ is given by a structural morphism 
\[
\gamma_\mathfrak{g}: \prod_{n \geq 0} \widehat{\Omega}^{\mathrm{s.a}}\ucom^*(n) ~\widehat{\otimes}_{\mathbb{S}_n} ~ \mathfrak{g}^{\otimes n} \longrightarrow \mathfrak{g}~,
\]
satisfying Conditions \ref{pdg condition}, \ref{associativity condition} and \ref{curved condition}. In particular, for any infinite sum 
\[
\sum_{n\geq 0} \sum_{\omega \geq 1} \sum_{\tau \in \mathrm{CRT}_n^\omega} \lambda_\tau \tau(g_{i_1}, \cdots, g_{i_n})~,
\]
there is a well-defined image 
\[
\gamma_\mathfrak{g}\left(\sum_{n\geq 0} \sum_{\omega \geq 1} \sum_{\tau \in \mathrm{CRT}_n^\omega} \lambda_\tau \tau(g_{i_1}, \cdots, g_{i_n})\right)~,
\]
in the pdg module $\mathfrak{g}$. Thus the graded module $\mathfrak{g}$ admits infinite sums in "two direction": both \textit{arity-wise} infinite sums and \textit{weight-wise} infinite sums are well-defined. If, starting with a curved $\mathcal{L}_\infty$-algebra structure $\{l_n\}_{n \geq 0}$, one wants to "lift it" into a curved absolute $\mathcal{L}_\infty$-algebra structure, one therefore needs to endow $\mathfrak{g}$ with \textit{two} underlying filtrations that make these two types of sums converge. 

\begin{Definition}[Bifiltered graded module]
A \textit{bifiltered graded module} $(A,\mathrm{F}_{(\bullet,\bullet)} A)$ is the data of a graded module $A$ together with a degree-wise decreasing bifiltration:
\[
\begin{tikzcd}[column sep=0.5pc,row sep=0.5pc]
A \arrow[r,phantom,"\supseteq"] \arrow[d,phantom,"\supseteq" rotate=-90]
&\mathrm{F}_{(1,1)}A \arrow[r,phantom,"\supseteq"] \arrow[d,phantom,"\supseteq" rotate=-90]
&\mathrm{F}_{(2,1)}A \arrow[r,phantom,"\supseteq"] \arrow[d,phantom,"\supseteq" rotate=-90]
&\cdots \arrow[r,phantom,"\supseteq"] 
&\mathrm{F}_{(n,1)}A \arrow[r,phantom,"\supseteq"] \arrow[d,phantom,"\supseteq" rotate=-90]
&\cdots \\
\mathrm{F}_{(0,2)}A \arrow[r,phantom,"\supseteq"] \arrow[d,phantom,"\supseteq" rotate=-90] \arrow[d,phantom,"\supseteq" rotate=-90]
&\mathrm{F}_{(1,2)}A \arrow[r,phantom,"\supseteq"] \arrow[d,phantom,"\supseteq" rotate=-90]
&\mathrm{F}_{(2,2)}A \arrow[r,phantom,"\supseteq"] \arrow[d,phantom,"\supseteq" rotate=-90]
&\cdots \arrow[r,phantom,"\supseteq"] 
&\mathrm{F}_{(n,2)}A \arrow[r,phantom,"\supseteq"] \arrow[d,phantom,"\supseteq" rotate=-90]
&\cdots \\
\mathrm{F}_{(0,3)}A \arrow[r,phantom,"\supseteq"] \arrow[d,phantom,"\supseteq" rotate=-90]
&\mathrm{F}_{(1,3)}A \arrow[r,phantom,"\supseteq"] \arrow[d,phantom,"\supseteq" rotate=-90]
&\mathrm{F}_{(2,3)}A  \arrow[r,phantom,"\supseteq"] \arrow[d,phantom,"\supseteq" rotate=-90]
&\cdots \arrow[r,phantom,"\supseteq"] 
&\mathrm{F}_{(n,3)}A \arrow[r,phantom,"\supseteq"] 
&\cdots \\
\vdots \arrow[d,phantom,"\supseteq" rotate=-90]
&\vdots \arrow[d,phantom,"\supseteq" rotate=-90]
&\vdots \arrow[d,phantom,"\supseteq" rotate=-90]
&\ddots\\
\mathrm{F}_{(0,\omega)}A \arrow[d,phantom,"\supseteq" rotate=-90] \arrow[r,phantom,"\supseteq"]
&\mathrm{F}_{(1,\omega)}A \arrow[d,phantom,"\supseteq" rotate=-90] \arrow[r,phantom,"\supseteq"]
&\mathrm{F}_{(2,\omega)}A \arrow[d,phantom,"\supseteq" rotate=-90] 
&
&\mathrm{F}_{(n,\omega)}A \\
\vdots
&\vdots
&\vdots
& 
&
&\ddots
\end{tikzcd}
\]
for $n \geq 0$ and $\omega \geq 1$ such that $\mathrm{F}_{(0,1)} A = A$.
\end{Definition}

A bifiltered pdg module can be seen as a filtered filtered graded module. We call the filtration $\mathrm{F}_{(\bullet,\omega)}$ with respect to $n$ the \textit{arity filtration}, and the filtration $\mathrm{F}_{(n,\bullet)}$ with respect to $\omega$ the \textit{weight filtration}.

\begin{Definition}[Bicomplete graded module]
Let $(A,\mathrm{F}_{(\bullet,\bullet)}A)$ be a bifiltered graded module. It is a \textit{bicomplete} graded module if the canonical morphism 
\[
\pi_A: A \longrightarrow \lim_{(n,\omega) \geq (0,1)} A/\mathrm{F}_{(n,\omega)}A \cong \lim_{n \geq 0} \lim_{\omega \geq 1}A/\mathrm{F}_{(n,\omega)}A \cong \lim_{\omega \geq 1}\lim_{n \geq 0} A/\mathrm{F}_{(n,\omega)}A
\]
is an isomorphism of bifiltered graded modules. 
\end{Definition}

\begin{Definition}[Bicomplete curved $\mathcal{L}_\infty$-algebra]
A \textit{bicomplete curved} $\mathcal{L}_\infty$-\textit{algebra} $(\mathfrak{g},\{l_n\}_{n \geq 0},\mathrm{F}_{(\bullet,\bullet)})$ amounts to the data of a bicomplete graded module $(\mathfrak{g},\mathrm{F}_{(\bullet,\bullet)})$ together with a curved $\mathcal{L}_\infty$-algebra structure $\{l_n\}_{n \geq 0}$ on the graded module $\mathfrak{g}$ that is compatible with the bifiltration in the following way.
\begin{enumerate}
\item The operations $\{l_n\}_{n \geq 0}$ preserve the arity filtration
\[
l_n(\mathrm{F}_{(i_1 ,\omega)}\mathfrak{g} \odot \cdots \odot \mathrm{F}_{(i_n,\omega)}\mathfrak{g}) \subseteq \mathrm{F}_{(i_1 + \cdots + i_n,\omega)}\mathfrak{g}~.
\]
\item The operations $\{l_n\}_{n \geq 0}$ raise the weight filtration by one 
\[
l_n(\mathrm{F}_{(n,\omega_1)}\mathfrak{g} \odot \cdots \odot \mathrm{F}_{(n,\omega_n)}\mathfrak{g}) \subseteq \mathrm{F}_{(n,\omega_1 + \cdots + \omega_n + 1)}\mathfrak{g}~.
\]
\end{enumerate}
\end{Definition}

\begin{Remark}
Notice that for any $n \geq 0$, the graded module 
\[
\mathfrak{g}/\mathrm{F}_n\mathfrak{g} \coloneqq \lim_{\omega}\mathfrak{g}/\mathrm{F}_{(n,\omega)}\mathfrak{g}
\]
does not form a complete curved $\mathcal{L}_\infty$-algebra. While, on the other hand, for any $\omega \geq 1$, the graded module 
\[
\mathfrak{g}/\mathrm{F}_\omega \mathfrak{g} \coloneqq \lim_{n} \mathfrak{g}/\mathrm{F}_{(n,\omega)}\mathfrak{g}
\]
does form a complete curved $\mathcal{L}_\infty$-algebra. The weight filtration is \textit{operadic} in nature, while the arity filtration is \textit{ad hoc}, simply there to make certain sums converge.
\end{Remark}

\begin{Proposition}\label{prop: if it converges it converges}
There is a faithful forgetful functor 
\[
\mathrm{U}: \mathsf{curv}~\mathcal{L}_\infty\text{-}\mathsf{alg}^{\mathsf{bicomp}} \longrightarrow \mathsf{curv}~\mathsf{abs}~\mathcal{L}_\infty\text{-}\mathsf{alg}^{\mathsf{comp}}~,
\]
from the category of bicomplete curved $\mathcal{L}_\infty$-algebras to the category of complete curved absolute $\mathcal{L}_\infty$-algebras. 
\end{Proposition}

\begin{proof}
Let $\mathfrak{g}$ be a bicomplete curved $\mathcal{L}_\infty$-algebra. We set 
\[
\gamma_\mathfrak{g}\left(\sum_{n\geq 0} \sum_{\omega \geq \omega_0} \sum_{\tau \in \mathrm{CRT}_n^\omega} \lambda_\tau \tau(g_{i_1}, \cdots, g_{i_n})\right) \coloneqq \sum_{n\geq 0} \sum_{\omega \geq \omega_0} \sum_{\tau \in \mathrm{CRT}_n^\omega} \lambda_\tau \gamma_\mathfrak{g}(\tau(g_{i_1}, \cdots, g_{i_n})) ~,
\]
where $\gamma_\mathfrak{g}(\tau(g_{i_1}, \cdots, g_{i_n}))$ is given by replacing each vertex of $\tau$ by an operation $l_n$ in the following way: if $v$ is a vertex with $k$ incoming edges, then it is replaced by $l_k$. The infinite sum on the right-hand side is well-defined since we basically forced it to be so. The axioms of a curved $\mathcal{L}_\infty$-algebra ensure that $\gamma_\mathfrak{g}$ satisfies Conditions \ref{pdg condition}, \ref{associativity condition} and \ref{curved condition}, since all the sums are "split". 

\medskip 

Furthermore, let us show that the resulting curved absolute $\mathcal{L}_\infty$-algebra is complete. Observe that 

\[
\mathrm{W}_\omega \mathfrak{g} \subseteq \mathrm{F}_{(0,\omega)}\mathfrak{g}~,
\]

for all $\omega \geq 1$, therefore 

\[
\bigcap_{\omega \geq 0} \mathrm{W}_\omega \mathfrak{g} \subseteq \bigcap_{\omega \geq 0} \mathrm{F}_{(0,\omega)}\mathfrak{g} = \{0\}~,
\]

since $\mathfrak{g}$ is bicomplete. This implies that the canonical morphism 

\[
\varphi_\mathfrak{g}: \mathfrak{g} \twoheadrightarrow \lim_{\omega} \mathfrak{g}/\mathrm{W}_\omega \mathfrak{g}
\]

is an isomorphism, since it is always an epimorphism. 
\end{proof}

\begin{Remark}
The content of Proposition \ref{prop: if it converges it converges} could be subsumed as follows: if a curved $\mathcal{L}_\infty$-algebra admits a topology in which all infinite sums arity and weight-wise converge, then it can be thought as a curved absolute $\mathcal{L}_\infty$-algebra. 
\end{Remark}

\section*{Appendix B: Choosing a model}\label{Section: Appendix B}
In order to work with simpler algebraic structures, we recall the \textit{semi-augmented} Bar construction introduced in \cite[Section 3.3]{HirshMilles12}, which provides us with smaller but non-functorial cofibrant resolutions for non-augmented dg operads. We introduce its dual construction for dg counital partial cooperads in the sense of Chapter \ref{Chapter 1}.

\begin{Definition}[Semi-augmented dg operad]
A \textit{semi-augmented} dg operad $(\PP,\gamma,\eta,d_\PP,\epsilon)$ is the data of a dg operad $(\PP,\gamma,\eta,d_\PP)$ together with a degree $0$ morphism
\[
\epsilon: \PP \longrightarrow \I~,
\]
such that $\epsilon \circ \eta = \mathrm{id}_\I~.$ 
\end{Definition}

Let $\overline{\PP}$ be the graded $\mathbb{S}$-module given by $\overline{\PP} \coloneqq \mathrm{Ker}(\epsilon)$. 

\begin{Definition}[Semi-augmented Bar construction]
Let $(\PP,\gamma,\eta,d_\PP,\epsilon)$ be a semi-augmented dg operad. Its \textit{semi-augmented Bar construction}, denoted by $\mathrm{B}^{\mathrm{s.a}}\PP$, is given by
\[
\mathrm{B}^{\mathrm{s.a}}\PP \coloneqq \left(\mathscr{T}^c(s\overline{\PP}),d_{\mathrm{bar}} \coloneqq d_1 + d_2 ,\Theta_{\mathrm{bar}}\right)~.
\] 
Here $\mathscr{T}^c(s\overline{\PP})$ denotes the cofree conilpotent coaugmented cooperad generated by the suspension of the graded $\mathbb{S}$-module $\overline{\PP}$. The pre-differential $d_{\mathrm{bar}}$ is given by the sum of two terms. The first term $d_1$ is given by the unique coderivation extending the map
\[
\begin{tikzcd}[column sep=4pc,row sep=0.5pc]
\mathscr{T}^c(s\overline{\PP}) \arrow[r,twoheadrightarrow]
&s\overline{\PP} \arrow[r,"sd_{\overline{\PP}}"]
&s\overline{\PP}~.
\end{tikzcd}
\]
The second term $d_2$ is given by the unique coderivation extending the map
\[
\begin{tikzcd}[column sep=4pc,row sep=0.5pc]
\mathscr{T}^c(s\overline{\PP}) \arrow[r,twoheadrightarrow]
&s^2(\overline{\PP} \circ_{(1)} \overline{\PP}) \arrow[r,"s^{-1}\overline{\gamma}_{(1)}"]
&s\overline{\PP}~,
\end{tikzcd}
\]
where $\overline{\gamma}_{(1)}$ is given by the composition 
\[
\begin{tikzcd}[column sep=4pc,row sep=0.5pc]
\overline{\PP} \circ_{(1)} \overline{\PP} \arrow[r,"\gamma_{(1)}"]
&s\PP \arrow[r,twoheadrightarrow]
& s\overline{\PP}~.
\end{tikzcd}
\]
The curvature $\Theta_{\mathrm{bar}}$ is given by 
\[
\begin{tikzcd}[column sep=3.5pc,row sep=0.5pc]
\Theta_{\mathrm{bar}}: \mathscr{T}^c(\overline{\PP}) \arrow[r,twoheadrightarrow]
&s^2(\overline{\PP} \circ_{(1)} \overline{\PP}) \arrow[r,"s^{-2}\gamma_{(1)}"]
&s\PP \arrow[r,"\epsilon"]
&\I~.
\end{tikzcd}
\]
The semi-augmented Bar construction $\mathrm{B}^{\mathrm{s.a}}\PP$ forms a conilpotent curved cooperad.
\end{Definition}

It provides smaller but non-functorial resolutions. Indeed, it is only functorial with respect to morphisms of \textit{semi-augmented} dg operads.

\begin{theorem}[{\cite[Theorem 3.4.4]{HirshMilles12}}]
Let $(\PP,\gamma,\eta,d_\PP,\epsilon)$ be a semi-augmented dg operad. There is a quasi-isomorphism of dg operads
\[
\varepsilon_{\mathcal{P}}: \Omega \mathrm{B}^{\mathrm{s.a}}\PP \qi \PP~.
\]
\end{theorem}

Furthermore, it coincides in certain cases with the Boardmann-Vogt construction. 

\begin{theorem}\label{thm: Boardman-Vogt and cellular chains}
Let $\PP$ an operad in the category of cellular topological spaces. There is an isomorphism of dg operads
\[
\Omega \mathrm{B}^{\mathrm{s.a}} C_*^c(\PP,\kk) \cong C_*^c( \mathrm{W} \PP,\kk)~,
\]
\vspace{0.5pc}

where $C_*^c(-,\kk)$ denotes the cellular chain functor and where $\mathrm{W}(-)$ denotes the Boardmann-Vogt construction. 
\end{theorem}

\begin{proof}
This result is shown for reduced dg operads in \cite{BergerMoerdijk06}. See \cite{grignou2022mapping} for the extension to semi-augmented dg operads using this resolution. 
\end{proof}

\begin{Example}
Let $\mathrm{uCom}$ be the operad in the category of sets defined by 
\[
\mathrm{uCom}(n) \coloneqq \{*\}~,
\]
for all $n \geq 0$, together with the obvious action of $\mathbb{S}_n$ and the obvious operad structure. Then
\[
\Omega \mathrm{B}^{\mathrm{s.a}} \ucom \cong \Omega \mathrm{B}^{\mathrm{s.a}} C_*^c(\mathrm{uCom},\kk) \cong  C_*^c( \mathrm{W} \mathrm{uCom},\kk)~,
\]
where $\mathrm{uCom}$ is viewed as an operad in the category of cellular topological spaces by endowing it with the discrete topology. 
\end{Example}

\medskip

We introduce the dual version for semi-coaugmented dg counital partial cooperads. For the definition of complete curved absolute (partial) operads, we refer to Appendix \ref{Appendix B}.

\begin{Definition}[Semi-coaugmented dg counital partial cooperad]
A \textit{semi-coaugmented} dg counital partial cooperad $(\C,\{\Delta_i\},\epsilon,d_\C,\eta)$ amounts to the data of a dg counital partial cooperad $(\C,\{\Delta_i\},\epsilon,d_\C)$ together with a degree $0$ morphism
\[
\eta: \I \longrightarrow \C~,
\]
such that $\epsilon \circ \eta = \mathrm{id}_\I~.$ 
\end{Definition}

Let $\overline{\C}$ be the graded $\mathbb{S}$-module given by $\overline{\C} \coloneqq \mathrm{Ker}(\epsilon)$. 

\begin{Definition}[Semi-coaugmented complete Cobar construction]
Let $(\C,\{\Delta_i\},\epsilon,d_\C,\eta)$ be a semi-coaugmented dg counital partial cooperad. Its \textit{semi-coaugmented complete Cobar construction}, denoted by $\widehat{\Omega}^{\mathrm{s.c}}\C$, is given by
\[
\widehat{\Omega}^{\mathrm{s.c}}\C \coloneqq \left(\mathscr{T}^\wedge(s^{-1}\overline{\C}),d_{\mathrm{cobar}} \coloneqq d_1 -d_2 ,\Theta_{\mathrm{cobar}}\right)~.
\] 
Here $\mathscr{T}^\wedge(s^{-1}\overline{\C})$ denotes the completed tree monad applied to the desuspension of the graded $\mathbb{S}$-module $\overline{\C}$. The pre-differential $d_{\mathrm{cobar}}$ is given by the difference of two terms. The first term $d_1$ is given by the unique derivation extending the map
\[
\begin{tikzcd}[column sep=4pc,row sep=0.5pc]
s^{-1}\overline{\C} \arrow[r,"sd_{\overline{\C}}"]
&s^{-1}\overline{\C} \arrow[r,hookrightarrow]
&\mathscr{T}^\wedge(s^{-1}\overline{\C})~.
\end{tikzcd}
\]
The second term $d_2$ is given by the unique derivation extending the map
\[
\begin{tikzcd}[column sep=4pc,row sep=0.5pc]
s^{-1}\overline{\C} \arrow[r,"s\overline{\Delta}_{(1)}"]
&s^{-2}(\overline{\C} \circ_{(1)} \overline{\C}) \arrow[r,hookrightarrow]
&\mathscr{T}^\wedge(s^{-1}\overline{\C})
\end{tikzcd}
\]
where $\overline{\Delta}_{(1)}$ is given by the composition
\[
\begin{tikzcd}[column sep=4pc,row sep=0.5pc]
\overline{\C} \arrow[r,"\Delta_{(1)}"]
&\C \circ_{(1)} \C \arrow[r,twoheadrightarrow]
&\overline{\C} \circ_{(1)} \overline{\C}~.
\end{tikzcd}
\]
The curvature $\Theta_{\mathrm{cobar}}$ is given by 
\[
\begin{tikzcd}[column sep=3.5pc,row sep=0.5pc]
\Theta_{\mathrm{cobar}}: \I \arrow[r,"\epsilon"]
&\C \arrow[r,"s^{-2}\overline{\Delta}_{(1)}"]
&s^{-2}(\overline{\C} \circ_{(1)} \overline{\C}) \arrow[r,hookrightarrow]
&\mathscr{T}^\wedge(s^{-1}\overline{\C})~.
\end{tikzcd}
\]
The semi-coaugmented complete Cobar construction $\widehat{\Omega}^{\mathrm{s.c}}\C$ forms a complete curved augmented absolute operad. 
\end{Definition}

\begin{Proposition}\label{lemma: Bucom dual lin}
Let $(\PP,\gamma,\eta,d_\PP,\epsilon)$ be a semi-augmented dg operad which is arity and degree-wise finite dimensional. There is an isomorphism of complete curved augmented absolute operads
\[
\left(\mathrm{B}^{\mathrm{s.a}}\PP \right)^* \cong \widehat{\Omega}^{\mathrm{s.c}}\PP^*~.
\]
\end{Proposition}

\begin{proof}
We first view $\PP$ as a semi-augmented dg unital partial operad. Since it is arity and degree-wise finite dimensional, its linear dual is a semi-coaugmented dg counital partial cooperad. The result follows from a straightforward computation.
\end{proof}

Let us first make explicit what a curved algebra over $\widehat{\Omega}^{\mathrm{s.c}}\ucom$ is.
Recall that any complete curved augmented absolute operad structure gives a curved operad structure by applying a forgetful functor. 

\begin{Proposition}\label{Prop: Curved algebras over OmegauCom}
Let $(V,d_V)$ be a pdg module. The data of a curved $\widehat{\Omega}^{\mathrm{s.c}}\ucom$-algebra structure on $V$ is equivalent to a family of symmetric operations 
\[
\left\{ l_n: V^{\odot n} \longrightarrow V \right\}
\]
of degree $-1$ for $n \neq 1$. When setting $l_1 \coloneqq d_V$, the family $\{l_n\}_{n \geq 0}$ forms a curved $\mathcal{L}_\infty$-algebra structure in the sense of Definition \ref{def: classical curved linfty alg}.
\end{Proposition}

\begin{proof}
This follows from a straightforward computation.
\end{proof}

\begin{Proposition}
Let $(V,\mathrm{F}_\bullet, d_V)$ be a complete pdg module. The data of a morphism of complete curved absolute operads 
\[
\varphi: \widehat{\Omega}^{\mathrm{s.c}}\ucomd \longrightarrow \mathrm{End}^{(\geq 1)}_V
\]
is equivalent to a family of symmetric operations 
\[
\left\{l_n \quad \mathrm{in} \quad \mathrm{F}_1 \mathrm{hom}(V^{\odot n},V) \right\}
\]
of degree $-1$ for $n \neq 1$. When setting $l_1 \coloneqq d_V$, the family $\{l_n\}_{n \geq 0}$ forms a complete curved $\mathcal{L}_\infty$-algebra in the sense of Definition \ref{def: complete curved L infinity algebra}.
\end{Proposition}

\begin{proof}
See Section \ref{Section: curved HTT}.
\end{proof}

\begin{Remark}[Mixed curved $\mathcal{L}_\infty$-algebras] 
One may wonder what type of algebraic structure would have appeared had we chosen to use the Bar-Cobar constructions that appear in Chapter \ref{Chapter 2} instead of their "semi-(co)augmented" counterparts. The type of algebraic structure that appears is not substantially different. A \textit{mixed curved} $\mathcal{L}_\infty$-algebra $(\mathfrak{g},\{l_n\}_{n \geq 0},d_\mathfrak{g})$ amounts to the data of a pdg module $(\mathfrak{g},d_\mathfrak{g})$ together with a family $\{l_n\}_{n \geq 0}$ of symmetric morphisms $l_n: \mathfrak{g}^{\odot n} \longrightarrow \mathfrak{g}$ of degree $-1$ such that 
\[
\sum_{p+q= n+1} \sum_{\sigma \in \mathrm{Sh}^{-1}(p,q)}(l_{p} \circ_1 l_q)^{\sigma} = -\partial(l_n)~,
\]
for all $n \geq 0$. In particular, for $n=1$, the relation satisfied is
\[
d_\mathfrak{g}^2 - \partial(l_1) = l_1 \circ_1 l_1 + l_2 \circ_1 l_0~.
\]
It is immediate that $d_\mathfrak{g} - l_1$ together with the other structure maps does define a curved $\mathcal{L}_\infty$-algebra structure. Furthermore, this gives a forgetful functor
\[
\begin{tikzcd}[column sep=4.5pc,row sep=0.3pc]
\mathsf{mix} ~ \mathsf{curv}\text{-}s\mathcal{L}_\infty \text{-}\mathsf{alg} \arrow[r,"\mathrm{blend}"]
&\mathsf{curv}\text{-}\mathcal{L}_\infty \text{-}\mathsf{alg}~.
\end{tikzcd}
\]
This type of structures can be found in \cite{calaque2021lie}. Considering curved algebras over the curved cooperad $\mathrm{B}\ucom$ would result in "mixed curved absolute" $\mathcal{L}_\infty$-algebras.
\end{Remark}

\bibliographystyle{alpha}
\bibliography{bribs}

\end{document}